\documentclass[11pt]{amsart}

\usepackage[mathscr]{eucal}
\usepackage{amsmath,amssymb,amsfonts,amsthm,enumerate,mathtools, setspace}

\usepackage{color}

\newtheorem{theorem}{Theorem}[section]
\newtheorem{mtheorem}[theorem]{Main Theorem}
\newtheorem{lemma}[theorem]{Lemma}
\newtheorem{corollary}[theorem]{Corollary}
\newtheorem{proposition}[theorem]{Proposition}

\newtheorem*{definition}{Definition}

\newcommand{\Z}{\mathbb Z}
\newcommand{\R}{\mathbb R}
\newcommand{\Q}{\mathbb Q}

\DeclareMathOperator{\ord}{ord}
\DeclareMathOperator{\lcm}{lcm}

\DeclareMathOperator{\supp}{Supp}

\newcommand{\la}{\langle}
\newcommand{\ra}{\rangle}
\newcommand{\be}{\begin{equation}}
\newcommand{\ee}{\end{equation}}
\newcommand{\und}{\;\mbox{ and }\;}
\newcommand{\nn}{\nonumber}
\newcommand{\ber}{\begin{eqnarray}}
\newcommand{\eer}{\end{eqnarray}}

\newcommand{\Sum}[2]{\underset{#1}{\overset{#2}{\sum}}}
\newcommand{\Summ}[1]{\underset{#1}{\sum}}

\newcommand{\LK}{\,[\![}
\newcommand{\RK}{]\!]}

\newcommand{\wtilde}{\widetilde}

\newcommand{\D}{\mathsf D}

\newcommand{\Fc}{\mathcal F}
\newcommand{\vp}{\mathsf v}

\DeclareMathOperator{\im}{im}
\newcommand{\x}{\mathbf x}
\newcommand{\y}{\mathbf y}
\newcommand{\z}{\mathbf z}
\newcommand{\C}{\mathsf C}
\newcommand{\Co}{{{\mathsf C}^\circ}}
\DeclareMathOperator{\wt}{wt}
\newcommand{\darrow}{{\downarrow}}
\newcommand{\uarrow}{{\uparrow}}

\DeclareSymbolFont{goo}{OMS}{cmsy}{b}{n}
\DeclareMathSymbol{\gooT}{\mathalpha}{goo}{"1}
\newcommand{\bdot}{\mathbin{\gooT}}
\numberwithin{equation}{section}

\begin{document}

\title[The Characterization of  Finite Elasticities]{The Characterization of Finite Elasticities}
\author{David J. Grynkiewicz}
\email{diambri@hotmail.com}
\address{Department of Mathematical Sciences\\ University of Memphis\\ Memphis, TN 38152, USA}

\begin{abstract}
\begin{singlespacing}
Our main motivating goal is the study of factorization in Krull Domains $H$ with finitely generated class group $G$. While factorization into irreducibles, called atoms,  generally fails to be unique, there are various measures of how badly this can fail. One of the most important is the elasticity $\rho(H)=\lim_{k\rightarrow \infty}\rho_k(H)/k$, where $\rho_k(H)$ is the maximal number of atoms in any re-factorization of a product of $k$ atoms. Having finite elasticity is a key indicator that factorization, while not unique, is not completely wild. The elasticities, as well as many other arithmetic invariants, are the same as those of an associated  combinatorial monoid $\mathcal B(G_0)$ of zero-sum sequences, where $G_0\subseteq G$ are the classes containing height one primes.

We characterize when finite elasticity holds for any Krull Domain with finitely generated class group $G$. Indeed, our results are valid for the more general class of Transfer Krull Monoids (over a subset $G_0$ of a finitely generated abelian group $G$).  Moreover, we show there is a minimal $s\leq (d+1)m$, where $d$ is the torsion free rank and $m$ is the exponent of the torsion subgroup, such that $\rho_s(H)<\infty$ implies $\rho_k(H)<\infty$ for all $k\geq 1$. This ensures  $\rho(H)<\infty$ if and only if  $\rho_{(d+1)m}(H)<\infty$. Our characterization is in terms of a simple combinatorial obstruction to infinite elasticity:  there existing a subset $G_0^\diamond\subseteq G_0$ and global bound $N$ such that there are no nontrivial zero-sum sequences with terms from $G_0^\diamond$, and every minimal zero-sum sequence has at most $N$ terms from $G_0\setminus G_0^\diamond$. We give an explicit description of $G_0^\diamond$ in terms of the Convex Geometry of $G_0$  modulo the torsion subgroup $G_T\leq G$, and show that finite elasticity is equivalent to there being no positive $\mathbb R$-linear combination of the elements of this explicitly defined subset $G^\diamond_0$ equal to $0$ modulo $G_T$. Additionally, we use our results to show finite elasticity implies the set of distances $\Delta(H)$, the catenary degree $\mathsf c(H)$ (for Krull Monoids) and  a weakened form of the tame degree (for Krull Monoids)  are all also finite,  and that the Structure Theorem for Unions holds---four of the most commonly used measurements of structured factorization, after the elasticity.

Our results for factorization in Transfer Krull Monoids are accomplished by developing  an extensive theory in Convex Geometry generalizing positive bases. The convex cone generated by $X\subseteq \mathbb R^d$ is $\mathsf C(X)=\{\sum_{i=1}^{n}\alpha_ix_i:\; n\geq 0,\; x_i\in X,\,\alpha_i\in\mathbb R,\,\alpha_i\geq 0\}$. A positive basis for $\mathbb R^d$ is a minimal by inclusion subset $X\subseteq \mathbb R^d$ such that $\mathsf C(X)=\mathbb R^d$. Positive bases were first introduced and studied in the mid 20th century, and  the structural work initiated  by Reay led to a simplified proof and strengthening  of  Bonnice and Klee's celebrated generalization of Carath\'eordory's Theorem. They have since found increasing importance in areas of applied mathematics. We show that the structural result of Reay can be extended to special types of complete simplicial fans, which we term Reay systems. We extend these results to a general theory dealing with infinite sequences of Reay systems, as well as their limit structures.
The latter, while more complex, avoid the introduction of linear dependencies into the limit structure  that were not originally present, circumventing the general obstacle that a limit of linearly independent sets can degenerate into linear dependence. The resulting theory is used to  study a broad family of infinite subsets $G_0$ of a lattice $\Lambda\subseteq \mathbb R^d$ that exhibit  various types of finite-like behaviour, and which generalize the class of  subsets having  finite elasticity.

\end{singlespacing}
\end{abstract}

\maketitle

\section{Introduction}
This work comprises extensive developments combining the areas of Convex Geometry, the Combinatorics of infinite subsets of lattice points in $\R^d$, and the arithmetic of Krull Domains and,  more generally, of Transfer Krull Monoids. Our main and original motivation concerns the study of  factorization in the very general setting of Transfer Krull Monoids. This is  accomplished by developing an extended theory significantly generalizing prior work in Convex Geometry, then linked to the algebraic study of factorization via intermediary Combinatorics involving infinite subsets of lattice points in $\R^d$.

That the geometry of $\R^d$ could be used to  study similar algebraic  questions has precedent  in the study of Primitive Partition Identities \cite{Diaconis-Graham}, closely related to   Toric Varieties \cite{sturmfels-book} \cite{Buchstaber-Panov}, with such geometric ideas employed in alternative combinatorial language in
 \cite{G-precursor} \cite{Freeze-schmid} \cite{Ge-Yu13a}. However, we will be concerned with more subtly defined algebraic questions, for whom a direct connection with Convex Geometry is much less evident. Nonetheless, we will show that such a connection does exist, though it will require considerable preparatory work in Convex Geometry.

\subsection*{Convex Geometry:}
Linear Algebra over $\R^d$ is concerned with subspaces and linear combinations. Imposing the seemingly benign condition that only non-negative scalars be considered leads to the subfield of Convex Geometry concerned with Convex Cones $\C(X)=\{\Sum{i=1}{n}\alpha_ix_i:\;n\geq 0,\,x_i\in X,\,\alpha_i\in \R,\,\alpha_i\geq 0\}$ and positive linear combinations, where $X\subseteq \R^d$. This natural line of geometric inquiry was initiated and studied in the mid 20th century by various authors \cite{Davis}  \cite{McKinney} \cite{Bonnice-Klee}  \cite{Reay-memoir} \cite{Reay-gen}  \cite{Reay-art}  \cite{hansen-klee} \cite{Bonnice-reay} \cite{shephard}.
Part of their motivation lay in the relationship with Linear and Integer Programming, Game Theory and  classical combinatorial aspects of Convex Geometry  \cite{Reay-gen} \cite{50years}, particularly the widely studied subfield having the theorems of Helly, Rado, Steinitz and Carath\'eodory as core tenets (see the surveys \cite{Handbook-eckhoff} \cite{New-Carath-survey} and the hundreds of references there listed).
In this setting, the natural analog of a linear basis is a positive basis for $\R^d$, which is a minimal by inclusion subset $X\subseteq \R^d$ such that $\C(X)=\R^d$. Unlike linear bases, the cardinality of a positive basis is not determined by the dimension, instead satisfying the basic bounds $d+1\leq |X|\leq 2d$ for $d\geq 1$.
While their cardinality is not unique, Reay \cite{Reay-memoir} gave a basic structural result for positive bases. The structural description of Reay closely tied positive bases to the combinatorial  foundations of Convex Geometry by first giving a simplified proof \cite{Reay-gen} of Bonnice and Klee's common generalization of the theorems of Carath\'eordory and Steinitz \cite{Bonnice-Klee} \cite{Handbook-eckhoff}, and then significantly generalizing the Bonnice-Klee Theorem itself \cite{Bonnice-reay} \cite{Handbook-eckhoff}. Positive bases have since found renewed interest also in applied areas of mathematics, particularly in Derivative-Free Optimization   \cite{DFO-book} \cite{Regis}.

The first half of this work is devoted to an expansive extension of the basic structural theory initiated by Reay for positive bases. We begin by first extending the basic theory of positive bases to specialized complete simplicial fans.
A fan is a finite collection of polyhedral cones $\C(X)\subseteq \R^d$, so $X\subseteq \R^d$ is finite,  each containing no nontrivial subspace. The faces of $\C(X)$ are sub-cones $\C(Y)$ with $Y\subseteq X$ obtained by intersecting  $\C(X)$ with a hyperplane defining a closed  half-space that contains $\C(X)$, and it is required that each face of a cone from the fan be an element of the fan, and that  the intersection of any two
cones in  the fan be a face of each. The fan is complete if the union of all its cones is the entire space, and it is simplicial if each cone $\C(X)$ has $X$ linearly independent, in which case the  faces of $\C(X)$ correspond to its subsets $Y\subseteq X$. From a purely geometric perspective,
 complete simplicial fans are in bijective correspondence with starshaped spheres (specialized polyhedral complexes defined on the unit sphere), and  rational fans (those whose vertices $X$ come from  a lattice) are central to the definition of Toric Varieties \cite{Buchstaber-Panov}, which constitute an entire subfield of Algebraic Geometry. The aforementioned structural partitioning result  of Reay implicitly gives each positive basis the structure of a complete simplicial fan associated to the partitioning.
 In Section \ref{sec-reay}, we methodically extend all basic theory of positive bases, including the structural result of Reay and several additional properties crucial to the continued development of the theory,  more generally to a type of  complete simplicial fan which we term a Reay system.

 One of our main later goals is to study infinite subsets $G_0$ of a lattice in $\R^d$. We wish to better understand the radial ``directions'' in which the set $G_0$ escapes to infinity. For a simple example like $G_0=\{(-1,y):\; y\in \Z,\,y\geq 1\}\cup \{(x,0):\; x\in \Z,\,x\geq 1\}\cup \{(0,-1)\}$, it is intuitively clear that the positive $x$ and $y$ axes constitute the unbounded  ``directions'' of the  set $G_0$. However, more general subsets $G_0\subseteq \R^d$ can exhibit much more complicated behaviour. Even for $d=2$, we could replace the elements $(-1,y)$ in the previously defined $G_0\subseteq \Z^2$ with the elements $(-f(y),y)$, for some sub-linear, monotonically increasing and unbounded  function $f(y)$, and still obtain a set having the positive $x$ and $y$ axes as its unbounded ``directions'', yet with  a marked  2nd order unbounded drift in the negative $x$-axis direction  occurring as the set  escapes to infinity  along the positive $y$-axis.
 As even this simple example illustrates, to make such a vague notion of unbounded direction precise, we will need very careful definitions. The framework we adopt for making this precise is lain out in Section 3, where we define the notion of an asymptotically filtered sequence
  with limit $\vec u$.
  With this framework in place, we can then use asymptotically filtered sequences of terms from $G_0$, as well as their associated limits $\vec u$, as a means of studying the directions $\vec u$ for which the set  $G_0$ is unbounded.

Our much more expansive extension of the structural work of Reay involves showing the fundamental properties of Reay systems can be extended to a general notion of convergent families of Reay systems, which we term a virtual Reay system, whose component Reay systems are defined by collections of asymptotically filtered sequences. The resulting theory is presented in Section \ref{sec-virtual}. A key obstacle to accomplishing this is the fact that our convergent family of Reay systems will generally not converge to another Reay system. This is exemplary of the basic fact that a convergent family of linearly independent sets can have a linearly dependent set as limit. For instance, the linearly independent sets $X_i=\{(-1,0), (1,1/i)\}$ converge to the linearly dependent set $X=\{(-1,0), (0,1)\}$ as $i\rightarrow \infty$. We overcome this fundamental obstacle by extending the notion of Reay system to a specialized variation of
 polyhedral complex that has highly constrained boundary requirements for its faces despite allowing them to be neither fully open nor closed, as required for an ordinary polyhedral complex. We term this more general limit structure an Oriented Reay System, and show in Section \ref{sec-oriented} how all  basic properties and associated theory for  ordinary Reay systems carries over into the more complete setting of oriented Reay systems, which is then pre-requisite for the later work in the setting of virtual Reay systems.

 In total, Sections \ref{sec-asym-seq}--\ref{sec-virtual} constitute the mostly self-contained extension of the basic theory of Reay and others from the original context of positive bases to the more expansive setting of virtual and oriented Reay systems, developed to the point where we are able to  define and deal with the basic quantities essential to our later applications for Krull Domains and Transfer Krull Monoids.

\subsection*{Krull Domains, Transfer Krull Monoids and Factorization:} Krull Domains are one of the most ubiquitous classes of rings in Commutative Algebraic, being the higher dimensional analog of Dedekind Domains. Each Krull Domain $D$ has a class group $G$, and this work is primarily concerned with the case when $G$ is finitely generated. In this setting, it is well-known that unique factorization into primes corresponds to having a trivial class group, and there is  unique factorization of (divisorial) ideals into (height one) prime ideals \cite{fossum} \cite{bourbaki} \cite{bourbakie-eng-1-7} \cite{HK98}. When $G$ is nontrivial, unique factorization fails, and there may be multiple ways to factor a non-unit $a=u_1\cdot\ldots\cdot u_k$ into irreducibles $u_i$, called atoms, with $k$  the length of the factorization. In such case, the degree of wildness of factorization is measured by various arithmetic invariants, the most common of which we introduce in more general context momentarily. The most fundamental question, then, is whether these invariants are finite, as this is a key indication that factorization is not completely wild.

Worth noting,  most of the arithmetic invariants controlling factorization  depend only on the distribution of height-one prime ideals in the class group, meaning the subset $G_0 \subseteq G$ consisting of all  classes containing height-one prime ideals is of prime concern.  If $G_0$ is finite  or $D$ is a tame domain, then sets of potential factorization lengths are well-structured and all invariants describing their structure are finite \cite{Ge-Gr09b} \cite{Ge-Ka10a}. It is easy to see that $G_0$ generates $G$ as a semigroup, while  realization theorems tell us this is essentially the only restriction for which subsets $G_0$ can occur \cite[Theorems 2.5.4 and 3.7.8]{alfredbook} \cite{Ea-He73} \cite{Gi-He-Sm96}. To overview the most commonly studied arithmetic invariants for factorizations, we do so in the fairly broad class of  unit-cancellative semigroups.

Let $H$ be a multiplicatively written, commutative semigroup with identity element.
Then $H$ is said to be  cancellative if $ab=ac$ implies $b=c$, whenever $a,b,c\in H$, and unit cancellative if $ua=a$ implies $u$ is a unit, whenever  $u,a\in H$. If $D$ is a commutative domain,  then the semigroup of nonzero elements of $D$ is cancellative, and if $D$ is noetherian, then the monoids of all nonzero ideals and of all invertible ideals are unit-cancellative (with usual ideal multiplication) but not cancellative in general.
While it may not immediately come to mind, an important example of factorization occurs for semigroups of (isomorphism classes of) modules, in which case the semigroup operation is the direct-sum. In this setting, factorization of a module  corresponds to a direct-sum decomposition $M=M_1\oplus\ldots\oplus M_k$ into a finite number of irreducible submodules. Many  classical unique decomposition theorems for modules are then equivalent to  factorization in the associated semigroup of modules being unique. Here, unit-cancellativity means  $M \cong M \oplus N$ implies $N=0$.
Thus unit-cancellativity states that all modules have to be directly finite or, in other words,  Dedekind-finite. This is a  frequent property (e.g., valid for all finitely generated modules over commutative rings), that is weaker than cancellativity  \cite{Go79a} \cite{Ha-Gu-Ki04a}.

Among the oldest and most important arithmetic  invariants of factorization are the elasticities. They were first studied (using alternative terminology) for rings of integers in algebraic fields (see \cite{Ge-Le90} for references to the old literature), with the now standard term elasticity introduced by Valenza \cite{elast1-valenza} in 1990. They have since been  studied by numerous authors. To list a few examples (focussing on surveys and more recent papers), see  \cite{An-An92} \cite{An-An94} \cite{elast999-anderson} \cite{An97a}
\cite{elast8-baeth}
\cite{elast2-baron} \cite{Ba-Co16a}
\cite{Ca-Ch95} \cite{Ch-Cl05}
\cite{elast99-chapman}  \cite{elast993-chapman} \cite{elast995-chapman} \cite{Ge-Zh18a}
\cite{elast4-fuj} \cite{elast5-gao}
\cite{elast7-garcia}
\cite{Gerold-lambert-rankone} \cite{alfredbook}  \cite{Ge-Ka10a} \cite{Go-ON20} \cite{Ha02b}
 \cite{elast996-halt} \cite{elast997-halt} \cite{Ka05a} \cite{Ka16b}
 \cite{elast2-oneil} \cite{Sc16a} \cite{Zh19a} \cite{Zh20a}. It can be defined
  as  $\rho(H)=\sup\{\sup \mathsf L(a)/\min \mathsf L(a):\;a\in H \mbox{ a non-unit}\}$, where $\mathsf L(a)=\{k\in \mathbb N:\; a=u_1\cdot \ldots\cdot u_k\mbox{ from some atoms $u_i\in H$}\}$ denotes the length set of $a$, that is, the set of all possible factorization lengths of $a$ written as a product of atoms, or it can be defined equivalently \cite[Proposition 1.4.2]{alfredbook} as $\rho(H)=\lim_{k\rightarrow \infty}\rho_k(H)/k$, where the $k$-th elasticity $\rho_k(H)$ denotes the maximum number of atoms in any re-factorization of a product of $k$-atoms.
The set of distances $\Delta(H)$, defined as the minimum difference of two consecutive elements of $\mathsf L(a)$ as we range over all non-units $a\in H$, and the Catenary degree $\mathsf c(H)$ (see Section \ref{sec-intro-factorization}) are two other of the most oft studied arithmetic invariants \cite{alfredbook}.
Yet another  measure of well-behaved factorization are  structural results for $\mathcal U_k(H)$, which is the set of all $\ell$ for which there are atoms $u_1,\ldots,u_k,v_1,\ldots,v_\ell\in H$ with $u_1\cdot \ldots\cdot u_k=v_1\cdot \ldots\cdot v_\ell$, so the possible re-factorization lengths of some product of $k$ atoms. For many semigroups, it is known that these sets must be highly structured in the following sense.
A finite set $X\subseteq \Z$ is said to be an almost arithmetic progression with difference $d\geq 1$ and bound $N\geq 0$ if $X=P\setminus Y$, where $P$ is an arithmetic progression with difference $d$ and $Y\subseteq P$ is a subset contained in the union of the first $N$ terms from $P$ and the last $N$ terms from $P$.  If there exists a constant $N\geq 0$ and difference $d\geq 1$ such that $\mathcal U_k(H)$ is an almost arithmetic progression with difference $d$ and bound $N$ for all sufficiently large $k$, then  $H$ is said to satisfy the Structure Theorem for Unions.
As shown in \cite[Theorem 4.2]{Gao-Ger-phok}, if $\Delta(H)$ is finite and there is a constant $M\geq 0$ such that $\rho_{k+1}(G_0)-\rho_k(G_0)\leq M$ for all $k\geq 1$, then the Structure Theorem for Unions holds for $H$, meaning this additional structure is implied by sufficiently strong finiteness results for the elasticities and  set of distances.

\medskip

The general class of semigroup treated in this paper are Transfer Krull Monoids. We remark that they include  all Krull Domains, Krull Monoids, and many examples of natural
semigroups of monoids, with a non-cancellative monoid of modules over Bass
rings that is Transfer Krull studied in \cite{Ba-Sm21}. We defer the formal definitions to Section \ref{sec-intro-factorization}, but continue with a detailed list illustrating the broad extent of the class of semigroups covered by our results (see also \cite[pp 972]{Ge16c} and \cite[Example 4.2]{Ge-Zh20a}).

\emph{Commutative Domains.} A Noetherian domain $D$ is Krull if and only if it is integrally closed, and  the integral closure of a Noetherian domain is Krull by the  Mori-Nagata Theorem. If $D$ is Krull and finitely generated over $\Z$, then its class group $G$ is finitely generated \cite[Chapter 2, Corollary 7.7]{La83}. Moreover,  $D$ is a Krull Domain if and only if its multiplicative monoid of nonzero elements is a Krull Monoid.

 \emph{Submonoids of Commutative Domains.} Regular congruence monoids defined in Krull Domains are Krull Monoids \cite[Section 2.11]{alfredbook}. Let $D$ be a factorial domain with quotient field $K$. Then the ring of integer-valued polynomials $\mathsf{Int} (D) = \{ f \in K[X] :\; f (D) \subseteq D\}$ is not Krull, but the divisor closed submonoid $\LK f \RK$ is a Krull Monoid for any polynomial $f \in \mathsf{Int} (D)$ \cite{Re14a, Fr16a}. Moreover, if $f \in D[X]$, then the class group of $\LK f \RK$ is finitely generated by \cite[Theorem 4.1]{Re17a}.

 \emph{Finitely Generated Krull Monoids.} A finitely generated monoid is a Krull Monoid if and only if it is root closed, and  finitely generated Krull Monoids have finitely generated class groups \cite[Theorem 2.7.14]{alfredbook}. The rank of the class group is studied in geometric terms in \cite[Corollary 1]{Le88}.  Two interesting special cases are the following.

Normal Affine Monoids. These are found in  combinatorial Commutative Algebra, and are equivalent to reduced finitely generated Krull Monoids with torsion free quotient group. For the class group  of a finite normal monoid algebra $D[M]$, where $D$ is a Noetherian Krull Domain and $M$ a normal affine monoid, its class group is the direct sum of the class group of $D$ and that of $M$ \cite[Theorem 4.60]{Br-Gu09a}, meaning it is finitely generated whenever the class group of $D$ is finitely generated.

Diophantine Monoids. A Diophantine monoid is an additive monoid of non-negative solutions to  a system of linear Diophantine equations, with the rank of the class group, in terms of the defining matrix, studied in \cite{Ch-Kr-Oe02}.

\emph{Semigroups of Modules:} A semigroup $\mathcal V$ of modules over a ring $R$ (as described earlier) closed under finite direct sums, direct summands, and isomorphisms is a reduced commutative semigroup. If the endomorphism rings $\mathsf{End}_R (M)$ are semilocal for all modules $M$ in the semigroup, then  $\mathcal V$ will be a Krull Monoid \cite[Theorem 3.4]{Fa02}. Lists of modules having this property may be found in
  \cite{Fa06b},  and examples when the class group is finitely generated are listed
  in \cite{Ba-Ge14b}. Conversely, every reduced Krull Monoid is isomorphic to a monoid of modules \cite[Theorem 2.1]{Fa-Wi04} with details available in the monograph \cite{Fa19a}.

\emph{Normalizing Krull Monoids.} A  semigroup $H$ is  normalizing if $aH=Ha$ for all $a \in H$. Normalizing Krull Monoids occur when studying  Noetherian semigroup algebras \cite{Je-Ok07a} and are Transfer Krull by \cite[Theorems 4.13 an 6.5]{Ge13a}.

\emph{Noncommutative rings.} Any bounded hereditary Noetherian prime ring $D$ for which every stably free left $D$-ideal is free is a Transfer Krull Domain \cite[Theorem 4.4]{Sm19a}.  Results giving other noncommutative families of Transfer Krull Monoids may be found in  \cite{Sm13a, Ba-Ba-Go14, Ba-Sm15, Sm16a}.

\emph{Commutative domains  close to a Krull Domain.} In many cases, if a  domain $R$ is ``close'' to a Krull Domain $D$, then $R$ will be  a Transfer Krull Domain. Results giving examples of this type may be found in \cite[Section 5]{Ge-Zh20a} and \cite[Proposition 4.6 and Theorem 5.8]{Ge-Ka-Re15a}. 

\subsection*{Zero-Sum Sequences:}
The topics discussed above under the headings for Convex Geometry and for Transfer Krull Monoids may seem so extremely disparate as to be completely unrelated. It is one of our main goals to show that the opposite is true, with there in fact being an extremely close connection  between these seemingly unrelated topics. In order to establish this connection, we will need an intermediary combinatorial structure involving subsets $G_0$ of a lattice $\Lambda\subseteq\R^d$.

Let $G$ be an abelian group and let $G_0\subseteq G$ be a subset. Following the tradition of Combinatorial Number Theory, a sequence over $G_0$ is a finite, unordered string of elements from $G_0$. The collection $\mathcal B(G_0)$ consisting of all sequences over $G_0$ whose terms sum to zero, equipped with the concatenation operation, makes $\mathcal B(G_0)$ into a monoid which is Krull. More importantly, every Krull Monoid $M$ has a transfer homomorphism to a monoid $\mathcal B(G_0)$ of zero-sum sequences that directly translates nearly all factorization properties of the original monoid $M$ into the corresponding ones for $\mathcal B(G_0)$. We go into more detail in Section \ref{sec-intro-factorization}. For our purposes, this means all factorization questions considered in this paper for a general Krull Monoid reduce to the study of the combinatorial object $\mathcal B(G_0)$ by established machinery, allowing us to focus solely on $\mathcal B(G_0)$. A similar statement is true for  Transfer Krull Monoids, though the increased broadness of this class means the corresponding transfer homomorphism  between $M$ and $\mathcal B(G_0)$ is necessarily  weaker. For our purposes, this means that, while  all  our main  results regarding  factorization will be valid in the highly general context of Transfer Krull Monoids, a small number of consequences will only be valid for  Krull Monoids. Thus  we  work entirely with the combinatorial object $\mathcal B(G_0)$ in this work with the corresponding results transferred automatically to Transfer Krull Monoids or Krull Monoids by established machinery (referenced in Section \ref{sec-intro-factorization}) that we need only use in passing. Note, arithmetic invariants for $\mathcal B(G_0)$ are generally abbreviated as $\rho(G_0):=\rho(\mathcal B(G_0))$, etc.

While the connection between factorization in Transfer Krull Monoids and the combinatorial object $\mathcal B(G_0)$ is well-known, any further connections with Convex Geometry were much more limited. The behavior of factorization when $G$ is not finitely generated is generally quite unrestricted. As such, the study of $\mathcal B(G_0)$ for  subsets $G_0$ of a finitely generated abelian group $G$ is the
the general framework for the foundational finiteness questions for the arithmetic invariants of Transfer Krull Monoids. In this regard, little in the way of characterization was known apart from the basic case when $G_0$ is finite or the case when $G$ has rank one \cite{elast999-anderson} \cite{Gerold-lambert-rankone}.
Since all factorization invariants are  finite when $G$ is finite (by basic arguments), it is natural to expect the chief obstacles for characterizing their finiteness would already be present for torsion-free abelian groups. As it will turn out (though not initially clear), this is precisely the case, allowing us to focus almost exclusively on the torsion-free case, later adapting these argument applied to $G/G_T$, where $G_T$ is the torsion subgroup. As such, we are reduced to considering infinite subsets $G_0\subseteq \Z^d$. The abelian group $\Z^d$ is the prototypical lattice in $\R^d$, meaning it is a full rank discrete subgroup of $\R^d$. While any full rank lattice $\Lambda \subseteq \R^d$ is isomorphic to $\Z^d$, it will be convenient to expand consideration to subsets $G_0$ of a general  full rank lattice $\Lambda\subseteq \R^d$, allowing the use of the geometry of $\R^d$ for studying $G_0$.

\subsection*{Main Results:}
Our driving goal is to characterize finite elasticity for $\mathcal B(G_0)$ (and thus for Transfer Krull Monoids in general) for any subset $G_0$ of a finitely generated abelian group $G$. One natural way to prevent infinite elasticity is if there is a subset $G_0^\diamond \subseteq G_0$
and global bound $N$ such that there are no nontrivial zero-sum sequences with terms from $G_0^\diamond$, and every minimal zero-sum sequence has at most $N$ terms from $G_0\setminus G_0^\diamond$. Such a condition trivially implies the elasticity  $\rho(G_0)$ is finite, in turn implying all refined elasticities $\rho_k(G_0)$ are also finite (see Proposition \ref{prop-pre-rho-char}). In the spirit of results like Hall's Matching Theorem, which shows that the simple combinatorial obstruction to a perfect matching characterizes when they exist,  we will show that this basic combinatorial obstruction characterizes finite elasticity. Indeed, it characterizes when $\rho_{(d+1)m}(G_0)$ is finite, where $m$ is the exponent of the torsion subgroup of $G$ and $d$ is the torsion-free rank, meaning there is a minimal $s\leq (d+1)m$ such that $\rho_s(G_0)<\infty$ implies $\rho_k(G_0)<\infty$ for all $k$. Moreover,
we give an explicit description of $G_0^\diamond$ in terms of the Convex Geometry of $G_0$   modulo the torsion subgroup $G_T\leq G$, and show that finite elasticity is equivalent to there being no  positive  $\R$-linear combination of the elements of this explicitly defined subset $G^\diamond_0$ equal to $0$ modulo $G_T$.
What this means is that the initial question of finite elasticity, involving equations of positive  $\Z$-linear combinations of lattice points, is  equivalent to one involving positive  $\R$-linear combinations of lattice points. Additionally, we obtain  a weak  structural description of the atoms in $\mathcal B(G_0)$, assuming finite elasticities, and use this to show  finite elasticity implies that the set of distances  $\Delta(G_0)$, the catenary degree $\mathsf c(G_0)$ and  a weakened form of the tame degree  are all also finite, and that  there are no arbitrarily large gaps in the sequence $\{\rho_k(G_0)\}_{k=1}^\infty$, which implies that the Structure Theorem for Unions  also holds.

The key means of accomplishing our characterization and many of its consequences  is by an in-depth study in Section \ref{sec-finitary} of an ample class of subsets $G_0\subseteq \Lambda\subseteq \R^d$, defined using our generalized theory of Reay systems given in Sections \ref{sec-asym-seq}--\ref{sec-virtual}, that possess several finite-like properties despite being (in general) infinite. We term such sets finitary, study their finite characteristics in detail, and show that finite elasticity implies $G_0$ is finitary (with the converse failing). The broader class of finitary sets, as it will turn out,
 shares most of the same structural  properties as sets $G_0$ with finite elasticities,
 while at the same time behaving better with respect to inductive arguments using quotients.

 We conclude the introduction by outlining the main results and content section by section. As a general remark, while somewhat subjective, we have labeled results as lemmas when they are more technical, often  with highly restricted hypotheses, and generally needed as part of a larger proof, as  propositions when they encode basic or fundamental properties of the concepts being explored, even when the proof may be quite involved, and as theorems when we wish to emphasize that it is one of our culminating results. Theorems that we wish to especially highlight will be referred to as main theorems.

 \begin{itemize}
\item[Section \ref{sec-prelim}:] We introduce the basic notation and preliminaries for Convex Geometry, partially ordered sets (Posets), lattices, zero-sum sequences, and Transfer Krull Monoids, as well as the requisite asymptotic notation. A more general notion of rational sequence over $G_0$ is introduced, where the multiplicities of terms are allowed to be non-negative rational numbers rather than integers, which will play a crucial role in later parts of Section \ref{sec-fact}. Precursors to this innovation may be found in \cite{G-precursor}.

\item[Section \ref{sec-asym-seq}.1:] The concept of an asymptotically filtered sequence with limit $\vec u$ and related notation and definitions are  introduced along with the basic properties of these definitions. This will be our main tool for measuring the directions in which a set escapes to infinity.
\item[Section \ref{sec-asym-seq}.2:]
    The companion concepts of encasement and minimal encasement of $\vec u$, as well as
    several basic properties, are given in the context of finite unions of polyhedral cones.
\item[Section \ref{sec-asym-seq}.3:]
    The notion of  a set $X$ being bound to another set $Y$, meaning every point in $X$ is within some globally bounded distance of some point from $Y$, is introduced. The section culminates with Theorem \ref{thm-nearness-characterization}, which gives a characterization of $X$ being bound to $Y$, assuming $Y$ is a finite union of polyhedral cones, in terms of encasement of limits $\vec u$ of asymptotically filtered sequences. We will mostly need the results of Sections \ref{sec-asym-seq}.2 and \ref{sec-asym-seq}.3  when $Y$ is a polyhedral cone.

\item[Section \ref{sec-reay}.1:] We characterize some basic non-degeneracy assumptions for arithmetic properties of $G_0$ in terms of the geometry of $\R^d$.
\item[Section \ref{sec-reay}.2:] The concept of an elementary atom, which is a slight modification to the definition as presented in \cite{G-precursor}, is given. It corresponds to the notion of atom when factoring using rational sequences rather than ordinary sequences, so using multiplicities from $\Q$ rather than $\Z$. A version of Carath\'eordory's Theorem, formulated for rational sequences, is  given in Theorem \ref{thm-carahtheodory-elm-atom}. Positive bases are introduced, with a lengthy list of equivalent defining conditions given in Proposition \ref{prop-char-minimal-pos-basis}, including equivalent formulations involving elementary atoms.
  \item[Section \ref{sec-reay}.3:]  The notion of a Reay system and related concepts are  introduced, and the basic theory of positive bases extended to this context. Proposition \ref{prop-reays-structureresult} corresponds to a refined statement of the structural result of Reay for positive bases, extended to Reay systems. The subsection concludes with Proposition \ref{prop-FanStability}, which contains the technical stability properties of complete simplicial fans, needed for handling some delicate arguments in Proposition \ref{prop-VReay-SupportSet}.2 and Lemma \ref{Lemma-VReay-RidRemainders}.2, and thus in turn for Theorem \ref{thm-keylemmaII} as well.

\item[Section \ref{sec-oriented}:] The notion of Reay system is generalized to that of oriented Reay systems, with all related notation introduced. The theory and results for Reay systems presented in Section \ref{sec-reay}.3 are extended to this more general context.

\item[Section \ref{sec-virtual}:] The notion of oriented Reay system is now extended to that of virtual Reay systems  with all related notation introduced. The theory and results for ordinary and oriented  Reay systems presented in Sections \ref{sec-reay}.3 and \ref{sec-oriented} are now extended to this more general context.

\item[Section \ref{sec-finitary}.1:] Three equivalent definitions for the subset $G_0^\diamond \subseteq G_0$, which  plays the  crucial role in the characterization of finite elasticities, are given in Proposition \ref{prop-G_0diamond-1st-easy-equiv}. The crucial notion of a finitary set is defined using the language of Sections \ref{sec-asym-seq}--\ref{sec-virtual}. The remainder of the section is devoted to an in-depth study of finitary sets. Theorem \ref{thm-finitary-diamond-containment} gives a condition on $G_0^\diamond$ that implies $G_0$ is finitary. In particular, $0\notin \C^*(G_0^\diamond)$ implies $G_0$ is finitary. Combining this result with Proposition \ref{prop-pre-rho-char} gives us the means to use finitary sets, along with all their associated properties developed in Section \ref{sec-finitary}, during  our characterization of finite elasticities. Theorem \ref{thm-neg-char} gives a 4th characterization of the set $G_0^\diamond$ in terms of the multiplicities of elements in elementary atoms, valid for finitary sets. This characterizes the elements of $G_0^\diamond$ in terms of positive $\Q$-linear primitive partition identities equal to zero. The subsection concludes with Theorem \ref{thm-keylemmaII}, giving geometric restrictions for the distribution of elements in  a finitary set $G_0$. In particular, it ensures that a finitary set $G_0$ has a linearly independent subset $X\subseteq G_0^\diamond\subseteq G_0$ such that $G_0$ is bound to $-\C(X)$, meaning the set $G_0$ must be concentrated around the simplicial cone $-\C(X)$.

\item[Section \ref{sec-finitary}.2:]  The goal of this subsection and the next two is to derive more structural information regarding   certain virtual Reay systems defined over a finitary set. To do this, the notion of a series decomposition  of a purely virtual Reay system is introduced (inspired by the Jordan-H\"older Theorem and other similar results in Algebra), along with detailed  notation and  concepts.

 \item[Section \ref{sec-finitary}.3:]    Series decompositions are used to study the structure of the sets $X(G_0)$ and $\mathfrak X(G_0)$ associated to a finitary set in Section \ref{sec-finitary}.2.
     Four of the key finiteness properties of finitary sets are given in Theorems \ref{thm-finitary-FiniteProps-I}.1, \ref{thm-finitary-FiniteProps-I}.2, \ref{thm-finitary-FiniteProps-II} and \ref{thm-finitary-FiniteProps-III}.
     Included is  the striking property given in Proposition \ref{prop-finitary-FiniteDeletion} that a finitary set $G_0$ has a finite subset $Y\subseteq G_0$ such that every zero-sum sequence with terms from $G_0$ must use at least one term from $Y$.

\item[Section \ref{sec-finitary}.4:]
The notions of lattice types and minimal types, as well as  related notation and concepts, are introduced. Their crucial finiteness and  interchangeability properties, needed for many results in Section \ref{sec-fact}, are given in sole result of the subsection: Proposition  \ref{prop-finitary-mintype}.

\item[Section \ref{sec-fact}.1:] The core of this subsection is Theorem \ref{thm-structural-Lambert}. Buried within the statement of Theorem \ref{thm-structural-Lambert} is a multi-dimensional generalization of a result of Lambert \cite{Lambert} \cite[Lemma 4.3]{Gerold-lambert-rankone}. The result of Lambert handles the case of dimension $d=1$, where a novel variation for the argument giving the upper bound $\mathsf D(G)\leq |G|$ for the Davenport constant yields the result. Though, in the end, the proof of Theorem \ref{thm-structural-Lambert} was derived as a natural extension of the theory developed in Sections \ref{sec-asym-seq}--\ref{sec-finitary} rather than as a generalization of the argument of Lambert,  this innocuous one-dimensional result nonetheless served as the inspiration for a good portion of the material contained in this work. Theorem \ref{thm-structural-Lambert} plays the central role in our characterization of finite elasticities, and its  proof involves a detailed algorithm utilizing the theory developed for finitary sets in Section \ref{sec-finitary}. The beginning of the section contains several basic results involving elasticity and $G_0^\diamond$. Corollary \ref{cor-G_0diamond-Zequiv} contains a 5th  characterization of the set $G_0^\diamond$ in terms of the multiplicities of elements in ordinary atoms, valid for finitary sets. This characterizes the elements of $G_0^\diamond$ in terms of positive $\Z$-linear primitive partition identities equal to zero. Theorem \ref{thm-rho-char} is our characterization result in the  torsion-free setting. Corollary \ref{cor-rhoelem} shows that finite elasticity is equivalent to a variation on elasticity defined using the more basic elementary atoms rather than ordinary atoms.

\item[Section \ref{sec-fact}.2:] This subsection contains all the additional consequences for factorization when $G$ is torsion-free. Theorem \ref{thm-structural-char} gives a weak structural characterization of the atoms over $G_0$ under the assumption of finite elasticities. Its principal use is to simulate globally bounded support for an atom. The number of distinct elements appearing in an atom over $G_0$ can generally be an arbitrarily large number. Omitting some details, Theorem \ref{thm-structural-char} shows that there are a finite number of ``types'' of elements, with elements of the same type behaving in the same, well-controlled manner as given by Proposition \ref{prop-swapping}, so that if one equates all elements of the same type in an atom (with some elements belonging to no type and left alone), then the total number of distinct  elements/types in an atom is globally bounded. While $\mathcal B(G_0)$ may not be tame in general, we introduce a weaker notion of tameness in Theorem \ref{thm-tame-pseudo} which is finite (assuming the elasticities are finite). While this notion of tameness is weaker than what is generally studied in the literature, it is nonetheless sufficient, under an assumption of finite elasticities, to showing the set of distances  and catenary degree are finite, as well as that the Structure Theorem for Unions holds. This is done in Theorem \ref{thm-Delta-finite}. At the end of Section \ref{sec-fact}.2, a more detailed summary of all the information derived for $G_0$ assuming finite elasticities is compiled.

\item[Section \ref{sec-fact}.3:] This final (comparatively short) subsection extends the results of Section \ref{sec-fact} from the  case $G$ torsion-free to general finitely generated abelian groups, and sets results in the more general framework of (Transfer) Krull Monoids. It is set separate from other sections so that the main ideas, already quite involved, can be presented in the more natural torsion-free setting without the few additional technical issues that must be handled for the general case. Proposition \ref{prop-fg-rho-transfer} contains the basic observations that allow finiteness of arithmetic invariants to be transferred between $G$ and $G/G_T$, where $G_T\leq G$ is the torsion subgroup. In Main Theorem \ref{thm-tame-pseudo-tor}, a technical modification to the torsion-free argument  extends the validity of Theorem \ref{thm-tame-pseudo} from the torsion-free case to the general case (for Krull Monoids), establishing the finiteness of the weak tame degree. Main Theorem \ref{thm-Delta-finite-tor} likewise extends Theorem \ref{thm-Delta-finite} to the general case,  establishing the finiteness of the set of distances  (for Transfer Krull Monoids),  the finiteness of the Catenary degree (for Krull Monoids), and that the Structure Theorem for Unions holds (for Transfer Krull Monoids).
    Proposition \ref{prop-localtame} contains the basic argument showing that having finite elasticities implies a Krull Monoid is locally tame.
    Proposition \ref{prop-fg-diamondZ} extends the two alternative characterizations of $G_0^\diamond$ given in Corollary \ref{cor-G_0diamond-Zequiv} to the general case.  Proposition \ref{prop-tor-deletion} extends Proposition \ref{prop-finitary-FiniteDeletion}, showing that finite elasticities ensures there is finite subset $Y\subseteq G_0$ such that every zero-sum sequence with terms from $G_0$ must use at least one term from $Y$. Finally, Main Theorem \ref{thm-rho-char-tor} extends Theorem \ref{thm-rho-char} to the general case (for Transfer Krull Monoids), giving our main characterization result for finite elasticities for general $G$. The subsection concludes with a more detailed summary of all the information derived for $G_0$ assuming finite elasticities.

\end{itemize}

\section{Preliminaries and General Notation}\label{sec-prelim}

\subsection{Convex Geometry}
 Variables introduced with inequalities, e.g., $x\geq 1$, are generally assumed to be integers unless otherwise stated. The letters $\alpha$, $\beta$, $\gamma$, $\delta$, $a$, $b$ and $c$ generally indicate real/rational numbers. Intervals are also discrete unless otherwise stated, so $[x,y]=\{z\in \Z:\;y\leq x\leq z\}$ for $x,\,y\in \R$. We use $\subset $ to denote proper inclusion.

For $x\in \R^d$, we let $\|x\|$ denote the usual Euclidian $L_2$-norm. Then $\mathsf d(x,y)=\|x-y\|$ is the distance between two points $x,\,y\in \R^d$. Given two nonempty sets $X,\,Y\subseteq \R^d$, we let $$\mathsf d(X,Y)=\inf\{\mathsf d(x,y):\; x\in X,\,y\in Y\}$$ and define $X+Y=\{x+y:\;x\in X,\,y\in Y\}$, $-X=\{-x:\;x\in X\}$ and $X-Y=X+(-Y)$.  For a real number $\epsilon>0$, let $B_\epsilon(x)$ denote an open ball of radius $\epsilon$ centered at $x$.

Let $X\subseteq \R^d$ be a subset of the $d$-dimensional Euclidian space, where $d\geq 0$. Then
\begin{align*}&\Z\la X\ra =\{\Sum{i=1}{n}\alpha_ix_i:\; n\geq 0,\; x_i\in X, \; \alpha_i\in \Z\}, \quad \Q\la X\ra=\{\Sum{i=1}{n}\alpha_ix_i:\; n\geq 0, \; x_i\in X, \; \alpha_i\in \Q\},
\\&\und \quad\R\la X\ra =\{\Sum{i=1}{n}\alpha_ix_i:\; n\geq 0,\; x_i\in X, \; \alpha_i\in \R\}, \end{align*} so $\R\la X\ra$ denotes the linear subspace spanned by $X$. Note $\R\la \emptyset\ra=\{0\}$.  For $x_1,\ldots,x_n\in \R^d$, we let $\R\la x_1,\ldots,x_n\ra=\R\la \{x_1,\ldots,x_n\}\ra$, and likewise define $\Q\la x_1,\ldots,x_n\ra$ and $\Z\la x_1,\ldots,x_n\ra$.
For a subspace $\mathcal E\subseteq \R^d$, we let $\dim \mathcal E$ denote the dimension of $\mathcal E$ and let  $\mathcal E^\bot\subseteq \R^d$ denote the orthogonal complement to $\mathcal E$.
 A subset $X=\{x_1,\ldots,x_n\}\subseteq \R^d$ of size $n\geq 0$ is said to be \textbf{linearly independent modulo $\mathcal E$} if $\Sum{i=1}{n}\alpha_ix_i\in \mathcal E$ with all $\alpha_i\in \R$ implies every $\alpha_i=0$. Equivalently, $\pi(X)$ is a linearly independent set of size $|X|=|\pi(X)|$ (so $\pi$ is injective on $X$), where $\pi:\R^d\rightarrow \mathcal E^\bot$ denotes the orthogonal projection.


Let $\Z_+$, $\Q_+$ and $\R_+$ denote the set of \emph{non-negative} integer, rational and real numbers, respectively.  A \textbf{positive linear combination} of some $x_1,\ldots,x_n\in \R^n$ is an expression of the form $\Sum{i=1}{n}\alpha_ix_i$ with all $\alpha_i\in \R_+$. A positive linear combination is \textbf{nontrivial} if some $\alpha_i>0$, and it is \textbf{strictly positive} if every  $\alpha_i>0$.
A  set $X\subseteq \R^d$  is \textbf{convex} if $x,\,y\in X$ implies $\alpha x+(1-\alpha)y\in X$ for all real numbers $\alpha\in [0,1]$, that is, all points on the line segment between $x$ and $y$ lie in $X$. A set $C\subseteq \R^d$ is a \textbf{cone} if $x\in C$ implies $\alpha x\in C$ for all real numbers $\alpha>0$, that is, the entire ray $\R_+x$ is contained in $C\cup \{0\}$. A \textbf{convex cone} is a set $C\subseteq \R^d$ that is both a cone and convex, i.e.,  $C$ is closed under nontrivial positive linear combinations: $x,\,y\in C$ implies that $\alpha x+\beta y\in C\cup \{0\}$ for all $\alpha,\,\beta\in \R_+$. The convex cones spanned by a subset $X\subseteq \R^d$ are defined as
\begin{align*}&\C_\Z( X) =\{\Sum{i=1}{n}\alpha_ix_i:\; n\geq 0,\; x_i\in X, \; \alpha_i\in \Z_+\}, \quad \C_\Q(X)=\{\Sum{i=1}{n}\alpha_ix_i:\; n\geq 0,  x_i\in X,  \alpha_i\in \Q_+\},
\\&\und \quad\C(X)=\C_\R(X) =\{\Sum{i=1}{n}\alpha_ix_i:\; n\geq 0,\; x_i\in X, \; \alpha_i\in \R_+\}. \end{align*} If $X$ is finite, then $\C(X)$ is a closed convex cone. We likewise define
$$\C^*(X)=\{\Sum{i=1}{n}\alpha_ix_i:\; n\geq 1, \; x_i\in X, \; \alpha_i\in \R_+, \; \alpha_i>0\}.$$ Note $\C^*(X)$ differs from $\C(X)$ only in that $0\in \C(X)$ is trivial while $0\in \C^*(X)$ is not. As special cases, we have $\C(\emptyset)=\{0\}$ and $\C^*(\emptyset)=\emptyset$. For $x_1,\ldots, x_n\in \R^d$, we use the abbreviations $\C(x_1,\ldots,x_n)=\C(\{x_1,\ldots,x_n\})$, $\C_\Q(x_1,\ldots,x_n)=\C_\Q(\{x_1,\ldots,x_n\})$, $\C_\Z(x_1,\ldots,x_n)=\C_\Z(\{x_1,\ldots,x_n\})$ and $\C^*(x_1,\ldots,x_n)=\C^*(\{x_1,\ldots,x_n\})$.

Given a subset $X\subseteq \R^d$, we let $X^\circ$ denote the relative interior of $X$, which is the interior of $X$ relative to the topological subspace $\R\la X-X\ra+X$ (called the affine subspace spanned by $X$), let $\partial(X)$ denote the relative boundary of $X$, which is the boundary of $X$ relative to the topological subspace $\R\la X-X\ra+X$, let $\overline X$ denote the closure of  $X$ (which is the same in both $\R\la X-X\ra+X$ and $\R^d$), and  let $\mathsf{int}(X)$ denote the interior of $X$ in $\R^d$. Note, if $0\in \overline X$ (e.g., if $X$ is a convex cone), then $\R\la X-X\ra+X=\R\la X\ra$ is simply the linear subspace spanned by $X$.
Note $$\R_+^\circ=\{\alpha\in \R:\;\alpha>0\}$$ is the set of strictly positive real numbers.
We set $$\C^\circ(X):=\C(X)^\circ \quad\und\quad \C^\circ(x_1,\ldots,x_n):=\C(x_1,\ldots,x_n)^\circ.$$
A cone $C\subseteq \R^d$ of the form $C=\C(X)$ for some finite subset $X\subseteq \R^d$ is called a finitely generated or \textbf{polyhedral cone}.
If $X$ is linearly independent, then $\C(X)$ is a \textbf{simplicial cone}, in which case $\C^\circ(X)=\{\Sum{i=1}{n}\alpha_ix_i:\;\alpha_i>0\}$, where $X=\{x_1,\ldots,x_n\}$ \cite[Theorem 4.17]{convexbookI}.
A  co-dimension one subspace $\mathcal H\subset \R^d$  naturally divides $\R^d$ into  two closed half spaces $\mathcal H_+$ and $\mathcal H_-$ whose common boundary is $\mathcal H$.

The texts \cite{50years} \cite{Fenchel} \cite{ConvexII} \cite{Gr07a} \cite{HandbookConvex} \cite{convexbookI}   contains many of the basic properties of convexity that we will regularly use with little further reference. The following are several important highlights. We remark that the slightly more general version of Carath\'eodory's Theorem given below is rather difficult to find so stated despite encapsulating what is actually proved, particularly the portion regarding the representation of $0$. The latter can be derived with ease from the argument used to prove  Carath\'eodory's Theorem or  from the basic theory of positive bases (see Section \ref{sec-reay}). It  is also a special consequence of Theorem \ref{thm-carahtheodory-elm-atom}, whose  brief proof we provide.  Let $X\subseteq \R^d$.

\begin{itemize}
\item (Carath\'eodory's Theorem) If  $x\in \C^*(X)$, then there exists a subset $Y\subseteq X$ with $|Y|\leq d+1$ and $x\in \C^*(Y)$. Moreover, if  $x\neq 0$, then $Y$ is linearly independent, and if $x=0$, then any proper subset of $Y$ is linearly independent  \cite[Theorem 4.27]{convexbookI}.

    \item  $\C(X)\cap -\C(X)$, called the \textbf{lineality space} of $\C(X)$, is the unique maximal linear subspace contained in $\C(X)$ \cite[Theorem 4.15]{convexbookI}, and is nontrivial if and only if $0\in \C^*(X\setminus \{0\})$.

 \item (Minkowski-Weyl Theorem) $\C(X)$ is a polyhedral cone if and only if $\C(X)$ is the intersection of a finite number of closed half-spaces \cite[Theoerem 11.9]{50years}.

\item (Duality) A closed convex set $X\subset \R^d$   is either the empty set or the intersection of all closed  half-spaces that contain $X$ \cite[Theorem 2.7]{ConvexII}.

    \item If    $\C(X)\neq \R\la X\ra$, then
    $\C(X)$ lies in a closed half-space of $\R\la X\ra$ \cite[Corollary 1]{Fenchel}.

    \item (Relative Interior and Closure) If $X\subseteq \R^d$ is  convex, then $X^\circ=\overline{X}^\circ$ and $\overline{X}=\overline{X^\circ}$, with both these sets convex \cite[Theorem 2.38]{convexbookI} \cite[Theorem 2.35]{convexbookI} \cite[Corollary 2.22]{convexbookI}.
    \item If $X\subseteq\R^d$ is  convex and $X\neq \emptyset$, then $X^\circ\neq \emptyset$ \cite[Corollary 2.18]{convexbookI}.

\end{itemize}

\subsection{Lattices and Partially Ordered Sets}\label{sec-poset-lattices}

A \textbf{lattice} is a discrete subgroup $\Lambda\leq \R^d$, where $d\geq 0$, meaning \emph{any bounded subset of $\R^d$ contains only finitely many lattice points} \cite{Cassels} \cite{Handbodk-Gruber} \cite{Handbook-gritzmann-wills}. In particular, any convergent sequence of lattice points stabilizes, that is, if  $\{x_i\}_{i=1}^\infty$ is a convergent sequence of lattice points $x_i\in \Lambda$ with $\lim_{i\rightarrow \infty} x_i= x$, then $x\in \Lambda$ and $x_i=x$ for all sufficiently large $i$.  As is well-known \cite[Theorem VI]{Cassels} \cite[Theorem 1]{Handbodk-Gruber}, being a lattice is equivalent to there existing a linearly independent generating subset $X\subseteq \Lambda$, so $\Z\la X\ra=\Lambda$ with $X\subseteq \R^d$ linearly independent. Such a subset $X\subseteq \Lambda$ is called a \textbf{lattice basis}. The rank of the lattice $\Lambda$ is $n=\dim \R\la \Lambda\ra=|X|$, and we say $\Lambda$ has \textbf{full rank} in $\R^d$ if $n=d$, that is, $\R\la \Lambda\ra=\R^d$. In particular, $\Lambda\cong \Z^n$ for a rank $n$ lattice $\Lambda$.

The following basic consequence of the Smith Normal Form \cite[Theorem III.7.8]{lang} \cite[Theorem 2]{Handbodk-Gruber} will be very important for us. We remark that the hypothesis that the kernel be generated by a subset of lattice points is quite necessary.

\begin{proposition}\label{Prop-lattice-homoIm}
Let $\Lambda\subseteq \R^d$ be a full rank lattice and let $\pi:\R^d\rightarrow \R^d$ be a linear transformation with $\ker \pi=\R\la X\ra$ for some $X\subseteq \Lambda$. Then $\pi(\Lambda)$ is a lattice having full rank in the subspace $\mathsf{im}\,\pi$.
\end{proposition}

\begin{proof}  By passing to a subset of $X$, we can w.l.o.g. assume the elements $x_1,\ldots,x_s\in X\subseteq \Lambda$ are linearly independent lattice points, and thus form a linear basis for $\ker \pi=\R\la X\ra$ with $\Z\la X\ra\leq \Lambda$ a sublattice of rank $s$. Via the Smith Normal form, we can find a lattice basis $e_1,\ldots,e_d\in \Lambda$ and integers $m_1\mid \ldots\mid m_d$, where $m_i>0$ for $i\leq s$ and $m_i=0$ for $i>s$, such that $y_i=m_ie_i$ for all $i\leq s$ with $\Z\la y_1,\ldots,y_s\ra=\Z\la X\ra=\Z\la x_1,\ldots,x_s\ra$.
Since $\ker \pi=\R\la x_1,\ldots,x_s\ra=\R\la y_1,\ldots,y_s\ra=\R\la m_1e_1,\ldots,m_se_s\ra$, it is now clear that $\pi(e_i)=0$ for $i\leq s$ while $\pi(e_{s+1}),\ldots,\pi(e_{d})\in \pi(\Lambda)$ are distinct linearly independent elements which generate the subgroup $\pi(\Lambda)$. This shows  $\pi(\Lambda)\leq \mathsf{im}\,\pi$ is a lattice, and since $\Lambda$ has full rank in $\R^d$, it follows that $\pi(\Lambda)$ has full rank in $\mathsf{im}\; \pi$.
\end{proof}

A \textbf{partially ordered set (poset)} $(P,\preceq)$ is a set $P$ together with a partial order $\preceq$, so $\preceq$ is a transitive, reflexive, anti-symmetric relation on $P$ \cite{Roman}. In such case, we use $x\prec y$ to indicate $x\preceq y$ but $x\neq y$. Any subset $X\subseteq P$ is also a partially ordered set using the partial order $\preceq $ inherited from $P$,
 and we implicitly consider subsets of posets to be posets using the inherited partial order.
 A maximal element $x\in  P$ is one for which there is no $y\in P$ with $x\prec y$. Likewise, a minimal element $x\in  P$ is one for which there is no $y\in P$ with $y\prec x$. We let $$\mathsf{Min}( P)\subseteq  P\quad\und\quad\mathsf{Max}( P)\subseteq  P$$ denote the set of minimal and maximal elements of $P$, respectively. Given a subset $X\subseteq P$, we let $$\darrow X=\{y\in \ P:\; y\preceq x\mbox{ for some $x\in X$}\}\quad\und\quad \uarrow X=\{y\in P:\; x\preceq y\mbox{ for some $x\in X$}\}$$ donate the down-set and up-set generated by $X$, respectively. We likewise define $\darrow x=\darrow \{x\}$ and $\uarrow x=\uarrow \{x\}$ for $x\in X$. An \textbf{anti-chain} is a subset $X\subseteq P$ such that no two distinct elements of $X$ are comparable. A \textbf{chain} is a totally ordered subset of $P$.  An \textbf{ascending chain} is a sequence $x_1\preceq x_2\preceq \ldots\preceq x_n\preceq\ldots$ with $x_i\in P$. Likewise, a \textbf{descending chain} is a sequence $x_1\succeq x_2\succeq \ldots \succeq x_n\succeq \ldots$ with $x_i\in P$. In both cases, the chain can either be finite (stopping at $n$) or infinite, and the chain is \textbf{strict} if each $\preceq$ or $\succeq$ in the chain  is always a strict inclusion  $\prec$ or $\succ$.

 If $(P,\preceq_1 )$ and $(Q,\preceq_2)$ are both posets, then $P\times Q$ is also a poset using the product partial order: $(x,y)\preceq (x',y')$ when $x\preceq_1 x'$ and $y\preceq_2 y'$. An important example of a partially ordered set is $P=\Z_+^d$ equipped with the product partial order: $(x_1,\ldots,x_d)\preceq (y_1,\ldots,y_d)$ when $x_i\leq y_i$ for all $i$.
 It is easily seen that the poset $\Z_+^d$ has no infinite strictly descending chain. Indeed, there are at most $x_1+\ldots+x_d+1$ elements in any strictly descending chain whose first element is  $(x_1,\ldots,x_d)$. A well-known consequence of Hilbert's Basis Theorem \cite[Theorem 4.1]{lang} is that the poset $\Z_+^d$ contains no infinite anti-chain: The points of $\Z_+^d$ correspond naturally to the monomials  in $\mathbb F_2[x_1,\ldots,x_d]$, and  if $X\subseteq \Z_+^d$ were an infinite anti-chain, then the monomial ideal generated by the monomials corresponding to the elements from $X$ would not be finitely generated, contrary to Hilbert's Basis Theorem.

We will be interested in posets which neither contain an infinite anti-chain nor an infinite strictly descending chain. Such posets are called \textbf{well-quasi-orderings} \cite{Kruskall}. We include the following basic propositions about such posets (mentioned in \cite{Kruskall} without proof) for completeness, which together give an alternative proof (as opposed to Hilbert's Basis Theorem) that any subset $X\subseteq \Z_+^d$ has only a finite number of minimal points (a result known as Dickson's Lemma \cite[Theorem 1.5.3]{alfredbook}).

\begin{proposition}\label{prop-poset-conseq}
If $(P,\preceq)$ is a poset that contains neither infinite anti-chains nor infinite strictly descending chains, then $\mathsf{Min}(X)$ is finite  and $\uarrow \mathsf{Min}(X)\cap X=X$ for any subset $X\subseteq P$.\end{proposition}

\begin{proof}
The set of minimal points in a poset is an anti-chain and the hypothesis that $P$ contains no infinite anti-chain inherits to $X$. Thus $\mathsf{Min}(X)$ is finite. If
$Y=X\setminus \uarrow \mathsf{Min}(X)$ is nonempty, then $Y=\darrow Y\cap X\subseteq X$ can contain no minimal point, and we can recursively select elements from $Y$ to form an infinite descending chain in $X$, and thus also in $P$, contrary to hypothesis.
\end{proof}

\begin{proposition}\label{prop-poset-prod}
If both $(P,\preceq)$ and $(Q,\preceq')$ are posets that contain neither infinite anti-chains nor infinite strictly descending chains, then the poset $P\times Q$, equipped with the product partial order, also contains no infinite anti-chains nor infinite strictly descending chains.
\end{proposition}

\begin{proof}
Assume to the contrary that  $\{(x_i,y_i)\}_{i=1}^\infty$ is a strictly descending chain. Then, per definition of the product partial order, both $\{x_i\}_{i=1}^\infty$ and $\{y_i\}_{i=1}^\infty$ must be descending chains. Since neither $P$ nor $Q$ contains infinite strictly descending chains, it follows that there exists some index $M$ such that  $x_i=x_j$ for all $i,\,j\geq M$, and likewise some index $N$ such that  $y_i=y_j$ for all $i,\,j\geq N$. Hence $(x_i,y_i)=(x_j,y_j)$ for all $i,\,j\geq \max\{N,M\}$, contradicting that $\{(x_i,y_i)\}_{i=1}^\infty$ is a strictly descending chain. This shows that $P\times Q$ contains no infinite strictly descending chain.

Next assume to the contrary that $Z=\{(x_1,y_1),(x_2,y_2),\ldots,\}$ is an infinite anti-chain.
Let $X=\{x_1,x_2,\ldots,\}$ and $Y=\{y_1,y_2,\ldots\}$. If $X$ is finite, then the pigeonhole principle ensures that there is some $x\in X$ that occurs as the first coordinate of an infinite number of pairs from $Z$. Thus, by passing to a subset, we can w.l.o.g. assume $Z=\{(x,y_1),(x,y_2),\ldots,\}$, in which case $Z$ can only be an infinite anti-chain, per definition of the product partial order, if $\{y_1,y_2,\ldots,\}$ is an infinite anti-chain in $Q$, which is assumed to not exist by hypothesis. Therefore instead assume that $X\subseteq P$ is infinite. Consequently, since any infinite poset contains either an infinite chain or an infinite anti-chain, and since any infinite chain contains either an infinite strictly descending or infinite strictly ascending chain \cite[Theorem 1.14]{Roman}, our hypotheses ensure that $X\subseteq P$ contains an infinite strictly  ascending chain. Hence, by passing to a subset of $Z$ and re-indexing the elements of $Z$, we can w.l.o.g. assume $\{x_i\}_{i=1}^\infty$ is an infinite strictly ascending chain in $P$.
Repeating the same argument using $Y$ instead of $X$, we conclude that $Y$ also contains an infinite strictly ascending chain, say $\{y_{i_j}\}_{j=1}^\infty$. However, since $\{x_i\}_{i=1}^\infty$ is an ascending chain but $Z=\{(x_1,y_1),(x_2,y_2),\ldots,\}$ is an anti-chain, the definition of the product partial order forces $i_{j}>i_{j+1}$ for all $j$. Thus $i_1>i_2>\ldots$ is an infinite strictly descending chain in $\Z_+$, which is not possible, completing the proof.
\end{proof}

\subsection{Sequences and Rational Sequences}\label{sec-intro-factorization}

Let $G\cong \Z^d\oplus G_T$ be a finitely generated  abelian group with  torsion subgroup $G_T$.  We let $\exp(G_T)$ denote the exponent of $G_T$, which is the minimal integer $m\geq 1$ such that $mg=0$ for all $g\in G_T$. Let $G_0\subseteq G$ be a subset.
Regarding sequences and sequence subsums over $G_0$, we follow the standardized notation from Factorization Theory  \cite{gao-ger-survey} \cite{Alfred-Ruzsa-book} \cite{alfredbook} \cite{Gbook}.  The key parts are summarized  here.

A \textbf{sequence} $S$ of terms from $G_0$ is viewed formally as an element of the free abelian monoid  with basis $G_0$, denoted  $\mathcal F(G_0)$. Context will always distinguish between a sequence $S\in \Fc(G_0)$ and a sequence $\{x_i\}_{i=1}^\infty$ of terms $x_i\in G_0$.  A sequence $S\in \Fc(G_0)$ is written as a finite multiplicative string of terms, using the boldsymbol  dot operation $\bdot$ to concatenate terms, and with the order irrelevant:
$$S=g_1\bdot\ldots\bdot g_\ell={\prod}^\bullet_{g\in G_0}g^{[\vp_g(S)]}$$ with $g_i\in G_0$ the terms of  $S$, with $\vp_g(S)=|\{i\in [1,\ell]:\;g_i=g\}|\in \Z_+$ the multiplicity of the term $g\in G_0$, and with $$|S|:=\ell=\Summ{g\in G_0}\vp_g(S)\geq 0\quad\mbox{ the \textbf{length} of $S$}\quad\und\quad \vp_X(S)=\Summ{x\in X}\vp_x(S),$$ where $X\subseteq G_0$.
Here
  $g^{[n]}={\underbrace{g\bdot\ldots\bdot g}}_n$ denotes the sequence consisting of the element $g$ repeated $n$ times, for $g\in G_0$ and $n\geq 0$. The notation is extended to sequences as well: $S^{[n]}={\underbrace{S\bdot\ldots\bdot S}}_n$.
If $S,\,T\in \Fc(G_0)$ are sequences, then  $S\bdot T\in \Fc(G_0)$ is  the sequence obtained by concatenating the terms of $T$ after those of $S$.
 We use $T\mid S$ to indicate that $T$ is a subsequence of $S$ and let ${T}^{[-1]}\bdot S$ or $S\bdot {T}^{[-1]}$ denote the sequence obtained by  removing the terms of $T$ from $S$.
Then
\begin{align*}
 &\supp(S)=\{g\in G_0:\; \vp_g(S)>0\}\subseteq G_0 \quad\mbox{ is the \textbf{support} of $S$},\quad\und\\
 &\sigma(S)=\Sum{i=1}{\ell}g_i=\Summ{g\in G_0}\vp_g(S)g\in G\quad\mbox{ is the sum of terms from $S$}.
\end{align*}
Given two sequences $S,\,T\in \Fc(G_0)$, we let $\gcd(S,T)\in \Fc(G_0)$ denote the maximal length sequence diving both $S$ and $T$.

Given a map
$\varphi \colon G_0 \to G'_0$, we let $\varphi(S)=\varphi(g_1)\bdot\ldots\bdot \varphi(g_\ell)\in \Fc(G'_0)$. The sequence $S$ is called \textbf{zero-sum} if $\sigma(S)=0$, and we let $\mathcal B(G_0)\subseteq \Fc(G_0)$ be the set of all zero-sum sequences over $G_0$. A \emph{nontrivial} sequence $S\in \Fc(G_0)$ is called an \textbf{atom} or \textbf{minimal zero-sum sequence} if $\sigma(S)=0$ but $\sigma(S')\neq 0$ for all proper, nontrivial subsequences $S'\mid S$. We let $\mathcal A(G_0)\subseteq \mathcal B(G_0)$ denote the set of all atoms over $G_0$. The \textbf{Davenport constant} for $G_0$ is $$\D(G_0)=\sup\{|U|:\;U\in \mathcal A(G_0)\}\in \Z_+\cup \{\infty\},$$ namely, the maximal length of an atom. Assuming $G$ is finitely generated and every $g\in G_0$ is contained in some zero-sum sequence over $G_0$, it is known that $\D(G_0)$ is finite if and only if $G_0$ is finite \cite[Theorem 3.4.2]{alfredbook}, and we have the general upper bound $\mathsf D(G)\leq |G|$ \cite[Propositoin 5.1.4.4]{alfredbook} \cite[Theorem 10.2]{Gbook}.

A sequence $S\in \Fc(G_0)$ has the form $S=\prod^\bullet_{g\in G_0}g^{[\vp_g(S)]}$ with $\vp_g(S)\in \Z_+$ the multiplicity of $g$ in $S$ and $\vp_g(S)>0$ for only finitely many $g\in G_0$. Thus $\Fc(G_0)$ can viewed as tuples from the monoid $(\Z_+^{G_0},+)$ with finite support, with the sequence $S=\prod^\bullet_{g\in G_0}g^{[\vp_g(S)]}$ corresponding to the tuple $(\vp_g(S))_{g\in G_0}$.  It will be useful to sometimes allow more general  exponents for terms in a sequence. This was done based on ideas from Matroid Theory in \cite{G-precursor}. We continue with a natural framework utilizing Convex Geometry instead.  The setup is as follows, which we present only in the context when $G$ is torsion-free (though a similar framework holds for general abelian groups provided one restricts to rational  multiplicities whose denominators are relatively prime to every $\ord(g)<\infty$ for $g\in G$).
We let $\Fc_{\mathsf{rat}}(G_0)$ denote the multiplicative monoid whose elements have the form $S=\prod^\bullet_{g\in G_0}g^{[\vp_g(S)]}$ with $\vp_g(S)\in \Q_+$ the multiplicity of $g$ in $S$ and $\vp_g(S)>0$ for only finitely many $g\in G_0$. We refer to the elements $S\in \Fc_{\mathsf{rat}}(G_0)$ as \textbf{fractional} or \textbf{rational sequences}. Thus the rational sequences $\Fc_{\mathsf{rat}}(G_0)$ can be viewed as tuples from the monoid $(\Q_+^{G_0},+)$ with finite support, with the sequence $S=\prod^\bullet_{g\in G_0}g^{[\vp_g(S)]}$ corresponding to the tuple $(\vp_g(S))_{g\in G_0}$. Any $S\in \Fc_{\mathsf{rat}}(G_0)$ has some positive integer $N>0$ such that $S^{[N]}\in \Fc(G_0)$ (simply take the least common multiple of all denominators of the nonzero $\vp_g(S)$).
We extend all notation regarding $\Fc(G_0)$ to $\Fc_{\mathsf{rat}}(G_0)$. In particular, $\supp(S)=\{g\in G_0:\;\vp_g(S)>0\}$, which is always a finite set by definition, and $|S|=\Summ{g\in G_0}\vp_g(S)$.
We would also like to speak of the sum of a rational sequence $S\in \Fc(G_0)$. Since $G$ is torsion-free, we have $G=\Z\la X\ra\cong \Z^{|X|}$ for some independent set $X$, and we can then embed $G$ into $Q=\Q^{|X|}$ (thus we embed $G$ into a divisible abelian group $Q$ \cite{hungerford}).
This means that the sum $\sigma(S)=\Summ{g\in G_0}\vp_g(S)g\in Q$ of a rational sequence $S\in \Fc_{\mathsf{rat}}(G_0)$ is now well-defined.
Let  $\mathcal B_{\mathsf{rat}}(G_0)$ consist of all rational sequences $S\in \Fc_{\mathsf{rat}}(G_0)$ whose sum is zero $\sigma(S)=\Summ{g\in G_0}\vp_g(S)g=0$.
This means $\mathcal B_{\mathsf{rat}}(G_0)$ is a rational generalized Krull Monoid \cite[Proposition 2]{Ch-HK-Kr02}.
 For $S\in \Fc_{\mathsf{rat}}(G_0)$, we define
$$\lfloor S\rfloor: ={\prod}^\bullet_{g\in G_0}g^{[\lfloor \vp_g(S)\rfloor]}\quad\und\quad \{S\}:={\prod}^\bullet_{g\in G_0}g^{[\{ \vp_g(S)\}]},$$ where $\{\vp_g(S)\}:=\vp_g(S)-\lfloor \vp_g(S)\rfloor$ denotes the fractional part. Thus $S=\lfloor S\rfloor \bdot \{S\}$ with $\lfloor S\rfloor \mid S$ the  unique maximal length subsequence  with $\lfloor S\rfloor \in \Fc(G_0)$, and with $\vp_g(\{S\})<1$ for all $g\in G_0$. It may be helpful to view the elements in $\{S\}$ as `damaged', and they will need to be avoided in some arguments as much as possible. Note $g\in \supp(\{S\})$ if and only if $\vp_g(S)\notin \Z$, while $g\in \supp(\lfloor S\rfloor)$ if and only if $\vp_g(S)\geq 1$.
The following properties for $S,\,T\in \Fc_{\mathsf{rat}}(G_0)$ are immediate from the definitions but quite useful:
\begin{align}
\label{RatSeq1}&|\{S\}|\leq |\supp(\{S\})|,\\\label{RatSeq2}
&|\supp(S)|,\,|S|\leq |\lfloor S\rfloor |+|\supp(\{S\})|,\\\label{RatSeq3} &\lfloor S\rfloor \bdot \lfloor T\rfloor\mid \lfloor S\bdot T\rfloor \quad\und\quad
\supp(\lfloor S\rfloor)\cup \supp(\lfloor T\rfloor)\subseteq \supp(\lfloor S\bdot T\rfloor),
\\\label{RatSeq4}
&\{S\bdot T\}\mid \{S\}\bdot \{T\}\quad\und\quad\supp(\{S\bdot T\})\subseteq \supp(\{S\})\cup \supp(\{T\}),\end{align}\begin{align}\label{RatSeq5}
 &\supp(\{S\bdot T^{[-1]}\})\subseteq \supp(\{S\})\cup \supp(\{T\})\quad\mbox{ for $T\mid S$},\\\label{RatSeq6}
 &\lfloor T\rfloor \mid \lfloor S\rfloor\quad\und\quad \lfloor S\bdot T^{[-1]}\rfloor\mid \lfloor S\rfloor \bdot \lfloor T\rfloor^{[-1]}\quad \mbox{ for } T\mid S.
\end{align}

While generally the order of terms in a sequence will be irrelevant, there are times when we will need to view our sequences as indexed, so that two distinct terms of a sequence $S\in \Fc(G_0)$ that are equal as elements can still be viewed as  distinct terms in the sequence $S$. For such occasions, we fix an indexing of the terms of $S$, say $$S=g_1\bdot\ldots\bdot g_\ell$$ with $g_i\in G_0$. Then for $I\subseteq [1,\ell]$,  we let $$S(I)={\prod}^\bullet_{i\in I}g_i\in \Fc(G_0)$$ denote the subsequence of terms indexed by $I$.

\subsection{Arithmetic Invariants for Transfer Krull Monoids}\label{sec-transferkrull}

Let $H$ be a unit-cancellative (associative and multiplicatively written) semigroup with identity $1_H$. A semigroup homomorphism is map $\varphi:H\rightarrow H'$ between semigroups with identity such that $\varphi(1_H)=\varphi(1_{H'})$ and $\varphi(xy)=\varphi(x)\varphi(y)$ for $x,\,y\in H$. Let $H^\times\subseteq H$ denote the subgroup of units (invertible elements) in $H$.
 If a non-unit $a\in H\setminus H^\times$ has a factorization   $a=u_1\cdot\ldots\cdot u_k$ with each $u_i\in H$ a non-unit, then  we call $k$ the length of the factorization. We are  generally interested in the case of factorizations into irreducible elements of $H$, called \textbf{atoms}, so $u_i\in \mathcal A(H)$ for all $i$ with $\mathcal A(H)$ the set of all atoms (irreducible elements) in $H$.  Let $\mathsf L(a)$ denote the set of all $k\geq 1$ for which $a$ has a factorization \emph{into atoms} of length $k$ and let $\mathcal L(H)=\{\mathsf L(a):\; a\in H\setminus H^\times\}$ denote the system of sets of lengths. Since, for purposes of factorization, elements of $H$ equal up to units are essentially identical, it is often convenient to replace $H$ (for $H$ commutative) with the reduced quotient monoid $H_{\mathsf{red}}=H/H^\times$ where this has been formally implemented.

A \textbf{Krull Monoid} is a commutative, cancellative semigroup $H$ with identity $1_H$ such that there is a \textbf{divisor homomorphism} $\varphi:H\rightarrow \mathcal F(P)$ into a free abelian monoid such that every $p\in P$ has some nonempty $X\subseteq H$ with $p=\gcd(\varphi(X))$. That $\varphi$ is divisor homomorphism means $\varphi$ is a semigroup homomorphism such that  $x\mid y$ in $H$ if and only if $\varphi(x)\mid \varphi(y)$ in $F$, where $x,\,y\in H$ \cite[Theorem 2.4.8]{alfredbook} \cite[Definition 4.]{Ge-Zh20a}. The monoid $\mathcal B(G_0)$ of zero-sum sequences is always a Krull Monoid, for any subset $G_0$ of an abelian group $G$, and the reduced monoid $H_{\mathsf{red}}=H/H^\times$of a Krull Monoid $H$ is also always a Krull Monoid.
A \textbf{Transfer Krull Monoid} is a unit-cancellative semigroup $H$ with identity $1_H$ such that there is a \textbf{(weak) transfer homomorphism} $\theta:H\rightarrow \mathcal B(G_0)$ for some subset $G_0$ of an abelian group $G$. The latter means $\theta$ is a surjective semigroup homomorphism  such that $\theta(z)$ is the trivial sequence if and only if $z\in H^\times$ is a unit, for $z\in H$,  and whenever $\theta(z)=U_1\bdot\ldots\bdot U_k$ with $U_1,\ldots,U_k\in \mathcal A(G_0)$, for $z\in H$, then there is permutation $\tau:[1,k]\rightarrow [1,k]$ and factorization  $z=x_1\cdot\ldots\cdot x_k$ with  $x_1,\ldots,x_k\in \mathcal A(H)$ such that  $\theta(x_i)=U_{\tau(i)}$ for all $i\in [1,k]$ (see \cite{Ge16c} \cite{Ba-Ba-Go14} \cite{Ba-Sm15} \cite{Sm16a} \cite{Sm19a}). Then $z\in H$ is an atom if and only if $\theta(z)\in \mathcal B(G_0)$ is an atom, and the permutation $\tau$ may be taken to be the identity when $H$ is commutative. Moreover, $\mathcal L(H)=\mathcal L(\mathcal B(G_0))$, meaning  any arithmetic invariant regarding the factorization of elements from $H$ that depends only on $\mathcal L(H)$ reduces to the study of the corresponding invariant for $\mathcal B(G_0)$.

We now  briefly review the basic theory of Krull Monoids. To this end, consider a Krull Monoid $H$ with divisor homomorphism $\varphi:H\rightarrow \mathcal F(P)$. The elements of $\mathcal F(P)$ are sequences $S=\prod_{p\in P}^\bullet p^{[\vp_p(S)]}$ with $\vp_p(S)\in \Z_+$. The set $P$ does not live in an abelian group apriori, but this can be resolved by formally defining $\wtilde G_P$ to be a  free abelian group with basis $P$, in which case $P\subseteq\wtilde G_P\cong \Z^P$. Now let $\wtilde H_P=\la \sigma(\varphi(x)):\;x\in H\ra\leq \wtilde G_P$,
let $G=\wtilde G_P/\wtilde H_P$, and let $[\cdot]:\wtilde G_P\rightarrow G=\wtilde G_P/\wtilde H_P$ be the natural homomorphism. For $S=\prod_{p\in P}^\bullet p^{[\vp_p(S)]}\in \mathcal F(P)$, we have $[S]=\prod_{p\in P}^\bullet [p]^{[\vp_p(S)]}\in \mathcal F(G_0)$, where $G_0=\{[p]:\; p\in P\}\subseteq G$, referred to as  the set of classes containing primes.
We remark that our notation is slightly different from the presentation given in \cite{alfredbook} \cite{Alfred-Ruzsa-book}. There, the formal quotient group $q(\Fc(P))$ of $\mathcal F(P)$ is used in place of $\wtilde G_P$, with $q(\varphi(H))$ used in place of $\wtilde H_P$, at times leading to confusion since $[S]$, for $S\in \mathcal F(P)$,  can both represent a sequence $[S]\in \mathcal F(G)$ as well as an element $[S]\in G$, the latter corresponding to the sum of terms in $[S]$ when considered as a sequence. The reformulation followed here avoids this potential confusion.

Any sequence $S=\prod_{p\in P}^\bullet p^{[\vp_p(S)]}\in \varphi(H)\subseteq \Fc(P)$ will have $[S]\in \mathcal B(G_0)$ by definition of $[\cdot]$.
This means the composition map $\theta:H\rightarrow \mathcal B(G_0)$, defined by $\theta(x)=[\varphi(x)]$, is a semigroup homomorphism taking each element of $H$ to a zero-sum sequence over the subset $G_0$ of the abelian group $G$. Now further assume $H$ is reduced, that is, assume we reduce our original Krull Monoid $H$ by replacing it with the reduced Krull Monoid $H_{\mathsf{red}}=H/H^\times$. Under this assumption,  the map $\theta$ is a transfer homomorphism, showing that the Krull Monoid $H_{\mathsf{red}}$ is a commutative Transfer Krull Monoid \cite[Theorem 1.3.4]{Alfred-Ruzsa-book}. Thus (up to units) every Krull Monoid is also a Transfer Krull Monoid. Moreover \cite[Proposition 2.4.2.4, Corollary 2.4.3.2]{alfredbook} (see also \cite[Theorem 1.3.4.1]{Alfred-Ruzsa-book}), \be\label{magic-krull-divisors} H_{\mathsf{red}}\cong \varphi(H)\quad\und\quad \varphi(H)=\{S\in \Fc(P):\;[S]\in \mathcal B(G_0)\}.\ee We then say $H$ is a Krull Monoid over the subset $G_0$ of the abelian group $G$.

We continue by describing some of the most commonly studied arithmetic invariants for factorizations that depend only on $\mathcal L(H)$.  For an integer $k\geq 1$, we let $\mathcal U_k(H)$ denote the set of all $\ell\geq 1$ for which there are atoms $U_1,\ldots,U_k,V_1,\ldots,V_\ell\in \mathcal A(H)$ with \be\label{factsample}U_1\bdot\ldots\bdot U_k=V_1\bdot\ldots\bdot V_\ell.\ee Thus $\mathcal U_k(H)$ is the union of all sets from $\mathcal L(H)$ containing $k$. We then define
\begin{align*}&\rho_k(H)=\sup\, \mathcal U_k(H) \quad\mbox{ for $k\geq 1$}, \und\\ &\rho(H)=\sup\Big\{\frac{\rho_k(H)}{k}:\;k\geq 1\Big\}=\sup\Big\{\frac{\sup \mathsf L(S)}{\min \mathsf L(S)}:\;S\in H\Big\}. \end{align*} The equality in the definition of $\rho(H)$ follows by a simple argument \cite[Proposition 1.4.2.3]{alfredbook}. Note  $\ell\leq \rho_k(H)$ must hold for any re-factorization of $k$ atoms as in \eqref{factsample}. A basic argument \cite[Proposition 1.4.2.1]{alfredbook} shows \be\label{rho-ascend-chain}\rho_1(H)\leq \rho_2(H)\leq \ldots.\ee Moreover, the inequalities are strict when the $\rho_k(H)$ are finite. The constant $\rho(H)$ is the \textbf{elasticity} of $H$, and the $\rho_k(H)$ are its refinements, which we  call the \textbf{elasticities} of $H$.
Worth noting, if $H$ is a Transfer Krull Monoid, then it is know that, for every $q\in \Q$ with $1<q<\rho(H)$, there is some $L\in \mathcal L(H)$ with $q=\max L/\min L$ \cite[Theorem 3.1]{Ge-Zh19a}. When $H=\mathcal B(G_0)$, we abbreviate $\rho(G_0)=\rho(\mathcal B(G_0))$, and likewise for all other arithmetic invariants.

Given a set $X\subseteq \Z$, we let $\Delta(X)\subseteq \Z_+$ denote the set of successive distances in $X$, so $\Delta(X)$ consists of all elements of the form $k_2-k_1$, where $k_2,\,k_2\in X$ and $k_1<k_2$, such that $X$ contains no elements from $[k_1+1,k_2-1]$. Then the set of successive distances for $H$ is $$\Delta(H)=\bigcup_{L\in \mathcal L(H)}\Delta(L)=\bigcup_{a\in H\setminus H^\times} \Delta(\mathsf{L}(a))\subseteq \Z_+.$$

A stronger measure of well-behaved factorization is indicated by structural results for $\mathcal U_k(H)$.  A finite set $X\subseteq \Z$ is said to be an almost arithmetic progression with difference $d\geq 1$ and bound $N\geq 0$ if $X=P\setminus Y$, where $P$ is an arithmetic progression with difference $d$ and $Y\subseteq P$ is subset contained in the union of the first $N$ terms from $P$ and the last $N$ terms from $P$.  If there exists a constant $N\geq 0$ and difference $d\geq 1$ such that $\mathcal U_k(H)$ is an almost arithmetic progression with difference $d$ and bound $N$ for all sufficiently large $k$, then we say that $H$ satisfies the \textbf{Structure Theorem for Unions}. As shown in \cite[Theorem 4.2]{Gao-Ger-phok}, if $\Delta(H)$ is finite and there is a constant $M\geq 0$ such that $\rho_{k+1}(H)-\rho_k(H)\leq M$ for all $k\geq 1$, then the Structure Theorem for Unions holds for $H$, providing additional motivation for the study of the elasticities.

The next invariant does not depend solely on $\mathcal L(H)$ but is known to nonetheless reduce to the study of $\mathcal B(G_0)$ for a large class of monoids, including Krull Monoids. Suppose we have  two factorizations $S=U_1\bdot\ldots\bdot U_k$ and $S=V_1\bdot\ldots\bdot V_\ell$ of the same $S\in H$. Re-index the $U_i$ and $V_j$ such that $I\subseteq [1,\min\{k,\ell\}]$ is a maximal subset with $U_i=V_i$ for all $i\in I$, meaning $U_i\neq V_j$ for all $i\in [1,k]\setminus I$ and $j\in [1,\ell]\setminus I$. Then we say that  the factorizations $V_1\bdot\ldots\bdot V_\ell$ and $U_1\bdot\ldots\bdot U_k$ \textbf{differ by $\max\{\ell-|I|,\,k-|I|\}$ factors}. The \textbf{catenary degree} $\mathsf c(H)$ is defined as the minimal integer $N\geq 0$ such that, given any two factorizations $U_1\bdot\ldots\bdot U_k=S=U'_1\bdot\ldots\bdot U'_{k'}$ of the same $S\in H$ into \emph{atoms} $U_i,\,U'_j\in \mathcal A(H)$, there is a sequence of factorizations $S=U_1^{(1)}\bdot\ldots\bdot U_{k_j}^{(j)}$ of $S$ into atoms $U_i^{(j)}\in \mathcal A(H)$, for $j=0,1,\ldots, \ell$, such that \begin{itemize}
\item $k_0=k$ and $U_i^{(0)}=U_i$ for all $i\in [1,k]$,
\item  $k_\ell=k'$ and $U_i^{(\ell)}=U'_i$ for all $i\in [1,k']$, and
\item each factorization $U_1^{(j)}\bdot\ldots\bdot U_{k_j}^{(j)}$ differs from the previous factorization $U_1^{(j-1)}\bdot\ldots\bdot U_{k_{j-1}}^{(j-1)}$, for $j\in [1,\ell]$, by at most $N=\mathsf c(H)$ factors.
\end{itemize}
If no such integer $N$ exists, then $\mathsf c(H)=\infty$. Having finite catenary degree indicates that any factorization can be transformed into any other factorization of the same element by a sequence of small modifications, and we have  $\mathsf c(G_0)\leq \mathsf c(H)\leq \max\{\mathsf c(G_0),2\}$ when $H$ is a Krull Monoid over the subset $G_0$ \cite[Theorem 3.4.10.5]{alfredbook}.

Finally, we define one further invariant. Let $U\in \mathcal A(H)$. If $U$ is prime, we set $\mathsf t(H,U)=0$. Otherwise, we let $\mathsf t(H,U)$ be the minimal integer $N\geq 1$ such that, whenever $U_1,\ldots,U_k\in \mathcal A(H)$ with $U\mid U_1\cdot\ldots\cdot U_k$, then there exists $I\subseteq [1,k]$ and $V_2,\ldots,V_{r}\in \mathcal A(H)$ with $\prod_{i\in I}U_i=U\cdot V_2\cdot\ldots\cdot V_{r}$ and $\max\{|I|,\,r\}\leq N$. If no such $N$ exists, we set  $\mathsf t(H,U)=\infty$. The constant $\mathsf t(H,U)$ is called the \textbf{local tame degree}, and $H$ is called \textbf{locally tame} if $\mathsf t(H,U)<\infty$ for all $U\in \mathcal A(H)$ \cite{alfredbook}.

\subsection{Asymptotic Notation}
Given sequences  $\{a_i\}_{i=1}^\infty$ and $\{b_i\}_{i=1}^\infty$ of positive real numbers $a_i,\,b_i\in \R_+$, we use the following notation for gauging  their comparative asymptotic growth: $a_i\in o(b_i)$ means $a_i/b_i\rightarrow 0$; $a_i\in O(b_i)$ means $\{a_i/b_i\}_{i=1}^\infty$ is a bounded sequence, and thus $a_i/b_i\rightarrow M$ for some real number $M\geq 0$ by passing to appropriate subsequences; $a_i\in \Theta(b_i)$ means there exist \emph{positive} real numbers $\alpha,\beta\in \R_+$ such that $\alpha b_i\leq a_i\leq \beta b_i$ for all $i$, i.e., $a_i\in O(b_i)$ and $b_i\in O(a_i)$, in turn ensuring $a_i/b_i\rightarrow M$ for some positive real number $M>0$ by passing to appropriate subsequences; and $a_i\sim b_i$ means $a_i/b_i\rightarrow 1$. We sometimes use $o(b_i)$ to represent some existent sequence $\{c_i\}_{i=1}^\infty$ with $|c_i|\in o(b_i)$, and likewise with the other asymptotic notation. In general, with only a few clear exceptions, all asymptotics will be with regard to the variable $i$ in this work.


\section{Asymptotically Filtered Sequences, Encasement and Boundedness}\label{sec-asym-seq}

\subsection{Asymptotically Filtered Sequences}\label{sec-asym-seq-1}
Given a subset $X\subseteq \R^d$, we will need a precise way to describe the directions in which $X$ escapes to infinity, allowing us to characterize when $X$ remains within bounded distance of another subset $Y\subseteq \R^d$. As a first approximation, we say that a sequence $\{x_i\}_{i=1}^\infty$ of nonzero points $x_i\in \R^d$ is \textbf{radially convergent} with \textbf{limit} $u$, where $u\in \R^d$ is a unit vector, if
$\lim_{i\rightarrow \infty}\|x_i\|\in \R_+\cup \{\infty\}$ exists and
the sequence of unit vectors $\{x_i/\|x_i\|\}_{i=1}^\infty$ is convergent with limit $\lim_{i\rightarrow \infty}x_i/\|x_i\|= u$. Since the unit sphere is a compact metric space, any sequence of nonzero points contains a radially convergent subsequence.

\begin{definition}For $X\subseteq \R^d$, let $X^\infty$ denote all unit vectors which are a limit of an unbounded radially convergent sequence of terms from $X$.\end{definition}

The set $X^\infty$ is  a crude notion of the ``directions'' in which the set $X$ escapes to infinity satisfying the following closure property.

\begin{lemma}\label{lemma-radpoints-closed} Let $X\subseteq \R^d$, where $d\geq 0$. Then $X^\infty$ is  a closed subset of the unit sphere in $\R^d$.
\end{lemma}

\begin{proof}For $d=0$, $X^\infty$ is empty, and there is nothing to show, so assume $d\geq 1$.
To show $X^\infty$ is closed, it suffices to show it contains all its limit points. To this end, let $\{u_i\}_{i=1}^\infty$ be a sequence of terms $u_i\in X^\infty$ with $u_i\rightarrow u$. Then $u$ is another unit vector. By passing to a subsequence, we can w.l.o.g. assume $\mathsf d(u_j,u)<1/2^j$ for all $j$.
For each $u_j\in X^\infty$, there exists a radially convergent sequence $\{x_{ij}\}_{i=1}^\infty$ with $x_{ij}\in X$, \  $\|x_{ij}\|\rightarrow \infty$, and  $\lim_{i\rightarrow \infty} x_{ij}/\|x_{ij}\|=u_j$. Thus we may take $y_j=x_{ij}$ for some fixed sufficiently large $i$ such that $\|y_j\|>2^j$ (possible as $\|x_{ij}\|\rightarrow \infty$) and $\mathsf d(y_j/\|y_j\|,u_j)<1/2^j$. Now consider the sequence $\{y_j\}_{j=1}^\infty$ of terms $y_j\in X$. By construction, $\|y_j\|\rightarrow \infty$ while the triangle inequality ensures that $\mathsf d(y_j/\|y_j\|,u)\leq \mathsf d(y_j/\|y_j\|,u_j)+\mathsf d(u_j,u)<1/2^j+1/2^j=1/2^{j-1}$, which tends to $0$ as $j\rightarrow \infty$. Thus $\{y_j\}_{j=1}^\infty$ is an unbounded radially convergent sequence of terms from $X$ with limit $u$, showing that $u\in X^\infty$. As $u$ was an arbitrary limit point, it follows that $X^\infty$ is closed.
\end{proof}

 While the set $X^\infty$ provides a crude notion of the directions in which $X$ escapes to infinity, it will turn out to  be insufficient for our needs, leading us to the following more refined notion.
\begin{definition}
Let  $\vec u=(u_1,\ldots,u_t)$ be a tuple of $t\geq 0$ orthonormal vectors in $\R^d$. A sequence $\{x_i\}_{i=1}^\infty$ of elements $x_i\in \R^d$ is called an \textbf{asymptotically filtered sequence} with \textbf{limit} $\vec u$ if
$$x_i= a_i^{(1)}u_1+\ldots+a^{(t)}_iu_{t}+y_i\quad\mbox{ for all $i\geq 1$,}$$
for some real numbers   $a_i^{(j)}>0$  and vectors $y_i\in \R\la u_1,\ldots,u_t\ra^\bot$ such that
\begin{itemize}
\item $\lim_{i\rightarrow \infty} a_i^{(j)}\in \R_+\cup \{\infty\}$ exists for each $j\in [1,t]$,
\item  $\|y_i\|\in o(a_i^{(t)})$, and $a_i^{(j+1)}\in o(a_i^{(j)})$  for all $j\in [1,t-1]$.
\end{itemize}
\end{definition}
The first bulleted condition above is added mostly for convenience and is not essential to the definition.
 We say that the limit $\vec u=(u_1,\ldots,u_t)$ is \textbf{fully unbounded} if $a_i^{(t)}\rightarrow \infty$ (and thus $a_i^{(j)}\rightarrow \infty$ for all $j\in [1,t]$). Note this requires $t\geq 1$. The empty tuple, corresponding to when $t=0$, is referred to as the \textbf{trivial} tuple. We call the limit $\vec u=(u_1,\ldots,u_t)$ \textbf{anchored} if $t=0$ or $\{a_i^{(t)}\}_{i=1}^\infty$ is instead a bounded sequence. The limit $\vec u$ is \textbf{complete} if $y_i=0$ for all $i$, and $\vec u$ is a \textbf{complete fully unbounded} limit if $\vec u$ is fully unbounded but $\{y_i\}_{i=1}^\infty$ is bounded.
 Given any limit $\vec u=(u_1,\ldots,u_t)$ and $j\leq t$, we refer to $(u_1,\ldots,u_j)$ as a \textbf{truncation} of $\vec u$, which is \textbf{strict} when $j<t$. For $t\geq 1$, we let $$\vec u^\triangleleft=(u_1,\ldots,u_{t-1})$$ denote the principal truncation of $\vec u=(u_1,\ldots,u_t)$. The choice of using tuples of orthonormal vectors to represent the limit of an asymptotically filtered sequence is purely a matter of canonical representation. Indeed, the tuple $(u_1,\ldots,u_t)$ really represents an ascending chain of half-spaces $\{0\}\subset \R^\circ_+u_1\subset \R\la u_1\ra+\R^\circ_+u_2\subset \R\la u_1,u_2\ra+\R_+^\circ u_3\subset\ldots\subset \R\la u_1,\ldots,u_{t-1}\ra+\R_+^\circ u_t$ (or more compactly, the set $\bigcup_{j=0}^t(\R\la u_1,\ldots,u_{j-1}\ra+\R_+^\circ u_j)$, where $u_0=0$, from which this chain can be recovered), and any tuple $\vec v=(v_1,\ldots,v_t)$ of vectors $v_1,\ldots,v_t\in \R^d$ such that $\R\la v_1\ldots,v_{j-1}\ra+\R_+^\circ v_j=\R\la u_1\ldots,u_{j-1}\ra+\R_+^\circ u_j$ for all $j\in [1,t]$ is considered an \textbf{equivalent tuple}.
 \begin{proposition}\label{prop-filt-char}
 Let $\vec u=(u_1,\ldots,u_t)$ be a tuple of $t\geq 0$ orthonormal vectors in $\R^d$ and let $\{x_i\}_{i=1}^\infty$ be a sequence of elements $x_i\in\R^d$ with $$x_i=a_i^{(1)}u_1+\ldots+a_i^{(t)}u_t+y_i \quad\mbox{ for all $i\geq 1$},$$ for some real numbers   $a_i^{(j)}>0$  and vectors $y_i\in \R\la u_1,\ldots,u_t\ra^\bot$. For $j\in [0,t]$, let $\pi_j:\R^d\rightarrow \R\la u_1,\ldots,u_j\ra^\bot$ denote the orthogonal projection. Then $\{x_i\}_{i=1}^\infty$ is an asymptotically filtered sequence with limit $\vec u$ if and only if $\{\pi_{j-1}(x_i)\}_{i=1}^\infty$ is a radially convergent sequence with limit $u_j$  for all $j\in [1,t]$. Moreover, when this is the case, we have $\lim_{i\rightarrow\infty}\|\pi_{j-1}(x_i)\|=\lim_{i\rightarrow \infty}a_i^{(j)}$ for every $j\in [1,t]$.
\end{proposition}

\begin{proof}
For $j\in [0,t]$, let $y_i^{(j)}=a_i^{(j+1)}u_{j+1}+\ldots+a_i^{(t)}u_t+y_i$ and observe that  $$\pi_{j-1}(x_i)=a_i^{(j)}u_j+y_i^{(j)}=y_i^{(j-1)}$$ for $j\in [1,t]$. If $\{x_i\}_{i=1}^\infty$ is asymptotically filtered with limit $\vec u$, then $\|y_i^{(j)}\|\in o(a_i^{(j)})$ for $j\in [1,t]$.
In view of the triangle inequality and  $\|y_i^{(j)}\|\in o(a_i^{(j)})$, we have   $\|a_i^{(j)}u_j+y_i^{(j)}\|=a_i^{(j)}+ o(a_i^{(j)})$.
Thus $\|y_i^{(j)}\|\in o(a_i^{(j)})$ further implies that $$\|y_i^{(j)}\|/\|a_i^{(j)}u_j+y_i^{(j)}\|= \|y_i^{(j)}\|/(a_i^{(j)}+o(a_i^{(j)}))=(\|y_i^{(j)}\|/a_i^{(j)})/(1+o(1))
 \rightarrow 0,$$ ensuring that $y_i^{(j)}/\|a_i^{(j)}u_j+y_i^{(j)}\|\rightarrow 0$ and $$\lim_{i\rightarrow \infty}\pi_{j-1}(x_i)/\|\pi_{j-1}(x_i)\|
 =\lim_{i\rightarrow \infty}(a_i^{(j)}u_j+y_i^{(j)})/\|a_i^{(j)}u_j+y_i^{(j)}\|=\lim_{i\rightarrow \infty}a_i^{(j)}u_j/\|a_i^{(j)}u_j+y_i^{(j)}\|=u_j.$$ Also, $\lim_{i\rightarrow \infty}\|\pi_{j-1}(x_i)\|=\lim_{i\rightarrow \infty}\|a_i^{(j)}u_j+y_i^{(j)}\|=\lim_{i\rightarrow \infty}(a_i^{(j)}+o(a_i^{(j)}))=\lim_{i\rightarrow \infty}a_i^{(j)}\in \R_+\cup \{\infty\}$ exists, showing that $\{\pi_{j-1}(x_i)\}_{i=1}^\infty$ is radially convergent with limit $u_j$.

Next assume that $\{\pi_{j-1}(x_i)\}_{i=1}^\infty$ is radially convergent with limit $u_j$ for each $j\in [1,t]$, so $(a_i^{(j)}u_j+y_i^{(j)})/\|a_i^{(j)}u_j+y_i^{(j)}\|\rightarrow u_j$ for all $j\in [1,t]$. In particular, since $u_j$ and $y_i^{(j)}$ are linearly independent, $\|y_i^{(j)}\|/\|a_i^{(j)}u_j+y_i^{(j)}\|\rightarrow 0$ and $\|a_i^{(j)}u_j+y_i^{(j)}\|\rightarrow a_i^{(j)}$. Thus
$$0=\lim_{i\rightarrow \infty} \|y_i^{(j)}\|/\|a_i^{(j)}+y_i^{(j)}\|=\lim_{i\rightarrow \infty} \|y_i^{(j)}\|/a_i^{(j)},$$ ensuring that $\|y_i^{(j)}\|\in o(a_i^{(j)})$ for $j\in [1,t]$. In particular, $\|y_i\|=\|y_i^{(t)}\|\in o(a_i^{(t)})$, while $$\|y_i^{(j)}\|=\|a_i^{(j+1)}u_{j+1}+y_i^{(j+1)}\|\in a_i^{(j+1)}+O(\|y_i^{(j+1)}\|)\subseteq
a_i^{(j+1)}+o(a_i^{(j+1)})$$  for $j\in [1,t-1]$. Thus $\|y_i^{(j)}\|\in o(a_i^{(j)})$ for $j\in [1,t]$ ensures that $a_i^{(j+1)}\in o(a_i^{(j)})$ for $j\in [1,t-1]$. Furthermore, $\lim_{i\rightarrow \infty} a_i^{(j)}=\lim_{i\rightarrow \infty}(a_i^{(j)}+o(a_i^{(j)}))=\lim_{i\rightarrow \infty} (\|a_i^{(j)}u_j+y_i^{(j)}\|)=\lim_{i\rightarrow \infty} \|\pi_{j-1}(x_i)\|\in \R_+\cup \{\infty\}$ exists as  $\{\pi_{j-1}(x_i)\}_{i=1}^\infty$ is radially convergent, completing the proof.
\end{proof}

Suppose $\{x_i\}_{i=1}^\infty$ is a sequence of elements $x_i\in \R^d$ that is not eventually the constant zero sequence. By passing to a subsequence, we can assume all $x_i$ are nonzero with $\{x_i\}_{i=1}^\infty$ a radially convergent sequence with limit (say) $u_1$.  Moreover, if $\{x_i\}_{i=1}^\infty$ is unbounded, then we can assume $\lim_{i\rightarrow \infty} \|x_i\|=\infty$. Write each $x_i=a_i^{(1)}u_1+y_i^{(1)}$ with $y_i^{(1)}\in \R\la u_1\ra^\bot$.  Since $x_i/\|x_i\|\rightarrow u_1$, we can assume $a_i^{(1)}>0$ for all $i$ by discarding the first few terms.
  Proposition \ref{prop-filt-char} implies that the  resulting subsequence
is asymptotically filtered with limit $u_1$ and $\lim_{i\rightarrow \infty}a^{(1)}_i=\lim_{i\rightarrow
  \infty}\|x_i\|$. If $\{y_i^{(1)}\}_{i=1}^\infty$ is eventually zero, then discarding the first few terms allows us to assume $x_i=a_i^{(1)}u_1$ for all $i\geq 1$. Otherwise, we can repeat the above procedure and, passing to a yet more refined subsequence, conclude that  $x_i=a_i^{(1)}u_1+a_i^{(2)}u_2+y_i^{(2)}$ is asymptotically filtered with limit $(u_1,u_2)$. Continuing to iterate the procedure, we find that any sequence of terms from $\R^d$ contains an asymptotically filtered subsequence with complete limit. Likewise, truncating appropriately,  any unbounded sequence in $\R^d$ must contain an asymptotically filtered subsequence with complete fully unbounded limit. Note Proposition \ref{prop-filt-char} also ensures that if $\{x_i\}_{i=1}^\infty$ is an asymptotically filtered sequence both with limit $\vec u$ and with limit $\vec v$, then either $\vec v$ is a truncation of $\vec u$ or $\vec u$ is a truncation of $\vec v$.

  The following proposition is routine but important and requires the following notation. Let $\vec u=(u_1,\ldots,u_t)$ be a tuple of orthonormal vectors $u_i\in \R^d$, let $\mathcal E\subseteq \R^d$ be a subspace,  and let  $\pi:\R^d\rightarrow \mathcal E^\bot$ be the orthogonal projection. Then $$\pi(\vec u):=(\overline u_1,\ldots,\overline u_\ell),$$ where the $\overline u_i$ are defined  as follows. Recursively define indices \be\label{tuple-pi}0=r_0<r_1<\ldots<r_\ell<r_{\ell+1}=t+1\ee by letting $r_j\in [r_{j-1}+1,t]$ (for $j\in [1,\ell]$) be the minimal index such that $\pi_{j-1}(u_{r_j})\neq 0$, where $\pi_{j-1}:\R^d\rightarrow (\mathcal E+\R\la u_{1},u_2\ldots,u_{r_{j-1}}\ra)^\bot$ is the orthogonal projection and $\ell\in [0,t]$ is the first index such that $\pi_\ell(u_i)=0$ for all $i$.  In particular, $\pi_0=\pi$ and  $r_1\in [1,t]$ is the first index such that $\pi(u_{r_1})\neq 0$ (unless $\ell=0$, in which case no index with this property exists). Then $\overline u_j:=\pi_{j-1}(u_{r_j})/\|\pi_{j-1}(u_{r_j})\|$, ensuring that $\pi(\vec u)$ is a tuple of orthonormal vectors from $\mathcal E^\bot$. Equivalently, the indices $r_j$ are those with   $\R\la \pi(u_1),\ldots,\pi(u_{r_j})\ra\neq \R\la \pi(u_1),\ldots,\pi(u_{r_j-1})\ra$, and  $\pi(\vec u)=(\overline u_1,\ldots,\overline u_\ell)$ is simply the canonical tuple of orthonormal vectors equivalent to the tuple $(\pi(u_{r_1}),\ldots,\pi(u_{r_\ell}))$. From this viewpoint, it is clear that, if $\pi':\R^d\rightarrow (\mathcal E')^\bot$ is an orthogonal projection with $\mathcal E\subseteq \mathcal E'$, then $\pi'(\vec u)=\pi'(\pi(\vec u))$.

\begin{proposition}\label{prop-infinite-limits-proj}
Let $X\subseteq \R^d$, where $d\geq 0$, let $\mathcal E\subseteq \R^d$ be a subspace, and let $\pi:\R^d\rightarrow \mathcal E^\bot$ be the orthogonal projection.
\begin{itemize}
\item[1.] Let $\{x_i\}_{i=1}^\infty$ be an asymptotically  filtered  sequence of terms $x_i\in X$ with limit $\vec{u}=(u_1,\ldots,u_t)$, say $x_i=a_i^{(1)}u_1+\ldots+a_i^{(t)}u_t+y_i$. Then the sufficiently large index terms in $\{\pi(x_i)\}_{i=1}^\infty$ form an asymptotically filtered  sequence with limit $\pi(\vec{u})=(\overline u_1,\ldots, \overline u_\ell)$ with  $$\pi(x_i)=b_i^{(1)}\overline u_{1}+\ldots+b_i^{(\ell)}\overline u_{\ell}+y'_i,$$   $b_i^{(j)}\in \Theta( a_i^{(r_j)})$ for all $j\in [1,\ell]$, and $\|y'_i\|\in O(\|y_i\|)$, where the indices $r_1<\ldots <r_\ell$ are those given by \eqref{tuple-pi}.
\item[2.] If $\{y_i\}_{i=1}^\infty$ is an asymptotically filtered sequence of terms $y_i\in \pi(X)$ with limit $\vec{v}$, then there is an asymptotically  filtered sequence $\{x_i\}_{i=1}^\infty$ of terms $x_i\in X$ with limit $\vec{u}$, such that, replacing $\{y_i\}_{i=1}^\infty$ with an appropriate subsequence, we have $\pi(x_i)=y_i$ for all $i$, \  $\pi(\vec{u})=\vec{v}$, and $\pi(\vec u^\triangleleft)=\pi(\vec u)^\triangleleft=\vec v^\triangleleft$ (if $\vec v$ is nontrivial).
\end{itemize}

\end{proposition}

\begin{proof}
1. Let the $r_j$ and $\pi_j$ be as given in the definition of $\pi(\vec u)$. For $j\in [1,\ell]$,
 we have $\pi_{j-1}(x_i)=a_i^{(1)}\pi_{j-1}(u_1)+\ldots+a_i^{(t)}\pi_{j-1}(u_t)+\pi_{j-1}(y_i)$. Since $a_i^{(j)}\in o(a_i^{(j-1)})$ for all $j\geq 2$ and $\|\pi_{j-1}(y_i)\|\in O(\|y_i\|)\subseteq o(a_i^{(t)})$ (the first inclusion follows as any linear transformation $\pi_{j-1}:\R^d\rightarrow \R^d$ is a bounded linear operator), it follows that
$\|\pi_{j-1}(x_i)\|\sim a_i^{(r_j)}\|\pi_{j-1}(u_{r_j})\|$. Now $\pi(x_i)=b_i^{(1)}\overline u_{1}+\ldots+b_i^{(\ell)}\overline u_{\ell}+y'_i$ with $$b_i^{(j)}=a_i^{(r_j)}\|\pi_{j-1}(u_{r_j})\|\pm a_i^{(r_j+1)}\|\pi_{j}^\bot(u_{r_j+1})\|\pm\ldots\pm
a_i^{(t)}\|\pi_{j}^\bot(u_{t})\|\pm\|\pi^\bot_{j}(y_i)\|$$ and $y'_i=\pi_\ell(y_i)$, where $\pi_{j}^\bot:\R^d\rightarrow \R \overline u_{j}$ is the orthogonal projection onto $\R \overline u_{j}$. For each $j\in [1,\ell]$, we have  $a_i^{(r_j)}\|\pi_{j-1}(u_{r_j})\|>0$ for all $i$,   $\|\pi_j^\bot(y_i)||\in O(\|y_i\|)\subseteq o(a_i^{(r_j)})$ and  $a_i^{(k)}\in o(a_i^{(r_j)})$ for all $k>r_j$. Thus $b^{(j)}_i\in \Theta(a_i^{(r_j)})$ with $b_i^{(j)}>0$ for all sufficiently large $i$. Also, $\|y'_i\|=\|\pi_\ell(y_i)\|\in O(\|y_i\|)\subseteq o(a_i^{(r_\ell)})=o(b_i^{(\ell)})$, and Item 1  now follows.

2.  Since each $y_i\in \pi(G_0)$, there exists some $x_i\in G_0$ such that $\pi(x_i)=y_i$, for $i\geq 1$. By passing to a subsequence, we can assume $\{x_i\}_{i=1}^\infty$ is asymptotically filtered with complete limit $\vec{u}'$, in which case Item 1 implies that the sufficiently large index terms in the sequence $\{\pi(x_i)\}_{i=1}^\infty=\{y_i\}_{i=1}^\infty$ are asymptotically filtered with complete limit $\pi(\vec{u}')$. Since $\{y_i\}_{i=1}^\infty$ is also asymptotically filtered with limit $\vec{v}$, we conclude that $\vec{v}$ is a truncation of $\pi(\vec{u}')$ (cf. Proposition \ref{prop-filt-char}), whence $\pi(\vec{u})=\vec{v}$ for an appropriate truncation $\vec{u}$ of $\vec{u}'$. Moreover, if $\vec v$ is nontrivial, then so is $\vec u$, and choosing a truncation $\vec u=(u_1,\ldots,u_t)$ with $t\geq 1$ minimal such that $\pi(\vec u)=\vec v$, we obtain that $\pi(\vec u^\triangleleft)=\pi(\vec u)^\triangleleft$. We may consider $\{x_i\}_{i=1}^\infty$ as an asymptotically  filtered sequence with the truncated limit $\vec{u}$, and Item 2 now follows.
\end{proof}

\begin{definition} Given a set $X\subseteq \R^d$, we let $X^{\mathsf{lim}}$ denote the set of all fully unbounded limits $\vec u=(u_1,\ldots,u_t)$ of an asymptotically filtered sequence of terms from $X$.
\end{definition}

Each element of $X^\infty$ occurs as a singleton tuple in $X^{\mathsf{lim}}$, and we view the fully unbounded limits from $X^{\mathsf{lim}}$ as the more complete set of ``directions'' in which the set $X$ escapes to infinity.

\subsection{Encasement and Boundedness}
  The following definition associates a family of cones that ``hug''  the boundary of each  potential direction from $X^{\mathsf{lim}}$.
\begin{definition}
Let $\vec u=(u_1,\ldots,u_t)$ be a tuple of linearly independent vectors $u_1,\ldots,u_t\in \R^d$, where $t\in [0,d]$.
A cone $C$ \textbf{encases} $\vec u$ if, for each $j\in [1,t]$, there are
 $\alpha_{1,j},\ldots,\alpha_{j,j}\in \R_+$ with $$\alpha_{j,j}>0,\quad z_j:= \alpha_{j,j}u_j+\alpha_{j-1,j}u_{j-1}+\ldots+\alpha_{1,j}u_1\quad \und\quad \C(z_1,\ldots,z_t)\subseteq C.$$ In the context of convex cones, we will often also say a subset  $X\subseteq \R^d$ \textbf{encases} $\vec u$ when $\C(X)$ does, and that   $X$ \textbf{minimally encases} $\vec u$ if, additionally, no proper subset of $X$ encases $\vec u$.
\end{definition}


 When $t=1$, we often speak of encasing the element $u_1$ rather than the tuple  $(u_1)$.
Removing the requirement that  $\alpha_{1,j},\ldots,\alpha_{j-1,j}\in \R_+$ in the definition of encasement results in a seemingly weaker definition that is nonetheless equivalent. Indeed, if $\C(y_1,\ldots,y_t)\subseteq C$ with each $y_j=\beta_{j,j}u_j+\beta_{j-1,j}u_{j-1}+\ldots+\beta_{1,j}u_1$ for some $\beta_{i,j}\in \R$ with $\beta_{j,j}>0$, then $y'_j=y_j+\max\{0,-\frac{\beta_{j-1,j}}{\beta_{j-1,j-1}}\} y_{j-1}\in \C(y_1,\ldots,y_t)$ is an element of the form $y'_j=\beta'_{j,j}u_j+\beta'_{j-1,j}u_{j-1}+\ldots+\beta'_{1,j}u_1$ for some $\beta'_{i,j}\in \R$ with $\beta'_{j,j}>0$ and $\beta'_{j-1,j}\geq 0$. Repeating this process sequentially for $y_{j-1},\ldots,y_1$ results in an element $z_j=\alpha_{j,j}u_j+\alpha_{j-1,j}u_{j-1}+\ldots+\alpha_{1,j}u_1\in \C(y_1,\ldots,y_t)\subseteq C$, where $\alpha_{i,j}\in \R_+$ and $\alpha_{j,j}>0$, showing that $C$ encases $(u_1,\ldots,u_t)$. Thus $C$ encasing $(u_1,\ldots,u_t)$ is equivalent to there existing a convex cone $C'\subseteq C$ that intersects $\R_+^\circ u_j+\R\la u_1,\ldots,u_{j-1}\ra$ for each $j\in [1,t]$. In particular, if $C$ encases $\vec u$, then $C$ encases all equivalent tuples to $\vec u$ too.

\begin{lemma}\label{lem-encasement-interseciton}
Let $(u_1,\ldots,u_t)$ be a tuple of linearly independent vectors $u_1,\ldots,u_t\in \R^d$, where $t\in [0,d]$. If  the cones  $C,\,C'\subseteq \R^d$ both encase  $(u_1,\ldots,u_t)$, then $C\cap C'$ encases  $(u_1,\ldots,u_t)$.
\end{lemma}

\begin{proof}Since $C,\,C'\subseteq \R^d$ each encase  $(u_1,\ldots,u_t)$, there are subsets $X=\{x_1,\ldots,x_t\}\subseteq C$ and $Y=\{y_1,\ldots,y_t\}\subseteq C'$ with $\C(X)\subseteq C$ and  $\C(Y)\subseteq C'$ such that, for each $j\in [1,t]$,
$x_j=\alpha_{1,j}u_1+\ldots+\alpha_{j,j}u_j\in \C(X)$ and $y_j=\beta_{1,j}u_1+\ldots+\beta_{j,j}u_j\in \C(Y)$ for some $\alpha_{i,j},\,\beta_{i,j}\geq 0$ with $\alpha_{j,j},\,\beta_{j,j}>0$. Let $A=(\alpha_{i,j})_{i,j}$ be the upper $t\times t$ triangular matrix given by the $\alpha_{i,j}$ with $i,\,j\in [1,t]$, so $\alpha_{i,j}=0$ whenever $i>j$, and likewise define the upper $t\times t$ triangular matrix $B=(\beta_{i,j})_{i,j}$.
For $j\in [1,t]$, we aim to construct $z_j=\gamma_{1,j}u_1+\ldots+\gamma_{j,j}u_j\in \C(X)\cap \C(Y)\subseteq C\cap C'$ with $\gamma_{i,j}\geq 0$ for all $1\leq i\leq j\leq t$ and each $\gamma_{j,j}>0$.

Fix $j\in [1,t]$ arbitrary. Since the $u_i$ are linearly independent,  $z_j=\gamma_{1,j}u_1+\ldots+\gamma_{j,j}u_j$  lies in $\C(X)$ precisely when there exists a vector $x=(r_1,\ldots,r_j,0,\ldots,0)$ with all entries non-negative such that $Ax=y$, where $y=(\gamma_{1,j},\ldots,\gamma_{j,j},0,\ldots,0)$. By well known back substitution formulas, \be\label{backsub-1}r_i=\frac{\gamma_{i,j}-\sum_{k=i+1}^{j}\alpha_{i,k}r_k}{\alpha_{i,i}}\quad\mbox{ for $i\in [1,j]$}.\ee
Likewise,  $z_j=\gamma_{1,j}u_1+\ldots+\gamma_{j,j}u_j$  lies in $\C(Y)$ precisely when there is an $x'=(r'_1,\ldots,r'_j,0,\ldots,0)$ with all entries non-negative such that $Ax'=y$. As before, we have
\be\label{backsub-2}r'_i=\frac{\gamma_{i,j}-
\sum_{k=i+1}^{j}\beta_{i,k}r'_k}{\beta_{i,i}}\quad\mbox{ for $i\in [1,j]$}.\ee Taking $\gamma_{j,j}=1$, we observe that both \eqref{backsub-1} and \eqref{backsub-2} ensure that $r_j>0$ and $r'_j>0$. But now we can recursively construct the $\gamma_{i,j}\geq 0$ for $i=j,j-1,\ldots,1$ with $r_i,\,r'_i\geq 0$ by simply choosing $\gamma_{i,j}\geq 0$ sufficiently large. Indeed, taking $\gamma_{i,j}=\max\{\sum_{k=i+1}^{j}\alpha_{i,k}r_k,\,\sum_{k=i+1}^{j}\beta_{i,k}r'_k\}$ for $i\in [1,j-1]$ suffices, completing the proof.
\end{proof}

One of the main goals of this section is to give a local containment characterization of the more subtle notion of being bound to $Y$, which we will introduce in the next subsection. Finding a satisfying such characterization for a general set $Y\subseteq \R^d$ is more difficult, so we will instead restrict our attention  to a suitably broad class of subsets of $\R^d$, namely, the class consisting of all \emph{finite unions of polyhedral cones}. We formally allow the empty set to be considered a finite union of polyhedral cones, viewing it as an empty union. Since polyhedral cones are closed, a finite union of polyhedral cones is also closed.  This class has several useful closure properties, summarized in the following lemma.


\begin{lemma}\label{lem-polyhedral-cone}
 Suppose $X,\,X_1,\ldots,X_s\subseteq \R^d$  are finite unions of  polyhedral cones.
\begin{itemize}
\item[1.] $\bigcup_{i=1}^{s}X_i$ and
$\bigcap_{i=1}^{s}X_i$  are both  finite unions of  polyhedral cones.
\item[2.] If $V\subseteq \R^d$ is a subspace containing $\R\la X\ra$, then
 $\overline{V\setminus X}$ is a finite union of polyhedral cones  with $\partial(\overline{V\setminus X})\subseteq \partial(X)$.
 \item[3.] If  $C_1,\ldots,C_s\subseteq \R^d$ are polyhedral cones with $\R\la C_i\ra=V$ for all $i\in [1,s]$, then $\partial(X)=\partial(\overline{V\setminus X})$, where $X=\bigcup_{j=1}^s C_j$.
\end{itemize}
\end{lemma}

\begin{proof}
1. That $\bigcup_{i=1}^{s}X_i$ is a finite union of polyhedral cones is immediate, as this is the case for each $X_i$. That the intersection of a finite number of  polyhedral cones is itself a polyhedral cone  follows from the characterization of a polyhedral cone as the intersection of a finite number of half spaces. By hypothesis, each $X_j$ is a finite union of polyhedral cones, say $X_j=\bigcup_{i=1}^{t_j}C^{(j)}_i$. Thus $\bigcap_{j=1}^{s}X_j=\bigcap_{j=1}^s\bigcup_{i=1}^{t_j}C_i^{(j)}=\underset{i_1\in [1,t_1],\ldots,i_s\in [1,t_s]}{\bigcup}\bigcap_{j=1}^s C_{i_j}^{(j)}$ with each $\bigcap_{j=1}^s C_{i_j}^{(j)}$ a polyhedral cone as previously remarked.

2. By hypothesis, $X\subseteq V\subseteq \R^d$ is a finite union of polyhedral cones, say $X=\bigcup_{j=1}^sC_j$ with each $C_j$ a polyhedral cone, so each $C_j=\bigcap_{i=1}^{t_j}H_{i,j}$ for some closed half spaces $H_{i,j}\subseteq V$. We may w.l.o.g. assume $V=\R^d$. Then \begin{align*}&V\setminus X=V\setminus (\bigcup_{j=1}^sC_j)=V\setminus (\bigcup_{j=1}^s\bigcap_{i=1}^{t_j}H_{i,j})=\bigcap_{j=1}^s(V\setminus \bigcap_{i=1}^{t_j}H_{i,j})=\bigcap_{j=1}^s\bigcup_{i=1}^{t_j}V\setminus H_{i,j}=\\&\bigcup_{i_1\in [1,t_1],\ldots,i_s\in [1,t_s]} \bigcap_{j=1}^s V\setminus H_{i_j,j}=\bigcup_{i_1\in [1,t_1],\ldots,i_s\in [1,t_s]} \bigcap_{j=1}^s \mathsf{int}(-H_{i_j,j}),\end{align*} where each $\mathsf{int}(-H_{i_j,j})$ is an open half space in $V$. Thus \begin{align}\nn&\overline{V\setminus X}=\overline{\bigcup_{i_1\in [1,t_1],\ldots,i_s\in [1,t_s]} \bigcap_{j=1}^s \mathsf{int}(-H_{i_j,j})}=
\bigcup_{i_1\in [1,t_1],\ldots,i_s\in [1,t_s]} \overline{\bigcap_{j=1}^s \mathsf{int}(-H_{i_j,j})}\\&=\bigcup_{i_1\in [1,t_1],\ldots,i_s\in [1,t_s]} \overline{\mathsf{int}\Big(\bigcap_{j=1}^s-H_{i_j,j}\Big)}.\label{gasping}
\end{align}
Each $\bigcap_{j=1}^s-H_{i_j,j}$ is a finite intersection of closed half spaces, and thus a polyhedral cone. In particular, $\bigcap_{j=1}^s-H_{i_j,j}$ is convex.
If $\mathsf{int}\Big(\bigcap_{j=1}^s-H_{i_j,j}\Big)\neq \emptyset$, then $\mathsf{int}\Big(\bigcap_{j=1}^s-H_{i_j,j}\Big)=\Big(\bigcap_{j=1}^s-H_{i_j,j}\Big)^\circ$ (as any nonempty open set contains a ball of positive radius, which is a full dimensional subset), whence $$\overline{\mathsf{int}\Big(\bigcap_{j=1}^s-H_{i_j,j}\Big)}=
\overline{\Big(\bigcap_{j=1}^s-H_{i_j,j}\Big)^\circ}=
\overline{\bigcap_{j=1}^s-H_{i_j,j}}=\bigcap_{j=1}^s-H_{i_j,j},$$ with the penultimate equality above in view of the convexity of $\bigcap_{j=1}^s-H_{i_j,j}$, and the final equality in view of $\bigcap_{j=1}^s-H_{i_j,j}$ being an intersection of closed subsets, and thus itself closed. Thus $\overline{\mathsf{int}\Big(\bigcap_{j=1}^s-H_{i_j,j}\Big)}=\bigcap_{j=1}^s-H_{i_j,j}$ is a polyhedral cone in such case. On the other hand, if
$\mathsf{int}\Big(\bigcap_{j=1}^s-H_{i_j,j}\Big)= \emptyset$, then $\overline{\mathsf{int}\Big(\bigcap_{j=1}^s-H_{i_j,j}\Big)}=\emptyset$ too. As a result, \eqref{gasping} shows that $\overline{V\setminus X}$ is a finite union of polyhedral cones, as desired.

Let $\mathscr E=\R\la X\ra\subseteq V=\R^d$. If $V$ properly contains $\mathscr E$, then $\overline{V\setminus X}=\overline{V\setminus \mathscr E}=V$ and $\partial(\overline{V\setminus X})=\partial(V)=\emptyset\subseteq \partial(X)$. On the other hand, if $V=\mathscr E$, then  \begin{align*}&\partial(X)=\partial(V\setminus X)={(\overline{V\setminus X})\setminus (V\setminus X)}^\circ\quad\und\\&\partial(\overline{V\setminus X})=(\overline{\overline{V\setminus X}})\setminus {(\overline{V\setminus X})}^\circ=(\overline{V\setminus X})\setminus {(\overline{V\setminus X})}^\circ,\end{align*} where the equality $\partial(X)=\partial(V\setminus X)$ follows in view of $X$ being closed (so that $V\setminus X$ lying in a proper subspace of $V$ is only possible if $X=V$). Consequently, since ${(V\setminus X)}^\circ\subseteq {(\overline{V\setminus X})}^\circ$ (as both $V\setminus X$ and $\overline{V\setminus X}$ span the same subspace), it follows that $\partial(\overline{V\setminus X})\subseteq \partial(X)$. This completes the proof of Item 2.

3.
Each $C_j$ has full dimension in $V$, and we can w.l.o.g. assume $V=\R^d$.  Since $X$ is closed and spans the subspace $V$, it follows that $V\setminus X$ is open in $V$ and is thus either empty or has full dimension. If $V\setminus X=\emptyset$, then $\overline{V\setminus X}=\emptyset$ too, whence $\partial(\overline{V\setminus X})=\partial(X)=\partial(V)=\emptyset$, as desired. So instead assume $V\setminus X$ has full dimension, whence $V\setminus X$ and $X$ both span the entire space $V$. As a result,  $\partial(X)=\partial(V\setminus X)=(\overline{V\setminus X})\setminus {(V\setminus X)}^\circ$ and $\partial(\overline{V\setminus X})=(\overline{\overline{V\setminus X}})\setminus {(\overline{V\setminus X})}^\circ=(\overline{V\setminus X})\setminus {(\overline{V\setminus X})}^\circ$. Letting $Y=V\setminus X$, it remains to show ${\overline{Y}}^\circ=Y^\circ$. The inclusion $Y^\circ\subseteq {\overline{Y}}^\circ$ is trivial (as $Y$ and $\overline{Y}$ span the same subspace). Since $Y=V\setminus X$ is open in $V$, we have $Y^\circ =Y$. Thus we need to show ${\overline{Y}}^\circ\subseteq Y$. Let $x\in {\overline{Y}}^\circ$ be arbitrary. Then $x\in \overline{Y}$ and there exists an open ball $B_\epsilon(x)\subseteq \overline{Y}$ for some $\epsilon>0$. Assume by contradiction that $x\notin Y=V\setminus X$, whence $x\in X=\bigcup_{j=1}^sC_j$, so w.l.o.g. $x\in C_1$, which is a polyhedral cone spanning $V$ by hypothesis.
 Since $C_1$ is a convex set with $x\in C_1$, we may take a point $y\in C_1^\circ$ and consider the line segment between $y$ and $x$. As $C_1$ is convex, all points on this line segment (apart from $x$) lie in $C_1^\circ$, ensuring that
  $B_\epsilon(x)\cap C_1$ contains some point $x_0\in C_1^\circ$. Thus $x_0\in C_1^\circ \cap B_\epsilon(x)$. However, since $B_\epsilon(x)\subseteq \overline{Y}=\overline{V\setminus X}$, we also have $x_0\in \overline{V\setminus X}\subseteq \overline{V\setminus C_1}$.
   Now $C_1$ is a closed set which spans the space $V$, which implies $V\setminus C_1$ is open in $V$, and thus also spans $V$ (it is nonempty else its closure could not contain $x_0$). Since $C_1$ and $V\setminus C_1$ span the same subspace $V$, it follows that $\overline{C_1}\cap \overline{V\setminus C_1}=\partial(C_1)=\partial(V\setminus C_1)$. But now $x_0\in \overline{V\setminus C_1}\cap C_1^\circ\subseteq \overline{V\setminus C_1}\cap \overline{C_1}=\partial(C_1)=\overline{C_1}\setminus C_1^\circ$, which contradicts that $x_0\in C_1^\circ$, completing the proof.
\end{proof}

Next, we need a notion of what it means for a set $X$ to be bounded relative to another set $Y$, extending the notion of a bounded set (which is the case $Y=\{0\}$).

\begin{definition}
Let $X,\,Y\subseteq \R^d$ be subsets. We say that $X$ is \textbf{bound} to $Y$ if there is some real number $\epsilon>0$ such that every $x\in X$ is distance less than $\epsilon$ from $Y$, i.e., $X\subseteq Y+B_\epsilon(0)$.   A sequence $\{x_i\}_{i=1}^\infty$ is bound to $Y$ if the set $X=\{x_i:\; i\geq 1\}$ is bound to $Y$.
\end{definition}

We continue with a basic connection between boundedness and encasement.

\begin{proposition} \label{prop-encasementcones-contain-aprox-seq} Let $\{x_i\}_{i=1}^\infty$ be an asymptotically  filtered sequence  of terms  $x_i\in \R^d$ with limit $\vec u=(u_1,\ldots,u_t)$ and let $\pi:\R^d\rightarrow \R\la u_1,\ldots,u_t\ra$ be the orthogonal projection.
Suppose  the cone  $C\subseteq \R^d$ encases  $\vec u$.
 \begin{itemize}
 \item[1.] $\{\pi(x_i)\}_{i=1}^\infty$ is  an asymptotically filtered sequence in $\R\la u_1,\ldots,u_t\ra\cong \R^t$ with limit $\vec u$.
 \item[2.] $\pi(x_i)\in C$ for all sufficiently large $i$.
 \item[3.]  If $\vec u$ contains a complete fully unbounded limit of $\{x_i\}_{i=1}^\infty$ as a truncation, then the sequence $\{x_i\}_{i=1}^\infty$ is bound to $C$. Indeed, $\mathsf d(x_i,C)\leq \sup_i \|\pi_t(x_i)\|<\infty$ for all sufficiently large $i$, where $\pi_t:\R^d\rightarrow \R\la u_1,\ldots,u_t\ra^\bot$ is the orthogonal projection.
 \end{itemize}
\end{proposition}

\begin{proof} For $j\in [0,t]$, let $\pi_j:\R^d\rightarrow \R\la u_1,\ldots,u_j\ra^\bot$ be the orthogonal projection. Since the cone $C$ encases   $(u_1,\ldots,u_t)$, there is some  $Y=\{y_1,\ldots,y_t\}\subseteq C$  with $\C(Y)\subseteq C$ such that, for every $j\in [1,t]$,  $y_j=\alpha_{1,j}u_1+\alpha_{2,j}u_2+\ldots+\alpha_{j,j}u_j$ for some $\alpha_{i,j}\geq 0$ with all $\alpha_{j,j}>0$.  Item 1 follows immediately from the definitions.

2. In view of Item 1, we may w.l.o.g. assume $x_i\in \R\la u_1,\ldots,u_t\ra$ for all $i$, so that each $\pi(x_i)=x_i=a_i^{(1)}u_1+\ldots+a_i^{(t)}u_t$ for some $a_i^{(j)}>0$ with $a_i^{(j)}\in o(a_i^{(j-1)})$ for $j\geq 2$. Let $A=(\alpha_{i,j})_{i,j}$ be the $t\times t$ upper triangular matrix given by the $\alpha_{i,j}$ with $i,\,j\in [1,t]$, so $\alpha_{i,j}=0$ whenever $i>j$. Then $x_i\in \C(Y)$ precisely when there exists a vector $x=(r_i^{(1)},\ldots,r_i^{(t)})$ with all $r_i^{(j)}\geq 0$ such that $Ax=z_i$, where $z_i=(a_i^{(1)},\ldots,a_i^{(t)})$. As in the proof of Lemma \ref{lem-encasement-interseciton}, by well known back substitution formulas, \be\label{backsub-3}r_i^{(j)}=\frac{a_i^{(j)}-\sum_{k=j+1}^{t}\alpha_{j,k}r_i^{(k)}}{\alpha_{j,j}}
\quad\mbox{ for $j\in [1,t]$}.\ee
Consequently, since $a_i^{(k)}\in o(a_i^{(j)})$ for all $k>j$, we find that $r_i^{(j)}\in \Theta(a_i^{(j)})$ for all $j\in [1,t]$. In particular, $r_i^{(k)}\in o(a_i^{(j)})$ for $k\geq j+1$, whence \eqref{backsub-3} ensures that $r_i^{(j)}>0$ for all sufficiently large $i$. But this means  $x_i\in \C(Y)\subseteq C$ for all sufficiently large $i$, completing Item 2.

3. By Item 2,   $\pi(x_i)\in C$ for all sufficiently large $i$, say all $i\geq N$. Since $\{x_i\}_{i=1}^\infty$ is asymptotically filtered with limit $(u_1,\ldots,u_t)$ containing a complete fully unbounded limit of $\{x_i\}_{i=1}^\infty$,  it follows that $\{\pi_t(x_i)\}_{i=1}^\infty$ is bounded, implying $\sup_i\|\pi_t(x_i)\|<\infty$.
We have $x_i=\pi(x_i)+\pi_t(x_i)$, whence
$\mathsf d(x_i,\pi(x_i))=\|\pi_t(x_i)\|\leq \sup_i\|\pi_t(x_i)\|$ for all $i$. Thus, for $i\geq N$,  we see $\pi(x_i)\in C$ is a point with distance at most $\sup_i\|\pi_t(x_i)\|<\infty$ from $x_i$, showing that the sequence $\{x_i\}_{i=1}^\infty$ is bound to $C$ with the desired bounds.
\end{proof}

For finite unions of polyhedral cones, the next proposition provides the key link between asymptotically  filtered sequences and encasement, helping explain why we have restricted attention to this class of sets.

\begin{proposition}\label{prop-finite-union--convergence-encasement}
Let $Y\subseteq \R^d$ be a finite union of polyhedral cones. If $\{y_i\}_{i=1}^\infty$ is an asymptotically filtered  sequence of terms $y_i\in Y$ with limit $\vec u$, then $Y$ encases $\vec u$.
\end{proposition}

\begin{proof}Let $\vec u =(u_1,\ldots,u_t)$ and  extend $(u_1,\ldots,u_t)$ to an orthonormal basis $(u_1,\ldots,u_d)$ of $\R^d$.
Since $Y$ is a \emph{finite} union of polyhedral cones, by replacing $\{y_i\}_{i=1}^\infty$  with an appropriate subsequence, we can w.l.o.g. assume $y_i\in C\subseteq Y$ for all $i$ for some polyhedral cone $C$. Let $\pi:\R^d\rightarrow \R\la u_1,\ldots,u_t\ra$ be the orthogonal projection. Observing that $\C(u_1,\ldots,u_t)$ encases $\vec u$,  we can apply Proposition \ref{prop-encasementcones-contain-aprox-seq}.2 to conclude $\pi(y_i)\in \C(u_1,\ldots,u_t)$ for all sufficiently large $i$. Thus, by passing to a subsequence, we can w.l.o.g. assume all $y_i$ lie in the polyhedral cone $C'=\C(u_1,\ldots,u_t,\pm u_{t+1},\ldots,\pm u_{d})$. As a result, since $C\cap C'\subseteq Y$ is also a polyhedral cone, with $y_i\in C\cap C'$ for all $i$, we see that we can replace $C$ with $C\cap C'$ and thereby assume $C\subseteq \C(u_1,\ldots,u_t,\pm u_{t+1},\ldots,\pm u_{d})$.

 Since $C$ is a polyhedral cone, so too is $\pi_j(C)$ for each $j\in [0,t]$, where $\pi_j:\R^d\rightarrow \R\la u_1,\ldots,u_j\ra^\bot$ is the orthogonal projection.
 Since $\lim_{i\rightarrow \infty} \pi_{j-1}(x_i)/\|\pi_{j-1}(x_i)\|=u_{j}$, for $j\in [1,t]$, is a limit of points $\pi_{j-1}(x_i)/\|\pi_{j-1}(x_i)\|$ from the closed cone  $\pi_{j-1}(C)$, it follows that $u_{j}\in \pi_{j-1}(C)$ for all $j\in [1,t]$.
 Consequently, since $\ker \pi_{j-1}=\R\la u_1,\ldots,u_{j-1}\ra$, there must be some $z_j\in C$ with $\pi_{j-1}(z_j)=u_{j}$, meaning  $z_j=\alpha_{1,j}u_1+\alpha_{2,j}u_2+\ldots+\alpha_{j,j}u_j$ with  $\alpha_{i,j}\in \R$ and $\alpha_{j,j}=1$. However, since each $z_j\in C\subseteq \mathsf C(u_1,\ldots,u_t,\pm u_{t+1},\ldots,\pm u_{d})$, we must have $\alpha_{i,j}\geq 0$ for all $i$ and $j$. This shows that $C$ encases $\vec u$, and since $C\subseteq Y$, it follows that $Y$ encases  $\vec u$, as desired.
\end{proof}

\begin{lemma}\label{lem-weak-nearness-implication}
Let $\{x_i\}_{i=1}^\infty$ be an asymptotically filtered sequence of points in $\R^d$ with fully unbounded limit $\vec u$ and let $Y\subseteq \R^d$. If $\{x_i\}_{i=1}^\infty$ is bound to $Y$, then there is an asymptotically filtered  sequence $\{y_i\}_{i=1}^\infty$ of points $y_i\in Y$ having  fully unbounded limit $\vec u$.
\end{lemma}

\begin{proof} Let $\vec u=(u_1,\ldots,u_t)$ and let   each $x_i=a_i^{(1)}u_1+\ldots+a_i^{(t)}u_t+x'_i$ with $x'_i\in \R\la u_1,\ldots,u_t\ra^\bot$.
Since $\vec u$ is  fully unbounded, we have $a_i^{(j)}\rightarrow \infty$ for every $j\in [1,t]$. Since $\{x_i\}_{i=1}^\infty$ is bound to $Y$, there exists some finite $\epsilon>0$ such that, for  each term $x_i$, there is some $y_i\in Y$ with $\mathsf d(x_i,y_i)\leq \epsilon$.
Consider the sequence $\{y_i\}_{i=1}^\infty$. For each $i\geq 1$, write $y_i=x_i+z_i$ and write $z_i=b^{(1)}_iu_1+\ldots +b_i^{(t)}u_t+z'_i$ with $z'_i\in \R\la u_1,\ldots,u_t\ra^\bot$. Then $\{z_i\}_{i=1}^\infty$ is a bounded sequence in view of  $\mathsf d(x_i,y_i)\leq \epsilon$, ensuring that $\|z_i\|\in O(1)$, and thus also $b_i^{(j)}\in O(1)$ for $j\in [1,t]$ and $z'_i\in O(1)$. Since $a_i^{(t)}\rightarrow \infty$, we have $O(1)\subseteq o(a_i^{(t)})\subseteq o(a_i^{(j)})$ for all $j\in [1,t]$.
In consequence, since  $\|x'_i\|\in o(a_i^{(t)})$, it follows that $y_i=(a_i^{(1)}+b_i^{(1)})u_1+\ldots+(a_i^{(t)}+b_i^{(t)})u_t+(x'_i+z'_i)$ gives an asymptotically filtered representation of the $y_i$ having fully unbounded limit $\vec u$ (passing to sufficiently large index terms to ensure $a_i^{(j)}+b_i^{(j)}>0$), as desired.
\end{proof}

We now come to the main result of this section, giving a local containment characterization for a set to be bound to a finite union of polyhedral cones. We remark that, if $Y\subseteq \R^d$ is only a  cone, not a finite union of polyhedral cones, then the argument below shows $4.\Rightarrow 1.\Rightarrow 2.\Rightarrow 3$. It is the implication $3.\Rightarrow 4.$ that requires $Y$ to be a finite union of polyhedral cones.

\begin{theorem}\label{thm-nearness-characterization}
Let $X,\,Y\subseteq \R^d$ be subsets with $Y\neq \emptyset$ a finite union of polyhedral cones. Then the following are equivalent.
\begin{itemize}
\item[1.] $X$ is bound to $Y$.
\item[2.] Every asymptotically filtered sequence  $\{x_i\}_{i=1}^\infty$ of terms $x_i\in X$ with fully unbounded limit  is bound to $Y$.
\item[3.] $X^{\mathsf{lim}}\subseteq Y^{\mathsf{lim}}$.
\item[4.]  $Y$ encases  every $\vec u\in X^{\mathsf{lim}}$.
\end{itemize}
\end{theorem}

\begin{proof}The implication $1.\Rightarrow 2.$ is trivial. The implication $2.\Rightarrow 3.$  follows from Lemma \ref{lem-weak-nearness-implication}.
As $Y$ is a finite union of polyhedral cones, the implication  $3.\Rightarrow 4.$ follows from Proposition \ref{prop-finite-union--convergence-encasement}. It remains to establish the implication $4.\Rightarrow 1.$
To this end, assume by contradiction that  $Y$ encases  every $\vec u\in X^{\mathsf{lim}}$ but $X$ is not bound to $Y$. Then there exists a sequence $\{x_i\}_{i=1}^\infty$ of points $x_i\in X$ with $\mathsf d(x_i,Y)\rightarrow \infty$. In particular, $\{\|x_i\|\}_{i=1}^\infty$ is unbounded (as $Y\neq \emptyset$). Thus, as discussed at the beginning of the section, there exists a subsequence which is asymptotically filtered with complete fully unbounded limit. Replacing $\{x_i\}_{i=1}^\infty$ with this subsequence, we can w.l.o.g. assume $\{x_i\}_{i=1}^\infty$ is itself an asymptotically filtered sequence of terms  $x_i\in X$ with complete fully unbounded limit $\vec u=(u_1,\ldots,u_t)$ and $\mathsf d(x_i,Y)\rightarrow \infty$.
Since  $Y$ encases $\vec u=(u_1,\ldots,u_t)$ by hypothesis,  Proposition  \ref{prop-encasementcones-contain-aprox-seq}.3 implies that $\{x_i\}_{i=1}^\infty$ is bound to $Y$. However, this contradicts that $\mathsf d(x_i,Y)\rightarrow \infty$, completing the proof.
\end{proof}

\begin{corollary}\label{cor-nearness-intersection-prop}
Let $X,\, Y,\, Z\subseteq \R^d$ be subsets with $X$ and $Y$ each finite unions of polyhedral cones. If $Z$ is bound to $X$ as well as to $Y$, then $Z$ is bound to $X\cap Y$.
\end{corollary}

\begin{proof} By Lemma \ref{lem-polyhedral-cone}, $X\cap Y$ is also a finite union of polyhedral ones. Thus the corollary follows from  Theorem \ref{thm-nearness-characterization}.4 and Lemma \ref{lem-encasement-interseciton}.
\end{proof}

\section{Elementary Atoms, Positive Bases and Reay Systems}\label{sec-reay}

\subsection{Basic Non-degeneracy Characterizations}

We begin with a few basic properties regarding atoms and the representation of $0$ as a positive linear combination.

\begin{lemma}\label{lemma-rational-approx}
Let $\Lambda\leq \R^d$ be a full rank lattice, where $d\geq 0$, let $x_1,\ldots,x_r\in \Lambda$ and suppose $\alpha_1x_1+\ldots+\alpha_rx_r=0$ for some $\alpha_i\in \R$.
Then, for all $\epsilon>0$, there exists $\alpha'_1,\ldots,\alpha'_r\in \Q$ with $\alpha'_1x_1+\ldots+\alpha'_rx_r=0$ and  $|\alpha_i-\alpha'_i|<\epsilon$ for all $i$.
\end{lemma}

\begin{proof}Let $E=\{e_1,\ldots,e_d\}\subseteq \Lambda$ be a lattice basis for $\Lambda$. The lemma is vacuous for $r=0$ and trivial if $\alpha_i=0$ for all $i$, so assume $r\geq 1$ and that not all $\alpha_i$ are zero.
 Let $x$ be the nonzero column vector whose $i$-th entry is $\alpha_i$.  Then, letting $M$ be the $d\times r$ integer matrix whose $j$-th column is $x_j$ expressed using the basis $E$, we see that  $x$ lies in  the null space of $M$. Since the entries of $M$ are integers, the null space of $M$ is generated by integer vectors, say $y_1,\ldots,y_s\in \Z^r$ with $s\geq 1$. Thus $x=\beta_1y_1+\ldots+\beta_sy_s$ for some $\beta_i\in \R$. Since the rational numbers approximate the reals, we can find rational numbers $\beta'_1,\ldots,\beta'_s\in \Q$ with $|\beta_i-\beta'_i|$ sufficiently small so as to guarantee that, for every $i\in [1,r]$, the $i$-th coordinate $\alpha'_i$ of the vector $x'=\beta'_1y_1+\ldots+\beta'_sy_s$ satisfies $|\alpha'_i-\alpha_i|<\epsilon$. Moreover,  $\alpha'_1x_1+\ldots+\alpha'_rx_r=0$ (as $x'$ lies in the null space of $M$) with  $\alpha'_i\in \Q$ for all $i$ (as each $y_i$ is an integer  vector and each $\beta'_i\in \Q$), as desired.\end{proof}

\begin{proposition}\label{prop-rational-atoms}
Let $\Lambda\leq \R^d$ be a full rank lattice, where $d\geq 0$, and let  $G_0\subseteq \Lambda$. Then $$\mathcal A(G_0)\neq \emptyset\quad\mbox{ if and only if  }\quad 0\in \C^*(G_0).$$ Moreover, if there exist $x_1,\ldots,x_r\in G_0$, where $r\geq 1$, and real numbers $\alpha_i>0$ such that $\alpha_1x_1+\ldots+\alpha_rx_r=0$, then there is a zero-sum subsequence $S\in \mathcal B(G_0)$ with $\supp(S)=\{x_1,\ldots,x_r\}$.
\end{proposition}

\begin{proof}

If $U\in \mathcal A(G_0)$ is an atom, then  $\Summ{g\in G_0}\vp_g(U)g=\sigma(U)=0$ shows that $0\in   \mathsf C^*(G_0)$. On the other hand, if $0\in \mathsf C^*(G_0)$, then there are $x_1,\ldots,x_r\in G_0$, where $r\geq 1$, and real numbers $\alpha_i>0$ such that $\alpha_1x_1+\ldots+\alpha_rx_r=0$.  Applying Lemma \ref{lemma-rational-approx} with $\epsilon=\min_i \alpha_i>0$, we conclude that  $\alpha'_1x_1+\ldots+\alpha'_rx_r=0$  for some $\alpha'_i\in \Q$ with $\alpha'_i>0$ for all $i\in [1,r]$. By multiplying by a common denominator, we can w.l.o.g. assume the $\alpha'_i\in \Z$ for all $i$, and now  $S=x_1^{[\alpha'_1]}\bdot\ldots\bdot x_r^{[\alpha'_r]}\in \mathcal B(G_0)$ is a nontrivial zero-sum sequence with $\supp(S)=\{x_1,\ldots,x_r\}\subseteq G_0$, showing that $\mathcal A(G_0)\neq \emptyset$.
\end{proof}

\begin{proposition}\label{prop-no-zs-char}
A subset $X\subseteq \R^d$, where $d\geq 0$, has $0\notin \C^*(X)$ if
and only if $0\notin X$ and there exist a sequence of subspaces $\mathcal
H^0\subseteq \mathcal H^1\subseteq \ldots\subseteq \mathcal H^d$ with $\dim
\mathcal H^i=i$ such that $X\cap \mathcal H^j$ is contained in the closed
half space $\mathcal H^{j-1}_+\subseteq \mathcal H^{j}$ for $j=1,\ldots,d$.
\end{proposition}

\begin{proof}The case $d=0$ is trivial, so we assume $d\geq 1$.
A simple inductive argument shows that  any $X\subseteq \R^d$ satisfying the
stated conditions has $0\notin\C^*(X)$. On the other hand, if
$0\notin \C^*(X)$, then $0\notin X$ and $\C(X)\cap -\C(X)\neq \R^d$ (as $d\geq 1$). Thus $\C(X)\neq \R^d$, ensuring that  $\C(X)$, and thus also $X\subseteq \C(X)$,
must be contained in a closed half space of some subspace $\mathcal H^{d-1}$
with $\dim \mathcal H^{d-1}=d-1$. Repeating this argument for $X\cap \mathcal H^{d-1}$
and iterating then yields the desired result.
\end{proof}

When studying zero-sum subsequences, it natural to focus on those subsets $G_0$ such that every $g\in G_0$ is contained in some atom. Otherwise, we could simply pass to a subset of $G_0$ having this property. The next proposition characterizes such sets $G_0\subseteq \Lambda$.

\begin{proposition}\label{prop-notrivialG_0}
Let $\Lambda\leq \R^d$ be a full rank lattice, where $d\geq 0$, and let  $G_0\subseteq \Lambda$.  Then every $g\in G_0$ has some $U\in \mathcal A(G_0)$
with $g\in \supp(U)$ if and only if $\C(G_0)=\R\la G_0\ra$.
\end{proposition}

\begin{proof}We may w.l.o.g. assume $\R\la G_0\ra=\R^d$ (by passing to the lattice $\Z\la G_0\ra\leq \R\la G_0\ra$). If $d=0$, then $G_0\subseteq \{0\}$ and the result is clear. Therefore we may assume $d\geq
1$. Recall that $G_0$ is contained in a closed half-space if and only if $\C(G_0)\neq \R\la G_0\ra= \R^d$.
First assume every $g\in G_0$ has some $U\in \mathcal A(G_0)$ with $g\in
\supp(U)$.
 If $G_0$ were contained in a closed half space $\mathcal H_+$, then it is clear by reduction modulo $\mathcal H$ that no atom
can contain a term outside of $\mathcal H$. Thus, as each $g_0\in G_0$ lies
in some atom by hypothesis,  we must have $G_0\subseteq \mathcal H\neq
\R^d$, contradicting that $\R\la G_0\ra=\R^d$. It remains to prove the other direction, for which we proceed by
induction on $d$. The cases $d\leq 1$ are easily verified, so
we assume $d\geq 2$.

Let $G_0'\subseteq G_0$ be the subset of all $g\in G_0$ that are contained in
the support of some atom. Assuming $G_0'\neq G_0$, we need to show that $G_0$
is contained in a closed half space of $\R^d$. For this, we can
w.l.o.g. assume $0\in G_0$, and thus that $G_0'$ is nonempty. Let $\mathcal
E=\R\la G_0'\ra$. First suppose that $\dim \mathcal
E=0$, i.e., that $G'_0=\{0\}$. Then $\mathcal A(G_0\setminus \{0\})=\mathcal A(G_0\setminus G'_0)=\emptyset$, whence Propositions \ref{prop-rational-atoms} and \ref{prop-no-zs-char} imply that $G_0\setminus \{0\}$ is contained in a closed half space of $\R^d$, and thus so too $G_0$, as desired.
 So we now assume $\dim \mathcal E\geq 1$.

Since every $g\in G'_0$ is contained in some atom over $G_0$ but no element
in $G_0\setminus G'_0$ is contained in an atom over $G_0$, we actually have
that every $g\in G'_0$ is contained in some atom over $G'_0$. Thus by the
already established direction of the proof, we conclude that $\C(G'_0)=\mathcal E$.
Consequently, any sequence $T\in\Fc(G_0)$ with $\sigma(T)\in \mathcal E$ can
be completed to a zero-sum sequence $T'\in \mathcal B(G_0)$ with $\supp(T)\subseteq \supp(T')$ in view of Proposition  \ref{prop-rational-atoms}.
This means no element of $G_0\setminus G'_0\neq \emptyset$ is contained in a zero-sum
modulo $\mathcal E$. Thus, by the induction hypothesis, we conclude that
$\pi(G_0)$ is contained in a closed half space $\mathcal
H_+$ of ${\mathcal E}^\bot$, where $\pi:\R^d\rightarrow {\mathcal E}^\bot$ is the orthogonal projection, which implies $G_0$ is
contained in the closed half space $\pi^{-1}(\mathcal H_+)$ of $\R^d$, as
desired.
\end{proof}

As a corollary to Propositions \ref{prop-notrivialG_0} and \ref{prop-rational-atoms}, we have the following.

\begin{corollary}\label{cor-nondegen}
Let $\Lambda\leq \R^d$ be a full rank lattice, where $d\geq 0$, let  $G_0\subseteq \Lambda$, let $\mathcal E=\C(G_0)\cap -\C(G_0)$, and let $\wtilde G_0=\{g\in G_0:\; g\in \supp(U)\mbox{ for some $U\in \mathcal A(G_0)$}\}$. Then $\wtilde G_0=G_0\cap \mathcal E$ and $\C(\wtilde G_0)=\mathcal E$.
\end{corollary}

\begin{proof}
By definition of $\wtilde G_0$, any $U\in \mathcal A(G_0)$ must have $\supp(U)\subseteq \wtilde G_0$. Thus Proposition \ref{prop-notrivialG_0} applied to $\wtilde G_0$ (contained in the full rank lattice $\Lambda\cap \R\la \wtilde G_0\ra$ of $\R\la \wtilde G_0\ra$) implies that $$\C(\wtilde G_0)=\R\la \wtilde G_0\ra.$$ If $g\in G_0\cap \R\la \wtilde G_0\ra$, then $\C(\wtilde G_0\cup\{g\})=\C(\wtilde G_0)=\R\la \wtilde G_0\ra$, in which case Proposition \ref{prop-notrivialG_0} applied to $\wtilde G_0\cup \{g\}$ (contained in the full rank lattice $\Lambda\cap \R\la \wtilde G_0\ra$ of $\R\la \wtilde G_0\ra$) implies that $g$ is contained in the support of some atom $U\in \mathcal A(G_0)$, whence $g\in \wtilde G_0$ by definition of $\wtilde G_0$. This shows that  $$\wtilde G_0=G_0\cap \R\la \wtilde G_0\ra.$$
Let $\pi:\R^d\rightarrow \R\la \wtilde G_0\ra^\bot$ be the orthogonal projection, so $\wtilde G_0=\{g\in G_0:\; \pi(g)=0\}$. Clearly, $\C(\wtilde G_0)=\R\la \wtilde G_0\ra$ is contained in the lineality space $\mathcal E$ for $\C(G_0)$. If $0\in \C^*(\pi(G_0)\setminus \{0\})=\C^*(\pi(G_0\setminus \wtilde G_0))$, then there will be a nontrivial positive linear combination of elements from $G_0\setminus \wtilde G_0$ contained in $\R\la \wtilde G_0\ra=\C(\wtilde G_0)$. But then Proposition \ref{prop-rational-atoms} ensures that there is some zero-sum $U\in \mathcal B(G_0)$ whose support contains elements from $G_0\setminus \wtilde G_0$, and thus an atom as well, contradicting the definition of $\wtilde G_0$.
Therefore $0\notin \C^*(\pi(G_0)\setminus \{0\})$, ensuring that $\C(\pi(G_0))$ has trivial lineality space, meaning $\R\la \wtilde G_0\ra=\ker \pi=\mathcal E$ is the lineality space for $\C(G_0)$. Hence $\C(\wtilde G_0)=\R\la \wtilde G_0\ra=\mathcal E$ and $\wtilde G_0= G_0\cap \R\la \wtilde G_0\ra=G_0\cap \mathcal E$, as desired.
\end{proof}

\subsection{Elementary Atoms and Positive Bases}

The set of atoms $\mathcal A(G_0)$ is central to the study of  factorizations over $\mathcal B(G_0)$, being the basic building block of all zero-sum sequences.  However,  some atoms are more elementary than others.

\begin{definition}
Let $G_0\subseteq G$ be a subset of a torsion-free abelian group. An atom $U\in\mathcal A(G_0)$ is called \textbf{elementary} if $\mathcal A(X)=\emptyset$ for every proper subset $X\subset \supp(U)$.  We let $\mathcal
A^{\mathsf{elm}}(G_0)\subseteq\mathcal A(G_0)$ denote the set of all elementary atoms.
\end{definition}

The notion of an elementary atom was introduced to the context of factorization theory in \cite{G-precursor} (with a  variation in the definition), but has independent and much older origins in both Convex Geometry and Matroid Theory \cite{Rockafellar}. The approach taken in \cite{G-precursor} followed the Matroid theoretic branch. However, the vein from  Convex Geometry provides a natural framework for adapting the arguments, and one which we will generalize quite extensively in this work.
We begin with the following theorem explaining why elementary atoms can, in some sense, be considered the basic building blocks for all atoms. Theorem \ref{thm-carahtheodory-elm-atom} is essentially Carath\'eordory's Theorem translated into the language of zero-sum sequences. At the very least, the key idea used in the proof of Carath\'eordory's Theorem is the same  used to prove Theorem \ref{thm-carahtheodory-elm-atom}.  The details are given  in   \cite[Theorem 3.7]{G-precursor} albeit using an alternative definition for an elementary atom. To avoid any confusion, we include the short proof below.

\begin{theorem}\label{thm-carahtheodory-elm-atom} Let $G_0\subseteq G$ be a subset of a torsion-free abelian group $G$. Any zero-sum $S\in \mathcal B_{\mathsf{rat}}(G_0)$ can be written as a product of rational powers of elementary atoms, i.e., there are elementary atoms $U_1,\ldots,U_\ell\in \mathcal A^{\mathsf{elm}}(G_0)$ and positive rational numbers  $\alpha_1,\ldots,\alpha_\ell\in \Q_+$ such that
$$S={\prod}_{i\in [1,\ell]}^\bullet U_i^{[\alpha_i]}\quad\mbox{ with } \quad \ell\leq |\supp(S)|.$$
\end{theorem}

\begin{proof}
 We show the theorem holds for any rational zero-sum $S\in \mathcal B_{\mathsf{rat}}(G_0)$ by induction on $|\supp(S)|$. When $|\supp(S)|=1$, then $S$ is itself a rational power of an elementary atom, and the theorem holds trivially. This completes the base of the induction. Given $S\in \mathcal B_{\mathsf{rat}}(G_0)$, there exists an elementary atom
 $U\in \mathcal A^{\mathsf{elm}}(G_0)$ with $\supp(U)\subseteq \supp(S)$. Indeed, simply consider an atom $U\in \mathcal A(G_0)$ with $\supp(U)\subseteq \supp(S)=\supp(S^{[N]})$ and $\supp(U)$ minimal subject to this constraint, where $S^{[N]}\in \mathcal B(G_0)$. Let $\alpha=\min\{\vp_{x}(S)/\vp_x(U):\; x\in\supp(U)\}$, which is a positive rational number as $\supp(U)\subseteq \supp(S)$ with $U,\,S\in \mathcal B_{\mathsf{rat}}(G_0)$.  Then $\vp_x(U^{[\alpha]})\leq \vp_x(S)$ for all $x\in \supp(S)$ with equality holding for any $x$ attaining the minimum in the definition of $\alpha$. Thus $U^{[\alpha]}\mid S$ with
 $S\bdot U^{[-\alpha]}\in \mathcal B_{\mathsf{rat}}(G_0)$ a rational zero-sum having $|\supp(S\bdot U^{[-\alpha]})|<|\supp(S)|$. Applying the induction hypothesis to $S\bdot U^{[-\alpha]}$ now completes the proof.
\end{proof}

Given a convex cone $C\subseteq \R^d$, a   subset $X\subseteq C$ such that
$\mathsf C(X)=C$ but $\mathsf C(Y)\neq C$ for all
proper subsets $Y\subset X$ is called a frame  of $C$ \cite{Davis}. A frame for a \emph{positive dimensional} subspace
$\mathcal E\subseteq \R^d$ is called a \textbf{positive basis}. This is the natural extension of  a linear basis to Convex Geometry.
Unlike ordinary linear bases, positive bases can exhibit complex algebraic structure.  Clearly, any positive basis  $X$ of an $n$-dimensional subspace $\mathcal E$ must have
$|X|\geq n+1$.   If equality holds, $X$ is called a \textbf{minimal positive basis} for the subspace $\mathcal E$, though the name is somewhat misleading.
We do not allow  positive bases of zero dimensional subspaces for technical reasons.
 While  a general positive basis can be complex, minimal positive bases are easily described and closely related to elementary atoms, as the following proposition shows. 
In particular, Proposition \ref{prop-char-minimal-pos-basis} shows $U$ is an elementary atom precisely when $\supp(U)$ is a minimal positive basis or $\{0\}$.

\begin{proposition}\label{prop-char-minimal-pos-basis}
For $X\subseteq \R^d$ with  $\mathcal E=\R\la X \ra$ a \emph{nontrivial} space, the following are equivalent.
\begin{itemize}
\item[1.] $X$ is a minimal positive basis for $\mathcal E$.
\item[2.] $0\in \C^*(X)$ and any proper subset $Y\subset X$ is linearly independent.
        \item[3.]  $X\setminus \{x\}$ is linearly independent with $-x\in \C^\circ(X\setminus \{x\})$ for every $x\in X$.
\item[4.] $X\setminus \{x\}$ is linearly independent with $-x\in \C^\circ(X\setminus \{x\})$ for some $x\in X$.
\item[5.] If $M$ is the $d\times |X|$ matrix whose columns are the elements from $X$, then the null space of $M$ has dimension $1$ and is generated by a vector having all  coordinates strictly positive.
\item[6.] $0\in \mathsf C^*(X)$ but $0\notin \mathsf C^*(Y)$ for all proper subsets $Y\subset X$.
\end{itemize}
If $X\subseteq \Lambda$ with $\Lambda\leq \R^d$ a full rank lattice, then the above are equivalent to each of   the following.
\begin{itemize}
\item[7.] $X=\supp(U)$ for some elementary atom $U$.
\item[8.] $|\mathcal A(X)|=1$ and $\mathcal A(Y)=\emptyset $ for all proper $Y\subset X$.
\end{itemize}
\end{proposition}

\begin{proof}


1. $\Rightarrow$ 2. Since $X$ is a positive basis,  $\C(X)=\mathcal E$, thus containing the subspace $\mathcal E$ of dimension $\dim \mathcal E\geq 1$, whence $0\in \C^*(X)$, showing that $X$ is linearly dependent. Since $\C(X)=\mathcal E$, it follows that $\R\la X\setminus \{x\}\ra=\mathcal E$ for all $x\in X$, for otherwise $\C(X)$ would be contained in the proper half space  $\R\la X\setminus\{x\}\ra +\R_+x\subset \mathcal E$, contradicting that $\C(X)=\mathcal E$.
   If a proper subset $Y\subset X$ were linearly dependent, then $X\setminus Y$ is nonempty and $\dim \mathcal E\leq (|Y|-1)+|X\setminus Y|= |X|-1=\dim \mathcal E$, with the last equality following as $X$ is a \emph{minimal} positive basis. This  forces $X\setminus Y$ to be a basis for $\mathcal E$ modulo the subspace $\mathcal E'=\R\la Y\ra$. In particular, $\R\la X\setminus \{x\}\ra\neq \mathcal E$ for any $x\in X\setminus Y$, contrary to what we established above. This shows each proper subset is linearly independent.

2. $\Rightarrow$ 3. For $x\in X$, we have $X\setminus \{x\}$ linearly independent by Item 2.
Since $0\in \C^*(X)$, there must be a strictly positive linear combination of \emph{all} elements from $X$ equal to zero, since any proper subset is linearly independent. Thus $-x$ can be written as a strictly positive linear combination of the linearly independent elements from $X\setminus \{x\}$, showing that $-x\in \C^\circ(X\setminus \{x\})$.

3. $\Rightarrow$ 4. This is immediate.

4. $\Rightarrow$ 5. The dimension of the null space of the $d\times |X|$ matrix $M$ is equal to $|X|-\dim \R\la X\ra$, which equals $1$ by Item 4, since $X\setminus\{x\}$ is linearly independent but $-x\in \C^\circ(X\setminus \{x\})\subseteq \R\la X\setminus \{x\}\ra$. Since $X\setminus \{x\}$ is linearly independent, $\C^\circ(X\setminus \{x\})$ consists of all elements which are a strictly positive linear combination of \emph{all} element from $X\setminus \{x\}$, in which case $-x\in \C^\circ (X\setminus\{x\})$ ensures that the null space of $M$ must contain a vector whose coordinates are all strictly positive, which must then be a generator.

5. $\Rightarrow$ 6.  Let $x_1,\ldots,x_r\in X$ be the distinct elements of $X$. Let $z=(\alpha_1,\ldots,\alpha_r)$ be the generator for the null space of $M$. Then $\alpha_i>0$ for all $i$ by item 5 and $\alpha_1x_1+\ldots+\alpha_rx_r=0$ showing that $0\in \C^*(X)$. On the other hand, if $0\in \mathsf C^*(Y)$ for some proper subset $Y\subset X$, then there would be $\beta_1,\ldots,\beta_r\in \R_+$ with $\beta_1x_1+\ldots+\beta_rx_r=0$, not all $\beta_i$  zero, and some $\beta_j=0$. But then $(\beta_1,\ldots,\beta_r)$ would be a nonzero element of  the null space of $M$  and clearly not a multiple of $z$ (as any non-zero multiple of $z$ has all non-zero coordinates), contradicting Item 5.

6. $\Rightarrow$ 1. Since $0\in \C^*(\{0\})$ and $\R\la X\ra$ is nontrivial, the hypotheses of Item 6 force $0\notin X$. Let $x_1,\ldots,x_r\in X$ be the distinct elements of $X$. Observe that $\C(Y)=\R\la Y\ra$, for a nonempty subset $Y\subseteq X$,  implies
$\C(Y)$ contains the positive dimensional subspace $\R\la Y\ra$, yielding $0\in \C^*(Y)$. By assumption of Item 6, the only nonempty subset $Y\subseteq X$ for which this can hold is $Y=X$. Thus to show $X$ is a positive basis, we only need to show $\C(X)=\R\la X\ra$. Since $0\in \C^*(X)$ but $0\notin \C^*(Y)$ for all proper subsets $Y\subset X$, it follows that there is a strictly positive linear combination  $\alpha_1x_1+\ldots+\alpha_rx_r=0$ using all elements from $X$. Let $y\in \C(X)$ be nonzero. Then $y=\beta_1x_1+\ldots+\beta_rx_r$ for some $\beta_1,\ldots,\beta_r\in \R_+$. Let $\beta=\max_i \beta_i>0$ and let $\alpha=\min_i \alpha_i>0$. Then $-y=(-\beta_1x_1-\ldots-\beta_rx_r)+\frac{\beta}{\alpha}
(\alpha_1x_1+\ldots+\alpha_rx_r)=
\Sum{i=1}{r}(-\beta_i+\frac{\beta}{\alpha}\alpha_i)x_i,$ with $-\beta_i+\frac{\beta}{\alpha}\alpha_i\geq 0$ for all $i$, showing $-y\in \C(X)$. As $y\in \C(X)$ was an arbitrary non-zero element, this means the convex cone $\C(X)$ is in fact a subspace, which is only possible if $\C(X)=\R\la X\ra=\mathcal E$. This shows that $X$ is a positive basis for $\mathcal E$.

Assuming $X$ is not a \emph{minimal} positive basis, we obtain $r=|X|\geq \dim \mathcal E+2$. In such case, there must be a linear combination of a \emph{proper} subset of $X$ equal to $0$. Thus, by re-indexing the $x_i$, we can w.l.o.g. assume $\gamma_1x_1+\ldots+\gamma_sx_s=0$ for some $s\in [1,\dim \mathcal E+1]\subseteq [1,r-1]$ and some $\gamma_i\in \R$ with not all $\gamma_i=0$. If $\gamma_i\geq 0$ for all $i$, then this contradicts the hypothesis $0\notin \C^*(Y)$ for $Y=\{x_1,\ldots,x_s\}\subset X$. Therefore we can define $\gamma=\min \{\alpha_i/|\gamma_i|:\;\gamma_i<0,\, i\in [1,s]\}>0$, which ensures that $\alpha_i+\gamma \gamma_i\geq 0$ for all $i\in [1,s]$, with equality holding for any $i$ attaining the minimum in the definition of $\gamma$. But now $(\alpha_1+\gamma\gamma_1)x_1+\ldots+(\alpha_s+\gamma\gamma_s)x_s+\alpha_{s+1}x_{s+1}
+\ldots+\alpha_{r}x_r=0$ is  a positive linear combination equal to zero. Moreover, the coefficient of $x_r$ is $\alpha_r>0$ as $r>s$, and the coefficient of some $x_i$ with $i\in [1,s]$ is zero (namely, any $i\in [1,s]$ attaining the minimum in the definition of $\gamma$), contradicting that $0\notin \C^*(Y)$ for all proper subsets $Y\subset X$.
 This completes the equivalence of the first 6 items.

Now assume $X\subseteq \Lambda$ with $\Lambda\leq \R^d$ a full rank lattice.  Then, in view of Proposition \ref{prop-rational-atoms}, we see that Items 6 and 7 are equivalent reformulations of each other. The implication 8. $\Rightarrow$ 7. is also immediate. To complete the proof, we show that 5. $\Rightarrow$ 8.   Since $X\subseteq \Lambda$, each vector $x\in X$ can be expressed as an integer linear combination of the elements from some lattice basis $E$ for $\Lambda$, which is then also a linear basis for $\R^d$. If we use $E$ to express the column vectors in $M$, it follows that $M$ is an integer matrix, whence
the generator $z$ of the null space of $M$ from Item 5 will have each of its $i$-th coordinates  $\alpha_i>0$ being a rational number. By clearing denominators, we can further assume $\alpha_i\in \Z$ with $\gcd(\alpha_1,\ldots,\alpha_r)=1$ and $\alpha_i\geq 1$ for all $i$. It is then clear that any positive multiple of $z=(\alpha_1,\ldots,\alpha_r)$ having all of its coordinates integers must be an integer multiple of $z$. This means  every zero-sum sequence $S\in \mathcal B(X)$ has the form $S=U^{[n]}$ for some integer $n\geq 1$, where $U=\prod_{i\in [1,r]}^\bullet x_i^{[\alpha_i]}\in \mathcal B(X)$, which implies that $U\in \mathcal A(X)$ is the unique atom with support contained in $X$. Hence $|\mathcal A(X)|=1$, and since $\supp(U)=X$, it follows that $\mathcal A(Y)=\emptyset$ for all proper $Y\subset X$, completing the proof.
\end{proof}

\subsection{Reay Systems}

Much of the material from this subsection can be extracted by a careful examination, variation and  reformulation of key ideas from the proofs found in the early works  \cite{Bonnice-Klee}  \cite{Davis} \cite{Reay-memoir} \cite{Reay-gen} \cite{McKinney} \cite{shephard}. The material in the format and generality we  require is not readily available,  so this section will serve as the foundation  for more extensive generalizations in later sections. The key definition is the following, which may be found buried in a proof of Reay in the specialized case when $X=X_1\cup \ldots\cup X_s$ is a positive basis \cite{Reay-memoir}. Since the authors of the aforementioned references were rather focussed on the study of positive bases,   their notions were not developed beyond this context in the fuller generality needed for this paper.

\begin{definition}
Let $X_1,\ldots,X_s\subseteq \R^d$ be nonempty subsets, where $s\geq 0$ and $d\geq 0$. For $j\in [0,s]$,     let  $\pi_j:\R^d\rightarrow \R\la X_1\cup \ldots\cup X_j\ra^\bot$ be the orthogonal projection.  We say that $\mathcal R=(X_1,\ldots,X_s)$ is a \textbf{Reay (coordinate) system} for the subspace $\R\la X_1\cup \ldots\cup X_s\ra\subseteq \R^d$ if
  $$\pi_{j-1}(X_j)\mbox{ is a minimal positive basis of size $|\pi_{j-1}(X_j)|=|X_j|$, \ for every $j\in [1,s]$}.$$
 \end{definition}

 We view the empty tuple as a Reay system for the trivial space.  Let $\mathcal R=(X_1,\ldots,X_s)$ be a Reay  system for the subspace $\mathcal E\subseteq \R^d$ and let $X=X_1\cup\ldots\cup X_s$.   We say that a subset $Y\subseteq \R^d$ \textbf{contains} the Reay system $\mathcal R$ if $X\subseteq Y$.  If the Reay system $\mathcal R$ has $X=X_1\cup\ldots\cup X_s$ being a positive basis, then we call $\mathcal R$ a \textbf{Reay basis}. We call $s\geq 0$ the \textbf{depth} of the Reay system, in which case a minimal positive basis  is just a Reay system of depth $1$, and  an element $x\in X$ with  $x\in X_j$ is said to be at depth $j$.
  The basic existence result for Reay systems is the following.

 \begin{proposition}\label{prop-reay-basis-exists}
 Let $X\subseteq\R^d$, where $d\geq 0$,  and  let $\mathcal E= \mathsf C(X)\cap-\mathsf C(X)$.
Then $X$ contains a Reay system for $\mathcal E$. Moreover,
  if $Y\subseteq X$ is a minimal positive basis, then $X$ contains a Reay system $(X_1,\ldots,X_s)$ for  $\mathcal E$ with $X_1=Y$.
\end{proposition}

\begin{proof} Recall that $\mathcal E$ is the maximal subspace contained in $\C(X)$.
If $\dim \mathcal E=0$, then $X$ contains no positive basis, and the empty tuple gives the desired Reay system. Therefore  assume $\dim \mathcal E\geq 1$ and  proceed by induction on $\dim \mathcal E$, the base case having just been completed.  Since $\C(X)$ contains the positive dimension subspace $\mathcal E$, we have
$0\in \C^*(X\setminus \{0\})$, and thus there must be a minimal (by inclusion)  subset $X_1\subseteq X\setminus \{0\}$ with  $0\in \C^*(X_1)$. In view of Proposition \ref{prop-char-minimal-pos-basis}.6, such a subset $X_1\subseteq X$ is a  minimal positive basis contained in $X$.  Let $\pi_1:\R^d\rightarrow \R\la X_1\ra^\bot$ be the orthogonal projection.  Since there is nothing otherwise special about $X_1$, we can w.l.o.g. assume $X_1$ is equal to \emph{any} minimal positive basis $Y\subseteq X$.
In view of the maximality of $\mathcal E$, we have $\C(X_1)=\R\la X_1\ra\subseteq \mathcal E$. If $\R\la X_1\ra=\mathcal E$, we are done. So instead assume $\dim \R\la X_1\ra<\dim \mathcal E$. Then $\pi_1(\mathcal E)$ will be the maximal subspace contained in $\pi_1(\C(X))=\C(\pi_1(X))$, so by  induction hypothesis $\pi_1(\mathcal E)$ has a Reay system $(\pi_1(X_2),\ldots,\pi_1(X_s))$ with $X_i\subseteq X$ for all $i$, and by discarding elements with equal images under $\pi_1$, we can w.l.o.g. assume $|X_j|=|\pi_1(X_j)|$ for all $j$. It now follows that $(X_1,\ldots,X_s)$ will be a Reay system for $\mathcal E$.
\end{proof}

We continue with some basic observations regarding Reay systems.

\begin{proposition}\label{prop-reay-basis-properties}
Let $\mathcal R=(X_1,\ldots,X_s)$ be a Reay system for a subspace $\mathcal E\subseteq \R^d$ of dimension $n$ and let $X=X_1\cup \ldots\cup X_s$. For $j\in [0,s]$, let $\mathcal E_j=\R\la X_1\cup\ldots\cup X_j\ra$ and  let  $\pi_j:\R^d\rightarrow \mathcal E_j^\bot$ be the orthogonal projection.
\begin{itemize}
\item[1.] $\C(X)=\R\la X\ra$, and $\bigcup_{i=1}^s X_i\setminus \{x_i\}$ is a linear basis for $\mathcal E$ for any $x_i\in X_i$.
\item[2.] $X=X_1\cup\ldots\cup X_s$ is a disjoint union with $|X|=\dim \mathcal E+s=n+s\leq 2n$.
\item[3.] $(X_1,\ldots,X_j)$ is a Reay system for $\mathcal E_j$ and $(\pi_{j-1}(X_j),\ldots,\pi_{j-1}(X_s))$ is a Reay system for $\pi_{j-1}(\mathcal E)$, for any $j\in [1,s]$.
\item[4.] If $X$ is a positive basis, then any Reay system $(Y_1,\ldots,Y_s)$ for $\mathcal E$ contained in $X$ is a Reay basis and has $X=Y_1\cup \ldots\cup Y_s$. Moreover, the Reay systems given in Item 3 are all Reay bases.
\end{itemize}
\end{proposition}

\begin{proof}
1. A quick inductive argument on $j=0,1,\ldots,s$ shows that $\C(X_1\cup \ldots\cup X_j)=\mathcal E_j$.
  Using Proposition \ref{prop-char-minimal-pos-basis}.3 and an inductive argument on $j=1,2,\ldots,s$ shows that $\bigcup_{i=1}^{j}X_i\setminus \{x_i\}$ is a linear basis for $\R\la X_1\cup\ldots\cup X_j\ra$ for any $x_i\in X_i$. The case $j=s$ yields Item 1.

2. That the elements in $X_1\cup \ldots\cup X_j$ are distinct follows by a simple inductive argument on $j=1,\ldots,s$ utilizing that $\pi_{j-1}$ is injective on $X_j$ with all elements in a minimal positive basis $\pi_{j-1}(X_j)$ nonzero. Hence $n=\dim \mathcal E=\Sum{i=1}{s}(|X_i|-1)=|X|-s$ by Item 1.  Moreover, since $1\leq \dim \mathcal E_1<\dim \mathcal E_2<\ldots<\dim \mathcal E_s= \dim \mathcal E=n$ (as each minimal positive basis  $\pi_{j-1}(X_j)$ must span a \emph{nontrivial} subspace), we have $s\leq n$, completing Item 2.

3. This follows immediately from the recursive definition of a Reay system.

4. That any Reay system $(Y_1,\ldots,Y_s)$ for $\mathcal E$ contained in the positive basis $X$  must have $X=Y_1\cup \ldots\cup Y_s$, and therefore be a Reay Basis, follows from Item 1, for otherwise $\C(Y)=\mathcal E$ for the proper subset $Y=Y_1\cup \ldots\cup Y_s$, contradicting that $X$ is a positive basis.
If $Y\subseteq X_j\cup \ldots\cup X_s$ is a subset with $\C(\pi_{j-1}(Y))=\pi_{j-1}(\mathcal E)$, then $\C(X_1\cup\ldots \cup X_{j-1}\cup Y)=\mathcal E$ follows in view of $\C(X_1\cup \ldots\cup X_{j-1})=\mathcal E_{j-1}=\ker \pi_{j-1}$ (which holds by Item 1 applied to the Reay system $(X_1,\ldots,X_{j-1})$). Likewise, if $Y\subseteq X_1\cup \ldots\cup X_j$ with $\C(Y)=\mathcal E_j$, then $\C(Y\cup X_{j+1}\cup \ldots\cup X_{s})=\mathcal E$ follows in view of $\C(\pi_j(X_{j+1})\cup \ldots\pi_j(X_s))=\pi_j(\mathcal E)$ (which holds in view of Item 1 applied to the Reay system $(\pi_{j}(X_{j+1}),\ldots,\pi_{j}(X_s))$).
Thus $\mathcal R$ being a Reay basis implies that  $(X_1,\ldots,X_j)$ and $(\pi_{j-1}(X_j),\ldots,\pi_{j-1}(X_s))$ are Reay bases too, for any $j\in [1,s]$.
\end{proof}

The key property of Reay systems is that they allow for a certain type of unique expression, a fact not highlighted in the original work of Reay.

\begin{proposition}\label{prop-reay-basis-unique-expression}
Let $(X_1,\ldots,X_s)$ be a Reay system for the subspace $\mathcal E\subseteq\R^d$.
Then every $z\in \mathcal E$ has a  unique expression as $$z=\Sum{j=1}{s}\Summ{x\in X_j}\alpha_{x}x$$ with all $\alpha_x\in\R_+$  but, for every
    $j\in [1,s]$, not all  $\alpha_x$ with $x\in X_j$ are nonzero.
\end{proposition}

\begin{proof}
If $\mathcal E$ is trivial, then $x=s=0$ and the unique expression for $x$ is the empty sum. Therefore assume $\dim \mathcal E\geq 1$. We proceed by induction on $s$. The first nontrivial case is $s=1$, when we have only the minimal positive basis $X_1$. In view of the characterization given in Proposition \ref{prop-char-minimal-pos-basis}.3, one can apply a linear transformation $\varphi$ mapping the first $|X_1|-1$ elements of $X_1$ to the standard basis vectors in $\R^d$, and then the remaining element of $X_1$ will map to an element with all coordinates strictly negative, a case for which the uniqueness of expression for $\varphi(z)$ is easily verified, and one which implies the same property for the original element $z$. This completes the base case. Now assume $s\geq 2$ and that we have unique expression for all smaller values of $s$. For $j\in [0,s]$, let $\mathcal E_j=\R\la X_1\cup \ldots\cup X_j\ra$ and let $\pi_j:\R^d\rightarrow \mathcal E_j^\bot$ be the orthogonal projection, so $\mathcal E_s=\mathcal E$. By Proposition \ref{prop-reay-basis-properties}.3, $(X_1,\ldots,X_{s-1})$ is a Reay system for $\mathcal E_{s-1}$ and $\pi_{s-1}(X_s)$ is a minimal positive basis for $\pi_{s-1}(\mathcal E)$.
Since $\mathcal E_{s-1}\subseteq \mathcal E_s=\mathcal E$, every $z\in \mathcal E$ has a unique expression as $z=a+b$ with $a\in \ker \pi_{s-1}=\mathcal E_{s-1}$ and $b\in \mathcal E_{s-1}^\bot\cap \mathcal E=\pi_{s-1}(\mathcal E)$.
Applying the base case to the minimal positive basis $\pi_{s-1}(X_s)$ for  $\pi_{s-1}(\mathcal E)$, we find there is a unique expression $b=\Summ{x\in X_s}\alpha_x\pi_{s-1}(x)=\pi_{s-1}\left(\Summ{x\in X_s}
\alpha_xx\right)$ with $\alpha_x\in \R_+$ for all $x\in X_s$ and not all $\alpha_x$  nonzero. Thus there is a unique expression $b=u_z+\Summ{x\in X_s}\alpha_xx$ with $u_z\in \ker \pi_{s-1}=\mathcal E_{s-1}$, with $\alpha_x\in \R_+$, and with not all $\alpha_x$  nonzero (since $b\in \pi_{s-1}(\mathcal E)$ with $\pi_{s-1}$ a projection, we have $\pi_{s-1}(b)=b$).
Since $b$ is uniquely determined by $z$, the element $u_z\in \mathcal E_{s-1}$ is uniquely determined by $z$.  But now there is a unique expression $z=a'+\Summ{x\in X_s}\alpha_xx$ with $a'\in \mathcal E_{s-1}$, with $\alpha_x\in \R^+$, and with not all $\alpha_x$  nonzero (namely, $a'=a+u_z$). Applying the induction hypothesis to the element $a'$ and  Reay system $(X_1,\ldots,X_{s-1})$  for $\mathcal E_{s-1}$ now yields the desired  unique expression for $z$.
\end{proof}

There is another way to view Proposition \ref{prop-reay-basis-unique-expression}. If $\mathcal R=(X_1,\ldots,X_s)$ is a Reay system for the subspace $\mathcal E\subseteq \R^d$ with $X=X_1\cup\ldots\cup X_s$, then define $$\mathfrak B=\{\C(Y):\; Y\subseteq X \und X_i\nsubseteq Y\mbox{ for every $i\in [1,s]$}\}.$$ For every choice of elements $x_i\in X_i$, for $i\in [1,s]$, the set $Y=\bigcup_{i=1}^s X_i\setminus \{x_i\}$ is a linear basis for $\mathcal E$ by Proposition \ref{prop-reay-basis-properties}.1. Thus  Proposition \ref{prop-reay-basis-unique-expression} implies that $\mathcal E=\bigcup_{C\in \mathfrak B} C^\circ$ is the disjoint union of the relative interiors of the cones in $\mathfrak B$.  Combining these observations, we see that if we intersect each cone from $\mathfrak B$ with the unit sphere in $\mathcal E$, we obtain an object homeomorphic to a  \emph{simplicial complex} of dimension $\dim \mathcal E-1$ whose union is the unit sphere in $\mathcal E$. In the parlance of topologists, such an object is called a \emph{simplicial sphere} or \emph{triangulated sphere}. Indeed, owing to the method of construction, we obtain a more restricted class of simplicial sphere known as a \emph{starshaped sphere} (somewhat surprisingly, not every simplicial sphere can be constructed this way). We direct the reader to \cite{Buchstaber-Panov} for a more comprehensive account, including more details of what follows below.

A \textbf{fan} is a \emph{finite} collection $\mathfrak B$ of polyhedral cones $\C(Y)\subseteq \R^d$  each having trivial lineality space $\C(Y)\cap -C(Y)=\{0\}$.
The \textbf{faces} of $\C(Y)$ are the sub-cones $\C(Z)$ with $Z\subseteq Y$ obtained by intersecting  $\C(Y)$ with a hyperplane defining a closed  half-space that contains $\C(Y)$, and it is required that each face of a cone from the fan be an element of the fan, and that  the intersection of any two
cones in  the fan be a face of each.
 It is a \textbf{simplicial} fan if every $\C(Y)\in \mathfrak B$ is generated by a linearly independent set $Y$, in which case there is a face for each subset $Z\subseteq Y$, and it is a \textbf{complete} fan if the union of all cones in $\mathfrak B$ equals an entire subspace. Thus, the set $\mathfrak B$ defined above from a Reay system is a special type of \textbf{complete simplicial fan} for the subspace $\mathcal E\subseteq \R^d$. Complete simplicial fans are in bijective correspondence with starshaped spheres, and  rational fans (those whose vertices come from  a lattice) are central to the definition of Toric Varieties, though we will only need their more basic properties.

One easily notes that a Reay system (with depth $s\geq 2$) is not always stable even under small perturbations of its defining vectors. This forces us to work with the more general concept of a complete simplicial fan, which maintains many of the essential features of a Reay system while gaining the important property of being stable under small perturbations. Worth noting, any complete fan must be pure (that is, all maximal cones must have the same dimension). Indeed, if $\mathfrak B$ is a complete fan for $\R^d$, then the sub-collection of all $d$-dimensional cones from $\mathfrak B$, together with all their faces (sub-cones), would also be a fan. Their union must be all of $\R^d$, for if it were not, then its complement would be a $d$-dimensional subset of $\R^d$, which clearly cannot be written as a union of a finite number of lower dimensional objects. But once we know their union is all of $\R^d$, it then follows that any lower dimensional cone must lie in one of the $d$-dimensional cones, and thus only $d$-dimensional cones are maximal.

Let $\mathfrak B$ be a simplicial fan in $\R^d$. For an integer $k\in [0,d]$, we use $\mathfrak B_k$ to denote the subset of $\mathfrak B$ consisting of all $k$-dimensional cones in $\mathfrak B$ (generated by $k$ elements). Thus the elements of $\mathfrak B_k$ when intersected with the unit sphere give rise (up to homeomorphism) to $(k-1)$-dimensional simplices.
A vertex set $V$ for  $\mathfrak B$ is a collection of nonzero elements, one chosen from each $C\in \mathfrak B_1$. For instance, the set $V=\bigcup_{C\in \mathfrak B_1}\big(C\cap \partial(B_1(0))\big)$ consisting of  the elements of the unit sphere contained in some $B\in \mathfrak B_1$ (i.e., the $0$-dimensional vertices of the associated starshaped sphere) is one possible set of vertices for $\mathfrak B$. Note that $|V(\mathfrak B)|=|\mathfrak B_1|$ is finite. Whenever dealing with a simplicial fan $\mathfrak B$, we will fix a vertex set $V(\mathfrak B)$, and if none is explicitly mentioned,  the unit sphere representatives are assumed to be  the vertices.
Every $x\in \bigcup_{C\in \mathfrak B}C$ has a unique cone  $C\in \mathfrak B$ for which $x\in C^\circ$. The cone $C$ is generated by a  set of linearly independent vectors $B_C\subseteq V(\mathfrak B)$, so $C=\C(B_C)$ and $x$ is a strictly positive linear combination of the elements of $B_C$ (note: if $x=0$, then $B_C=\emptyset$ and $C=\{0\}$).
If $\mathfrak B$ is defined using a Reay system $(X_1,\ldots,X_s)$, then we take $V(\mathfrak B)=X_1\cup\ldots\cup X_s$ for the vertices. The coefficients in the linear combination correspond to the baricentric coordinates for points inside the simplicial cone $C$. We define $\supp_{\mathfrak B}(x)=B_C\subseteq V(\mathfrak B)$ to be the  \textbf{support} set of the element $x$ with respect to the simplicial fan $\mathfrak B$, which is then the unique subset $\supp_\mathfrak B(x)\subseteq V$ such that $x\in \C^\circ(\supp_\mathfrak B(x))$. If $\mathfrak B$ arises from a Reay system $\mathcal R$, then we let $$\supp_{\mathcal R}(x)=\supp_{\mathfrak B}(x)$$ and call this set the Reay support of the element $x$.  By Proposition \ref{prop-reay-basis-properties}.1, $\supp_\mathcal R(x)$ is always a linearly independent set.  We carry on with some basic properties about Reay systems and positive bases.

\begin{proposition}\label{prop-reay-RayAlg}
Let $\mathcal R=(X_1,\ldots,X_s)$ be a Reay system for a subspace $\mathcal E\subseteq \R^d$ and  let $\mathcal E_j=\R\la X_1\cup\ldots\cup X_j\ra$ for $j\in [0,s]$.
For every $k\in [1,s]$, $$\C(X_k)\cap \mathcal E_{k-1}=\R_+u_k$$ for some (possibly zero) $u_k\in \mathcal E_{k-1}$. Moreover,  $$X'_k=\supp_{\mathcal R}(-u_k)\cup X_k\subseteq X_1\cup\ldots\cup X_k$$ is a minimal positive basis for some subspace $\mathcal E'_k\subseteq \mathcal E_k$.
\end{proposition}

\begin{proof}
Let
 $k\in [1,s]$ be arbitrary.
Since $X_k\setminus \{x_k\}$ is linearly independent for any $x_k\in X_k$ by Proposition \ref{prop-reay-basis-properties}.1, it follows that $\dim \R\la X_k\ra=|X_k|$ or $|X_k|-1$, depending on whether the elements of $X_k$ are linearly independent or linearly dependent.
If they are linearly dependent, so $\dim \R\la X_k\ra=|X_k|-1$, then, since $\mathcal E_k=\mathcal E_{k-1}+\R\la X_k\ra$, we have $\dim (\mathcal E_{k-1}\cap \R\la X_k\ra)=\dim \mathcal E_{k-1}+\dim \R\la X_k\ra-\dim \mathcal E_k=\Big(\Sum{i=1}{k-1}(|X_i|-1)\Big)+|X_k|-1-\Big(\Sum{i=1}{k}(|X_i|-1)\Big)=0$, meaning $\R\la X_k\ra\cap \mathcal E_{k-1}=\{0\}$. In this case, $\C(X_k)\cap \mathcal E_{k-1}=\R_+u_k$ for $u_k=0$. On the other hand, if they are linearly independent, so $\dim \R\la X_k\ra=|X_k|$, then  we instead have  $\dim (\mathcal E_{k-1}\cap \R\la X_k\ra)=1$, meaning $\R\la X_k\ra\cap \mathcal E_{k-1}$ is a one-dimensional subspace. Since $0\in\C^*(\pi_{k-1}(X_k))$ (in view of Proposition \ref{prop-char-minimal-pos-basis}.6), it follows that there is a nontrivial positive linear combination of   elements from $X_k$ that lies in $\mathcal E_{k-1}$. Since the elements of $X_k$ are linearly independent, this linear combination must be a \emph{nonzero} element of $\mathcal E_{k-1}$, say $u_k\in \mathcal E_{k-1}\cap \C^*(X_k)$. If  $-u_k$ were also contained in $\C(X_k)$, then $0\in \C^*(X_k)$, contradicting that the elements of $X_k$ are linearly independent. Therefore we instead conclude  that $\mathcal E_{k-1}\cap \C(X_k)=\R_+u_k$ in this case as well.

If $u_k=0$, then $X_k$ is linearly dependent and $\supp(-u_k)=\emptyset$, so $X'_k=X_k$. By Proposition \ref{prop-reay-basis-properties}.1, every proper subset of $X_k$ is linearly independent, in which case $X'_k=X_k$  is a minimal positive basis by Proposition \ref{prop-char-minimal-pos-basis}.2, as desired. Therefore we now assume $u_k\neq 0$ with $X_k$ linearly independent.

Since $-u_k\in \mathcal E_{k-1}$, we have $\supp_{\mathcal R}(-u_k)\subseteq X_1\cup\ldots\cup X_{k-1}$.
Since $u_k\in  \C^*(X_k)$ and  $-u_k\in \C(\supp_{\mathcal R}(-u_k))$, we have
 $0\in  \C^*(X'_k)$. Consequently, to show $X'_k$ is a minimal positive basis, it suffices by Proposition \ref{prop-char-minimal-pos-basis}.6 to show $0\notin \mathsf C^*(Y)$ for all proper subsets $Y\subset X'_k$. To this end, consider an arbitrary subset $Y\subseteq X'_{k}$ with $0\in \C^*(Y)$. Since $\supp_{\mathcal R}(-u_k)$ is a linearly independent subset, $Y\subseteq \supp_{\mathcal R}(-u_k)$ is not possible, implying $Y\cap X_k\neq \emptyset$. However, since $0\notin \C^*(\pi_{k-1}(Z))$ for any proper subset $Z\subseteq X_k$ (per Proposition \ref{prop-char-minimal-pos-basis}.6), any strictly positive linear combination of a proper subset of terms
 from $X_k$ lies outside the subspace $\mathcal E_{k-1}$, and thus cannot be combined with any linear combination of terms from $\mathcal E_{k-1}$ to yield $0$.
 In consequence, we conclude that $X_k\subseteq Y$. Thus $0\in \C^*(Y)$ ensures that $0=a+b$ with $a\in \C^*(X_k)$ and $b\in \C(\supp_{\mathcal R}(-u_k)\cap Y)\subseteq \mathcal E_{k-1}$. However, since $\mathcal E_{k-1}\cap \C(X_k)=\R_+u_k$, we must have $a=\alpha u_k$ for some positive $\alpha>0$, and by re-scaling we may w.l.o.g. assume $\alpha=1$.
 Hence $-u_k=-a=b\in \C(\supp_{\mathcal R}(-u_k)\cap Y)$. However, by definition of $\supp_{\mathcal R}(-u_k)$, there is no proper subset of $\supp_{\mathcal R}(-u_k)$ that contains $-u_k$ in its positive span, so we must have $\supp_{\mathcal R}(-u_k)\subseteq Y$, which together with $X_k\subseteq Y$ implies that $Y=X_k$. As $Y\subseteq X'_k$ was an arbitrary subset with $0\in \C^*(Y)$, we conclude that no proper subset  $Y\subset X'_k$ has $0\in \C^*(Y)$, completing the proof.
\end{proof}

Let $\mathcal R=(X_1,\ldots,X_s)$ be a Reay system for $\R^d$ and let $X=X_1\cup \ldots\cup X_s$. Then each $z\in \R^d$ corresponds via Proposition \ref{prop-reay-basis-unique-expression} uniquely to a tuple $\alpha(z)=(\alpha_x(z))_{x\in X}\in \R_+^{|X|}$ with $\Summ{x\in X}\alpha_x(z)x=z$ such that, for every $j\in [1,s]$, not all $\alpha_x(z)$ for $x\in X_j$ are non-zero.
Unlike Euclidean coordinates, if $y,\,z\in \R^d$, we may not have $\alpha(y+z)=\alpha(y)+\alpha(z)$. The problem is that, in one (or more) of  the $s$ groupings of coordinates in $\alpha(y)+\alpha(z)$ corresponding to the $X_j$, all  coordinates may be strictly positive, which is not allowed. Proposition \ref{prop-reay-RayAlg} gives a means to quickly transform $\alpha(y)+\alpha(z)$ into  $\alpha(y+z)$. For each $k\in [1,s]$, there is an expression $\Summ{x\in X_k}(\alpha_x^{(k)})x=u_k=\Sum{i=1}{k-1}\Summ{x\in X_i}\beta_{x}^{(k)}x$ with all $\alpha_x^{(k)}> 0$ and $\beta_{x}^{(k)}\geq 0$. We have $\beta_{x}^{(k)}>0$ precisely when $x\in \supp_{\mathcal R}(u_k)$. Let $\mathbf a_k=(\mathbf a_k(x))_{x\in X}\in \R^{|X|}$ be the vector with  $\mathbf a_k(x)=-\alpha_x^{(k)}$ for $x\in X_k$, with $\mathbf a_k(x)=\beta_x^{(k)}$ for $x\in X_1\cup \ldots\cup X_{k-1}$, and will all other coordinates zero.
Then, if the $k$-th grouping in  $\alpha(y)+\alpha(z)$ has all its coordinates positive, there will be a unique multiple of $\mathbf a_k$, namely $b_k \mathbf a_k$ with $b_k=\min_{x\in X_k} \frac{\alpha_x(y)+\alpha_x(y)}{|\mathbf a_k(x)|}$, such that $\alpha(y)+\alpha(z)+b_k\mathbf a_k$ has all coordinates non-negative but has at least one zero in the $k$-th grouping $X_k$.  Since only coordinates in grouping $k$ and lower change by adding $b_k \mathbf a_k$, we can sequentially apply the relations $b_k\mathbf a_k$ as needed for $k=s,s-1,\ldots,1$ until we reduce $\alpha(y)+\alpha(z)$ to $\alpha(y+z)$ in at most $s\leq d$ steps.

The following proposition is a refined statement of one of the main goals in the original work of Reay \cite{Reay-memoir}.

\begin{proposition}\label{prop-reays-structureresult}
Let $X\subseteq \R^d$ be a positive basis for $\R^d$ with $d\geq 1$. Then there exists a Reay basis $\mathcal R=(X_1,\ldots,X_s)$ for $\R^d$ with $X=X_1\cup\ldots\cup X_s$ and
$$|X_1|\geq |X'_{2}|\geq |X_2|\geq |X'_{3}|\geq |X_3|\geq\ldots\geq |X'_{s}|\geq |X_s|\geq 2,$$ where each $X'_j=(\supp_{\mathcal R}(-u_j)\setminus \mathcal E_{j-2})\cup X_j$ with  $\mathcal E_{j-2}=\R\la X_1\cup \ldots\cup X_{j-2}\ra$ and  $u_j$ as given in Proposition \ref{prop-reay-RayAlg}, for $j\in [2,s]$.
\end{proposition}

\begin{proof}
If $X$ is a minimal positive basis, then $X=X_1$ is itself a Reay basis of depth $s=1$, in which case the proposition is trivial. Therefore we can assume otherwise. In particular, the proposition is true when $|X|\leq 2$, allowing us to proceed by induction on $|X|$.
 By Proposition \ref{prop-reay-basis-exists}, $X$ contains a minimal positive basis $X_1$ for some subspace, so we may w.l.o.g. assume $X_1\subseteq X$ is a \emph{maximal cardinality} minimal positive basis contained in $X$.
Let $\pi_1:\R^d\rightarrow \R\la X_1\ra^\bot$ be the orthogonal projection. We can take any Reay system $(\pi_1(X_2),\ldots,\pi_1(X_s))$, where $X_i\subseteq X$ are subsets with $|\pi_1(X_i)|=|X_i|$, and then $(X_1,X_2,\ldots,X_s)$ will be a Reay system.
By  Proposition \ref{prop-reay-basis-exists}, there \emph{is} some Reay system $(X_1,\ldots,X_s)$ with $X=X_1\cup\ldots\cup X_s$. Now   $\pi_1$ is injective on $X\setminus X_1$ (by definition of a Reay system) while Proposition \ref{prop-reay-basis-properties}.4 ensures that   $\pi_1(X\setminus X_1)$ is a positive basis with $|\pi_1(X\setminus X_1)|=|X\setminus X_1|<|X|$.
 Apply the induction hypothesis to $\pi_1(X\setminus X_1)$ to find a Reay system $\mathcal R'=(\pi_1(X_2),\ldots,\pi_1(X_s))$  satisfying the conclusion of the proposition with $|\pi_1(X_i)|=|X_i|$ for all $i\geq 2$ and $X_2\cup\ldots\cup X_s=X\setminus X_1$.
 Then $\mathcal R=(X_1,\ldots,X_s)$ is a Reay basis with $X=X_1\cup\ldots\cup X_s$ as noted above.
 Let $X'_j$ and $u_j$, for $j\in [2,s]$, be as defined by the proposition for the Reay Basis $\mathcal R$, and let $\pi_1(X_j)'$ and $u'_j$, for $j\in [3,s]$, be the corresponding quantities for the Reay Basis $\mathcal R'$.  Then $\pi_1(u_j)=u'_j$ and  $\pi_1(X'_j)=\pi_1(X_j)'$ for $j\geq 3$.
Thus the induction hypothesis and injectivity of $\pi_1$ on $X\setminus X_1$ yield $$|X_2|\geq |X'_3|\geq |X_3|\geq \ldots|X'_s|\geq |X_s|\geq 2.$$
By definition, $|X'_2|\geq |X_2|$, while Proposition \ref{prop-reay-RayAlg} implies that $X'_2$ is minimal positive basis, so that the maximality of $X_1$ ensures $|X_1|\geq |X'_2|$, completing the proof.
\end{proof}

Reay Bases can be used to help better understand positive bases. However, most of the useful properties of Reay Bases hold for the more general class of Reay systems or even complete simplicial fans. If one needs more refined structure for a positive basis, there is a geometric interpretation  given by Shephard \cite{shephard} involving the Gale diagram of a linear representation of the positive basis.
We conclude the subsection with some important  properties of complete simplicial fans.

\begin{proposition}\label{prop-FanStability} Let $\mathfrak B$ be a complete simplicial fan for   $\R^d$ with $d\geq 1$, and let $\{x_1,\ldots,x_s\}=V(\mathfrak B)$  be the distinct vertices of $\mathfrak B$. Let $x\in\R^d$ and  let $\mathfrak B_d(x)\subseteq \mathfrak B_d$ consist of all cones $\C(B)\in \mathfrak B_d$ with $\supp_\mathfrak B(x)\subseteq B\subseteq V(\mathfrak B)$.
\begin{itemize}
\item[1.] If $\C(B)\in \mathfrak B$ with  $\supp_{\mathfrak B}(x)\nsubseteq B\subseteq V(\mathfrak B)$,  then $\C(B)\cap \R_+x=\{0\}$.
\item[2.]
 $x\in \mathsf{Int}(\bigcup_{C\in \mathfrak B_d(x)}C)$ with $x\notin \mathsf{Int}(\bigcup_{C\in Y}C)$ for any proper subset $Y\subset \mathfrak B_d(x)$.
 \item[3.]$\bigcap_{C\in \mathfrak B_d(x)}C=\C(\supp_\mathfrak B(x))$.
\item[4.] There is a sufficiently small $\epsilon>0$ (dependent on $\mathfrak B$) such that, for any  $y_1,\ldots,y_s\in \R^d\setminus \{0\}$ with $\mathsf d\big(x_i/\|x_i\|,y_i/\|y_i\|\big)<\epsilon$ for all $i$, the map $\varphi: \mathfrak B\rightarrow \mathfrak B'$ given by  $\varphi(C)=\C(\{y_i\}_{i\in I})$ for $C=\C(\{x_i\}_{i\in I})\in \mathfrak B$, where $I\subseteq [1,s]$, is a simplicial isomorphism between $\mathfrak B$ and $\mathfrak B':=\{\varphi(C):\; C\in \mathfrak B\}$. In particular, $\mathfrak B'$ is a complete simplicial fan for $\R^d$.

\item[5.] If $X=\{x_1,\ldots,x_{s}\}$ is a minimal positive basis for $\mathcal E\subseteq \R^d$, then there is an $\epsilon>0$ such that any set $\{y_1,\ldots,y_{d+1}\}\subseteq \mathcal E$ with $\mathsf d(x_i/\|x_i\|,y_i/|y_i\|)<\epsilon$ for all $i$ is also a minimal positive basis for $\mathcal E$.

    \item[6.] Let  $\mathfrak B'$ and $\varphi$ be  defined as in Item 4. Then, for  all sufficiently small $\epsilon>0$ (dependent on $\mathfrak B$ and $x$), $\varphi\big(\supp_{\mathfrak B}(x)\big)\subseteq \supp_{\mathfrak B'}(x)$.
\end{itemize}
\end{proposition}

\begin{proof}
1. Per definition of a complete simplicial fan, we have a disjoint decomposition of $\R^d$ given by $\R^d=\biguplus_{C\in \mathfrak B}C^\circ$. Each polyhedral cone $C\in \mathfrak B$ corresponds uniquely to some linearly independent subset of vertices $B_C\subseteq V(\mathfrak B)$ with $\C(B_C)=C$, with the faces of the cone $C$ corresponding to the subsets of $B_C$. Thus $C'\subseteq C$ is a face when $C'=\C(B_{C'})$ with $B_{C'}\subseteq B_C$. By definition, the cone $\C(B_x)\in\mathfrak B$, where $B_x:=\supp_{\mathfrak B}(x)$, is the unique cone in $\mathfrak B$ that contains $\R^\circ_+x$ in its relative interior. If $\C(B)\in \mathfrak B$ is a cone that  contains $\R_+x$, then $x$ will be contained in the relative interior of some face of $\C(B)$. Consequently, since all such faces of $\C(B)$ lie in $\mathfrak B$,  the uniqueness of $\C(B_x)$ ensures that this face must be $\C(B_x)$, i.e., $B_x\subseteq B$.

2. If $x=0$, then $\supp_\mathfrak B(0)=\emptyset$ and $\mathfrak B_d(0)=\mathfrak B_d$,  in which case the statement follows trivially in view of $\mathfrak B$ being complete.  Otherwise, in view of Item 1, it follows that there is a small  neighborhood around each cone $C\in \mathfrak B_d\setminus \mathfrak B_d(x)$ with the property that the closure of this neighborhood does not contain  $x$. Thus this is also true for the finite union $\bigcup_{C\in \mathfrak B_d\setminus \mathfrak B_d(x)}C$, implying $x\in \mathsf{Int}(\bigcup_{C\in \mathfrak B_d(x)}C)$ in view of  $\bigcup_{C\in \mathfrak B_d}C=\R^d$ (as $\mathfrak B$ is complete). Recall that the relative interiors of the cones from $\mathfrak B$ form a disjoint partition of $\R^d$ with $x$ contained in the relative interior of the cone $C_x:=\C(\supp_\mathfrak B(x))$.
Let $C'\in \mathfrak B_d$ be arbitrary. By definition of $\mathfrak B_d$, we have $C_x\subseteq C'$, so  $x\in C_x\subseteq \overline{C'}=\overline{(C')^\circ}$, with the latter equality following since $C'$ is a convex set. It follows that $B_\epsilon(x)\cap (C')^\circ\neq \emptyset$ for all $\epsilon>0$, and since $(C')^\circ$ is disjoint from $\bigcup_{C\in \mathfrak B_d(x)\setminus \{C'\}}C$ (as the relative interiors of the cones in $\mathfrak B$ form a disjoint partition of $\R^d$), it follows that $B_\epsilon(x)\nsubseteq \bigcup_{C\in \mathfrak B_d(x)\setminus \{C'\}}C$ for all $\epsilon>0$, showing that  $x\notin \mathsf{Int}(\bigcup_{C\in \mathfrak B_d(x)\setminus \{C'\}}C)$. Since $C'\in \mathfrak B_d(x)$ was arbitrary, this implies $x\notin \mathsf{Int}(\bigcup_{C\in Y}C)$ for any proper subset $Y\subset \mathfrak B_d(x)$.

3. Let $C_x=\C(\supp_\mathfrak B(x))$. By definition of $\mathfrak B_d(x)$, we have $C_x\subseteq \bigcap_{C\in \mathfrak B_d(x)}C$. If the reverse inclusion fails, then there must be some vertex $z\notin B_x:=\supp_\mathfrak B(x)$ contained in  every $C\in \mathfrak B_d(x)$. In particular, $B_x\cup \{z\}$ is linearly independent (as the generating vertices of each cone are linearly independent) and linearly spans some subspace $\mathcal E$. For any $C=\C(B)\in \mathfrak B_d(x)$, we have $B_x\cup \{z\}\subseteq B$ with the vertices in $B\subseteq V(\mathfrak B)$ linearly independent. Thus $C\cap \mathcal E=\C(B_x\cup \{z\})$. As this is true for every $C\in \mathfrak B_d(x)$, we conclude that $\left(\bigcup_{C\in \mathfrak B_d(x)}C\right)\cap \mathcal E=\C(B_x\cup \{z\})$.
Consequently, since $B_x\cup \{z\}$ is linearly independent with $x\in \C^\circ(B_x)$, it follows that   $\left(\bigcup_{C\in \mathfrak B_d(x)}C\right)\cap \C(-z,x)=\C(B_x\cup\{z\})\cap \C(-z,x)=\R_+x$, ensuring that $x$ is \emph{not} contained in the interior of $\bigcup_{C\in \mathfrak B_d(x)}C$ (else all points $x-\alpha z\in \C(-z,x)$ with $\alpha>0$ sufficiently small would  be contained in $\bigcup_{C\in \mathfrak B_d(x)}C$), contrary to  Item 2.

4. For each vertex $x_j$, there are only a finite number of  $\C(B)\in \mathfrak B$ with $x_j\in B\subseteq V(\mathfrak B)$, defining a finite number of subspaces linearly spanned by the sets $B\setminus \{x_j\}$. Each such subspace $\mathcal H$ does not contain $x_j$, as the elements of $B$ are linearly independent, so there is some finite positive distance (along the unit sphere) between $x_j/\|x_j\|$ and any such subspace $\mathcal H$ intersected with the unit sphere. Let $\epsilon '>0$ be the minimum such distance, where the minimum runs over all possible vertices  $x_j\in V(\mathfrak B)$ and all possible subspace pairings. Given a subspace $\mathcal H$ generated by linearly independent unit vectors, if we perturb the generators of $\mathcal H$ by a small amount, replacing each with a new unit vector generator some sufficiently small distance $\epsilon>0$ from the original vector, then the resulting set of perturbed generators will generate a perturbed subspace $\mathcal H'$ of equal dimension which, when intersected with the unit sphere, has all its points some small distance away from the original hyperplane  $\mathcal H$ intersected with the sphere. Thus, by this continuity property, we may choose $\epsilon>0$ sufficiently small with respect to $\epsilon'>0$ to ensure that each perturbed point $y_j/\|y_j\|$ remains disjoint from each perturbed paired subspace $\mathcal H'$ (ensuring that each cone in $\mathfrak B'$ is still generated by linearly independent elements) and on the same side (in the case of the maximal co-dimension $1$ subspaces) as the original vector $x_j/\|x_j\|$. By doing so, we ensure that the simplicial structure of $\mathfrak B$ is preserved in $\mathfrak B'$, and the result follows.
%
%
See also  \cite[Section 5.2]{Buchstaber-Panov}.

5. This is a special case of Item 4 in view of Proposition \ref{prop-char-minimal-pos-basis}.

6. If $x=0$, we have $\supp_\mathfrak B(x)=\emptyset=\supp_{\mathfrak B'}(x)$, and the result is clear. Therefore we may assume $x$ is nonzero. Now, if $x$ is a positive scaler multiple of a vertex of $\mathfrak B$, then we can assume by rescaling $x$ that $x\in V(\mathfrak B)$. On the other hand, if $x$ is not a positive scaler multiple of any vertex of $\mathfrak B$, then we can perform a baricentric subdivision at $x$ in $\mathfrak B$ to create a new complete simplicial fan $\mathfrak C$ having $V(\mathfrak C)=V(\mathfrak B)\cup \{x\}$
(so we remove each  $\C(B)$  with $\supp_\mathfrak B(x)\subseteq B\subseteq V(\mathfrak B)$ and replace it with the collection of cones of the form $\C(B\setminus \{y\}\cup \{x\})$ for $y\in \supp_\mathfrak B(x)$).
In the former case (when $x\in V(\mathfrak B)$), set $\mathfrak C=\mathfrak B$. By item 4 applied to $\mathfrak C$, for sufficiently small $\epsilon>0$, replacing each vertex of $\mathfrak C$ with a new vertex at radial distance at most $\epsilon$ from the original vector results in a new complete simplicial fan $\mathfrak C'$ isomorphic to $\mathfrak C$. Let $\mathfrak B'\subseteq \mathfrak C'$ be the complete simplicial fan associated to the image of the original vertex set $\varphi\big(V(\mathfrak B)\big)$, in which case $\varphi(C)\in \mathfrak B'$ for $C\in \mathfrak B$. In view of Item 2 applied to $\mathfrak B$, we have $x\in \mathsf{Int}(\bigcup_{C\in \mathfrak B_d(x)}C)$, whence  \be\label{starpower}x\in \mathsf{Int}(\bigcup_{C\in \mathfrak B_d(x)}\varphi(C))\ee in view of our choice of $\epsilon>0$. By Item 2 applied to $\mathfrak B'$, we know that $x\in\mathsf{Int}(\bigcup_{\varphi(C)\in \mathfrak B'_d(x)}\varphi(C))$ with this \emph{failing} for any proper subset of $\mathfrak B'_d(x)$. Combining this with \eqref{starpower} and Item 1,
we conclude that $$\mathfrak B'_d(x)\subseteq \varphi(\mathfrak B_d(x)).$$ Consequently, since every element of $\mathfrak B_d(x)$ contains the set $\supp_\mathfrak B(x)$, it follows that \be\label{grob}\C\Big(\varphi(\supp_\mathfrak B(x))\Big)\subseteq \bigcap_{\varphi(C)\in \mathfrak B'_d(x)}\varphi(C)=\C\Big(\supp_{\mathfrak B'}(x)\Big),\ee with the final equality above in view of  Item 3 applied to $\mathfrak B'$. However, \eqref{grob} is equivalent to  $\varphi(\supp_\mathfrak B(x))\subseteq \supp_{\mathfrak B'}(x)$, which completes the proof.
\end{proof}

\subsection{$\mathcal F$-filtered Sequences, Minimal Encasement and Reay Systems} Next, we extend the concept of an asymptotically filtered sequence, which we  introduced in Section \ref{sec-asym-seq}. Note,  if we specialize below to the case when $w_i^{(j)}=0$ for all  $j\in [1,\ell]$ and $i\geq 1$ with $\mathcal E_j=\R\la u_1,\ldots,u_j\ra$ for $j\in [1,\ell]$, then we recover the notion of an asymptotically filtered sequence. We extend much of the terminology introduced in Section \ref{sec-asym-seq} from asymptotically filtered sequences to $\mathcal F$-filtered sequences.
\begin{definition}
Let $\vec u=(u_1,\ldots,u_\ell)$ be a tuple of $\ell\geq 0$ orthonormal vectors in $\R^d$ and let $\{0\}=\mathcal E_0\subset \mathcal E_1\subset\ldots\subset \mathcal E_\ell\subseteq \R^d$ be a chain of subspaces such that $u_j\in \mathcal E_j\cap \mathcal E_{j-1}^\bot$ for all $j\in [1,\ell]$. A sequence  $\{x_i\}_{i=1}^\infty$ of terms $x_i\in \R^d$ is an $\mathcal F$-\textbf{filtered} sequence with  filter $\mathcal F=(\mathcal E_1,\ldots,\mathcal E_\ell)$ and limit $\vec u$ if
$$x_i= ({a}_i^{(1)}u_1+w_{i}^{(1)})+\ldots+(a_i^{(\ell)} u_{\ell}+w_i^{(\ell)})+y_i\quad\mbox{ for all $i\geq 1$,}$$
for some real numbers  $a_i^{(j)}>0$,  vectors  $u_j,\,w_i^{(j)}\in \mathcal E_j\cap \mathcal E_{j-1}^\bot$,  and  $y_i\in \mathcal E_\ell^\bot$ \ such that
\begin{itemize}
\item $\lim_{i\rightarrow \infty} a_i^{(j)}\in \R_+\cup \{\infty\}$ exists for each $j\in [1,\ell]$,
\item  $\|y_i\|,\,\|w_i^{(j)}\|\in o(a_i^{(j)})$ for all $j\in [1,\ell]$, and $a_i^{(j+1)}\in o(a_i^{(j)})$ for all $j\in [1,\ell-1]$.
\end{itemize}
\end{definition}

%
%

 Let $\vec u=(u_1,\ldots,u_t)$ be a tuple of orthonormal vectors $u_i\in \R^d$ and let $\mathcal F=(\mathcal E_1,\ldots,\mathcal E_\ell)$ be a tuple of subspaces with  $$\{0\}=\mathcal E_0\subset \mathcal E_1\subset \ldots\subset \mathcal E_\ell\quad\und\quad\R\la u_1,\ldots,u_t\ra\subseteq\mathcal E_\ell.$$ Then,  for each $i\in [1,t]$,  there is a unique $j\in [1,\ell]$ with $u_i\in \mathcal E_j\setminus \mathcal E_{j-1}$.
 Let $J\subseteq [1,\ell]$ consist of all indices $j$ for which $\mathcal E_j\setminus \mathcal E_{j-1}$ contains some $u_i$ with $i\in [1,t]$, and  for each $j\in J$, let $r_j\in [1,t]$ be the minimal index with $u_{r_j}\in \mathcal E_j\setminus \mathcal E_{j-1}$. If $J=[1,\ell]$ and $r_i\leq r_j$ holds whenever $i\leq j$, for $i,\,j\in [1,\ell]$,
 then we say that $\mathcal F$ is a \textbf{compatible filter} for $\vec u$,
 define $$\mathcal F(\vec u)=(\overline u_{1},\ldots, \overline u_{\ell}),$$ where $\overline u_j=\pi_{j-1}(u_{r_j})/\|\pi_{j-1}(u_{r_j})\|$ with $\pi_{j-1}:\R^d\rightarrow \mathcal E_{j-1}^\bot$ the orthogonal projection, and set $r_{\ell+1}=t+1$,
 in which case $$1=r_1<r_2<\ldots<r_\ell<r_{\ell+1}=t+1,$$ which refer to as the \textbf{associated indices} for $\mathcal F(\vec u)$. Note $\ell\geq 1$ except when $\vec u$ is the empty tuple, in which case $\mathcal F(\vec u)$ is the empty tuple.

\begin{proposition}\label{prop-compat-filter}
Let $\vec u=(u_1,\ldots,u_t)$ be a tuple of $t\geq 0$ orthonormal vectors $u_i\in \R^d$ and  let  $\mathcal F=(\mathcal E_1,\ldots,\mathcal E_\ell)$ be a compatible filter for $\vec u$  with  $1=r_1<\ldots<r_\ell<r_{\ell+1}=t+1$ the associated indices. Then  $$u_i\in \mathcal E_{j-1}\quad\mbox{ for all $i<r_{j}$ and $j\leq \ell+1$}.$$ Moreover, if
 $\{x_i\}_{i=1}^\infty$ is an asymptotically   filtered sequence with limit $\vec u$, say with $x_i=a_i^{(1)}u_1+\ldots+a_i^{(t)}u_t+y_i$, then $\{x_i\}_{i=1}^\infty$ is an $\mathcal F$-filtered sequence with limit $\mathcal F(\vec u)$, say with $$x_i=(\alpha_i^{(1)}\overline u_1+w_i^{(1)})+\ldots+(\alpha_i^{(\ell)}\overline u_\ell+w_i^{(\ell)})+y'_i,$$ where $\alpha_i^{(j)}\in \Theta(a_i^{(r_j)})$ for $j\in [1,\ell]$ and $\|y'_i\|\in O(\|y_i\|)$.
\end{proposition}

\begin{proof}
If $u_i\notin\mathcal E_{j-1}$ with $i<r_{j}$, then $u_i\in \mathcal E_k\setminus \mathcal E_{k-1}$ for some $k\geq j$ (as $\R\la u_1,\ldots,u_t\ra\subseteq \mathcal E_\ell$), whence $r_k\leq i<r_{j}$, contradicting that $\mathcal F$ is a compatible filter for $\vec u$ in view of $k\geq j$. Therefore we instead conclude that $u_i\in\mathcal E_{j-1}$ for $i<r_{j}$.
 For $j\in [0,\ell]$, let $\pi_j:\R^d\rightarrow \mathcal E_j^\bot$ and $\pi^\bot_j:\R^d\rightarrow \mathcal E_j$ be the orthogonal projections, where $\mathcal E_0=\{0\}$. We may assume $\ell\geq 1$, for $\ell=0$ implies $t=0$, in which case $x_i=y_i=y'_i$ is trivially an $\mathcal F$-filtered sequence.
 By definition of the $\overline u_j$, we have $\overline u_j\in \mathcal E_j\cap \mathcal E_{j-1}^\bot$ for $j\in [1,\ell]$.
 For $j\in [1,\ell]$, let $y_i^{(j)}=\pi_j^\bot(x_i)-\pi_{j-1}^\bot(x_i)$, and set $y'_i=\pi_\ell(x_i)=\pi_\ell(y_i)\in \mathcal E_\ell^\bot$ (since $\R\la u_1,\ldots,u_t\ra\subseteq \mathcal E_\ell=\ker \pi_\ell$). Then $x_i=y_i^{(1)}+\ldots+y_i^{(\ell)}+y'_i$ with each $y_i^{(j)}\in \mathcal E_j\cap \mathcal E_{j-1}^\bot$ and $y'_i\in \mathcal E_\ell^\bot$.
  Let $\alpha_i^{(j)}=a_i^{(r_j)}\|\pi_{j-1}(u_{r_j})\|>0$ (since $u_{r_j}\notin \mathcal E_{j-1}$) and let $$w_i^{(j)}=y_i^{(j)}-\alpha_i^{(j)}\overline u_{j}=y_i^{(j)}-\pi_{j-1}(a_i^{(r_j)}u_{r_j})=
  \pi_{j-1}(y_i^{(j)}-a_i^{(r_j)}u_{r_j})
  \in \mathcal E_j\cap \mathcal E_{j-1}^\bot$$ (since $y_i^{(j)}\in \mathcal E_j\cap \mathcal E_{j-1}^\bot$ and  $u_{r_j}\in \mathcal E_j$ by definition of $r_j$). Thus $$x_i=(\alpha_i^{(1)}\overline u_1+w_i^{(1)})+\ldots+(\alpha_i^{(\ell)}\overline u_\ell+w_i^{(\ell)})+y'_i,$$ with $\overline u_j,\,w_i^{(j)}\in \mathcal E_j\cap \mathcal E_{j-1}^\bot$ and $y'_i\in \mathcal E_\ell^\bot$.
  By definition of $r_j$, we are assured that $\|\pi_{j-1}(u_{r_j})\|>0$. Thus $\alpha_i^{(j)}\in \Theta(a_i^{(r_j)})$ and $\|y'_i\|=\|\pi_\ell(y_i)\|\in O(\|y_i\|)$ (since linear operators between finite dimensional spaces are bounded).
  Consequently, since $x_i=a_i^{(1)}u_1+\ldots+a_i^{(t)}u_t+y_i$ is asymptotically filtered, we have $\alpha_i^{(j)}\in o(a_i^{(r_j-1)})\subseteq o(a_i^{(r_{j-1})})=o(\alpha_i^{(j-1)})$ and $\|y'_i\|\in O(\|y_i\|)\subseteq o(a_i^{(t)})\subseteq o(a_i^{(r_\ell)})=o(\alpha_i^{(\ell)})$.
Since $u_i\in \mathcal E_{j-1}\subset \mathcal E_j$ for $i<r_j$, and since $u_{r_j}\in \mathcal E_j\setminus \mathcal E_{j-1}$, it follows that   $$w_i^{(j)}=\pi_{j-1}(y_i^{(j)}-a_i^{(r_j)}u_{r_j})=\pi_{j-1}\pi_j^\bot(x_i-a_i^{(r_j)}u_{r_j})=
\pi_{j-1}\pi_j^\bot(a_i^{(r_j+1)}u_{r_j+1}+\ldots+
  a_i^{(t)}u_t+y_i).$$
  In consequence, if $r_j<t$, then $\|w_i^{(j)}\|\in O(a_i^{(r_j+1)})\subseteq o(a_i^{(r_j)})=o(\alpha_i^{(j)})$, and if $r_j=t$, then $j=\ell$ and $\|w_i^{(j)}\|\in O(\|y_i\|)\subseteq o(a_i^{(t)})=o(a_i^{(r_j)})=o(\alpha_i^{(j)})$, and the proof is complete in either case.
\end{proof}

\begin{lemma}\label{lemma-matrix-assymp-solutions}
Let $X\subseteq \R^d$ be a linearly independent subset, let $\{x_i\}_{i=1}^\infty$ be a sequence of terms  $x_i\in \R\la X\ra$, and let $x_i=\Summ{x\in X}\alpha_i^{(x)}x$ with $\alpha_i^{(x)}\in \R$ for $i\geq 1$. Then $|\alpha_i^{(x)}|\in O(\|x_i\|)$ for all $x\in X$.
\end{lemma}

\begin{proof} Let $|X|=s\leq d$ (as $X$ is linearly independent). Let $M'$ be the $d\times s$ matrix with column vectors the elements from $X$. Since $X$ is linearly independent, the matrix $M'$ has rank $s$, allowing us to add an additional $d-s$ columns to the right of $M'$ to create an invertible $d\times d$ matrix $M$.  Since $x_i\in \R\la X\ra$, for each $i$, there is a  vector $y_i=(\alpha_{i,1},\ldots,\alpha_{i,s},0,\ldots,0)$ with $My_i=x_i$. Then $\|y_i\|=\|M^{-1}x_i\|\leq \|M^{-1}\|\; \|x_i\|$, where $\|M^{-1}\|$ is the matrix operator norm induced by the Euclidean  $L_2$-norm. This shows $\alpha_{i,j}^2\leq \alpha_{i,1}^2+\ldots+\alpha_{i,s}^2=\|y_i\|^2\leq C^2 \|x_i\|^2$ for each $j\in [1,s]$, where $C=\|M^{-1}\|>0$, implying $|\alpha_{i,j}|\leq C\|x_i\|$ for all $i$ and $j\in [1,s]$, as desired.
\end{proof}

The next proposition links  minimal encasement, Reay systems and $\mathcal F$-filtered sequences.

\begin{proposition}\label{prop-min-encasement-char}
Let $X\subseteq\R^d$, and let $\vec u=(u_1,\ldots,u_t)$ be a tuple of $t\geq 0$ orthonormal  vectors $u_i\in \R^d$, where $d\geq 0$. Then $-X$ minimally encases $\vec u$ if and only if
\begin{itemize}
\item[1.]
there exists a disjoint partition  $X=\bigcup_{i=1}^\ell X_i$ such that $\mathcal F=(\mathcal E_1,\ldots,\mathcal E_\ell)$ is a compatible filter for $\vec u$, where $\mathcal E_j=\R\la X_1\cup \ldots\cup  X_{j}\ra$ for $j\in [1,\ell]$, and
\item[2.] $(X_1\cup \{u_{r_1}\},\ldots ,X_\ell\cup \{u_{r_\ell}\})$ is a Reay system, where  $1=r_1<\ldots<r_\ell<r_{\ell+1}=t+1$ are the associated indices for $\mathcal F(\vec u)$.
\end{itemize}
Moreover, the  Reay system $(X_1\cup \{u_{r_1}\},\ldots,X_\ell\cup \{u_{r_\ell}\})$ satisfying Items 1 and 2 is unique.
\end{proposition}

\begin{proof}
If $t=0$, then only the empty set $X=\emptyset$  minimally encases the empty tuple, and the empty partition with $\ell=0$ satisfies the desired conditions. Likewise, if such a partition exists for a set $X$ when $t=0$, then $\ell=0$ follows, implying that $X=\emptyset$. Thus we may assume $t\geq 1$.

Suppose Items 1 and 2 hold.  Item 2 allows us to apply Proposition \ref{prop-reay-basis-properties} to conclude that  $-u_{r_j}\in \R\la X_1\cup \ldots\cup X_j\ra=\C(X_1\cup \{u_{r_1}\}\cup\ldots \cup X_j\cup \{u_{r_j}\})$ for every $j\in [1,\ell]$,
whence $-u_{r_j}=z+a_1u_{r_1}+\ldots+a_{j-1}u_{r_{j-1}}$ for some $z\in \C(X_1\cup \ldots\cup X_j)$ and $a_i\geq 0$, implying that $u_{r_j}+a_1u_{r_1}+\ldots+a_{j-1}u_{r_{j-1}}\in \C(-X_1\cup\ldots\cup -X_j)\subseteq \C(-X)$. By item 1 and Proposition \ref{prop-compat-filter}, we have $u_i\in \mathcal E_{j-1}$ whenever $i<r_j$ and $j\leq \ell+1$. Thus, for any $i\in [1,t]$ with $r_{j}<i<r_{j+1}$ and $j\in [1,\ell]$, we have  $u_i\in \R \la X_1\cup\ldots\cup X_{j}\ra=\C(X_1\cup\{u_{r_1}\}\ldots\cup X_{j}\cup \{u_{r_j}\})$, so that a similar argument yields $u_i+a_1u_{r_1}+\ldots+a_ju_{r_j}\in \C(-X_1\cup\ldots\cup -X_j)\subseteq \C(-X)$ for some $a_i\geq 0$. This shows that $-X$ encases $\vec u=(u_1,\ldots,u_{t})$. It remains to establish the minimality of the encasement. To this end, it suffices to show, for an arbitrary $x\in X$, that $\C(-X\setminus\{-x\})$ does not encase $(u_1,\ldots,u_{t})$. Let $j\in [1,\ell]$ be the index such that $x\in X_j$ and let $\pi_{j-1}:\R^d\rightarrow \mathcal E_{j-1}^\bot$ be the orthogonal projection. Then $\mathcal R=(\pi_{j-1}(X_j\cup \{u_{r_j}\}),\ldots,\pi_{j-1}(X_\ell\cup \{u_{r_\ell}\}))$ is a Reay system and $\supp_\mathcal R(-\pi_{j-1}(u_{r_j}))=\pi_{j-1}(X_j)$. Assume by contradiction that  $-X\setminus \{-x\}$ encases $(u_1,\ldots,u_{t})$. Then, since $\pi_{j-1}(u_i)=0$ for all $i<r_j$ (by Proposition \ref{prop-compat-filter}), it follows that  $-\pi_{j-1}(u_{r_j})\in\C(\pi_{j-1}(X\setminus \{x\}))=\C(\pi_{j-1}(X_j\setminus \{x\})\cup\ldots\cup \pi_{j-1}(X_\ell))$, implying that $\pi_{j-1}(x)\notin \supp_\mathcal R(-\pi_{j-1}(u_{r_j}))$. However, this contradicts that $\supp_\mathcal R(-\pi_{j-1}(u_{r_j}))=\pi_{j-1}(X_j)$ with $x\in X_j$. So we conclude that $-X$ minimally encases $\vec u$, as desired.

Next suppose $-X$ minimally encases $\vec u$. Then  $X$ must be finite (by the minimality of $-X$), and  there are vectors $v_1,\ldots,v_{t}\in \C(-X)$ with each \be\label{treewindow} v_j=\alpha_{1,j}u_1+\ldots+\alpha_{j,j}u_j\in \C(-X)\ee for some real numbers $\alpha_{i,j}\geq 0$ with $\alpha_{j,j}>0$.
In particular, $u_1\in \C(-X)$, so there must be  a subset $X_1\subseteq X$ such that $X_1\cup \{u_1\}$ is a minimal positive basis (in view of Proposition \ref{prop-char-minimal-pos-basis}.4 and Carth\'eordory's Theorem).
Let $r_1=1$, let $\mathcal E_0=\{0\}$, and let $\mathcal E_1=\R\la X_1\ra$.
 We proceed to recursively construct, for $j=1,2,\ldots,$  nonempty subsets $X_1,\ldots,X_j\subseteq X$, subspaces $\mathcal E_j=\R\la X_1\ldots\cup X_j\ra$, and indices $1=r_1<r_2<\ldots<r_j\leq t$ such that $r_j\in [1,t]$ is the minimal index with $u_{r_j}\notin \mathcal E_{j-1}$, and $(X_1\cup \{u_{r_1}\},\ldots, X_j\cup \{u_{r_j}\})$ is a Reay system. We have just shown this is possible for $j=1$, so assume $j\geq 2$ and that the sets $X_1,\ldots,X_{j-1}\subseteq X$ and indices $1=r_1<\ldots<r_{j-1}$ have already been found such that $(X_1\cup \{u_{r_1}\},\ldots,X_{j-1}\cup \{u_{r_{j-1}}\})$ is a Reay system  and each $r_{i}\in [1,t]$ is the minimal index with $u_{r_i}\notin \mathcal E_{i-1}=\R\la X_1\cup \ldots\cup X_{i-1}\ra$ for $i\leq j-1$. Let $r_j\in [1,t]$ be the minimal index such that $u_{r_j}\notin \mathcal E_{j-1}=\R\la X_1\cup\ldots\cup X_{j-1}\ra$ (if such an index exists), and otherwise set $r_j=t+1$. By Proposition \ref{prop-reay-basis-properties}.1, we have $$\mathcal E_{j-1}=\R\la X_1\cup \ldots\cup X_{j-1}\ra=\C(X_1\cup \{u_{r_1}\}\cup \ldots\cup X_{j-1}\cup \{u_{r_{j-1}}\}).$$ Thus $r_j>r_{j-1}$. If $r_j=t+1$, then $-u_i\in \mathcal E_{j-1}=\C(X_1\cup \{u_{r_1}\}\cup \ldots\cup X_{j-1}\cup \{u_{r_{j-1}}\})$ for all $i\in [1,t]$. In such case, for each $i>r_{j-1}$, we have $u_i+\alpha_{i,1}u_{r_1}+\ldots \alpha_{i,j-1}u_{r_{j-1}}\in \C(-X_1\cup \ldots\cup -X_{j-1})$ for some $\alpha_{i,r_1},\ldots,\alpha_{i,r_{j-1}}\in \R_+$.

 If no such index $r_2$ exists (so $r_2=t+1$), then $-u_i\in \mathcal E_1=\C(X_1\cup \{u_1\})$ for all $i\geq 1$, meaning each $-u_i=b_i+a_iu_1$ for some $a_i\geq 0$ and $b_i\in \C(X_1)$.
But then $u_i+a_iu_1\in \C(-X_1)$ for all $i\geq 2$, implying that $\C(-X_1)$ encases $(u_1,\ldots,u_{t})$.  In this case, the minimality of $X$ ensures that $X=X_1$, and the desired partition of $X$ follows with $\ell=1$ and $\mathcal F=(\mathcal E_1)$ compatible with $\vec u$ (as $u_1\in \mathcal E_1\setminus \{0\}$ with $\R\la u_1,\ldots,u_t\ra\subseteq \mathcal E_1$). So we may now assume the index $r_2$ exists.  Moreover, the previous argument ensures that $\C(-X_1)$ encases $(u_1,\ldots,u_{r_2-1})$.
Now  let $\pi_1:\R^d\rightarrow \mathcal E_1^\bot$ be the orthogonal projection. Since $u_1,\ldots,u_{r_2-1}\in \mathcal E_1$, we have  $-\pi_1(u_{r_2})\in \C(\pi_1(X))$ by \eqref{treewindow}. Thus, as before, we can find a subset $X_2\subseteq X\setminus X_1$ such that $|\pi_1(X_2)|=|X_2|$ and $\pi_1(X_2\cup \{u_{r_2}\})$ is a minimal positive basis, in which case $(X_1\cup \{u_{r_1}\},X_2\cup \{u_{r_2}\})$ is a Reay system.
Since $\pi_1(X_2)\cup \{\pi_1(u_{r_2})\}$ is a minimal positive basis, we have  $u_{r_2}+b\in \C^\circ(-X_2)$ for some $b\in \ker \pi_1=\mathcal E_1=\C(X_1\cup \{u_{r_1}\})$, in turn implying $u_{r_2}+a_{r_2}u_{r_1}\in \C(-X_1\cup -X_2)$ for some $a_{r_2}\geq 0$. Thus $-X_1\cup -X_2$ encases $(u_1,\ldots,u_{r_2})$.
Let $\mathcal E_2=\R\la X_1\cup X_2\ra=\C(X_1\cup X_2\cup \{u_{r_1},u_{r_2}\})$ (in view of Proposition \ref{prop-reay-basis-properties}.1) and let $r_3$ be the minimal index such that $u_{r_3}\notin \mathcal E_2$.
If no such index $r_3$ exists, then $-u_i\in \mathcal E_2=\C(X_1\cup X_2\cup \{u_{r_1},u_{r_2}\})$ for all $i\geq 1$, meaning each $-u_i=c_i+b_iu_{r_2}+a_iu_{r_1}$ for some $a_i,\,b_i\geq 0$ and $c_i\in \C(X_1\cup X_2)$. But then $u_i+a_iu_{r_1}+b_iu_{r_2}\in \C(-X_1\cup -X_2)$ for all $i>r_2$, implying that $\C(-X_1\cup -X_2)$ encases $(u_1,\ldots,u_{t})$. As before, the minimality of $X$ then ensures that $X=X_1\cup X_2$, and the desired partition follows with $\ell=2$. So we may now assume the index $r_3$ exists. Moreover, $\C(-X_1\cup -X_2)$ encases $(u_1,\ldots,u_{r_3-1})$. Continuing to iterate these arguments (as in Proposition \ref{prop-reay-basis-exists}) now leads to the desired partition of $X$ after $\ell\leq t$ steps.

Finally, it remains to show $\mathcal R=(X_1\cup \{u_{r_1}\},\ldots,X_\ell\cup \{u_{r_\ell}\})$ is unique, which we handle by induction on $\ell$.  To this end, suppose $\mathcal R'=(X'_1\cup \{u_{r'_1}\},\ldots, X'_{\ell'}\cup \{u_{r'_{\ell'}}\})$ is another Reay system satisfying Items 1 and 2, so $1=r'_1<\ldots<r'_{\ell'}<r'_{\ell'+1}=t+1$. Thus $-u_1=-u_{r_1}\in \C^\circ(X_1)$ and $-u_1=-u_{r'_1}\in \C^\circ(X'_1)$ by Proposition \ref{prop-char-minimal-pos-basis}.3. Let $\pi_1:\R^d\rightarrow \mathcal E_1^\bot$ be the orthogonal projection. Since $\mathcal R$ is a Reay system, it follows that  $\pi_1(X_2)\cup \ldots\cup \pi_1(X_\ell)$ is a linearly independent set of size $|X_2|+\ldots+|X_\ell|$ (per definition of a Reay system and Proposition \ref{prop-reay-basis-properties}). Consequently, any linear combination of elements from $X$ equal to $-u_1\in \mathcal E_1=\ker \pi_1$ can only involve terms from $X_1$. Thus  $-u_1=-u_{r'_1}\in \C^\circ(X'_1)$ implies  $X'_1\subseteq X_1$. Exchanging the roles of $\mathcal R$ and $\mathcal R'$ and repeating the argument, we find that $X_1\subseteq X'_1$, whence $X_1=X'_1$. As a result, $u_{r_2}=u_{r'_2}$ now follows from Item 1 (since $\R\la X_1\ra=\R\la X'_1\ra$). This completes the base case when $\ell=1$, so we can assume $\ell\geq 2$.  Now $(\pi_1(X_2)\cup \{\pi_1(u_{r_2})\},\ldots,\pi_1(X_\ell)\cup \{\pi_1(u_{r_\ell})\})$ and $(\pi_1(X'_2)\cup \{\pi_1(u_{r'_2})\},\ldots,\pi_1(X'_{\ell'})\cup \{\pi_1(u_{r_{\ell'}})\})$ are both Reay systems (per Proposition \ref{prop-reay-basis-properties}.3) showing that $-\pi_1(X\setminus X_1)$ minimally encases $\pi_1(\vec u)$. Thus, by induction hypothesis, $\ell=\ell'$ and $\pi_1(X_j)=\pi_1(X'_j)$ for all $j\in [2,\ell]$, implying $X_j=X'_j$ for all $j\in [2,\ell]$ as $\pi_1$ is injective on $X\setminus X_1$. This shows  $\mathcal F=(\mathcal E_1,\ldots,\mathcal E_\ell)$ is uniquely defined, where $\mathcal E_j=\R\la X_1\cup \ldots\cup X_j\ra$, in which case  Item 1 ensures that  the indices $1=r_1<\ldots<r_\ell<r_{\ell+1}=t+1$ are also  uniquely defined,  i.e., $r_j=r'_j$ for all $j\in [1,\ell]$. \end{proof}

If $-X$ minimally encases $\vec u$,  we have a unique Reay system $\mathcal R=(X_1\cup \{u_{r_1}\},\ldots,X_\ell\cup \{u_{r_\ell}\})$ and  indices $1=r_1<\ldots<r_\ell<r_{\ell+1}=t+1$ given by Proposition \ref{prop-min-encasement-char}, which we refer to as the Reay system and indices \textbf{associated} to the minimal encasement of $\vec u$ by $-X$.

\begin{proposition}\label{prop-min-encasement-minposbasis}
 Let $\vec u=(u_1,\ldots,u_t)$ be a tuple of  orthonormal  vectors $u_i\in \R^d$, where $t,\,d\geq 1$, let $X\subseteq \R^d$ be a subset  minimally encasing $-\vec u$, and let $X=\bigcup_{i=1}^\ell X_i$ and $1=r_1<\ldots<r_\ell<r_{\ell+1}=t+1$ be the Reay system and indices associated to the minimal encasement.

If $\{x_i\}_{i=1}^\infty$ is an asymptotically  filtered sequence with limit $\vec u$, say $x_i=a_i^{(1)}u_1+\ldots+a_i^{(t)}u_t+y_i,$ and $y_i\in \R\la X\ra$ for all $i$, then   $X\cup \{x_i\}$ is a minimal positive basis for $\R\la X\ra$ for all sufficiently large $i$. Moreover, letting $-x_i=\Summ{x\in X}\alpha_i^{(x)}x$ be the unique positive linear combination with $\alpha_i^{(x)}>0$ (for $i$ sufficiently large), we have $\alpha_i^{(x)}\in \Theta(a_i^{(r_j)})$ for $x\in X_j$.
\end{proposition}

\begin{proof}
For  $j\in [1,\ell]$, let $\mathcal E_j=\R\la X_1\cup\ldots\cup X_j\ra$ and let $\pi_j:\R^d\rightarrow \mathcal E_j^\bot$ be the orthogonal projection.  Then $X=X_1\cup\ldots\cup X_\ell$ is linearly independent by Proposition \ref{prop-reay-basis-properties}.1, so  $X\cup\{x_i\}$ is a minimal positive basis if and only if  $x_i\in -\C^\circ(X)$ (by Proposition \ref{prop-char-minimal-pos-basis}.4). By Proposition \ref{prop-encasementcones-contain-aprox-seq}, $\{x_i-y_i\}_{i=1}^\infty$ is  asymptotically filtered  with limit $\vec u$ and $x_i-y_i\in -\C(X)$ for all sufficiently large $i$.

Let us first show that, if the proposition holds for the sequence $\{x_i-y_i\}_{i=1}^\infty$, then it holds for the sequence $\{x_i\}_{i=1}^\infty$ as well. To this end, suppose $x_i-y_i=-\Summ{x\in X}\alpha_i^{(x)}x$ for all sufficiently large $i$, for some $\alpha_i^{(x)}>0$ with $\alpha_i^{(x)}\in \Theta(a_i^{(r_j)})$ for all $x\in X_j$ and $j\in [1,\ell]$. Since $y_i\in \R\la X\ra$ with $X$ linearly independent, Lemma \ref{lemma-matrix-assymp-solutions} implies that $y_i=-\Summ{x\in X}\beta_i^{(x)}x$ for some $\beta_i^{(x)}\in \R$ with $|\beta_i^{(x)}|\in O(\|y_i\|)\subseteq o(a_i^{(t)})\subseteq o(a_i^{(r)})$ for all $x\in X$ and $r\in [1,t]$. Thus $\alpha_i^{(x)}+\beta_i^{(x)}\in \Theta(a_i^{(r_j)})$ for all $x\in X_j$ and $j\in [1,\ell]$. Moreover,  $\alpha_i^{(x)}+\beta_i^{(x)}>0$ for all sufficiently large $i$. Thus, since $x_i=(x_i-y_i)+y_i=-\Summ{x\in X}(\alpha_i^{(x)}+\beta_i^{(x)})x$, it follows that $x_i\in -\C^\circ(X)$ for all sufficiently large $i$, ensuring that the proposition holds for $\{x_i\}_{i=1}^\infty$, as desired.

We proceed by induction on the depth $\ell$ to show the proposition holds when $y_i=0$ for all $i$, which will complete the proof by what was just shown. If $\ell=1$, then $x_i=a_i^{(1)}u_1$ with $u_1\in -\C^\circ(X_1)$ (as $X_1\cup \{u_{r_1}\}$ is a minimal positive basis and $r_1=1$). Then there is a unique strictly positive linear combination $-u_1=\Summ{x\in X}\alpha_xx$, implying that $-x_i=a_i^{(1)}(-u_1)=\Summ{x\in X}(a_i^{(1)}\alpha_x)x$. Since $\alpha_x>0$ and $\alpha_i^{(1)}>0$, we have $a_i^{(1)}\alpha_x\in \Theta(a_i^{(1)})=\Theta(a_i^{(r_1)})$ with $a_i^{(1)}\alpha_x>0$, ensuring $-x_i\in \C^\circ(X_1)=\C^\circ(X)$, as desired. This completes the induction base, so now assume $\ell\geq 2$.

Let $z_i=a_i^{(r_{\ell-1}+1)}u_{r_{\ell-1}+1}+\ldots+a_i^{(t)}u_{t}$. Write $z_i=z_i^*+z_i^\bot$ with $z_i^*\in \mathcal E_{\ell-1}$ and $z_i^\bot\in \mathcal E_{\ell-1}^\bot$, so   $$x_i=a_i^{(1)}u_1+\ldots+a_i^{(r_{\ell-1})}u_{r_{\ell-1}}+z_i^*+z_i^\bot.$$
Since $\pi_{\ell-1}(X_\ell)\cup \{\pi_{\ell-1}(u_{r_\ell})\}$ is a minimal positive basis, Proposition \ref{prop-infinite-limits-proj} allows us to apply the base case to $\{\pi_{\ell-1}(-x_i)\}_{i=1}^\infty$ yielding that  $-z_i^\bot=\pi_{\ell-1}(-x_i)=\Summ{x\in X_\ell}\alpha_i^{(x)}\pi_{\ell-1}(x)$ for some $\alpha_i^{(x)}>0$ with $\alpha_i^{(x)}\in \Theta(a_i^{(r_\ell)})\subseteq o(a_i^{(r_{\ell-1})})$, for all $x\in X_\ell$ and all sufficiently large $i$. Thus \be\label{combine}\Summ{x\in X_\ell}\alpha_i^{(x)}x=-z_i^\bot+\xi_i\ee for some $\xi_i\in \mathcal E_{\ell-1}$.
 Since $z_i^*$ and $z_i^\bot$ are orthogonal, we have $\|z_i^*\|,\,\|z_i^\bot\|\leq \|z_i\|$, ensuring
 $\|z_i^*\|,\,\|z_i^\bot\|\in O(\|z_i\|)=O(a_i^{(r_{\ell-1}+1)})\subseteq o(a_i^{(r_{\ell-1})})$.
 Thus $z_i^*+\xi_i\in \mathcal E_{\ell-1}$ with $\|z_i^*+\xi_i\|\in o(a_i^{(r_{\ell-1})})$, allowing us to apply the induction hypothesis to the sequence $\{x_i-z_i^\bot+\xi_i\}_{i=1}^\infty$ to conclude
 $-x_i+z_i^\bot-\xi_i=\Summ{x\in X\setminus X_\ell}\alpha_i^{(x)}x$ for some $\alpha_i^{(x)}>0$ with $\alpha_i^{(x)}\in \Theta(a_i^{(r_j)})$, for all $x\in X_j$, all $j\in [1,\ell-1]$ and all sufficiently large $i$. Combined with \eqref{combine}, it follows that $-x_i=\Summ{x\in X}\alpha_i^{(x)}x$ with $\alpha_i^{(x)}>0$ and $\alpha_i^{(x)}\in \Theta(a_i^{(r_j)})$, for $x\in X_j$ and  sufficiently large $i$, showing that $x_i\in -\C^\circ(X)$, which completes the induction and the proof.
 \end{proof}

\section{Oriented Reay Systems}
\label{sec-oriented}

If $\mathbf x\subseteq \R^d$ is a half-space inside the subspace $\R\la \mathbf x\ra$ with partial boundary (meaning $\x$ is obtained from a closed half-space in $\R\la\x\ra$ by removing elements from  the boundary subspace), so $\mathbf x^\circ \subseteq \mathbf x\subseteq \overline{\mathbf x}$ with $\overline{\mathbf x}=\overline {\x^\circ}$ a closed half-space in $\R\la \mathbf x\ra$ and $\mathbf x^\circ$  an open half space, then we call $\mathbf x$ a  \textbf{(relative) half-space}, though we will henceforth simply refer to such sets as  half-spaces for brevity. With regards to Convex Geometry, there is little difference between a nonzero point $x\in \R^d$ and the one-dimensional ray $\R_+ x$ that it defines, which is a one-dimensional half-space. In this way, we informally view a higher dimensional half-space $\x\subseteq \R^d$ as a type of higher-dimensional element. We will later see that, in many ways, half-spaces share similar behavior with ordinary elements. An element $x\in \mathbf x^\circ$ is called a \textbf{representative} for the relative half-space $\mathbf x$. Thus $\overline \x=\partial(\x)+\R_+x$ for any representative $x$.
Recall that a simplicial cone is set of the form $\C(X)$ with $X\subseteq \R^d$ a linearly independent set.
If $\pi:\R^d\rightarrow \R^d$ is a linear transformation with $\ker \pi\cap \R\la \mathbf x\ra\subseteq \partial(\mathbf x)$, or  equivalently, with $\mathbf x\nsubseteq \ker \pi+\partial(\mathbf x)$ (as both are equivalent to $\pi(x)\neq 0$ for all $x\in \x^\circ$),  then $\pi(\mathbf x)$ will also be a relative half-space with \be\label{halfspace-lintrans}\partial(\pi(\mathbf x))=\pi(\partial(\mathbf x)),\quad\pi(\x)\cap \partial(\pi(\x))=\pi(\x\cap \partial(\x))\quad\und\quad \pi(\x^\circ)=\pi(\x)^\circ.\ee
 We define a \textbf{blunted simplicial cone} to be a set of the form $\mathcal E+\C(X)$ with $\mathcal E\subseteq \R^d$ a subspace and $X\subseteq \R^d$ a subset for which   $\pi(X)$ is a linearly independent subset of $|X|$ elements, where $\pi:\R^d\rightarrow \mathcal E^\bot$ is the orthogonal projection. In order to avoid excessive use of dummy variables, given a set $\mathcal X$ whose elements $\mathbf x\in \mathcal X$ are subsets $\x\subseteq \R^d$, we let
$$\R^\cup\la \mathcal X\ra=\R\la {\bigcup}_{\mathbf x\in \mathcal X}\mathbf x\ra,\quad\C^\cup(\mathcal X)=\C({\bigcup}_{\mathbf x\in \mathcal X}\mathbf x),\quad\und \quad \C^\cup(\mathcal X)^\circ=\Co(\bigcup_{\mathbf x\in \mathcal X}\mathbf x).$$

We adapt the convention that relative half-spaces will be denoted in boldface, e.g. $\x$, and that the corresponding non-boldface symbol denotes a fixed representative for the half-space, e.g., $x\in \x^\circ$. Likewise, a collection of relative half-spaces will be denoted in calligraphic script, e.g. $\mathcal X$, with the corresponding non-calligraphic symbol denoting a set obtained by replacing each half-space $\x\in \mathcal X$ with a fixed representative $x\in \x^\circ$, e.g., $X$ denotes a set of representatives for $\mathcal X$. Generally, we will use a representative set $X$ only in contexts where it is irrelevant which representative is chosen for each $\x\in \mathcal X$, and the representative sets will be fixed.

\begin{definition}
For $j\in [1,s]$, let $\mathcal X_j\cup\{\mathbf v_j\}$ be a subset of relative half-spaces in $\R^d$ with distinguished element $\mathbf v_j\notin \mathcal X_j$, where $d\geq 0$ and  $s\geq 0$. For $j\in [0,s]$, let $\mathcal E_j=\R^\cup \la \mathcal X_1\cup \{\mathbf v_1\}\cup \ldots \cup \mathcal X_j\cup\{\mathbf v_j\}\ra$ and let  $\pi_j:\R^d\rightarrow \mathcal E_j^\bot$ be the orthogonal projection.  Suppose, for each $j\in [1,s]$, that the following hold.
\begin{itemize}
\item[(OR1)]  For every $\mathbf x\in \mathcal X_j\cup \{\mathbf v_j\}$, we have $\partial(\mathbf x)=\R^\cup\la\mathcal B_{\mathbf x}\ra$
    and $\partial(\mathbf x)\cap \x=\C^\cup(\mathcal B_\x)$
    for some $\mathcal B_{\mathbf x}\subseteq \mathcal X_1\cup\ldots\cup \mathcal X_{j-1}$ (which we will later denote by $\partial(\{\x\})=\mathcal B_\x$).
\item[(OR2)] $\pi_{j-1}( X_j\cup \{v_j\})$ is a minimal positive basis with $|\pi_{j-1}(X_j\cup\{v_j\})|=|\mathcal X_j|+1$.
\end{itemize}
Then we call $\mathcal R=(\mathcal X_1\cup \{\mathbf v_1\},\ldots,\mathcal X_s\cup\{\mathbf v_s\})$ an \textbf{orientated Reay system} for the subspace $\mathcal E_s$.
\end{definition}

Note, in view of (OR1), that $\pi_{j-1}(\mathcal X_j\cup\{\mathbf v_j\})$ is a set of rays in $\mathcal E_{j-1}^\bot$, so (OR2) does not depend on the choice of representatives.   Using (OR1) and (OR2) and a recursive argument for $j=1,2,\ldots,s$, it follows that  $(X_1\cup\{v_1\},\ldots,X_j\cup\{v_j\})$ is a Reay system for $\mathcal E_j$ with $\mathcal E_j=\R\la X_1\cup\ldots\cup X_j\ra=\R^\cup\la \mathcal X_1\cup\ldots\cup \mathcal X_j\ra$ for all $j\in [0,s]$.
Also, (OR2) ensures that  $\pi_{j-1}(\mathbf x)\neq 0$ for every $\mathbf x\in  \mathcal X_j\cup\{\mathbf v_j\}$, whence $\mathbf x\nsubseteq \mathcal E_{j-1}= \mathcal E_{j-1}+\partial(\mathbf x)$, and thus $\mathbf x\nsubseteq \mathcal E_{i-1}+\partial(\mathbf x)$ for any
$i\leq j$ as well. Consequently, for any $j\in [1,s]$, $$\pi_{j-1}(\mathcal R):=\big(\pi_{j-1}(\mathcal X_j)\cup\{\pi_{j-1}(\mathbf v_j)\},\ldots,\pi_{j-1}(\mathcal X_s)\cup\{\pi_{j-1}(\mathbf v_s)\}\big)$$ is also an orientated Reay system in view of \eqref{halfspace-lintrans}, while it
is clear from the recursive nature of the  definition that $(\mathcal X_1\cup\{\mathbf x_1\},\ldots,\mathcal X_{j}\cup \{\mathbf x_{j}\})$ is an orientated Reay system for any $j\in [0,s]$.

If $\mathcal B\subseteq \mathcal X_1\cup\ldots\cup\mathcal X_s$ and $\x\subseteq \C^\cup(\mathcal B)$ with $\x\in \mathcal X_j\cup\{\mathbf v_j\}$ for $j\in [1,s]$, then \be\label{faqueen}\x\in \mathcal X_j\mbox{ and there is some $\y\in \mathcal B$ with $\x\subseteq \y$},\ee which can be seen by the following short inductive proof on $s$ using (OR1) and (OR2). When $j=s$,  applying $\pi_{s-1}$ and using (OR2) yields the desired result, which covers the case $s=1$. Let $\mathcal B_s=\mathcal B\cap \mathcal X_s$. Observe that (OR2) ensures there is no nontrivial linear combination of elements from $B_s$ lying in $\mathcal E_{s-1}$, whence $\mathcal E_{s-1}\cap \C^\cup (\mathcal B_s)=\C\big(\bigcup_{\y\in \mathcal B_s}(\mathcal E_{s-1}\cap \y)\big)=\C\big(\bigcup_{\y\in \mathcal B_s}(\partial(\y)\cap \y)\big)$, with the latter equality in view of (OR1) and (OR2).
As a result, when $j<s$, we have $\x\subseteq \mathcal E_{s-1}\cap \C^\cup(\mathcal B)=\C^\cup(\mathcal B\setminus \mathcal B_s)+(\mathcal E_{s-1}\cap \C^\cup(\mathcal B_s))=\C^\cup(\mathcal B\setminus
\mathcal B_s)+\C\big(\bigcup_{\y\in \mathcal B_s}(\partial(\y)\cap \y)\big)$.
In view of (OR1), there is a subset $\mathcal X'\subseteq \mathcal X_1\cup\ldots\cup \mathcal X_{s-1}$ with $\C^\cup(\mathcal X')=\C\big(\bigcup_{\y\in \mathcal B_s}(\partial(\y)\cap \y)\big)$, namely $\mathcal X'=\bigcup_{\y\in \mathcal B_s}\mathcal B_\y$.
Moreover, if $\y\in \mathcal B_s$ and $\z\in \mathcal B_\y$, then $\z\subseteq \C^\cup(\mathcal B_\y)= \partial(\y)\cap \y\subseteq \y$, meaning every $\z\in \mathcal X'$ has some $\y\in \mathcal B$ with $\z\subseteq \y$.  Applying the induction hypothesis to $\mathcal B\setminus \mathcal B_s\cup\mathcal X'\subseteq \mathcal X_1\cup\ldots\cup \mathcal X_{s-1}$ now yields the desired result.

We define a partial order on the elements $\x,\,\y\in \mathcal X_1\cup\{\mathbf v_1\}\cup\ldots\cup\mathcal X_s\cup\{\mathbf v_s\}$ by declaring $\x\preceq \y$ when $\x\subseteq \y$. If $\x,\,\y\in \mathcal X_j\cup\{\mathbf v_j\}$, where $j\in [1,s]$, then $\x\preceq \y$ is only possible if $\x=\y$ (which can be seen by applying the map $\pi_{j-1}$ and using (OR2)). If $\x\in \mathcal X_{j'}\cup \{\mathbf v_{j'}\}$ and $\y\in \mathcal X_{j}\cup \{\mathbf v_{j}\}$ with  $\x\prec \y$, then $j'<j$ (which can be seen by applying the map $\pi_{j}$ to $\x$ and $\y$ to conclude $j'\leq j$ and then using the previous observation). In such case, we have $\x\subseteq\mathcal E_{j-1}\cap \y=\partial(\y)\cap \y$ (the equality follows in view of (OR1) and (OR2) as before), whence (OR1) and \eqref{faqueen} ensure that $\x\in \mathcal X_{j'}$. Thus each $\mathbf v_j$ is a maximal element.
If $\mathcal B\subseteq \mathcal X_1\cup\{\mathbf v_1\}\cup\ldots\cup\mathcal X_s\cup\{\mathbf v_s\}$, we let
$$\darrow \mathcal B=\{\x\in \mathcal X_1\cup\{\mathbf v_1\}\cup\ldots\cup\mathcal X_s\cup\{\mathbf v_s\}:\; \x\preceq \y \mbox{ for some $\y\in \mathcal B$}\}$$ denote the down-set generated by $\mathcal B$. Likewise, we let $\darrow B$ denote the set of representatives for $\darrow \mathcal B$, where $B\subseteq X_1\cup\{v_1\}\cup\ldots\cup X_s\cup\{v_s\}$ is the set of representatives for $\mathcal B$. Indeed, since (OR2) ensures there is a bijective correspondence between $\mathcal X_1\cup\{\mathbf v_1\}\cup\ldots\cup\mathcal X_s\cup\{\mathbf v_s\}$ and  $X_1\cup\{ v_1\}\cup\ldots\cup X_s\cup\{v_s\}$, the partial order defined above inherits to one on $X_1\cup\{ v_1\}\cup\ldots\cup X_s\cup\{v_s\}$. We let $\mathcal B^*\subseteq \mathcal B$ denote the subset of all maximal elements of $\mathcal B$, that is, all $\x\in\mathcal B$ for which there is no $\y\in \mathcal B$ with $\x\prec \y$. Clearly, $\darrow (\mathcal B^*)=\darrow \mathcal B$, \ $(\darrow \mathcal B)^*=\mathcal B^*$,  \be\label{rboow}\C^\cup(\mathcal B^*)=\C^\cup(\mathcal B)=\C^\cup(\darrow \mathcal B)\quad\und\quad \R^\cup(\mathcal B^*)=\R^\cup(\mathcal B)=\R^\cup(\darrow \mathcal B).\ee

A short argument now shows there  is a uniquely defined subset $\mathcal B_\x$ satisfying (OR1) with the additional property that $\mathcal B_\x^*=\mathcal B_\x$. Indeed, the existence of such a set follows in view of $\C^\cup(\mathcal B^*_\x)=\C^\cup(\mathcal B_\x)$. On the other hand, if $\mathcal C_\x\subseteq \mathcal X_1\cup\ldots\cup X_{j-1}$ is another set satisfying (OR1) with $\mathcal C_\x^*=\mathcal C_\x$, then we have $\C^\cup(\mathcal C_\x)=\partial(\mathbf x)\cap \x=\C^\cup(\mathcal B_x)$. Consequently, if $\y\in \mathcal B_\x$ is arbitrary, then $\y\subseteq \C^\cup(\mathcal B_\x)=\C^\cup(\mathcal C_\x)$, whence \eqref{faqueen} implies $\y\subseteq \z$ for some $\z\in \mathcal C_\x$. Thus $\mathcal B_\x\subseteq \darrow \mathcal C_\x$, implying $\darrow \mathcal B_\x\subseteq \darrow \mathcal C_\x$. Swapping the roles of $\mathcal B_\x$ and $\mathcal C_\x$ and repeating this argument shows $\darrow \mathcal C_\x\subseteq \darrow \mathcal B_\x$. As a result, we find that $\darrow \mathcal B_\x= \darrow \mathcal C_\x$, in turn implying $\mathcal B_\x=\mathcal B_\x^*=(\darrow \mathcal B_\x)^*=(\darrow \mathcal C_\x)^*=\mathcal C_\x^*=\mathcal C_\x$, establishing  the uniqueness of $\mathcal B_\x$.
We now henceforth
use $\partial(\{\x\}):=\mathcal B_\x$ to denote the unique set satisfying (OR1) with $\partial(\{\x\})^*=\partial(\{\x\})$, and let  $\partial(\{x\})$ denoting the set of representatives for $\partial(\{\x\})$.
Note, if $\x,\,\y\in \mathcal X_1\cup\{\mathbf v_1\}\cup\ldots\cup\mathcal X_s\cup\{\mathbf v_s\}$ with $\y\prec \x$ and $\x\in \mathcal X_j\cup \{\mathbf v_j\}$, then $\y\subseteq \mathcal E_{j-1}\cap \x=\partial(\x)\cap \x=\C^\cup (\partial(\{\x\}))$, whence $\y\in \darrow \partial(\{\x\})$ by \eqref{faqueen}. In consequence, since $\C^\cup (\darrow \partial(\{\x\}))=\C^\cup (\partial(\{\x\}))=\partial(\x)\cap \x\subseteq \x$, we find that  $$\darrow \partial(\{\x\})=\darrow \x\setminus \{\x\}.$$ Also, if $\partial(\x)=\{0\}$, then $\partial(\{\x\})=\emptyset$, which will be the case for any $\x\in \mathcal X_1\cup \{\mathbf v_1\}$.

In view of (OR1), it follows that, for any subset $\mathcal B\subseteq \mathcal X_1\cup\{\mathbf v_1\}\ldots\cup \mathcal X_j\cup \{\mathbf v_j\}$, where $j\in [1,s]$, there exists a subset $\partial(\mathcal B)\subseteq \mathcal X_1\cup\ldots\cup \mathcal X_{j-1}$ with
\be\label{leprechaumy}\R\la {\bigcup}_{\x\in \mathcal B}\partial(\x)\ra=\R^\cup\la \partial(\mathcal B)\ra\quad\und\quad\C\Big({\bigcup}_{\x\in \mathcal B}(\partial(\x)\cap \x)\Big)=\C^\cup(\partial(\mathcal B)).\ee For instance, we could take $\partial(\mathcal B)=(\bigcup_{\x\in \mathcal B}\mathcal B_{\x})^*$. Moreover, as argued in the previous paragraph, if we set $$\partial(\mathcal B):=\big(\bigcup_{\x\in \mathcal B}\mathcal B_{\x}\big)^*=\big(\bigcup_{\x\in \mathcal B}\partial(\{\x\})\big)^*,$$ then $\partial(\mathcal B)\subseteq \mathcal X_1\cup\ldots\cup \mathcal X_{j-1}$ will be the unique set satisfying \eqref{leprechaumy} with the additional property that $\partial(\mathcal B)^*=\partial(\mathcal B)$, which we henceforth assume is the case. We let $\partial(B)\subseteq X_1\cup\ldots\cup X_{j-1}$ denote the set of representatives for $\partial(\mathcal B)$. Note, from its definition, if $\mathcal A\subseteq \mathcal B$, then $\partial(\mathcal A)\subseteq \darrow \partial(\mathcal B)$.
Indeed, $\darrow \partial(\mathcal B)=\bigcup_{\x\in \mathcal B}\darrow \partial(\{\x\}),$ ensuring \be\label{delta-union}\darrow \partial(\mathcal A\cup \mathcal B)=\darrow \partial(\mathcal A)\cup \darrow \partial(\mathcal B)\quad\mbox{ for $\mathcal A,\,\mathcal B\subseteq\mathcal X_1\cup\{\mathbf v_1\}\cup\ldots\cup\mathcal X_s\cup \{\mathbf v_s\}$}.\ee

We extend the  partial order $\preceq$ to the subsets of $\mathcal X_1\cup\{\mathbf v_1\}\cup\ldots\cup \mathcal X_s\cup \{\mathbf v_s\}$ as follows. Given any subset $\mathcal B\subseteq \mathcal X_1\cup\{\mathbf v_1\}\cup\ldots\cup \mathcal X_s\cup \{\mathbf v_s\}$, we define the immediate predecessors to $\mathcal B$ to be the sets $\mathcal B'=\mathcal B\setminus \{\mathbf x\}\cup \partial(\{\x\})$ for $\mathbf x\in \mathcal B$.
 Since the tuple $(|\mathcal B\cap (\mathcal X_s\cup \{\mathbf v_s\})|,\ldots,|\mathcal B\cap (\mathcal X_1\cup \{\mathbf v_1\})|)$ associated to $\mathcal B$ strictly decreases in the lexicographic order under this operation, extending this relation transitively then defines the partial order $\preceq$. From its definition and a short inductive argument on $(|\mathcal B\cap ( \mathcal X_s\cup\{\mathbf v_s\})|,\ldots,|\mathcal B\cap (\mathcal X_1\cup\{\mathbf v_1\})|)$, we find that \be\label{lassie1}\mathcal A\preceq \mathcal B\quad\mbox{implies}\quad\mathcal A\subseteq \darrow \mathcal B,\ee thus ensuring that the partial order on the subsets of $\mathcal X_1\cup\{\mathbf v_1\}\cup\ldots\cup\mathcal X_s\cup\{\mathbf v_s\}$ is compatible with the pre-order induced from the partial order for the elements of  $\mathcal X_1\cup\{\mathbf v_1\}\cup\ldots\cup\mathcal X_s\cup\{\mathbf v_s\}$.
 As a result, \eqref{rboow} implies that, if  $\mathcal A\preceq \mathcal B$, then $\C^\cup(\mathcal A)\subseteq \C^\cup(\mathcal B)$.
Also, $$\mathcal A\subseteq \mathcal B\quad\mbox{implies}\quad\mathcal A\preceq \mathcal B,$$ as the following argument shows. Note, it suffices to show $\mathcal B\setminus\{\x\}\preceq \mathcal B$ for $\x\in \mathcal B$, as this can then be iterated. To see this,  we observe that  we may simply replace each $\x$ with the half-spaces from $\partial(\{\x\})$, and then replace each $\y\in \partial(\{\x\})\setminus \mathcal B$ with the half-spaces from $\partial(\{\y\})$, and so forth, until all such elements and their successors are replaced either using the empty set or an element already in $\mathcal B\setminus \{\x\}$.

Next, $$\mathcal A\subseteq \mathcal B\quad\mbox{implies}\quad(\mathcal B\setminus \mathcal A)\cup\partial(\mathcal A)\preceq \mathcal B,$$ which can be seen as follows. Sequentially replacing each $\x\in \mathcal A$ with $\partial(\{\x\})$, always choosing the next $\x\in \mathcal A$ in the sequence to be an element minimal among the remaining elements of $\mathcal A$ with respect to $\preceq$ (for instance, we could first take all $\x\in \mathcal A\cap (\mathcal X_1\cup\{\mathbf v_1\})$, then all $\x\in\mathcal A\cap (\mathcal X_2\cup\{\mathbf v_2\})$, and so forth), shows $(\mathcal B\setminus \mathcal A)\cup \bigcup_{\x\in \mathcal A}\partial(\{\x\})\preceq \mathcal B$. Then, since $(\mathcal B\setminus \mathcal A)\cup\partial(\mathcal A)=(\mathcal B\setminus \mathcal A)\cup \big(\bigcup_{\x\in \mathcal A}\partial(\{\x\})\big)^*\subseteq (\mathcal B\setminus \mathcal A)\cup \bigcup_{\x\in \mathcal A}\partial(\{\x\})\preceq\mathcal B$, the claimed result $(\mathcal B\setminus \mathcal A)\cup\partial(\mathcal A)\preceq \mathcal B$ follows. In particular, $\partial(\mathcal A)\preceq \mathcal B$ when $\mathcal A\subseteq \mathcal B$.

We remark that \be\label{lassie2}\{\x\in \mathcal A:\;\mathcal A\preceq \mathcal B\}=\darrow \mathcal B,\ee as can be seen by an inductive argument on $s$. Indeed, the inclusion $\{\x\in \mathcal A:\;\mathcal A\preceq \mathcal B\}\subseteq \darrow \mathcal B$ follows from \eqref{lassie1}, while $\x\in \darrow \mathcal B$ implies $\x\preceq \y$ for some $\y\in \mathcal B$. If $\x=\y\in \mathcal B\preceq \mathcal B$, the reverse inclusion holds. Otherwise, $\x\subseteq \mathcal E_{j-1}\cap \y=\partial(\y)\cap \y=\C^\cup(\partial(\{\y\}))$, where $\y\in \mathcal X_j\cup \{\mathbf v_j\}$. In view of  \eqref{faqueen}, we have $\x\preceq \z$ for some $\z\in \partial(\{\y\})\preceq \mathcal B$, i.e., $\x\in \darrow \partial(\{\y\})$, and now applying the induction hypothesis to $\partial(\{\y\})\subseteq \mathcal X_1\cup\ldots\cup\mathcal X_{s-1}$ yields the reverse inclusion.

Let $\partial^n(\mathcal B)={\underbrace{\partial(\partial(\ldots\partial}}_n(\mathcal B))\ldots))$ for $n\geq 0$, so $\partial ^0(\mathcal B)=\mathcal B$.  A similar inductive argument   on $s$ yields \be\label{newunion}\darrow \mathcal B=\bigcup_{n=0}^{s-1}\partial^n(\mathcal B).\ee
Note $\partial^s(\mathcal B)=\emptyset$ for any $\mathcal B\subseteq \mathcal X_1\cup \{\mathbf v_1\}\cup \ldots\cup \mathcal X_s\cup \{\mathbf v_s\}$ and that the case $s=1$ is clear since $\partial(\{\x\})=\emptyset$ for all $\x\in \mathcal B$ in this case.  The inclusion $\bigcup_{n=0}^{s-1}\partial^n(\mathcal B)\subseteq \bigcup_{n=0}^{s-1}\darrow\partial^n(\mathcal B)\subseteq \darrow \mathcal B$ follows in view of $\partial^{s-1}(\mathcal B)\preceq \partial^{s-2}(\mathcal B)\preceq\ldots\preceq \partial^0(\mathcal B)=\mathcal B$ and \eqref{lassie1}. On the other hand, if $\x\in \darrow \mathcal B$, then either $\x\in \mathcal B=\partial^0(\mathcal B)\subseteq \bigcup_{n=0}^{s-1}\partial^n(\mathcal B)$, or else
 $\x\subseteq \mathcal E_{j-1}\cap \y=\partial(\y)\cap \y=\C^\cup(\partial(\{\y\}))$, for some $\y\in \mathcal B$ with  $\y\in \mathcal X_j\cup \{\mathbf v_j\}$, in which case \eqref{faqueen} yields $\x\in \darrow \partial(\{\y\})$, and now applying the induction hypothesis to $\partial(\mathcal B)\subseteq \mathcal X_1\cup\ldots\cup\mathcal X_{s-1}$ yields $\x\in \darrow \partial(\{\y\})\subseteq \bigcup_{n=1}^{s-1}\darrow \partial^n(\mathcal B)=\darrow \partial(\mathcal B)=\bigcup_{n=1}^{s-1}\partial^n(\mathcal B)\subseteq \bigcup_{n=0}^{s-1}\partial^n(\mathcal B)$, establishing the reverse inclusion.


It may be helpful to view the half-spaces from $\mathcal X_1\cup \{\mathbf v_1\}\cup\ldots\cup \mathcal X_s\cup \{\mathbf v_s\}$ as vertices in a directed graph with each half-space $\x$ connected to the half-spaces from $\partial(\{\x\})$ by a directed edge.
Let $\mathcal B \subseteq \mathcal X_1\cup\{\mathbf v_1\}\cup\ldots\cup\mathcal X_s\cup\{\mathbf v_s\}$ and let $\mathcal Y$ denote the subset of vertices reachable from some $\x\in \mathcal B$, so $\y\in \mathcal Y$ means there is a sequence $\y_0,\ldots,\y_r$ with $\y_i\in \partial(\{\y_{i-1}\})$ for $i\geq 1$, $\y_0\in \mathcal B$ and $\y_r=\y$. Note \eqref{newunion} and \eqref{delta-union} imply $\darrow\partial(\darrow \mathcal A)=
\darrow \partial(\bigcup_{n=0}^{s-1}\partial^n(\mathcal A))=\bigcup_{n=1}^{s-1}\darrow \partial^n(\mathcal A)=\darrow \darrow \partial(\mathcal A)=\darrow \partial(\mathcal A)$ for any $\mathcal A\subseteq \mathcal X_1\cup \{\mathbf v_1\}\cup\ldots\cup \mathcal X_s\cup \{\mathbf v_s\}$.
This can be used in an inductive argument on $i=0,1,\ldots,r$ to show $\y_i\in\darrow \partial^i(\mathcal B)$ for $i\in [0,r]$. Indeed, $\y_0\in \mathcal B$ ensures $\y_0\in \darrow \mathcal B=\darrow \partial^0(\mathcal B)$, while $\y_i\in \partial(\{\y_{i-1}\})$ and $\y_{i-1}\in \darrow \partial^{i-1}(\mathcal B)$ then ensure
$\y_i\in \darrow \partial(\darrow \partial^{(i-1)}(\mathcal B))=\darrow \partial^i(\mathcal B)$, completing the induction.
In particular, $\y=\y_r\in \darrow\partial^r(\mathcal B)$, so that \eqref{newunion} yields $\mathcal Y\subseteq \darrow \mathcal B$. The reverse inclusion follows more directly from \eqref{newunion}, meaning  $\darrow \mathcal B=\mathcal Y$ consists of all vertices which can be reached via a directed path starting at some vertex from the subset  $\mathcal B$. The grading condition that $\partial(\mathcal B)\subseteq \mathcal X_1\cup \ldots\cup \mathcal X_{j-1}$ for $\mathcal B\subseteq \mathcal X_1\cup\{\mathbf v_1\}\cup \ldots\cup \mathcal X_{j}\cup\{\mathbf v_j\}$ ensures there are no directed cycles.

As already remarked, $\mathcal A\preceq \mathcal B$ implies $\mathcal A\subseteq \darrow\mathcal B$. The partial converse to this is \be\label{lassie3}\mathcal A\subseteq \darrow \mathcal B\quad\mbox{implies}\quad \mathcal A^*\preceq \mathcal B.\ee To see this, we modify the argument for showing $\mathcal A\subseteq \mathcal B$ implies $\mathcal A\preceq \mathcal B$. Order the elements of $\mathcal B\setminus \mathcal A$. Beginning with the first $\x\in\mathcal B\setminus \mathcal A$,  replace $\x$  with the half-spaces from $\partial(\{\x\})$, and then replace each $\y\in \partial(\{\x\})\setminus  \mathcal A$ with the half-spaces from $\partial(\{\y\})$, and so forth, until all such elements and their successors are replaced either using the empty set or an element already in $\mathcal A$. Let $\mathcal B'\subseteq \mathcal B\cup \mathcal A$ be the resulting set. Take the next $\x'\in \mathcal B'\cap (\mathcal B\setminus \mathcal A)=\mathcal B'\setminus \mathcal A$ and repeat the procedure. Continue until all elements from $\mathcal B\setminus \mathcal A$ have been exhausted and let $\mathcal C\subseteq \mathcal A$ be the resulting set. By construction, $\mathcal C\preceq \mathcal B$. Since $\mathcal A^*$ is the set of maximal elements in $\mathcal A$, it follows from $\mathcal A\subseteq \darrow \mathcal B$ and  \eqref{newunion} that $\mathcal A^*\subseteq \mathcal C$, whence $\mathcal A^*\preceq \mathcal C\preceq \mathcal B$, as desired.

If $\mathcal A\subseteq \darrow \mathcal B$, then $\partial(\mathcal A)\subseteq \darrow \partial(\mathcal B)$ can be seen as follows. Let $\x\in \partial(\mathcal A)$. Then $\x\in \partial(\{\x'\})$ for some $\x'\in \mathcal A\subseteq \darrow \mathcal B$, in which case there is some $\y\in \mathcal B$ with $\x\preceq \x'\preceq \y$. If $\x'=\y$, then $\x\in \partial(\{\x'\})\subseteq \darrow \partial(\mathcal B)$. Otherwise, $\x\subseteq \mathcal E_{j-1}\cap \y=\partial(\y)\cap
\y=\C^\cup(\partial(\{\y\}))$, where $\y\in \mathcal X_j\cup\{\mathbf v_j\}$, in which case \eqref{faqueen} implies $\x\preceq \y'$ for some $\y'\in \partial(\{\y\})$, i.e., $\x\in \darrow \partial(\{\y\})\subseteq \darrow \partial(\mathcal B)$, as desired. Combining this with \eqref{lassie3}, we conclude that $$\mathcal A\subseteq \darrow \mathcal B\quad\mbox{implies}\quad \partial(\mathcal A)\preceq \partial(\mathcal B).$$

\begin{definition}Let $\mathcal R=(\mathcal X_1\cup\{\mathbf v_1\},\ldots,\mathcal X_s\{\mathbf v_s\})$ be an oriented Reay system for the subspace $\mathcal E_s\subseteq \R^d$ and let $\mathcal B\subseteq \mathcal X_1\cup\{\mathbf v_1\}\cup\ldots\cup\mathcal X_s\cup \{\mathbf v_s\}$.
We say that $\mathcal B$ is a \textbf{support set} for $\mathcal R$ if $$\mathcal B^*=\mathcal B\quad\und\quad\mathcal X_i\cup \{\mathbf v_i\}\nsubseteq \darrow \mathcal B\quad\mbox{ for all $i\in [1,s]$}.$$ We say that $\mathcal B$ is a \textbf{virtual independent} set if
$$\mathcal B^*=\mathcal B\quad\und\quad\mathcal \darrow B\;\mbox{ is linearly independent}.$$
\end{definition}

   We remark that Proposition \ref{prop-orReay-BasicProps}.3 ensures that the definition of a virtual independent set does not depend on the choice of representative set $\darrow B$.

\begin{proposition}\label{prop-orReay-BasicProps}
Let $\mathcal R=(\mathcal X_1\cup\{\mathbf v_1\},\ldots,\mathcal X_s\cup\{\mathbf v_s\})$ be a an orientated Reay system for the subspace $\mathcal E_s\subseteq \R^d$ and let $\mathcal B\subseteq \mathcal X_1\cup\{\mathbf v_1\}\cup\ldots\cup \mathcal X_{s}\cup \{\mathbf v_s\}$. For $j\in [1,s]$, let $\mathcal E_{j-1}=\R\la \mathcal X_1\cup\ldots\cup \mathcal X_{j-1}\ra$, let $\mathcal E'_{j-1}=\mathcal E_{j-1}+\R^\cup \la \partial(\mathcal B)\ra$, let $\pi'_{j-1}:\R^d\rightarrow (\mathcal E'_{j-1})^\bot$ be the orthogonal projection, and  let $\mathcal B_j=\mathcal B\cap (\mathcal X_j\cup \{\mathbf v_j\})$.

\begin{itemize}
\item[1.]   $\R\la \darrow B\ra=\R^\cup\la\mathcal B\ra$.
\item[2.] $\C^\cup(\mathcal B)=\z_1+\ldots+\z_\ell$ is a convex cone containing $0$, and $\overline{\C^\cup(\mathcal B)}=\overline \z_1+\ldots+\overline \z_\ell$ is a polyhedral cone, where $\z_1,\ldots,\z_\ell\in \mathcal B$ are the distinct half-spaces in $\mathcal B$.
\item[3.] $\mathcal B$ is a virtual independent set if and only if $\pi(B)$ is a linearly independent set of size $|\mathcal B|$, where $\pi:\R^d\rightarrow \R^\cup \la\partial(\mathcal B)\ra^\bot$ is the orthogonal projection.
 \item[4.] If $\mathcal B$ is virtual independent, then $\overline{\C^\cup(\mathcal B)}=\overline \z_1+\ldots+\overline \z_\ell$ is a blunted simplicial cone with lineality space $\R^\cup \la \partial(\mathcal B)\ra$, \ $\C^\cup(\mathcal B)$ has trivial lineality space, and $\C^\cup(\mathcal B)^\circ=\z_1^\circ+\ldots+\z_\ell^\circ$, where $\z_1,\ldots,\z_\ell\in \mathcal B$ are the distinct half-spaces in $\mathcal B$.
\item[5.] If $\mathcal B$ is a support set for $\mathcal R$, then $\mathcal B$ is a virtual independent set.
\item[6.] If $\mathcal B$ is a support set for $\mathcal R$,  then $\pi'_{j-1}\big(\bigcup_{i=j}^sB_i\big)$ is a linearly independent subset of $\Sum{i=j}{s}|\mathcal B_i|\geq 0$ distinct elements, for any $j\in [1,s]$.
\item[7.] If $\mathcal B$ is virtual independent   and   $\mathcal A\subseteq \mathcal B$, then  $\mathcal C=((\mathcal B\setminus \mathcal A)\cup\partial(\mathcal A))^*$ is virtual independent  with $\mathcal B\setminus \mathcal A\subseteq \mathcal C$. Moreover, if $\mathcal B$ is a support set, then so is $\mathcal C$.
\item[8.] If $\mathcal B\subseteq \mathcal X_1\cup\ldots\cup \mathcal X_s$, then $\mathcal B^*$ is a support set. In particular, $\partial(\mathcal B)$ is always a support set for any $\mathcal B\subseteq \mathcal X_1\cup\{\mathbf v_1\}\cup\ldots\cup\mathcal X_s\cup\{\mathbf v_s\}$.
\item[9.] If $\mathcal B\subseteq \mathcal X_1\cup\ldots\cup \mathcal X_s$, \ $\x\in\mathcal X_1\cup\ldots\cup \mathcal X_s$ and $\x\subseteq \R^\cup\la \mathcal B\ra$, then $\x\in \darrow \mathcal B$.
    \end{itemize}
\end{proposition}

\begin{proof}

1.
We proceed by induction on the depth $s$.
In view of \eqref{newunion}, we have \be\label{larkst}\darrow \mathcal B=
\mathcal B_s\cup\darrow \partial(\mathcal B)\cup\darrow (\mathcal B\setminus \mathcal B_s)\quad\und\quad  \darrow  B=
B_s\cup\darrow \partial(B)\cup \darrow (B\setminus B_s).\ee
By induction hypothesis,  $\R\la \darrow \partial( B)\cup \darrow (B\setminus B_s)\ra=\R^\cup\la \partial(\mathcal B)\cup (\mathcal B\setminus \mathcal B_{s})\ra$ (note $\darrow\big(\partial(\mathcal B)\cup (\mathcal B\setminus \mathcal B_s)\big)=\emptyset$ for $s=1$).
 But now \eqref{larkst} implies that
  $$\R\la \darrow B\ra=\R\la B_s\ra+\R\la \darrow \partial(B)\cup \darrow (B\setminus  B_s)\ra=\R\la B_s\ra+\R^\cup\la\partial(\mathcal B)\ra
  +\R^\cup\la \mathcal B\setminus  \mathcal B_{s}\ra=\R^\cup\la  \mathcal B\ra,$$ as desired.

2. Suppose $\mathcal B\subseteq \mathcal X_1\cup \{\mathbf v_1\}\cup \ldots\mathcal X_j\cup\{\mathbf v_j\}$, where $j\in [0,s]$, and let $\z_1,\ldots,\z_\ell\in \mathcal B$ be the distinct half-spaces in $\mathcal B$.
For $j=0$, we have $\mathcal B=\emptyset$ and  $\C^\cup(\mathcal B)=\overline{\C^\cup(\mathcal B)}=\{0\}$, which is the sum of an empty number of half-spaces. For $j=1$, we have  $\partial(\{\z_i\})=\emptyset$, $\partial(\z_i)=\{0\}$ and  $\z_i=\overline \z_i=\R_+z_i$ for all $i$ (by (OR2)), in which case  $\C^\cup (\mathcal B)=\C(z_1,\ldots,z_r)$ with the result clear. Thus we assume $j\geq 2$ and proceed by induction on $j$. By induction hypothesis applied to $\partial(\{\z_i\})$, we have $0\in \C^\cup (\partial(\{\z_i\}))=\z_i\cap \partial(\z_i)$, ensuring each $\z_i$ is a convex cone with $0\in (\z_i\cap \partial(\z_i))\subseteq \z_i$, for $i\in [1,\ell]$, in turn implying
 $0\in \C^\cup (\mathcal B)=\z_1+\ldots+\z_\ell$ with $\C^\cup(\mathcal B)$ a convex cone.

 Note $\overline \z_i=\R_+z_i+\partial(\z_i)=\C(Y_i\cup \{z_i\})$, where $Y_i$ is any positive basis for the subspace $\partial(\z_i)$.
Thus $\overline \z_1+\ldots+\overline \z_\ell=\C\big(\bigcup_{i=1}^\ell(Y_i\cup \{z_i\}\big)\big)$, which is polyhedral cone, and thus closed. Hence $\overline{\C^\cup(\mathcal B)}=\overline{\z_1+\ldots+\z_\ell}\subseteq \overline{ \overline \z_1+\ldots+\overline \z_\ell}\subseteq \overline \z_1+\ldots+\overline \z_\ell$.
 On the other hand, since $0\in \z_i$ for all $i$ with  $\C^\cup (\mathcal B)$, and thus also $\overline{\C^\cup(\mathcal B)}$, a convex cone, it follows that
$\overline \z_1+\ldots+\overline \z_\ell\subseteq \overline{\z_1+\ldots+\z_\ell}=\overline{\C^\cup(\mathcal B)}$, establishing the reverse inclusion, showing $\overline{\C^\cup(\mathcal B)}=\overline\z_1+\ldots+\overline\z_\ell=\C\big(\bigcup_{i=1}^\ell(Y_i\cup \{z_i\}\big)$ is a polyhedral cone.

3. Suppose $\mathcal B^*=\mathcal B$ and $\darrow B$ is linearly independent. Since $\mathcal B^*=\mathcal B$, it follows that $\mathcal B$ is disjoint from $\darrow \partial(\mathcal B)$. By Item 1, $\R\la \darrow \partial(B)\ra=\R^\cup \la \partial(\mathcal B)\ra$. Thus $\ker \pi$ is generated by the linearly independent subset $\darrow \partial(B)\subseteq \darrow B$ (by Proposition \ref{prop-reay-basis-properties}.1), ensuring  $\pi(B)$ is a linearly independent set of size $|\mathcal B\setminus \darrow \partial(\mathcal B)|=|\mathcal B|$, with the equality since $\mathcal B$ is disjoint from $\darrow \partial(\mathcal B)$.
Next instead suppose  $\pi(B)$ is linearly independent of size $|\mathcal B|$. Since $\pi(B)$ is linearly independent,  no half-space in $\mathcal B$ is contained in $\R^\cup\la \partial(\mathcal B)\ra$. If $\x\prec \y$ with $\x,\,\y\in \mathcal B$, then we must have $\x\in \darrow \partial(\{\y\})$, contradicting that $\x\nsubseteq \R^\cup\la\partial(\mathcal B)\ra=\R^\cup\la\darrow \partial( \mathcal B)\ra$. Therefore  we instead conclude that $\mathcal B^*=\mathcal B$. Since $\darrow\partial(\mathcal B)\subseteq \mathcal X_1\cup\ldots\cup\mathcal X_{s-1}$, Proposition \ref{prop-reay-basis-properties}.1 ensures  that $\darrow\partial(B)$ is always a linearly independent set.
Consequently, since $\pi(B)$ is a linearly independent set of size $|\mathcal B|$ with $\ker \pi=\R\la \darrow \partial(B)\ra$ (by Item 1), it follows that $\darrow B=B\cup \darrow \partial(B)$ is linearly independent, as desired.

4.  In view of Items 2 and 3, $\overline{\C^\cup(\mathcal B)}=\overline \z_1+\ldots+\overline \z_\ell$ is a blunted simplicial cone with lineality space $\R^\cup \la \partial(\mathcal B)\ra$, and  $\C^\cup(\mathcal B)^\circ=\z_1^\circ+\ldots+\z_\ell^\circ$. In particular, if $\partial(\mathcal B)=\emptyset$, then $\C^\cup(\mathcal B)$ has trivial lineality space. Assuming Items 5 and 8 have been established, $\partial(\mathcal B)$ will be virtual independent. We can then use an inductive argument on $|\darrow \mathcal B|$ to  show $\C^\cup(\mathcal B)$ has trivial lineality space for a general virtual independent set $\mathcal B$. Indeed, in view of the established portions of Item 4, the lineality space of $\C^\cup(\mathcal B)$ must be contained in $(\z_1\cap \partial(\z_1))+\ldots+(\z_\ell\cap \partial(\z_\ell))=\C^\cup(\partial(\{\z_1\}))+\ldots+\C^\cup(\partial(\{\z_\ell\}))
=\C^\cup(\partial(\mathcal B))$, and so applying the induction hypothesis to $\partial(\mathcal B)$, completes the proof.

5. This follows from Proposition \ref{prop-reay-basis-properties}.1.

6. By definition, $\pi'_{j-1}(\mathbf x)$ is either zero or a ray for any $\mathbf x\in \mathcal B$ and $j\in [1,s]$, so $\pi_{j-1}(x)$ is a representative for the half-space $\pi_{j-1}(\x)$ whenever $\pi_{j-1}(\x)\neq \{0\}$.
We proceed by induction on $j=s,s-1,\ldots,1$.
If Item 6 fails, then there must be a nontrivial linear combination of the elements of $\bigcup_{i=j}^sB_i$ equal to an element of $\mathcal E'_{j-1}$, say \be\label{stay}\Summ{x
\in B_{j+1}\cup\ldots\cup B_s}\alpha_{x}x+\Summ{y\in B_j}\beta_yy\in \mathcal E'_{j-1},\ee where the $\alpha_x,\,\beta_y\in \R$ are not all zero. But then $\Summ{x
\in B_{j+1}\cup\ldots\cup B_s}\alpha_{x}x\in \mathcal E'_{j-1}+\mathcal E_j=\mathcal E'_j=\ker \pi'_{j}$. By induction hypothesis, $\pi'_j(\bigcup_{i=j+1}^sB_i)$ is a linearly independent set of $\Sum{i=j+1}{s}|B_i|\geq 0$ elements, meaning this is only possible if $\alpha_x=0$ for all $x\in B_{j+1}\cup\ldots\cup B_s$ (note this is trivially true when $j=s$). Hence \be\label{stardreamt}\Summ{y\in B_j}\beta_yy\in \mathcal E'_{j-1}=\mathcal E_{j-1}+\R^\cup \la \partial(\mathcal B)\ra.\ee In view of Item 1, we have $\R^\cup \la \partial(\mathcal B)\ra=\R\la \darrow \partial(B)\ra$.
Since $\mathcal B^*=\mathcal B$ for the support set $\mathcal B$, it follows that  $\darrow \partial(B)$ and $B$ are disjoint subsets of $\darrow B$, which is linearly independent by Item 5 as $\mathcal B$ is a support set. Also, all elements from  $\darrow \partial(B)\cap (X_1\cup \ldots\cup X_{j-1})$ are contained in $\mathcal E_{j-1}$.  Thus \eqref{stardreamt} implies that there is a nontrivial linear combination \be\label{weent}\Summ{y\in B_j}\beta_yy+\Summ{z\in  Z}\gamma_zz\in \mathcal E_{j-1},\ee for some $\gamma_z\in \R$, where $Z:= \darrow \partial(B)\cap (X_j\cup \ldots\cup X_s)$. Since $B_j\cup Z\subseteq \darrow B$ is a disjoint union with  $\mathcal B$ a support set, Propositions \ref{prop-reay-basis-properties}.3 and \ref{prop-reay-basis-properties}.1 imply that $\pi_{j-1}(B_j\cup Z)$ is a linearly independent set of size $|B_j\cup Z|=|B_j|+|Z|$. Thus, applying $\pi_{j-1}$ to \eqref{weent}, it follows that $\beta_y=0$ for all $y\in B_j$, and $\gamma_z=0$ for all $z\in Z$, contradicting that the linear combination was nontrivial, which completes the proof of Item 6.

7. By definition of $\mathcal C$, we have $\mathcal C^*=\mathcal C$. Since $\mathcal A\subseteq \mathcal B$ (implying $\partial(\mathcal A)\preceq \mathcal B$ and then $\partial(\mathcal A)\subseteq \darrow \mathcal B$), we have $\darrow \mathcal C=\darrow (\mathcal B\setminus \mathcal A)\cup \darrow \partial(\mathcal A)\subseteq \darrow \mathcal B$. If $\mathcal B$ is a support set, then  $\mathcal X_j\cup \{\mathbf v_j\}\nsubseteq \darrow \mathcal B$ for every $j\in [1,s]$. As a result, since $\darrow \mathcal C\subseteq \darrow \mathcal B$, we also have $\mathcal X_j\cup \{\mathbf v_j\}\nsubseteq \darrow \mathcal C$ for every $j\in [1,s]$, which shows that $\mathcal C$ is a support set. If $\mathcal B$ is virtual independent, then $\darrow B$ is linearly independent, whence $\darrow C\subseteq \darrow  B$ is also linearly independent, showing that $\mathcal C$ is virtual independent.

Suppose $\mathcal B\setminus \mathcal A\nsubseteq \mathcal C$.  Then, in view of the definition of $\mathcal C$, there must be some $\x\in \mathcal B\setminus \mathcal A$ and $\y\in (\mathcal B\setminus \mathcal A)\cup \partial(\mathcal A)$ with $\x\prec \y$.
If $\y\in \mathcal B\setminus \mathcal A$, then $\x\prec \y$ contradicts that $\mathcal B^*=\mathcal B$ for the virtual independent set $\mathcal B$. Therefore we must have $\y\in\partial(\mathcal A)$, implying that $\x\prec \y\prec \z$ for some $\z\in \mathcal A\subseteq \mathcal B$, which again contradicts  that $\mathcal B^*=\mathcal B$ for the virtual independent  set $\mathcal B$. So we instead conclude that $\mathcal B\setminus \mathcal A\subseteq \mathcal C$, completing the proof of Item 7.

8.  Since $\mathcal B\subseteq \mathcal X_1\cup\ldots\cup \mathcal X_s$ implies $\darrow \mathcal B\subseteq \mathcal X_1\cup\ldots\cup \mathcal X_s$, we have $\mathcal X_j\cup \{\mathbf v_j\}\nsubseteq \darrow \mathcal B=\darrow(\mathcal B^*)$ for all $j\in [1,s]$, while $\mathcal B^{**}=\mathcal B^*$, which shows that $\mathcal B^*$ is a support set.

9. In view of Item 1, we have $\R^\cup\la \mathcal B\ra=\R\la \darrow B\ra$. Thus $x\in \x\subseteq \R^\cup\la \mathcal B\ra=\R\la \darrow B\ra$. Since $\mathcal B\subseteq \mathcal X_1\cup\ldots\cup \mathcal X_s$, we also have $\darrow B\subseteq X_1\cup\ldots\cup X_s$, while $x\in X_1\cup\ldots\cup X_s$ since $\x\in \mathcal X_1\cup\ldots\cup \mathcal X_s$. Proposition \ref{prop-reay-basis-properties}.1 implies that $X_1\cup\ldots\cup X_s$ is a linearly independent set, whence $x\in \R\la \darrow B\ra$ with $x\in X_1\cup\ldots\cup X_s$ and $\darrow B\subseteq X_1\cup\ldots\cup X_s$ is only possible if $x\in \darrow B$, implying $\x\in \darrow \mathcal B$.
\end{proof}

Next, we show that we have a well-behaved quotient oriented Reay system defined modulo $\mathcal R^\cup \la \mathcal B\ra$, for any $\mathcal B\subseteq \mathcal X_1\cup\ldots\cup \mathcal X_s$. We will use the notation $\mathcal C^\pi=\pi(\mathcal C)\setminus \{\{0\}\}$,
$\pi(\mathcal R)$ and $\pi^{-1}(\mathcal D)$ defined in Proposition \ref{prop-orReay-modulo} for the remainder of this work.

\begin{proposition}\label{prop-orReay-modulo}
Let $\mathcal R=(\mathcal X_1\cup\{\mathbf v_1\},\ldots,\mathcal X_s\cup\{\mathbf v_s\})$ be an orientated Reay system for the subspace $\mathcal E_s\subseteq \R^d$, let $\mathcal B\subseteq \mathcal X_1\cup\ldots\cup \mathcal X_{s}$ be a subset, let  $\mathcal E=\R^\cup\la  \mathcal B\ra$, and let $\pi:\R^d\rightarrow \mathcal E^\bot$ be the orthogonal projection. For a subset $\mathcal C\subseteq \mathcal X_1\cup\{\mathbf v_1\}\cup \ldots\cup \mathcal X_s\cup \{\mathbf v_s\}$, set $$\mathcal C^\pi:=\pi(\mathcal C)\setminus \{\{0\}\}.$$ For $j\in [1,s]$, let $\mathcal X'_j\subseteq \mathcal X_j$ be all those $\x\in \mathcal X_j$ with $\pi(\x)\neq \{0\}$, and let  $J\subseteq [1,s]$ be all those indices $j$ with $\mathcal X_j^\pi\neq \emptyset$, say $J=\{j_1,\ldots,j_t\}$ with $1\leq j_1<\ldots<j_t\leq s$.
For a subset $\mathcal D\subseteq \bigcup_{i\in J}(\mathcal X_i^\pi\cup \{\pi(\mathbf v_i)\})$, set $$\pi^{-1}(\mathcal D):=\Big\{\x\in \mathcal X_1\cup\ldots\cup\mathcal X_s\cup \{\mathbf v_i:\;i\in J\}:\; \pi(\x)\in \mathcal D\Big\}\subseteq \bigcup_{i\in J}(\mathcal X'_i\cup \{\mathbf v_i\}).$$
  \begin{itemize}
%
\item[1.]   $$\pi(\mathcal R):=\Big(\mathcal X^\pi_j\cup \{\pi(\mathbf v_j)\}\Big)_{j\in J}=\big(\mathcal X^\pi_{j_1}\cup\{\pi(\mathbf v_{j_1})\},\ldots,\mathcal X^\pi_{j_{t}}\cup\{\pi(\mathbf v_{j_{t}})\}\big)$$ is an oriented Reay system for $\pi(\mathcal E_s)$ with $\pi$ injective on $\bigcup_{i\in J}(\mathcal X'_i\cup \{\mathbf v_i\})$ and $\pi(\mathbf v_i)\neq \{0\}$ for all $i\in J$.
\item[2.]
If $\mathcal C\subseteq \mathcal X_1\cup\ldots\cup \mathcal X_s\cup \{\mathbf v_j:\; j\in J\}$, then $$(\mathcal C^\pi)^*=(\mathcal C^*)^\pi, \quad \darrow (\mathcal C^\pi)=(\darrow \mathcal C)^\pi\quad\und\quad \partial(\mathcal C^\pi)=\partial(\mathcal C)^\pi.$$

\item[3.] \begin{itemize}
\item[(a)] If $\mathcal C_1,\,\mathcal C_2\subseteq \mathcal X_1\cup\ldots\cup \mathcal X_s\cup \{\mathbf v_j:\; j\in J\}$ with $\mathcal C_1\prec \mathcal C_2$, then $\mathcal C_1^\pi\preceq \mathcal C_2^\pi$. Moreover, if  in a sequence of replacing  elements $\x$ by $\partial(\{\x\})$ showing that $\mathcal C_1\prec \mathcal C_2$ there is some $\x$ with $\pi(\x)\neq \{0\}$, then $\mathcal C_1^\pi\prec \mathcal C_2^\pi$.
    \item[(b)] If $\mathcal D_1,\,\mathcal D_2\subseteq \bigcup_{j\in J}(\mathcal X_j^\pi\cup \{\pi(\mathbf v_j)\})$ with $\mathcal D_1\prec \mathcal D_2$, then $\pi^{-1}(\mathcal D_1)\prec \pi^{-1}(\mathcal D_2)$.
    \end{itemize}

\item[4.]  If $\mathcal D\subseteq \bigcup_{i\in J}(\mathcal X_i^\pi\cup \{\pi(\mathbf v_i)\})$, then $$\pi\big((\pi^{-1}(\mathcal D))^*\big)=\mathcal D^*,\quad (\mathcal \pi^{-1}(\mathcal D))^*\subseteq (\pi^{-1}(\mathcal D)\cup \mathcal B)^*\quad\und\quad \big((\pi^{-1}(\mathcal D)\cup \mathcal B)^*\big)^\pi=\mathcal D^*.$$ In particular, $\mathcal D^*=\mathcal D$ if and only if $(\pi^{-1}(\mathcal D))^*=\pi^{-1}(\mathcal D)$.
\item[5.]
\begin{itemize}
 \item[(a)]If $\mathcal D\subseteq \bigcup_{i\in J}(\mathcal X_i^\pi\cup \{\pi(\mathbf v_i)\})$ is a virtual independent set, then $(\pi^{-1}(\mathcal D)\cup \mathcal B)^*$ and $\pi^{-1}(\mathcal D)\subseteq (\pi^{-1}(\mathcal D)\cup \mathcal B)^*$ are virtual independent sets.
     \item[(b)] If $\mathcal D\subseteq \bigcup_{i\in J}(\mathcal X_i^\pi\cup \{\pi(\mathbf v_i)\})$ is a support set for $\pi(\mathcal R)$, then  $(\pi^{-1}(\mathcal D)\cup \mathcal B)^*$ and $\pi^{-1}(\mathcal D)$ are support sets for $\mathcal R$.
        \end{itemize}
        \noindent In either case, if  $\mathcal B^*=\mathcal B$, then  $(\pi^{-1}(\mathcal D)\cup \mathcal B)^*= \pi^{-1}(\mathcal D)\cup \big(\mathcal B\setminus \darrow \partial\big(\pi^{-1}(\mathcal D)\big)\big)$.
\item[6.] \begin{itemize}\item[(a)] If $\mathcal C\subseteq \mathcal X_1\cup\ldots\cup \mathcal X_s\cup \{\mathbf v_j:\; j\in J\}$ is a virtual independent set with $\mathcal B\subseteq \darrow \mathcal C$, then $\mathcal C^\pi$ is a virtual independent set.

   \item[(b)] If $\mathcal C\subseteq \mathcal X_1\cup \{\mathbf v_1\}\cup \ldots\cup\mathcal X_s\cup\{\mathbf v_s\}$ is a support set with $\mathcal B\subseteq \darrow \mathcal C$, then $\mathcal C^\pi$ is a support set with $\mathcal C\subseteq \mathcal X_1\cup\ldots\cup \mathcal X_s\cup \{\mathbf v_j:\; j\in J\}$.

    \item[(c)] If $\mathcal C\subseteq \mathcal X_1\cup \ldots\cup\mathcal X_s$ is a support set, then $\mathcal C^\pi$ is a support set.\end{itemize}
    \end{itemize}

\end{proposition}

\begin{proof}
For $j\in [1,s+1]$, let $\mathcal E_{j-1}=\R^\cup \la \mathcal X_1\cup\ldots\cup X_{j-1}\ra$ and  let $\pi_{j-1}:\R^d\rightarrow \mathcal E_{j-1}^\bot$ and   $\varpi_{j-1}:\R^d\rightarrow (\mathcal E_{j-1}+\mathcal E)^\bot$ be the orthogonal projections.
Note $\varpi_0=\pi$. Since $\ker \pi_{j-1}\leq \ker \varpi_{j-1}$, we have  $\varpi_{j-1}=\varpi_{j-1}\pi_{j-1}$ for all $j\in [1,s]$.
In view of Proposition \ref{prop-orReay-BasicProps}.1, $$\R\la \darrow B\ra=\R^\cup\la \mathcal B\ra=\mathcal E$$ with $\darrow B\subseteq X_1\cup\ldots\cup X_s$ (since $\mathcal B\subseteq \mathcal X_1\cup\ldots\cup \mathcal X_s$).
Thus $\mathcal E'_{j-1}=\mathcal E_{j-1}+\R\la \darrow B\ra$, while $X_j\subseteq \mathcal E_j$ and $\mathcal E_{j-1}\subseteq \mathcal E_j$, implying $\mathcal E_{j-1}+\R\la \darrow B\cap X_j\ra\subseteq \mathcal E'_{j-1}\cap \mathcal E_j$. Any element  $x\in\mathcal E'_{j-1}\cap \mathcal E_j$ has $x=y+z$ with $y\in \mathcal E_{j-1}$, $z$  a linear combination of elements from $\darrow B$ and $x=y+z\in \mathcal E_j$. Then $\pi_j(z)=\pi_j(x)=0$. Applying $\pi_j$ to the linear combination representing $z$ and using  Propositions \ref{prop-reay-basis-properties}.1 and \ref{prop-reay-basis-properties}.3 shows that $x=y+z\in \mathcal E_{j-1}+\R\la \darrow B\cap X_j\ra$. Hence
\be\mathcal E'_{j-1}\cap \mathcal E_j=\mathcal E_{j-1}+\R\la \darrow B\cap X_j\ra\label{dunntee}.\ee
Since $\varpi_{j-1}=\varpi_{j-1}\pi_{j-1}$, we have $\ker \varpi_{j-1}\cap \pi_{j-1}(\mathcal E_j)=\pi_{j-1}(\ker \varpi_{j-1})\cap \pi_{j-1}(\mathcal E_j)$, meaning \be\label{dunntee2}\ker \varpi_{j-1}\cap \pi_{j-1}(\mathcal E_j)=\pi_{j-1}(\mathcal E'_{j-1})\cap \pi_{j-1}(\mathcal E_j).\ee
Any element $x\in \pi_{j-1}(\mathcal E'_{j-1})\cap \pi_{j-1}(\mathcal E_j)$ is a linear combination of terms $\pi_{j-1}(b)\in \pi_{j-1}(\darrow B)$  (with each $b\in \darrow B$)  that lies in $\pi_{j-1}(\mathcal E_j)$. Applying $\pi_j$ to this linear combination and using Propositions \ref{prop-reay-basis-properties}.1 and \ref{prop-reay-basis-properties}.3, we find that only those $b\in \darrow B\cap \mathcal E_j$ have non-zero coefficients, meaning $x\in \pi_{j-1}(\mathcal E'_{j-1}\cap \mathcal E_j)$. Thus $\pi_{j-1}(\mathcal E'_{j-1})\cap \pi_{j-1}(\mathcal E_j)=\pi_{j-1}(\mathcal E'_{j-1}\cap \mathcal E_j)$ (the reverse inclusion is trivial), which combined with \eqref{dunntee} and \eqref{dunntee2} yields
$$\ker \varpi_{j-1}\cap \pi_{j-1}(\mathcal E_j)=\R\la \pi_{j-1}(\darrow B\cap X_j)\ra.$$  Thus the kernel of $\varpi_{j-1}$ restricted to $\pi_{j-1}(\mathcal E_j)$ is generated by a subset of the linearly independent set $\pi_{j-1}(X_j)$.
Consequently, since $\pi_{j-1}(X_j)\cup \{\pi_{j-1}(v_j)\}$ is a minimal positive basis for $\pi_{j-1}(\mathcal E_j)$ of size $|X_j|+1$, it follows that $(X_j\cup\{v_j\})^{\varpi_{j-1}}$ is either empty or  a minimal positive basis for $\varpi_{j-1}(\mathcal E_j)$ of size \be\label{size1}|(X_j\cup\{v_j\})^{\varpi_{j-1}}|=|(X_j\cup \{v_j\})\setminus (\mathcal E_{j-1}+\mathcal E)|=|X_j|+1-|\darrow B\cap X_j|.\ee

Let $j\in J$.
 If $\varpi_{j-1}(x)=0$ for some $x\in X_j\cup \{v_j\}$, so  $x\in \mathcal E'_{j-1}\cap \mathcal E_j=\mathcal E_{j-1}+\R\la \darrow B\cap X_j\ra$ (the equality follows from \eqref{dunntee}), then  $\pi_{j-1}(x)\in \R\la \pi_{j-1}(\darrow B\cap X_j)\ra$ with $\pi_{j-1}(\darrow B\cap X_j)\subseteq \pi_{j-1}(X_j)$. Thus,
since $\pi_{j-1}(X_j)\cup \{\pi_{j-1}(v_j)\}$ is a minimal positive basis of size $|X_j|+1$ and $x\in X_j\cup \{v_j\}$, we must either  have
$\darrow B\cap X_j=X_j$ or $x\in \darrow B$. In the former case, $\mathcal X_j\subseteq \darrow \mathcal B$, ensuring that $\mathcal X_j^\pi=\emptyset$, which is contrary to the definition of $j\in J$. Therefore we must instead have the latter case, $x\in \darrow B$, meaning $\pi(x)=0$, $\x\in \darrow \mathcal B$ and $\pi(\x)=\{0\}$.
As a result, we find that \be\label{tukt}\varpi_{j-1}(\x)\neq \{0\}\quad\mbox{ if and only if } \quad \pi(\x)\neq \{0\},\quad\mbox{ for all $\x\in \mathcal X_j\cup \{\mathbf v_j\}$ and $j\in J$}.\ee Moreover, since each $\mathbf v_j\notin \darrow \mathcal B\subseteq \mathcal X_1\cup\ldots\cup \mathcal X_s$, the above work yields \be\label{5goop}\varpi_{j-1}(\mathbf v_j)\neq \{0\}\quad\mbox{ for any $j\in J$}.\ee

If $\pi(\x)\neq \{0\}$ for some $\x\in \mathcal X_j\cup \{\mathbf v_j\}$ with $j\in J$, then $\varpi_{j-1}(\x)\neq \{0\}$, meaning $\x\nsubseteq \mathcal E_{j-1}+\mathcal E=\ker \varpi_{j-1}+\mathcal \partial(\x)$ for $\x\in \mathcal X'_j\cup \{\mathbf v_j\}$ and $j\in J$, and thus also $\x\nsubseteq \ker \varpi_i+\partial(\x)$ for any $i<j$. This ensures that, given any $\x\in \bigcup_{i\in J}(\mathcal X_i\cup \{\mathbf v_i\})$ and $j\in [1,s]$, we either have $\varpi_{j-1}(\x)=\{0\}$ or else $\varpi_{j-1}(\x)$ remains a relative half-space. In particular, $\pi(\x)$ is a relative half-space for all $\x\in \bigcup_{i\in J}(\mathcal X'_i\cup \{\mathbf v_i\})$.

Suppose $\x,\,\y\in \bigcup_{i\in J}(\mathcal X_i\cup \{\mathbf v_i\})$ and $j\in [1,s]$ with $\varpi_{j-1}(\x)=\varpi_{j-1}(\y)\neq \{0\}$. We aim to show $\x=\y$. Now  w.l.o.g. $\x\in  \mathcal X_i\cup \{\mathbf v_i\}$ and $\y\in  \mathcal X_{i'}\cup \{\mathbf v_{i'}\}$ with  $i,\,i'\in J$ and $i,\,i'\geq j$.
Moreover, if $i>i'$, then $\varpi_{i-1}(\x)=\varpi_{i-1}\varpi_{j-1}(\x)=\varpi_{i-1}\varpi_{j-1}(\y)=\varpi_{i-1}(\y)=\{0\}$. Thus, since $i\in J$, it follows from \eqref{tukt} that $\varpi_{i-1}(\x)=\pi(\x)=0$, implying $\varpi_{j-1}(\x)=\varpi_{j-1}\pi(\x)=\{0\}$, contrary to assumption. Hence $i\leq i'$, and likewise $i'\leq i$, whence $i=i'$. This argument also shows that $\varpi_{i-1}(\x)=\varpi_{i-1}(\y)\neq \{0\}$. Since $\partial(\x),\,\partial(\y)\subseteq \mathcal E_{i-1}\subseteq \ker \varpi_{i-1}$, it follows from $\varpi_{i-1}(\x)=\varpi_{i-1}(\y)\neq \{0\}$ that  $\R\varpi_{i-1}(x)=\R\varpi_{i-1}(y)\neq \{0\}$. However, since $x,\,y\in X_i\cup \{v_i\}$, this contradicts that $(X_i\cup\{v_i\})^{\varpi_{i-1}}$ is a minimal positive basis of size $|(X_i\cup \{v_i\})\setminus (\mathcal E_{i-1}+\mathcal E)|$ (cf. \eqref{size1}) unless $x=y$, in which case $\x=\y$. In summary, we have just shown that $\varpi_{j-1}(\x)=\varpi_{j-1}(\y)\neq \{0\}$, for $\x,\,\y\in \bigcup_{i\in J}(\mathcal X_i\cup \{\mathbf v_i\})$, implies $\x=\y$. In particular, $\pi=\varpi_0$ is injective on  $\bigcup_{i\in J}(\mathcal X'_i\cup \{\mathbf v_i\})$.

From \eqref{tukt},  \eqref{5goop} and the injectivity of   $\pi$ just established, we conclude that   \be\label{size2}|(X_j\cup \{v_j\})\setminus (\mathcal E_{j-1}+\mathcal E)|=|(X_j\cup \{v_j\})\setminus \mathcal E|=|X_j^\pi|+1\geq 2\quad\mbox{ for any $j\in J$}.\ee
%
%
In view of \eqref{size1}, \eqref{5goop} and \eqref{size2}, it follows that $X_j^{\varpi_{j-1}}\cup \{\varpi_{j-1}(v_j)\}$ is a
 minimal positive basis for $\varpi_{j-1}(\mathcal E_j)$ of size $|X_j^\pi|+1$, for $j\in J$. Hence $\pi(\mathcal R)$ satisfies (OR2) in the definition of an orientated Reay system (note $\varpi_{j-1}\pi=\omega_{j-1}$ since $\ker \pi\leq \ker \varpi_{j-1}$). Also, if  $\x\in \mathcal X'_j\cup \{\mathbf v_j\}$ and $j\in J$, then $\pi(\x)\neq \{0\}$ is a relative half-space with $\varpi_{j-1}(\x)\neq \{0\}$ as shown earlier, so  $\mathbf x\nsubseteq \mathcal E_{j-1}+\ker \pi$, in which case  (OR1) holding for $\mathcal R$ together with \eqref{halfspace-lintrans} implies (OR1) holds for $\pi(\mathcal R)$ with \be\label{witchwhich}\mathcal B_{\pi(\x)}=\mathcal B_\x^\pi=\partial(\{\x\})^\pi,\ee showing that $\pi(\mathcal R)$ is an oriented Reay system for $\pi(\mathcal E_s)$. 
This establishes Item 1.

2. We begin by showing
 \be\label{partial-pit} (\mathcal C^*)^\pi=(\mathcal C^\pi)^* \quad\mbox{ for any $\mathcal C\subseteq \mathcal X_1\cup \ldots\cup \mathcal X_s$}.\ee
Let $\pi(\x)\in (\mathcal C^*)^\pi$ be arbitrary, so $\x\in \mathcal C^*$ with $\pi(\x)\neq \{0\}$. If $\pi(\x)\notin (\mathcal C^\pi)^*$, then there must be some $\y\in \mathcal C$ with $\{0\}\neq \pi(\x)\prec \pi(\y)$. Thus $\x\subseteq \y+\mathcal E\subseteq \R^\cup \la \{\y\}\cup \mathcal B\ra$, in which case Proposition \ref{prop-orReay-BasicProps}.9 ensures that $\x\in \darrow \{\y\}\cup \darrow \mathcal B$ (since $\x,\,\y\in \mathcal C\cup \mathcal B\subseteq \mathcal X_1\cup \ldots\cup \mathcal X_s$).
Since $\pi(\x)\neq \{0\}$, we have $\x\notin \darrow \mathcal B$, forcing $\x\in \darrow \{\y\}$, i.e., $\x\preceq \y$. However, since $\x\in \mathcal C^*$ and $\y\in \mathcal C$, this is only possible if $\x=\y$. Hence $\pi(\x)=\pi(\y)$, contradicting that $\pi(\x)\prec \pi(\y)$. This shows that $(\mathcal C^*)^\pi\subseteq (\mathcal C^\pi)^*.$

Next let $\pi(\x)\in (\mathcal C^\pi)^*$ be arbitrary, so $\x\in \mathcal C$ and $\pi(\x)\neq \{0\}$. If $\x\notin \mathcal C^*$, then there is some $\y\in \mathcal C$ with $\x\prec \y$, implying $\{0\}\neq \pi(\x)\preceq \pi(\y)$.
Thus, since $\pi(\x)\in (\mathcal C^\pi)^*$, we must have $\{0\}\neq \pi(\x)=\pi(\y)$, in which case the injectivity of $\pi$ given in Item 1 implies $\x=\y$, contradicting that $\x\prec \y$.
Therefore we instead conclude that $\x\in \mathcal C^*$, whence $\pi(\x)\in (\mathcal C^*)^\pi$, establishing the reverse inclusion
$(\mathcal C^\pi)^*\subseteq (\mathcal C^*)^\pi$, which establishes \eqref{partial-pit}.

 In view of \eqref{partial-pit}, we have $(\partial(\{\x\})^\pi)^*=(\partial(\{\x\})^*)^\pi=\partial(\{\x\})^\pi$ for any $\x\in \bigcup_{i\in j}(\mathcal X'_i\cup \{\mathbf v_i\})$, which combined with \eqref{witchwhich} and the definition of $\partial(\{\pi(\x)\})$ implies \be\label{trepteeze}\partial(\{\pi(\x)\})=\partial(\{\x\})^\pi\quad\mbox{ for any $\x\in \bigcup_{i\in j}(\mathcal X'_i\cup \{\mathbf v_i\})$}.\ee Let $\mathcal C\subseteq \mathcal X_1\cup \ldots\cup \mathcal X_s\cup \{\mathbf v_j:\;j\in J\}$. Then \begin{align}\nn\partial(\mathcal C^\pi)=&\Big(\bigcup_{\x\in \mathcal C,\,\pi(\x)\neq \{0\}} \partial(\{\pi(\x)\})\Big)^*=\Big(\bigcup_{\x\in \mathcal C,\,\pi(\x)\neq \{0\}} \partial(\{\x\})^\pi\Big)^*=\Big(\Big(\bigcup_{\x\in \mathcal C,\,\pi(\x)\neq \{0\}} \partial(\{\x\})\Big)^\pi\Big)^*\\=&\Big(\Big(\bigcup_{\x\in \mathcal C} \partial(\{\x\})\Big)^\pi\Big)^*=\Big(\Big(\bigcup_{\x\in \mathcal C} \partial(\{\x\})\Big)^*\Big)^\pi=\partial(\mathcal C)^\pi,\label{constantant}\end{align}
 where the first equality follows by definition of $\partial(\mathcal C^\pi)$, the second by \eqref{trepteeze} in view of $\mathcal C\subseteq \mathcal X_1\cup \ldots\cup \mathcal X_s\cup \{\mathbf v_j:\;j\in J\}$, the third is a trivial identity, the fourth follows since $\pi(\x)=\{0\}$ ensures $\pi(\y)=\{0\}$ for all $\y\in \darrow \x$, and thus for all $\y\in \partial(\{\x\})\subseteq \darrow \x$ as well, the fifth  follows by \eqref{partial-pit} applied to $\Big(\bigcup_{\x\in \mathcal C} \partial(\{\x\})\Big)$, and  the sixth follows by definition of $\partial(\mathcal C)$.
 The identity $(\darrow \mathcal C)^\pi=\darrow(\mathcal C^\pi)$ now follows from repeated application of \eqref{constantant} to \eqref{newunion}: \be\label{darroho}(\darrow \mathcal C)^\pi=\Big(\bigcup_{n=0}^{s-1}\partial^n(\mathcal C)\Big)^\pi=\bigcup_{n=0}^{s-1}\big(\partial^n(\mathcal C)\big)^\pi=\bigcup_{n=0}^{s-1}\partial^n(\mathcal C^\pi)=\darrow (\mathcal C^\pi).\ee

 Finally, it remains to establish $(\mathcal C^*)^\pi=(\mathcal C^\pi)^*$  in the case when $\mathcal C\subseteq \mathcal X_1\cup \ldots\cup \mathcal X_s\cup\{\mathbf v_j:\;j\in J\}$. We have $(\mathcal C^\pi)^*\subseteq (\mathcal C^*)^\pi$ by the argument used to establish this inclusion for \eqref{partial-pit}. To see the reverse inclusion $(\mathcal C^*)^\pi\subseteq (\mathcal C^\pi)^*$, let $\pi(\x)\in (\mathcal C^*)^\pi$ be arbitrary, so $\x\in \mathcal C^*$ with $\pi(\x)\neq \{0\}$. If $\pi(\x)\notin (\mathcal C^\pi)^*$, then there must be some $\y\in\mathcal C$ with $\{0\}\neq \pi(\x)\prec \pi(\y)$. Thus \eqref{darroho} implies  $\pi(\x)\in \darrow ( \{\y\}^\pi)=(\darrow \{\mathbf y\})^\pi$. As a result, since $\pi(\x)$ and $\pi(\y)$ are both nonzero, it follows from the injectivity of $\pi$ established in Item 1 that $\x\preceq \y$. However, since $\x\in \mathcal C^*$ and $\y\in \mathcal C$ by hypothesis, this is only possible if $\x=\y$. Thus $\pi(\x)=\pi(\y)$, contradicting that $\pi(\x)\prec \pi(\y)$. This shows that $(\mathcal C^*)^\pi\subseteq (\mathcal C^\pi)^*,$ whence $(\mathcal C^*)^\pi= (\mathcal C^\pi)^*$, which completes the proof of Item 2.

3. Let $\mathcal C\subseteq \mathcal X_1\cup\ldots\cup \mathcal X_s\cup \{\mathbf v_j:\;j\in J\}$ and let $\x\in \mathcal C$. If $\pi(\x)=\{0\}$, then we have $(\mathcal C\setminus \{\x\}\cup \partial(\{\x\}))^\pi=\mathcal C^\pi$. If $\pi(\x)\neq \{0\}$,  then Item 2 and the injectivity of $\pi$ from Item 1 imply  that $(\mathcal C\setminus \{\x\}\cup \partial(\{\x\}))^\pi=\mathcal C^\pi\setminus \{\pi(\x)\}\cup \partial(\{\pi(\x)\})$. This shows $\mathcal C_1\prec \mathcal C_2$ implies $\mathcal C_1^\pi\preceq \mathcal C_2^\pi$ when $\mathcal C_1$ is the immediate predecessor of $\mathcal C_2$, and iterating then yields the first part of Item 3.
Moreover, if $\pi(\x)\neq \{0\}$, then $\mathcal C_1^\pi\prec \mathcal C_2^\pi$.
 Next let $\mathcal D\subseteq \bigcup_{j\in J}(\mathcal X_j^\pi\cup \{\pi(\mathbf v_j)\})$ and $\x\in \pi^{-1}(\mathcal D)$. Then $\pi^{-1}\big(\mathcal D\setminus \{\pi(\x)\}\cup \partial(\{\pi(\x)\})\big)\subseteq \pi^{-1}(\mathcal D)\setminus\{\x\}\cup \partial(\{\x\})\prec \pi^{-1}(\mathcal D)$ in view of Item 2 and the injectivity of $\pi$ given in Item 1, which ensures $\pi^{-1}\big(\mathcal D\setminus \{\pi(\x)\}\cup \partial(\{\pi(\x)\})\big)\prec \pi^{-1}(\mathcal D)$. This shows that $\mathcal D_1\prec \mathcal D_2$ implies $\pi^{-1}(\mathcal D_1)\prec \pi^{-1}(\mathcal D_2)$ when $\mathcal D_1$ is the immediate predecessor of $\mathcal D_2$, and iterating then yields the remaining part of Item 3.

4. Let $\mathcal D\subseteq \bigcup_{i\in J}(\mathcal X_i^\pi\cup \{\pi(\mathbf v_i)\})$ and note that $\pi^{-1}(\mathcal D)\subseteq \bigcup_{j\in J}(\mathcal X'_j\cup \{\mathbf v_j\})$ by definition and Item 1. Thus, Item 2 implies that
\be\label{taiti}\pi\big((\pi^{-1}(\mathcal D))^*\big)=
\big((\pi^{-1}(\mathcal D))^*\big)^\pi=\big((\pi^{-1}(\mathcal D))^\pi\big)^*=(\pi(\pi^{-1}(\mathcal D)))^*=\mathcal D^*.\ee
If $(\pi^{-1}(\mathcal D))^*=\pi^{-1}(\mathcal D)$, then \eqref{taiti} yields $\mathcal D^*=\pi\big((\pi^{-1}(\mathcal D))^*\big)=\pi\big(\pi^{-1}(\mathcal D)\big)=\mathcal D$.  On the other hand,
if $\mathcal D^*=\mathcal D$, then \eqref{taiti}  combined with the injectivity of $\pi$ established in Item 1 implies that $(\pi^{-1}(\mathcal D))^*=\pi^{-1}(\mathcal D)$. This shows that $\mathcal D^*=\mathcal D$ if and only if $(\pi^{-1}(\mathcal D))^*=\pi^{-1}(\mathcal D)$.

Let $\x\in (\pi^{-1}(\mathcal D))^*$ be arbitrary, so $\x\in \bigcup_{j\in J}(\mathcal X'_j\cup \{\mathbf v_j\})$ with $\pi(\x)\in \mathcal D$. If $\x\notin (\pi^{-1}(\mathcal D)\cup \mathcal B)^*$, then there must be some $\y\in \pi^{-1}(\mathcal D)\cup \mathcal B$ with $\x\prec \y$. Since $\x\in (\pi^{-1}(\mathcal D))^*$, we cannot have $\y\in \pi^{-1}(\mathcal D)$, forcing $\x\prec \y\in \mathcal B$. However, since $\ker \pi=\R^\cup \la \mathcal B\ra$, this implies $\pi(\x)=\pi(\y)=\{0\}$, contradicting that $\x\in \bigcup_{j\in J}(\mathcal X'_j\cup \{\mathbf v_j\})$. This  establishes that \be\label{joyl}(\pi^{-1}(\mathcal D))^*\subseteq (\pi^{-1}(\mathcal D)\cup \mathcal B)^*.\ee By Item 2,
$\big(\big(\pi^{-1}(\mathcal D)\cup \mathcal B\big)^*\big)^\pi=
\big(\big(\pi^{-1}(\mathcal D)\cup \mathcal B\big)^\pi\big)^*=\big((\pi^{-1}(\mathcal D))^\pi\cup \mathcal B^\pi\big)^*=
\mathcal D^*$, with the latter equality since $\mathcal B^\pi=\emptyset$, which completes Item 4.

5.  Since $\mathcal D$ is virtual independent, we have $\mathcal D^*=\mathcal D$, whence Item 4 implies $$\pi^{-1}(\mathcal D)=(\pi^{-1}(\mathcal D))^*\subseteq (\pi^{-1}(\mathcal D)\cup \mathcal B)^*.$$ We have $\darrow \mathcal (\pi^{-1}(\mathcal D)\cup \mathcal B)^*= \darrow \pi^{-1}(\mathcal D)\cup \darrow \mathcal B$. By Item 1 and Proposition \ref{prop-orReay-BasicProps}.9,  $\pi$ is injective on $\darrow \pi^{-1}(\mathcal D)\setminus \darrow \mathcal B$, and thus on $\darrow \pi^{-1}(D)\setminus\darrow B$.
By Item 2, $(\darrow \pi^{-1}(\mathcal D))^\pi=\darrow \mathcal D$.
Consequently, if we have a nontrivial linear combination of terms from $\darrow \pi^{-1}(D)\cup \darrow B$ equal to zero, then applying $\pi$ to this linear combination and using that $\mathcal D$ is virtual independent shows that only terms from $\darrow B$ occur in the linear combination, contradicting that $\darrow B$ is linearly independent in view of $\mathcal B\subseteq \mathcal X_1\cup\ldots\cup \mathcal X_s$. This shows that $(\pi^{-1}(\mathcal D)\cup \mathcal B)^*$ is virtual independent. As a result, any subset $\mathcal C\subseteq \darrow(\pi^{-1}(\mathcal D)\cup \mathcal B)^*$ with $\mathcal C^*=\mathcal C$ is also virtual independent. In particular, $\pi^{-1}(\mathcal D)$ is virtual independent.

Next suppose $\mathcal D$ is a support set for $\pi(\mathcal R)$. To show $\pi^{-1}(\mathcal D)$ and $(\pi^{-1}(\mathcal D)\cup \mathcal B)^*$ are support sets, we need  to show  $\mathcal X_j\cup\{\mathbf v_j\}\nsubseteq \darrow\mathcal B\cup \darrow \pi^{-1}(\mathcal D)$, for all $j\in [1,s]$. If $j\notin J$, then $\mathbf v_j\notin \darrow \pi^{-1}(\mathcal D)\subseteq \mathcal X_1\cup \ldots\cup \mathcal X_s\cup \{\mathbf v_i:\;i\in J\}$, while  $\mathcal B\subseteq \mathcal X_1\cup\ldots\cup \mathcal X_s$ ensures that $\mathbf v_j\notin \darrow \mathcal B$, as desired. On the other hand, if  $j\in J$,   then,  since $\mathcal D$ is a support set, we have  $$\mathcal X_j^\pi\cup\{\pi(\mathbf v_j)\}\nsubseteq  \darrow \mathcal D=\darrow \pi(\pi^{-1}(\mathcal D))=(\darrow \pi^{-1}(\mathcal D))^\pi,$$ with the last equality above in view of Item 2, which show that there is some  $\y\in \mathcal X_j\cup \{\mathbf v_j\}$ with $\pi(\y)\neq \{0\}$ and $\y\notin \darrow \pi^{-1}(\mathcal D)$. Note that $\darrow \mathcal B\subseteq \R^\cup \la \darrow \mathcal B\ra=\R^\cup \la \mathcal B\ra=\mathcal E$, implying $\pi(\x)=\{0\}$ for every $\x\in \darrow \mathcal B$. Thus $\y\notin \darrow \mathcal B$ as well, whence  $\mathcal X_j\cup\{\mathbf v_j\}\nsubseteq \darrow\mathcal B\cup \darrow \pi^{-1}(\mathcal D)$ follows, as desired. This  establishes that
$\pi^{-1}(\mathcal D)$ and $(\pi^{-1}(\mathcal D)\cup \mathcal B)^*$ are support sets.

To complete Item 5, now suppose $\mathcal B^*=\mathcal B$. By Item 4, we have $\pi^{-1}(\mathcal D)=\pi^{-1}(\mathcal D)^*\subseteq (\pi^{-1}(\mathcal D)\cup \mathcal B)^*$, implying $(\pi^{-1}(\mathcal D)\cup \mathcal B)^*= \pi^{-1}(\mathcal D)\cup \mathcal B'$ for some $\mathcal B'\subseteq \mathcal B$.
If $\x\in \mathcal B\setminus \mathcal B'$, then we have $\x\prec \y$ for some $\y\in \pi^{-1}(\mathcal D)$ (as $\mathcal B^*=\mathcal B$), implying $\x\subseteq \R^\cup \la\partial(\pi^{-1}(\mathcal D))\ra$, and then $\x\in \darrow \partial(\pi^{-1}(\mathcal D))$ by Proposition \ref{prop-orReay-BasicProps}.9. Hence $\mathcal B\setminus \mathcal B'\subseteq \darrow \partial(\pi^{-1}(\mathcal D))$. On the other hand, if $\y\in \mathcal B\cap \darrow \partial(\pi^{-1}(\mathcal D))$, then $\y\prec \x$ for some $\x\in \pi^{-1}(\mathcal D)$, ensuring that $\y\notin \mathcal B'$ (since $(\pi^{-1}(\mathcal D)\cup \mathcal B)^*=\pi^{-1}(\mathcal D)\cup \mathcal B'$).
 Thus $\mathcal B'\cap \darrow \partial(\pi^{-1}(\mathcal D))=\emptyset$.
  Combining this  with $\mathcal B\setminus \mathcal B'\subseteq \darrow \partial(\pi^{-1}(\mathcal D))$,  we conclude that $\mathcal B'=\mathcal B\setminus \darrow \partial(\pi^{-1}(\mathcal D))$, which completes Item 5.

6. Suppose $\mathcal C\subseteq  \mathcal X_1\cup \{\mathbf v_1\}\cup\ldots\cup \mathcal X_s\cup \{\mathbf v_s\}$ is a support set with $\mathcal B\subseteq \darrow \mathcal C$.
 If $\mathbf v_j\in \mathcal C$ for some $j\notin J$, then $\x\subseteq \ker \pi=\R^\cup \la \mathcal B\ra$ for all $\x\in \mathcal X_j$, whence  Proposition \ref{prop-orReay-BasicProps}.9 ensures that $\mathcal X_j\subseteq \darrow \mathcal B\subseteq \darrow \mathcal C$, in which case $\mathcal X_j\cup \{\mathbf v_j\}\subseteq \darrow \mathcal C$, contradicting that $\mathcal C$ is a support set. Therefore we may assume otherwise, in which case $(\mathcal C^\pi)^*=(\mathcal C^*)^\pi=\mathcal C^\pi$ by Item 2 (as $\mathcal C^*=\mathcal C$ is a support set). But now, if $\mathcal X_j^\pi\cup \{\pi(\mathbf v_j)\}\subseteq \darrow\mathcal C^\pi=(\darrow \mathcal C)^\pi$ for some $j\in J$, with the equality form Item 2, then the injectivity of $\pi$ established in Item 1  ensures that $\mathcal X'_j\cup \{\mathbf v_j\}\subseteq \darrow \mathcal C$. As before, since $\x\subseteq \ker \pi=\R^\cup \la \mathcal B\ra$ for all $\x\in \mathcal X_j\setminus \mathcal X'_j$ (by definition of $\mathcal X'_j$), Proposition \ref{prop-orReay-BasicProps}.9 ensures that $\mathcal X_j\setminus \mathcal X'_j\subseteq \darrow \mathcal B\subseteq \darrow \mathcal C$, in which case $\mathcal X_j\cup \{\mathbf v_j\}\subseteq \darrow \mathcal C$, contradicting that $\mathcal C$ is a support set. Therefore we instead find that $\mathcal X_j^\pi\cup \{\pi(\mathbf v_j)\}\nsubseteq \mathcal C^\pi$ for $j\in J$, which combined with $(\mathcal C^\pi)^*=\mathcal C^\pi$ implies that $\mathcal C^\pi$ is a support set.

Suppose $\mathcal C\subseteq \mathcal X_1\cup \ldots\cup \mathcal X_s$ is a support set. Then $\mathcal C^*=\mathcal C$. Thus  Item 2 implies that $(\mathcal C^\pi)^*=(\mathcal C^*)^\pi=\mathcal C^\pi$ and  $\darrow \mathcal C^\pi=(\darrow \mathcal C)^\pi\subseteq \bigcup_{j\in J}\mathcal X_j^\pi$, ensuring that $\mathcal C^\pi$ is a support set.

 Finally, suppose $\mathcal C\subseteq \mathcal X_1\cup \ldots\cup \mathcal X_s\cup \{\mathbf v_j:\;j\in J\}$ is virtual independent with $\mathcal B\subseteq \darrow \mathcal C$. Then $\mathcal C^*=\mathcal C$, whence Item 2 implies $(\mathcal C^\pi)^*=(\mathcal C^*)^\pi=\mathcal C^\pi$. Also, $\darrow C$ is linearly independent. Thus, since $\mathcal B\subseteq \darrow \mathcal C$ ensures that $\darrow B\subseteq \darrow C$ with $\ker \pi=\R\la \darrow B\ra$ by Proposition \ref{prop-orReay-BasicProps}.1, it follows that $\pi(\darrow C)\setminus \{0\}=\pi(\darrow C\setminus \darrow B)$ is linearly independent, while item 2 ensures $\darrow (\mathcal C^\pi)=(\darrow \mathcal C)^\pi=\pi(\darrow \mathcal C)\setminus \{\{0\}\}$. Hence item 1 implies that $\darrow (C^\pi)=\pi(\darrow  C)\setminus \{0\}$ is linearly independent, implying $\mathcal C^\pi$ is virtual independent.
\end{proof}

If $\mathcal R=(\mathcal X_1\cup \{\mathbf v_1\},\ldots,\mathcal X_s\cup \{\mathbf v_s\})$ is an oriented Reay system,  $\mathcal B\subseteq \mathcal X_1\cup\ldots\cup \mathcal X_s$, and $\pi:\R^d\rightarrow \R^\cup\la \mathcal B\ra^\bot$ is the orthogonal projection,  Proposition \ref{prop-orReay-modulo}.1 implies $\pi(\mathcal R)=(\mathcal X_j^\pi\cup \{\pi(\mathbf v_j)\})_{j\in J}$ is also an oriented Reay system. In particular, either  $\pi(\x)=\{0\}$ or $\pi(\x)$ is a relative half-space, for any $\x\in \mathcal X_1\cup\ldots\cup \mathcal X_s\cup \{\mathbf v_j:\;j\in J\}$; thus $\pi(x)=0$ if and only if $\pi(\x)=\{0\}$, for $\x\in \mathcal X_1\cup\ldots\cup \mathcal X_s\cup \{\mathbf v_j:\;j\in J\}$. Most of the well-behaved properties of $\pi(\mathcal R)$ that pull-back to $\mathcal R$ require $\x\in \mathcal X_1\cup\ldots\cup \mathcal X_s\cup \{\mathbf v_j:\;j\in J\}$. This is a technical restriction that will stay with us throughout  the remainder of this work.

Let $\mathcal D\subseteq \bigcup_{j\in J}(\mathcal X_j^\pi\cup \{\pi(\mathbf v_j)\})$. In view of Proposition \ref{prop-orReay-modulo}.1, $\pi$ is injective on all half-spaces of $\mathcal X_1\cup\ldots\cup\mathcal X_1\cup \{\mathbf v_j:\;j\in J\}$ not mapped to $\{0\}$, so $|\pi^{-1}(\mathcal D)|=|\mathcal D|$. The set $\pi^{-1}(\mathcal D)\subseteq \bigcup_{i\in J}(\mathcal X'_i\cup \{\mathbf v_i\})$ is then the unique subset $\mathcal C\subseteq \mathcal X_1\cup\ldots\cup\mathcal X_s\cup \{\mathbf v_j:\;j\in J\}$ with $\mathcal C^\pi=\pi(\mathcal C)=\mathcal D$. We call $\pi^{-1}(\mathcal D)$ the \textbf{pull-back} of $\mathcal D$ and $(\pi^{-1}(\mathcal D)\cup \mathcal B)^*$ the \textbf{lift} of $\mathcal D$. These sets feature in Proposition \ref{prop-orReay-modulo} and will reoccur often. When $\mathcal D$ is a support/virtual independent set for $\pi(\mathcal R)$, Propositions \ref{prop-orReay-modulo}.5 and \ref{prop-orReay-modulo}.4 ensure that both the pull-back and lift of $\mathcal D$ are support/virtual independent  sets for $\mathcal R$ with $$\pi^{-1}(\mathcal D)=(\pi^{-1}(\mathcal D))^*\subseteq (\pi^{-1}(\mathcal D)\cup \mathcal B)^*\quad\und\quad((\pi^{-1}(\mathcal D)\cup \mathcal B)^*)^\pi=\mathcal D^*=\mathcal D.$$

We continue with the generalization of Proposition \ref{prop-reay-basis-unique-expression} to oriented Reay systems.

\begin{proposition}\label{Prop-orReay-coord}
Let $\mathcal R=(\mathcal X_1\cup\{\mathbf v_1\},\ldots,\mathcal X_s\cup\{\mathbf v_s\})$ be an orientated Reay system for the subspace $\mathcal E_s\subseteq \R^d$. Then $\mathcal E_s=\bigcup_{\mathcal B}\C^\cup( \mathcal B)^\circ$ is a disjoint union, where the union runs over all support sets $\mathcal B$ for $\mathcal R$.
\end{proposition}

\begin{proof}
For $j\in [1,s]$, let $\mathcal E_{j-1}=\R^\cup \la \mathcal X_1\cup\ldots\cup \mathcal
X_{j-1}\ra$ and let $\pi_{j-1}:\R^d\rightarrow \mathcal E_{j-1}^\bot$ be the orthogonal projection.
We proceed by induction on $s$. Let $x\in \mathcal E_s$ be arbitrary. We need to show there is a unique support set $\mathcal B$ with $x\in \C^\cup (\mathcal B)^\circ$.

Since $\pi_{s-1}(X_s\cup \{v_s\})$ is a minimal positive basis of size $|X_s|+1$
(by definition of an oriented Reay system),
it follows from Proposition \ref{prop-reay-basis-unique-expression} that there is a uniquely defined proper subset
$\mathcal B_s\subset \mathcal X_s\cup \{\mathbf v_s\}$ with $\pi_{s-1}(x)\in \C^\cup(\pi_{s-1}(\mathcal B_s))^\circ$, which completes the proof in the base case when $s=1$. Thus we may assume $s\geq 2$. Let $\z_1,\ldots,\z_\ell\in \mathcal B_s$ be the distinct elements of $\mathcal B_s$, let $\mathcal E=\R^\cup \la \partial(\mathcal B_s)\ra$, and let $\pi:\R^d\rightarrow \mathcal E^\bot$ be the orthogonal projection. Since $\partial(\mathcal B_s)\subseteq \mathcal X_1\cup\ldots\cup \mathcal X_{s-1}$, we have  $\ker \pi=\mathcal E\subseteq \mathcal E_{s-1}=\ker \pi_{s-1}$, ensuring $\pi_{s-1}\pi=\pi_{s-1}$.
 By Proposition \ref{prop-char-minimal-pos-basis}, $\pi_{s-1}(B_s)$ is linearly independent, so there is a unique linear combination of the elements of $\pi_{s-1}(B_s)$ equal to $\pi_{s-1}(x)$, say $\Sum{i=1}{\ell}\alpha_i\pi_{s-1}(z_i)=\pi_{s-1}(x)$ with $\alpha_i\in \R$, and this linear combination is strictly positive since $\pi_{s-1}(x)\in \C^\cup(\pi_{s-1}(\mathcal B_s))^\circ$, so $\alpha_i>0$ for all $i$. Thus \be\label{uniquedefin}\Sum{i=1}{\ell}\alpha_i\pi(z_i)=\pi(x)-\xi=\pi(x-\xi)\quad\mbox{ for some $\xi\in \mathcal E_{s-1}\cap \mathcal E^\bot$}.\ee
 If there were some other linear combination $\Sum{i=1}{\ell}\alpha'_i\pi(z_i)=\pi(x)-\xi'$ for some $\xi'\in \mathcal E_{s-1}\cap \mathcal E^\bot$ and $\alpha'_i\in \R$, then, applying $\pi_{s-1}$ to this linear combination and using that $\pi_{s-1}\pi=\pi_{s-1}$, we conclude from the uniqueness of the linear combination $\Sum{i=1}{\ell}\alpha_i\pi_{s-1}(z_i)=\pi_{s-1}(x)$ that $\alpha'_i=\alpha_i$ for all $i$, whence $\xi=\xi'$ as well. In consequence, the $\alpha_i$ and $\xi\in \mathcal E_{s-1}\cap \mathcal E^\bot$ in \eqref{uniquedefin} are uniquely defined.

 Note $\mathcal B_s\subset \mathcal X_s\cup \{\mathbf v_s\}$ is a support set for $\mathcal R$ since $\x\prec \y$ is impossible for half-spaces from the same level $\mathcal X_s\cup \{\mathbf v_s\}$. Thus Proposition \ref{prop-orReay-BasicProps} ensures that $\C^\cup (\mathcal B_s)^\circ = \z_1^\circ+\ldots+\z_\ell^\circ$ with $\overline{\C^\cup (\mathcal B_s)} = \overline \z_1+\ldots+\overline \z_\ell$ a blunted simplicial cone having lineality space $\mathcal E=\partial(\z_1)+\ldots+\partial(\z_\ell)=\R^\cup\la \partial(\mathcal B_s)\ra\subseteq \mathcal E_{s-1}$. As a result, since $\alpha_i>0$ for all $i$,  \eqref{uniquedefin} implies \be\label{reju}x-\xi+\mathcal E\subseteq \C^\cup(\mathcal B_s)^\circ=\z_1^\circ+\ldots+\z_\ell^\circ.\ee Moreover, $\xi=\pi(\xi)$ is the unique element from $\mathcal E_{s-1}\cap \mathcal E^\bot$ such that $x-\xi\in \C^\cup(\mathcal B_s)^\circ$.

 In view of Proposition \ref{prop-orReay-modulo}.1, $$\pi(\mathcal R)=\Big(\mathcal X^\pi_{j}\cup\{\pi(\mathbf v_{j})\}\Big)_{j\in J}$$ is an orientated Reay system for $\pi(\mathcal E_s)$ for some $J\subseteq [1,s]$.
 Moreover, since $\mathcal E\subseteq \mathcal E_{s-1}$, we have $s\in J$. Since $\xi\in \mathcal E_{s-1}\cap \mathcal E^\bot=\pi(\mathcal E_{s-1})$, we can apply the induction hypothesis to the oriented Reay system
 $\Big(\mathcal X_{j}^\pi\cup\{\pi(\mathbf v_{j})\}\Big)_{j\in J\setminus \{s\}}$ to conclude there exists a unique support set $\mathcal D\subseteq \bigcup_{j\in J\setminus\{s\}}(\mathcal X_j^\pi\cup \{\pi(\mathbf v_j)\})$ with \be\label{jumpcliff}\xi\in \C^\cup (\mathcal D)^\circ=\pi(\y_1)^\circ+\ldots+\pi(\y_r)^\circ=\pi(\y_1^\circ)+\ldots+\pi(\y_r^\circ),\ee where $\y_1,\ldots,\y_r\in \pi^{-1}(\mathcal D)$ are the distinct half-spaces from $\pi^{-1}(\mathcal D)$, with the first equality in \eqref{jumpcliff} in view of Proposition \ref{prop-orReay-BasicProps}.4 applied to the support set $\mathcal D$, and the second holding since each $\pi(\y_i)$ remains a relative half-space in $\pi(\mathcal R)$.

 Since $\mathcal D$ is a support set, Proposition \ref{prop-orReay-modulo}.5 ensures that $\pi^{-1}(\mathcal D)$ and $(\pi^{-1}(\mathcal D)\cup \mathcal \partial(\mathcal B_s))^*$ are both support sets with $\pi^{-1}(\mathcal D)\subseteq (\pi^{-1}(\mathcal D)\cup \mathcal \partial(\mathcal B_s))^*$. Set $$\mathcal B:=\pi^{-1}(\mathcal D)\cup \mathcal B_s.$$ If $\mathcal B^*\neq \mathcal B$, then there are $\y_1,\,\y_2\in \pi^{-1}(\mathcal D)\cup \mathcal B_s$ with $\y_1\prec \y_2$. Since $\pi^{-1}(\mathcal D)^*=\pi^{-1}(\mathcal D)\subseteq \mathcal X_1\cup\{\mathbf v_1\}\cup\ldots\cup\mathcal X_{s-1}\cup \{\mathbf v_{s-1}\}$ (as $\pi^{-1}(\mathcal D)$ is a support set) and $\mathcal B_s\subset \mathcal X_s\cup \{\mathbf v_s\}$, this is only possible if $\y_1\in \pi^{-1}(\mathcal D)$ and $\y_2\in \mathcal B_s$, in which case $\y_1\in \darrow \partial(\mathcal B_s)$, whence $\pi(\y_1)=\{0\}$ as $\mathcal E=\R^\cup\la \partial(\mathcal B_s)\ra$, contradicting that $\y_1\in \pi^{-1}(\mathcal D)$ with $\mathcal D$ a set of non-zero half-spaces. Therefore we instead conclude that $\mathcal B^*=\mathcal B$. Since $\pi^{-1}(\mathcal D)\subseteq \mathcal X_1\cup\{\mathbf v_1\}\cup\ldots\cup\mathcal X_{s-1}\cup \{\mathbf v_{s-1}\}$ and $\mathcal B_s\subset \mathcal X_s\cup \{\mathbf v_s\}$, we have $\mathcal X_s\cup\{\mathbf v_s\}\nsubseteq \darrow \mathcal B$. Since $\mathcal B_s\subseteq \mathcal X_s\cup \{\mathbf v_s\}$, any   $\y\in \mathcal X_j\cup\{\mathbf v_j\}$  with $j<s$ and  $\y\in \darrow \mathcal B$ must have  $\y\in \darrow \pi^{-1}(\mathcal D)$ or $\y\in \darrow \partial(\mathcal B_s)$. Thus $\mathcal X_j\cup \{\mathbf v_j\}\subseteq \darrow \mathcal B$, for $j<s$, would imply $\mathcal X_j\cup \{\mathbf v_j\}\subseteq \darrow (\pi^{-1}(\mathcal D)\cup \partial(\mathcal B_s))^*$, contradicting that $(\pi^{-1}(\mathcal D)\cup \partial(\mathcal B_s))^*$ is a support set. This shows that $\mathcal B=\mathcal B^*$ is itself a support set.

 Since $\mathcal B=\pi^{-1}(\mathcal D)\cup \mathcal B_s$ is a support set, Proposition \ref{prop-orReay-BasicProps}.4 implies \be\label{coffeacup}\C^\cup (\mathcal B)^\circ=\C^\cup (\pi^{-1}(\mathcal D))^\circ+\C^\cup (\mathcal B_s)^\circ.\ee
  Since $\pi^{-1}(\mathcal D)$ is a support set, Proposition \ref{prop-orReay-BasicProps}.4 implies $\C^\cup(\pi^{-1}(\mathcal D))^\circ=\y_1^\circ+\ldots+\y_r^\circ$.
  Thus \eqref{jumpcliff} implies  that    $\xi=\pi(\xi)\in \C^\cup (\pi^{-1}(\mathcal D))^\circ+\mathcal E$, which combined with \eqref{reju} and \eqref{coffeacup} yields   \be\label{gusterbust} x\in \C^\cup(\pi^{-1}(\mathcal D))^\circ+\C^\cup( \mathcal B_s)^\circ=\C^\cup (\mathcal B)^\circ.\ee  It remains to establish the uniqueness of $\mathcal B$.

Now suppose $\mathcal B'$ were any support set for $\mathcal R$ with $x\in \C^\cup(\mathcal B')^\circ$. We need to show it equals the support set $\mathcal B$ constructed above. Note   $\pi_{s-1}(x)\in \C^\cup \Big(\pi_{s-1}\Big( \mathcal B'\cap (\mathcal X_s\cup \{\mathbf v_s\})\Big)\Big)^\circ$ in view of Proposition \ref{prop-orReay-BasicProps}.4, in which case we must have $\mathcal B_s=\mathcal B'\cap (\mathcal X_s\cup \{\mathbf v_s\})$ by the uniqueness property established with the  existence of $\mathcal B_s$. Let $\mathcal C=\mathcal B'\setminus \mathcal B_s$. We need to show $\mathcal C=\pi^{-1}(\mathcal D)$.

Since $\partial(\mathcal B_s)\subseteq \darrow \mathcal B_s\subseteq \darrow \mathcal B'$, Proposition \ref{prop-orReay-modulo}.6 implies that $(\mathcal B')^\pi=(\mathcal C\cup \mathcal B_s)^\pi=\mathcal C^\pi\cup \mathcal B_s^\pi$ is a support set for $\pi(\mathcal R)$ with $\mathcal C\subseteq \mathcal X_1\cup\ldots\cup \mathcal X_{s-1}\cup \{\mathbf v_j:\;j\in J\setminus\{s\}\}$.
If there were some $\y\in \mathcal C$ with $\pi(\y)=\{0\}$, then, since $\pi(\mathbf v_j)\neq \{0\}$ for all $j\in J$ by Proposition \ref{prop-orReay-modulo}.1, we must have $\y\in \mathcal X_1\cup \ldots\cup \mathcal X_{s-1}$ and $\y\subseteq \ker \pi=\R^\cup \la \partial(\mathcal B_s)\ra$, whence Proposition \ref{prop-orReay-BasicProps}.9 implies that $\y\in \darrow \partial(\mathcal B_s)\subseteq \darrow \mathcal B_s$. However this contradicts that  $(\mathcal C\cup \mathcal B_s)^*=(\mathcal C\cup \mathcal B_s)$ for the support set $\mathcal B'=\mathcal C\cup \mathcal B_s$. Therefore, we conclude that  $\mathcal C^\pi=\pi(\mathcal C)$, ensuring $\mathcal C=\pi^{-1}(\pi(\mathcal C))$ in view of the injectivity of $\pi$ given in Proposition \ref{prop-orReay-modulo}.1. It remains to show $\mathcal C^\pi=\pi(\mathcal C)=\mathcal D$.

Since $\mathcal B'=\mathcal C\cup \mathcal B_s$ and $(\mathcal B')^\pi=\mathcal C^\pi\cup \mathcal B_s^\pi$ are  support sets,  $\mathcal C$ and $\mathcal C^\pi$ are also support sets.
We have $x\in \C^\cup (\mathcal C\cup \mathcal B_s)^\circ=\C^\cup (\mathcal C)^\circ+\C^\cup (\mathcal B_s)^\circ$, ensuring $\pi(x)\in \C^\cup (\mathcal C^\pi)^\circ +\C^\cup(\mathcal B_s^\pi)^\circ$, both in view of Proposition \ref{prop-orReay-BasicProps}.4 (since $\mathcal C\cup \mathcal B_s$ and $\mathcal C^\pi\cup \mathcal B_s^\pi$ are both  support sets). Thus $\pi(x)-\xi'\in \C^\cup(\mathcal B_s^\pi)$ for some $\xi'\in \C^\cup (\mathcal C^\pi)^\circ
\subseteq \pi(\mathcal E_{s-1})=\mathcal E_{s-1}\cap \mathcal E^\bot$, in which case the uniqueness of $\xi$ given in \eqref{uniquedefin} ensures that $\xi'=\xi$. Hence $\xi=\xi'\in \C^\cup (\mathcal C^\pi)^\circ$ with $\mathcal C^\pi$ a support set, in which case the uniqueness property established with the  existence of $\mathcal D$ ensures that $\mathcal  D=\mathcal C^\pi=\pi(\mathcal C)$, completing the proof.
\end{proof}

In view of Proposition \ref{Prop-orReay-coord}, given any $x\in \mathcal E_s$, there is a uniquely defined support set $\mathcal B$ for $\mathcal R$ with $x\in \C^\cup( \mathcal B)^\circ$, which we denote by $\mathcal B=\supp_{\mathcal R}(x)$.  Given a support set $\mathcal B$ for $\mathcal R$, we define \be\label{weight-deff}\wt(\mathcal B)=|\mathcal B\cap \{\mathbf v_1,\ldots,\mathbf v_s\}|.\ee Likewise, for $x\in\mathcal E_s$, we let $\wt(x)=\wt(\supp_\mathcal R(x))$.

\begin{definition} Let $\vec u=(u_1,\ldots,u_t)$ be a tuple of orthonormal vectors from $\R^d$, let $\mathcal R=(\mathcal X_1\cup\{\mathbf v_1\},\ldots,\mathcal X_s\cup\{\mathbf v_s\})$ be an oriented ray system, and let $\mathcal B\subseteq \mathcal X_1\cup\{\mathbf v_1\}\cup\ldots\cup \mathcal X_s\cup\{\mathbf v_s\}$. We say that $\mathcal B$ \textbf{encases} $\vec u$ provided $\C^\cup(\mathcal B)$ encases $\vec u$. We say that $\mathcal B$ \textbf{minimally encases} $\vec u$ if $\mathcal B$ encases $\vec u$ but no  proper $\mathcal B'\prec \mathcal B$ encases $\vec u$.
\end{definition}

Since $\C^\cup (\mathcal B^*)=\C^\cup (\mathcal B)$ and $\mathcal B^*\preceq \mathcal B$, it is clear that $\mathcal B$ minimally encasing $\vec u$ is only possible if $\mathcal B^*=\mathcal B$. Thus any subset $\mathcal B\subseteq \mathcal X_1\cup\ldots\cup \mathcal X_s$ which minimally encases $\vec u$ must be a support set. The following lemma deals with  minimal encasement for $t=1$.

\begin{lemma}\label{lem-orReay-mincase-t=1}Let
$\mathcal R=(\mathcal X_1\cup\{\mathbf v_1\},\ldots,\mathcal X_s\cup\{\mathbf v_s\})$ be an oriented ray system for a subspace $\mathcal E_s\subseteq \R^d$, let $\mathcal B\subseteq \mathcal X_1\cup\{\mathbf v_1\}\cup\ldots\cup \mathcal X_s\cup\{\mathbf v_s\}$, and let $u_1\in \R^d$ be a unit vector.
 \begin{itemize}
 \item[1.] If $\mathcal B$ minimally encases $-u_1$, then $-u_1\in \z_1^\circ+\ldots+\z_\ell^\circ$, where $\z_1,\ldots,\z_\ell\in \mathcal B$ are the distinct half-spaces from $\mathcal B$.
 \item[2.] If $\mathcal B$ is virtual independent, then $\mathcal B$ minimally encases $-u_1$ if and only if $-u_1\in \C^\cup (\mathcal B)^\circ$.

 \end{itemize}
\end{lemma}

\begin{proof}1.
Suppose $\mathcal B$ minimally encases $-u_1$.   Then $-u_1\in \C^\cup(\mathcal B)= \z_1+\ldots+\z_\ell$ (by Proposition \ref{prop-orReay-BasicProps}.2). If $-u_1\in (\partial(\z_1)\cap \z_1)+\z_2+\ldots+\z_\ell=\C(\partial(\{\z_1\}))+\z_2+\ldots+\z_\ell$, then $-u_1\in \C(\mathcal B')$, where $\mathcal B'=(\mathcal B\setminus\{\z_1\})\cup \partial(\{\z_1\})\prec \mathcal B$, contradicting the minimality of $\mathcal B$. A similar argument can be used for $\z_2,\ldots,\z_\ell$, and we conclude that  $-u_1\in \z_1^\circ+\ldots+\z_\ell^\circ$, as desired.

2. Since $\mathcal B$ is virtual independent, Proposition \ref{prop-orReay-BasicProps}.4 implies that $\C^\cup(\mathcal B)^\circ=\z_1^\circ+\ldots+\z_\ell^\circ$. Thus, if  $\mathcal B$ minimally encases  $-u_1$, then Item 1 implies $-u_1\in \C^\cup(\mathcal B)^\circ$. Next suppose that $-u_1\in \C^\cup(\mathcal B)^\circ=\z_1^\circ+\ldots+\z_\ell^\circ$.  Let $\pi:\R^d\rightarrow \R^\cup \la \partial(\mathcal B)\ra^\bot$ be the orthogonal projection. Since  $-u_1\in \z_1^\circ+\ldots+\z_\ell^\circ$, it follows that $-\pi(u_1)$ can be written as a strictly positive linear combination of the distinct linearly independent elements $\pi(z_1),\ldots,\pi(z_\ell)\in \R^\cup\la \partial(\mathcal B)\ra^\bot$ (in view of Proposition \ref{prop-orReay-BasicProps}.3). Since the $\pi(z_i)$ are distinct and linearly independent, this is then the unique way to write $-\pi(u_1)$ as a linear combination of the elements $\pi(z_1),\ldots,\pi(z_\ell)\in \R^\cup\la \partial(\mathcal B)\ra^\bot$.
  If $\mathcal B$ does not minimally encase $-u_1$, then w.l.o.g. $-u_1\in \C^\cup\big(\mathcal B\setminus\{\z_1\}\cup \partial(\{\z_1\})\big)= \C^\cup(\partial(\{\z_1\}))+\z_2+\ldots+\z_\ell= (\partial(\z_1)\cap \z_1)+\z_2+\ldots+\z_\ell$. Thus $-\pi(u_1)$ can be also be written as a positive linear combination of the  elements $\pi(z_2)\ldots,\pi(z_\ell)\in \R^\cup\la \partial(\mathcal B)\ra^\bot$, contradicting that the unique way to express $-\pi(u_1)$ has the coefficient of $\pi(z_1)$ being nonzero.
\end{proof}

Any support set $\mathcal B$ is always virtual independent (by Proposition \ref{prop-orReay-BasicProps}.5). Consequently, since $-u_1\in \C(\supp_{\mathcal R}(-u_1))^\circ$ by definition of $\supp_{\mathcal R}$, we conclude via Lemma \ref{lem-orReay-mincase-t=1} that $\mathcal B=\supp_{\mathcal R}(-u_1)$ is always a support set which minimally encases $-u_1\in \R^\cup \la \mathcal X_1\cup \ldots\cup \mathcal X_s\ra$, while  Proposition \ref{Prop-orReay-coord} and Lemma \ref{lem-orReay-mincase-t=1} ensure  that $\supp_{\mathcal R}(-u_1)$ is the unique
  \emph{support} set for $\mathcal R$ which minimally encases $-u_1$, though other non-support sets may also do the same.
 Indeed, if  $\mathcal B$ is a virtual independent set  and
  $\pi:\R^d\rightarrow \R^\cup\la \partial(\mathcal B)\ra^\bot$ is the orthogonal projection,
then Lemma \ref{lem-orReay-mincase-t=1} and  Proposition \ref{prop-orReay-BasicProps} (Items 3--4) ensure that $\mathcal B$ minimally encasing the \emph{nonzero} element  $-u_1$ is equivalent to $\pi(B)\cup \{\pi(u_1)\}$ being a minimal positive basis.

As the above discussion shows, there is always a unique support set $\mathcal B=\supp_{\mathcal R}(-u_1)$ which minimally encases the element $u_1\in \R^\cup \la \mathcal X_1\cup \ldots\cup \mathcal X_s\ra$. However, if $\vec u=(u_1,\ldots,u_t)$ is a tuple of  orthonormal vectors $u_1,\ldots,u_t\in \R^\cup \la \mathcal X_1\cup \ldots\cup \mathcal X_s\ra$, there is no guarantee that $-\vec u$ will be minimally encased by some  support set from $\mathcal R$ when $t\geq 2$. We will later show that this problem does not occur if we  impose additional conditions on $\mathcal R$. However, until we can achieve this, we will have need of the following definition.

\begin{definition} Let $\mathcal R=(\mathcal X_1\cup\{\mathbf v_1\},\ldots,\mathcal X_s\cup\{\mathbf v_s\})$ be an oriented Reay system in $\R^d$,  let $\mathcal B\subseteq \mathcal X_1\cup\{\mathbf v_1\}\cup\ldots\cup \mathcal X_s\cup\{\mathbf v_s\}$, and let $\vec u=(u_1,\ldots,u_t)$ be a tuple of orthonormal vectors from $\R^d$. Suppose $\mathcal B$ minimally encases $\vec u$. If $\mathcal B\neq \emptyset$, then $t\geq 1$ and there will be a maximal index  $t'\in [0,t-1]$  such that $\mathcal B$ does \emph{not} minimally encase $(u_1,\ldots,u_{t'})$. Note $\mathcal B$ still encases $(u_1,\ldots,u_{t'})$, so there is some $\mathcal A\prec \mathcal B$ such that $\mathcal A$ minimally encases $(u_1,\ldots,u_{t'})$. If $\mathcal B=\emptyset$ or \begin{align*}\mathcal B\subseteq \mathcal X_1\cup \ldots\cup \mathcal X_s\cup \{\mathbf v_j:\;j\in J\} \mbox{ is virtual independent}\quad\und \quad  \quad \mathcal A\subseteq \mathcal X_1\cup\ldots\cup \mathcal X_s,\end{align*} where $\pi(\mathcal A)=(\mathcal X_i^\pi\cup \{\pi(\mathbf v_i)\})_{i\in J}$ and $\pi:\R^d\rightarrow \R^\cup\la \mathcal A\ra^\bot$ is the orthogonal projection, then we say that $\mathcal B$ minimally encases $\vec u$  \textbf{urbanely}.\end{definition}

We remark that,
 when $\mathcal B$ is a support set, the condition $\mathcal B\subseteq \mathcal X_1\cup \ldots\cup \mathcal X_s\cup \{\mathbf v_j:\;j\in J\}$ in the above definition holds automatically  in view of Proposition \ref{prop-orReay-modulo}.6, and so can be dropped.  We will later see in Proposition \ref{prop-orReay-minecase-char}.1   that the $\mathcal A$ occurring in the definition of urbane minimal encasement  is uniquely defined.
If $\mathcal B\subseteq \mathcal X_1\cup\ldots\cup \mathcal X_s$ and $\mathcal B$ minimally encases $-\vec u$, then it must always do so urbanely and be a support set in view of $\mathcal A\subseteq \darrow \mathcal B\subseteq \mathcal X_1\cup \ldots\cup \mathcal X_s$. However, for more general subsets $\mathcal B\subseteq \mathcal X_1\cup\{\mathbf v_1\}\cup\ldots\cup \mathcal X_s\cup \{\mathbf v_s\}$, it is possible for $\mathcal B$ to minimally encase $-\vec u$ non-urbanely. The following proposition contains the basic properties regarding urbane minimal encasement and is the analogue of Propositions \ref{prop-min-encasement-char} and \ref{prop-min-encasement-minposbasis} for oriented Reay systems.

 \begin{proposition}\label{prop-orReay-minecase-char}
 Let $\mathcal R=(\mathcal X_1\cup\{\mathbf v_1\},\ldots,\mathcal X_s\cup\{\mathbf v_s\})$ be an oriented ray system in $\R^d$, let $\vec u=(u_1,\ldots,u_t)$ be a tuple of $t\geq 0$ orthonormal vectors in $\R^d$,  and let $\mathcal B\subseteq \mathcal X_1\cup \mathcal \{\mathbf v_1\}\cup\ldots\cup \mathcal X_s\cup \{\mathbf v_s\}$. 

 \begin{itemize}
 \item[1.] If $\mathcal B,\,\mathcal C\subseteq  \mathcal X_1\cup \mathcal \{\mathbf v_1\}\cup\ldots\cup \mathcal X_s\cup \{\mathbf v_s\}$ are support sets that both minimally encase $\vec u$ urbanely, then $\mathcal B=\mathcal C$.  In particular, if $\mathcal B,\,\mathcal C\subseteq  \mathcal X_1\cup \ldots\cup \mathcal X_s$
      with $\mathcal B$ minimally encasing $\vec u$ and $\mathcal C$ encasing $\vec u$, then $\mathcal B\preceq \mathcal C$.
 \item[2.]  $\mathcal B$  minimally encases $-\vec u$ urbanely if and only if  there are indices $1=r_1<\ldots<r_\ell<r_{\ell+1}=t+1$ and virtual independent sets $\mathcal C_i$, for $i=0,1,\ldots,\ell$, satisfying  $$\emptyset=\mathcal C_0\prec\mathcal C_1\prec \ldots\prec\mathcal C_{\ell-1}\subseteq \mathcal X_1\cup \ldots\cup \mathcal X_s\quad\und \quad \mathcal C_{\ell-1}\prec \mathcal C_\ell=\mathcal B$$ such that either $\ell=0$ or else
 \begin{itemize}
 \item[(a)]  $\mathcal C_\ell^{\pi_{\ell-1}}=\mathcal D_\ell$ for some virtual independent set $\mathcal D_\ell$ from $\pi_{\ell-1}(\mathcal R)$ that minimally encases $-\pi_{\ell-1}(u_{r_\ell})$ and $\mathcal C_\ell$ contains no $\mathbf v_i$ with $\mathcal X_i^{\pi_{\ell-1}}=\emptyset$ (the latter which necessarily holds when $\mathcal C_\ell$ is a support set),
 \item[(b)] $\mathcal C_j^{\pi_{j-1}}=\mathcal D_j:=\supp_{\pi_{j-1}(\mathcal R)}(-\pi_{j-1}(u_{r_j}))$ for every $j\in [1,\ell-1]$, and
 \item[(c)] $u_i\in \R^\cup\la \mathcal C_j\ra$ for all $i<r_{j+1}$ and $j\in [1,\ell]$,
 \end{itemize}
 where $\pi_{j-1}:\R^d\rightarrow \R^\cup\la \mathcal C_{j-1}\ra^\bot$  is the orthogonal projection for $j\in [1,\ell]$.

 Moreover, $\mathcal B$ is a support set if and only if $\ell=0$ or $\mathcal D_\ell=\supp_{\pi_{\ell-1}(\mathcal R)}(-\pi_{\ell-1}(u_{r_\ell}))$.
 \end{itemize}
 Now assume the conditions of Item 2 hold along with the relevant notation.
 \begin{itemize}
 \item[3.] $\mathcal F=(\R^\cup \la \mathcal C_1\ra,\ldots,\R^\cup\la \mathcal C_\ell\ra)$ is a compatible filter for $\vec u$ with associated indices $1=r_1<\ldots<r_\ell<r_{\ell+1}=t+1$, $\mathcal F(\vec u)=(\overline u_{1},\ldots,\overline u_{\ell})$, and for every $j\in [1,\ell]$, the following hold:
 \begin{itemize}
 \item[(a)] $\mathcal C_j=(\pi_{j-1}^{-1}(\mathcal D_j)\cup \mathcal C_{j-1})^*$ is the lift of $\mathcal D_j$,
 \item[(b)] $\C^\cup(\mathcal C_j)+\C(u_{r_1},\ldots,u_{r_j})=\C^\cup(\mathcal C_j)+\C(\overline u_{1},\ldots,\overline u_{j})=\R^\cup\la \mathcal C_j\ra$, and
     \item[(c)] $\mathcal B$ minimally encases $-\mathcal F(\vec u)$ and  $-(u_{r_1},\ldots,u_{r_\ell})$ urbanely.
 \end{itemize}

 \item[4.]
     Let $\mathcal A\subseteq \mathcal X_1\cup\ldots\cup \mathcal X_s$, let $\mathcal E=\R^\cup\la \mathcal A\ra$, let $\pi:\R^d\rightarrow \mathcal E^\bot$ be the orthogonal projection, let $J\subseteq [1,\ell]$ be all those indices $j\in [1,\ell]$ with $\mathcal C^\pi_{j-1}\prec \mathcal C_j^\pi$, and
     let $\mathcal F^\pi=(\R^\cup\la \mathcal C_j^\pi\ra)_{j\in J}$. For $j\in [1,\ell]$, let $\tau_{j-1}:\R^d\rightarrow (\mathcal E+\R^\cup\la \mathcal C_{j-1}\ra)^\bot$ be the orthogonal projection.
     Suppose \begin{align*}&\mathcal B\subseteq \mathcal X_1\cup \ldots\cup \mathcal X_s\quad \mbox{ or } \quad\mathcal A\subseteq \darrow \mathcal B,\quad\mbox{ and either }\\ &\mathcal B\mbox{ is a support set } \quad \mbox{or }\quad \mbox{ $\mathcal B$ contains no $\mathbf v_i$ with $\mathcal X_i^{\tau_{\ell-1}}=\emptyset$ for $i\in [1,s]$}.\end{align*}
     Then
     $\mathcal F^\pi$ is a compatible filter for $\pi(\vec u)$ with
     $$\mathcal F^\pi(\pi(\vec u))=(\overline {u^*_{i}})_{i\in J},$$
      where $ \overline{u^*_{i}}=\tau_{i-1}(u_{r_i})/\|\tau_{i-1}(u_{r_i})\|$, and the virtual independent set  $\mathcal B^\pi$ for $\pi(\mathcal R)$ minimally encases $\pi(\vec u)$ urbanely with the  sets $\mathcal C_j^\pi$ for $j\in J$ those satisfying  Item 2 for $\pi(\vec u)$. \end{itemize}
 \end{proposition}

\begin{proof}
1. We begin by proving Item 1 in the case when $\mathcal B,\,\mathcal C\subseteq \mathcal X_1\cup \ldots\cup \mathcal X_s$. For this, it suffices to consider subsets $\mathcal B,\,\mathcal C\subseteq \mathcal X_1\cup\ldots\cup \mathcal X_s$ which minimally encase $\vec u$ and show $\mathcal B=\mathcal C$, in which case $\mathcal B$ and $\mathcal C$ are both support sets.   Let $\mathcal E= \R^\cup\la \mathcal B\ra$, let $\pi:\R^d\rightarrow \mathcal E^\bot$ be the orthogonal projection, and  let $\mathcal C'\subseteq \mathcal C$ be all those half-spaces $\x\in \mathcal C$  with $\pi(\x)\neq \{0\}$.
By Proposition \ref{prop-orReay-modulo}.6, $\mathcal C^\pi=\pi(\mathcal C')$ is a support set for $\pi(\mathcal R)$ with $\pi$ injective on $\mathcal C'$ (by Proposition \ref{prop-orReay-modulo}.1), while Proposition \ref{prop-orReay-modulo}.2 implies that $\partial(\mathcal C^\pi)=\partial(\mathcal C)^\pi$.
Thus Proposition \ref{prop-orReay-BasicProps} (Items 3 and 5) applied to the support set $\mathcal C^\pi$ for $\pi(\mathcal R)$ implies that the elements of $C'$ are linearly independent modulo $\mathcal E+\R\la \partial(\mathcal C)\ra$. Consequently,  \be\label{drysal}\C^\cup(\mathcal C)\cap \mathcal E=(\C^\cup(\mathcal C')\cap \mathcal E)+\C^\cup(\mathcal C\setminus \mathcal C')\subseteq \C^\cup(\partial(\mathcal C'))+\C^\cup(\mathcal C\setminus \mathcal C').\ee
Since $\mathcal C$ encases $(u_1,\ldots,u_t)$ with $\R\la u_1,\ldots, u_t\ra\subseteq \R^\cup \la \mathcal B\ra=\mathcal E$ (as $\mathcal B$ encases $(u_1,\ldots,u_t)$), it follows that $\C^\cup(\mathcal C)\cap \mathcal E$ also encases $(u_1,\ldots,u_t)$, whence \eqref{drysal} ensures that $(\mathcal C\setminus \mathcal C')\cup \partial(\mathcal C')\preceq \mathcal C$ encases $(u_1,\ldots,u_t)$. Consequently, since $\mathcal C$ \emph{minimally} encases $(u_1,\ldots,u_t)$, we conclude that  $(\mathcal C\setminus \mathcal C')\cup \partial(\mathcal C')= \mathcal C$, which is only possible if $\mathcal C'=0$, that is, if all half-spaces in $\mathcal C\subseteq \mathcal X_1\cup\ldots\cup \mathcal X_s$ are contained in $\mathcal E=\R^\cup\la \mathcal B\ra$. Hence Proposition \ref{prop-orReay-BasicProps}.9 implies that $\mathcal C\subseteq \darrow \mathcal B$, further implying $\mathcal C=\mathcal C^*\preceq \mathcal B$ (the equality follows since $\mathcal C$ is a support set). But now, since $\mathcal B$ \emph{minimally} encases $(u_1,\ldots,u_t)$, we must have $\mathcal C=\mathcal B$, completing the proof of Item 1 in the case $\mathcal B,\,\mathcal C\subseteq \mathcal X_1\cup \ldots\cup \mathcal X_s$. We will complete the more general case  for Item 1 later.

2--3. If $t=0$, then Item 2 holds trivially with $\ell=0$, since only $\mathcal B=\emptyset$ can minimally encase the empty tuple. If $\ell=0$, then $1=r_1=r_{\ell+1}=t+1$ implies $t=0$. Therefore we may assume $t,\,\ell\geq 1$.
Suppose the virtual independent  sets $\mathcal C_i$ exist with the prescribed properties 2(a), 2(b) and 2(c). Since $\mathcal C_j\subseteq \mathcal X_1\cup \ldots\cup \mathcal X_s$ for $j<\ell$, it follows that each $\mathcal C_j$ with $j<\ell$ is not just a virtual independent set, but also a support set. We must show that $\mathcal B$ minimally encases $-(u_1,\ldots,u_t)$ and that the conditions in Item 3 all hold.

If there were some $\mathbf v_i\in \mathcal C_j$ with $\mathcal X_i^{\pi_{j-1}}=\emptyset$, for some $i\in [1,s]$ and $j\in [1,\ell]$,  then Proposition \ref{prop-orReay-BasicProps}.9 implies $\mathcal X_i\subseteq \darrow \mathcal C_{j-1}\subseteq \darrow \mathcal C_j$, with the latter inclusion in view of $\mathcal C_{j-1}\prec \mathcal C_j$. Thus $\mathcal X_i\cup \{\mathbf v_i\}\subseteq \darrow \mathcal C_j$, ensuring that $\mathcal C_j$ is \emph{not} a support set, which is only possible when $j=\ell$, in which case 2(a) gives us the contradiction $\mathbf v_i\notin \mathcal C_\ell=\mathcal C_j$ by hypothesis. Therefore we instead conclude there is no $\mathbf v_i\in \mathcal C_j$ with $\mathcal X_i^{\pi_{j-1}}=\emptyset$ for any $i\in [1,s]$ and $j\in [1,\ell]$. 

Suppose $\pi_{j-1}(u_{r_j})=0$ for some $j\in [1,\ell]$. Then, since $\mathcal D_j$ minimally encases $-\pi_{j-1}(u_{r_j})=0$, it follows that   $\mathcal D_j=\emptyset$, whence 2(a) or 2(b)  implies $\mathcal C_j^{\pi_{j-1}}=\emptyset$, i.e., $\pi_{j-1}(\x)=\{0\}$ for all $\x\in \mathcal C_j$. Combining this with the conclusion of the previous paragraph and Proposition \ref{prop-orReay-modulo}.1, we find $\mathcal C_j\subseteq \mathcal X_1\cup\ldots\cup \mathcal X_s$, and now Proposition \ref{prop-orReay-BasicProps}.9 implies  $\mathcal C_j\subseteq \darrow \mathcal C_{j-1}$. Thus  $\mathcal C_j=\mathcal C_j^*\preceq \mathcal C_{j-1}$ as $\mathcal C_j$ is virtual independent, which contradicts that $\mathcal C_{j-1}\prec \mathcal C_j$.
 So we instead conclude that $\pi_{j-1}(u_{r_j})\neq 0$ for all $j\in [1,\ell]$, whence  2(c) ensures that $\mathcal F=(\R^\cup \la \mathcal C_1\ra,\ldots,\R^\cup\la \mathcal C_\ell\ra)$ is a compatible filter for $\vec u$ with $\mathcal F(\vec u)=(\overline u_{1},\ldots,\overline u_{\ell})$, where $\overline u_{j}=\pi_{j-1}(u_{r_j})/\|\pi_{j-1}(u_{r_j})\|$ for $j\in [1,\ell]$, as required for Item 3.

 We proceed by induction on $j\in[0,\ell]$ to show that $\mathcal C_j$ minimally encases $-(u_1,\ldots,u_{r_{j+1}-1})$ as well as $-(u_{r_1},\ldots,u_{r_{j}})$ and $-(\overline u_{1},\ldots,\overline u_{j})$, all urbanely,  with $\C^\cup(\mathcal C_j)+\C(u_{r_1},\ldots,u_{r_j})=\C^\cup(\mathcal C_j)+\C(\overline u_{1},\ldots,\overline u_{j})=\R^\cup \la \mathcal C_j\ra$.  The case $j=\ell$ will then verify 3(b) and 3(c), and also show that $\mathcal B$ minimally encases $\vec u$ (as required for Item 2). During the course of the proof, we will also see that 3(a) holds.
 The base case is $j=0$, in which case the empty set $\mathcal C_0=\emptyset$ minimally encases the empty tuple trivially with $\R^\cup \la \emptyset\ra=\{0\}=\C^\cup(\emptyset)$. Therefore we assume $j\geq 1$ and that this has been shown for $j-1$. To lighten notation, let $$\pi=\pi_{j-1}.$$

  Since $\mathcal C_{j-1}\prec \mathcal C_j$ implies $\mathcal C_{j-1}\subseteq \darrow\mathcal C_j$, we have $\C^\cup(\mathcal C_{j-1})\subseteq \C^\cup (\mathcal C_j)$. As a result, since $\mathcal C_{j-1}$ minimally encases $-(u_1,\ldots,u_{r_j-1})$,  $-(u_{r_1},\ldots,u_{r_{j-1}})$ and $-(\overline u_{1},\ldots,\overline u_{j-1})$ by induction hypothesis, it follows that $\mathcal C_{j}$ encases these tuples as well.
  By induction hypothesis, we also have $\ker \pi=\R^\cup \la \mathcal C_{j-1}\ra=\C^\cup(\mathcal C_{j-1})+\C(u_{r_1},\ldots,u_{r_{j-1}})=\C^\cup(\mathcal C_{j-1})+\C(\overline u_{1},\ldots,\overline u_{j-1})$. Consequently, in view of 2(c), to show $\mathcal C_j$ encases
 $-(u_1,\ldots,u_{r_j})$ and  $-(u_{r_1},\ldots,u_{r_{j}})$, we
  we just need to know $-\pi(u_{r_j})\in \C^\cup(\mathcal C_j^\pi)$,
for if this is the case, then $u_{r_j}+x\in   -\C^\cup(\mathcal C_j)$ for some $x\in \R^\cup \la \mathcal C_{j-1}\ra= \C^\cup(\mathcal C_{j-1})+\C(u_{r_1},\ldots,u_{r_{j-1}})$, whence $(-x+\C(u_{r_1},\ldots,u_{r_{j-1}}))\cap  -\C^\cup(\mathcal C_j)\neq \emptyset$, and then  the desired conclusion $((u_{r_j}+x)-x+\C(u_{r_1},\ldots,u_{r_{j-1}}))\cap -\C^\cup (\mathcal C_j)\neq \emptyset$ follows in view of $\C^\cup(\mathcal C_j)$ being a convex cone.
Likewise, to show $\mathcal C_j$ encases  $-(\overline u_{1},\ldots,\overline u_{j})$,
  we just need to know $-\pi(\overline u_{j})\in \C^\cup(\mathcal C_j^\pi)$, and  since $\overline u_{j}=\pi(u_{r_j})/\|\pi(u_{r_j})\|$, this is equivalent to the previous condition  $-\pi(u_{r_j})\in \C^\cup(\mathcal C_j^\pi)$.
However, that $-\pi(u_{r_j})\in \C^\cup(\mathcal C_j^\pi)$ holds follows directly from 2(a) or 2(b) and the definition of $\mathcal D_j$. Indeed, Lemma \ref{lem-orReay-mincase-t=1} ensures that $\mathcal C_j^{\pi}=\mathcal D_j$ is a virtual independent  set which minimally encases $-\pi(u_{r_j})$.
Thus we now know $\mathcal C_j$ encases $-(u_1,\ldots,u_{r_j})$, $-(u_{r_1},\ldots,u_{r_j})$ and $-(\overline u_1,\ldots,\overline u_j)$.


Since $\ker \pi=\R^\cup \la \mathcal C_{j-1}\ra$ and $\mathcal C_j$ contains no $\mathbf v_i$ with $\mathcal X_i^{\pi}=\emptyset$, as shown earlier, it follows in view of $\mathcal C_j^{\pi}=\mathcal D_j$ and Propositions \ref{prop-orReay-BasicProps}.9 and \ref{prop-orReay-modulo}.1 that
$\pi^{-1}(\mathcal D_j)\subseteq  \mathcal C_j\subseteq \darrow (\pi^{-1}(\mathcal D_j)\cup \mathcal C_{j-1})$.
Since $\mathcal C_{j-1}\prec \mathcal C_j$, we have $\darrow \mathcal C_{j-1}\subseteq \darrow \mathcal C_j$, whence $\darrow \mathcal C_j
\subseteq \darrow \pi^{-1}(\mathcal D_j)\cup \darrow \mathcal C_{j-1}\subseteq \darrow \mathcal C_j\cup \darrow \mathcal C_{j-1}\subseteq \darrow \mathcal C_j$, implying \be\label{liftlook}\mathcal C_j=\mathcal C_j^*=(\darrow \pi^{-1}(\mathcal D_j)\cup \darrow \mathcal C_{j-1})^*=(\pi^{-1}(\mathcal D_j)\cup \mathcal C_{j-1})^*,\ee which shows that the virtual independent set $\mathcal C_j$ is the lift of the virtual independent set $\mathcal D_j$ as required for 3(a). Since $\mathcal C_{j-1}$ is a support set, Proposition \ref{prop-orReay-modulo}.5 implies that \be\label{Clineup}\mathcal C_j\setminus \pi^{-1}(\mathcal D_j)=\mathcal C_{j-1}\setminus \darrow \partial(\mathcal \pi^{-1}(\mathcal D_j))=\mathcal C_{j-1}\setminus \darrow \pi^{-1}(\mathcal D_j),\ee where the second equality follows since $\mathcal C_{j-1}$ and $\pi^{-1}(\mathcal D_j)$ are disjoint (as $\pi(\x)=\{0\}$ for $\x\in \mathcal C_{j-1}$ but $\pi(\x)\neq \{0\}$ for $\x\in \pi^{-1}(\mathcal D_j)$).

Now $\pi(\pi^{-1}(\mathcal D_j))=\mathcal D_j$  is a virtual independent  set  with $\pi$ injective on $\pi^{-1}(\mathcal D_j)$ (by Proposition \ref{prop-orReay-modulo}.1) while $\partial(\mathcal D_j)=\partial(\pi(\pi^{-1}(\mathcal D_j)))=\partial(\pi^{-1}(\mathcal D_j))^{\pi}$ by Proposition \ref{prop-orReay-modulo}.2. Thus Proposition
\ref{prop-orReay-BasicProps}.3 implies that $\pi^{-1}( D_j)$ is linearly independent modulo
$\ker \pi+\R^\cup\la \partial(\pi^{-1}(\mathcal D_j))\ra=\R^\cup\la \mathcal C_{j-1}\ra+\R^\cup\la \partial(\pi^{-1}(\mathcal D_j))\ra$, ensuring that \be\label{spillmilk}\C^\cup(\darrow \pi^{-1}(\mathcal D_j))\cap \R^\cup\la \mathcal C_{j-1}\ra=\C^\cup(\pi^{-1}(\mathcal D_j))\cap \R^\cup\la \mathcal C_{j-1}\ra\subseteq \C^\cup\big(\partial(\pi^{-1}(\mathcal D_j))\big).\ee

\subsection*{Claim A} $\mathcal C_j\setminus \pi^{-1}(\mathcal D_j)\cup \partial(\pi^{-1}(\mathcal D_j))$ encases the tuples $-(u_1,\ldots,u_{r_j-1})$, $-(u_{r_1},\ldots,u_{r_{j-1}})$ and $-(\overline u_{1},\ldots,\overline u_{j-1})$, but  $\mathcal C'_{j-1}\cup \partial(\pi^{-1}(\mathcal D_j))$ encases none of these    for any $\mathcal C'_{j-1}\prec\mathcal C_j\setminus \pi^{-1}(\mathcal D_j)$.

\begin{proof}
We know $\mathcal C_j$ encases  $-(u_1,\ldots,u_{r_j-1})$, $-(u_{r_1},\ldots,u_{r_{j-1}})$ and $-(\overline u_{1},\ldots,\overline u_{j-1})$ with $u_i,\,\overline u_{k}\in \R^\cup \la \mathcal C_{j-1}\ra$ for all $i<r_j$ and  $k<j$. In view of \eqref{liftlook} and \eqref{spillmilk}, we have \begin{align*}\C^\cup(\mathcal C_j)\cap \R^\cup \la \mathcal C_{j-1}\ra=&\;\C^\cup(\darrow \pi^{-1}(\mathcal D_j)\cup \mathcal C_{j-1})\cap \R^\cup \la \mathcal  C_{j-1}\ra\\\subseteq &\;\C^\cup(\partial(\pi^{-1}(\mathcal D_j)))+\C^\cup(\mathcal C_{j-1}\setminus\darrow  \pi^{-1}(\mathcal D_j)).\end{align*} Thus $\mathcal C_j\setminus \pi^{-1}(\mathcal D_j)\cup \partial(\pi^{-1}(\mathcal D_j))$ encases the tuples $-(u_1,\ldots,u_{r_j-1})$, $-(u_{r_1},\ldots,u_{r_{j-1}})$ and $-(\overline u_{1},\ldots,\overline u_{j-1})$  in view of \eqref{Clineup}.
Suppose by contradiction that $\mathcal C'_{j-1}\cup \partial(\pi^{-1}(\mathcal D_j))$ encases one of the tuples $-(u_1,\ldots,u_{r_j-1})$, $-(u_{r_1},\ldots,u_{r_{j-1}})$ or $-(\overline u_{1},\ldots,\overline u_{j-1})$ for some $\mathcal C'_{j-1}\prec\mathcal C_j\setminus \pi^{-1}(\mathcal D_j)=\mathcal C_{j-1}\setminus \darrow \pi^{-1}(\mathcal D_j)$ (in view of \eqref{Clineup}).
Note that $\mathcal C_j\setminus \pi^{-1}(\mathcal D_j)\subseteq \mathcal C_j$ is a virtual independent  set as it is a subset of the virtual independent  set $\mathcal C_j$. As a result, $\mathcal C'_{j-1}\prec\mathcal C_{j-1}\setminus \darrow \pi^{-1}(\mathcal D_j)$ implies that  there is some
\be\label{cotag}\y\in \mathcal C_{j-1}\setminus \darrow \pi^{-1}(\mathcal D_j)\quad\mbox{ with }\quad\y\notin\darrow \mathcal C'_{j-1}.\ee
Since both $\mathcal C_{j-1}$ and $\mathcal C'_{j-1}\cup \partial(\pi^{-1}(\mathcal D_j))$ are subsets of $\mathcal X_1\cup\ldots\cup \mathcal X_s$ that encase  $-(u_1,\ldots,u_{r_j-1})$, $-(u_{r_1},\ldots,u_{r_{j-1}})$ or $-(\overline u_{1},\ldots,\overline u_{j-1})$  with the encasement by $\mathcal C_{j-1}$ minimal by induction hypothesis, it follows from the already established case in Item 1 that $\mathcal C_{j-1}\preceq \mathcal C'_{j-1}\cup \partial(\pi^{-1}(\mathcal D_j))$, implying $\y\in \mathcal C_{j-1}\subseteq \darrow \mathcal C'_{j-1}\cup \darrow \partial(\pi^{-1}(\mathcal D_j))\subseteq \darrow \mathcal C'_{j-1}\cup \darrow \pi^{-1}(\mathcal D_j)$.  However, this contradicts \eqref{cotag}, and Claim A is established.
\end{proof}

We showed above that $\mathcal C_j$ encases the tuples
$-(u_1,\ldots,u_{r_j})$, $-(u_{r_1},\ldots,u_{r_{j}})$ and  $-(\overline u_{1},\ldots,\overline u_{j})$.
 Let us next show that is does so \emph{minimally}. To this end, its suffices to show that the immediate predecessor $\mathcal C'_j=\mathcal C_j\setminus \{\x\}\cup \partial(\{\x\})$ encases neither $-(u_1,\ldots,u_{r_j})$ nor $-(u_{r_1},\ldots,u_{r_j})$ nor $-(\overline u_{1},\ldots,\overline u_{j})$ for any $\x\in \mathcal C_j$. Suppose by contradiction that this fails for $\x\in \mathcal C_j$. If $\pi(\x)\neq \{0\}$, then Proposition \ref{prop-orReay-modulo}.3 implies that $(\mathcal C'_j)^{\pi}\prec \mathcal C_j^{\pi}$ with $-\pi(u_{r_j})\in \C^\cup((\mathcal C'_j)^{\pi})$ or $-\pi(\overline u_{r_j})\in \C^\cup((\mathcal C'_j)^\pi)$.
 Noting that $\overline u_{j}=\pi(u_{r_j})/||\pi(u_{r_j})||$, we see that the latter case implies the former, and now both cases  contradict that $\mathcal C_j^{\pi}=\mathcal D_j$ minimally encases the element $-\pi(u_{r_j})$ by definition of $\mathcal D_j$.
So instead suppose $\pi(\x)=\{0\}$, so $\x\in \mathcal C_{j}\setminus \pi^{-1}(\mathcal D_j)$. In this case, $$\mathcal C'_j=\pi^{-1}(\mathcal D_j)\cup \mathcal C'_{j-1}$$ with $\mathcal  C'_{j-1}\prec \mathcal C_j\setminus \pi^{-1}(\mathcal D_j)$. Since $\mathcal C'_j$  encases $-(u_1,\ldots,u_{r_j})$, $-(u_{r_1},\ldots,u_{r_{j}})$ or $-(\overline u_{1},\ldots,\overline u_{j})$, it also encases
$-(u_1,\ldots,u_{r_j-1})$, $-(u_{r_1},\ldots,u_{r_{j-1}})$ or $-(\overline u_{1},\ldots,\overline u_{j-1})$
 with $u_i, \overline u_{r_k}\in \R^\cup \la \mathcal C_{j-1}\ra$ for all $i<r_j$ and $k<j$. As argued in Claim A, it follows in view of \eqref{spillmilk} that $$\C^\cup (\mathcal C'_j)\cap \R^\cup \la \mathcal C_{j-1}\ra=\C^\cup (\pi^{-1}(\mathcal D_j)\cup \mathcal C'_{j-1})\cap \R^\cup \la\mathcal C_{j-1}\ra\subseteq \C^\cup(\partial(\pi^{-1}(\mathcal D_j)))+\C^\cup(\mathcal C'_{j-1}).$$
 Thus $\mathcal C'_{j-1}\cup \partial(\pi^{-1}(\mathcal D_j))$  also encases one of the tuples $-(u_1,\ldots,u_{r_j-1})$, $-(u_{r_1},\ldots,u_{r_{j-1}})$ or $-(\overline u_{1},\ldots,\overline u_{j-1})$, which is contrary to Claim A. This shows $\mathcal C_j$ minimally encases $-(u_1,\ldots,u_{r_j})$, $-(u_{r_1},\ldots,u_{r_{j}})$ and $-(\overline u_{1},\ldots,\overline u_{j})$.

Next, we show that $\R^\cup\la \mathcal C_j\ra=\C^\cup(\mathcal C_j)+\C(u_{r_1},\ldots,u_{r_j})=\C^\cup(\mathcal C_j)+\C(\overline u_{1},\ldots,\overline u_{j})$. By induction hypothesis,   $\R^\cup\la \mathcal C_{j-1}\ra=\C^\cup(\mathcal C_{j-1})+\C(u_{r_1},\ldots,u_{r_{j-1}})\subseteq \C^\cup(\mathcal C_j)+\C(u_{r_1},\ldots,u_{r_j})$ (the inclusion follows as $\mathcal C_{j-1}\prec \mathcal C_j$ implies $\mathcal C_{j-1}\subseteq \darrow \mathcal C_j$), and likewise $\R^\cup(\mathcal C_{j-1})\subseteq \C^\cup (\mathcal C_j)+\C(\overline u_{1},\ldots,\overline u_{j})$. Consequently, it suffices to show $\C^\cup(\mathcal C_j^{\pi})+\C(\pi(u_{r_j}))=\R^\cup \la\mathcal C_{j}^{\pi}\ra$ and $\C^\cup(\mathcal C_j^{\pi})+\C(\pi(\overline u_{j}))=\R^\cup \la\mathcal C_{j}^{\pi}\ra$.
Recalling that $\overline u_{j}=\pi(u_{r_j})/\|\pi(u_{r_j})\|$,  we find that the latter condition is equivalent to the former. As already remarked above, the virtual independent set $\mathcal C_j^{\pi}=\mathcal D_j$ minimally encases the element $-\pi(u_{r_j})$, which is equivalent to
$D_j\cup\{\pi(u_{r_j})\}$ being a minimal positive basis modulo  $\R^\cup \la \partial(\mathcal D_j)\ra$ (as remarked after Lemma \ref{lem-orReay-mincase-t=1}).
Thus there is a strictly positive linear combination $\Summ{\x\in \mathcal D_j}\alpha_\x x+\beta \pi(u_{r_j})=\xi\in \R^\cup \la \partial(\mathcal D_j)\ra$, so $\beta>0$ and $\alpha_\x>0$ for $\x\in \mathcal D_j$. But then $x+\partial(\x)\subseteq \x$ for each $\x\in \mathcal D_j$ ensures that $\R^\cup\la \partial(\mathcal D_j)\ra=\xi+\Summ{\x\in \mathcal D_j}\partial(\x)\subseteq \C^\cup(\mathcal D_j)+\C(\pi(u_{r_j}))$, and now $D_j\cup\{\pi(u_{r_j})\}$ being a minimal positive basis modulo  $\R^\cup \la \partial(\mathcal D_j)\ra$  ensures that $\C^\cup(\mathcal C_j^\pi)+\C(\pi(u_{r_j}))=\C^\cup(\mathcal D_j)+\C(\pi(u_{r_j}))=\R^\cup\la \mathcal D_j\ra=\R^\cup \la\mathcal C_{j}^{\pi}\ra$. This establishes $\R^\cup\la \mathcal C_j\ra=\C^\cup(\mathcal C_j)+\C(u_{r_1},\ldots,u_{r_j})=\C^\cup(\mathcal C_j)+\C(\overline u_{1},\ldots,\overline u_{j})$.

By 2(c) and the conclusion of the previous paragraph,  $-u_i\in \R^\cup\la \mathcal C_j\ra=\C^\cup(\mathcal C_j)+\C(u_{r_1},\ldots,u_{r_j})$ for $i<r_{j+1}$, meaning $(u_i+\C(u_{r_1},\ldots,u_{r_j}))\cap -\C^\cup(\mathcal C_j)\neq \emptyset$.
In particular, this is true for $i\in [r_j+1,r_{j+1}-1]$, which implies $\mathcal C_j$ encases not just $-(u_1,\ldots,u_{r_j})$ but also $-(u_1,\ldots,u_{r_{j+1}-1})$. As a result, since  we already know $\mathcal C_j$ minimally encases $-(u_1,\ldots,u_{r_j})$, we conclude that  $\mathcal C_j$ minimally encases $-(u_1,\ldots,u_{r_{j+1}-1})$.  Since  $\mathcal C_{j-1}$ minimally encases $-(u_1,\ldots,u_{r_j-1})$ by induction hypothesis with $\mathcal C_{j-1}\prec \mathcal C_j$, we cannot have $\mathcal C_{j}$ \emph{minimally} encasing  $-(u_1,\ldots,u_{r_j-1})$. Thus $t'_j=r_j-1\in [0,r_{j+1}-1]$ is the maximal index such that $\mathcal C_j$ does not minimally encase $-(u_1,\ldots,u_{t'_j})$. Consequently, since $\mathcal C_{j-1}\subseteq \mathcal X_1\cup \ldots\cup \mathcal X_s$, we see that $\mathcal C_j$ minimally encases $-(u_1,\ldots,u_{r_{j+1}-1})$ urbanely. Likewise,  the minimal encasement of $-(u_{r_1},\ldots,u_{r_{j}})$ and $-(\overline u_{1},\ldots,\overline u_{j})$ must also be urbane, and the induction is complete. As already noted, this  completes one of the implications in Item 2 and all parts of Item 3

\smallskip

To prove the other implication in Item 2, now assume that $\mathcal B$ minimally encases $-(u_1,\ldots,u_{t})$ urbanely, which ensures that  $\mathcal B$ is a  virtual independent set. Let $t'\in [0,t-1]$ be the maximal index such that $\mathcal B$ does not minimally encase $-(u_1,\ldots,u_{t'})$, and let $\mathcal A\prec \mathcal B$ be a subset which minimally encases $-(u_1,\ldots,u_{t'})$. Since $\mathcal B$ minimally encases $-\vec u$ urbanely, we have $\mathcal A\subseteq \mathcal X_1\cup \ldots\cup\mathcal X_s$. In view of the already established portion of Item 1, it follows that $\mathcal A$ is the unique subset of $\mathcal X_1\cup\ldots\cup \mathcal X_s$ which minimally encases $-(u_1,\ldots,u_{t'})$.

We first construct the support sets $\mathcal C_i$ and indices $r_i$ satisfying  2(b) and 2(c) for the set  $\mathcal A$ recursively. We will then show $\mathcal C_{\ell-1}=\mathcal A$ and $t'=r_\ell-1$ with $\mathcal C_\ell=\mathcal B$ also satisfying 2(a) and 2(c) afterwards to complete the proof of Item 2.
Suppose the sets $\mathcal C_i\subseteq \mathcal X_1\cup \ldots\cup \mathcal X_s$  have already been constructed for $i=0,1\ldots,j-1$, where $j\geq 1$ (we set $r_0=0$ and $\mathcal C_0=\emptyset$). If $\mathcal C_{j-1}=\mathcal A$, we are done with the initial construction, so assume otherwise.
Let $r_{j}\in [r_{j-1}+1,t']$ be the minimal index such that $u_{r_j}\notin \R^\cup\la \mathcal C_{j-1}\ra$, or set $r_j=t'+1$ if no such index exists. In view of the already completed implication in Item 2, we see that $\mathcal C_{j-1}$ minimally encases $-(u_1,\ldots,u_{r_j-1})$. If $\mathcal A$ also minimally encases $-(u_1,\ldots,u_{r_j-1})$, then both $\mathcal C_{j-1}$ and $\mathcal A$ minimally encase  $-(u_1,\ldots,u_{r_j-1})$, in which case the already completed portion of Item 1 implies that $\mathcal C_{j-1}=\mathcal A$, contrary to assumption.
Therefore we may assume $\mathcal A$ does not minimally encases $-(u_1,\ldots,u_{r_j-1})$. In particular, $r_j\leq t'$. Thus, since $\mathcal C_{j-1}\subseteq \mathcal X_1\cup \ldots\cup \mathcal X_s$ minimally encases $-(u_1,\ldots,u_{r_j-1})$ and $\mathcal A\subseteq \mathcal X_1\cup \ldots\cup \mathcal X_s$ encases $-(u_1,\ldots,u_{r_j-1})$, it again follows from the already completed portion of Item 1 that  $\mathcal C_{j-1}\prec \mathcal A$.
As before, let $$\pi=\pi_{j-1}$$ to simplify notation.
By Proposition \ref{prop-orReay-modulo}, $\pi(\mathcal R)$ is an oriented Reay system with $\mathcal A^{\pi}\subseteq \bigcup_{i=1}^{s}\mathcal X_i^{\pi}$ a support set. Since $\mathcal A$ encases $-(u_1,\ldots,u_{t'})$ but $u_{r_j}\notin\R^\cup\la \mathcal C_{j-1}\ra=\ker\pi$ and $u_1,\ldots,u_{r_j-1}\in \ker \pi$ (by induction hypothesis), it follows that $\mathcal A^\pi$
encases $-\pi(u_{r_j})$.
Thus there is some $\mathcal D_j\preceq \mathcal A^\pi$ that minimally encases $-\pi(u_{r_j})$, and by  Proposition \ref{prop-orReay-modulo}.3 we have  $\pi^{-1}(\mathcal D_j)\preceq \pi^{-1}(\mathcal A^\pi)\subseteq \mathcal A\subseteq \mathcal X_1\cup\ldots\cup \mathcal X_s$.
Since $\mathcal D_j\preceq \mathcal A^\pi\subseteq\bigcup_{j\in J}\mathcal X_j^\pi$ minimally encases $-\pi(u_{r_j})$, it follows by Lemma \ref{lem-orReay-mincase-t=1} that $\mathcal D_j$ is a support set with $-\pi(u_{r_j})\in \C^\cup (\mathcal D_j)^\circ$, whence $\mathcal D_j=\supp_{\pi(\mathcal R)}(-\pi(u_{r_j}))$.
 By Proposition \ref{prop-orReay-modulo}.5, $$\mathcal C_j:=(\pi^{-1}(\mathcal D_j)\cup \mathcal C_{j-1})^*=\pi^{-1}(\mathcal D_j)\cup \mathcal C'_{j-1}$$ is a support set for $\mathcal R$, where $\mathcal C'_{j-1}=\mathcal C_{j-1}\setminus \darrow \partial(\pi^{-1}(\mathcal D_j))$. Since $\mathcal C_{j-1}\subseteq \darrow (\pi^{-1}(\mathcal D_j)\cup \mathcal C_{j-1})=\darrow \mathcal C_j$, we have $\mathcal C_{j-1}=\mathcal C_{j-1}^*\preceq \mathcal C_j$. Proposition \ref{prop-orReay-modulo}.4 implies
 $\mathcal C_{j}^{\pi_{j-1}}=\mathcal C_j^\pi=\mathcal D_j$, whence  2(b) holds for $\mathcal C_j$. Also, since $\mathcal D_j$ is a nonempty set of nonzero elements while $\mathcal C_{j-1}^\pi=\emptyset$, we must have $\mathcal C_{j-1}\prec \mathcal C_j$. Letting $r_{j+1}\in [r_j+1,\ldots,t']$ be the minimal index such that $u_{r_{j+1}}\notin \R^\cup\la \mathcal C_j\ra$, or setting $r_{j+1}=t'+1$ if no such index exists, we see that 2(c) also holds. This defines the support sets $\mathcal C_i$ and indices $r_i$ for the set $\mathcal A$.
Since we cannot have an infinite ascending chain $\emptyset=\mathcal C_0\prec \mathcal C_1\prec \mathcal C_2\prec\ldots$ of subsets from the finite set $\mathcal X_1\cup\ldots\cup\mathcal X_s$, the process must eventually terminate with some index $r_{\ell}=t'+1$ with $\mathcal C_{\ell-1}=\mathcal A$. The remainder of the proof of Item 2 is similar to what we have just seen, with some important but subtle differences. We now set $j=\ell$, so $$\pi=\pi_{\ell-1} \quad\und \quad \ker \pi=\R^\cup\la \mathcal C_{\ell-1}\ra=\R^\cup \la \mathcal A\ra.$$

By definition of $t'$, $\mathcal B$ minimally encases $-(u_1,\ldots,u_{t'+1})=-(u_1,\ldots,u_{r_\ell})$.  Since $\mathcal A\prec \mathcal B$, this ensures that $\mathcal A$ does not minimally encase $-(u_1,\ldots,u_{r_\ell})$, and thus $\pi(u_{r_\ell})\neq 0$. In consequence, $r_\ell$ is the minimal index such that $u_{r_\ell}\notin \R^\cup\la \mathcal C_{\ell-1}\ra=\R^\cup \la \mathcal A\ra$.
In view of  Proposition \ref{prop-orReay-modulo}.1,
$\pi(\mathcal R)=
(\mathcal X_j^{\pi}\cup \{\pi(\mathbf v_j)\})_{j\in J}$ is an oriented Reay system.
 Since $\mathcal B$ is a virtual independent set with $\mathcal C_{\ell-1}=\mathcal A\prec \mathcal B$, and since $\mathcal B\subseteq \mathcal X_1\cup \ldots\cup \mathcal X_s\cup \{\mathbf v_j:\;j\in J\}$ per definition of urbane encasement,  Proposition \ref{prop-orReay-modulo}.6 implies that  $\mathcal B^{\pi}$ is a  virtual independent  set.

Since $\mathcal B$ encases $-(u_1,\ldots,u_{t})$ but $u_{r_\ell}\notin\R^\cup\la \mathcal C_{\ell-1}\ra=\ker\pi$ and $u_1,\ldots,u_{r_\ell-1}\in \ker \pi$ (since $\mathcal A=\mathcal C_{\ell-1}$ encases $-(u_1,\ldots,u_{r_\ell-1})$), it follows that  $\mathcal B^{\pi}$
encases $-\pi(u_{r_\ell})$.
Thus there is some $\mathcal D_\ell\preceq \mathcal B^\pi$ which minimally encases $-\pi(u_{r_\ell})$, and by Proposition \ref{prop-orReay-modulo}.3, we have $$\pi^{-1}(\mathcal D_\ell)\preceq \pi^{-1}(\mathcal B^\pi)\subseteq \mathcal B.$$
Since $\mathcal B^{\pi}$ is a virtual independent  set, and since $(\mathcal D_\ell)^*=\mathcal D_\ell$ holds by virtue of $\mathcal D_\ell$ minimally encasing $-\pi(u_{r_\ell})$,  it follows that  $\mathcal D_\ell\preceq \mathcal B^{\pi}$ is also a virtual independent  set, and one which minimally encases the element $-\pi(u_{r_\ell})$.
 By Proposition \ref{prop-orReay-modulo}.5, $$\mathcal C_\ell:=(\pi^{-1}(\mathcal D_\ell)\cup \mathcal C_{\ell-1})^*=\pi^{-1}(\mathcal D_\ell)\cup \mathcal C'_{\ell-1}$$ is a virtual independent set for $\mathcal R$, where $\mathcal C'_{\ell-1}=\mathcal C_{\ell-1}\setminus \darrow \partial(\pi^{-1}(\mathcal D_\ell))$. In view of  $\mathcal C_{\ell-1}\subseteq \darrow (\pi^{-1}(\mathcal D_\ell)\cup\mathcal C_{\ell-1})=\darrow \mathcal C_\ell$, we have $\mathcal C_{\ell-1}=\mathcal C_{\ell-1}^*\preceq \mathcal C_\ell$.
Proposition \ref{prop-orReay-modulo}.4 implies $\mathcal C_{\ell}^{\pi}=\mathcal D_\ell$, whence  2(a) holds for $\mathcal C_\ell=(\pi^{-1}(\mathcal D_\ell)\cup \mathcal C_{\ell-1})^*$. Note $\pi^{-1}(\mathcal D_\ell)$ contains no $\mathbf v_i$ with $\pi(\mathbf v_i)=\{0\}$ by definition, while this is also the case for $\mathcal C_{\ell-1}\subseteq \mathcal X_1\cup\ldots\cup \mathcal X_s$. Since $\mathcal D_\ell$ is a nonempty set of nonzero elements while $\mathcal C_{\ell-1}^\pi=\emptyset $, we must have $\mathcal C_{\ell-1}\prec \mathcal C_\ell$. Letting $r_{\ell+1}\in [r_\ell+1,\ldots,t]$ be the minimal index such that $u_{r_{\ell+1}}\notin \R^\cup\la \mathcal C_\ell\ra$, or setting $r_{\ell+1}=t+1$ if no such index exists, we see that 2(c) also holds.
It remains to show $\mathcal C_\ell=\mathcal B$. Since $\pi^{-1}(\mathcal D_\ell)\preceq \pi^{-1}(\mathcal B^\pi)\preceq \mathcal B$ and $\mathcal C_{\ell-1}=\mathcal A\prec \mathcal B$, we have $\mathcal C_\ell\subseteq \darrow \mathcal B$, which combined with $\mathcal C_\ell^*=\mathcal C_\ell$ ensures that $\mathcal C_\ell\preceq \mathcal B$. By the already completed direction in  Item 2, $\mathcal C_{\ell}$ minimally encases $-(u_1,\ldots,u_{t'+1})$. Thus, since $\mathcal C_\ell\preceq\mathcal B$, and since $\mathcal B$ \emph{minimally} encases $-(u_1,\ldots,u_{t'+1})$, it follows that $\mathcal C_\ell=\mathcal B$, completing the reverse implication in  Item 2.

It remains only to prove the moreover part of Item 2. If $\mathcal B=\mathcal C_\ell$ is a support set, then Proposition \ref{prop-orReay-modulo}.6 implies that $\mathcal C_\ell^{\pi_{\ell-1}}=\mathcal D_\ell$ is a support set, and one which minimally encases $-\pi_{\ell-1}(u_{r_\ell})$ (in view of 2(a)), forcing $\mathcal D_\ell=\supp_{\pi_{\ell-1}(\mathcal R)}(-\pi_{\ell-1}(u_{r_\ell}))$ by definition of $\supp_{\pi_{\ell-1}(\mathcal R)}$. Conversely, if $\mathcal D_\ell= \supp_{\pi_{\ell-1}(\mathcal R)}(-\pi_{\ell-1}(u_{r_\ell}))$, then $\mathcal D_\ell$ is a support set which minimally encases $-\pi_{\ell-1}(u_{r_\ell})$, in which case Proposition \ref{prop-orReay-modulo}.5 implies that $(\pi_{\ell-1}^{-1}(\mathcal D_\ell)\cup \mathcal C_{\ell-1})^*$ is also a support set. However, $\mathcal B=\mathcal C_\ell=(\pi_{\ell-1}^{-1}(\mathcal D_\ell)\cup \mathcal C_{\ell-1})^*$ was shown to be the lift of $\mathcal D_\ell$ during the proof (cf.  Item 3(a)), so we conclude that $\mathcal B$ is support set, completing the proof of Item 2.

1. Next, we establish Item 1 in the unrestricted case when $\mathcal B,\,\mathcal C\subseteq \mathcal X_1\cup \mathcal \{\mathbf v_1\}\cup\ldots\cup \mathcal X_s\cup \{\mathbf v_s\}$ are both \emph{support sets} that minimally encase $\vec u$ urbanely. If $t=0$, then $\mathcal B=\mathcal C=\emptyset$, as desired, so we may assume $t\geq 1$, ensuring that $\mathcal B$ and $\mathcal C$ are both nonempty.
%
%
Let $\emptyset=\mathcal C_0\prec \mathcal C_1\prec \ldots\prec \mathcal C_{\ell_B-1}\prec \mathcal C_{\ell_B}=\mathcal B$ be the support sets given by application of Item 2 to $\mathcal B$. Observe that $\mathcal C_1=\supp_{\mathcal R}(-u_1)$ depends only on $\vec u$ and not on $\mathcal B$, and thus by an iterative argument (using 2(a)--2(c) and 3(a)), none of the sets $\mathcal C_i$, for $i\in [0,\ell_B]$, depend on $\mathcal B$ at all
 (note the case $i=\ell_B$ requires the moreover statement in Item 1, which is available as $\mathcal B$ and $\mathcal C$ are support sets by hypothesis), meaning the only portion of Item 2 that is dependent on $\mathcal B$, and not $\vec u$, is the number of iterations $\ell_B$ that occur for $\mathcal B$. Applying Item 2 to $\mathcal C$, we arrive at the same conclusion. Thus, letting $\emptyset=\mathcal C'_0\prec \mathcal C'_1\prec \ldots\prec \mathcal C'_{\ell_C-1}\prec \mathcal C'_{\ell_C}=\mathcal C$  be resulting support sets, and w.l.o.g. assuming $\ell_C\leq \ell_B$, we find that $\mathcal C'_i=\mathcal C_i$ for all $i\in [0,\ell_C]$. If $\ell_B=\ell_C$, then $\mathcal C=\mathcal C'_{\ell_C}=\mathcal C_{\ell_B}=\mathcal B$, as desired. Otherwise, $\ell_B>\ell_C$, in which case $\mathcal C=\mathcal C'_{\ell_C}=\mathcal C_{\ell_C}\prec \mathcal C_{\ell_B}=\mathcal B$. However, since both $\mathcal C$ and $\mathcal B$ \emph{minimally} encase $-\vec u$, this is not possible, completing the proof of Item 1.

4. Since any virtual independent subset $\mathcal B\subseteq \mathcal X_1\cup\ldots\cup \mathcal X_s$ must be a support set, the hypotheses of Item 4 imply that \be\label{bequark}\mbox{$\mathcal B\subseteq \mathcal X_1\cup\ldots\cup \mathcal X_s$ is a support set}\quad\mbox{ or } \quad \mbox{
 $\mathcal A\subseteq  \darrow \mathcal B$}.\ee They also imply that  \be\label{nobadv}\mbox{each $\mathcal C_j$, for $j\in [1,\ell]$,  contains no $\mathbf v_i$ with $\mathcal X_i^{\tau_{\ell-1}}=\emptyset$}.\ee Note \eqref{nobadv} is trivially true when $\mathcal C_j\subseteq \mathcal X_1\cup\ldots\cup \mathcal X_s$, and thus for  $j<\ell$, while it holds directly by hypothesis for $\mathcal C_\ell=\mathcal B$ except when $\mathcal A\subseteq \darrow \mathcal B$ with $\mathcal B$ a support set.
 However, in this last remaining case in question,  $\mathcal X_j^{\tau_{\ell-1}}=\emptyset$ implies $\mathcal X_j\subseteq \darrow \mathcal A\cup \darrow \mathcal C_{\ell-1}\subseteq \darrow \mathcal C_\ell=\darrow \mathcal B$ in view of  Proposition \ref{prop-orReay-BasicProps}.9 and $\ker \tau_{\ell-1}=\R^\cup \la \mathcal A\cup \mathcal C_{\ell-1}\ra$ (note $\darrow \mathcal C_{\ell-1}\subseteq \darrow \mathcal C_\ell$ follows from $\mathcal C_{\ell-1}\prec \mathcal C_\ell$), so that the definition of support set instead ensures $\mathbf v_j\notin \darrow \mathcal B=\darrow \mathcal C_\ell$.
 Thus \eqref{nobadv} is established in all cases, which together with \eqref{bequark} allows us to  apply Proposition \ref{prop-orReay-modulo}.6 to each $\mathcal C_j$, for $j\in [1,\ell]$, to conclude  $\mathcal C_j^\pi$ is a support set for $j<\ell$, and thus virtual independent, and that $\mathcal C_\ell^\pi=\mathcal B^\pi$ is also virtual independent (note $\ker \pi\leq \ker \tau_{\ell-1}$).
 Proposition \ref{prop-orReay-modulo}.3 implies that $\mathcal C_0^\pi\preceq\mathcal C_1^\pi\preceq\ldots\preceq \mathcal C_\ell^\pi=\mathcal B^\pi$.

 If $i,\,j\in [0,\ell]$ with $i<j$ and $\mathcal C_i^\pi\prec \mathcal C_j^\pi$, then there must be some $\y\in \mathcal C_j$ with $\pi(\y)\neq \{0\}$ and $\y\notin \darrow \mathcal C_i$, for otherwise $\mathcal C_j^\pi\subseteq (\darrow \mathcal C_i)^\pi=\darrow \mathcal C_i^\pi$ (by Proposition \ref{prop-orReay-modulo}.2), yielding the contradiction $\mathcal C_j^\pi=(\mathcal C_j^\pi)^*\preceq \mathcal C_i^\pi$ (as $\mathcal C_j^\pi$ is virtual independent). Suppose $\y\subseteq \mathcal E+\R^\cup \la \mathcal C_i\ra=\R^\cup \la\mathcal A\cup \mathcal C_i\ra$. If $\mathcal A\subseteq \darrow \mathcal B$, then $y\in \y\subseteq \R^\cup \la\mathcal A\cup \mathcal C_i\ra=\R\la \darrow A\cup \darrow C_i\ra$
meaning $y$ can be written as a linear combination of the elements from $\darrow A\cup \darrow C_i$.
However, it follows in view of $\mathcal A\subseteq \darrow \mathcal B$, \ $\y\in \mathcal C_j$ and $\mathcal C_i\prec \mathcal C_j\preceq \mathcal C_\ell=\mathcal B$ that $\{y\}\cup \darrow A\cup \darrow C_i\subseteq \darrow B$ is a linearly independent subset (as $\mathcal B$ is a virtual independent set), meaning the only way $y$ can be written as a linear combination of the elements from $\darrow A\cup \darrow C_i$ is if $y\in \darrow A\cup \darrow C_i$, i.e., if $\y\in \darrow \mathcal A\cup \darrow \mathcal C_i$.
However, in view of $\pi(\y)\neq \{0\}$ and $\y\notin \darrow \mathcal C_i$, neither of these is possible. On the other hand, if $\mathcal A\nsubseteq \darrow \mathcal B$, then \eqref{bequark} ensures  $\mathcal B\subseteq \mathcal X_1\cup\ldots\cup X_s$ is a support set and $\y\in \mathcal C_j\subseteq \darrow \mathcal B\subseteq \mathcal X_1\cup \ldots\cup \mathcal X_s$. In this case, $\y\subseteq \R^\cup \la\mathcal A\cup \mathcal C_i\ra$ combined with Proposition \ref{prop-orReay-BasicProps}.9 implies that $\y\in \darrow \mathcal A\cup \darrow \mathcal
C_i$, and we obtain the same contradiction as before.
  So we instead conclude that any  $\y\in \mathcal C_j$ with $\pi(\y)\neq \{0\}$ and $\y\notin \darrow \mathcal C_i$ must satisfy $\y\nsubseteq \mathcal E+\R^\cup \la \mathcal C_i\ra$.

By hypothesis, $\mathcal F^\pi=(\R^\cup\la \mathcal C_i^\pi\ra)_{i\in J}$ with $J\subseteq [1,\ell]$ the subset of indices $j\in [1,\ell]$ with $\mathcal C_{j-1}^\pi\prec \mathcal C_{j}^\pi$.
  If $j_1,\,j_2\in \{0\}\cup J$ are consecutive elements with $j_1<j_2$, then $\mathcal C_{j_1}^\pi=\mathcal C_i^\pi$ for all $i\in [j_1,j_2-1]$  and $\mathcal C_{j_1}^\pi\prec \mathcal C_{j_2}^\pi$.
Thus $\mathcal E+\R^\cup \la \mathcal C_{j_1}\ra=\mathcal E+\R^\cup \la \mathcal C_{i}\ra$ for $i\in [j_1,j_2-1]$.
Since $\mathcal C_{j_1}^\pi\prec \mathcal C_{j_2}^\pi$, there must be some $\y\in \mathcal C_{j_2}$ with $\pi(\y)\neq \{0\}$ and $\y\notin \darrow \mathcal C_{j_1}$, implying $\y\nsubseteq \mathcal E+\R^\cup \la \mathcal C_{j_1}\ra$ as shown above. Thus $\mathcal E+\R^\cup \la \mathcal C_{j_1}\ra\subset \mathcal E+\R^\cup \la \mathcal C_{j_2}\ra$, and
in particular, for any $j\in J$, there must be some $\y\in\mathcal C_j$ with $\y\nsubseteq \mathcal E+\R^\cup \la \mathcal C_{j-1}\ra$.

Let $j\in J$ be arbitrary, let $\y\in\mathcal C_j$ with $\y\nsubseteq \mathcal E+\R^\cup \la \mathcal C_{j-1}\ra$ be arbitrary (at least one such $\y$ exists as just noted), and let  $$\tau:\R^d\rightarrow \R^\cup\la \mathcal A\cup \mathcal C_{j-1}\cup\partial(\mathcal C_j)\ra^\bot$$ be the orthogonal projection.
 By Proposition \ref{prop-orReay-modulo}.1, we know that $\tau(\mathcal R)$ is an oriented Reay system. Moreover, $\mathcal C_j$ is a virtual independent  set and either  $\mathcal C_j\subseteq \mathcal X_1\cup\ldots\cup\mathcal X_s$ or else  $j=\ell$ with $\mathcal C_\ell=\mathcal B$ and $\mathcal A\subseteq \darrow \mathcal B$. In the former case, Proposition \ref{prop-orReay-modulo}.6 ensures that $\mathcal C_j^\tau$ is a support set, and  thus $\tau(C_j)\setminus \{0\}$ is a linearly independent set of size $|C_j\setminus \ker \tau|$.
 In the latter case, $\mathcal C_j=\mathcal C_\ell$, and Proposition \ref{prop-orReay-BasicProps}.1 implies $\ker \tau=\R^\cup \la \mathcal A\cup \mathcal C_{\ell-1}\cup \partial(\mathcal B)\ra=\R\la \darrow A\cup \darrow C_{\ell-1}\cup \darrow \partial(B)\ra$  with $\darrow A\cup \darrow C_{\ell-1}\cup \darrow \partial(B)\subseteq \darrow B$.
  Thus, since $\darrow B$ is linearly independent (as $\mathcal B$ is virtual independent), it follows that  $\tau(\darrow B)\setminus \{0\}$ is a  linearly independent set of size $|\darrow B\setminus \ker \tau|$, and since  $B\subseteq \darrow B$, we then conclude that   $\tau(B)\setminus \{0\}=\tau(C_\ell)\setminus \{0\}$ is a linearly independent set of size $|C_\ell\setminus \ker \tau|$. In both cases, $\tau(C_j)\setminus \{0\}$ is a linearly independent set of size $|C_j\setminus \ker \tau|$.

Now the virtual independent set $\mathcal C_j^{\pi_{j-1}}=\mathcal D_j$ minimally encases $-\pi_{j-1}(u_{r_j})$ by 2(a) or 2(b).  Thus, since we also have $\ker \pi_{j-1}\subseteq \R^\cup\la \mathcal C_{j-1}\cup  \partial(\mathcal C_j)\ra\subseteq \ker \tau$, it follows from Lemma \ref{lem-orReay-mincase-t=1} that
 there exists a strictly positive linear combination of the elements from $C_j\setminus \ker \tau$ equal to $-u_{r_j}+\xi$ for some $\xi\in \ker \tau$.

Let $\mathcal C'_j=\mathcal C_j\setminus\{\y\}\cup \partial(\{\y\})$.
   Since $\y\nsubseteq \mathcal E+\R^\cup \la \mathcal C_{j-1}\ra=\R^\cup \la \mathcal A\cup \mathcal C_{j-1}\ra$, we have $\y\notin \darrow \mathcal A\cup \darrow \mathcal C_{j-1}$. In particular, $\y\notin \mathcal C_{j-1}$, ensuring that $\mathcal C_{j-1}\subseteq \darrow \mathcal C'_j$ (as $\mathcal C_{j-1}\prec \mathcal C_j$ implies $\mathcal C_{j-1}\subseteq \darrow \mathcal C_j$). Since $\mathcal \y\in \mathcal C_j=\mathcal C_j^*$ (as $\mathcal C_j$ is a virtual independent  set), we also have $\y\notin \darrow \mathcal C'_j$. In summary, \be\label{trollop} \y\notin \darrow \mathcal A\cup \darrow \mathcal C_{j-1}\cup \darrow \mathcal C'_j\quad\und\quad \mathcal C_{j-1}\subseteq \darrow \mathcal C'_j.\ee

 Suppose $\y\subseteq \mathcal E+\R^\cup \la \mathcal C'_j\ra=\R^\cup \la \mathcal A\cup \mathcal C'_j\ra$. Consequently, if $\mathcal C_j\subseteq \mathcal X_1\cup \ldots\cup \mathcal X_s$, then Proposition \ref{prop-orReay-BasicProps}.9 enures that $\y\in \darrow \mathcal A\cup \darrow \mathcal C'_j$, contradicting \eqref{trollop}.
 On the other hand, if  $\mathcal C_j\nsubseteq \mathcal X_1\cup \ldots\cup \mathcal X_s$, then $j=\ell$, \ $\mathcal C_j=\mathcal B$ and $\mathcal A\subseteq \darrow \mathcal B$.
 In this case, $y\in \y\subseteq \R^\cup \la \mathcal A\cup \mathcal C'_j\ra=\R\la \darrow A\cup \darrow C'_j\ra$ with $\{y\}\cup\darrow A\cup \darrow C'_j\subseteq \darrow B$.
 However, since $\mathcal B=\mathcal C_\ell=\mathcal C_j$ is a virtual independent set, it follows that $\darrow B$, and thus also $\{y\}\cup\darrow A\cup \darrow C'_j$, is a linearly independent set, contradicting that $y\in  \R\la \darrow A\cup \darrow C'_j\ra$  with  $y\notin \darrow A\cup \darrow C'_j$ (by \eqref{trollop}).
 So we instead conclude that
\be\label{dustrust} \y\nsubseteq \mathcal E+\R^\cup \la \mathcal C'_j\ra.\ee Since $\ker \tau=\mathcal E+\R^\cup \la \mathcal C_{j-1}\cup\partial(\mathcal C_j)\ra\subseteq \mathcal E+\R^\cup\la \mathcal C_{j-1}\cup \mathcal C'_j\ra=\mathcal E+\R^\cup\la \mathcal C'_j\ra$, with the final equality in view of \eqref{trollop},
and since $\tau(C_j)\setminus\{0\}$ is a linearly independent set of size $|C_j\setminus \ker \tau|$,
we conclude from \eqref{dustrust} and $\y\in \mathcal C_j$ that \be\label{lades}\R^\cup \la \mathcal C_j\ra\nsubseteq\mathcal E+\R^\cup \la \mathcal C_{j-1}\cup\partial(\mathcal C_j)\ra\quad\und\quad\tau(C'_j)\setminus \{0\}\subset\tau(C_j)\setminus \{0\}.\ee

Suppose  $u_{r_j}\in \mathcal E+\R^\cup \la\mathcal C_{j-1}
\cup\partial(\mathcal C_j)\ra=\ker \tau$. As remarked above, $\tau(C_j)\setminus \{0\}$ is a linearly independent set of size $|C_j\setminus \ker \tau|$ and there exists a strictly positive linear combination of the elements from $C_j\setminus \ker \tau$ equal to $-u_{r_j}+\xi$ for some $\xi\in \ker \tau$. However, since $u_{r_j}\in \ker \tau$, this contradicts that $\tau(C_j)\setminus \{0\}$ is a linearly independent set of size $|C_j\setminus \ker \tau|$ unless $C_j\subseteq \ker \tau$, i.e., $\R^\cup \la \mathcal C_j\ra\subseteq \mathcal E+\R^\cup \la\mathcal C_{j-1}\cup\partial(\mathcal C_j)\ra$, which is contrary to \eqref{lades}.
So we instead conclude that \be\label{lamplite}u_{r_j}\notin \mathcal E+\R^\cup \la\mathcal C_{j-1}\cup \partial(\mathcal C_j)\ra\quad\mbox{ for every $j\in J$}.\ee Thus $\tau(u_{r_j})\neq 0$ and, since  $\tau(C_j)\setminus \{0\}$ is linearly independent with a strictly positive linear combination of the elements of $\tau(C_j)\setminus \{0\}$ equal to $-\tau(u_{r_j})$, we see that
   $\tau(C_j)\setminus \{0\}\cup \{\tau(u_{r_j})\}$ is minimal positive basis of size $|C_j\setminus \ker \tau|+1$.


By definition, we have $\pi(\vec u)=(u^*_1,u^*_2,\ldots,u^*_{\ell_\pi})$ with the $u^*_i$ for $i\in [1,\ell_\pi]$ defined as follows. We recursively define the indices $0=s_0<s_1<\ldots<s_{\ell_\pi}<s_{\ell_\pi+1}=r_{\ell+1}=t+1$ by letting $s_i\in [s_{i-1}+1,t]$ be the minimal index such that $\pi^*_{i-1}(u_{s_i})\neq 0$, where $$\pi^*_{i-1}:\R^d\rightarrow (\mathcal E+\R\la u_1,u_2,\ldots,u_{s_{i-1}}\ra)^\bot$$ is the orthogonal projection, and then  $$u^*_i=\pi^*_{i-1}(u_{s_i})/\|\pi^*_{i-1}(u_{s_i})\|.$$ Note this ensures $$\mathcal E+\R\la u_1,u_2,\ldots,u_{s_{i-1}}\ra=\mathcal E+\R\la u_1,u_2,\ldots,u_{s_{i}-1}\ra.$$

 Recall that $j\in J\subseteq [1,\ell]$ is arbitrary. By 2(c), we have $u_i\in \R^\cup \la \mathcal C_{j-1}\ra\subseteq \mathcal E+\R^\cup\la \mathcal C_{j-1}\ra$ for all $i<r_j$.
 Hence, if $s_{i-1}<r_j<s_{i}$ for some $i\in [1,\ell_\pi+1]$, then $\ker \pi_{i-1}^*\subseteq \mathcal E+\R^\cup\la \mathcal C_{j-1}\ra$ and $\pi_{i-1}^*(u_{r_j})=0$, the latter  in view of the minimality in the definition of $s_i$, which contradicts  \eqref{lamplite}. Therefore, we instead conclude that, for each $j\in J$,  $$r_j=s_{j^*}\quad\mbox{ for some $j^*\in [1,\ell_\pi]$}.$$
 For $1\leq i<j^*$, we have $1\leq s_i<s_{j^*}=r_j$, ensuring $u_{s_i}\in \R^\cup \la \mathcal C_{j-1}\ra\subseteq \mathcal E+\R^\cup \la \mathcal C_{j-1}\ra$ and $\mathcal E\subseteq \ker \pi^*_{i-1}=\mathcal E+\R\la u_1,\ldots,u_{s_{i-1}}\ra\subseteq \mathcal E+\R^\cup \la \mathcal C_{j-1}\ra$
  by 2(c) for $\vec u$.
 Thus \be\label{blackwing1} u_i^*\in\R^\cup \la \mathcal C^\pi_{j-1}\ra\quad\mbox{ for all $i<j^*$},\ee where $j\in J$.
  Likewise, $u_k\in \R^\cup \la \mathcal C_{j-1}\ra$ for $k\leq s_{j^*-1}<s_{j^*}=r_j$ and $u_{s_{j*}}=u_{r_j}\in \R^\cup \la \mathcal C_j\ra\subseteq \mathcal E+\R^\cup \la \mathcal C_j\ra$,  while  \eqref{lamplite} implies that $u_{s_{j^*}}=u_{r_j}\notin \mathcal E+\R^\cup \la \mathcal C_{j-1}\ra$. Thus \be\label{merchantike}\mathcal E\subseteq \ker \pi^*_{j^*-1}=\mathcal E+\R\la u_1,\ldots,u_{s_{j^*-1}}\ra\subseteq \ker \tau_{j-1}=\mathcal E+\R^\cup \la \mathcal C_{j-1}\ra\subseteq \mathcal E+\R^\cup \la \mathcal C_j\ra\ee with \be\label{blackwing2}u_{j^*}^*\in \R^\cup \la \mathcal C_j^\pi\ra\setminus \R^\cup \la \mathcal C_{j-1}^\pi\ra\ee and
   $\tau_{j-1}\pi^*_{j^*-1}(u_{s_{j^*}})=\tau_{j-1}(u_{r_j})$, for $j\in J$. As a result,  \eqref{blackwing1} and \eqref{blackwing2} ensure $\mathcal F^\pi$ is  compatible with $\pi(\vec u)$ with $\mathcal F^\pi(\pi(\vec u))=(\overline{u^*_{j}})_{j\in J}$, where each $\overline{u^*_{j}}=\tau_{j-1}(u^*_{j^*})/\|\tau_{j-1}(u^*_{j^*})\|
   =\tau_{j-1}\pi^*_{j^*-1}(u_{s_{j^*}})/
 \|\tau_{j-1}\pi^*_{j^*-1}(u_{s_{j^*}})\|=\tau_{j-1}(u_{r_j})/\|\tau_{j-1}(u_{r_j})\|$.

Let $j\in J$ be arbitrary and let $j_+$ be the next consecutive
 element of $J$ after $j$, or set $j_+=\ell+1$ and $j_+^*=(\ell+1)^*:=\ell_\pi+1$ if $j$ is the final element of $J$.
 Since $j,\,j_+\in J\cup\{\ell+1\}$ are consecutive, the definition of $J$ implies  $\mathcal E+\R^\cup \la \mathcal C_i\ra=\mathcal E+\R^\cup\la \mathcal C_{j^+-1}\ra$ for all $i\in [j,j^+-1]$. In particular,
 \be\label{utype}\mathcal E+\R^\cup \la \mathcal C_j\ra=\mathcal E+\R^\cup\la \mathcal C_{j_+-1}\ra.\ee
 Let   $r_j=s_{j^*}$ and $r_{j_+}=s_{j^*_+}$ with $j^*,\,j^*_+\in [1,\ell_\pi+1]$  as shown above. If $J=\emptyset$, then $\pi(\vec u)$ is the trivial tuple, $\emptyset=\mathcal C_0^\pi=\mathcal C_\ell^\pi=\mathcal B^\pi$, and all parts of Item 4 hold trivially. Therefore we can assume $J$ is nonempty.
We now proceed  to show that 2(a), 2(b) and 2(c) hold for $\pi(\vec u)=(u^*_1,\ldots,u^*_{\ell_\pi})$ using the virtual independent sets $\mathcal C_j^\pi$ for $j\in \{0\}\cup J$, indices $j^*$ for $j\in J\cup \{\ell+1\}$, and elements $$u_{j^*}^*=\pi^*_{j^*-1}(u_{s_{j^*}})/\|\pi^*_{j^*-1}(u_{s_{j^*}})\|
=\pi^*_{j^*-1}(u_{r_j})/\|\pi^*_{j^*-1}(u_{r_j})\|\quad\mbox{ for $j\in J$}$$ in place of the virtual independent sets  $\mathcal C_j$ for $j\in [0,\ell]$, indices $r_j$ for $j\in[1,\ell+1]$, and elements $u_{r_j}$ for $j\in [1,\ell]$, which will imply that $\mathcal B^\pi$ minimally encases $\pi(\vec u)$ urbanely since $\mathcal B^\pi=\mathcal C_\ell^\pi=\mathcal C_j^\pi$ for the final element $j\in J$ (and thereby complete the proof).
Note the $\mathcal C_j^\pi$ for $j\in J$ are virtual independent sets for $\pi(\mathcal R)$ in view of Proposition \ref{prop-orReay-modulo}.6 (as remarked earlier).

If $j\in J$ is the final element, then  $u_{s_i}\in \mathcal E+\R^\cup \la \mathcal C_j\ra=\mathcal E+\R^\cup \la \mathcal C_\ell\ra$ and $\mathcal E\subseteq \ker \pi_{i-1}^*\subseteq
\mathcal E+\R^\cup \la \mathcal C_j\ra=\mathcal E+\R^\cup \la \mathcal C_\ell\ra$
for all $i<\ell_\pi+1=j^*_+$ follows by 2(c) for $\vec u$, in turn implying that  $u^*_i\in \R^\cup \la \mathcal C_j^\pi\ra$ for all $i<j_+^*=\ell_\pi+1$. Otherwise,
\eqref{blackwing1} (applied with $j=j_+$) and \eqref{utype} yield $u^*_i\in \mathcal \R^\cup\la \mathcal C^\pi_{j_+-1}\ra=\R^\cup \la \mathcal C_j^\pi\ra$ for all $i<j_+^*$. As a result,  2(c) holds for $\pi(\vec u)$ for all $j\in J$.

For each $j\in J$, $u^*_{j^*}$ is a positive multiple of $\pi^*_{j*-1}(u_{s_{j^*}})=\pi^*_{j^*-1}(u_{r_{j}})$ with $\ker \pi_{j^*-1}^*\subseteq \ker \tau_{j-1}$ by \eqref{merchantike}, and  thus $\tau_{j-1}(u^*_{j^*})$ is a positive multiple of $\tau_{j-1}\pi^*_{j^*-1}(u_{r_j})=\tau_{j-1}(u_{r_j})$. By \eqref{nobadv},  $\mathcal C_j$ contains no $\mathbf v_i$ with $\mathcal X_i^{\tau_{j-1}}=\emptyset$. Thus  Proposition \ref{prop-orReay-modulo}.6 and \eqref{bequark} ensure that $\mathcal C_j^{\tau_{j-1}}$ is a virtual independent  set,
which will be a support set when $j\in J$ is not the final element (as this ensures $\mathcal C_j\subseteq \mathcal X_1\cup\ldots\cup\mathcal X_s$). Consequently,  to establish 2(a) and 2(b) for $\pi(\vec u)$, we just need to show $\mathcal C_j^{\tau_{j-1}}$ minimally encases $-\tau_{j-1}(u^*_{j^*})$ for each $j\in J$, and since $\tau_{j-1}(u^*_{j^*})$ is a positive multiple of $\tau_{j-1}(u_{r_j})$ for $j\in J$ (as just shown above), this is  equivalent to showing $\mathcal C_j^{\tau_{j-1}}$ minimally encases $-\tau_{j-1}(u_{r_j})$.

Since $-\pi_{j-1}(u_{r_j})\in \C^\cup (\mathcal C_j^{\pi_{j-1}})$ by 2(a) or 2(b) for $\mathcal B$, and since $\ker \pi_{j-1}\leq \ker \tau_{j-1}$, we have $-\tau_{j-1}(u_{r_j})\in \C^\cup(\mathcal C_j^{\tau_{j-1}})$, showing that $\mathcal C_j^{\tau_{j-1}}$ encases $-\tau_{j-1}(u_{r_j})$. To show the encasement is minimal, we need to show $(\mathcal C'_j)^{\tau_{j-1}}$ does not encase $-\tau_{j-1}(u_{r_j})$ for any  $\mathcal C'_j=\mathcal C_j\setminus \{\y\}\cup \partial(\{\y\})$ with $\y\in \mathcal C_j$ and $\y\nsubseteq \mathcal E+\R^\cup\la \mathcal C_{j-1}\ra=\ker \tau_{j-1}$ (cf. Proposition \ref{prop-orReay-modulo}.3: note $\mathcal D\prec \mathcal C_j^{\tau_{j-1}}$ implies $\tau_{j-1}^{-1}(\mathcal D)\prec \tau_{j-1}^{-1}(\mathcal C_j^{\tau_{j-1}})\subseteq \mathcal C_j$).
 Recall that $\tau:\R^d\rightarrow (\mathcal E+\R^\cup\la \mathcal C_{j-1}\cup \partial(\mathcal C_j)\ra)^\bot$ is the orthogonal projection.
 Suppose by contradiction that $-\tau_{j-1}(u_{r_j})\in \C^\cup((\mathcal C'_j)^{\tau_{j-1}})$ for some $\mathcal C'_j=\mathcal C_j\setminus \{\y\}\cup \partial(\{\y\})$ with $\y\in \mathcal C_j$ and $\y\nsubseteq \mathcal E+\R^\cup\la \mathcal C_{j-1}\ra.$
 Consequently, since $\ker \tau_{j-1}\subseteq \ker \tau$, \eqref{lades} ensures $\tau(C'_j)\setminus \{0\}$ is a proper subset of $\tau (C_j)\setminus\{0\}$ that encases $-\tau(u_{r_j})$.  However, as  noted immediately after \eqref{lamplite}, $\tau(C_j)\setminus \{0\}\cup \{\tau(u_{r_j})\}$ is minimal positive basis of size $|C_j\setminus \ker \tau|+1$, meaning it is not possible for a proper subset $\tau(C'_j)\setminus \{0\}$ of $\tau(C_j)\setminus \{0\}$ to encase $-\tau(u_{r_j})$, and with this contradiction, we establish that $\mathcal B^\pi$ minimally encases $\pi(\vec u)$ as described in Item 4, completing the proof.
\end{proof}

\begin{proposition}\label{prop-orientedReay-ExtraUnique}
Let $\mathcal R=(\mathcal X_1\cup \{\mathbf v_1\},\ldots,\mathcal X_s\cup\{\mathbf v_s\})$ be an oriented Reay system in $\R^d$,  let $\vec u=(u_1,\ldots,u_t)$ be a tuple of $t\geq 0$ orthonormal vectors from $\R^d$, and let $\mathcal A,\,\mathcal B\subseteq \mathcal X_1\cup \ldots\cup \mathcal X_s$. Suppose $\mathcal B$ minimally encases $-\vec u$ and $u_1,\ldots,u_t\in \R^\cup \la \mathcal A\ra$. Then $\darrow \mathcal B \subseteq  \darrow \mathcal A$.
\end{proposition}

\begin{proof}Let $\pi:\R^d\rightarrow \R^\cup \la \mathcal A \ra^\bot$ be the orthogonal projection.
Since $\mathcal B\subseteq \mathcal X_1\cup \ldots\cup \mathcal X_s$, it follows that $\mathcal B$ minimally encases $-\vec u$ urbanely. By Proposition \ref{prop-orReay-minecase-char}.4, $\mathcal B^\pi$ minimally encases $-\pi(\vec u)$. However, since $u_1,\ldots,u_t\in \R^\cup \la \mathcal A\ra=\ker \pi$, we have $\pi(\vec u)$ equal to the empty tuple, and  now $\mathcal B^\pi$ minimally encasing the empty tuple forces $\mathcal B^\pi=\emptyset$. Hence $\x\in \ker \pi=\R^\cup \la\mathcal A\ra$ for all $\x\in \mathcal B\subseteq \mathcal X_1\cup \ldots\cup \mathcal X_s$. Consequently, since $\mathcal A\subseteq \mathcal X_1\cup \ldots\cup\mathcal X_s$,   Proposition \ref{prop-orReay-BasicProps}.9 implies $\x\in \darrow \mathcal A$ for all $\x\in \mathcal B$. Thus $\mathcal B\subseteq \darrow \mathcal A$, implying $\darrow \mathcal B\subseteq \darrow \mathcal A$, as desired.
\end{proof}

The next proposition shows that, if a half-space $\x$ occurs in two separate oriented Reay systems $\mathcal R$ and $\mathcal R'$, then the quantity $\partial(\{\x\})$ is the same for both, and is thus intrinsically defined by the half-space $\x$ itself.

The proof of Proposition \ref{prop-orientedReay-halfspace-uniquelydetermines} requires the translation invariant notion for the lineality subspace of a convex cone. Recall that, if $C\subseteq\R^d$ is convex cone, then $C\cap -C$ is the lineality subspace of $C$, which is the maximal subspace contained in $C$. For a general set $X\subseteq \R^d$, we define $$o(X)=\{x\in X:\;x+\R_+(y-x)\subseteq X\mbox{ for every $y\in X$}\}$$ to be the set of apex points of $X$ (see \cite{convexbookI}).
If $C\subseteq \R^d$ is a convex cone containing $0$, then \be\label{lin-apex}o(C)=C\cap -C,\ee as the following argument shows.  If $x\in C\cap -C$ and $y\in C$, then $-x\in C$ and $y\in C$ imply via the convexity of $C$ that $y-x\in C$, whence $x+\R_+(y-x)\subseteq C$ since $x\in C$ with $C$ a convex cone. This shows that $C\cap -C\subseteq o(C)$. On the other hand, if $x\in o(C)$, then since $0\in C$ by hypothesis, it follows from the definition of $o(C)$ that $x+\R_+(0-x)\subseteq C$. In particular, $-x=x-2x\in C$, whence $x\in C\cap -C$ (as $o(C)\subseteq C$ by definition). This establishes the reverse inclusion for \eqref{lin-apex} It is also routine to verify that $$o(z+C)=z+o(C)$$ for any $z\in \R^d$ and convex cone $C\subseteq \R^d$. Combined with \eqref{lin-apex}, it follows that \be\label{apex-lin}o(z+C)-o(z+C)=C\cap -C\ee for any convex cone $C\subseteq \R^d$ containing $0$ (recall that $C\cap -C$ is a subspace).

\begin{proposition}\label{prop-orientedReay-halfspace-uniquelydetermines}
Let $\mathcal R=(\mathcal X_1\cup \{\mathbf v_1\},\ldots,\mathcal X_s\cup\{\mathbf v_s\})$ and $\mathcal R'=(\mathcal X'_1\cup \{\mathbf v'_1\},\ldots,\mathcal X'_{s'}\cup\{\mathbf v'_{s'}\})$ be two oriented Reay systems in $\R^d$. Suppose $\x\in \mathcal X_1\cup \ldots\cup \mathcal X_s\cup\{\mathbf v_i:\;i\in [1,s]\}$ and $\x'\in \mathcal X'_1\cup \ldots\cup \mathcal X'_s\cup\{\mathbf v'_i:\;i\in [1,s']\}$ with $\x=\x'$. Then $\partial_\mathcal R(\{\x\})=\partial_{\mathcal R'}(\{\x'\})$, where $\partial_\mathcal R(\{\x\})=\partial(\{\x\})\subseteq \mathcal X_1\cup\ldots\cup \mathcal X_s$ and $\partial_{\mathcal R'}(\{\x\})=\partial(\{\x'\})\subseteq \mathcal X'_1\cup\ldots\cup \mathcal X'_s$ are the respective sets for $\x=\x'$  when considered as a half-space for the oriented Reay system $\mathcal R$, and when considered as a half-space for the oriented Reay system $\mathcal R'$.
\end{proposition}

\begin{proof}
Let $\x_1,\ldots,\x_r\in \partial_\mathcal R(\{\x\})$ and $\x'_1,\ldots,\x'_{r'}\in \partial_{\mathcal R'}(\{\x'\})$ be the distinct half-spaces from $\partial_\mathcal R(\{\x\})$ and $\partial_{\mathcal R'}(\{\x'\})$. Observe that $\overline{\C^\cup (\partial_\mathcal R(\{\x\}))}=\overline{\x_1+\ldots+\x_r}=\overline \x_1+\ldots+\overline \x_r$, and likewise $\overline{\C^\cup (\partial_{\mathcal R'}(\{\x'\}))}=\overline{\x'_1+\ldots+\x'_r}=\overline \x'_1+\ldots+\overline \x'_{r'}$ (recall Proposition \ref{prop-orReay-BasicProps}.2).
Since $\x=\x'$, we have $$\C^\cup (\partial_\mathcal R(\{\x\}))=\x\cap \partial(\x)=\x'\cap \partial(\x')=\C^\cup(\partial_{\mathcal R'}(\{\x'\})).$$
Let $\mathcal E\subseteq \overline{\C^\cup (\partial_\mathcal R(\{\x\}))}=\overline{\C^\cup (\partial_{\mathcal R'}(\{\x'\}))}$ be the lineality subspace. Since $\partial_\mathcal R(\{\x\})$ and $\partial_{\mathcal R'}(\{\x'\})$ are both support sets, and thus virtual independent,  Proposition \ref{prop-orReay-BasicProps}.4 implies that $$\mathcal E=\R^\cup\la\partial_\mathcal R^2(\{\x\})\ra=\partial(\x_1)+\ldots+\partial(\x_r)=
\partial(\x'_1)+\ldots+\partial(\x'_{r'})=\R^\cup\la \partial^2_{\mathcal R}(\{\x'\})\ra.$$ Moreover, $\pi(x_1),\ldots, \pi(x_r)$ are distinct, linearly independent elements by Proposition \ref{prop-orReay-BasicProps}.3, where $\pi:\R^d\rightarrow \mathcal E^\bot$ is the orthogonal projection. Likewise $\pi(x'_1),\ldots, \pi(x'_{r'})$ are distinct, linearly independent elements. We have $\C(\pi(x_1),\ldots,\pi(x_r))=\pi(\x_1)+\ldots+\pi(\x_r)=
\pi\Big(\C^\cup(\partial_\mathcal R(\{\x\}))\Big)=\pi\Big(\C^\cup(\partial_{\mathcal R'}(\{\x'\}))\Big)=
\pi(\x'_1)+\ldots+\pi(\x'_{r'})=\C(\pi(x'_1),\ldots,\pi(x'_{r'}))$. Thus, since $\{\pi(x_1),\ldots,\pi(x_r)\}$ and $\{\pi(x'_1),\ldots,\pi(x'_{r'})\}$ are both linearly independent subsets of distinct elements, we conclude that  $r=r'$ with (after re-indexing appropriately) $$\pi(\x_i)=\pi(\x'_i)\quad\mbox{ for all $i\in [1,r]$.}$$ It remains to show $\x_i=\x'_i$ for every $i\in [1,r]$. By replacing each $x'_i$ with an appropriate positive scalar multiple, we can w.l.o.g. also assume $\pi(x_i)=\pi(x'_i)$ for all $i\in [1,r]$.

For $j\in [1,r]$, let $Z_j\subseteq \C^\cup(\partial_\mathcal R(\{\x\}))=\C^\cup(\partial_{\mathcal R'}(\{\x'\}))$ consist of all $z\in \C^\cup(\partial_\mathcal R(\{\x\}))=\C^\cup(\partial_{\mathcal R'}(\{\x'\}))$ with $\pi(z)=\pi(x_j)=\pi(x'_j)$. Since $\{\pi(x_1),\ldots,\pi(x_r)\}=\{\pi(x'_1),\ldots,\pi(x'_{r})\}$ is a  linearly independent subset of distinct elements and $\x_1+\ldots+\x_r=\C^\cup (\partial_\mathcal R(\{\x\}))=\C^\cup (\partial_{\mathcal R'}(\{\x'\}))=\x'_1+\ldots+\x'_r$, it follows that \begin{align}\label{linlike}Z_j&=x_j+\partial(\x_j)+\C^\cup(\partial^2_\mathcal R(\{\x\}))=x'_j+\partial(\x'_j)+\C^\cup(\partial_{\mathcal R'}^2(\{\x'\})).
\end{align} Now $$C_j:=Z_j-x_j=\partial(\x_j)+\C^\cup(\partial_\mathcal R^2(\{\x\}))\quad\und\quad C'_j:=Z_j-x'_j=\partial(\x'_j)+\C^\cup(\partial_{\mathcal R'}^2(\{\x'\}))$$ are both convex cones containing $0$ by Proposition \ref{prop-orReay-BasicProps}.2. Consequently, in view of \eqref{apex-lin}, we have $$o(Z_j)-o(Z_j)=C_j\cap -C_j=C'_j\cap -C'_j.$$ Since $\partial(\x_j)$ is a subspace,  $\partial(\x_j)\subseteq C_j\cap -C_j$. If the inclusion were strict, then $\pi_j\Big(\C^\cup(\partial_\mathcal R^2(\{\x\}))\Big)=\pi_j(C_j)$ would have nontrivial lineality subspace, where $\pi_j:\R^d\rightarrow \partial(\x_j)^\bot=\R^\cup\la \partial(\{\x_j\})\ra^\bot$ is the orthogonal projection. Proposition \ref{prop-orReay-modulo} implies that $\pi_j(\mathcal R)$ is an oriented Reay system with $\pi_j\Big(\C^\cup(\partial_\mathcal R^2(\{\x\}))\Big)=\C^\cup(\partial_\mathcal R^2(\{\pi_j(\x)\}))$, and  Proposition \ref{prop-orReay-BasicProps}.4 ensures that
$\C^\cup(\partial_\mathcal R^2(\{\pi_j(\x)\}))$ has trivial lineality subspace. As a result, the inclusion $\partial(\x_j)\subseteq C_j\cap -C_j$ cannot be strict, forcing $\partial(\x_j)=C_j\cap -C_j$. An analogous argument shows $\partial(\x'_j)=C'_j\cap -C'_j$. Thus \be\label{baselined}\partial(\x'_j)=C'_j\cap -C'_j=C_j\cap -C_j=\partial(\x_j),\quad\mbox{ for every $j\in [1,r]$}.\ee


 Let us next show that \be\label{goldnug}(\x_1+\ldots+\x_r)\cap \partial(\x_j)=\x_j\cap \partial(\x_j), \quad\mbox{for each $j\in [1,r]$}.\ee The inclusion $\x_j\cap \partial(\x_j)\subseteq (\x_1+\ldots+\x_r)\cap \partial(\x_j)$ is trivial. Let   $y\in (\x_1+\ldots+\x_r)\cap \partial(\x_j)$ be arbitrary.
 Since $y\in \x_1+\ldots+\x_r$, it follows that $\{\x_1,\ldots,\x_r\}=\partial_{\mathcal R}(\{\x\})$ encases $y$, so  there is some $\mathcal Z\preceq \partial_{\mathcal R}(\{\x\})$ with   $\mathcal Z\subseteq \darrow \{\x_1,\ldots,\x_r\}=\darrow \partial_\mathcal R(\{\x\})$ that minimally encases $y$. In view of the minimality of $\mathcal Z$, there is a subset $Z\subseteq \R^d$ of representatives for the half-spaces from $\mathcal Z$ such that $y=\Summ{z\in Z}z$.
 Since $\mathcal Z\subseteq \darrow \partial_\mathcal R(\{\x\})$, we can extend the representative set $Z$ to a set of representatives $\darrow \partial_\mathcal R(\{x\})$ for $\darrow \partial_\mathcal R(\{\x\})\subseteq \mathcal X$, which will then be linearly independent in view of  Proposition \ref{prop-reay-basis-properties}.1. Thus, since $\Summ{z\in Z}z=y\in \partial(\x_j)=\R^\cup\la\partial_\mathcal R(\{\x_j\})\ra=\R\la \darrow \partial_\mathcal R(\{x_j\})\ra$ (by Proposition \ref{prop-orReay-BasicProps}.1) with $\darrow \partial_\mathcal R(\{x_j\})\subseteq \darrow \partial_\mathcal R(\{x\})$, \ $Z\subseteq \darrow \partial_\mathcal R(\{x\})$ and $\darrow \partial_\mathcal R(\{x\})$ linearly independent, it follows that $Z\subseteq \darrow\partial_\mathcal R(\{x_j\})$, whence $\mathcal Z\subseteq \darrow \partial_\mathcal R(\{\x_j\})$.
 But this means $y=\Summ{z\in Z}z\in \C^\cup(\mathcal Z)\subseteq \C^\cup (\darrow \partial_\mathcal R(\{\x_j\}))=\C^\cup(\partial_\mathcal R(\{\x_j\}))=\x_j\cap \partial(\x_j)$. Since $y\in (\x_1+\ldots+\x_r)\cap \partial(\x_j)$ was arbitrary, this shows the nontrivial inclusion $(\x_1+\ldots+\x_r)\cap \partial(\x_j)\subseteq \x_j\cap \partial(\x_j)$, and \eqref{goldnug} is established.

 By an analogous argument, we conclude that $$(\x'_1+\ldots+\x'_r)\cap \partial(\x'_j)=\x'_j\cap \partial(\x'_j),\quad\mbox{ for each $j\in [1,r]$}.$$ We have $\x_1+\ldots+\x_r=\C^\cup(\partial_\mathcal R(\{\x\}))=\C^\cup(\partial_{\mathcal R'}(\{\x'\}))=\x'_1+\ldots+\x'_r$, while \eqref{baselined} gives  $\partial(\x_j)=\partial(\x'_j)$. In consequence, $\x'_j\cap \partial(\x'_j)=(\x'_1+\ldots+\x'_r)\cap \partial(\x'_j)=(\x_1+\ldots+\x_r)\cap \partial(\x_j)=\x_j\cap \partial(\x_j)$, so \be\label{alittlecloser} \x_j\cap \partial(\x_j)=\x'_j\cap \partial(\x'_j),\quad\mbox{ for every $j\in [1,r]$}.\ee
 As a result, \begin{align*}&\C^\cup(\partial_\mathcal R^2(\{\x\}))=\C^\cup(\partial_\mathcal R(\{\x_1\}))+\ldots+\C^\cup(\partial_\mathcal R(\{\x_r\}))=(\x_1\cap \partial(\x_1))+\ldots+(\x_r\cap \partial(\x_r))\\&=(\x'_1\cap \partial(\x'_1))+\ldots+(\x'_r\cap \partial(\x'_r))=\C^\cup(\partial_{\mathcal R'}(\{\x'_1\}))+\ldots+\C^\cup(\partial_{\mathcal R'}(\{\x'_r\}))=\C^\cup(\partial_{\mathcal R'}^2(\{\x'\})).\end{align*}
 Combining the above equality with \eqref{baselined} and \eqref{linlike}, we conclude that $$x_j-x'_j+\partial(\x_j)+\C^\cup(\partial_\mathcal R^2(\{\x\}))=
 \partial(\x_j)+\C^\cup(\partial_\mathcal R^2(\{\x\})),$$ which readily implies that $x_j-x'_j$ is contained in the lineality subspace of $\partial(\x_j)+\C^\cup(\partial_\mathcal R^2(\{\x\}))=C_j$, which is equal to $\partial(\x_j)$ by \eqref{baselined}. This shows $x_j-x'_j\in \partial(\x_j)$, which combined with \eqref{baselined} ensures that $\x_j$ and $\x'_j$ linearly span the same subspace, and then $x_j-x'_j\in \partial(\x_j)$ further ensures that $x_j$ and $x'_j$ lie on the same side of the codimension one subspace $\partial(\x_j)=\partial(\x'_j)$. Thus $\overline \x_j=\overline \x'_j$, which combined with  \eqref{alittlecloser} yields the desired conclusion $\x_j=\x'_j$, for all $j\in [1,r]$.\end{proof}

\section{Virtual Reay Systems}\label{sec-virtual}

When trying to better understand the geometric properties of infinite subsets $G_0\subseteq \R^d$, one of the crucial problems encountered is that small perturbations of linearly dependent sets can result in linearly independent sets. This allows limits of linearly independent subsets to degenerate into linearly dependent ones. We aim to better understand $G_0$ by looking at the limiting behavior of sequences of terms from $G_0$. However, for this to be effective, we need to avoid introducing linear dependencies into the limiting structures that do not exist in the original sequences. This will be accomplished by a careful development of the following generalization of a Reay system.  It may be  helpful to think of a virtual Reay system, defined below,  as a convergent family of ordinary Reay systems whose limiting structure is an Oriented Reay system. This will be made more formal later.

\begin{definition}
Let $G_0\subseteq \R^d$ be a subset.
 If $\mathcal R=(\mathcal X_1\cup \{\mathbf v_1\},\ldots,\mathcal X_s\cup\{\mathbf v_s\})$ is an oriented Reay system in $\R^d$ such that
\begin{itemize}
 \item[(V1)] every $\x\in \mathcal X_1\cup \{\mathbf v_1\}\cup \ldots\cup\mathcal X_1\cup\{\mathbf v_s\}$ has an asymptotically filtered sequence $\{\x(i)\}_{i=1}^\infty$ of terms $\x(i)\in G_0$ with  limit $\vec u_\x=(u^\x_1,\ldots,u^\x_{t_\x})$
     and $t_\x\geq 1$ such that  $-\vec u_\x^\triangleleft=-(u^\x_1,\ldots,u^\x_{t_\x-1})$ is minimally encased by $\partial(\{\x\})$ and $\overline \x=\R^\cup\la \partial(\{\x\})\ra+\R_+ u^\x_{t_\x}$,
 \end{itemize}
     then we say that $\mathcal R$ is a \textbf{virtual Reay system} in (or over) $G_0$ for the subspace $\R^\cup\la\mathcal X_1\cup\ldots\cup \mathcal X_s\ra$.
 \end{definition}


 \begin{definition}
     If $\mathcal R$ is a virtual Reay system in $G_0\subseteq \R^d$ and
     \begin{itemize}
 \item[(V2)]
 $\vec u_\x$ is anchored for every $\x\in \mathcal X_1\cup\ldots\cup \mathcal X_s$,
 \end{itemize}
 then we call $\mathcal R$ an \textbf{anchored} virtual Reay system.  If $\mathcal R$ is a virtual Reay system in $G_0\subseteq \R^d$ and
 \begin{itemize}
 \item[(V3)]  $\vec u_{\mathbf v_j}$ is fully unbounded for every $j\in [1,s]$,\end{itemize}
 then we call $\mathcal R$ a \textbf{purely} virtual Reay system.
\end{definition}

 Let $\x(i)=a_i^{(1)}u_1+\ldots+a_i^{(t)}u_t+z_i$ be the representation of $\{\x(i)\}_{i=1}^\infty$ as an asymptotically filtered sequence with limit $(u_1,\ldots,u_t)$,  let $\pi:\R^d\rightarrow \R^\cup \la \partial(\{\x\})\ra^\bot$  and $\pi^\bot:\R^d\rightarrow \R^\cup \la \partial(\{\x\})\ra$ be the orthogonal projections, and let $\overline u_t=\pi(u_t)/\|\pi(u_t)\|$. Then $z_i=w'_i+a'_i\overline u_t+y_i$ for some $w'_i\in \R^\cup \la \partial(\{\x\})\ra$, $a'_i\in \R$ and $y_i\in \R \la \x\ra^\bot=(\R^\cup \la\partial(\{\x\})\ra+\R \overline u_t)^\bot$, and $a_i^{(t)}u_t=a_i^{(t)}\pi^\bot(u_t)+a_i^{(t)}\|\pi(u_t)\|\overline u_t$. Consequently, setting $b^{(t)}_i=a_i^{(t)}\|\pi(u_t)\|+a'_i$ and $w_i=w'_i+a_i^{(t)}\pi^\bot(u_t)$, we find that
 \be\label{x(i)-filtered-form}\x(i)=(a_i^{(1)} u_{1}+\ldots+a_i^{(t-1)} u_{t-1}+w_i)+b_i^{(t)}\overline u_t+y_i\ee with $w_i\in \R^\cup \la \partial(\{\x\})\ra$, $\|w_i\|\in o(a_i^{(t-1)})$
 (since $\|w'_i\|=\|\pi^\bot(z_i)\|\in O(\|z_i\|)\subseteq o(a_i^{(t)})$ and $a_i^{(t)}\in o(a_i^{(t-1)})$)
 and $w_i=0$ when $t=1$, with $b^{(t)}_i\in \Theta(a_i^{(t)})$ and $b_i^{(t)}>0$ for all sufficiently large $i$ (since $a'_i=\|a'_i\overline u_t\|\leq \|z_i\|\in o(a_i^{(t)})$ with $\pi(u_t)\neq 0$ by (V1)), and with $y_i\in \R \la \x\ra^\bot$ and $\|y_i\|\in o(a_i^{(t)})=o(b_i^{(t)})$ (since $\|y_i\|\leq \|z_i\|\in o(a_i^{(t)})$).  We will now generally assume, by discarding the first few terms in $\{\x(i)\}_{i=1}^\infty$, that any representative sequence $\{\x(i)\}_{i=1}^\infty$ must satisfy $b_i^{(t)}>0$ for all $i$ in any virtual Reay system.

As a matter of notation,
  we let \ber\nn \tilde \x(i)&=&(a_i^{(1)} u_{1}+\ldots+a_i^{(t-1)} u_{t-1}+w_i)+b_i^{(t)}\overline u_t\quad\und\\\nn \tilde \x^{(t-1)}(i)&=&(a_i^{(1)} u_{1}+\ldots+a_i^{(t-1)} u_{t-1}+w_i).\eer
  Note that $\{\tilde\x^{(t-1)}(i)\}_{i=1}^\infty$ is an asymptotically  filtered sequence of terms $\tilde\x^{(t-1)}(i)\in \R^\cup\la \partial(\{\x\})\ra$ with limit $(u_1,\ldots,u_{t-1})$ and that $\tilde \x(i)$ is always a representative for the half-space $\x$ for any $i\geq 1$ (assuming we have discarded terms to attain $b_i^{(t)}>0$ for all $i$). If $\mathcal B,\,\mathcal C\subseteq \mathcal X_1\cup\{\mathbf v_1\}\cup\ldots\cup \mathcal X_s\cup\{\mathbf v_s\}$ with $\mathcal B\subseteq \mathcal C$ and $k=(i_\x)_{\x\in \mathcal C}$ is a tuple of indices $i_\x\geq 1$, then we let $$B(k)=\{\x(i_\x):\; \x\in\mathcal B\}\quad\und\quad \tilde B(k)=\{\tilde \x(i_\x):\; \x\in\mathcal B\}\quad\mbox{ for $k=(i_\x)_{\x\in \mathcal C}$.}$$
  In view of Proposition \ref{prop-orReay-minecase-char}.1 and (V1), the limit $\vec u_\x$ uniquely determines the half-space $\x\in \mathcal X_j\cup \{\mathbf v_j\}$ (assuming all
 $\y\in \mathcal X_1\cup\{\mathbf v_1\}\cup \ldots\cup \mathcal X_{j-1}\cup \{\mathbf v_{j-1}\}$ have already been determined by their respective limits $\vec u_\y$).  It is important to realize  that the half-spaces $\x$ defined in (V1) depend only on the tuples $\vec u_\x$ and not on the representative sequences $\x(i)$.
  It is even possible  for distinct half-spaces $\x,\,\y\in \mathcal X_1\cup\{\mathbf v_1\}\cup\ldots\cup \mathcal X_s\cup \{\mathbf v_s\}$ to have their representative  sequences $\{\x(i)\}_{i=1}^\infty$ and $\{\y(i)\}_{i=1}^\infty$ being the same, since the same sequence may be considered as an asymptotically  filtered sequence with limit $(u_1,\ldots,u_t)$ and also one with limit $(u_1,\ldots,u_{t'})$ with $t<t'$. There are then many compatible sequences $\{\x(i)\}_{i=1}^\infty$ that can be used for the representative sequence associated to $\x$.
 Indeed, any asymptotically filtered sequence of terms from $G_0$ having limit $\vec u_\x=(u_1,\ldots,u_t)$ with $\vec u^\triangleleft_\x$  minimally encased by $\partial(\{\x\})$, and $u_t$ a representative for $\x$, will do once the  first few terms with  $b_i^{(t)}\leq 0$ are removed. Of course, it is important at least one such asymptotically filtered sequence exist. Besides the existence of an asymptotically filtered sequence $\{\x(i)\}_{i=1}^\infty$ of terms $\x(i)\in G_0$ with limit $\vec u_\x$, it is relevant to the properties (V2) and (V3) whether $\{\x(i)\}_{i=1}^\infty$ can be chosen to have $\vec u_\x$ as a fully unbounded or  anchored limit, as a different choice can change whether $\vec u_\x$ is considered as a fully unbounded  or anchored in $\mathcal R$.

 We continue with the analogue of Proposition \ref{prop-orReay-modulo} for virtual Reay systems, showing we once again have a well-behaved notion of a quotient virtual Reay system modulo $\R^\cup \la \mathcal A\ra$, for any $\mathcal A\subseteq \mathcal X_1\cup\ldots\cup \mathcal X_s$.

\begin{proposition}\label{prop-VReay-modulo} Let $G_0\subseteq \R^d$ be a subset.
Suppose $\mathcal R=(\mathcal X_1\cup \{\mathbf v_1\},\ldots,\mathcal X_s\cup\{\mathbf v_s\})$ is a virtual Reay system in  $G_0$.
Let $\mathcal A\subseteq \mathcal X_1\cup\ldots\cup\mathcal X_s$ and let $\pi:\R^d\rightarrow \R^\cup \la \mathcal A\ra^\bot$ be the orthogonal projection. Then $\pi(\mathcal R)=\big(\mathcal X_j^\pi\cup \{\pi(\mathbf v_j)\}\big)_{j\in J}$ is a virtual Reay system in $\pi(G_0)$ with $$\vec u_{\pi(\x)}=\pi(\vec u_\x),\quad \pi(\vec u_\x)^\triangleleft=\pi(\vec u_\x^\triangleleft),\quad\und\quad\partial(\{\pi(\x)\})=\partial(\{\x\})^\pi,$$ for each $\x\in \bigcup_{i\in J}(\mathcal X_i\cup\{\mathbf v_i\})$ with $\pi(\x)\neq \{0\}$, having representative sequence $\{\pi(\x)(i)\}_{i=1}^\infty$ the sufficiently large index terms in $\{\pi(\x(i))\}_{i=1}^\infty$. Moreover, if $\mathcal R$ is purely virtual, then so is  $\pi(\mathcal R)$, and if $\mathcal R$ is anchored, then so is $\pi(\mathcal R)$.
\end{proposition}

\begin{proof}
It follows from  Proposition \ref{prop-orReay-modulo} that  $\pi(\mathcal R)=\big(\mathcal X_j^\pi\cup \{\pi(\mathbf v_j)\}\big)_{j\in J}$ is an oriented Reay system with $\partial(\{\pi(\x)\})=\partial(\{\x\})^\pi$ a support set for all half-spaces $\x\in \bigcup_{j\in J}(\mathcal X_j\cup \{\mathbf v_j\})$ with $\pi(\x)\neq \{0\}$. Let $\x\in \bigcup_{j\in J}(\mathcal X_j\cup \{\mathbf v_j\})$ with $\pi(\x)\neq \{0\}$ be arbitrary, let $\vec u_\x=(u_1,\ldots,u_{t})$ and let $\mathcal E=\R^\cup\la \mathcal A\ra$.
By Proposition \ref{prop-orReay-minecase-char}.4 and (V1), the support set $\partial(\{\x\})^\pi$ minimally encases $-\pi(\vec u^\triangleleft_\x)$. By Proposition
\ref{prop-infinite-limits-proj}, the sufficiently large index terms in $\{\pi(\x(i))\}_{i=1}^\infty$ are an asymptotically filtered sequence with limit $\pi(\vec u_\x)$.
If $\pi(\vec u_\x^\triangleleft)=\pi(\vec u_\x)$, then this would mean $u_{t}\in \mathcal E+\R\la u_1,\ldots,u_{t-1}\ra\subseteq \mathcal E+\R^\cup \la \partial(\{\x\})\ra$, with inclusion following since $\partial(\{\x\})$ minimally encases $-\vec u_\x^\triangleleft$ by (V1).
However, since $\pi(\x)\neq \{0\}$ is a relative half-space by Proposition \ref{prop-orReay-modulo}.1, we must have $\x\nsubseteq \mathcal E+\partial(\x)=\mathcal E+ \R^\cup\la \partial(\{\x\})\ra$, contrary to what was just shown in view of $u_{t}\in \x^\circ$ (since  (V1) holds for $\mathcal R$).
Therefore we instead conclude that $\pi(\vec u_\x^\triangleleft)\neq \pi(\vec u_\x)$, whence $\pi(\vec u_\x^\triangleleft)=\pi(\vec u_\x)^\triangleleft$ with the last coordinate in $\pi(\vec u_\x)$ equal to  $\overline u_{t}:=\tau(u_{t})/\|\tau(u_{t})\|$,
 where $\tau:\R^d\rightarrow (\mathcal E+\R\la u_1,\ldots,u_{t-1}\ra)^\bot$ is the orthogonal projection.
 Thus (V1) for $\mathcal R$ and \eqref{halfspace-lintrans} give $\pi(\overline \x)=\R^\cup\la\partial(\{\x\})^\pi\ra+\R_+ \overline u_{t}$ and $\partial(\pi(\x))\cap \pi(\x)=\pi(\partial(\x)\cap \x)=\C^\cup(\partial(\{\x\})^\pi)$, which establishes (V1) for $\pi(\mathcal R)$. Since the last coordinate of $\pi(\vec u_\x)$ equals $\overline u_{t}=\tau(u_{t})/\|\tau(u_{t})\|$ (ensuring $r_\ell=t$ in Proposition \ref{prop-infinite-limits-proj} when applied to $\pi(\vec u_\x)$), Proposition \ref{prop-infinite-limits-proj}.1 ensures that the limit $\pi(\vec u_\x)$ is anchored if and only if $\vec u_\x$ is anchored, ensuring that (V2) holds for $\pi(\mathcal R)$ when it does for $\mathcal R$, and also that (V3) holds for $\pi(\mathcal R)$ when it does for  $\mathcal R$, which completes the proof.
\end{proof}

Let $\mathcal R=(\mathcal X_1\cup \{\mathbf v_1\},\ldots,\mathcal X_s\cup \{\mathbf v_s\})$ be a virtual Reay system in $G_0\subseteq \R^d$. Now suppose  $\mathcal R'=(\mathcal Y_1\cup \{\mathbf w_1\},\ldots,\mathcal Y_r\cup \{\mathbf w_r\})$ is another virtual Reay system in $G_0$ with $\mathcal Y_1\cup \ldots\cup \mathcal Y_{r-1}\subseteq \mathcal \bigcup_{i=1}^s\mathcal X_i$ and $\mathcal Y_r\subseteq \bigcup_{i=1}^s(\mathcal X_i\cup \{\mathbf v_i\})$. Then each $\x\in \mathcal Y_1\cup \ldots\cup \mathcal Y_{r}$ has a boundary neighborhood $\partial_{\mathcal R}(\{\x\})$ when considered as a half-space from $\mathcal R$ and a boundary neighborhood $\partial_{\mathcal R'}(\{\x\})$ when considered as a half-space from $\mathcal R'$. Proposition \ref{prop-orientedReay-halfspace-uniquelydetermines} ensures both these neighborhoods are equal: $\partial_{\mathcal R}(\{\x\})=\partial_{\mathcal R'}(\{\x\})$.
Thus the  partial order in $\mathcal R'$ agrees with that in $\mathcal R$ and means there is no need to distinguish whether $\x$ lies in $\mathcal R$ or $\mathcal R'$ when dealing with quantities like $\partial(\{\x\})$ or $\darrow \x$.

By Proposition \ref{prop-orReay-minecase-char}, urbane minimal encasement of a tuple $\vec u$ by a \emph{support} set from $\mathcal R$ corresponds to an ascending chain of support sets from $\mathcal R$. Our next major goal is to show this chain can be completed to an entire virtual Reay system, allowing us to use our machinery for virtual Reay systems when dealing with the support sets associated to urbane minimal encasement by a support set. We begin with the following lemma.
 Informally, we consider a fixed virtual Reay system $\mathcal R$ together with a ``sub-'' virtual Reay system $\mathcal R_A$ and quotient virtual Reay system $\mathcal R_C$, both of $\mathcal R$. Lemma \ref{lemma-VReay-extension} then shows that, under certain conditions, these two virtual Reay systems can be combined to create a new  sub- virtual Reay System $\mathcal R'$ of $\mathcal R$.

\begin{lemma}
\label{lemma-VReay-extension}
Let $\mathcal R=(\mathcal X_1\cup\{\mathbf v_1\},\ldots,\mathcal X_s\cup \{\mathbf v_s\})$ and $\mathcal R_A=(\mathcal A_1\cup\{\mathbf a_1\},\ldots,\mathcal A_t\cup \{\mathbf a_t\})$ be  virtual Reay systems in $G_0\subseteq \R^d$ with $\mathcal A=\bigcup_{i=1}^t\mathcal A_i\subseteq \bigcup_{i=1}^s\mathcal X_i$, and
let $\pi:\R^d\rightarrow \R^\cup\la\mathcal A\ra^\bot$ be the orthogonal projection  with  $\pi(\mathcal R)=(\mathcal X_i^\pi\cup \{\pi(\mathbf v_i)\})_{i\in J}$. Suppose $\mathcal R_C=(\pi(\mathcal B_1)\cup \{\mathbf c_1\},\ldots,\mathcal \pi(\mathcal B_r)\cup \{\mathbf c_r\})$ is a virtual Reay system in $\pi(G_0)$, for some
$\mathcal B_j\subseteq \bigcup_{i=1}^s\mathcal X_i$ for $j\in [1,r-1]$ and  $\mathcal B_r\subseteq \mathcal X_1\cup \ldots\mathcal X_s\cup \{\mathbf v_i:\;i\in J\}$ with $\vec u_{\pi(\x)}=\pi(\x)$ for each
 $\x\in \mathcal B_1\cup \ldots\cup \mathcal B_r$, and that  each $\mathbf c_j$ for $j\in [1,r]$ is defined by a limit $\pi(\vec u_j)$ with $\vec u_j$ the limit of an asymptotically  filtered sequence of terms from $G_0$,
  $\pi(\vec u_j^\triangleleft)=\pi(\vec u_j)^\triangleleft$ and  $-\vec u_j^\triangleleft$ encased by $\mathcal A\cup \mathcal C_{j}$ for some  $\mathcal C_{j}\subseteq \mathcal B_1\cup\ldots\cup \mathcal B_{j-1}$ with $\pi(\mathcal C_{j})=\partial(\{\mathbf c_j\})$. Then
$$\mathcal R'=(\mathcal A_1\cup\{\mathbf a_1\},\ldots,\mathcal A_s\cup \{\mathbf a_s\},\mathcal B_1\cup \{\mathbf b_1\},\ldots,\mathcal B_r\cup\{\mathbf b_r\})$$ is virtual Reay system in $G_0$  with each $\mathbf b_j$ for $j\in [1,r]$ defined by the limit $\vec u_j$, $\pi(\mathbf b_j)=\mathbf c_j$, $\partial(\{\mathbf b_j\})\preceq \mathcal A\cup \mathcal C_{j}$ and $\mathcal C_{j}\subseteq \mathcal \partial(\{\mathbf b_j\})$. Moreover, if $\mathcal R_A$ and $\mathcal R_C$ are  purely virtual, then so is $\mathcal R'$, and if $\mathcal R_A$ and $\mathcal R_C$ are  anchored, then so is $\mathcal R'$.
\end{lemma}

\begin{proof}
By Proposition \ref{prop-VReay-modulo}, $\pi(\mathcal R)$ is a virtual Reay system in $\pi(G_0)$. Let $j\in [1,r]$ be arbitrary. Since $\mathcal R_C$ is an oriented Reay system, we have $\pi(\x)\neq \{0\}$ for every $\x\in \mathcal B_j$, and thus also for every $\x\in \mathcal C_{j}$. Since $\mathcal R_A$ is an oriented Reay system, we must have $\darrow \mathcal A=\mathcal A$.

Since $-\vec u_j^\triangleleft$ is encased by $\mathcal A\cup \mathcal C_{j}$, there exists $\mathcal C'_{j}\preceq \mathcal A\cup \mathcal C_{j}$ which minimally encases $-\vec u_j^\triangleleft$, and since $\mathcal A\cup \mathcal C_{j}\subseteq \mathcal X_1\cup \ldots\cup \mathcal X_s$, it must do so urbanely and be a support set. If $\mathcal C_{j}\nsubseteq \mathcal C'_{j}$, then ${(\mathcal C'_{j})}^\pi\prec \mathcal C_{j}^\pi=\partial(\{\mathbf c_j\})$ follows in view of Proposition \ref{prop-orReay-modulo}.3 since $\pi(\x)\neq \{0\}$ for all $\x\in \mathcal C_{j}$. By Proposition \ref{prop-orReay-minecase-char}.4, $(\mathcal C'_{j})^\pi$ minimally encases $-\pi(\vec u_j^\triangleleft)=-\pi(\vec u_j)^\triangleleft$, with the equality holding by hypothesis, but then, in view of
$(\mathcal C'_{j})^\pi\prec \partial(\{\mathbf c_j\})$, we contradict  that $\partial(\{\mathbf c_j\})$ minimally encases $-\pi(\vec u_j)^\triangleleft$ by (V1). Therefore we conclude that $\mathcal C_{j}\subseteq \mathcal C'_{j}$. Hence, since $(\mathcal C'_{j})^*=\mathcal C'_{j}$ (as $\mathcal C'_{j}$ is a support set), let $\mathcal A'_j=\mathcal C'_{j}\setminus \mathcal C_{j}\subseteq \darrow \mathcal A=\mathcal A$. Then $\mathcal C'_{j}=\mathcal A'_j\cup \mathcal C_{j}\subseteq \mathcal A\cup \mathcal B_1\cup\ldots\cup \mathcal B_{j-1}\subseteq \mathcal X_1\cup \ldots\cup \mathcal X_s$ and $(\mathcal C'_{j})^\pi=\mathcal C_{j}^\pi$.

By hypothesis, $\mathcal R_A=(\mathcal A_1\cup \{\mathbf a_1\},\ldots,\mathcal A_s\cup \{\mathbf a_s\})$ is a virtual Reay system.  If $\x\in \mathcal B_{j}$ for $j\in [1,r]$, then we know (V1) and (OR1) hold for $\x$ with the set $\partial(\{\x\})$ when we consider $\x$ as part of the virtual Reay system $\mathcal R$. Thus, to show this is also the case when we consider $\x$ as part of $\mathcal R'$, we just need to know $\partial(\{\x\})\subseteq \mathcal A\cup  \mathcal B_1\cup\ldots\cup \mathcal B_{j-1}$. Since $\pi(\x)\in \pi(\mathcal B_j)$, we have $\pi(\x)\neq \{0\}$ and $\partial(\{\x\})^\pi=\partial(\{\pi(\x)\})\subseteq \pi(\mathcal B_1)\cup \ldots\cup \pi(\mathcal B_{j-1})$, where the first equality follows from proposition \ref{prop-orReay-modulo}.2. Thus, in view of the injectivity of $\pi$ (Proposition \ref{prop-orReay-modulo}.1) and Proposition \ref{prop-orReay-BasicProps}.9, we have $\partial(\{\x\})\subseteq \darrow \mathcal A\cup \mathcal B_1\cup\ldots\cup \mathcal B_{j-1}=\mathcal A\cup \mathcal B_1\cup\ldots\cup \mathcal B_{j-1}$, as desired. By Proposition \ref{prop-VReay-modulo} and hypothesis, we have $\pi(\vec u_\x)=\vec u_{\pi(\x)}$ (both in $\pi(\mathcal R)$ and $\mathcal R_C$). Thus, if $\mathcal R_\mathcal C$ is anchored, then $\pi(\vec u_\x)=\vec u_{\pi(\x)}$ is anchored, which is only possible if $\vec u_\x$ is anchored (cf.  Proposition \ref{prop-infinite-limits-proj}.1).

Let $\vec u_j=(u_{j,1},\ldots,u_{j,t_j})$ for $j\in [1,r]$.
 The hypothesis $\pi(\vec u_j^\triangleleft)=\pi(\vec u_j)^\triangleleft$ simply means $u_{j,t_j}\notin \R^\cup\la \mathcal A\ra+\R\la u_{j,1},\ldots,u_{j,t_j-1}\ra$. Thus,
since $\mathcal R_C$ is a virtual Reay system with $(\mathcal C'_{j})^\pi=\partial(\{\mathbf c_j\})$ and $\mathbf c_j$ defined by the limit $\pi(\vec u_j)$,
it follows that $u_{j,1},\ldots,u_{j,t_j-1}\in \R^\cup\la \mathcal A\ra+\R^\cup\la \mathcal C'_j\ra$ with $u_{j,t_j}\notin \R^\cup\la \mathcal A\ra+\R^\cup\la \mathcal C'_j\ra$.
Consequently, since $\mathcal C'_{j}\subseteq \mathcal A\cup \mathcal B_1\cup\ldots\cup \mathcal B_{j-1}$ minimally encases $-\vec u_j^\triangleleft$ with $\vec u_j$ the limit of an asymptotically  filtered sequence of terms from $G_0$, it follows that we can define a half-space $\mathbf b_j$ by the limit $\vec u_j$ and it will satisfy the needed requirements in (V1) and (OR1) with $\partial(\{\mathbf b_j\})=\mathcal C'_{j}$.
Since $\mathbf c_j$ is defined by the limit $\pi(\vec u_j)$ with $\pi(\vec u_j^\triangleleft)=-\pi(\vec u_j)^\triangleleft$, we have $\overline{\mathbf c}_j=\partial(\mathbf c_j)+\R_+\pi(u_{j,t_j})=\R^\cup \la \partial(\{\mathbf c_j\})\ra+\R_+\pi(u_{j,t_j})$ and $\mathbf c_j\cap \partial(\mathbf c_j)=\C^\cup (\partial(\{\mathbf c_j\}))$. Thus, since $$\partial(\{\mathbf b_j\})^\pi=(\mathcal C'_{j})^\pi=\mathcal C_j^\pi=\pi(\mathcal C_{j})=\partial(\{\mathbf c_j\})$$ with $\mathbf b_j$ defined by the limit $\vec u_j$, we have $\pi(\mathbf b_j)\cap \partial(\pi(\mathbf b_j))=\pi(\mathbf b_j\cap \partial(\mathbf b_j))=\C^\cup(\partial(\{\mathbf b_j\})^\pi)=\C^\cup(\partial(\{\mathbf c_j\}))=
\mathbf c_j\cap \partial(\mathbf c_j)$ and $\pi (\overline{\mathbf b}_j)=\pi(\partial(\mathbf b_j))+\R_+\pi(u_{j,t_j})=\R^\cup\la\partial(\{\mathbf b_j\})^\pi\ra+\R_+\pi(u_{j,t_j})=\R^\cup\la\partial(\{\mathbf c_j\})\ra+\R_+\pi(u_{j,t_j})=\overline{\mathbf c}_j$, whence $\pi(\mathbf b_j)=\mathbf c_j$.
Moreover, if $\mathcal R_C$ is purely virtual, then $\pi(\vec u_j)$ will be fully unbounded, implying $\vec u_j$ is fully unbounded in view of $\pi(\vec u_j^\triangleleft)=-\pi(\vec u_j)^\triangleleft$ and Proposition \ref{prop-infinite-limits-proj}.1. It remains to show (OR2) holds to complete the proof.

For $j\in [1,r]$, let $\pi_{j-1}:\R^d\rightarrow (\R^\cup\la \mathcal A\ra+\R^\cup\la \mathcal B_1\cup\ldots\cup \mathcal B_{j-1}\ra)^\bot$ be the orthogonal projection. Then $\pi_{j-1}(\mathcal B_j\cup \{\mathbf b_j\})=\pi_{j-1}\pi(\mathcal B_j\cup \{\mathbf b_j\})=\pi_{j-1}(\pi(\mathcal B_j))\cup \{\pi_{j-1}(\mathbf c_j)\}$. Hence it follows from the injectivity of $\pi$ (Proposition \ref{prop-orReay-modulo}.1) and  (OR2) holding for $\mathcal R_C$ that (OR2) holds for $\mathcal R'$.
\end{proof}

The proof of Proposition \ref{prop-VReay-modularCompletion} gives an algorithm by which  $\mathcal R'$ can be constructed. Also worth noting, if $\vec u$ is fully unbounded and the strict truncation of any limit defining a half-space from $\darrow \mathcal B$ in $\mathcal R$ is either trivial or fully unbounded, then $\mathcal R'$ will be purely virtual, and if every $\x\in \darrow \mathcal B$ has $\vec u_\x$ anchored (e.g., if $\mathcal B\subseteq \mathcal X_1\cup \ldots,\cup \mathcal X_s$ with $\mathcal R$ anchored),  then $\mathcal R'$ will be anchored.

\begin{proposition}
\label{prop-VReay-modularCompletion}
Let $G_0\subseteq \R^d$ be a subset.
Suppose $\mathcal R=(\mathcal X_1\cup \{\mathbf v_1\},\ldots,\mathcal X_s\cup\{\mathbf v_s\})$ is a virtual Reay system in  $G_0$ and  $\{x_i\}_{i=1}^\infty$ is an asymptotically filtered sequence of terms $x_i\in G_0$ with limit $\vec u=(u_1,\ldots,u_t)$ such that $-\vec u$ is minimally encased urbanely by a support set $\mathcal B\subseteq \mathcal X_1\cup\{\mathbf v_1\}\cup\ldots\cup \mathcal X_s\cup\{\mathbf v_s\}$.
Let $$\emptyset=\mathcal C_0\prec \mathcal C_1\prec \ldots\prec C_{\ell-1}\subseteq \mathcal X_1\cup \ldots\cup \mathcal X_s\quad\und\quad \mathcal C_{\ell-1}\prec \mathcal C_{\ell}=\mathcal B$$ be the support sets and let $1=r_1<\ldots <r_\ell<r_{\ell+1}=t+1$ be the indices  given by Proposition \ref{prop-orReay-minecase-char} for $\mathcal B$.
  Then there exists a virtual Reay system $$\mathcal R'=(\mathcal C^{(1)}_1\cup\{\mathbf v_1^{(1)}\},\ldots,\mathcal C^{(1)}_{s_1}\cup\{\mathbf v^{(1)}_{s_1}\},
  \mathcal C^{(2)}_1\cup\{\mathbf v_1^{(2)}\},\ldots,\mathcal C^{(2)}_{s_2}\cup\{\mathbf v^{(2)}_{s_2}\},
  \ldots,\mathcal C^{(\ell)}_1\cup\{\mathbf v_1^{(\ell)}\},\ldots,\mathcal C^{(\ell)}_{s_\ell}\cup\{\mathbf v^{(\ell)}_{s_\ell}\})$$ in  $G_0$ for the subspace $\R^\cup\la \mathcal B\ra$ such that the following hold for each $j\in [1,\ell]$.
    \begin{itemize}
    \item[(a)] $\mathcal C^{(j)}_{s_j}=\mathcal C_j$ and  $\bigcup_{k\in [1,s_j]}\mathcal C^{(j)}_k=\darrow \mathcal C_j\setminus \darrow \mathcal C_{j-1}$, whence $\darrow\mathcal B=\bigcup_{\alpha\in [1,\ell]}\bigcup_{\beta\in [1,s_{\alpha}]}\mathcal C^{(\alpha)}_\beta$. In particular, $\mathcal C_\beta^\alpha\subseteq \mathcal X_1\cup \ldots\cup\mathcal X_{s-1}$ for all $\beta$ and $\alpha$ except possibly when $\alpha=\ell$ and $\beta=s_\ell$.
    \item[(b)] $\mathbf v^{(j)}_{s_j}$ is the relative half-space defined by the limit $(u_1,\ldots,u_{r_{j}})$ taking  $\mathbf v^{(j)}_{s_j}(i)=x_i$ (for all sufficiently large $i$), with $-(u_1,\ldots,u_{r_{j+1}-1})$ minimally encased by $\mathcal C_{s_j}^{(j)}=\mathcal C_j$ and $\partial(\{\mathbf v_{s_j}^{(j)}\})=\mathcal C_{j-1}$.
\item[(c)] Each $\mathbf v^{(j)}_k$ with $k<s_j$ is defined by a strict truncation of a  limit associated to some $\y\in \darrow\mathcal C_j$ taking $\mathbf v^{(j)}_{k}(i)=\y(i)$ (for all sufficiently large $i$).
    \end{itemize}
\end{proposition}

\begin{proof}Since $\mathcal B$ is a support set, implying $\mathcal B^*=\mathcal B$, we have  $\darrow \mathcal B\setminus \mathcal B=\darrow \partial(\mathcal B)\subseteq \mathcal X_1\cup \ldots\cup\mathcal X_{s-1}$. Thus, the in particular statement in (a) follows in view of $\mathcal C^\ell_{s_\ell}=\mathcal C_\ell=\mathcal B$ and $\darrow\mathcal B=\bigcup_{\alpha\in [1,\ell]}\bigcup_{\beta\in [1,s_{\alpha}]}\mathcal C^{(\alpha)}_\beta$.

We define a directed graph with vertices $\mathcal X_1\cup \{\mathbf v_1\}\cup\ldots\cup \mathcal X_s\cup\{\mathbf v_s\}$ as follows. Each vertex $\y$ is defined by an associated limit $\vec u_\y=(u^\y_1,\ldots,u_{t_\y}^\y)$ such that $\partial(\{\y\})$ minimally encases $-\vec u_\y^\triangleleft$, and thus via Proposition \ref{prop-orReay-minecase-char} there is a uniquely defined associated sequence of support sets $\emptyset=\mathcal C_0^\y\prec \mathcal C_1^\y\prec\ldots\prec \mathcal C_{\ell_\y}^\y=\partial(\{\y\})\subset \darrow \{ \y\}$ and indices $1=r_1^\y<\ldots<r^\y_{\ell_\y}<r^\y_{\ell_\y+1}=t_\y+1$.
We define a directed edge between $\y$ and each half-space from $\bigcup_{j=1}^{\ell_\y}\mathcal C_j^\y$. Moreover, each $\z\in \bigcup_{j=1}^{\ell_\y}\mathcal C_j^\y$ lies in some $\mathcal C_j^\y$, and associating to $\z$ the minimal index $j$ for which this is true gives a way to  group the neighbors of $\y$ in way that places all vertices from $\mathcal C_1^\y$ before those from $\mathcal C_2^\y\setminus \mathcal C_1^\y$, which then come before those from $\mathcal C_3^\y\setminus (\mathcal C_2^\y\cup \mathcal C_1^\y)$, and so forth.
Since $\mathcal B$ is a support set by hypothesis, each $\mathcal C_j$ for $j\in [1,\ell]$ is also a support set.
If we start with the vertices from $\mathcal B$ and include all their neighbors, followed by all the neighbors of  their neighbors, and continue in this fashion, we obtain an induced subgraph on $\darrow \mathcal B$. We can add a vertex $\y_0$ to this subgraph and connect it to all vertices from $\mathcal B$ to obtain a graph rooted at $\y_0$. Set
 $\vec u_{\y_0}=(u_1,\ldots,u_t)$,  $\y_0(i)=x_i$ for $i\geq 1$, and
 $\mathcal C_j^{\y_0}=\mathcal C_j$ for $j\in [1,\ell]$, which minimally encases $-(u_1,\ldots,u_{r_{j+1}-1})$ by Proposition \ref{prop-orReay-minecase-char}.2.  
We now describe how the sets $\mathcal C^{(j)}_k\cup \{\mathbf v^{(j)}_k\}$ can be constructed via a depth-first search argument rooted at the vertex $\y_0$ which respects the ordering of neighbors described above.

Starting at  $\y_0$, consider a directed path, say with vertex sequence $\y_0,\y_1,\ldots,\y_{\alpha+1}$, such that the vertex $\y_j$ is always chosen from among the minimal available neighbors of $\y_{j-1}$, that is, $\y_j\in \mathcal C_1^{\y_{j-1}}$. Moreover, when possible, choose $\y_j\in \mathcal C_1^{\y_{j-1}}$ with $\partial(\{\y_j\})\neq \emptyset$. Choosing vertices this way ensures that $\y_{\alpha+1}\prec \y_\alpha\prec \y_{\alpha-1}\prec\ldots\prec \y_1$. Suppose some $\y_j$ were a neighbor of some $\y_i$ with $i<j-1$.
 Then $\y_j\prec \y_{i+1}\in \mathcal C_1^{\y_{i}}$ and  $\y_{j}\in \mathcal C_{k}^{\y_i}$ for some $k\in [1,\ell_{\y_i}]$. Since $\mathcal C_1^{\y_i}\preceq \mathcal C_k^{\y_i}$ implies  $\mathcal C_1^{\y_i}\subseteq \darrow\mathcal C_k^{\y_i}$, we see this would mean there  is some $\z\in \mathcal C_k^{\y_i}$ with $\y_j\prec \y_{i+1}\preceq \z$. Hence $(\mathcal C_k^{\y_i})^*\neq\mathcal C_k^{\y_i}$. However, this contradicts that  $\mathcal C_k^{\y_i}$ is a support set.
   So we instead conclude that each $\y_j$ is not in the neighborhood of any $\y_i$ with $i<j-1$. If  we consider such a directed path with non-extendable length, then (in view of our preference for choosing $\y_j$ with $\partial(\{\y_j\})\neq \emptyset$) we must have $\partial(\{\z\})=\emptyset$ for every $\z\in \mathcal C_1^{\y_{\alpha}}$. By definition, $\mathcal C_1^{\y_{\alpha}}$ minimally encases $-u_1^{\y_\alpha}$,
 with $u_1^{\y_\alpha}$  a truncation of the limit $(u_1^{\y_\alpha},\ldots,u_{t_{\y_\alpha}}^{\y_\alpha})$ which defines
 $\y_{\alpha}$. For $\alpha>0$, we have  $\y_{\alpha+1}\in \mathcal C_1^{\y_\alpha}\preceq \partial(\{\y_{\alpha}\})$, thus ensuring  $t_{\y_\alpha}>1$ (if $t_{\mathbf y_\alpha}=1$, then (V1) implies that $\partial(\{\y_\alpha\})=\emptyset$, contradicting that $\y_{\alpha+1}\in \mathcal C_1^{\y_\alpha}\preceq \partial(\{\y_{\alpha}\})$). Hence  $u_1^{\y_\alpha}$  is a \emph{strict} truncation for $\alpha>0$. Define the half-space $\mathbf v^{(1)}_1=\R_+ u_1^{\y_\alpha}$ via this truncated limit
 $u_1^{\y_\alpha}$
 and set $\mathbf v_1^{(1)}(i)=\y_{\alpha}(i)$ for all $i\geq 1$. Note this ensures that $(\mathcal C_1^{(1)}\cup \{\mathbf v_1^{(1)}\})$ is a virtual Reay system in $G_0$, where $\mathcal C_1^{(1)}=\mathcal C_1^{\y_\alpha}$. If $\alpha=0$, then $\y_{\alpha}=\y_0$, and we have $\mathcal C_1^{(1)}=\mathcal C_1$ and $u_1^{\y_\alpha}=u_1=u_{r_1}$. Thus, for $\ell=1$ and $\alpha=0$ (meaning $\partial(\{\z\})=\emptyset$ for all $\z\in \mathcal B=\mathcal C_\ell=\mathcal C_1=\mathcal C_1^{\mathbf y_\alpha}$), the proof is complete. So we can proceed by induction on $|\darrow \mathcal B|$ with the base case when $|\darrow \mathcal B|=1$ complete.

Note that the only way $\mathcal C_1^{(1)}=\mathcal C_1^{\y_\alpha}\subseteq \mathcal X_1\cup\ldots\cup\mathcal X_s$ can fail  is when $\alpha=0$ and $\ell=1$, which was the case covered in the base of the induction. Therefore we may assume $\mathcal C_1^{(1)}\subseteq \mathcal X_1\cup\ldots\cup \mathcal X_s$.
 Let $\mathcal E=\R^\cup\la \mathcal C_1^{(1)}\ra$ and let $\pi:\R^d\rightarrow \mathcal E^\bot$ be the orthogonal projection.
 In view of Proposition \ref{prop-VReay-modulo}, $\pi(\mathcal R)=\big(\mathcal X_i^\pi\cup \{\pi(\mathbf v_i)\}\big)_{i\in J}$ is a virtual Reay system in $\pi(G_0)$. By proposition \ref{prop-orReay-modulo}.6,  $\darrow \mathcal B\subseteq \mathcal X_1\cup\ldots\cup \mathcal X_s\cup \{\mathbf v_j:\;j\in J\}$ with $\mathcal B^\pi$ a support set.

 Let $\mathcal A\subseteq \darrow \mathcal B$ be arbitrary.  Since $\mathcal B\subseteq \mathcal X_1\cup\ldots\cup \mathcal X_s\cup\{\mathbf v_j:\; j\in J\}$ and  $\partial(\{\z\})=\emptyset$ for every $\z\in \mathcal C_1^{\y_\alpha}=\mathcal C_1^{(1)}$,  Propositions \ref{prop-orReay-BasicProps}.9  and \ref{prop-orReay-modulo}.1 ensure that  $\mathcal A^\pi=\pi(\mathcal A\setminus \mathcal C_1^{(1)})$, while  Proposition \ref{prop-orReay-modulo}.2 implies that $\darrow (\mathcal A^\pi)=(\darrow \mathcal A)^\pi$.
 Since $\mathcal B$ is a support set, it follows that $\mathcal A\subseteq \darrow \mathcal B$ is a support set if and only if $\mathcal A^*=\mathcal A$. When this is the case, Proposition \ref{prop-orReay-modulo}.2 implies that $(\mathcal A^\pi)^*=\mathcal A^\pi$, and thus $\mathcal A^\pi=\pi(\mathcal A\setminus \mathcal C_1^{(1)})$ is also a support set since $\mathcal B^\pi$ is a support set with $\mathcal A^\pi\subseteq (\darrow \mathcal B)^\pi=\darrow  \mathcal B^\pi$.

By proposition \ref{prop-infinite-limits-proj}.1, the sufficiently large index terms in $\{\pi(x_i)\}_{i=1}^\infty$ form an asymptotically  filtered sequence with limit $\pi(\vec u)$, and by discarding the first few terms, we can w.l.o.g. assume $\{\pi(x_i)\}_{i=1}^\infty$ is an asymptotically filtered sequence with limit $\pi(\vec u)$.

\subsection*{Case 1: $\alpha=0$}  In this case, $\ell>1$ and $t>1$ (else we fall in the already completed base of the induction),  $\mathcal C_1^{(1)}=\mathcal C_1$ and $u_1,\ldots,u_{r_2-1}\in \mathcal E=\R^\cup\la \mathcal C_1\ra$ (as $-(u_1,\ldots,u_{r_2-1})$ is minimally encased by $\mathcal C_1$ by Proposition \ref{prop-orReay-minecase-char}).
 If $\mathcal C_j\subseteq \mathcal X_1\cup\ldots\cup \mathcal X_s$ with $j\geq 2$, then Proposition \ref{prop-orReay-BasicProps}.9 together with $\mathcal C_{j-1}\prec \mathcal C_j$ ensures that there is some $\y\in \mathcal C_j$ with $\y\nsubseteq\R^\cup \la \mathcal C_{j-1}\ra$ (else $\mathcal C_j\subseteq \darrow \mathcal C_{j-1}$, yielding the contradiction $\mathcal C_j=\mathcal C_j^*\preceq \mathcal C_{j-1}$), and thus with $\pi(\y)\neq \{0\}$ as well, so that $\mathcal C_{j-1}^\pi\prec \mathcal C_j^\pi$ (by Proposition \ref{prop-orReay-modulo}.3). On the other hand, if $\mathcal C_j\nsubseteq \mathcal X_1\cup\ldots\cup \mathcal X_s$, then $j=\ell$ with $\mathcal C_\ell=\mathcal B$ and $\mathcal C_{\ell-1}\subseteq \mathcal X_1\cup \ldots\cup \mathcal X_s$ (as $\mathcal B$ minimally encases $-\vec u$ urbanely). In this case, $\mathcal C_\ell$ contains some $\mathbf v_i$ with $\mathbf v_i\notin \mathcal C_{\ell-1}\subseteq \mathcal X_1\cup\ldots\cup \mathcal X_s$.
  Since $\mathbf v_i\in \darrow \mathcal B\subseteq \mathcal X_1\cup\ldots\cup \mathcal X_s\cup \{\mathbf v_j:\;j\in J\}$, Proposition \ref{prop-orReay-modulo}.1 implies  $\pi(\mathbf v_i)\neq 0$ , and thus $\mathcal C_{\ell-1}^\pi\prec \mathcal C_\ell^\pi$ holds in this case as well (in view of Proposition \ref{prop-orReay-modulo}.3). Therefore, we find that
$\emptyset=\mathcal C_1^\pi\prec \mathcal C_2^\pi\prec \ldots\prec \mathcal C_\ell^\pi$. In such case, Proposition \ref{prop-orReay-minecase-char}.4 implies that  $\mathcal B^\pi$ minimally encases $-\pi(\vec u)$ with $\mathcal F^\pi(\pi(\vec u))=(\overline u_{r_2},\ldots,\overline u_{r_\ell})$, where $\mathcal F(\vec u)=(\overline u_{r_1},\ldots,\overline u_{r_\ell})$ and $\mathcal F=(\R^\cup (\mathcal C_1),\ldots,\R^\cup(\mathcal C_\ell))$. Moreover, $\mathcal C_2^\pi\prec \ldots\prec \mathcal C_\ell^\pi$ are the support sets given by Proposition \ref{prop-orReay-minecase-char}.1 for $\pi(\vec u)$, and $r_2<\ldots<r_\ell<r_{\ell+1}=t+1$ are the indices given by Proposition
\ref{prop-orReay-minecase-char}.1, ensuring that the encasement is urbane as it was for $\mathcal B=\mathcal C_\ell$. By induction hypothesis (and Proposition \ref{prop-orReay-modulo}), there exists a virtual Reay system (found by depth first search) $$\mathcal R''=(\mathcal C^\pi_{1,2}\cup\{\mathbf w_1^{(2)}\},\ldots,\mathcal C^\pi_{s_2,2}\cup\{\mathbf w^{(2)}_{s_2}\},
  \ldots,\mathcal C^\pi_{1,\ell}\cup\{\mathbf w_1^{(\ell)}\},\ldots,\mathcal C^\pi_{s_\ell,\ell}\cup\{\mathbf w^{(\ell)}_{s_\ell}\})$$ in  $\pi(G_0)$ for the subspace $\R^\cup \la \mathcal B^\pi\ra$ satisfying (a)--(c), with each $\mathcal C_{i,j}$ disjoint from $\mathcal C_1^{(1)}$ (so $\mathcal C_{i,j}^\pi=\pi(\mathcal C_{i,j})$), with each $\mathbf w_k^{(j)}$ with $k<s_j$ defined by a strict truncation of the  limit associated to some $\pi(\y)$ with $\y\in \darrow \mathcal C_j\setminus \mathcal C_1^{(1)}$,  and with each $\mathbf w_{s_j}^{(j)}$ defined by the limit $\pi((u_1,\ldots,u_{r_j}))$ with  $\mathcal F^\pi(\pi((u_1,\ldots,u_{r_{j}})))=(\overline u_{r_2},\ldots,\overline u_{r_j})$, and with $\mathcal C_{i,j}\subseteq \mathcal X_1\cup \ldots\cup \mathcal X_s$ for all $i$ and $j$ except possibly when $i=s_\ell$ and $j=\ell$. We proceed to show we can apply Lemma \ref{lemma-VReay-extension}.
 %
  %

  For $\mathbf w_k^{(j)}$ with $k<s_j$, we have the half-space $\mathbf w_k^{(j)}$ defined by a strict truncated limit associated to some $\pi(\y)$ with $\y\in \darrow \mathcal C_j\setminus \mathcal C_1^{(1)}$. If $\vec u_\y$ is the limit which defines $\mathbf y$, then Proposition \ref{prop-infinite-limits-proj}.1 ensures that $\vec u_\y$ has a strict truncation $\vec v_\y$ with $\pi(\vec v_\y)$ the limit defining $\mathbf w_k^{(j)}$, and by choosing such a truncation  of minimal length, we can assume $\pi(\vec v_\y)^\triangleleft=\pi(\vec v_\y^\triangleleft)$.
  Since $\partial(\{\y\})$ encases $-\vec u_\y^\triangleleft$, it also encases the truncation $-\vec v_\y^\triangleleft$, so there must be some $\mathcal B_{k,j}\preceq \partial(\{\y\})\subseteq \darrow \mathcal C_j$ which minimally encases $-\vec v_\y^\triangleleft$, and since $\mathcal B_{k,j}\preceq \partial(\{\y\})\subseteq \mathcal X_1\cup \ldots\cup\mathcal X_{s-1}$, it does so urbanely.
  Thus Proposition \ref{prop-orReay-minecase-char}.4 ensures that $\mathcal B_{k,j}^\pi$ minimally encases $-\pi(\vec v_\y^\triangleleft)=-\pi(\vec v_\y)^\triangleleft$.
   Since $\mathcal B_{k,j}\subseteq\darrow \mathcal C_j$, it follows from (a) holding for $\mathcal R''$ and Proposition \ref{prop-orReay-modulo}.2 that $\mathcal B_{k,j}^\pi\subseteq \mathcal C_{1,2}^\pi\cup \ldots\cup\mathcal C^\pi_{s_j,j}$, in which case $\mathcal B_{k,j}^\pi$ minimally encases $-\pi(\vec v_\y^\triangleleft)=-\pi(\vec v_\y)^\triangleleft$ urbanely in $\mathcal R''$.
  As a result, since the support set $\partial(\{\mathbf w_k^{(j)}\})$ in $\mathcal R''$ also minimally encases $-\pi(\vec v_\y)^\triangleleft$ (urbanely) by (V1), it follows from Proposition \ref{prop-orReay-minecase-char}.1 that $\mathcal B_{k,j}^\pi=\partial(\{\mathbf w_k^{(j)}\})$, and now we must have $\mathcal B_{k,j}^\pi=\partial(\{\mathbf w_k^{(j)}\})\subseteq \mathcal C_{1,2}^\pi\cup \ldots\cup\mathcal C^\pi_{k,j}$ (with $\mathcal B_{k,j}^\pi$  disjoint from $\mathcal C^\pi_{k,j}$), and thus $\mathcal B_{k,j}\subseteq \mathcal C_1^{(1)}\cup \mathcal C_{1,2}\cup \ldots\cup\mathcal C_{k,j}$ (with $\mathcal B_{k,j}$  disjoint from $\mathcal C_{k,j}$).

 Next consider a half-space $\mathbf w_{s_j}^{(j)}$ with $j\in [2,\ell]$, which is defined by the limit $\pi((u_1,\ldots,u_{r_{j}}))$ with $\pi((u_1,\ldots,u_{r_{j}}))^\triangleleft=\pi((u_1,\ldots,u_{r_{j}-1}))=
  \pi((u_1,\ldots,u_{r_j})^\triangleleft)$ (since $r_j$ occurs as the index of the last coordinate of $\mathcal F^\pi(\pi((u_1,\ldots,u_{r_j})))=(\overline u_{r_2},\ldots,\overline u_{r_j})$). Since $\mathcal B=\mathcal C_\ell$ minimally encases $-\vec u$, it also encases the truncation $-(u_1,\ldots,u_{r_j-1})$,
 so there must be some $\mathcal B_{s_j,j}\preceq \mathcal B\subseteq \darrow \mathcal C_\ell$ which minimally encases $-(u_1,\ldots,u_{r_j-1})$. Indeed, $\mathcal B_{s_j,j}=\mathcal C_{j-1}\subseteq \mathcal X_1\cup\ldots\cup \mathcal X_s$ by Proposition \ref{prop-orReay-minecase-char}.2, ensuring the encasement is urbane.
 Thus Proposition \ref{prop-orReay-minecase-char}.4 ensures that $\mathcal B_{s_j,j}^\pi$ minimally encases $-\pi((u_1,\ldots,u_{r_j-1}))$.
   Since $\mathcal B_{s_j,j}\subseteq\darrow \mathcal C_\ell$, it follows from (a) holding for $\mathcal R''$ and Proposition \ref{prop-orReay-modulo}.2 that $\mathcal B_{s_j,j}^\pi\subseteq \mathcal C_{1,2}^\pi\cup \ldots\cup\mathcal C^\pi_{s_\ell,\ell}$, in which case  $\mathcal B_{s_j,j}^\pi$ minimally encases $-\pi((u_1,\ldots,u_{r_j-1}))$ urbanely.
  As a result, since the support set $\partial(\{\mathbf w_{s_j}^{(j)}\})$ also minimally encases $-\pi(u_1,\ldots,u_{r_j-1})$ (urbanely) by (V1), it follows from Proposition \ref{prop-orReay-minecase-char}.1 that $\mathcal B_{s_j,j}^\pi=\partial(\{\mathbf w_{s_j}^{(j)}\})$, and now we must have $\mathcal B_{s_j,j}^\pi=\partial(\{\mathbf w_{s_j}^{(j)}\})\subseteq \mathcal C_{1,2}^\pi\cup \ldots\cup\mathcal C^\pi_{s_j,j}$ (with $\mathcal B_{s_j,j}^\pi$  disjoint from $\mathcal C^\pi_{s_j,j}$), and thus $\mathcal B_{s_j,j}\subseteq \mathcal C_1^{(1)}\cup \mathcal C_{1,2}\cup \ldots\cup\mathcal C_{s_j,j}$ (with $\mathcal B_{s_j,j}$  disjoint from $\mathcal C_{s_j,j}$). Additionally, Proposition \ref{prop-VReay-modulo} ensures that $\vec u_{\pi(\x)}=\pi(\vec u_\x)$ for $\x\in \mathcal C_{1,2}\cup\ldots\cup \mathcal C_{s_\ell,\ell}$.

The work of the above two paragraphs allows us to invoke Lemma \ref{lemma-VReay-extension} with $\mathcal R_\mathcal A=(\mathcal C_1^{(1)}\cup \{\mathbf v_1^{(1)}\})$ and $\mathcal R_C=\mathcal R''$, where $\mathbf v_1^{(1)}=\R_+ u_1$: here the $\mathcal B_j$  for $j\in [1,r]$  are taken to be the  $\mathcal C_{k,j}$ with $j\in [2,\ell]$ and  $k\in [1,s_j]$,  the $\mathbf c_j$ and $\mathcal C_j$ for $j\in [1,r]$ are taken to be the $\mathbf w_k^{(j)}$ and $\mathcal B_{k,j}$ with $j\in [2,\ell]$ and $k\in [1,s_j]$, and the $\vec u_j$ for $j\in [1,r]$ are taken to be the $\vec v_\y$ for   $\y\in \mathcal C_{k,j}$ when $k<s_j$ as well as the $(u_1,\ldots,u_{r_j})$ when $k=s_j$, for $j\in [2,\ell]$ (note $\mathcal B_{k,j}$ minimally encasing $-\vec v_\y^\triangleleft$ ensures that $\mathcal C_1^{(1)}\cup \mathcal B_{k,j}$ encases $-\vec v_\y^\triangleleft$). As a result, we find  that
$$\mathcal R'=(\mathcal C_1^{(1)}\cup \{\mathbf v_1^{(1)}\},\mathcal C_{1,2}\cup\{\mathbf v_1^{(2)}\},\ldots,\mathcal C_{s_2,2}\cup\{\mathbf v^{(2)}_{s_2}\},
  \ldots,\mathcal C_{1,\ell}\cup\{\mathbf v_1^{(\ell)}\},\ldots,\mathcal C_{s_\ell,\ell}\cup\{\mathbf v^{(\ell)}_{s_\ell}\})$$ is a virtual Reay system in $G_0$, where each $\pi(\mathbf v^{(j)}_k)=\mathbf w_k^{(j)}$ with $\mathbf v_k^{(j)}$ for $k<s_j$ defined by the limit $\vec v_\y$ (described above) and each $\mathbf v_{s_j}^{(j)}$ defined by the limit $(u_1,\ldots,u_{r_j})$. Since (a) held for $\mathcal R''$, it follows by construction that (a) holds for $\mathcal R'$. Since both $\partial(\{\mathbf v_{s_j}^{(j)}\})$ and $\mathcal C_{j-1}$ minimally encase $-(u_1,\ldots,u_{r_{j}-1})$ urbanely (the first from (V1) for $\mathcal R'$ and the second from Proposition \ref{prop-orReay-minecase-char}.2), Proposition \ref{prop-orReay-minecase-char}.1 ensures that $\partial(\{\mathbf v_{s_j}^{(j)}\})=\mathcal C_{j-1}$, and the remaining parts of (b)--(c) follow directly by construction, completing the case.

  \subsection*{Case 2: $\alpha>0$} In this case, $\y_\alpha\in \darrow \mathcal B$.
   To lighten the notation, we use the following abbreviations for $\beta\in [0,\alpha]$: $\mathcal C_j^{\y_\beta}=\mathcal C_j^\beta$, $\vec u_{\y_\beta}=\vec u_\beta$, $u^{\y_\beta}_j=u^\beta_j$, $\ell_{\y_\beta}=\ell_\beta$, $t_{\y_\beta}=t_\beta$, $r^{\y_\beta}_j=r^\beta_j$ and $u^{\y_\beta}_{r^{\y_\beta}_j}=u^\beta_{r_j}$, etc., for $j\in [1,\ell_\beta]$.
Since $\mathcal C^\beta_{j-1}$ minimally encases the tuple  $-(u^\beta_1,\ldots,u^\beta_{r_j-1})$ (by Proposition \ref{prop-orReay-minecase-char}.2), we have $u^\beta_1,\ldots,u^\beta_{r_j-1}\in \R^\cup \la \mathcal C^\beta_{j-1}\ra$ for any $\beta\in [0,\alpha]$.
We have $\mathcal C_1^\beta\preceq \mathcal C_{\ell_\beta}^\beta=\partial(\{\mathbf y_{\beta}\})\prec \{\mathbf y_{\beta}\}\subseteq \mathcal C_1^{\beta-1}$ for $\beta\in [1,\alpha]$. Thus $\mathcal C_1^{(1)}=\mathcal C_1^\alpha\prec \mathcal C_1^{\alpha-1}\prec\ldots\prec \mathcal C_1^{0}=\mathcal C_1$, ensuring  $\mathcal E=\R^\cup \la \mathcal C_1^{\alpha}\ra\subseteq \R^\cup \la\mathcal C^{\beta}_1\ra$ for $\beta\in [0,\alpha]$.
  Observe that $\mathcal C^{\beta}_1$ minimally encases $-u^{\beta}_1$ by Proposition \ref{prop-orReay-minecase-char}.2, implying that $C_1^\beta\cup \{u_1^{\beta}\}$ is a minimal positive basis modulo $\R^\cup \la \partial(\mathcal C_1^\beta)\ra$ (see the comments after Lemma \ref{lem-orReay-mincase-t=1}).
  For $\beta\in [0,\alpha-1]$, we have  $\mathcal C_1^\alpha\preceq\partial(\{\y_\alpha\})\prec \mathcal C_1^{\alpha-1}\preceq \partial(\{\mathbf y_{\alpha-1}\})\prec\ldots\preceq \partial(\{\mathbf y_{\beta+1}\})$ with $\y_{\beta+1}\in \mathcal C_1^\beta$, ensuring that $\mathcal C_1^\alpha\subseteq \darrow \partial(\mathcal C_1^\beta)$ and $\mathcal E=\R^\cup\la \mathcal C_1^\alpha\ra\subseteq \R^\cup\la \partial(\mathcal C_1^\beta)\ra$. Thus, since  $C_1^\beta\cup \{u_1^{\beta}\}$ is a minimal positive basis modulo $\R^\cup \la \partial(\mathcal C_1^\beta)\ra$, we must have $u_1^\beta\notin \mathcal E$ with $(\mathcal C_1^\beta)^\pi\neq \emptyset$.
  Likewise, for $j\geq 2$ and $\beta\in [0,\alpha]$, we have $\mathcal E\subseteq \R^\cup \la \mathcal C^{\beta}_{1}\ra\subseteq \R^\cup \la \mathcal C^{\beta}_{j-1}\ra$, while  $C_j^\beta\cup \{u_{r_j}^\beta\}$ is a minimal positive basis modulo  $\R^\cup \la \mathcal C^{\beta}_{j-1}\cup \partial(\mathcal C^{\beta}_j)\ra$ (cf. Proposition \ref{prop-orReay-minecase-char} and the comments after Lemma \ref{lem-orReay-mincase-t=1}).
  In particular, $u^{\beta}_{r_j}\notin \mathcal E+\R^\cup \la \mathcal C^{\beta}_{j-1}\ra=\R^\cup\la \mathcal C_{j-1}^\beta\ra$.
   If $\mathcal C_j^\beta$ contains some $\mathbf v_i$, then $\mathbf v_i\in \mathcal C_j^\beta\subseteq\darrow\mathcal B\subseteq \mathcal X_1\cup\ldots\cup \mathcal X_s\cup \{\mathbf v_j:\;j\in J\}$ ensures that $\{0\}\neq \pi(\mathbf v_i)\in \mathcal C_j^\beta$ with $\mathbf v_i\notin \mathcal C_{j-1}^\beta\subseteq \mathcal X_1\cup \ldots\cup \mathcal X_s$. Otherwise, if every $\y\in \mathcal C_j^\beta$ has $\y\subseteq \mathcal E+\R^\cup \la \mathcal C_{j-1}^\beta\ra=\R^\cup \la\mathcal C_{j-1}^\beta\ra$, then Proposition \ref{prop-orReay-BasicProps}.9 would imply $\mathcal C_j^\beta=(\mathcal C_j^\beta)^*\preceq \mathcal C_{j-1}^\beta$, contradicting that $\mathcal C_{j-1}^\beta\prec \mathcal C_j^\beta$.
   In either case, Proposition \ref{prop-orReay-modulo}.3 now ensures that $(\mathcal C_{j-1}^\beta)^\pi\prec (\mathcal C_j^\beta)^\pi$.
   In summary, the above works shows that, for every $j\geq 1$ and $\beta\in [0,\alpha-1]$, we have $u^{\beta}_{r_j}\notin \mathcal E+\R \la u^{\beta}_1,\ldots,u^{\beta}_{r_j-1}\ra$ and $ (\mathcal C_{j-1}^\beta)^\pi\neq (\mathcal C_j^\beta)^\pi$.

  As a result, Proposition \ref{prop-orReay-minecase-char}.4 implies that $(\mathcal C^\beta_{\ell_{\beta}})^\pi$ minimally encases $-\pi(\vec u_{\beta})$ urbanely with $\mathcal F^\pi(\pi(\vec u_{\beta}))=(\overline u^{\beta}_{r_1},\ldots,\overline u^{\beta}_{r_{\ell_{\beta}}})$, for each $\beta\in [0,\alpha-1]$. Moreover, $(\mathcal C^\beta_1)^\pi\prec \ldots\prec (\mathcal C_\ell^\beta)^\pi$ are the support sets given by Proposition \ref{prop-orReay-minecase-char}.1 for $\pi(\vec u_\beta)$, and $r^\beta_1<\ldots<r^\beta_{\ell_\beta}<r^\beta_{\ell_\beta+1}=t_\beta+1$ are the indices given by Proposition
\ref{prop-orReay-minecase-char}.1, ensuring that the encasement is urbane as it was for $\mathcal C_\ell^\beta$.  Additionally, since (as remarked at the start of the proof) no element of $\mathcal C^\alpha_{1}$ is in the neighborhood of any $\y_{\beta}$ with $\beta<\alpha$ (as any element of $\mathcal C_1^\alpha$ may be taken as $\y_{\alpha+1}$), it follows that each $\mathcal C_j^\beta$ is disjoint from $\mathcal C_1^\alpha=\mathcal C_1^{(1)}$ for $\beta\in [0,\alpha-1]$, and thus $\pi(\mathcal C_j^\beta)=(\mathcal C_j^\beta)^\pi$.


 The case $\beta=0$ above tells us $\mathcal C_\ell^\pi=\mathcal B^\pi$ minimally encases $-\pi(\vec u)$ urbanely. By induction hypothesis (and Proposition \ref{prop-orReay-modulo}), there exists a virtual Reay system (found by depth first search) $$\mathcal R''=(\mathcal C^\pi_{2,1}\cup\{\mathbf w_2^{(1)}\},\ldots,\mathcal C^\pi_{s_1,1}\cup\{\mathbf w^{(1)}_{s_1}\},
  \ldots,\mathcal C^\pi_{1,\ell}\cup\{\mathbf w_1^{(\ell)}\},\ldots,\mathcal C^\pi_{s_\ell,\ell}\cup\{\mathbf w^{(\ell)}_{s_\ell}\})$$ in  $\pi(G_0)$ for the subspace $\R^\cup \la \mathcal B^\pi\ra$ satisfying (a)--(c), with each $\mathcal C_{i,j}$ disjoint from $\mathcal C_1^{(1)}$ (so $\mathcal C_{i,j}^\pi=\pi(\mathcal C_{i,j})$), with each $\mathbf w_k^{(j)}$ with $k<s_j$ defined by a strict truncation of the  limit associated to some $\pi(\y)$ with $\y\in \darrow \mathcal C_j\setminus \mathcal C_1^{(1)}$,  with each $\mathbf w_{s_j}^{(j)}$ defined by the limit $\pi((u_1,\ldots,u_{r_j}))$ with  $\mathcal F^\pi(\pi((u_1,\ldots,u_{r_{j}})))=(\overline u_{r_1},\ldots,\overline u_{r_j})$, and with $\mathcal C_{i,j}\subseteq \mathcal X_1\cup \ldots\cup \mathcal X_s$ for all $i$ and $j$ except possibly when $i=s_\ell$ and $j=\ell$.
  In view of Proposition \ref{prop-orReay-modulo}.2 and the work above (which holds for all $\beta<\alpha$), the path $\y_0,\y_1,\ldots,\y_\alpha$ remains a path modulo $\mathcal E$ with the ordering of vertices in each neighborhood preserved. If $\ell_\alpha>1$ (so $\mathcal C_1^{(1)}=\mathcal C_1^\alpha\neq \mathcal C_{\ell_\alpha}^\alpha$), then $\partial(\{\pi(\y_\alpha)\})$ remains nonempty, and we can assume we used some continuation of the path $\y_0,\y_1,\ldots,\y_\alpha$ when constructing $\mathcal R''$ via a depth-first search. If $\ell_\alpha=1$, then $\partial(\{\pi(\y_\alpha)\})=\emptyset$, and we may need to use some other $\y'_\alpha\in \mathcal C_1^{\alpha-1}$ with $\partial(\{\pi(\y'_\alpha)\})\neq \emptyset$ (assuming such $\y'_\alpha$ exists) in place of $\y_\alpha$ when constructing $\mathcal R''$. Nonetheless, in either case, we can assume we used some continuation of the path $\y_0,\y_1,\ldots,\y_{\alpha-1},\y'_\alpha$, where $\y'_\alpha\in \mathcal C_1^{\alpha-1}$, when constructing $\mathcal R''$ via  a depth-first search.

We can apply Lemma \ref{lemma-VReay-extension}, taking $\mathcal R_\mathcal A=(\mathcal C_1^{(1)}\cup \{\mathbf v_1^{(1)}\})$ and $\mathcal R_C=\mathcal R''$, where $\mathbf v_1^{(1)}=\R_+ u^\alpha_1$, as can be seen by a near identical argument to that used in Case 1, the only difference being that the filtered limits begin with $\overline u_{r_1}$ instead of $\overline u_{r_2}$, and the $\mathbf c_j$ for $j\in [1,r]$ correspond to the $\mathbf w_{k}^{(j)}$ for $(j,k)\in \big([1,\ell]\times [1,s_j]\big)\setminus \{(1,1)\}$. Thus
$$\mathcal R'=(\mathcal C_1^{(1)}\cup \{\mathbf v_1^{(1)}\},\mathcal C_{2,1}\cup\{\mathbf v_2^{(1)}\},\ldots,\mathcal C_{s_1,1}\cup\{\mathbf v^{(1)}_{s_1}\},
  \ldots,\mathcal C_{1,\ell}\cup\{\mathbf v_1^{(\ell)}\},\ldots,\mathcal C_{s_\ell,\ell}\cup\{\mathbf v^{(\ell)}_{s_\ell}\})$$ is a virtual Reay system in $G_0$, where each $\pi(\mathbf v^{(j)}_k)=\mathbf w_k^{(j)}$ with $\mathbf v_k^{(j)}$ for $k<s_j$ defined by the limit $\vec v_\y$ (described above) and each $\mathbf v_{s_j}^{(j)}$ defined by the limit $(u_1,\ldots,u_{r_j})$. Recall that $\mathbf v_1^{(1)}$ is defined by $u_1^{\y_\alpha}$, which is a strict truncation of $\vec u_{\y_\alpha}$ when $\alpha>0$ (as remarked at the start of the proof) with $\y_\alpha\in \mathcal C_1^{\alpha-1}\subseteq \darrow \mathcal C_1$. As a result, since (a) held for $\mathcal R''$, it follows by construction (and Proposition \ref{prop-orReay-modulo}.2) that (a) holds for $\mathcal R'$, and (b)--(c) follow directly by construction as argued in Case 1, completing the case and proof.
\end{proof}

\begin{lemma}\label{lemma-minencase-mod}
Let $\vec u=(u_1,\ldots,u_s)$ a tuple of $s\geq 0$ orthonormal vectors in $\R^d$, let $t\in [0,s]$ be an index, let $X\subseteq \R^d$ be a subset which minimally encases $-(u_1,\ldots,u_t)$, and  let $\pi:\R^d\rightarrow \R\la X\ra^\bot$ be the orthogonal projection.
\begin{itemize}
\item[1.] $\R\la X\ra=\C(X\cup\{u_1,\ldots,u_t\})$.
\item[2.] If $Y\subseteq \R^d$ is subset such that 
 $\pi(Y)$ minimally encases $-\pi(\vec u)$ with $|\pi(Y)|=|Y|$, then $X\cup Y$ minimally encases $-(u_1,\ldots,u_s)$.
\end{itemize}
\end{lemma}

\begin{proof} 1. Since $X$ minimally encases $-(u_1,\ldots,u_t)$, Proposition \ref{prop-min-encasement-char} implies there is a Reay system  $\mathcal R=(X_1\cup \{u_{r_1}\},\ldots,X_\ell\cup \{u_{r_\ell}\})$ with $X=\bigcup_{i=1}^\ell X_i$ such that $\mathcal F=(\mathcal E_1,\ldots,\mathcal E_\ell)$ is a compatible filter for $-(u_1,\ldots,u_t)$ having associated indices $1=r_1<\ldots<r_\ell$, where $\mathcal E_j=\R\la X_1\cup \ldots\cup X_j\ra$ for $j\in [0,\ell]$. Since $\mathcal F$ is a compatible filter, we have   $u_i\in \mathcal E_\ell=\R\la X\ra$ for all $i$, and now applying Proposition  \ref{prop-reay-basis-properties}.1 to $\mathcal R$ yields $\C(X\cup\{u_{r_1},\ldots,u_{r_\ell}\})=\R\la X\ra=\C(X\cup\{u_1,\ldots,u_t\})$.


2. Since $X$ minimally encases $-(u_1,\ldots,u_t)$, it follows by Propositions \ref{prop-min-encasement-char} and \ref{prop-reay-basis-properties}.1 that $X$ is linearly independent. Since $\pi(Y)$ minimally encases $-\pi(\vec u)$ with $|\pi(Y)|=|Y|$, it likewise follows that $\pi(Y)$ and $X\cup Y$ are both linearly independent.

Let us first show that it suffices to know $X\cup Y$ simply encases $-\vec u$ to complete the proof. Assuming this is the case, there will be a subset $Z\subseteq X\cup Y$ that \emph{minimally} encases $-\vec u=-(u_1,\ldots,u_s)$. We need to show $Z=X\cup Y$. Since $Z$ encases   $-(u_1,\ldots,u_s)$, it also encases the truncation $-(u_1,\ldots,u_t)$.
Hence, since the elements of $Z\subseteq X\cup Y$ are linearly independent and $u_i\in \R\la X\ra$ for $i\leq t$ (as $X$ encases $-(u_1,\ldots,u_t)$), it follows that $Z\cap X\subseteq X$ encases $-(u_1,\ldots,u_t)$. But since $X$ \emph{minimally} encases $-(u_1,\ldots,u_t)$ by hypothesis, this is only possible if $Z\cap X=X$, that is, $X\subseteq Z$. Next, since $Z$ encases $-\vec u$,  and since $X\cup Y$ is linearly independent with $\ker \pi=\R\la X\ra$, it follows that $\pi(Z)\setminus\{0\}=\pi(Z\setminus X)\subseteq \pi(Y)$ encases $-\pi(\vec u)$ (cf. the equivalent version of encasement mentioned immediately after its definition and the commentary regarding equivalent tuples found there). But since $\pi(Y)$ \emph{minimally} encases $-\pi(\vec u)$ by hypothesis, this is only possible if $\pi(Y)=\pi(Z\setminus X)$, in turn implying $Z\setminus X=Y$ as $\pi$ is injective on $Y$ by hypothesis. Combined with $X\subseteq Z$, it follows that  $Z=X\cup Y$, as desired. Thus it remains to show $X\cup Y$ simply encases $-\vec u$ to complete the proof, as the above argument shows  any such encasement will necessarily be minimal.

If $t=s$, then $\pi(\vec u)$ is the empty tuple, in which case we must have $Y=\emptyset$ (as $\pi(Y)$ minimally encases $-\pi(\vec u)$), and now Item 2 is trivial.

If $t=s-1$, then $X$ minimally encases $-\vec u^\triangleleft =-(u_1,\ldots,u_{s-1})$, while $\pi(Y)$ encases $-\pi(\vec u)=-\pi(u_s)/\|\pi(u_s)\|$, and thus also $-\pi(u_s)$, meaning $-\pi(u_s)\in \C(\pi(Y))$ (a fact trivially true if $\pi(u_s)=0$). Hence $u_s+a\in -\C(Y)$ for some $a\in \ker \pi =\R\la X\ra=\C(X\cup \{u_1,\ldots,u_{s-1}\})$, with the final equality by Item 1, implying that $(u_s+\C(u_1,\ldots,u_{s-1}))\cap -\C(X\cup Y)\neq \emptyset$. Thus, since $X\subseteq X\cup Y$ encases $-(u_1,\ldots,u_{s-1})$, it now follows that $X\cup Y$ encases $-(u_1,\ldots, u_s)$, completing the proof as noted above. So we may assume $t\leq s-2$ and proceed by induction on $s-t$.

Let $Y'\subseteq Y$ be a subset such that $\pi(Y')$ minimally encases $-\pi((u_1,\ldots,u_{s-1}))$. Then $X\cup Y'$ minimally encases $-(u_1,\ldots,u_{s-1})$ by induction hypothesis. Letting $\tau:\R^d\rightarrow \R\la X\cup Y'\ra^\bot$ be the orthogonal projection, we find that $\tau (Y)\setminus \{0\}=\tau(Y\setminus Y')$ is a linearly independent set of size $|Y\setminus Y'|$ (since $X\cup Y$ is linearly independent with $X\cup Y'\subseteq X\cup Y$). Since $\pi(Y)$ encases $-\pi(\vec u)$, it follows that $\tau(\pi(Y))\setminus\{0\}=\tau (Y\setminus Y')$ encases $-\tau(\pi(\vec u))=-\tau(\vec u)$ (cf. the equivalent version of encasement mentioned immediately after its definition and the commentary regarding equivalent tuples found there). Thus there must be some subset $Z\subseteq Y\setminus Y'$ such that $\tau(Z)$ \emph{minimally} encases $-\tau(\vec u)$. Since $\tau(Y\setminus Y')$ is a linearly independent set of size $|Y\setminus Y'|$, it follows that $\tau$ is injective on $Y\setminus Y'$, and thus also on $Z\subseteq Y\setminus Y'$. As a result, we can apply the base of the induction using $X\cup Y'$, $Z$ and $\tau$ (in place of $X$, $Y$ and $\pi$, with $t=s-1$) to conclude that $X\cup Y'\cup Z\subseteq X\cup Y$ minimally encases $-\vec u$, implying that $X\cup Y$ encases $-\vec u$, which completes the proof as shown above.
\end{proof}

We now make precise the idea that the underlying oriented Reay system to a virtual Reay system is a limit structure.
Suppose for each tuple $k=(i_1,\ldots,i_t)\in \Z^t$  of indices $i_j\geq 1$ we have a set $X_k\subseteq \R^d$. Then we write $\lim_{k\rightarrow \infty} X_k=X$ if $X_k\subseteq X$ holds once  all $i_j$ are sufficiently large and, for every $x\in X$, there exists some bound $N_x>0$ such that $x\in X_k$ for all tuples $k$ having all coordinates $i_j\geq N_x$ for $j\in [1,t]$. We say that $\lim_{k\rightarrow \infty} X_k=X$ \textbf{order uniformly} if (additionally) there is a global constant  $N>0$ such that, for every integer $m\geq N$, there is a relative constant $N_m\geq N$ such that, for every tuple $k=(i_1,\ldots,i_t)$ of indices with all coordinates $i_j\geq N$ and $i_\alpha=m$ for some $\alpha\in [1,t]$ and every tuple $k'=(i'_1,\ldots,i'_t)$ of indices with $i'_j=i_j$ for $j\in [1,t]\setminus \{\alpha\}$ and $i'_\alpha\geq N_m$, we have $X_k\subseteq X_{k'}$. (The terms global and relative are introduced above to later make referencing these constants easier.) When this is the case, we can find an increasing subsequences of indices $N\leq \iota_1<\iota_2<\ldots$  so that, given any tuples $k=(\iota_{\alpha_1},\ldots,\iota_{\alpha_t})$ and $k'=(\iota_{\beta_1},\ldots,\iota_{\beta_t})$ having all coordinates equal except one which is larger in $k'$, we have $X_k\subseteq X_{k'}$.   Indeed, we could take $\iota_1=N$, $\iota_2=\max\{N_N,\,\iota_1+1\}$, $\iota_3=\max\{N_{(N_N)},\,\iota_2+1\}$, and so forth. However, it is easily derived from this property that, if $k=(\iota_{\alpha_1},\ldots,\iota_{\alpha_t})$ and $k'=(\iota_{\beta_1},\ldots,\iota_{\beta_t})$ are any pair of tuples with $\iota_{\beta_j}\geq \iota_{\alpha_j}$ for all $j\in [1,t]$, then $X_{k}\subseteq X_{k'}$. In such case, the natural partial order on the tuples $k=(\iota_{\alpha_1},\ldots,\iota_{\alpha_t})$ corresponds to set theoretic inclusion among the sets indexed by $k$. (In other words, these additional properties are  obtained by passing to appropriate subsequences.)
 Additionally, if $k=(i_1,\ldots,i_t)$ is any tuple with $i_j\geq N$ for all $j$, then there exists a constant $N_k$ such that $X_k\subseteq X_{k'}$ for any tuple $k'=(i'_1,\ldots,i'_t)$ with $i'_j\geq N_k$ for all $j$. Indeed, we may simply take $N_k=\max\{N_{i_1},\ldots,N_{i_t}\}$, and this will then follow from the definition of order uniform limit.

\begin{proposition}\label{prop-VReay-RepBasics} Let $G_0\subseteq \R^d$ be a subset and  let  $\mathcal R=(\mathcal X_1\cup \{\mathbf v_1\},\ldots,\mathcal X_s\cup\{\mathbf v_s\})$ be a virtual Reay system in  $G_0$.
 There is a bound $N>0$ such that the following hold for all  $\x\in \mathcal X_1\cup\{\mathbf v_1\}\cup\ldots\cup\mathcal X_1\cup\{\mathbf v_s\}$.

\begin{itemize}
\item[1.]
$\darrow \wtilde{\partial(\{x\})}(k)$ minimally encases $-\vec u_\x^\triangleleft$ for any tuple $k=(i_\z)_{\z\in \darrow \partial(\{\x\})}$ with all $i_\z\geq N$.

\item[2.]
$\lim_{ k\rightarrow \infty}\C\big((\darrow \tilde \x)(k)\big)=\x$, \ $\lim_{ k\rightarrow \infty}\C^\circ\big((\darrow \tilde \x)(k)\big)=\x^\circ$ and $\lim_{ k\rightarrow \infty}\Big(\mathsf C\big((\darrow \tilde \x)( k)\big)\cap \partial(\x)\Big)=\lim_{k\rightarrow \infty}\C(\darrow \wtilde{\partial(\{x\})}(k))=\partial(\x)\cap \x$, with all limits holding order uniformly.

\item[3.] $\darrow\wtilde{\partial(\{x\})}(k)\cup \{\tilde \x^{(t-1)}(i)\}$ is a minimal positive basis when $t>1$,  where $\vec u_\x=(u_1,\ldots,u_t)$, for any tuple $k=(i_\z)_{\z\in \darrow \partial(\{\x\})}$ with all $i_\z\geq N$ and any $i\geq N$.
\end{itemize}
\end{proposition}

\begin{proof}
Let
$\x\in \mathcal X_1\cup\{\mathbf v_1\}\cup\ldots\cup\mathcal X_1\cup\{\mathbf v_s\}$ be arbitrary,
let $\vec u_\x=(u_1,\ldots,u_t)$,  let $\mathcal B_\x=\partial(\{\x\})$  and let
 $\emptyset=\mathcal C_0\prec \mathcal C_1\prec\ldots\prec \mathcal C_\ell=\mathcal B_\x$ and $1=r_1<\ldots<r_\ell<r_{\ell+1}=t$ be the support sets and indices given by Proposition \ref{prop-orReay-minecase-char} applied to urbane minimal encasement of $-\vec u_\x^\triangleleft$ by $\mathcal B_\x=\partial(\{\x\})$.
 For $\y\in \mathcal B_\x=\partial(\{\x\})$, let $\mathcal B_\y=\partial(\{\y\})$ and $\vec u_\y=(u^\y_1,\ldots,u_{t_\y}^\y)$.
 Let $\mathcal B'_\x\subseteq \mathcal B_\x$ be the subset of all half-spaces $\y\in \mathcal B_\x$ with $\partial(\{\y\})\neq \emptyset$.
 Let $\hat k=(i_\y)_{\y\in \darrow \x}$ be a tuple of sufficiently large indices $i_\y\geq N_\x$ (the constant $N_\x>0$ will be determined during the course of the proof) and let $k=(i_\y)_{\y\in \mathcal \darrow \mathcal B_\x}$ be the sub-tuple consisting of all indices $i_\y$ with $\y\in \darrow \mathcal B_\x=\darrow \x\setminus\{\x\}$. For $\y\in \mathcal B_\x$, let $k_\y=(i_\z)_{\z\in \darrow\mathcal B_\y}$
 be the sub-tuple of $k$ consisting of coordinates indexed by $\darrow \mathcal B_\y$, and let $\hat k_\y=(i_\z)_{\z\in\darrow \y}$ be the sub-tuple of $k$ consisting of coordinates indexed by $\darrow \y=\darrow \mathcal B_\y\cup \{\y\}$.
 Since each $\mathcal C_j$ is a support set and replacing each half-space in an oriented Reay system with a representative yields an ordinary Reay system, Propositions \ref{prop-reay-basis-properties}.1  implies that $\darrow \tilde C_j(k')$ is a basis for $\R^\cup \la \mathcal C_j\ra=\R\la \darrow \tilde C_j(k')\ra$ for any tuple $k'$ indexed by $\darrow\mathcal C_{j}$ and $j\in [0,\ell]$.
 Likewise, $(\darrow \tilde \x)(\hat k)$ is a linearly independent set. Since there are only a finite number of $\x\in \mathcal X_1\cup\{\mathbf v_1\}\cup \ldots\mathcal X_s\cup \{\mathbf v_s\}$, if we can show Items 1--3 hold for our arbitrary $\x$ using some bound $N_\x$, with $N_\x$ also being the global constant required in the definition of order uniform limits,  then the proposition will follow by taking $N=\max_\x N_\x$.


We  first handle the case when $\partial(\{\z\})=\emptyset$ for all $\z\in
\mathcal X_1\cup\{\mathbf v_1\}\cup\ldots\cup\mathcal X_1\cup\{\mathbf v_s\}$. In this case, $t=1$, $\mathcal B_\x=\emptyset$ and $\C\big((\darrow \tilde \x)(k)\big)=\x$ for all tuples $k$, in which case Items 1--3 hold trivially with $N=1$. In particular, this completes the case $s=1$, and thus also the case $\dim \R^\cup \la \mathcal X_1\cup\ldots\cup \mathcal X_s\ra=\dim\R^\cup \la \mathcal X_1\ra=1$. We proceed by induction on $\dim \R^\cup \la \mathcal X_1\cup\ldots\cup \mathcal X_s\ra$.  Since $\dim \R^\cup \la\mathcal X_1\cup\ldots\cup \mathcal X_{s-1}\ra<\dim\R^\cup\la \mathcal X_1\cup\ldots\cup \mathcal X_s\ra$, it suffices by induction hypothesis and the base case to consider $\x\in \mathcal X_s\cup\{\mathbf v_s\}$ with $s\geq 2$. If $t=1$, then $\partial(\{\x\})=\emptyset$, and Items  1--3 hold trivially as before, so we can assume $t\geq 2$, and thus also $\ell\geq 1$.
\medskip

1. We first show  Item 1 holds, which we do in two cases based on whether $\ell>1$ or $\ell=1$.

\subsection*{Case 1:} Suppose that $\ell>1$. Let $\mathcal E=\R^\cup \la \mathcal C_{\ell-1}\ra$ and let $\pi:\R^d\rightarrow \mathcal E^\bot$ be the orthogonal projection.  By proposition \ref{prop-VReay-modulo}, $\pi(\mathcal R)$ is a virtual Reay system with $\vec u_{\pi(\x)}=\pi(\vec u_\x)$, $\pi(\vec u_\x)^\triangleleft=\pi(\vec u_\x^\triangleleft)$ and $\partial(\{\pi(\x)\})=\partial(\{\x\})^\pi=\mathcal B_\x^\pi$.
By proposition \ref{prop-VReay-modularCompletion}, there exists a virtual Reay system $\mathcal R'=(\mathcal Y_1\cup\{\mathbf w_1\},\ldots,\mathcal Y_{s'}\cup \{\mathbf w_{s'}\})$ in $G_0$ for the subspace $\R^\cup \la \mathcal B_\x\ra=\R^\cup \la \partial(\{\x\})\ra$ such that $\mathcal Y_{s'}=\mathcal B_\x$, \ $\darrow \mathcal B_\x=\bigcup_{i=1}^{s'}\mathcal Y'_i$, \ $\partial(\{\mathbf w_{s'}\})=\mathcal C_{\ell-1}$  and $\mathbf w_{s'}$ is defined by the limit $(u_1,\ldots,u_{r_\ell})$. By the induction hypothesis applied to $\mathcal R'$, $\darrow \tilde C_{\ell-1}(k)$ minimally encases $-(u_1,\ldots,u_{r_\ell-1})$ once all indices in the tuple $k$ are sufficiently large.  By Propositions \ref{prop-orReay-minecase-char}.4 and \ref{prop-orReay-modulo}.2, $\partial(\{\pi(\x)\})=\mathcal B_\x^\pi=\mathcal C_\ell^\pi$ minimally encases $-\pi(\vec u_\x^\triangleleft)=-\pi(\vec u_\x)^\triangleleft$.
Thus, by the induction hypothesis applied to $\pi(\mathcal R)$, $\pi(\darrow \tilde C_\ell(k))\setminus\{0\}$ minimally encases $-\pi(\vec u_\x^\triangleleft)$ once all indices in the tuple $k$ are sufficiently large.
Since $\mathcal C_\ell=\mathcal B_\x=\partial(\{\x\})\subseteq \mathcal X_1\cup\ldots\cup \mathcal X_s$,
Propositions \ref{prop-orReay-BasicProps}.9 and \ref{prop-orReay-modulo}.2 ensure that $\darrow \mathcal C_\ell^\pi=(\darrow \mathcal C_\ell)^\pi=\pi(\darrow \mathcal C_\ell\setminus \darrow \mathcal C_{\ell-1})$.
As noted above,  $\mathcal E=\R^\cup \la\mathcal C_{\ell-1}\ra=\R\la\darrow \tilde C_{\ell-1}(k)\ra$ with $\darrow \tilde C_\ell(k)$ a set of $|\darrow \mathcal C_\ell|$ linearly independent elements, and we have $\darrow \mathcal C_{\ell-1}\subseteq \darrow \mathcal C_\ell$ in view of $\mathcal C_{\ell-1}\prec \mathcal C_\ell$. Consequently, Lemma \ref{lemma-minencase-mod} (applied with $X=\darrow \tilde C_{\ell-1}(k')$ and $Y=(\darrow \tilde C_\ell\setminus \darrow \tilde C_{\ell-1})(k'')$, where $k'$ and $k''$ are the appropriate sub-tuples of $k$ indexed by $\darrow \mathcal C_{\ell-1}$ and $\darrow \mathcal C_\ell\setminus \darrow \mathcal C_{\ell-1}$, respectively, and with $\vec u=\vec u^\triangleleft_\x$ and $t=r_\ell-1$) implies that $\darrow \tilde B_\x(k)=\darrow \tilde C_\ell(k)$ minimally encases $-\vec u_\x^\triangleleft$, as desired. (Note $N_\x$ must be at least the values obtained inductively for $\mathcal R'$ and $\pi(\mathcal R)$.)

\subsection*{Case 2:} Next suppose that $\ell=1$. Then $\mathcal B_\x=\mathcal C_1$ minimally encases $-u_1$ with $u_i\in \R^\cup\la \mathcal B_\x\ra$ for all $i<t$.
In consequence, since $\mathcal B_\x$ is a support set, Lemma \ref{lem-orReay-mincase-t=1} and Proposition \ref{prop-orReay-BasicProps}.4 imply $-u_1\in \C^\cup(\mathcal B_\x)^\circ=\Summ{\y\in \mathcal B_\x}\y^\circ$.
Thus we must have a representation of the form $\Summ{\y\in \mathcal B_\x}b_\y=-u_1$ with each  $b_\y\in \y^\circ$. Applying Item 2 of the induction hypothesis to each $\y\in \mathcal B_\x$, we have $\y^\circ =\lim_{\hat k_\y\rightarrow \infty}\C^\circ\big((\darrow \tilde \y)(\hat k_\y))$, allowing us to assume $b_\y\in \C^\circ((\darrow \tilde \y)(\hat k_\y)))$ so long as all indices from the tuple $\hat k_\y$ are sufficiently large.  (Note $N_\x$ must be at least the value of the constant $N_{b_\y}$ in the definition of $b_\y\in \lim_{\hat k_\y\rightarrow \infty}\C^\circ\big((\darrow \tilde \y)(\hat k_\y)))$, for each $\y\in \mathcal B_\x$, and we do not require order uniformity.)
Consequently, since $\darrow \y=\mathcal B_\y\cup\{\y\}\subseteq \darrow \mathcal B_\x$, we see that $\darrow \tilde B_\x(k)$ encases $-u_1$.
If it does not do so minimally, then there must be some proper subset $\mathcal B\subset \darrow \mathcal B_\x$ such that $\tilde B(k)$ encases $-u_1$.
Let $\mathcal E=\R^
\cup \la\partial(\mathcal B_\x)\ra=\R^\cup\la \partial^2(\{\x\})\ra$, and let $\pi:\R^d\rightarrow \mathcal E^\bot$ and  $\pi^\bot: \R^d\rightarrow \mathcal E$ be the orthogonal projections.
By Proposition \ref{prop-orReay-BasicProps} (Items 3 and 5),  $\pi(\tilde B_\x(k))$ is a linearly independent set of size $|\mathcal B_\x|$ with $-\pi(u_1)\in \Summ{\y\in \mathcal B_\x}\pi(b_\y)\in \C^\cup(\pi(\tilde{\mathcal B}_\x(k)))^\circ$.
Thus Proposition \ref{prop-char-minimal-pos-basis}.4 implies that $\pi(\tilde B_\x(k))\cup \{\pi(u_1)\}$ is a minimal positive basis of size $|\tilde B_\x(k)|+1$, meaning $\pi(\tilde B_\x(k))$ minimally encases $-\pi(u_1)$.  Consequently, we must have $\mathcal B_\x\subseteq \mathcal B$, as otherwise  $\pi(\tilde B(k))\setminus \{0\}$ would be a proper subset of $\pi(\tilde B_\x(k))$ which encases $-\pi(u_1)$, contradicting the minimality of $\pi(\tilde B_\x(k))$. Therefore, there must be some element from $\darrow \mathcal B_\x\setminus \mathcal B_\x$ missing from the proper subset $\mathcal B\subset \darrow \mathcal B_\x$, which ensures that $\mathcal B'_\x\neq\emptyset$ (since $\mathcal B'_\x=\emptyset$ implies $\darrow \mathcal B_\x=\mathcal B_\x$), and thus that $\mathcal E$ is nontrivial.

For $\y\in \mathcal B_\x$, we have  $\tilde \y(i_\y)=(a_{i_\y,1}^\y  u_{1}^\y+\ldots+a^\y_{i_\y,t_\y-1} u_{t_\y-1}^\y+w_{i_\y}^\y)
+b^\y_{i_\y}\overline u_{t_\y}^\y$ with $t_\y\geq 1$,  $u_{1}^\y,\ldots, u_{t_\y-1}^\y,w_{i_\y}^\y\in \R^\cup \la \partial(\{\y\})\ra\subseteq \mathcal E$ and $b_{i_\y}^{\y},\,\|w_{i_\y}^\y\|\in o(a_{i_\y,t_\y-1}^\y$) (as given in \eqref{x(i)-filtered-form}). Thus \be\label{torc}\pi(\tilde \y(i_\y))=b_{i_\y}^\y\pi(\overline u_{t_\y}^\y)\neq 0,\ee in view of Proposition \ref{prop-orReay-BasicProps}.3, and  \ber\label{kinggo}\pi^\bot(\tilde \y(i_\y))&=&
(a_{i_\y,1}^\y  u_{1}^\y+\ldots+a^\y_{i_\y,t_\y-1} u_{t_\y-1}^\y+w_{i_\y}^\y)
+b^\y_{i_\y}\pi^\bot(\overline u_{t_\y}^\y)\\&=&\nn
\tilde \y^{(t_\y-1)}(i_\y)+b_{i_\y}^\y\pi^\bot (\overline u_{t_\y}^\y).\eer
Note $t_\y\geq 2$ holds precisely for those $\y\in \mathcal B'_\x$, while $\tilde \y^{(t_\y-1)}(i_\y)=0$ for all $i_\y$ and $\y\in \mathcal B_\x\setminus \mathcal B'_\x$.

Since $\pi(\tilde B_\x(k))\cup \{\pi(u_1)\}$ is a minimal positive basis of size $|B_\x(k)|+1$, so too is the re-scaled set $\{\pi(\overline u_{t_\y}^\y):\y\in\mathcal B_\x\}\cup \{\pi(u_1)\}$. As a result, there is a unique linear combination  \be\label{mullberry}\Summ{\y\in \mathcal B_\x}\alpha_\y\pi(\overline u_{t_\y}^\y)=-\pi(u_1),\ee and this linear combination has  $\alpha_\y>0$ for all $\y\in \mathcal B_\x$.  Consequently,
in view of  \eqref{kinggo}, \eqref{mullberry} and \eqref{torc},  we have
\be\label{strawberry}\Summ{\y\in \mathcal B_\x}\frac{\alpha_\y}{b^\y_{i_\y}}\tilde \y(i_\y)=-\pi(u_1)+\Summ{\y\in \mathcal B_\x}\Big(\frac{\alpha_\y}{b_{i_\y}^\y}\tilde \y^{(t_\y-1)}(i_\y)+\alpha_\y \pi^\bot(\overline u_{t_\y}^\y)\Big).\ee
Since $\tilde B(k)$ encases $-u_1$, there is a  positive linear combination of the elements from $\tilde B(k)$ equal to $-u_1$. Moreover, considering this linear combination modulo $\mathcal E$, we see (in view of the uniqueness of the linear combination in \eqref{mullberry}) that the coefficient of each $\y\in \mathcal B_\x\subseteq \mathcal B$ must be $\frac{\alpha_\y}{b^\y_{i_\y}}$.
Thus, in view of \eqref{strawberry}, we find that we have a positive linear combination of the elements from $(\tilde B\setminus \tilde B_\x)(k)$ equal to \be\label{elderberry}-z_k:=-\Summ{\y\in \mathcal B'_\x}\frac{\alpha_\y}{b^\y_{i_\y}}\tilde \y^{(t_\y-1)}(i_\y)-\xi,\ee where $\xi=\pi^\bot(u_1)+\Summ{\y\in \mathcal B_\x}\alpha_\y\pi^\bot(\overline u_{t_\y}^\y)\in \mathcal E$ is a fixed element (independent of $i_\y$), i.e., \be\label{blackberry}-z_k\in \C((\tilde B\setminus \tilde B_\x)(k)).\ee

Applying Item 1 of the induction hypothesis to each $\y\in\mathcal B_\x$, we find that $\darrow \tilde B_\y(k_\y)$ minimally encases $-(u_1^\y,\ldots,u^\y_{t_\y-1})$ for each $\y\in \mathcal B_\x$ so long as all coordinates in the tuple $k_\y$ are sufficiently large (recall that $\mathcal B_\y=\partial(\{\y\})$). (Note $N_\x$ must be at least the value obtained inductively for each $\y\in \mathcal B_\x\subseteq \mathcal X_1\cup\ldots\cup \mathcal X_{s-1}$.) For each $\y\in \mathcal B_\x$,
fix one particular tuple $\kappa_\y=(\iota_\z)_{\z\in \darrow \mathcal B_\y}$ such that $\darrow \tilde B_\y(\kappa_\y)$ minimally encases $-(u_1^\y,\ldots,u^\y_{t_\y-1})$. For convenience, we  can assume  that the coordinate $\iota_\z$ is the same among all tuples $\kappa_\y$  that contain a coordinate indexed by $\z\in \darrow \mathcal B_\y$. By Item 2 of the induction hypothesis applied to each $\y\in \mathcal B_\x$, we have $\lim_{k_\y\rightarrow \infty}\C(\darrow \tilde B_\y(k_\y))=\partial(\y)\cap \y$ order uniformly. Consequently, so long as we choose the coordinates in the tuples $\kappa_\y$ sufficiently large, we can be assured that $\C(\darrow \tilde B_\y(\kappa_\y))\subseteq \C(\darrow \tilde B_\y(k_\y))$ for any tuple $k_\y$ with all coordinates sufficiently large. (Note each coordinate $\iota_\z$ in $\kappa_\y$, for $\y\in \mathcal B_\x$ and $\z\in\darrow \mathcal B_\y$, must be at least the value of the global constant from the order uniform limit
$\lim_{k_\y\rightarrow \infty}\C(\darrow \tilde B_\y(k_\y))$ given by the inductive hypothesis applied to $\y\in \mathcal B_\x$, and then $N_\x$ must be at least $\max_{\y\in \mathcal B_\x,\,\z\in\darrow \mathcal B_\y} N_{\iota_\z}$, where $N_{\iota_\z}$ is the relative constant given in the definition of the order uniform limit $\lim_{k_\y\rightarrow \infty}\C(\darrow \tilde B_\y(k_\y))$.)

Observe that $$\bigcup_{\y\in \mathcal B'_\x}\darrow \mathcal B_\y=\bigcup_{\y\in \mathcal B_\x}\darrow \mathcal B_\y=\darrow\partial(\mathcal B_\x)=\darrow \mathcal B_\x\setminus \mathcal B_\x,$$ with the final equality since  $\mathcal B_\x^*=\mathcal B_\x=\partial(\{\x\})$ is a support set.
Thus \be\label{cloverdrag}\mathcal B\setminus \mathcal B_\x\subset \bigcup_{\y\in\mathcal B'_\x}\darrow \mathcal B_\y\quad\und\quad \mathcal E=\Summ{\y\in \mathcal B'_\x}\R^\cup \la \mathcal B_\y\ra=\Summ{\y\in \mathcal B'_\x}\R\la \darrow \tilde B_\y(k_\y)\ra.\ee (Recall that $\mathcal B_\y=\partial(\{\y\})=\emptyset$ for $\y\in\mathcal B_\x\setminus \mathcal B'_\x$.) Hence, since $\xi\in \mathcal E$, we may write $\xi=\Summ{\y\in \mathcal B'_\x}\xi_\y$ with each $\xi_\y\in \R^\cup \la \mathcal B_\y\ra=\R\la \darrow \tilde B_\y(k_\y)\ra$.

Since $b_{i_\y}^\y\in o(a_{i_\y,j}^\y)$, we have $a_{i_\y,j}^\y/b_{i_\y}^\y\rightarrow \infty$, for all $\y\in \mathcal B'_\x$ and $j\in [1,t_\y-1]$. Hence, in view of \eqref{kinggo},
we see that $\{\frac{\alpha_\y}{b^\y_{i_\y}}\tilde \y^{(t_\y-1)}(i_\y)+\xi_\y\}_{i_\y=1}^\infty$ is an asymptotically filtered sequence of terms from $\R^\cup \la\mathcal B_\y\ra=\R^\cup \la \darrow \tilde B_\y(k_\y)\ra$ with limit $(u^\y_1,\ldots,u_{t_\y-1}^\y)$ (once $i_\y$ is sufficiently large), for each $\y\in \mathcal B'_\x$.
Consequently, since $\darrow \tilde B_\y(\kappa_\y)$ minimally encases $-(u_1^\y,\ldots,u_{t_\y-1}^\y)$ (as noted above), it follows from Proposition \ref{prop-min-encasement-minposbasis} that $-\frac{\alpha_\y}{b^\y_{i_\y}}\tilde \y^{(t_\y-1)}(i_\y)-\xi_\y\in\C^\circ (\darrow \tilde B_\y(\kappa_\y)) \subseteq \C^\circ(\darrow \tilde B_\y(k_\y))$ when $i_\y$ is sufficiently large.
(Note this requires $N_\x$ to be at least the value resulting from the application of Proposition \ref{prop-min-encasement-minposbasis} to the \emph{fixed} set $\darrow \tilde B_\y(\kappa_\y)$, for $\y\in \mathcal B_\x$, and thus does not depend on the infinite possible values for   $k_\y$, which is a subtle but important point.)
Hence \eqref{elderberry} yields $-z_k\in \Summ{\y\in \mathcal B'_\x}\C^\circ(\darrow \tilde B_\y(k_\y))$.
As a result, since the elements of $\bigcup_{\y\in \mathcal B'_\x}\darrow \tilde B_\y(k_\y)\subseteq \darrow \tilde C_\ell(k)$ are linearly independent for any tuple $k$, it follows that $-z_k\in \C^\circ(\bigcup_{\y\in \mathcal B'_\x}\darrow \tilde B_\y(k_\y))$, implying that $\bigcup_{\y\in \mathcal B'_\x}\darrow \tilde B_\y(k_\y)$ \emph{minimally} encases $-z_k$. However, since $(\tilde B\setminus \tilde B_\x)(k)\subset \bigcup_{\y\in \mathcal B'_\x} \darrow \tilde B_\y(k_\y)$ is a proper subset (in view of \eqref{cloverdrag}), this contradicts \eqref{blackberry}. So we have now established that $\darrow \tilde B_\x(k)$ minimally encases $-u_1$ when all coordinates in the tuple $k$ are sufficiently large. Thus, since $u_i\in \R^\cup\la \mathcal B_\x\ra=\R\la \darrow \tilde B_\x(k)\ra=\C(\darrow \tilde B_\x(k)\cup \{u_1\})$ for all $i<t$ (the final equality follows by Lemma \ref{lemma-minencase-mod}.1 while the initial inclusion was remarked at the start of Case 2), it follows that $\darrow  \tilde B_\x(k)$ must also minimally encase $-(u_1,\ldots,u_{t-1})$,  and Item 1 is established for Case 2 as well.

2. We next show that Item 2 holds, for which we will implicitly make use of the fact that $(\darrow \tilde \x)(\hat k)$ is always a linearly independent subset (noted at the beginning of the proof).

 Since $\darrow \tilde B_\x(k)$ encases $-(u_1,\ldots,u_{t-1})$, it follows that $u_1,\ldots,u_{t-1}\in \R\la \darrow \tilde B_\x(k)\ra=\R^\cup\la \mathcal B_\x\ra=\partial(\x)$. Thus  $\mathsf C\big((\darrow \tilde \x)(\hat k)\big)$ is contained in the closed half space $\partial(\x)+\R_+ u_t=\overline \x$, and likewise $\mathsf C^\circ\big((\darrow \tilde \x)(\hat k)\big)\subseteq \x^\circ$.
 Moreover, \be\label{whispwillow}\mathsf C\big((\darrow \tilde \x)(\hat k)\big)\cap \partial(\x)=\mathsf C(\darrow \tilde B_\x(k))=\C\big(\bigcup_{\y\in \mathcal B_\x}(\darrow \tilde \y)(\hat k_\y)\big)=\Summ{\y\in\mathcal B_\x}\C\big((\darrow \tilde \y)(\hat k_\y)\big).\ee
 By the induction hypothesis applied to each $\y\in \mathcal B_\x$, it follows that $\y=\lim_{\hat k_\y\rightarrow \infty}\C\big((\darrow \tilde \y)(\hat k_\y)\big)$ order uniformly. Consequently, since $\partial(\x)\cap \x=\C^\cup(\partial(\{\x\}))=\C^\cup(\mathcal B_\x)=\Summ{\y\in \mathcal B_\x}\y$, it now follows in view of \eqref{whispwillow} that
 \be\label{hummingwasp}\lim_{\hat k\rightarrow \infty}\Big(\mathsf C\big((\darrow \tilde \x)(\hat k)\big)\cap \partial(\x)\Big)=\lim_{k\rightarrow \infty}\C(\darrow \tilde B_\x(k))=\partial(\x)\cap \x\ee also order uniformly. It remains to show $\lim_{\hat k\rightarrow \infty}\C^\circ\big((\darrow \tilde \x)(\hat k)\big)=\x^\circ$ order uniformly, as this together with \eqref{hummingwasp} implies that $\lim_{\hat k\rightarrow \infty}\C\big((\darrow \tilde \x)(\hat k)\big)=\x$ order uniformly as well, which will complete the proof. (Note this requires the global constant $N_\x$ in the order uniform limits to be at least the value for the global constant obtained inductively for each $\y\in \mathcal B_\x\subseteq \mathcal X_1\cup \ldots\mathcal X_{s-1}$.)
 Let $$\tilde \x(i_\x)=(a_{i_\x}^{(1)}u_1+\ldots+a_{i_\x}^{(t-1)}u_{t-1}+w_{i_\x})+b_{i_\x}\overline u_t=\tilde\x^{(t-1)}(i_\x)+b_{i_\x}\overline u_t
 $$ for $i_\x\geq 1$ be as given by \eqref{x(i)-filtered-form}. In particular, $w_{i_\x}\in \R^\cup \la\mathcal B_\x\ra$ with $b_{i_\x},\, \|w_{i_\x}\|\in o(a_{i_\x}^{(t-1)})$.

Since $\darrow \tilde B_\x(k)$ minimally encases $-(u_1,\ldots,u_{t-1})$ so long as all coordinates in $k$ are sufficiently large, we can fix one particular tuple $\kappa=(\iota_\y)_{\y\in \darrow \mathcal B_\x}$ such that $\darrow \tilde B_\x(\kappa)$ minimally encases $-(u_1,\ldots,u_{t-1})$. As shown above, $\lim_{k\rightarrow \infty}\C(\darrow \tilde B_\x(k))=\partial(\x)\cap \x$ order uniformly.
 Consequently, so long as we choose the coordinates in the tuple $\kappa$ sufficiently large, we can be assured that $\C(\darrow \tilde B_\x(\kappa))\subseteq \C(\darrow \tilde B_\x(k))$ for any tuple $k$ with all coordinates sufficiently large. (Note each coordinate $\iota_\z$ in $\kappa$, for $\z\in\mathcal \darrow \mathcal B_\x$, must be at least the value of the global constant from the order uniform limit
$\lim_{k\rightarrow \infty}\C(\darrow \tilde B_\x(k))$ obtained above inductively, and then $N_\x$ must be at least $\max_{\z\in\darrow \mathcal B_\x} N_{\iota_\z}$, where $N_{\iota_\z}$ is the relative constant given in the definition of the order uniform limit $\lim_{k\rightarrow \infty}\C(\darrow \tilde B_\x(k))$.)

To show $\lim_{\hat k\rightarrow \infty}\C^\circ\big((\darrow \tilde \x)(\hat k)\big)=\x^\circ$, let  $z=\alpha \overline u_{t}+w\in \x^\circ=(\partial(\x)+\R_+ \overline u_t)^\circ$ be an arbitrary element  with $\alpha>0$ and $w\in \partial(\x)=\R\la \darrow \tilde B_\x(k)\ra$.
Now \be\label{lockgo2}(\alpha/b_{i_\x}) \tilde \x(i_\x)= \alpha \overline u_{t}+(\alpha/b_{i_\x})\tilde \x^{(t-1)}(i_\x)=z-w+(\alpha/b_{i_\x})
\tilde \x^{(t-1)}(i_\x),\ee
and $\{w-(\alpha/b_{i_\x})\tilde \x^{(t-1)}(i_\x)\}_{i_\x=1}^\infty$ is an asymptotically  filtered sequence of terms (once $i_\x$ is sufficiently large) from $\partial(\x)=\R\la \darrow \tilde B_\x(k)\ra$ with limit $-(u_1,\ldots,u_{t-1})$ minimally encased by $(\darrow \tilde B_\x)(\kappa)$, in view of $b_{i_\x}\in o(a_{i_\x}^{(j)})$ for all $j<t$.
Thus, by Proposition \ref{prop-min-encasement-minposbasis}, $w-(\alpha/b_{i_\x})\tilde \x^{(t-1)}(i_\x)\in \C^\circ((\darrow \tilde B_\x)(\kappa))\subseteq \C^\circ((\darrow \tilde B_\x)(k))$ for all sufficiently large $i_\x$ (independent of $k$ since
$\C^\circ((\darrow \tilde B_\x)(\kappa))$ is a fixed set), which combined with \eqref{lockgo2} ensures that $z\in \C^\circ((\darrow\tilde \x)(\hat k))$ for any tuple $\hat k$  with all coordinates sufficiently large.
 This shows that $\lim_{\hat k\rightarrow \infty}\C^\circ\big((\darrow \tilde \x)(\hat k)\big)=\x^\circ$.
 It remains to show the limit holds order uniformly. For this, it suffices, in view of the order uniform limit in \eqref{hummingwasp}, to show that, for  each sufficiently large integer $m\geq N_\x$, there is a bound $N_m$ such that $\C^\circ(\darrow \tilde \x(\hat k))\subseteq \C^\circ(\darrow \tilde \x(\hat k'))$ whenever $\hat k=(i_\z)_{i_\z\in \darrow \x}$ and $\hat k'=(i'_\z)_{i'_\z\in \darrow \x}$ are tuples of sufficiently large indices $i_\z,\,i'_\z\geq N_\x$ with  $$i_\x=m,\quad i'_\x\geq N_m,\quad\und\quad i_\z=i'_\z\quad\mbox{ for $\z\in \darrow \mathcal B_\x$}.$$ (Note this reduction requires the global bound $N_\x$ be at least the global bound for the order uniform limit in \eqref{hummingwasp}.) In particular, we have $k=k'$ as only the entry indexed by  $\x$ differs between $\hat k$ and $\hat{k'}$.

Let $z+\alpha \overline u_{t}\in \x^\circ =(\partial(\x)+\R_+ \overline u_{t})^\circ$ be an arbitrary element, where $\alpha> 0$ and $z\in \partial(\x)$.
Since
 $\overline u_{t}\notin \partial(\x)$, any linear combination of elements from $\darrow \tilde \x(\hat k)$ equal to  $z+\alpha\overline u_t$ must have
 the coefficient of $\tilde \x(i_\x)$ being  $\alpha/b_{i_\x}$. In consequence,  $z+\alpha \overline u_{t}\in \C^\circ(\darrow \tilde \x(\hat k))$ is equivalent to $
z+\alpha\overline u_t-(\alpha/b_{i_\x})\tilde \x(i_\x)=z-(\alpha/b_{i_\x}) \tilde \x^{(t-1)}(i_\x)
\in \C^\circ(\darrow \tilde B_\x(k))$.
This means the elements $z+\alpha \overline u_{t}\in \C^\circ(\darrow \tilde \x(\hat k))$ are those with
\be\label{teest}z-(\alpha/b_{i_\x}) \tilde \x^{(t-1)}(i_\x)\in \C^\circ(\darrow \tilde B_\x( k)),\ee and we simply need to know  $$z-(\alpha/b_{i'_\x}) \tilde \x^{(t-1)}(i'_\x)\in \C^\circ(\darrow \tilde B_\x(k))$$ also holds for all sufficiently large $i'_\x$ (independent of $k=k'$) to show $z+\alpha u_t\in \C^\circ(\darrow \tilde \x(\hat k'))$ (since $i_\z=i'_\z$ for all $\z\in \darrow \mathcal B_\x$). To achieve this, in view of \eqref{teest}  and the fact that $\C(\darrow \tilde B_\x(k))$ is a convex cone, it suffices  to know \be\label{threepiecer}z_{i'_\x}:=(1/b_{i_\x})\tilde \x^{(t-1)}(i_\x)-(1/b_{i'_\x})\tilde \x^{(t-1)}(i'_\x)\in \C(\darrow B_\x(\kappa))\subseteq \C(\darrow \tilde B_\x(k))\ee for all sufficiently large $i'_\x$ (independent of $k$  and $z+\alpha \overline u_t$ but dependent on $m=i_\x$).
  However, each \ber\nn z_{i'_\x}&=&\Sum{j=1}{t-1}\Big(a_{i_\x}^{(j)}/b_{i_\x}-
  a_{i'_\x}^{(j)}/b_{i'_\x}\Big)u_j+(1/b_{i_\x})w_{i_\x}-(1/b_{i'_\x})w_{i'_\x}\\\nn &=&-\Sum{j=1}{t-1}(a_{i'_\x}^{(j)}/b_{i'_\x}) u_j+\Big(\Sum{j=1}{t-1}(a_m^{(j)}/b_m)u_j+(1/b_{m})w_m-(1/b_{i'_\x})w_{i'_\x}\Big)
  \eer
  with $a_{i'_\x}^{(j)}/b_{i'_\x}\rightarrow \infty$ in view of $b_{i'_\x}\in o(a_{i'_\x}^{(j)})$ and  $j<t$
  and with $\|(1/b_{i'_\x})w_{i'_\x}\|\in o(a^{(t-1)}_{i'_\x}/b_{i'_\x})$  (recall that $i_\x=m$ is fixed). Thus  the sequence $\{z_{i'_\x}\}_{i'_\x=N_\x}^\infty$ is an asymptotically filtered sequence of terms from $\R^\cup \la \mathcal B_\x\ra$ with limit $-(u_1,\ldots,u_{t-1})$ (once $i'_\x$ is sufficiently large).  As a result,  since $\darrow B_\x(\kappa)$ minimally encases $-(u_1,\ldots,u_{t-1})$, it follows from Proposition \ref{prop-min-encasement-minposbasis}  that $z_{i'_\x}\in  \C(\darrow\tilde B_\x(\kappa))\subseteq \C(\darrow \tilde B_\x(k))$\ for all sufficiently large $i'_\x\geq N_m$ (independent of $k$ and $z+\alpha\overline u_t$ but dependent on $i_\x=m$), as desired. This establishes \eqref{threepiecer} and completes the proof.
(Note this final step for Item 2 requires the relative constant $N_m$ be at least the value resulting from the application of Proposition \ref{prop-min-encasement-minposbasis} to the \emph{fixed} set $\darrow \tilde B_\x(\kappa)$).

\medskip

3. In view of the main part in Item 1, let $\kappa=(\iota_\z)_{\z\in \darrow \mathcal B_\x}$ be a \emph{fixed} tuple with all $\iota_\z$ sufficiently large that $\darrow \wtilde B_\x(\kappa)$ minimally encases $-(u_1,\ldots,u_{t-1})$.
Moreover, choose the $\iota_\z$ to each be at least the global constant from the order uniform limits given in Item 2.
Then, for any tuple $k=(i_\z)_{\z\in \darrow \mathcal B_\x}$ with all $i_\z$ sufficiently large, we have $\C^\circ(\darrow \tilde B_\x(\kappa))=\Summ{\y\in \mathcal B_\x}\C^\circ(\darrow \wtilde \y(\kappa))\subseteq \Summ{\y\in \mathcal B_\x}\C^\circ(\darrow \wtilde \y(k))=\C^\circ(\darrow \tilde B_\x(k))$ by Item 2 applied to each $\y\in \mathcal B_\x$.
By proposition \ref{prop-min-encasement-minposbasis}, $\darrow \tilde B_k(\kappa)\cup \{\tilde x^{(t-1)}(i_\x)\}$ is a minimal positive basis for all sufficiently large $i_\x$, implying $-\tilde x^{(t-1)}(i_\x)\in \C^\circ(\darrow \tilde B_\x(\kappa))\subseteq \C^\circ(\darrow \tilde B_\x(k))$, and thus ensuring that $\darrow \tilde B_k(k)\cup \{\tilde x^{(t-1)}(i_\x)\}$ is a minimal positive basis for all sufficiently large $i_\x$ (dependent only on the \emph{fixed} set $\darrow \tilde B_k(\kappa)$), which completes the proof.  (Note each coordinate $\iota_z$ for $\kappa$ must be at least the constant $N_\x$ needed to obtain the main part of Item 1 for $\x$ as well as the global constant from the order uniform limits $\lim_{k\rightarrow \infty}\C^\circ(\darrow \tilde{\y}(k))=\y^\circ$ for $\y\in \mathcal B_\x$, and then $N_\x$ must be at least $\max_{\y\in \mathcal B_\x,\,\z\in \darrow \y} N_{\iota_\z}$, where $N_{\iota_\z}$ is the relative constant given in the definition of the order uniform limits $\lim_{k\rightarrow \infty}\C^\circ(\darrow \tilde{\y}(k))=\y^\circ$ for $\y\in\mathcal B_\x$, as well as at least the value given by Proposition \ref{prop-min-encasement-minposbasis} applied to the fixed set $\darrow \tilde B_\x(\kappa)$).
\end{proof}

For a half-space $\x$ from a general virtual Reay system, we only have $\tilde \x(i)$ as a representative for $\x$, which is a truncated approximation of the actual element $\x(i)\in G_0$. The following proposition shows that,  when $G_0$ is a subset of a lattice, there is no need to truncate when each $\y\in \darrow \x$ has $\vec u_\y$ anchored.
We remark that the hypothesis in Proposition \ref{prop-VReay-Lattice} that $\vec u_\y$ be anchored for every $\y\in \darrow \x$ holds automatically when $\mathcal R$ is anchored and $\x\in \mathcal X_1\cup \ldots\cup \mathcal X_s$.

\begin{proposition}\label{prop-VReay-Lattice}
Let $\Lambda\subseteq \R^d$ be a full rank lattice, let $G_0\subseteq \Lambda$ be a subset, and let  $\mathcal R=(\mathcal X_1\cup \{\mathbf v_1\},\ldots,\mathcal X_s\cup\{\mathbf v_s\})$ be a virtual Reay system in  $G_0$. Suppose
 $\x\in \mathcal X_1\cup\{\mathbf v_1\}\cup\ldots\cup\mathcal X_s\cup\{\mathbf v_s\}$ with  $\vec u_\y$ anchored for every $\y\in \darrow \x$. Then $\x(i)=\tilde \x(i)$ for all sufficiently large $i$ and $\vec u_\x^\triangleleft$ is either trivial or fully unbounded. Moreover, if $\vec u_\x=(u_1,\ldots,u_t)$, then $\x(i)-\tilde \x^{(t-1)}(i)=\xi\neq 0$ is constant for all sufficiently large $i$. In particular, if $t=1$, $\x(i)=\xi$ for all sufficiently large $i$.
\end{proposition}

\begin{proof}
Let $\vec u_\x=(u_1,\ldots,u_t)$, where $t\geq 1$, be the limit associated to $\x$, let $\mathcal B_\x=\partial(\{\x\})$ and let $\pi:\R^d\rightarrow \R^\cup\la \mathcal B_\x\ra^\bot$ and $\pi^\bot:\R^d\rightarrow \R^\cup\la \mathcal B_\x\ra$ be the orthogonal projections. Let $$\x(i)=(a_i^{(1)}u_1+\ldots+a_i^{(t-1)}u_{t-1}+w_i)+b_i\overline u_t+y_i$$ be the asymptotically filtered sequence of lattice points  $\x(i)\in G_0\subseteq \Lambda$ with limit $(u_1,\ldots,u_t)$ in the form given by \eqref{x(i)-filtered-form}, where $\overline u_t=\pi(u_t)/\|\pi(u_t)\|$. In particular, $\pi(\overline u_t)=\overline u_t$ is a positive multiple of $\pi(u_t)$, $y_i\in \R\la \x\ra^\bot$, $\|y_i\|\in o(b_i)$, $w_i\in \R^\cup \la \partial(\{\x\})\ra$, and $b_i,\,\|w_i\|\in o(a_i^{(t-1)})$.
We proceed inductively on $s$. In particular, since $\vec u_\y$ is anchored for every $\y\in \darrow \mathcal B_\x=\darrow \partial(\{\x\})\subseteq \darrow \x$, and since $\darrow \mathcal B_\x\subseteq \mathcal X_1\cup\ldots\cup \mathcal X_{j-1}$ when $\x\in \mathcal X_{j}$,  it follows from the induction hypothesis that  $\R^\cup\la \mathcal B_\x\ra=\R \la \darrow \tilde B_\x(k)\ra=\R \la \darrow B_\x(k)\ra$ is a subspace generated by the lattice points $\darrow B_\x(k)\subseteq G_0\subseteq \Lambda$ (this is trivially true when $s=1$, in which case $\pi$ is the identity map and $\darrow \mathcal B_\x=\emptyset$) for any tuple $k$ with all coordinates sufficiently large. Thus, by Proposition \ref{Prop-lattice-homoIm}, $\pi(G_0)\subseteq \pi(\Lambda)$ is a subset of the full rank lattice $\pi(\Lambda)\subseteq \R^\cup \la \mathcal B_\x\ra^\bot$.
Hence $\{\pi(\x(i))\}_{i=1}^\infty$ is a bounded sequence of lattice points from $\R^\cup\la \mathcal B_\x\ra^\bot$ (in view of  $\vec u_\x$ being \emph{anchored} by hypothesis),
meaning there are only a \emph{finite} number of possibilities for the values of  $$\pi(\x(i))=\pi(\tilde \x(i))+\pi(y_i)=b_i\pi(\overline u_{t})+\pi(y_i)=b_i \overline u_{t}+\pi(y_i).$$ Since $b_i$ is bounded (as $\vec u_\x$ is anchored) and $\|y_i\|\in  o(b_i)$, it follows that $\pi(y_i)\rightarrow 0$ and \be\label{lim-help2}\lim_{i\rightarrow\infty} \pi(\x(i))=(\lim_{i\rightarrow \infty} b_i)\overline u_{t}=\xi,\ee for some $\xi\in  \R^\cup\la \mathcal B_\x\ra^\bot$ (recall that $b_i\in \Theta(a_i^{(t)})$ is a convergent sequence by definition of an asymptotically  filtered sequence $\x(i)$). Consequently, since any convergent sequence of lattice points must eventually stabilize, it follows, for all sufficiently large $i$, that $$b_i\overline u_i+\pi(y_i)=\pi(\x(i))=\xi\in \R^\cup\la \mathcal B_\x\ra^\bot\cap \pi(\Lambda).$$

 If $\xi=0$, then $b_i\|\overline u_{t}||-\|\pi(y_i)\|\leq \|\pi(\x(i))\|=0$ for all sufficiently large $i$ (by the triangle inequality), implying $1=\|\overline u_t\|\leq \|\pi(y_i)\|/b_i$ for all sufficiently large $i$. However, $\|\pi(y_i)\|/b_i\rightarrow 0$ in view of  $\|y_i\|\in o(b_i)$, making this impossible. Therefore we conclude that $\xi\neq 0$, in which case the limit of  $b_i$   must be  nonzero. Thus, since $b_i\in o(a_i^{(t-1)})$ when $t\geq 2$, we must have $a_i^{(t-1)}$ being unbounded for $t\geq 2$, meaning $\vec u_\x^\triangleleft=(u_1,\ldots,u_{t-1})$ is either trivial or fully unbounded. Also,  \eqref{lim-help2} now ensures that $\xi$ is a positive multiple of $\overline u_t$, so $\xi=\|\xi\|\overline u_t$ and $\R\la \x\ra=\R^\cup\la \mathcal B_\x\ra+\R \xi$.

 Since $\xi\neq 0$, it follows that the space $\R \overline u_t=\R\xi$ is linearly spanned by the lattice point $\xi\in \pi(\Lambda)$, in which case Proposition \ref{Prop-lattice-homoIm} ensures $\pi'(\Lambda)$ is a lattice, where $\pi':\R^d\rightarrow (\R^\cup\la \mathcal B_\x\ra+\R\xi)^\bot$ is the orthogonal projection.
Now $\pi'(\x(i))=\pi'(y_i)=y_i$, with the latter equality in view of $y_i\in \R\la\x\ra^\bot=(\R^\cup\la \mathcal B_\x\ra+\R\xi)^\bot$, and we also have $\pi'(y_i)\rightarrow 0$ (since $\pi(y_i)\rightarrow 0$). Hence, since any convergent sequence of lattice points must stabilize, we conclude $y_i=\pi'(y_i)=\pi'(\x(i))=0$ for all sufficiently large $i$, in turn implying $\tilde \x(i)=\x(i)$, and all parts of the proposition follow.
\end{proof}

When dealing with a virtual Reay system over a subset of lattice points $G_0$, we will now always  assume, by removing the first few terms, that the representative sequences $\{\x(i)\}_{i=1}^\infty$ for $\x\in \mathcal X_1\cup\{\mathbf v_1\}\cup \ldots\cup \mathcal X_s\cup \{\mathbf v_s\}$ when  $\vec u_\y$ is anchored for all $\y\in \darrow \x$ (in particular, for all $\x\in\mathcal X_1\cup\ldots\cup \mathcal X_s$ when $\mathcal R$ is anchored)
satisfy the conclusion of Proposition \ref{prop-VReay-Lattice} for all $i\geq 1$.
Then, combining Propositions \ref{prop-VReay-Lattice} and \ref{prop-VReay-RepBasics}, we can remove all $\sim$'s from the statement of Proposition \ref{prop-VReay-RepBasics} for such $\x$. Likewise,  if $\mathcal B\subseteq \mathcal X_1\cup \ldots\cup \mathcal X_s$ and   $\mathcal R$ is anchored, then
 Proposition \ref{prop-VReay-Lattice} implies that $\vec u^\triangleleft_\x$ is trivial or fully unbounded for every half-space $\x$ from $\darrow \mathcal B$, in which case Proposition \ref{prop-VReay-modularCompletion} outputs a virtual Reay system $\mathcal R'$ which is anchored, and also purely virtual provided the limit $\vec u$ is fully unbounded.  Thus we gain important simplifications when restricting to virtual Reay systems over a subset of lattice points $G_0\subseteq \Lambda$.

Our next goal is to provide the analog of Proposition \ref{prop-reay-RayAlg} for virtual Reay systems (and thus for oriented Reay systems as well), which will be done with some effort in Proposition \ref{prop-VReay-SupportSet}. However, we first need some lemmas (Lemma \ref{lem-urbane-encasement-guarantee} will be needed in the next section).

\begin{lemma}
\label{lemma-minencase-Rep}
 Let $G_0\subseteq \R^d$ be a subset, where $d\geq 1$, let  $u\in \R^d$ be a unit vector, and  let  $\mathcal R=(\mathcal X_1\cup \{\mathbf v_1\},\ldots,\mathcal X_s\cup\{\mathbf v_s\})$ be a virtual Reay system in  $G_0$. Suppose there is some virtual independent set $\mathcal A\subseteq \mathcal X_1\cup \{\mathbf v_1\}\cup\ldots\cup \mathcal X_s\cup \{\mathbf v_s\}$ which minimally encases $-u$. Then, for any tuple $k=(i_\y)_{\y\in \darrow \mathcal A}$ with all $i_\y$ sufficiently large, $\darrow \tilde A(k)$ minimally encases $-u$ and $\darrow \tilde A(k)\cup \{u\}$ is a minimal positive basis.
\end{lemma}

\begin{proof}
Since $\mathcal A$ is a virtual independent set,
 $\darrow \tilde A(k)$ is linearly independent for any tuple $k$. Let $\z_1,\ldots,\z_\ell\in \mathcal A$ be the distinct half-spaces in $\mathcal A$.
Since $\mathcal A$  minimally encases $-u$, Lemma \ref{lem-orReay-mincase-t=1} and Proposition \ref{prop-orReay-BasicProps}.4 imply  $-u\in \z_1^\circ+\ldots+\z_\ell^\circ=\C^\cup(\mathcal A)^\circ$. Hence, by Proposition \ref{prop-VReay-RepBasics}.2 applied to each $\z_j^\circ$, we find that $-u\in \C^\circ(\darrow\tilde  A(k))$ for any tuple $k$ with all $i_\y$ sufficiently large, and since $\darrow \tilde A(k)$ is linearly independent, it follows that $\darrow \tilde A(k)$ minimally encases $-u$ with $\darrow \tilde A(k)\cup \{u\}$ a minimal positive basis, completing the proof.
\end{proof}

\begin{lemma}
\label{lemma-MinEncase-Induct}
 Let $\mathcal R=(\mathcal X_1\cup\{\mathbf v_1\},\ldots,\mathcal X_s\cup\{\mathbf v_s\})$ be an oriented Reay system in $\R^d$, let $\vec u=(u_1,\ldots,u_t)$ be a tuple of $t\geq 1$ orthonormal vectors in $\R^d$,   let $\mathcal B\subseteq \mathcal X_1\cup\ldots\cup \mathcal X_s$ be a subset minimally encasing $-\vec u^\triangleleft$, let $\pi:\R^d\rightarrow \R^\cup \la \mathcal B\ra^\bot$ be the orthogonal projection, and let $\pi(\mathcal R)=(\mathcal X_j^{\pi}\cup \{\pi(\mathbf v_j)\})_{j\in J}$. Suppose $\mathcal D\subseteq \bigcup_{j\in J}(\mathcal X_j^\pi\cup \{\mathbf v_j\})$ is a virtual independent set that minimally encases $-\pi(\vec u)$.  Then $\mathcal C=(\pi^{-1}(\mathcal D)\cup \mathcal B)^*$ is virtual independent and minimally encases $-\vec u$ urbanely, \  $\mathcal C^{\pi}=\mathcal D$ and $\pi^{-1}(\mathcal D)\subseteq \mathcal C$.
\end{lemma}

\begin{proof}  By Proposition \ref{prop-orReay-modulo}.1, $\pi(\mathcal R)$ is an oriented Reay system. By Proposition \ref{prop-orReay-modulo}.5, $\mathcal C$ is a virtual independent set with $\pi^{-1}(\mathcal D)\subseteq \mathcal C$, and $\mathcal C^\pi=\mathcal D^*=\mathcal D$ follows by Proposition \ref{prop-orReay-modulo}.4 and the fact that $\mathcal D$ is  virtual independent.  Since $\mathcal B\subseteq \darrow \mathcal C$, it follows that $\mathcal B=\mathcal B^*\preceq \mathcal C$ (with the first equality in view of $\mathcal B$ minimally encasing $-\vec u^\triangleleft$). If $\mathcal B=\mathcal C$, then
$\mathcal D=\mathcal C^\pi=\mathcal B^\pi=\emptyset$,  whence $\pi(\vec u)$ must be the empty tuple.
In this case, $\mathcal B=\mathcal C$ minimally encases $-\vec u$ (cf. Proposition \ref{prop-orReay-minecase-char}), with the encasement urbane since $\mathcal B\subseteq \mathcal X_1\cup\ldots\cup \mathcal X_s$, as desired. Therefore, we may assume $\mathcal B\prec \mathcal C$. Thus, since $\mathcal B\subseteq \mathcal X_1\cup\ldots\cup \mathcal X_s$ minimally encases $-\vec u^\triangleleft$, since $\mathcal C^\pi=\mathcal D$ is a virtual independent set which minimally encases $-\pi(\vec u)=-\pi(u_t)/\|\pi(u_t)\|$, and since $\mathcal C\subseteq \mathcal X_1\cup \ldots\cup \mathcal X_s\cup \{\mathbf v_j:\;j\in J\}$ per definition of $\mathcal C$ and $\pi^{-1}(\mathcal D)$,  it follows from Proposition \ref{prop-orReay-minecase-char}.2 that $\mathcal C$ minimally encases $-\vec u$ urbanely, completing the proof.
\end{proof}

\begin{lemma}\label{lem-urbane-encasement-guarantee}
Let $G_0\subseteq \R^d$, where $d\geq 1$, let $\mathcal R=(\mathcal X_1\cup \{\mathbf v_1\},\ldots,\mathcal X_s\cup \{\mathbf v_s\})$ be a purely virtual Reay system over $G_0$, let $\vec u=(u_1,\ldots,u_t)$ be a tuple of orthonormal vectors $u_1,\ldots,u_t\in \R^d$, where $t\geq 1$, and let $\mathcal A\subseteq \mathcal X_1\cup \{\mathbf v_1\}\cup\ldots\cup \mathcal X_s\cup \{\mathbf v_s\}$ be a support set. Suppose $\mathcal A$ minimally encases $-\vec u$ and there is some $\mathbf x\in \mathcal A$ with $\vec u_{\x}$ fully unbounded. Then there exists some $t_0\in [1,t]$ and $\mathcal A'\preceq \mathcal A$ such that $\mathcal A'$ minimally encases $-(u_1,\ldots,u_{t_0})$ \emph{urbanely} and has some $\y\in \mathcal A'$ with $\vec u_{\y}$ fully unbounded.
\end{lemma}

\begin{proof}
 If $t=1$, then the \emph{support} set $\mathcal A$ minimally encases $-u_1$ and  $\mathcal A$ does not \emph{minimally} encase the empty tuple, while the empty set $\emptyset\subseteq \mathcal X_1\cup \ldots\cup \mathcal X_s$ does. Thus $\mathcal A$ minimally encases $-u_1$ urbanely (as the regularity condition required for urbane encasement holds automatically for support sets, as remarked after the definition of urbane encasement), and the lemma follow taking $\mathcal A'=\mathcal A$, \ $t_0=t=1$ and $\y=\x$.
 Therefore we may assume $t\geq 2$ and proceed by induction on $t$.
 Let $t'\in [0,t-1]$ be the minimal index such that $\mathcal A$ does not minimally encase $-(u_1,\ldots,u_{t'})$ and let $\mathcal A'\prec \mathcal A$ be such that $\mathcal A'$ minimally encases $-(u_1,\ldots,u_{t'})$.
 If $\mathcal A'\subseteq \mathcal X_1\cup\ldots\cup \mathcal X_s$, then the \emph{support} set $\mathcal A$ minimally encases $-\vec u$ urbanely, and the lemma follows taking $\mathcal A'=\mathcal A$, \ $t_0=t$ and $\y=\x$.
 Otherwise, $t'\geq 1$ and there must be some $\mathbf v_j\in \mathcal A'$, and since $\mathcal R$ is purely virtual, $\vec u_{\mathbf v_j}$ is fully unbounded. Since $\mathcal A'$ minimally encases $-(u_1,\ldots,u_{t'})$, we must have $(\mathcal A')^*=\mathcal A'$, while $\darrow \mathcal A'
\subseteq \darrow  \mathcal A$
 ensures $\mathcal X_j\subseteq \{\mathbf v_j\}\nsubseteq \darrow \mathcal A'$ for all $j\in [1,s]$  as $\mathcal X_j\subseteq \{\mathbf v_j\}\nsubseteq \darrow \mathcal A$ for the support set $\mathcal A$. In consequence,  $\mathcal A'$ is also a support set. Thus we can apply the induction hypothesis to $\mathcal A'$ to find some $\mathcal A''\preceq \mathcal A'\prec \mathcal A$ and some $t_0\in [1,t']\subseteq [1,t]$ such that $\mathcal A''$ minimally encases $-(u_1,\ldots,u_{t_0})$ urbanely and contains some $\y\in \mathcal A''$ with $\vec u_{\y}$ fully unbounded, as desired.
\end{proof}

Proposition \ref{prop-VReay-SupportSet}.2(c) ensures that $\mathcal A_\x$ is the lift modulo $\partial(\x)$ of the unique support set minimally encasing $-\vec u_\x$ modulo $\partial(\x)$, and   $\mathcal A_j$ is simply the pull-back of this same set union $\x$ with it not mattering which half-space $\x\in X_j\cup \{\mathbf v_j\}$ is used to perform this construction by Proposition \ref{prop-VReay-SupportSet}.1. Furthermore, the other half-spaces $\y\in\mathcal A_j$ behave  symmetrically in this regard, with the exception that we do \emph{not} have $\mathcal A_\y$ and $\mathcal A_\y^{\pi_\y}$ being support sets;  instead, they are  simply  virtual independent sets, which is only slightly weaker.

\begin{proposition} \label{prop-VReay-SupportSet} Let $G_0\subseteq \R^d$ be a subset, where $d\geq 1$, let $\mathcal R=(\mathcal X_1\cup \{\mathbf v_1\},\ldots,\mathcal X_s\cup\{\mathbf v_s\})$ be a virtual Reay system in $G_0$. For  $\x\in \mathcal X_1\cup\{\mathbf v_1\}\cup \ldots\cup \mathcal X_s\cup \{\mathbf v_s\}$,  let $\pi_\x:\R^d\rightarrow \R^\cup \la \partial(\{\x\})\ra^\bot$ denote  the orthogonal projection and  $\mathcal D_\x=\supp_{\pi_\x(\mathcal R)}(-\pi_\x(\vec u_\x))$. Let $j\in [1,s]$.

\begin{itemize}\item[1.] We have $$\mathcal A_j:=\pi_\x^{-1}(\mathcal D_\x)\cup \{\x\}=\pi_{\z}^{-1}(\mathcal D_{\z})\cup \{\z\} \quad\mbox{ for every $\x,\,\z\in \mathcal X_j\cup \{\mathbf v_j\}$}.$$
\item[2.] For every   $\x\in \mathcal X_j\cup \{\mathbf v_j\}$ and every $\y\in \mathcal A_j$, the following hold.
\begin{itemize}
\item[(a)]  $\mathcal X_j\cup \{\mathbf v_j\}\subseteq \mathcal A_j\subseteq \mathcal X_1\cup \{\mathbf v_1\}\cup \ldots\cup \mathcal X_j\cup \{\mathbf v_j\}$ and $\mathcal A_j^*=\mathcal A_j$.
\item[(b)] $\mathcal A_\y:=\big(\mathcal A_j\setminus \{\y\}\cup \partial(\{\y\})\big)^*$ and $\mathcal A_\y^{\pi_\y}=\pi_\y(\mathcal A_j\setminus \{\y\})$ are  virtual independent sets minimally encasing $-\vec u_\y$ and $-\pi_\y(\vec u_\y)$, respectively, with $\mathcal A_j\setminus \{\y\}\subseteq \mathcal A_\y$.
\item[(c)] $\mathcal A_\x=\big(\pi_\x^{-1}(\mathcal D_\x)\cup \partial(\{\x\})\big)^*$ and $\mathcal A_\x^{\pi_\x}=\mathcal D_\x$ are support sets.
\item[(d)] $\mathcal A_j=\pi_\y^{-1}(\mathcal A_\y^{\pi_\y})\cup\{\y\}$ and $\pi_\y^{-1}(\mathcal A_\y^{\pi_\y})=\mathcal A_j\setminus\{\y\}$.

\item[(e)] $\mathcal A_j\subseteq \mathcal X_1\cup\ldots\cup \mathcal X_s\cup \{\mathbf v_i:\; i\in J\}$ and $A_j^\pi$ is a minimal positive basis of size $|\mathcal A_j|$, where $\pi:\R^d\rightarrow \R^\cup \la\partial(\mathcal A_j)\ra^\bot$ is the orthogonal projection and $\pi(\mathcal R)=(\mathcal X_i^\pi\cup \{\pi(\mathbf v_i)\})_{i\in J}$.
\end{itemize}
\item[3.] If $\C(G_0)=\R^d$, then there exists a finite set $Y\subseteq G_0$ such that, for any $\x\in \mathcal A_j$ and any tuple $k=(i_\z)_{\z\in \darrow \mathcal A_\x}$  with all $i_\z$ sufficiently large,  there is a subset $Y_{k}\subseteq Y$ such that
    $$
    ( A_j\setminus \{x\})(k)\cup \darrow  \wtilde{\partial(A_j)}(k) \cup Y_k$$ minimally encases $-\vec u_\x$. Moreover, if $(A_j\setminus \{x\})(k)\subseteq \R^\cup \la \mathcal A_j\ra$, then $Y_k=\emptyset$.
\end{itemize}
\end{proposition}

\begin{proof} If $\x\in \mathcal X_j\cup \{\mathbf v_j\}$ with $\vec u_\x=(u_1,\ldots,u_t)$, then $\partial(\{\x\})$ encases $\vec u^\triangleleft_\x=(u_1,\ldots,u_{t-1})$, ensuring that $\pi_\x(\vec u_\x)=\pi_\x(u_t)/\|\pi_\x(u_t)\|$, both in view of (V1). Thus (V1) further ensures that $\pi_\x(\vec u_\x)$ is a unit vector contained in the subspace $\R^\cup\la \mathcal X^{\pi_\x}_1\cup \ldots\cup \mathcal X^{\pi_\x}_j\ra$, in which case Lemma \ref{lem-orReay-mincase-t=1} (and the comments after its statement) imply that $\mathcal D_\x\subseteq \bigcup_{i\in J_\x,\,i\leq j}\big(\mathcal X_i^{\pi_\x}\cup \{\mathbf \pi_\x(\mathbf v_i)\}\big)$, where
$\pi_\x(\mathcal R)=(\mathcal X_i^{\pi_\x}\cup \{\pi_\x(\mathbf v_i)\})_{i\in J_\x}$ is the virtual Reay system given by Proposition \ref{prop-VReay-modulo}. As a result, $\pi_\x^{-1}(\mathcal D_\x)\cup \{\x\},\,\mathcal A_\y\subseteq \mathcal X_1\cup \{\mathbf v_1\}\cup\ldots\cup \mathcal X_j\cup \{\mathbf v_j\}$ for every $\x\in \mathcal X_j\cup \{\mathbf v_j\}$ and $\y\in \mathcal A_j$. It is now clear that all sets and half-spaces occurring in the proposition are found in $\mathcal X_1\cup \{\mathbf v_1\}\cup\ldots\cup\mathcal X_j\cup \{\mathbf v_j\}$, so we may  w.l.o.g. assume $j=s$, freeing the variable $j$ for other use later.

1. and 2. Let $\x\in \mathcal X_s\cup \{\mathbf v_s\}$ be arbitrary and let $\vec u_\x=(u_1,\ldots,u_t)$.  By (V1), we know that the support set $\partial(\{\x\})\subseteq \mathcal X_1\cup\ldots\cup \mathcal X_{s-1}$ minimally encases $-\vec u^\triangleleft_\x=-(u_1,\ldots,u_{t-1})$.   By proposition \ref{prop-VReay-modulo}, $\pi_\x(\mathcal R)=\big(\mathcal X_i^{\pi_\x}\cup \{\pi_\x(\mathbf v_i)\}\big)_{i\in J_\x}$ is a virtual Reay system with $\vec u_{\pi_\x(\x)}=\pi_\x(\vec u_\x)$. Note $\pi_\x(u_i)=0$ for $i<t$ as
$\partial(\{\x\})$ encases $\vec u^\triangleleft_\x$, but $\pi_\x(u_t)\neq 0$ and $$\pi_\x(\x)=\R_+\pi_\x(u_t)$$ in view of (V1).
By definition, $\mathcal D_\x$ is the unique support set minimally encasing $-\pi_\x(u_t)$, and $\pi_\x^{-1}(\mathcal D_\x)$ is a support set by Proposition \ref{prop-orReay-modulo}.5.
In view of $\pi_\x(\x)=\R_+\pi_\x(u_t)$, we cannot have $\pi_\x(\x)\in \mathcal D_\x$ (by definition of $\mathcal D_\x$), whence $\x\notin \pi_\x^{-1}(\mathcal D_\x)$.
Let \be\label{AD-defined} \mathcal A_\x:=\big(\pi_\x^{-1}(\mathcal D_\x)\cup \partial(\{\x\})\big)^*.\ee
Thus
Proposition \ref{prop-orReay-modulo}.5 and Proposition \ref{prop-orReay-modulo}.4  imply that $\mathcal A_\x$ is a support set with
\be\label{asparagus} \mathcal A_\x^{\pi_\x}=\mathcal D_\x \quad \und\quad \pi_\x^{-1}(\mathcal D_\x)\subseteq \mathcal A_\x.\ee
By Lemma \ref{lemma-MinEncase-Induct} (applied to $\mathcal B=\partial(\{\x\})$ and $\vec u=\vec u_\x$), $\mathcal A_\x$ minimally encases $-\vec u_\x$ urbanely.

Let $\pi_{s-1}:\R^d\rightarrow \R^\cup \la \mathcal X_1\cup \ldots\cup \mathcal X_{s-1}\ra^\bot$ be the orthogonal projection. Note that $\ker \pi_\x\leq \ker \pi_{s-1}$.
  Since $-\pi_\x(u_t)\in \C^\cup(\mathcal D_\x)$ (by definition of $\mathcal D_\x$) and $\ker \pi_\x\subseteq \ker \pi_{s-1}$, it follows that \be\label{guenn}-\pi_{s-1}(u_t)\in \C^\cup(\pi_{s-1}( \mathcal D_\x))=\C \big(\pi_{s-1}(D_\x)\cap \big(\pi_{s-1}( X_s)\cup \{\pi_{s-1}( v_s)\}\big)\big).\ee
By (OR2), $\pi_{s-1}( X_s)\cup \{\pi_{s-1}(v_s)\}$ is a minimal positive basis of size $|\mathcal X_s|+1$ which contains (a positive multiple of) $\pi_{s-1}(u_t)$ with $\R_+\pi_{s-1}(u_t)=\pi_{s-1}(\x)\neq \{0\}$  (recall that $\R_+\pi_\x(u_t)=\pi_\x(\x)$). Thus, since $\pi_\x(\mathbf x)\notin \mathcal D_\x$, we conclude from \eqref{asparagus} and \eqref{guenn} that
\be\label{XisIn}(\mathcal X_s\cup \{\mathbf v_s\})\setminus \{\x\}\subseteq \pi_\x^{-1}(\mathcal D_\x)\subseteq \mathcal A_\x.\ee

 \bigskip

 Now fix $\x\in \mathcal X_s\cup\{\mathbf v_s\}$ and let $$\mathcal A_s=(\mathcal A_\x\setminus \partial(\{\x\}))\cup \{\x\}=\pi_\x^{-1}(\mathcal D_\x)\cup \{\x\},$$ where the second equality follows from \eqref{AD-defined}, \eqref{asparagus} and the fact that all half-spaces from $\pi_\x^{-1}(\mathcal D_\x)$ are nonzero modulo $\ker \pi_\x=\R^\cup \la \partial(\{\x\})\ra$. In view of \eqref{XisIn}, we see that \be\label{trunkeli}\mathcal X_s\cup \{\mathbf v_s\}\subseteq \mathcal A_s.\ee Since $\mathcal A_\x^*=\mathcal A_\x$ (as $\mathcal A_\x$ is a support set) and there is no half-space $\z\in \mathcal A_s$ with $\x\prec \z$ (as this would require $\z\in \mathcal X_j\cup\{\mathbf v_j\}$ for some $j\geq s+1$) nor any half-space $\z\in \mathcal A_\x\setminus \partial(\{\x\})=\pi^{-1}_\x(\mathcal D_\x)$ (the equality follows from \eqref{AD-defined} and \eqref{asparagus}) with $\z\prec \x$ (as this would imply $\pi_\x(\z)=0$, contrary to the definition of $\pi_\x^{-1}(\mathcal D_\x)$), we conclude that $\mathcal A_s^*=\mathcal A_s$. Thus 2(a) holds.
By definition and \eqref{XisIn}, $\mathcal A_s\setminus \{\x\}\subseteq \mathcal A_\x$. Also, $$(\mathcal A_s\setminus \{\x\}\cup\partial(\{\x\}))^*=(\mathcal A_\x\cup \mathcal \partial(\{\x\}))^*=\mathcal A^*_\x=\mathcal A_\x,$$ where the second equality follows in view of $\partial(\{\x\})\subseteq \darrow (\pi_\x^{-1}(\mathcal D_\x)\cup \partial(\{\x\}))=\darrow \mathcal A_\x$.
 Let $\y\in \mathcal A_s$ be arbitrary, let $$\mathcal A_\y=(\mathcal A_s\setminus \{\y\}\cup \partial(\{\y\}))^*,$$  and let $\vec u_\y=(u'_1,\ldots, u'_{t'})$. Note the case $\y=\x$ agrees with the previous definition for $\mathcal A_\x$ as just shown. Since $\y\in \mathcal A_s=\mathcal A_s^*$, the  definition of $\mathcal A_\y$ implies
  \be\label{peppep} \darrow \mathcal A_s=\darrow \mathcal A_\y\cup \{\y\},\ee with the union disjoint. Let $\pi(\mathcal R)=(\mathcal X_i^{\pi}\cup\{\pi(\mathbf v_i)\})_{i\in J}$, where $\pi:\R^d\rightarrow \R^\cup \la \partial(\mathcal A_s)\ra^\bot$ is the orthogonal projection.

\subsection*{Claim A} $\mathcal A_s\subseteq \mathcal X_1\cup \ldots\cup \mathcal X_s\cup \{\mathbf v_j:\; j\in J\}$ with $\pi(\y)\neq \{0\}$ for all $\y\in \mathcal A_s$.

\begin{proof}Consider an arbitrary index $j\in [1,s]\setminus J$, so   $\mathcal X_j^{\pi}=\emptyset$. Then, since $\partial(\mathcal A_s)\subseteq \mathcal X_1\cup\ldots\cup \mathcal X_{s-1}$, it follows from (OR2) that $j<s$, and it follows from Proposition \ref{prop-orReay-BasicProps}.9 that $\mathcal X_j\subseteq  \darrow \partial(\mathcal A_s)\subseteq \darrow \mathcal A_\x$, with the latter inclusion in view of \eqref{peppep} (using $\y=\x$) and $\mathcal A_s^*=\mathcal A_s$. Thus, since $\mathcal A_\x$ is a support set, we must have $\mathbf v_j\notin \darrow \mathcal A_\x$.
 Since $j<s$ and $\x\in \mathcal X_s\cup \{\mathbf v_s\}$, we also have $\mathbf v_j\neq \x$, in which case  $\mathbf v_j\notin \darrow \mathcal A_\x\cup \{\x\}=\darrow \mathcal A_s$, with  the equality in view of \eqref{peppep} (used in the case $\y=\x$), which shows that $\mathcal A_s\subseteq \mathcal X_1\cup \ldots\cup \mathcal X_s\cup \{\mathbf v_j:\; j\in J\}$.

Next suppose that $\pi(\y)=\{0\}$ for some $\y\in \mathcal A_s$. In consequence, if $\y\in \mathcal X_1\cup\ldots\cup \mathcal X_s$, then Proposition \ref{prop-orReay-BasicProps}.9 implies that $\y\in \darrow \partial(\mathcal A_s)$, implying that $\y\prec \z$ with $\z,\,\y\in \mathcal A_s$, which contradicts that $\mathcal A_s^*=\mathcal A_s$. Therefore we must instead have $\y=\mathbf v_j\in \mathcal A_s$ for some $j\in [1,s]$, which as just shown implies $j\in J$. However Proposition \ref{prop-orReay-modulo}.1 implies that $\pi(\mathbf v_j)\neq \{0\}$ for all $j\in J$, contrary to assumption, which completes the claim.
\end{proof}

From Claim A and Proposition \ref{prop-orReay-modulo}.1, we conclude that  $\pi_\y$ is injective on $\mathcal A_s$ with $\pi_\y(\mathcal A_s)=\mathcal A_s^{\pi_\y}$ for all $\y\in \mathcal A_s$ (since $\ker \pi_\y\subseteq \ker \pi$).
Since $\mathcal A_s^*=\mathcal A_s$ and
$\y\in \mathcal A_s$, ensuring that $\darrow\partial(\{\y\})$ is disjoint from $\mathcal A_s$, we have  \be\label{singingbird}\mathcal A_s\setminus \{\y\}\subseteq \mathcal A_\y=(\mathcal A_s\setminus \{\y\}\cup \partial(\{\y\}))^*\subseteq \mathcal A_s\setminus \{\y\}\cup \partial(\{\y\}).\ee Thus \be\label{catchertae}\mathcal A_\y^{\pi_\y}=(\mathcal A_s\setminus \{\y\})^{\pi_\y}=\pi_\y(\mathcal A_s\setminus\{\y\})\quad\und\quad \pi_\y^{-1}(\mathcal A_\y^{\pi_\y})=\mathcal A_s\setminus \{\y\},\ee with the second and third  equalities in view of Claim A and the injectivity of $\pi_\y$ on $\mathcal A_s$, establishing 2(d).
Consequently,  \eqref{catchertae},  Proposition \ref{prop-orReay-modulo}.2 and $\ker \pi_\y=\R^\cup\la \partial(\{\y\})\ra$ imply \be\label{laotat}\partial(\mathcal A_\y^{\pi_\y})=\partial((\mathcal A_s\setminus \{\y\})^{\pi_\y})=\partial(\mathcal A_s\setminus\{\y\})^{\pi_\y}=\partial(\mathcal A_s)^{\pi_\y}.\ee

In view of Claim A and Proposition \ref{prop-orReay-modulo}.1, we see that $\pi$ is injective on $\mathcal A_s$ with $\mathcal A_s^\pi=\pi(
\mathcal A_s)$. Thus, since $\ker \pi_\x\subseteq \ker \pi$, \eqref{catchertae} (applied with $\y=\x$) yields $\pi(A_\x^{\pi_\x})=\pi\pi_\x (A_s\setminus\{x\})=\pi( A_s\setminus \{x\})$.
Since $\mathcal A_\x^{\pi_\x}=\mathcal D_\x$ is a support set that minimally encases $-\pi_\x(u_t)$, Proposition \ref{prop-orReay-BasicProps} (Items 3 and 5) applied to $\mathcal A_\x^{\pi_\x}$, combined with \eqref{laotat} (used with  $\y=\x$) and Lemma \ref{lem-orReay-mincase-t=1}, ensures that
  $$\pi(A_\x^{\pi_\x})=\pi( A_s\setminus \{x\})$$ is a linearly independent set  that minimally encases $-\pi\pi_\x(u_t)=-\pi(u_t)$ and $\R_+\pi(u_t)=\pi(\x)\neq \{0\}$. Consequently,
   $ A_s^\pi=\pi (A_s)=\pi( A_s\setminus \{x\})\cup \{\pi(x)\}=\pi( A_\x^{\pi_x})\cup \{\pi(x)\}$ is a minimal positive basis of size $|\mathcal A_s|$, showing 2(e) holds.


Since $\y\notin \darrow \mathcal A_\y$, two applications of \eqref{peppep} (once taking $\y=\y$ and once taking $\y=\x$) gives $\darrow A_\y=\darrow A_s\setminus \{y\}=(\darrow A_\x\cup \{x\})\setminus \{y\}$. Since $\mathcal A_\x$ is a support set, and thus also a virtual independent set, $\darrow A_\x$ is  linearly independent. Thus, if $\darrow A_\y$ were not linearly independent, we could write the representative $x$ for $\x$ as a linear combination of the elements from $\darrow A_s\setminus \{x,y\}$, which when considered modulo $\R^\cup \la \darrow\partial(A_s)\ra$ would contradict that $\pi(A_s)$ is a minimal positive basis of size $|\mathcal A_s|$ with $x,\,y\in A_s$. Therefore $\darrow A_\y$ is linearly independent, which combined with $\mathcal A^*_\y=\mathcal A_\y$ implies that $\mathcal A_\y$ is a virtual independent set. Thus, in view of Claim A, Proposition \ref{prop-orReay-modulo}.6 and $\partial(\{\y\})\subseteq \darrow (\mathcal A_s\setminus \{\y\}\cup \partial(\{\y\}))^*=\darrow \mathcal A_\y$, it follows that $\mathcal A_\y^{\pi_\y}$ is also a virtual independent set.

Let us show that $\mathcal A_\y^{\pi_\y}=\pi_\y(\mathcal A_s\setminus \{\y\})$  (the equality follow from \eqref{catchertae}) minimally encases $-\pi_\y(\vec u_\y)=-\pi_\y(u'_{t'})/\|\pi_\y(u'_{t'})|$. Note $\R_+\pi_\y(u'_{t'})=\pi_\y(\y)$ by (V1). As just seen,  $A_s^\pi=\pi( A_s)$ is a minimal positive basis of size $|\mathcal A_s|$ with  $\pi$ is injective on $\mathcal A_s$. Thus $\pi$ is also injective on $\pi_\y(\mathcal A_s\setminus\{\y\})=\mathcal A_\y^{\pi_\y}$ (as $\pi\pi_\y=\pi$) and  \be\label{cattrick}-\pi(u'_{t'})\in \C^\circ(\pi(\mathcal A_s\setminus\{\y\}))=\C^\circ(\pi(\mathcal A_\y^{\pi_\y})).\ee
In view \eqref{cattrick} and \eqref{laotat}, we find that  $-\pi_\y(u'_{t'})+a\in \C^\circ(\mathcal A_\y^{\pi_\y})$ for some $a\in \R^\cup\la \partial(\mathcal A_s)^{\pi_\y}\ra=\R^\cup \la\partial(\mathcal A_\y^{\pi_\y})\ra$, which implies $-\pi_\y(u'_{t'})\in \C^\circ(\mathcal A_\y^{\pi_\y})$ by  Proposition \ref{prop-orReay-BasicProps}.4. Thus the virtual independent set $\mathcal A_\y^{\pi_\y}$ minimally encases $-\pi_\y(u'_{t'})=-\pi_\y(\vec u_\y)$ by Lemma \ref{lem-orReay-mincase-t=1}, while $\partial(\{\y\})$ minimally encases $-\vec u^\triangleleft_\y$ by (V1). As a result,  Lemma \ref{lemma-MinEncase-Induct}, Claim A and \eqref{catchertae} now imply that  $(\pi_\y^{-1}(\mathcal A_\y^{\pi_\y})\cup \partial(\{\y\}))^*=(\mathcal A_s\setminus \{\y\}\cup \partial(\{\y\}))^*=\mathcal A_\y$ minimally encases $-\vec u_\y$. Consequently, since we have already shown that $\mathcal A_\y$ and $\mathcal A_\y^{\pi_\y}$ are  virtual independent sets,  2(b) now follows in view of \eqref{singingbird}.

Let $\z\in \mathcal X_s\cup \{\mathbf v_s\}$ be arbitrary. By  definition of $\mathcal A_\z$ and $\mathcal A_s$, we have \be\label{twigrig}\mathcal A_{\z}=\big(\mathcal A_s\setminus \{\z\}\cup\partial(\{\z\})\big)^*=\Big(\big(\pi_\x^{-1}(\mathcal D_\x)\cup \{\x\}\big)\setminus \{\z\}\cup \partial(\{\z\})\Big)^*.\ee Since $\mathcal D_\x$ is a support set, Proposition \ref{prop-orReay-modulo}.5 implies that  $(\pi_\x^{-1}(\mathcal D_\x)\cup \partial(\{\x\}))^*$ is a  support set.  We have already seen that $\mathcal A_\x$ is a support set.
If $\z\neq \x$, then $\z\in \pi_\x^{-1}(\mathcal D_\x)$ by \eqref{XisIn}, whence \eqref{twigrig} implies $\darrow \mathcal A_{\z}\subseteq \darrow (\pi_\x^{-1}(\mathcal D_\x)\cup \partial(\{\x\}))^*\cup \{\x\}$. Consequently, since $(\pi_\x^{-1}(\mathcal D_\x)\cup \partial(\{\x\}))^*$ is a  support set and $\x\in \mathcal X_s\cup\{\mathbf v_s\}$, the only way $\mathcal A_{\z}=\mathcal A_{\z}^*$ can fail to be a support set is if $\mathcal X_s\cup \{\mathbf v_s\}\subseteq \darrow \mathcal A_{\z}$. However this contradicts that $\z\in \mathcal X_s\cup \{\mathbf v_s\}$ but $\z\notin \darrow \mathcal A_{\z}=\darrow (\mathcal A_s\setminus \{\z\}\cup \partial(\{\z\}))$ (as \eqref{peppep} using  is a disjoint union, using $\y=\z$). In summary, we now conclude that $\mathcal A_{\z}$ is a support set for every $\z\in  \mathcal X_s\cup \{\mathbf v_s\}$, in which case Proposition \ref{prop-orReay-modulo}.6 implies that $\mathcal A_{\z}^{\pi_{\z}}$ is also a support set. By the established Item 2(b), the support set $\mathcal A_{\z}^{\pi_{\z}}$ minimally encases $-\pi_{\z}(\vec u_{\z})$, and must then equal $\mathcal D_{\z}=\supp_{\pi_{\z}(\mathcal R)}(-\pi_{\z}(\vec u_{\z}))$, so \be\label{pharz}\mathcal A_{\z}^{\pi_{\z}}=\mathcal D_{\z}\quad\mbox{ for
$\z\in \mathcal X_s\cup \{\mathbf v_s\}$}.\ee Since
 \eqref{catchertae} and \eqref{pharz} imply that  $\mathcal A_\z=\big(\mathcal A_s\setminus \{\z\}\cup \partial(\{\z\})\big)^*=\big(\pi_\z^{-1}(\mathcal A_\z^{\pi_\z})\cup \partial(\{\z\})\big)^*=\big(\pi_\z^{-1}(\mathcal D_\z)\cup \partial(\{\z\})\big)^*$ for $\z\in \mathcal X_s\cup \{\mathbf v_s\}$, all parts of 2(c) are now established, completing  Item 2. Additionally, Item 2(a), \eqref{catchertae} (applied with $\y=\z$) and \eqref{pharz} imply $\mathcal A_s=\pi_\z^{-1}(\mathcal A_\z^{\pi_\z})\cup \{\z\}=\pi_\z^{-1}(\mathcal D_\z)\cup \{\z\}$ for $\z\in \mathcal X_s\cup \{\mathbf v_s\}$, establishing Item 1.

\medskip

 3. Let $\x\in \mathcal A_s$ be arbitrary and let $k=(i_\z)_{\z\in \darrow \mathcal A_\x}$ be an arbitrary tuple of indices with all $i_\z$ sufficiently large (as will be determined during the course of the proof). We maintain the general notation used in Items 1 and 2 and let $$\mathcal B_\x=\partial(\{\x\}).$$
  By Proposition \ref{prop-orReay-BasicProps}.1, we have $$\R^\cup\la \mathcal A_s\ra=\R\la \darrow\wtilde A_s(k)\ra\quad\und\quad \R^\cup \la \mathcal B_\x\ra=\R\la \darrow \wtilde B_\x(k)\ra$$ for any tuple $k$. Let $$\varpi:\R^d\rightarrow \R^\cup \la \mathcal A_s\ra^\bot$$ be the orthogonal projection. Since $\C(G_0)=\R^d$, we also have $\C(\varpi(G_0))=\R^\cup \la \mathcal A_s\ra^\bot$. Thus, by Proposition \ref{prop-reay-basis-exists}, we can find a finite subset $Y=Y_1\cup \ldots\cup Y_\ell\subseteq G_0$ such that $(\varpi(Y_1),\ldots,\varpi(Y_\ell))$ is a Reay system for $\R^\cup \la \mathcal A_s\ra^\bot$ with $\varpi$ injective on $Y$.
We  will show Item 3 holds for an arbitrary such $Y$.

 In view of Proposition \ref{prop-VReay-RepBasics}.1, once all $i_\z$ are sufficiently large, we can assume $\darrow \tilde B_\x(k)$ minimally encases $-\vec u_\x^\triangleleft$. In view of Item 2(b)(e), Proposition \ref{prop-orReay-modulo}.1 and Proposition \ref{prop-orReay-BasicProps}.9, it follows that the half-spaces from $\darrow \mathcal A_\x$ mapped to $\{0\}$ by $\pi_\x$ are precisely those in $\darrow \mathcal B_\x$. Consequently, in view of Lemma \ref{lemma-minencase-mod}.2 (applied with $X=\darrow \tilde B_\x(k)$ and $Y=(A_s\setminus\{x\})(k)\cup
 \big(\darrow\wtilde{ \partial(A_s)}\setminus \darrow \tilde B_\x\big)(k)\cup Y_k$), in order to complete the proof, it suffices to show $$\pi_\x(A_s\setminus\{x\})(k)\cup
 \pi_\x\big(\darrow\wtilde{ \partial(A_s)}\setminus \darrow \tilde B_\x\big)(k)\cup \pi_\x(Y_k)$$ minimally encases $-\pi_\x(u_t)$ with $\pi_\x$ injective on $(A_s\setminus\{x\})(k)\cup
 \big(\darrow\wtilde{ \partial(A_s)}\setminus \darrow \tilde B_\x\big)(k)\cup Y_k$,  for some $Y_k\subseteq Y$ (as well as the additional moreover statement). We begin with the injectivity of $\pi_\x$.

\subsection*{Claim B} $\pi_\x$ is injective on $(A_s\setminus\{x\})(k)\cup
 \big(\darrow\wtilde{ \partial(A_s)}\setminus \darrow \tilde B_\x\big)(k)\cup Y$ for any tuple $k=(i_\z)_{\z\in \darrow \mathcal A_s\setminus\{\x\}}$ with all $i_\z$ sufficiently large.

\begin{proof}
 Item  2(b) implies that  \be\label{dessertgrease}\mathcal A_\x^{\pi_\x}=\pi_\x(\mathcal A_s\setminus \{\x\})\ee minimally encases $-\pi_\x(u_t)$ with $\mathcal A_\x^{\pi_\x}$ a virtual independent set. Hence, by Lemma \ref{lemma-minencase-Rep} (applied to $\mathcal A_\x^{\pi_\x}$) and Proposition \ref{prop-orReay-modulo}.2, $\darrow  \tilde A_\x^{\pi_\x}(k)=(\darrow \tilde A_\x)^{\pi_\x}(k)$ minimally encases $-\pi_\x(u_t)$ (once all $i_\z$ are sufficiently large).
 In view of Item 2(e),  Proposition \ref{prop-orReay-BasicProps}.9 and Proposition \ref{prop-orReay-modulo}.1, it follows  that $(\darrow \tilde A_\x)^{\pi_\x}(k)=\pi_\x\big(\darrow \tilde A_\x\setminus \darrow \tilde B_\x\big)(k)$ with  $\pi_\x$ injective on $\darrow \mathcal A_\x\setminus \darrow \mathcal B_\x$ and  $(\darrow \tilde  A_\x\setminus \darrow  \tilde B_\x)(k)$. As a result,
$\pi_\x(\darrow \tilde A_\x\setminus \darrow \tilde B_\x)(k)$ minimally encases $-\pi_\x(u_t)$ and it follows that $$\pi_\x(\darrow \tilde A_\x\setminus \darrow \tilde B_\x)(k)\cup \{\pi_\x(u_t)\}\quad\mbox{ is a minimal positive basis of size $|\darrow \mathcal A_\x\setminus \darrow \mathcal B_\x|+1$}$$ for the fixed (independent of $k$) subspace $\R^\cup \la\mathcal A_\x^{\pi_\x}\cup\{\pi_\x(u_t)\}\ra$ (by Proposition \ref{prop-orReay-BasicProps}.1). Recall that
$\mathcal A_\x^{\pi_\x}=\pi_\x(\mathcal A_s\setminus \{\x\})$ and $\R_+\pi_\x(\x)=\R_+\pi_\x(u_t)$. As a result,  $$\R^\cup \la\mathcal A_\x^{\pi_\x}\cup\{\pi_\x(u_t)\}\ra=\R^\cup \la \pi_\x(\mathcal A_s)\ra.$$
Thus, since $\ker \pi_\x\leq \ker \varpi$,  Proposition  \ref{prop-reay-basis-exists} implies that,  for all $k$ with every $i_\z$ sufficiently large,
  $$\mathcal R_k=(\pi_\x(\darrow \tilde A_\x\setminus \darrow \tilde B_\x)(k)\cup \{\pi_\x(u_t)\},\pi_\x(Y_1),\ldots,\pi_\x(Y_\ell))$$ is an ordinary Reay system for $\R^\cup \la\mathcal B_\x\ra^\bot$.
  In view of $\ker \pi_\x\leq \ker \varpi$, $\pi_\x$ must be injective on $Y$ as $\varpi$ is injective on $Y$. Thus, since $\mathcal R_k$ is a Reay system,   $\pi_\x$ is injective on $$( \darrow \tilde A_\x\setminus  \darrow \tilde B_\x)(k)\cup Y=(\tilde A_s\setminus\{\tilde x\})(k)\cup \big(\darrow\wtilde{ \partial(A_s)}\setminus \darrow \tilde B_\x\big)(k)\cup Y$$ (as it is injective on each component set in $\mathcal R_k$ as already noted). We still must show that $\pi_\x$ is injective on $(A_s\setminus\{x\})(k)\cup \big(\darrow\wtilde{ \partial(A_s)}\setminus \darrow \tilde B_\x\big)(k)\cup Y$ for $k$ with all $i_\z$ sufficiently large.

 Observe by Item 2(b) that $\R^\cup\la \partial(\mathcal A_s)\ra=\R^\cup \la \mathcal B_\x\cup \partial(\mathcal A_s\setminus \{\x\})\ra=\R^\cup\la \mathcal B_\x\cup \partial(\mathcal A_\x)\ra$.
Thus $\pi:\R^d\rightarrow \R^\cup\la \partial(\mathcal A_s)\ra^\bot=\R^\cup \la\mathcal B_\x\cup \partial(\mathcal A_\x)\ra^\bot$ is the orthogonal projection with
 $\pi \pi_\x=\pi$ (as $\ker \pi_\x\leq \pi$).
 Moreover, $\partial(\mathcal A_\x)^{\pi_\x}=\partial(\mathcal A_\x^{\pi_\x})\subseteq \darrow \mathcal A_\x^{\pi_\x}=(\darrow\mathcal A_\x)^{\pi_\x}=\pi_\x(\darrow \mathcal A_\x\setminus \darrow \mathcal B_\x)$ in view of Item 2(b)(e), Proposition \ref{prop-orReay-modulo}.2 and Proposition \ref{prop-orReay-BasicProps}.9. As a result,
 since an ordinary Reay system may be considered as an oriented Reay system (for which all half-spaces have trivial boundary) by replacing each element
  with the ray it defines, Propositions \ref{prop-orReay-modulo}.1  and \ref{prop-orReay-BasicProps}.1 ensure that $\pi(\mathcal R_k)$ is a Reay system for $\R^\cup\la \partial(\mathcal A_s)\ra^\bot$.
 Now  $\mathcal A_\x^{\pi_\x}$ is a virtual independent set which minimally encases $-\pi_\x(u_t)$ (by Item 2(b)). Thus Lemma \ref{lem-orReay-mincase-t=1}.1 ensures that $-\pi_\x(u_t)$ can be written as a sum of representatives from all the half-spaces in  $\mathcal A_\x^{\pi_\x}$. Since $\pi(\mathcal A_\x^{\pi_\x})=\pi(\mathcal A_s\setminus \{\x\})$ by \eqref{dessertgrease}, applying $\pi$ to this sum shows that $-\pi(u_t)=-\pi\pi_\x(u_t)$ is a sum of representatives of all the half-spaces from $\pi(\mathcal A_\x^{\pi_\x})=\pi(\mathcal A_s\setminus \{\x\})$, with these representatives being linearly independent in view of Item 2(e). In particular, since each half-space from $\pi(\mathcal A_s)$ has trivial boundary, it follows that  $\pi(\mathcal A_\x^{\pi_\x})=\pi(\mathcal A_s\setminus \{\x\})$ is a virtual independent set in $\pi(\mathcal R_k)$ (identifying this set of $1$-dimensional rays with the corresponding set of representatives in $\pi(\mathcal R_k)$), and now
   Lemma \ref{lem-orReay-mincase-t=1}.2 applied to the virtual independent set $\pi(\mathcal A_s\setminus \{\x\})$  shows that
 $\pi(\mathcal A_\x^{\pi_\x})=\pi(\mathcal A_s\setminus \{\x\})$ minimally encases $-\pi(u_t)$. Note  $\pi$ is injective on $\mathcal A_s\setminus \{\x\}$ in view of Item 2(e). However, by definition of $\pi$, we have $\R_+\pi(\tilde \y(i_\y))=\pi(\y)$ for all $\y\in \mathcal A_s\setminus\{\x\}$. Thus $\pi(\tilde A_s\setminus \{\tilde x\})(k)\cup \{\pi(u_t)\}$ is a minimal positive basis of size $|\mathcal A_s|$, and $\pi(\mathcal R_k)=(\pi(\tilde A_s\setminus \{\tilde x\})(k)\cup \{\pi(u_t)\},\pi(Y_1),\ldots,\pi(Y_\ell))$ is a Reay system with $\pi$ injective on $Y$  in view of $\ker \pi\leq \ker \varpi$. For each $\y\in \mathcal A_s\setminus \{\x\}$, we have $\pi(\y(i_\y))/\|\pi(\y(i_\y))\|\rightarrow v_\y$, where $v_\y=\pi(\tilde \y(i_\y))/\|\pi(\tilde \y(i_\y))\|$ is the constant unit vector pointing in the direction given by the $1$-dimensional ray $\pi(\y)$
(in view of   $\R_+\pi(\tilde \y(i_\y))=\pi(\y)$). Let $\mathfrak V^\pi_k$ be the associated complete simplicial fan with $V(\mathfrak V^\pi_k)=\pi(\tilde A_s\setminus \{\tilde x\})(k)\cup \{\pi(u_t)\}\cup\pi(Y)$.
Note $\mathfrak V^\pi_k$ does not depend on $k$ apart from the choice of vertices used to represent the one-dimensional rays from $\mathfrak V^\pi_k$ (since $\R_+\pi(\tilde \y(i_\y))=\pi(\y)$).
Then, by Proposition \ref{prop-FanStability}.4, for any tuple $k$ with all $i_\z$ sufficiently large, there is a simplicial isomorphism between $\mathfrak V_k^\pi$ and the simplicial fan $(\mathfrak V'_k)^\pi$ with vertices $\pi(A_s\setminus \{x\})(k)\cup \{\pi(u_t)\}\cup\pi(Y)$. In particular, $|\pi(A_s\setminus \{x\})(k)\cup \{\pi(u_t)\}\cup\pi(Y)|=|\pi(\tilde A_s\setminus \{\tilde x\})(k)\cup \{\pi(u_t)\}\cup\pi(Y)|=|\mathcal A_s|+|Y|$, which shows that $\pi$, and thus also $\pi_\x$ (in view of $\ker\pi_\x\leq \ker \pi$), is injective on $$(A_s\setminus \{x\})(k)\cup Y,$$ mapping all such elements to distinct non-zero elements. Thus, since   all half-spaces from $\darrow \partial(\mathcal A_s)$ are mapped to zero by $\pi$, it follows that  $$\pi_\x(A_s\setminus\{x\})(k)\cap \pi_\x(\wtilde{\darrow \partial(A_s)}\setminus \darrow \tilde B_\x)(k)=\emptyset.$$
Since $\pi_\x$ is injective on  $(\tilde A_s\setminus\{\tilde x\})(k)\cup\big(\wtilde{\darrow \partial(A_s)}\setminus \darrow \tilde B_\x\big)(k)\cup Y$, it is, in particular,  injective on $\big(\wtilde{\darrow \partial(A_s)}\setminus \darrow \tilde B_\x\big)(k)\cup Y$. Combining the last three conclusions, it follows that $\pi_\x$ is  injective on
$(A_s\setminus\{x\})(k)\cup \big(\wtilde{\darrow \partial(A_s)}\setminus \darrow \tilde B_\x\big)(k)\cup Y$, completing the claim.
\end{proof}

In view of claim B, it remains to  show  $\pi_\x(A_s\setminus\{x\})(k)\cup
 \pi_\x\big(\darrow\wtilde{ \partial(A_s)}\setminus \darrow \tilde B_\x\big)(k)\cup \pi_\x(Y_k)$ minimally encases $-\pi_\x(u_t)$ for some $Y_k\subseteq Y$ (as noted earlier), as well as the additional moreover statement. By Proposition \ref{prop-VReay-modulo}, $\pi_\x(\mathcal R)$ is a virtual Reay system.

\subsection*{Claim C} $\pi_\x(\mathcal A_s)$ is the set given by Proposition \ref{prop-VReay-SupportSet}.1 for $\mathcal X_s^{\pi_\x}\cup \{\pi_\x(\mathbf v_s)\}$. In particular,  $\mathcal A_\x^{\pi_\x}=\mathcal A_{\pi_\x(\x)}$.

\begin{proof} Let $\y\in \mathcal X_s$ and $\mathcal B_\y=\partial(\{\y\})$.
Let
 $\pi_{\x\y}:\R^d\rightarrow \R^\cup \la \mathcal B_\x\cup \mathcal B_\y\ra^\bot$ be the orthogonal projection. Now $\mathcal A_\y^{\pi_\y}=\mathcal D_\y=\supp_{\pi_\y(\mathcal R)}(-\pi_\y(\vec u_\y))$ is the unique support set minimally encasing $-\pi_\y(\vec u_\y)$ by Item 2(c), and we  have \be\label{petal-bloom}\darrow \mathcal B_\x^{\pi_\y}=(\darrow \mathcal B_\x)^{\pi_\y}\subseteq \darrow \partial(\mathcal A_s)^{\pi_{\y}}=(\darrow \partial(\mathcal A_s))^{\pi_{\y}}\subseteq (\darrow \mathcal A_\y)^{\pi_\y}=\darrow \mathcal A_\y^{\pi_\y}\ee in view of Item 2(b)(e) and Proposition \ref{prop-orReay-modulo}.2, in which case Proposition \ref{prop-orReay-modulo}.6 ensures that $\mathcal A_\y^{\pi_{\x\y}}=(\mathcal A_\y^{\pi_\y})^{\pi_{\x\y}}$ is a support set.
 Since $-\pi_\y(\vec u_\y)$ consists of a single element, its minimal encasement by the support set $\mathcal A_\y^{\pi_\y}$ is urbane. Thus Proposition \ref{prop-orReay-minecase-char}.4 and \eqref{petal-bloom} ensure that the support set
 $\mathcal A_\y^{\pi_{\x\y}}=(\mathcal A_\y^{\pi_\y})^{\pi_{\x\y}}$
  minimally encases $-\pi_{\x\y}(\vec u_\y)$, meaning $\supp_{\pi_{\x\y}(\mathcal R)}(-\pi_{\x\y}(\vec u_\y))=\mathcal A_\y^{\pi_{\x\y}}$.
  In view of Item 2(b), we have $\mathcal A_\y^{\pi_\y}=\pi_\y(\mathcal A_s\setminus \{\y\})$, while Item 2(e) ensures that $\pi$, and thus also $\pi_{\x\y}$, is injective on $\mathcal A_s\setminus \{\y\}$ mapping no half-space to $0$, whence
 $\pi_{\x\y}^{-1}(\mathcal A_\y^{\pi_{\x\y}})=\mathcal A_s\setminus \{\y\}$, which implies  that $(\mathcal A_s\setminus \{\y\})^{\pi_\x}$ is the  pull-back of
  $\supp_{\pi_{\x\y}(\mathcal R)}(-\pi_{\x\y}(\vec u_\y))=\mathcal A_\y^{\pi_{\x\y}}$  to $\pi_\x(\mathcal R)$. Thus, in view of Item 1, we find that $(\mathcal A_s\setminus \{\y\})^{\pi_\x}\cup \{\pi_\x(\y)\}$ is the set given by Proposition \ref{prop-VReay-SupportSet}.1 for $\mathcal X_s^{\pi_\x}\cup \{\pi_\x(\mathbf v_s)\}$. Since $\pi_\x(\mathcal A_s)=\pi_\x(\mathcal A_s\setminus \{\y\})\cup  \{\pi_\x(\y)\}=(\mathcal A_s\setminus \{\y\})^{\pi_\x}\cup \{\pi_\x(\y)\}$ by Item 2(e), the main part of the claim is  complete. To establish the
   in particular statement, note the main part together with Item 2(e) and Proposition \ref{prop-orReay-modulo}.2 yields  $\mathcal A_{\pi_\x(\x)}=\big(\pi_\x(\mathcal A_s)\setminus \{\pi_{\x}(\x)\}\cup \partial(\{\pi_\x(\x)\})\big)^*=\big((\mathcal A_s\setminus \{\x\})^{\pi_\x}\cup \partial(\{\x\})^{\pi_\x}\big)^*=\big(\big((\mathcal A_s\setminus \{\x\})\cup \partial(\{\x\})\big)^*\big)^{\pi_\x}=\mathcal A_\x^{\pi_\x}$.
\end{proof}

If we knew Item 3 held for any $\x\in\mathcal A_s$ whenever $\mathcal B_\x=\emptyset$, that is, when $\pi_\x:\R^d\rightarrow \R^d$ is the identity map, $t=1$ and $\mathcal A_\x=\mathcal A_s\setminus\{\x\}$, then applying this case to $\pi_\x(\x)\in\pi_\x(\mathcal A_s)$ in the virtual Reay system $\pi_\x(\mathcal R)$  would yield (in view of Claim C) that $$\pi_\x(A_s\setminus \{x\})(k)\cup \darrow \wtilde{\partial(\pi_\x( A_s))}(k)\cup \pi_\x(Y_k)$$ minimally encases $-\pi_\x(u_t)=-\vec u_{\pi(\x)}$ (by Proposition \ref{prop-VReay-modulo}) for some $Y_k\subseteq Y$, whenever all $i_\z$ are sufficiently large.
Moreover, if $(\mathcal A_s\setminus \{x\})(k)\subseteq \R^\cup\la \mathcal A_s\ra$, then $\pi_\x(\mathcal A_s\setminus \{x\})(k)\subseteq \R^\cup\la \pi_\x(\mathcal A_s)\ra$, and Item 3 would further give $\pi_\x(Y_k)=\emptyset$, forcing $Y_k=\emptyset$.
However,  Item 2(e), Proposition \ref{prop-orReay-modulo}.2 and Proposition \ref{prop-orReay-BasicProps}.9 imply $\darrow \wtilde{\partial(\pi_\x( A_s))}(k)=\pi_\x(\darrow \wtilde{\partial(\mathcal A_s)}\setminus \darrow\tilde B_\x)(k)$, and thus we obtain the needed conclusion that  $\pi_\x(A_s\setminus \{x\})(k)\cup \pi_\x(\darrow \wtilde{\partial(\mathcal A_s)}\setminus \darrow\tilde B_\x)(k)\cup \pi_\x(Y_k)$ minimally encases $-\pi_\x(u_t)$, which would complete the proof as already remarked.
Thus it suffices to handle this case, so we now assume that $\partial(\{\x\})=\mathcal B_\x=\emptyset$, so that $\pi_\x:\R^d\rightarrow \R^d$ is the identity map, $t=1$, \be\label{gulftrack}\mathcal A_\x=\mathcal A_s\setminus\{\x\}\quad\und\quad  \partial(\mathcal A_\x)=\partial(\mathcal A_s).\ee This will simplify notation. In particular, we now have
\be\label{eletrunk}\mathcal R_k=(\darrow \tilde A_\x(k)\cup \{u_1\},Y_1,\ldots,Y_\ell),\ee which is an ordinary Reay system for $\R^d$, for any tuple $k$ with all $i_\z$ sufficiently large (as argued in Claim B).
Note that \be\label{supperdinner}\supp_{\mathcal R_k}(-u_1)=\darrow \tilde A_\x(k).\ee
We also have $-u_1\in \C^\cup(\mathcal A_\x)^\circ$ (cf. Lemma \ref{lem-orReay-mincase-t=1}) as the virtual independent set $\mathcal A_\x=\mathcal A_s\setminus\{\x\}=\mathcal A_\x^{\pi_\x}$ minimally encases $-u_1=-\pi_\x(u_t)$.

Now let $\kappa=(\iota_\z)_{\z\in \darrow \mathcal A_\x}$ be a \emph{fixed} tuple with all $\iota_\z$ sufficiently large that all statements derived above (for the arbitrary tuple $k$ with all $i_\z$ sufficiently large)  are applicable for $\kappa$ and such that  all coordinates $\iota_\z$  are also greater than the global constant for the order uniform limit $\lim_{k\rightarrow \infty}\C((\darrow \tilde \y)(k))=\y$
given by Proposition \ref{prop-VReay-RepBasics}.2 for each $\y \in \mathcal A_\x$.
In view of Proposition \ref{prop-VReay-RepBasics}.2 and the definition of order uniform limits, by increasing how large each $i_\z$ must be (relative to $\kappa$), we find that
\be\label{umbrello}
 \C((\darrow \tilde \y)(\kappa))\subseteq \C((\darrow \tilde \y)(k))\ee for each $\y\in \mathcal A_\x$ and any tuple $k$ with all $i_\z$ sufficiently large.

Let $\y\in \mathcal A_\x$ be arbitrary and let $\overline v_\y\in \mathcal E_\y^\bot:=\R^\cup \la \partial(\{\y\})\ra^\bot$ be the unit vector such that $\overline \y=\mathcal E_\y+\R_+ \overline v_\y$.
By (V1) and \eqref{x(i)-filtered-form}, we know that $$\y(i_\y)=-x_{i_\y}^\y+\alpha_{i_\y}^\y\overline v_\y+\epsilon_{i_\y}^\y\quad\und\quad \tilde \y(i_\y)=-x_{i_\y}^\y+\alpha_{i_\y}^\y\overline v_\y$$
for some $x_{i_\y}^\y\in \R^\cup (\partial(\{\y\}))= \mathcal E_\y$,  \ $\alpha_{i_\y}^\y>0$ and $\epsilon_{i_\y}^\y\in \R \la\y\ra^\bot=(\mathcal E_\y+\R \overline v_\y)^\bot$, with \be\label{by5}\alpha_{i_\y}^\y\in o(\|x_{i_\y}^\y\|)\quad \mbox{if $\partial(\{\y\})\neq \emptyset$ }\quad\und\quad \|\epsilon_{i_\y}^\y\|\in o(\alpha_{i_\y}^\y).\ee   In particular,
$$\tilde \y(\iota_\y)=-x_{\iota_\y}^\y+\alpha_{\iota_\y}^\y\overline v_\y=-\chi_\y+
a_\y\overline v_\y,$$ where $\chi_\y:=x_{\iota_\y}^\y\in \mathcal E_\y$ and $a_\y:=\alpha_{\iota_\y}^\y>0$, is a fixed, nonzero vector.

In view of \eqref{umbrello}, there is a positive linear combination of the elements from $\darrow\tilde \y(k)$ equal to $\tilde \y(\iota_\y)=-\chi_\y+
a_\y\overline v_\y\in \C( (\darrow \tilde \y) (\kappa) )\subseteq \C( (\darrow \tilde \y) (k) )$ (the element inclusion holds as $\tilde \y(\iota_\y)\in \darrow \tilde \y(\kappa)$). Since  $\tilde \y(i_\y)\in \darrow
\tilde\y(k)$ is the only representative for a half-space not contained in $\mathcal E_\y$, the coefficient of $\y(i_\y)$ in this linear combination must be $a_\y/\alpha_{i_\y}^\y>0$. Replacing $\tilde \y(i_\y)$ with $\y(i_\y)$ in this linear combination, we find that
\be\label{y-con}z_{i_\y}^\y:=\tilde \y(\iota_\y)+(a_\y/\alpha_{i_\y}^\y)\epsilon^\y_{i_\y}=-\chi_\y+a_\y\overline v_\y+(a_\y/\alpha_{i_\y}^\y)\epsilon^\y_{i_\y}\in \R_+^\circ\y(i_\y)+\C\Big(\wtilde{\darrow\partial(\{\y\})}(k)\Big).\ee
Since $\|\epsilon_{i_\y}^\y\|\in o(\alpha_{i_\y}^\y)$ (by \eqref{by5}) with $a_\y>0$, we have $\|(a_\y/\alpha_{i_\y}^\y)\epsilon^\y_{i_\y}\|\rightarrow 0$, meaning $\|(a_\y/\alpha_{i_\y}^\y)\epsilon^\y_{i_\y}\|\in o(\|\tilde\y(\iota_\y)\|)=o(1)$ (as $\tilde \y(\iota_\y)\neq 0$).
In consequence $$\{z_{i_\y}^\y\}_{i_\y=1}^\infty=\{\tilde \y(\iota_\y)+(a_\y/\alpha_{i_\y}^\y)\epsilon^\y_{i_\y}\}_{i_\y=1}^\infty$$ is a radially convergent   sequence of terms $z_{i_\y}^\y\in\R^\circ_+\y(i_\y)+\C\Big(\wtilde{\darrow\partial(\{\y\})}(k)\Big)$ with limit $\tilde \y(\iota_\y)/\|\tilde \y(\iota_\y)\|$. Moreover, if $(A_s\setminus \{x\})(k)\subseteq \R^\cup\la \mathcal A_s\ra$, then $\y(i_\y)\in \R^\cup \la \mathcal A_s\ra$, whence $\epsilon_{i_\y}^\y\in \R^\cup \la \mathcal A_s\ra$ and $z_{i_\y}^\y\in \R^\cup \la \mathcal  A_s\ra$.

Let $\mathfrak F$ be the complete simplicial fan associated to the Reay system $\mathcal R_\kappa$ with vertex set $V(\mathfrak F)=\{\tilde \y(\iota_\y):\;\y\in \mathcal A_\x\}\cup
\{\tilde \y(\iota_\y):\;\y\in \darrow \partial(\mathcal A_\x)\} \cup Y\cup \{u_1\}$. Since each $\{z_{i_\y}^\y\}_{i_\y=1}^\infty$ is a radially convergent sequence with limit $\tilde\y(\iota_\y)/\|\tilde\y(\iota_\y)\|$,
Proposition \ref{prop-FanStability}.4 implies that, for any tuple $k$ with all $i_\z$ sufficiently large,  the map $\varphi_k: \mathfrak F\rightarrow \mathfrak F^{(k)}$, which replaces each vertex
$\tilde \y(\iota_\y)$ for $\y\in \mathcal A_\x$ with the slightly perturbed vertex $z^\y_{i_\y}$
(and leaves all other vertices fixed), is a simplicial isomorphism of $\mathfrak F$ with $\mathfrak F^{(k)}$, where $\mathfrak F^{(k)}$ denotes the resulting complete simplicial fan.
Moreover, further assume all coordinates $i_\y$ are sufficiently large that Proposition \ref{prop-FanStability}.6 can be applied to  $\mathfrak F^{(k)}$ with $x=-u_1$ (possible as the sequences  $\{z_{i_\y}^\y\}_{i_\y=1}^\infty$ are radially convergent). Then, by  proposition \ref{prop-FanStability}.6 and \eqref{supperdinner}, we have \be\label{trackstill}\{z_{i_\y}^\y:\;\y\in \mathcal A_\x\}\cup \{\tilde \y(i_\y):\; \y\in\darrow \partial(\mathcal A_\x)\} =\varphi_k\big(\mathcal \darrow  \tilde A_\x(k)\big)\subseteq \supp_{\mathfrak F^{(k)}}(-u_1)\ee for any  tuple $k$ with all $i_\z$ sufficiently large.
Consequently,  $$\supp_{\mathfrak F^{(k)}}(-u_1)=\{z_{i_\y}^\y:\;\y\in \mathcal A_\x\}\cup \{\tilde \y(i_\y):\; \y\in \darrow \partial(\mathcal A_\x)\}\cup Y_k$$ for some subset $Y_k\subseteq Y$ (note $u_1$ can never be in the support set $\supp(-u_1)$ in view of the definition of $\supp$). Since $\mathfrak F$ is a complete simplicial fan associated to the Reay system $\mathcal R_\kappa$, any support set for $\mathfrak F$ cannot contain all elements from $Y_j$, for any $j\in [1,\ell]$. Thus, since the simplicial isomorphism $\varphi_k$ gives a correspondence between  support sets in $\mathfrak F$ and support sets in $\mathfrak F^{(k)}$, we have  \be\label{toocan}Y_j\nsubseteq Y_k\quad\mbox{for all $j\in [1,\ell]$}.\ee
If $(A_s\setminus\{x\})(k)\subseteq \R^\cup \la \mathcal A_s\ra$, then we can perform all the above using the complete simplicial fan $\mathfrak F'$ with  $V(\mathfrak F')=\{\tilde \y(\iota_\y):\;\y\in \mathcal A_\x\}\cup
\{\tilde \y(\iota_\y):\;\y\in \darrow \partial(\mathcal A_\x)\} \cup \{u_1\}$ in place of $\mathfrak F$ (since all $z_{i_\y}^\y\in \R^\cup \la \mathcal A_s\ra$ in this case) and thereby conclude that $Y_k=\emptyset$.
Now by definition of $\supp_{\mathfrak F^{(k)}}(-u_1)$, we have $$-u_1\in \C^\circ\big(\{z_{i_\y}^\y:\;\y\in \mathcal A_\x\}\big)+\C^\circ\big(\darrow \wtilde{\partial(A_\x)}(k)\big)+\C^\circ(Y_k)$$
with $\supp_{\mathfrak F^{(k)}}(-u_1)$ linearly independent. Hence, in view of \eqref{y-con} and the definition of the $z_{i_\y}^\y$, we find that $$-u_1\in \C^\circ\big(A_\x(k)\cup \darrow\wtilde{\partial(  A_\x)}(k)\cup Y_k\big)$$ will follow once we know $A_\x(k)\cup \darrow\wtilde{\partial( A_\x)}(k)\cup Y_k$ is linearly independent, which will then, in turn, also ensure that $A_\x(k)\cup \darrow\wtilde{\partial(  A_\x)}(k)\cup Y_k
=(A_s\setminus \{x\})(k)\cup \darrow \wtilde{\partial(A_s)}(k)\cup Y_k$ (cf. \eqref{gulftrack})
minimally encases $-u_1=\pi_\x(\vec u)$. Thus, to complete the proof, it remains to show  the following claim.

\subsection*{Claim D} $A_\x(k)\cup \darrow\wtilde{\partial(  A_\x)}(k)\cup Y_k$ is linearly independent for $k$ with all $i_\z$ sufficiently large.

\begin{proof}
The proof is a variation on the argument used to establish the injectivity of $\pi_\x$ in Claim B.
Since $\darrow\wtilde{\partial( A_\x)}(k)\cup Y_k$ is a subset of the linearly independent set $\supp_{\mathfrak F^{(k)}}(-u_1)$, it follows that \be\label{tonicA}\darrow\wtilde{\partial(  A_\x)}(k)\cup Y_k\quad\mbox{is linearly independent.}\ee
Since $\mathcal A_\x=\mathcal A_s\setminus\{\x\}$ with $\partial(\{\x\})=\emptyset$, we have $\R^\cup\la \partial(\mathcal A_s)\ra=\R^\cup\la \partial(\mathcal A_\x)\ra$. Thus
  $\pi:\R^d\rightarrow \R^\cup \la \partial(\mathcal A_\x)\ra^\bot$ is the orthogonal projection. Since an ordinary Reay system may be considered as an oriented Reay system (for which all half-spaces have trivial boundary) by replacing each element
  with the ray it defines, Proposition \ref{prop-orReay-modulo}.1 ensures that $\pi(\mathcal R_k)$ is a Reay system.
  Since $\mathcal A_\x=\mathcal A_s\setminus\{\x\}=\mathcal A_\x^{\pi_\x}$ is a virtual independent set which minimally encases $-u_1=-\pi_\x(u_t)$, it follows from Lemma \ref{lem-orReay-mincase-t=1} and Proposition \ref{prop-orReay-BasicProps}.3  that $\pi(\mathcal A_\x)$ minimally encases $-\pi(u_1)$  with $\pi$ injective on $\mathcal A_\x$. However, by definition of $\pi$, we have $\R_+\pi(\tilde \y(i_\y))=\pi(\y)$ for all $\y\in \mathcal A_\x$. Thus $\pi(\tilde A_\x(k))\cup \{\pi(u_1)\}$ is a minimal positive basis of size $|\mathcal A_\x|+1$ for the subspace $\R^\cup\la \pi(\tilde A_\x(k))\cup \pi(u_1)\ra=\R^\cup\la \pi(\mathcal A_\x)\cup \pi(\x)\ra=\R^\cup\la \pi(\mathcal A_s)\ra$, and $$\mathcal \pi(\mathcal R_k)=(\pi(\tilde A_\x(k))\cup \{\pi(u_t)\},\pi(Y_1),\ldots,\pi(Y_\ell))$$ is a Reay system with $\pi$ injective on $\tilde A_\x(k)\cup Y\cup \{u_t\}$  (since   $\ker \pi\leq \ker \varpi$).

  For each $\y\in \mathcal A_\x$, we have $\pi(\y(i_\y))/\|\pi(\y(i_\y))\|\rightarrow v_\y$, where $v_\y=\pi(\tilde \y(i_\y))/\|\pi(\tilde \y(i_\y))\|$ is a constant unit vector pointing in the same direction as $\pi(\tilde \y(i_\y))$ for $\y\in\mathcal A_\x$
(in view of   $\R_+\pi(\tilde \y(i_\y))=\pi(\y)$).
Let $\mathfrak F^\pi_k$ be the associated complete simplicial fan with $V(\mathfrak
F^\pi_k)=\pi(\tilde A_\x(k))\cup \{\pi(u_t)\}\cup\pi(Y)$ associated to the Reay system $\pi(\mathcal R_k)$.
Note $\mathfrak F_k^\pi$ is fixed (and independent of $k$) apart from some possible variation in which positive scalar multiples are used for the vertices.
Then, by Proposition \ref{prop-FanStability}.4, for any tuple $k$ with all $i_\z$ sufficiently large, there is a simplicial isomorphism $\varphi$ between $\mathfrak F_k^\pi$ and a simplicial fan ${\mathfrak F'_k}^\pi$ with vertices $\pi(A_\x(k))\cup \{\pi(u_t)\}\cup\pi(Y)$. In particular,  support sets map to support sets. Since $\mathfrak F_k^\pi$ is a simplicial fan associated to the Reay system $\pi(\mathcal R_k)$, \eqref{toocan} ensures that $\pi(\tilde A_\x(k))\cup \pi(Y_k)$ is a support set for $\mathfrak F^\pi_k$, whence its image  $\varphi\big(\pi(\tilde A_\x(k))\cup \pi(Y_k)\big)=\pi(A_\x(k))\cup\pi(Y_k)$ is also a support set, and thus linearly independent (as support sets of a simplicial fan are linearly independent by definition).
Furthermore, since the simplicial isomorphism $\varphi$ is injective on vertices, and since $\pi$ is injective on $\tilde A_\x(k)\cup Y$, it follows that $|\mathcal A_\x|+|Y_k|=|\tilde A_\x(k)\cup Y_k|=|\varphi\big(\pi(\tilde A_\x(k))\cup \pi(Y_k)\big)|=|\pi(A_\x(k))\cup\pi(Y_k)|$, forcing $\pi$ to be injective on $A_\x(k)\cup Y_k$. Thus $A_\x(k)\cup Y_k$ is linearly independent modulo $$\ker \pi=\R^\cup\la \partial(\mathcal A_\x)\ra=\R\la \darrow \wtilde{\partial(A_\x)}(k)\ra,$$ with the second  equality above in view of  Proposition \ref{prop-orReay-BasicProps}.1. As a result,  \eqref{tonicA} now implies that  $A_\x(k)\cup \darrow\wtilde{\partial(  A_\x)}(k)\cup Y_k$ is linearly independent, completing the claim and the proof.\end{proof}
\end{proof}

\section{Finitary Sets}\label{sec-finitary}

\subsection{Core Definitions and Properties} We have now developed the asymptotic framework generalizing the notion of positive basis to the point where we can define our main object of study later in the section. We begin by giving three equivalent definitions for the special subset $G_0^\diamond\subseteq G_0\subseteq \R^d$, which plays a key role  in our characterization result for finite elasticities. Note the three equivalent conditions defining $G_0^\diamond$ are all dependent only on notions from Convex Geometry---involving linear combinations over $\R_+$---rather than combinatorial properties dealing with $\mathcal A(G_0)$ and  linear combinations over $\Z_+$. When we impose additional conditions on $G_0$, we will later be able to give two further equivalent definitions for $G_0^\diamond$, first in terms of $\mathcal A^{\mathsf{elm}}(G_0)$ and linear combinations over $\Q_+$, and then in terms of $\mathcal A(G_0)$ and linear combinations over $\Z_+$.
Recall that $G_0^{\mathsf{lim}}$ was defined in Section \ref{sec-asym-seq-1}. Note the equality $\C(G_0\cap \mathcal E)=\mathcal E$ mentioned in the definition below follows by considering $\C(\pi(G_0))$, where $\pi:\R^d\rightarrow \mathcal E^\bot$ is the orthogonal projection, which must have trivial lineality space as $\mathcal E$ is the maximal subspace contained in $\C(G_0)$, thus ensuring that $0\notin \C^*(\pi(G_0)\setminus \{0\})$, which means any positive linear combination of elements from $G_0$ that lies in $\mathcal E$ must have all it elements lying in $\mathcal E$.

\begin{definition}
Let $G_0\subseteq \R^d$ be a subset with lineality space  $\mathcal E=\C(G_0)\cap -\C(G_0)$, so $\C(G_0\cap \mathcal E)=\mathcal E$.  Then we let   $G_0^\diamond\subseteq G_0\cap \mathcal E\subseteq G_0$ denote the subset of elements $g\in G_0\cap \mathcal E$ satisfying the equivalent conditions of Proposition \ref{prop-G_0diamond-1st-easy-equiv} (applied to  $G_0\cap \mathcal E$ in place of $G_0$).
\end{definition}

\begin{proposition}\label{prop-G_0diamond-1st-easy-equiv}
Let $G_0\subseteq \R^d$ be a subset with $\C(G_0)=\mathcal E\subseteq \R^d$ a subspace, where $d\geq 0$, and let $g\in G_0$. The following are equivalent.
 \begin{itemize}
 \item[1.]  There exists a subset $X\subseteq G_0$ and $\vec u\in G_0^{\mathsf{lim}}$   with $g\in X$ and $X$ minimally encasing  $-\vec u$.
 \item[2.] There exists a linearly independent subset $X\subseteq G_0$ with $g\in X$ and a sequence $\{x_i\}_{i=1}^\infty$ of terms $x_i\in G_0\cap \C^\circ(-X)$ such that $-x_i=\Summ{x\in X}\alpha_i^{(x)}x$ with $\alpha_i^{(x)}>0$ and $\alpha_i^{(g)}\rightarrow \infty$.
 \item[3.] There exists a linearly independent set $X\subseteq G_0$ with $g\in X$  and   $\C(X)\cap -G_0$ not bound to $\C(X\setminus \{g\})$.
 \end{itemize}
\end{proposition}

\begin{proof}
We may w.l.o.g. assume $\mathcal E=\R^d$.

 1. $\Rightarrow$ 2. Suppose there exists an asymptotically filtered sequence $\{x_i\}_{i=1}^\infty$ of terms $x_i\in G_0$ with fully unbounded limit $\vec u=(u_1,\ldots,u_t)$ and a subset $X\subseteq G_0$ with $g\in X$ such that $X$ minimally encases $-\vec u$. Write each $$x_i=a_1^{(1)}u_1+\ldots+a_i^{(t)}u_t+y_i$$ with $a_i^{(j)}>0$ and $y_i\in \R\la u_1,\ldots,u_t\ra^\bot$ such that $\|y_i\|\in o(a_i^{(t)})$, and let $x_i^{(t)}=x_i-y_i=a_1^{(1)}u_1+\ldots+a_i^{(t)}u_t$ be the truncated terms for $i\geq 1$.  By hypothesis, $t\geq 1$ and $a_i^{(t)}\rightarrow \infty$. Let $\pi:\R^d\rightarrow \R\la X\ra^\bot$ be the orthogonal projection.
 In view of Proposition \ref{prop-min-encasement-minposbasis}, by removing the first few terms, we can w.l.o.g. assume $X\cup \{x^{(t)}_i\}$ is a minimal positive basis for all $i$.
 Since $\C(G_0)=\R^d$, we can find a subset $Y\subseteq G_0$ such that $|\pi(Y)|=|Y|$ and $\pi(Y)$ is a positive basis for $\R\la X\ra^\bot$ corresponding to a Reay system $\mathcal R_{\pi(Y)}$.
 For each $i\geq 1$, let $Y_i\subseteq Y$ be the subset with $\pi(Y_i)=\supp_{\mathcal R_\pi(Y)}(-\pi(x_i))$, so $\pi(Y_i)\cup \{\pi(x_i)\}$ is a minimal positive basis or $Y_i=\emptyset$ with $\pi(x_i)=0$.
 By passing to a subsequence, we can assume all $Y_i$ are equal, say $Y_i=Y_1$ for all $i\geq 1$.
 By Proposition \ref{prop-char-minimal-pos-basis}.3, $X$ and $\pi(Y_1)$ are both  linearly independent. As a result, since $\pi$ is injective on $Y_1\subseteq Y$ with $\ker \pi=\R\la X\ra$, it follows that $X\cup Y_1$ is a linearly independent subset.
 Since $\pi(Y_1)\cup \{\pi(x_i)\}$ is a minimal positive basis with $\pi(x_i)=\pi(y_i)$, Lemma \ref{lemma-matrix-assymp-solutions} implies  \be\label{oneone}\Summ{y\in Y_1}\alpha_i^{(y)}y=-y_i+z_i\ee for some $z_i\in \R\la X\ra$ and $\alpha_i^{(y)}>0$ with $\alpha_i^{(y)}\in O(\|\pi(y_i)\|)\subseteq O(\|y_i\|)\subseteq o(a_i^{(t)})$ ($z_i=y_i$ if $Y_1=\emptyset$ with $\pi(x_i)=0$). Hence $\|z_i\|\in o(a_i^{(t)})$ and $x_i-y_i+z_i=x_i^{(t)}+z_i\in \R\la X\ra$.
 Thus $x_i-y_i+z_i=x_i^{(t)}+z_i=a_1^{(1)}u_1+\ldots+a_i^{(t)}u_t+z_i$ is an asymptotically filtered sequence with limit $(u_1,\ldots,u_t)$ (once all $i$ are sufficiently large).
 Applying Proposition \ref{prop-min-encasement-minposbasis} to $\{x_i-y_i+z_i\}_{i=1}^\infty$, we conclude that
 \be\label{twotwo}\Summ{x\in X}\alpha_i^{(x)}x=-x_i+y_i-z_i\ee for some $\alpha_i^{(x)}>0$ (passing to sufficiently large index terms). Moreover, for each $x\in X$, we have $\alpha_i^{(x)}\in \Theta(a_i^{(j)})$ for some $j\in [1,t]$.
 Since  $g\in X$ and  $a_i^{(j)}\rightarrow \infty$ for all $j\leq t$, we have $\alpha_i^{(g)}\rightarrow \infty$, whence Item  2 follows from \eqref{oneone} and \eqref{twotwo} taking $X$ to be $X\cup Y_1$.

2. $\Rightarrow$ 3.  Let $X\subseteq G_0$ and $\{x_i\}_{i=1}^\infty$ be as given by Item 2, so $x_i\in G_0$ for all $i\geq 1$. Then each $-x_i=\Summ{x\in X}\alpha_i^{(x)}x\in\C^\circ(X)\cap -G_0$  with $\alpha_i^{(x)}>0$ and $\alpha_i^{(g)}\rightarrow \infty$.  Assuming by contradiction that $\{-x_i\}_{i=1}^\infty$ is bound to $\C(X\setminus \{g\})$, then there is a bound $M>0$ such that each  $-x_i$ has some $z_i\in \C(X\setminus \{g\})$ with $\|-x_i-z_i\|\leq M$. Let $T:\R^d\rightarrow \R$ be a linear transformation that sends $g$ to $1$ and $\R\la X\setminus \{g\}\ra$ to $0$,  which exists since $g\in X$ with $X$ linearly independent. Then $T(-x_i-z_i)=T(-x_i)=\alpha_i^{(g)}$ with $ \alpha_i^{(g)}=\|T(-x_i-z_i)\|\leq C_T \|-x_i-z_i\|\leq C_TM<\infty $, where $C_T$ is the operator norm of $T$ with respect to the Euclidean metric,  contradicting that $a_i^{(g)}\rightarrow \infty$.

3. $\Rightarrow$ 1. Let $X\subseteq G_0$ be as given by Item 3. Since $\C(X)\cap -G_0$ is not bound to $\C(X\setminus \{g\})$, we can find a sequence $\{-x_i\}_{i=1}^\infty$ of terms $-x_i\in \C(X)\cap -G_0$ such that $\mathsf d(-x_i,\C(X\setminus \{g\}))\rightarrow \infty$.  By passing to a subsequence, we can assume $\{x_i\}_{i=1}^\infty$ is an asymptotically filtered sequence with  complete fully unbounded limit $\vec u=(u_1,\ldots,u_t)$, where $t\geq 1$. Thus $\vec u\in G_0^{\mathsf{lim}}$ since $\vec u$ is fully unbounded and $x_i\in G_0$ for all $i\geq 1$. Since $-x_i\in \C(X)$ for all $i\geq 1$,  Proposition \ref{prop-finite-union--convergence-encasement} implies that $\C(X)$ encases $-\vec u$, so let $Y\subseteq X$ be a subset for which $Y$ minimally encases $-\vec u$. If $g\in Y$, then Item 1 follows. Otherwise, $X\setminus \{g\}$ encases $-\vec u$, in which case Proposition \ref{prop-encasementcones-contain-aprox-seq}.3 implies that $\{-x_i\}_{i=1}^\infty$ is bound to  $\C(X\setminus\{g\})$, contradicting that $\mathsf d(-x_i,\C(X\setminus \{g\}))\rightarrow \infty$.
\end{proof}

We continue with the following basic inclusion for $G_0^\diamond$ for subsets $G_0$ of lattice points.

\begin{proposition}\label{prop-diamond-basic-containment}
Let $\Lambda\subseteq \R^d$ be a full rank lattice, where $d\geq 0$, and let
$G_0\subseteq \Lambda$ be a subset with $\C(G_0)=\R^d$.  Then
$$G_0^\diamond\subseteq\{g\in G_0:\; \sup\{\vp_{g}(U):\; U\in \mathcal A^{\mathsf{elm}}(G_0)\}=\infty\}.$$
\end{proposition}

\begin{proof}

Let $g\in G_0^\diamond$ be arbitrary.  By Proposition \ref{prop-G_0diamond-1st-easy-equiv}.2, there exists a linearly independent subset $X\subseteq G_0$ with $g\in X$ and a sequence $\{x_i\}_{i=1}^\infty$ of terms $x_i\in G_0\cap -\C^\circ(X)$ such that \be\label{eaglet}-x_i=\Summ{x\in X}\alpha_i^{(x)}x\ee with $\alpha_i^{(x)}>0$ and $\alpha_i^{(g)}\rightarrow \infty$.
By Proposition \ref{prop-char-minimal-pos-basis} (Items 4 and 7), each $X\cup\{x_i\}$ is a minimal positive basis and there exists an elementary atom $U_i$ with $\supp(U_i)=X\cup \{x_i\}$. Since the vector
$(\vp_x(U_i))_{x\in X\cup \{x_i\}}$  has $\Summ{x\in X\cup \{x_i\}}\vp_x(U_i)x=\sigma(U_i)=0$, it follows from
 Proposition \ref{prop-char-minimal-pos-basis}.5 and \eqref{eaglet} that $\vp_x(U_i)=\vp_{x_i}(U_i)\alpha_i^{(x)}$ for all $x\in X$. In particular, $\vp_g(U_i)=\vp_{x_i}(U_i)\alpha_i^{(g)}\geq \alpha_i^{(g)}$ (the inequality follows in view of $x_i\in \supp(U_i)$) with $\alpha_i^{(g)}\rightarrow \infty$, showing that $ \sup\{\vp_{g}(U):\; U\in \mathcal A^{\mathsf{elm}}(G_0)\}=\infty$. This establishes the desired inclusion.
\end{proof}

Lemma \ref{Lemma-VReay-RidRemainders} links  the set $\mathcal A_j$ from Proposition \ref{prop-VReay-RepBasics} with the diamond subset $G_0^\diamond$.

\begin{lemma}
\label{Lemma-VReay-RidRemainders} Let $\Lambda\subseteq \R^d$ be a full rank lattice, where $d\geq 1$, let $G_0\subseteq \Lambda$ be a subset of lattice points with $\C(G_0)=\R^d$,  let $\mathcal R=(\mathcal X_1\cup \{\mathbf v_1\},\ldots,\mathcal X_s\cup\{\mathbf v_s\})$ be a virtual Reay system in $G_0$,  let $\mathcal A_s$ be the subset given by Proposition \ref{prop-VReay-SupportSet}.1 for $j=s$, suppose   that the virtual Reay system $\mathcal R'=(\mathcal X_1\cup \{\mathbf v_1\},\ldots,\mathcal X_{s-1}\cup \{\mathbf v_{s-1}\})$ is anchored, and let $\x\in \mathcal A_s$.

\begin{itemize}

\item[1.] If $\vec u_\x$ is fully unbounded and $\y\in \mathcal A_s\setminus \{\x\}$,  then there exists  an asymptotically filtered sequence $\{y_i\}_{i=1}^\infty$ of terms $y_i\in \C(G_0^\diamond)\cap \R^\cup\la \mathcal A_s\ra$ with limit $\vec u_\y$.  Indeed, there is a finite subset $Z\subseteq G_0^\diamond$ such that $y_i\in \y(i)+\C(Z)$ for all sufficiently large $i$, and $\z(i)\in  G_0^\diamond$ for every $\z\in \darrow \mathcal A_s\setminus \{\x\}$ and all sufficiently large $i$.
\item[2.] Suppose, for each $\y\in \mathcal A_s\setminus\{\x\}$, that $\{y_{i_\y}^\y\}_{i_\y=1}^\infty$ is an asymptotically filtered sequence of terms $y_{i_\y}^\y\in \R^\cup \la \mathcal A_s\ra$ with limit $\vec u_\y$. Then, for any tuple $k=(i_\z)_{\z\in \darrow \mathcal A_s\setminus \{\x\}}$ with all $i_\z$ sufficiently large, the set $\darrow \partial(A_s)(k)\cup \{y^\y_{i_\y}:\; \y\in \mathcal A_s\setminus \{\x\}\}$ minimally encases $-\vec u_\x$ and  $\R^\cup\la \mathcal A_s\ra=\R\la \darrow \partial(A_s)(k)\cup \{y^\y_{i_\y}:\; \y\in \mathcal A_s\setminus \{\x\}\}\ra$.
    \end{itemize}
\end{lemma}

\begin{proof}
1. Let $\vec u_\x=(u_1^\x,\ldots,u_{t_\x}^\x)$ and $\vec u_\y=(u_1^\y,\ldots,u_{t_\y}^\y)$.
Let $\pi:\R^d\rightarrow \R\la u^\y_1,\ldots,u^\y_{t_\y}\ra^\bot$ be the orthogonal projection. Since $\{\y(i)\}_{i=1}^\infty$ is an asymptotically filtered sequence with limit $\vec u_\y=(u^\y_1,\ldots,u^\y_{t_\y})$, we have $$\y(i)=a_i^{(1)}u_1^\y+\ldots+a_i^{(t_\y)}u^\y_{t_\y}+\pi(\y(i))\quad\mbox{for all $i$},$$ where \be\label{gooseduck}a_i^{(j)}\in o(a_i^{(j-1)})\quad\mbox{ for $j\in [2,t_\y]$}\quad\und\quad \|\pi(\y(i))\|\in o(a_i^{(t_\y)}).\ee

Our general strategy is as follows. We will partition the  terms in $\y(i)$ into a finite number of infinite subsequences, say $I_1\cup \ldots\cup I_r=\Z_+\setminus\{0\}$ with this union disjoint and each $I_j$ infinite. We will  show Item 1 holds for the sufficiently large index terms in each $\{\y(i)\}_{i\in I_j}$ and, additionally,  it  does so with each $y_i=\tilde \y(i)+\xi_i=\y(i)-\pi(\y(i))+\xi_i$, for $i\in I_j$, such that  $\|\xi_i\|\in o(a_i^{(t_\y)})$, \ $\xi_i\in \R^\cup \la \mathcal A_s\ra$ and $-\pi(\y(i))+\xi_i\in \C(Z_j)$ for some fixed, finite subset $Z_j\subseteq G_0^\diamond$. Then Item 1 will hold for the sufficiently large index terms in $\{\y(i)\}_{i=1}^\infty$ by setting $Z=\bigcup_{j=1}^r Z_j$ and using the sequence $\{y_i\}_{i=1}^\infty$, where each $y_i$ for $i\in I_j$ was the term defined for the subsequence $\{\y(i)\}_{i\in I_j}$.

Let $Y\subseteq G_0$ be the subset given by Proposition \ref{prop-VReay-SupportSet}.3 for $j=s$.
Let $k=(i_\z)_{\z\in \darrow \mathcal A_\x}$ be a tuple of indices. In view of Proposition \ref{prop-VReay-SupportSet}.3, once all $i_\z$ are sufficiently large, then we can assume that $( A_s\setminus\{x\})(k)\cup \darrow \wtilde{\partial(A_s)}(k)\cup Y_k$ minimally encases $-\vec u_\x$ for some $Y_k\subseteq Y$. Since $\mathcal R'$ is anchored, Proposition \ref{prop-VReay-Lattice} implies that $\tilde \z(i_\z)=\z(i_\z)$ for all $\z\in \darrow \partial(\mathcal A_s)$ once $i_\z$ is sufficiently large. Thus $( A_s\setminus\{x\})(k)\cup \darrow \partial(A_j)(k)\cup Y_k$ minimally encases $-\vec u_\x$ (once all $i_\z$ are sufficiently large).
As a result, since  $\vec u_\x$ is fully unbounded, it follows that  $\z(i_\z)\in G_0^\diamond\subseteq G_0$ for all $\z\in (\mathcal A_s\setminus \{\x\})\cup \darrow \partial(\mathcal A_s)=\darrow \mathcal A_s\setminus \{\x\}$, once all $i_\z$ are sufficiently large.
Now fix all indices $i_\z$ with $\z\neq \y$ (chosen sufficiently large that Proposition \ref{prop-VReay-SupportSet}.3 and Proposition \ref{prop-VReay-Lattice} are  applicable) and consider $i=i_\y\rightarrow \infty$. Since $Y$ is finite, there are only a finite number of possibilities for the $Y_k$ as $i\rightarrow\infty$.
Thus (as described at the start of the proof), by passing to a subsequence of $\{\y(i)\}_{i=1}^\infty$, we can w.l.o.g. assume the same set $Y_k$ occurs for every $i$, in which case $\{\y(i)\}\cup Z$ minimally encases $-\vec u_\x$ for all $i$, where $Z:=(A_s\setminus \{x,\,y\})(k)\cup\darrow \partial(A_j)(k)\cup Y_k\subseteq G_0$ is a fixed subset. Since $\vec u_\x$ is fully unbounded, we have  $$\{\y(i)\}\cup Z\subseteq G_0^\diamond$$ for all $i$.
By Proposition \ref{prop-min-encasement-char}, for each $i$, there are indices $1=r_1<\ldots<r_{\ell+1}=t_\x+1$ and a disjoint partition $\{\y(i)\}\cup Z=Z_1\cup\ldots\cup Z_{\ell}$ such that $(Z_1\cup \{u^\x_{r_1}\},\ldots,Z_\ell\cup \{u^\x_{r_{\ell}}\})$ is a Reay system.
By Proposition \ref{prop-reay-RayAlg}, for each $i$, there is some subset $Z'\subseteq Z\cup \{u^\x_1,\ldots,u^\x_{t_\x}\}$ such that $Z'\cup \{\y(i)\}$ is a minimal positive basis. Since there are only a finite number of possibilities for $Z'$, by passing to a subsequence of $\{\y(i)\}_{i=1}^\infty$, we can assume the same set $Z'$ occurs for every $i$ (using the same partitioning argument as before). Let $J\subseteq [1,t_\x]$ be the set of indices $j$ for which $u_j^\x\in Z'$ and let $U=\{u^\x_j:\;j\in J\}$. Then $$Z'\setminus U\subseteq Z\subseteq G_0^\diamond.$$

  In view of (VR1), we have $u_1^\y,\ldots,u_{t_\y}^\y\in \R\la \y\ra\subseteq \R^\cup \la \mathcal A_s\ra$. Hence, if $\pi(\y(i))=0$ for an infinite number of $i$, then Item 1 follows (on the partition of indices in $I$ having this property) taking $\{y_i\}_{i=1}^\infty$ to be a subsequence of $\{\y(i)\}_{i=1}^\infty$ with all terms having $\xi_i=\pi(\y(i))=0$, in view of $\y(i)\in G_0^\diamond $ for all $i$.
 Assuming this does not occur (so we restrict to  a sub-partition of  $I$ where $\pi(\y(i))= 0$ holds only for a finite number of indices  $i$), by discarding the first few terms in the sequence $\{\y(i)\}_{i=1}^\infty$, we can assume $\pi(\y(i))\neq 0$ for all $i$. Thus,
since $Z'\cup \{\y(i)\}$ is a minimal positive basis, it follows by Carath\'eordory's Theorem that $\pi(Z'')\cup \{\pi(\y(i))\}$ is a minimal positive basis of size $|Z''|+1$ for some subset $Z''\subseteq Z'$. As there are only a finite number of possibilities for $Z''$, by once more passing to a subsequence of $\{\y(i)\}_{i=1}^\infty$, we can assume the same set $Z''$ occurs for every $i$ (by the same partitioning argument as before). As a result, $-\pi(\y(i))\in \C^\circ(\pi(Z''))$ for all $i$, in which case Lemma \ref{lemma-matrix-assymp-solutions} implies that $$\Summ{z\in Z''\setminus U}\alpha_i^{z}\pi(z)+\Summ{z\in Z''\cap U}\beta_i^{z}\pi(z)=-\pi(\y(i))$$ for some $\alpha_i^{z}>0$ and $\beta_i^{z}>0$ with $\alpha_i^{z},\,\beta_i^{z}\in O(\|\pi(\y(i))\|)$ for all  $z\in Z''$. Thus \be\nn\Summ{z\in Z''\setminus U}\alpha_i^{z}z=-\pi(\y(i))+\xi_i\ee for some $\xi_i\in\
\ker \pi+\R\la U\ra\subseteq \R\la u^\x_1,\ldots,u^\x_{t_\x},u_1^{\y},\ldots,u_{t_\y}^\y\ra$ with $\|\xi_i\|\in O(\|\pi(\y(i))\|)\subseteq o(a_i^{(t_\y)})$, with the latter inclusion in view of \eqref{gooseduck}.
Setting $y_i=\Summ{z\in Z''\setminus U}\alpha_i^{z}z+\y(i)$, we have $$y_i=a_i^{(1)}u_1^\y+\ldots+a_i^{(t_\y)}u^\y_{t_\y}+\xi_i$$ and, by discarding the first few terms in $\{y_i\}_{i=1}^\infty$, we find that   $\{y_i\}_{i=1}^\infty$ is an asymptotically  filtered sequence with limit $\vec u_\y$ in view of $\|\xi_i\|\in o(a_i^{(t_\y)})$. Since $Z''\setminus U\subseteq Z'\setminus U\subseteq G_0^\diamond$ and $\y(i)\in G_0^\diamond$, it follows that $y_i\in \C(G_0^\diamond)$ for all $i$. By (VR1), we know $u_1^\y,\ldots,u_{t_\y}^\y,u_1^\x,\ldots,u_{t_\x}^{t_\x}\in \R\la \x\cup \y\ra\subseteq \R^\cup\la \mathcal A_s\ra$. Consequently, since $\xi_i\in \R\la u^\x_1,\ldots,u^\x_{t_\x},u_1^{\y},\ldots,u_{t_\y}^\y\ra$, it follows that $y_i\in \R\la \x\cup\y\ra\subseteq \R^\cup \la \mathcal A_s\ra$ for all $i$, and the sequence $\{y_i\}_{i=1}^\infty$ now has the desired properties (using $Z''\setminus U$ for $Z$), completing the proof of Item 1.

2. Define a new virtual Reay system $\wtilde{\mathcal R}$ on the set  $$G_0\cup \{y^\y_{i_\y}:\: \y\in \mathcal A_s\setminus \{\x\}\und i_\y\geq 1\}\subseteq \R^d$$ identical to $\mathcal R$ except that we replace each sequence $\{\y(i_\y)\}_{i_\y=1}^\infty$ with the  sequence $\{y_{i_\y}^\y\}_{i_\y=1}^\infty$ for $\y\in \mathcal A_s\setminus\{\x\}$. Since both these sequences are asymptotically  filtered sequences with the same limit for every $\y\in \mathcal A_s\setminus\{\x\}$,  this leaves all half-spaces and their boundaries unchanged. In particular, $\mathcal A_s$ is also the set given by Proposition \ref{prop-VReay-SupportSet}.1 for $\wtilde{\mathcal R}$ as well as $\mathcal R$. Consequently, since $y^\y_{i_\y}\in \R^\cup\la \mathcal A_s\ra$ for all $i_\y$ and $\y\in \mathcal A_s\setminus \{\x\}$ by hypothesis, Proposition \ref{prop-VReay-SupportSet}.3 applied to $\wtilde{\mathcal R}$ implies that $\wtilde{\darrow \partial(A_s)}(k)\cup \{y^\y_{i_\y}:\; \y\in \mathcal A_s\setminus \{\x\}\}$ minimally encases $-\vec u_\x$ once all $i_\z$ in the tuple $k=(i_\z)_{\z\in \darrow \mathcal A_s\setminus \{\x\}}$ are sufficiently large. However, since  $G_0\subseteq \Lambda$ and $\mathcal R'$ is anchored, Proposition \ref{prop-VReay-Lattice} implies that  $\darrow \wtilde {\partial(A_s)}(k)=\darrow  \partial(A_s)(k)$ once all $i_\z$ are sufficiently large, and thus $\darrow \partial(A_s)(k)\cup \{y^\y_{i_\y}:\; \y\in \mathcal A_s\setminus \{\x\}\}$ minimally encases $-\vec u_\x$ once all $i_\z$ are sufficiently large.

 By definition, each $y_{i_\y}^\y=z_{i_\y}^\y+\epsilon_{i_\y}^\y$ with $\{z_{i_\y}^\y\}_{i_\y=1}^\infty$ a complete asymptotically  filtered sequence with limit $\vec u_\y=(u_1^\y,\ldots,u_{t_\y}^\y)$ and $\epsilon_{i_\y}^\y\in \R^\cup \la \mathcal A_s\ra$ the remainder term, so $z_{i_\y}^\y=a_{i_\y,\y}^{(1)}u_1^\y+\ldots+a_{i_\y,\y}^{(t_\y)}u^\y_{t_\y}$ with $\|\epsilon_{i_\y}^\y\|\in o(a_{i_\y,\y}^{(t_\y)})$.
Proposition \ref{prop-orReay-BasicProps}.1 and Proposition \ref{prop-VReay-SupportSet}.2(e) imply that $$\R^\cup \la \partial(\mathcal A_s)\ra=\R^\cup\la \darrow \partial(A_s)(k)\ra\quad\und\quad \R^\cup\la \mathcal A_s\ra=\R\la \darrow \partial(A_s)(k)\cup \{z^\y_{i_\y}:\; \y\in \mathcal A_s\setminus \{\x\}\}\ra.$$ Thus, letting $\pi:\R^d\rightarrow \R^\cup\la \partial(\mathcal A_s)\ra^\bot$ be the orthogonal projection, it follows that $$\R^\cup \la \pi(\mathcal A_s)\ra=\R\la \{\pi(z^\y_{i_\y}):\; \y\in \mathcal A_s\setminus \{\x\}\}\ra.$$
Each $z_{i_\y}^\y $ is  a representative for the half-space $\y\in \mathcal A_s\setminus\{\x\}$, so Proposition \ref{prop-VReay-SupportSet}.2(e) implies that   $$\{\pi(z^\y_{i_\y}):\; \y\in \mathcal A_s\setminus \{\x\}\}$$ is linearly independent, and thus a basis of size $|\mathcal A_s\setminus\{\x\}|$ for the  subspace $\R^\cup \la \pi(\mathcal A_s)\ra$.
In view of Proposition \ref{prop-VReay-SupportSet}.2(e), the value $\pi(z_{i_\y}^\y)/\|\pi(z_{i_\y}^\y)\|=\pi(u_{t_\y}^\y)/\|\pi(u_{t_\y}^\y)\|$ is constant, where  $\vec u_\y=(u_1^\y,\ldots,u_{t_\y}^\y)$, so $\|\pi(z_i)\|=\Theta(a_{i_\y,\y}^{(t_\y)})$.
In consequence, each $\{\pi(y_{i_\y}^\y)\}_{i_\y=1}^\infty$ for $\y\in \mathcal A_s\setminus\{\x\}$ is a sequence of terms from $\R^\cup\la \pi(\mathcal A_s)\ra$ that radially converges to $\pi(z_{i_\y}^\y)/\|\pi(z_{i_\y}^\y)\|$,
and thus Proposition \ref{prop-FanStability}.5 implies that, once all $i_\y$ are sufficiently large,  $\{\pi(y^\y_{i_\y}):\; \y\in \mathcal A_s\setminus \{\x\}\}$ is also a basis of size $|\mathcal A_s\setminus\{\x\}|$ for the  subspace $\R^\cup \la \pi(\mathcal A_s)\ra$, which, combined with $\R^\cup \la \partial(\mathcal A_s)\ra=\R^\cup\la \darrow \partial(A_s)(k)\ra$,  gives the desired conclusion $\R^\cup\la \mathcal A_s\ra=\R\la \darrow \partial(A_s)(k)\cup \{y^\y_{i_\y}:\; \y\in \mathcal A_s\setminus \{\x\}\}\ra$, completing the proof.
\end{proof}

We are now in position to give the \emph{key} definition that will be the subject of the remainder of this section. As we will later see, finitary sets form  an ample class of subsets of lattice points that behave well under inductive arguments and include infinite subsets that nonetheless behave like finite sets, particularly in relation to combinatorial properties related to $\mathcal A(G_0)$.

\begin{definition}
Let $\Lambda\leq \R^d$ be a full rank lattice, where $d\geq 0$, and let $G_0\subseteq \Lambda$. If every purely virtual Reay system over $G_0$ is anchored, we say that $G_0$ is \textbf{finitary}.

Let $t\in [0,d]$ be a integer. If every purely virtual Reay system $\mathcal R=(\mathcal X_1\cup \{\mathbf v_1\},
\ldots,\mathcal X_s\cup \{\mathbf v_s\})$ over $G_0$ with $\dim \R^\cup\la \mathcal X_1\cup \ldots\cup \mathcal X_s\ra\leq t$ is anchored, then we say that $G_0$ is \textbf{finitary  up to dimension $t$}. (We will need this refined definition for several inductive arguments.)
\end{definition}

As a point of clarification, if $\mathcal R=(\mathcal X_1\cup \{\mathbf v_1\},\ldots,\mathcal X_s\cup \{\mathbf v_s\})$ is a virtual Reay system in $G_0$, then each $\x\in \mathcal X_1\cup\ldots\cup \mathcal X_s$ has various possibilities for the asymptotically filtered sequence $\{\x(i)\}_{i=1}^\infty$ having limit $\vec u_\x$. It is possible for some of these sequences to have $\vec u_\x$ as fully unbounded limit, and some to have $\vec u_\x$ as an anchored limit. In such case, it is possible to consider $\vec u_\x$ either as fully unbounded or anchored, meaning the virtual Reay system $\mathcal R$ can be considered anchored or not, depending on which sequences  $\{\x(i)\}_{i=1}^\infty$ are chosen. The definition of finitary requires that all choices of sequences $\{\x(i)\}_{i=1}^\infty$ result in $\mathcal R$ being anchored. A similar fact holds regarding purely virtual Reay systems, though here we only need each $\mathbf v_j$ to have \emph{some}  sequence $\{\mathbf v_j(i)\}_{i=1}^\infty$ having $\vec u_{\mathbf v_j}$ as fully unbounded limit in order to consider $\mathcal R$ as purely virtual.

We now make the observation that if $\mathcal R=(\mathcal X_1\cup\{\mathbf v_1\},\ldots,\mathcal X_s\cup \{\mathbf x_s\})$ is a virtual Reay system in $G_0$, then so too is $\mathcal R=(\mathcal X_1\cup\{\mathbf v_1\},\ldots,\mathcal X_{s-1}\cup \{\mathbf v_{s-1}\},((\mathcal X_s\cup \{\mathbf x_s\})\setminus\{\y\})\cup \{\y\})$ for any $\y\in \mathcal X_s\cup \{\mathbf v_s\}$. This observation will allow us to apply many of the propositions that require $\mathcal B\subseteq \mathcal X_1\cup\ldots\cup \mathcal X_s$ as a hypothesis when we only have $\mathcal B\subseteq (\mathcal X_1\cup\ldots\cup \mathcal X_{s})\cup \{\mathbf v_s\}$ with $\mathcal X_s\cup\{\mathbf v_s\}\nsubseteq \mathcal B$, at least in special circumstances (in theory, many of these propositions could have been stated in this more general form, but the trick above means the added generality is implied by applying the more limited form to a modified oriented Reay system).

We continue by giving some important properties of finitary sets $G_0$. Recall that $\wt(-g)$ was defined in \eqref{weight-deff}.

\begin{proposition}\label{prop-finitary-basics}
Let $\Lambda\subseteq \R^d$ be a full rank lattice, where $d\geq 0$, let
$G_0\subseteq \Lambda$ be a subset with $\C(G_0)=\R^d$, and let  $\mathcal R=(\mathcal X_1\cup \{\mathbf v_1\},\ldots,\mathcal X_s\cup\{\mathbf v_s\})$ be a  purely virtual Reay system in $G_0$,  where $s\geq 0$. Suppose $G_0$ is finitary up to dimension $\dim\R^\cup\la \mathcal X_1\cup\ldots\cup\mathcal X_s\ra$.
\begin{itemize}
\item[1.] $\x(i)\in G_0^\diamond$ for every $\x\in \mathcal X_1\cup\ldots\cup \mathcal X_s$ and  all  sufficiently large $i$.
\item[2.] If $\mathcal A_j$ is the subset given by Proposition \ref{prop-VReay-SupportSet} for $j\in [1,s]$, then $\mathcal A_j\setminus \{\mathbf v_j\}\subseteq \mathcal X_1\cup\ldots\cup \mathcal X_j$.
\item[3.] If  $\vec u=(u_1,\ldots,u_t)\in G_0^{\mathsf{lim}}$  with $u_1,\ldots,u_t\in \R^\cup \la \mathcal X_1\cup\ldots\cup \mathcal X_s\ra$, then $-\vec u$ is  encased by $\mathcal X_1\cup \ldots\cup \mathcal X_s$.
\item[4.] If  $\vec u=(u_1,\ldots,u_t)$ is an anchored limit of an asymptotically filtered sequence of terms from $G_0$ with $u_1,\ldots,u_t\in \R^\cup \la \mathcal X_1\cup \ldots\cup \mathcal X_s\ra$, then $-\vec u$ is minimally encased urbanely by a support subset $\mathcal B\subseteq \mathcal X_1\cup \{\mathbf v_1\}\cup \ldots\cup \mathcal X_s\cup\{\mathbf v_s\}$ with $|\mathcal B\cap \{\mathbf v_1,\ldots,\mathbf v_s\}|\leq 1$. In particular, $\mathsf{wt}(-g)\leq 1$ for all $g\in G_0\cap \R^\cup \la \mathcal X_1\cup\ldots\cup \mathcal X_s\ra$.
\end{itemize}
\end{proposition}

\begin{proof}
Let $\mathcal X=\mathcal X_1\cup \ldots\cup \mathcal X_s$.
 If  $\dim \R^\cup \la \mathcal X\ra=0$, then $s=0$ and all items hold vacuously. So we may assume $\dim \R^\cup \la \mathcal X\ra>0$ (implying $s\geq 1$) and proceed by induction on $\dim \R^\cup \la \mathcal X\ra$.  Let $\mathcal R'_0=(\mathcal X_1\cup \{\mathbf v_1\},\ldots,\mathcal X_{s-1}\cup \{\mathbf v_{s-1}\})$. Applying the induction hypothesis to $\mathcal R'_0$, we find that $(X_1\cup\ldots\cup X_{s-1})(k)\subseteq G_0^\diamond$ for any tuple $k=(i_\z)$ with all $i_\z$ sufficiently large, $\mathcal A_j\setminus \{\mathbf v_j\}\subseteq \mathcal X_1\cup \ldots\cup\mathcal X_j$ for all $j\in [1,s-1]$, and Items 3 and 4 both hold for $\mathcal R'_0$.

1. Since $G_0$ is finitary up to dimension $\dim\R^\cup\la \mathcal X_1\cup\ldots\cup\mathcal X_s\ra$, it follows that $\mathcal R$ is anchored, meaning $\vec u_\x$ is anchored for every $\x\in \mathcal X_1\cup\ldots\cup \mathcal X_{s}$.
Since $\mathcal R$ is purely virtual,   $\vec u_{\mathbf v_s}$ is fully unbounded.
Thus Lemma \ref{Lemma-VReay-RidRemainders}.1 implies that $\x(i)\in G_0^\diamond$ for every $\x\in \darrow \mathcal A_s\setminus \{\mathbf v_s\}$ and all sufficiently large $i$. Since $\mathcal X_s\subseteq \mathcal A_s$, Item 1 now follows.

2. Suppose Item 2 fails.
Then $\mathcal A_s\setminus \{\mathbf v_s\}$ contains some $\mathbf v_j$ with $j\in [1,s-1]$.
Since $\mathcal R$ is purely virtual, $\vec u_{\mathbf v_j}$ is fully unbounded.
Let $\mathcal A_{\mathbf v_s}=(\mathcal A_s\setminus \{\mathbf v_s\}\cup \partial(\{\mathbf v_s\}))^*$.
Since $\mathbf v_j\in\mathcal A_s\setminus \{\mathbf v_s\}$ and every $\mathbf v_i$ is a maximal element in $\mathcal X_1\cup\ldots\cup \mathcal X_s\cup \{\mathbf v_1,\ldots,\mathbf v_s\}$, it follows that  $\mathbf v_j\in \mathcal A_{\mathbf v_s}$. Proposition \ref{prop-VReay-SupportSet}.2(b) implies that $\mathcal A_{\mathbf v_s}$ minimally encases $-\vec u_{\mathbf v_s}$, which is fully unbounded since $\mathcal R$ is purely virtual.
Proposition \ref{prop-VReay-SupportSet}.2(c) implies that $\mathcal A_{\mathbf v_s}$ is a support set.
Let $\vec u_{\mathbf v_s}=(u_1,\ldots,u_t)$.
Applying Lemma \ref{lem-urbane-encasement-guarantee} to $\mathcal A_{\mathbf v_s}$ and $\vec u_{\mathbf v_s}$, we find there is some $\mathcal A\preceq \mathcal A_{\mathbf v_s}$ and $t_0\in [1,t]$ such that $\mathcal A$ minimally encases $-(u_1,\ldots,u_{t_0})$ urbanely and contains some $\y\in \mathcal A$ with $\vec u_{\y}$ fully unbounded.
Since $\mathcal A\preceq \mathcal A_{\mathbf v_s}$ with $\mathcal A_{\mathbf v_s}$ a support set, it follows that $\mathcal A=\mathcal A^*$ is also a support set (note $\mathcal A^*\mathcal =\mathcal A$ since $\mathcal A$ minimally encases $-(u_1,\ldots,u_{t_0})$).
Apply Proposition \ref{prop-VReay-modularCompletion} to $\mathcal A$ and let $\mathcal R'=(\mathcal Y_1\cup \{\mathbf w_1\},\ldots,\mathcal Y_{s'}\cup \{\mathbf w_{s'}\})$ be the resulting virtual Reay system with $\mathcal Y_1\cup \ldots\cup \mathcal Y_{s'}=\darrow \mathcal A$.
Since $\vec u_{\mathbf v_s}=(u_1,\ldots,u_t)$ is fully unbounded and $t_0\geq 1$, it follows that $(u_1,\ldots,u_{t_0})$ is also fully unbounded.
Since $G_0\subseteq \Lambda$,
Proposition \ref{prop-VReay-Lattice} implies that every $\x\in \darrow \mathcal A$ has $\vec u_\x^\triangleleft$ either trivial or fully unbounded (since each $\vec u_{\mathbf v_i}$ is fully unbounded, this is trivially true when $\x=\mathbf v_i$, while all $\x\in \mathcal X_1\cup \ldots\cup \mathcal X_s$ have $\vec u_\x$ anchored). Thus $\mathcal R'$ is purely virtual as noted before Proposition \ref{prop-VReay-modularCompletion}. As a result, since $\darrow \mathcal A\subseteq \darrow \mathcal A_{\mathbf v_s}\subseteq \mathcal X_1\cup\{\mathbf v_1\}\cup\ldots\cup \mathcal X_s\cup \{\mathbf v_s\}$ implies $\dim \R^\cup\la \mathcal Y_1\cup\ldots\cup \mathcal Y_{s'}\ra\leq \dim \R^\cup\la \mathcal X_1\cup\ldots\cup \mathcal X_s\ra$, and since $G_0$ is finitary up to dimension $\dim \R^\cup\la \mathcal X_1\cup\ldots\cup \mathcal X_s\ra$, it follows that $\mathcal R'$ must be anchored. However, this contradicts that $\y\in \mathcal A\subseteq \darrow \mathcal A=\mathcal Y_1\cup\ldots\cup \mathcal Y_{s'}$ with $\vec u_{\y}$ fully unbounded, and Item 2 is now established.

3. To prove Item 3, suppose $\vec u=(u_1,\ldots,u_t)$ is a  fully unbounded limit of an asymptotically filtered sequence $\{x_i\}_{i=1}^\infty$  of terms $x_i\in G_0$ with $u_1,\ldots,u_t\in \R^\cup \la \mathcal X\ra=\C^\cup (\mathcal X\cup \{\mathbf v_1,\ldots,\mathbf v_s\})$ (the equality follows by Proposition \ref{Prop-orReay-coord}) such that $-\vec u$ is not encased by $\mathcal X$. By choosing such a counter-example with $t$ minimal, we can assume $-\vec u^\triangleleft$ is encased by $\mathcal X$.
Thus there is some $\mathcal A\subseteq \mathcal X$ which minimally encases $-\vec u^\triangleleft$ urbanely.
Let $\pi:\R^d\rightarrow \R^\cup\la\mathcal A\ra^\bot$ be the orthogonal projection.
Since $\mathcal X$ does not encase $-\vec u$, we must have $u_t\notin \R^\cup \la\mathcal A\ra$ (by Proposition \ref{prop-orReay-minecase-char}).
Let $\mathcal C\subseteq \mathcal X_1\cup \{\mathbf v_1\}\cup\ldots\cup \mathcal X_s\cup \{\mathbf v_s\}$ be the lift of $\supp_{\pi(\mathcal R)}(-\pi(u_t))$, which is a support set by Proposition \ref{prop-orReay-modulo}.5. Then Propositions \ref{prop-orReay-minecase-char}.2 and \ref{prop-orReay-minecase-char}.3(a) imply that $\mathcal C$ minimally encases $-\vec u$ urbanely (since $\mathcal A\subseteq \mathcal X$). Apply Proposition \ref{prop-VReay-modularCompletion} to $\mathcal C$ and let $\mathcal R'=(\mathcal Y_1\cup\{\mathbf w_1\},\ldots,\mathcal Y_r\cup \{\mathbf w_r\})$ be the resulting virtual Reay system in $G_0$ for the subspace $\R^\cup \la \mathcal C\ra\subseteq \R^\cup \la \mathcal X\ra$ with $\mathcal Y_r=\mathcal C$, $\vec u_{\mathbf w_r}=\vec u$ and $\mathcal Y_1\cup\ldots\cup \mathcal Y_r=\darrow \mathcal C$.
Since $G_0\subseteq \Lambda$,
Proposition \ref{prop-VReay-Lattice} implies that every $\x\in \darrow \mathcal C$ has $\vec u_\x^\triangleleft$ either trivial or fully unbounded (since each $\vec u_{\mathbf v_i}$ is fully unbounded, this is trivially true when $\x=\mathbf v_i$, while all $\x\in \mathcal X_1\cup \ldots\cup \mathcal X_s$ have $\vec u_\x$ anchored).
Thus, since $\vec u=(u_1,\ldots,u_t)$ is fully unbounded, it follows that  $\mathcal R'$ is purely virtual as noted before Proposition \ref{prop-VReay-modularCompletion}.
Hence, since $G_0$ is finitary up to dimension $\dim \R^\cup \la \mathcal X_1\cup\ldots\cup\mathcal X_s\ra$, and since $\mathcal Y_1\cup \ldots\cup \mathcal Y_s=\darrow \mathcal C\subseteq \mathcal X_1\cup\{\mathbf v_1\}\ldots\cup\mathcal X_s\cup \{\mathbf v_s\}$, it follows that $\mathcal R'$ is anchored.
 In particular, $\vec u_\x$ is anchored for every $\x\in \mathcal Y_r=\mathcal C\subseteq \mathcal X_1\cup \{\mathbf v_1\}\cup \ldots\cup \mathcal X_s\cup \{\mathbf v_s\}$. However, since $\mathcal R$ is purely virtual, this means $\x$ cannot equal any $\mathbf v_j$, and thus $\mathcal C\subseteq \mathcal X$, contrary to our assumption that $\vec u$ was a counter-example to Item 3. Thus Item 3 is now established.

4. To prove Item 4, suppose $\vec u=(u_1,\ldots,u_t)$ is an anchored  limit of an asymptotically filtered sequence $\{x_i\}_{i=1}^\infty$  of terms $x_i\in G_0$ with $u_1,\ldots,u_t\in \R^\cup \la \mathcal X\ra=\C^\cup (\mathcal X\cup \{\mathbf v_1,\ldots,\mathbf v_s\})$. If $t=0$, then Item 4 follows taking  $\mathcal B=\emptyset$, so we may assume $t\geq 1$. Let $t''\in [0,t-1]$ be the maximal index such that $\vec u''=(u_1,\ldots,u_{t''})$ is fully unbounded or trivial. If $t''=0$, then $\{x_i\}_{i=1}^\infty$ is a bounded sequence of lattice points. Since $\{x_i\}_{i=1}^\infty$ is an asymptotically filtered sequence of lattice points, this is only possible if $\{x_i\}_{i=1}^\infty$ is eventually constant, in which case $t=1$ and $\vec u^\triangleleft$ is trivial. On the other hand, if $t''>0$, then Item 3 ensures that there is a subset $\mathcal A\subseteq \mathcal X$ such that $\mathcal A$ minimally encases $-\vec u''$ urbanely.
This is also true for $t''=0$ with $\mathcal A=\emptyset$.
Let $\pi:\R^d\rightarrow \R^\cup \la \mathcal A\ra^\bot$ be the orthogonal projection and let $t'\in [t''+1,t]$ be the minimal index such that $\pi(u_{t'})\neq 0$.
Note $t'$ exists else $\mathcal A$ minimally encases $-\vec u$ as well (by Proposition \ref{prop-orReay-minecase-char}), in which case  Item 4 follows taking $\mathcal B=\mathcal A$. Let $\mathcal C\subseteq \mathcal X_1\cup \{\mathbf v_1\}\cup\ldots\cup \mathcal X_s\cup \{\mathbf v_s\}$ be the lift of $\supp_{\pi(\mathcal R)}(-\pi(u_{t'}))$, which is a support set by Proposition \ref{prop-orReay-modulo}.5.
Then Propositions \ref{prop-orReay-minecase-char}.2 and \ref{prop-orReay-minecase-char}.3(a) imply that $\mathcal C$ minimally encases $-\vec u':=-(u_1,\ldots,u_{t'})$ urbanely (since $\mathcal A\subseteq \mathcal X$).
Apply Proposition \ref{prop-VReay-modularCompletion} to $\mathcal C$ and let $\mathcal R'=(\mathcal Y_1\cup\{\mathbf w_1\},\ldots,\mathcal Y_r\cup \{\mathbf w_r\})$ be the resulting virtual Reay system in $G_0$ for the subspace $\R^\cup \la \mathcal C\ra\subseteq \R^\cup \la \mathcal X\ra$ with $\mathcal Y_r=\mathcal C$, $\vec u_{\mathbf w_r}=\vec u'$ and $\mathcal Y_1\cup\ldots\cup \mathcal Y_r=\darrow \mathcal C$.
Moreover, in the notation of Proposition \ref{prop-VReay-modularCompletion}, $\mathcal C_{s_{\ell-1}}^{(\ell-1)}=\partial(\{\mathbf w_r\})=\mathcal A\subseteq \darrow \mathcal C$.
Note $\vec u_{\mathbf w_r}=\vec u'=(u_1,\ldots,u_{t'})$ is anchored in view of the maximality of $t''$.
Since $\mathcal Y_1\cup\ldots\cup \mathcal Y_{r-1}\subseteq \mathcal X_1\cup \ldots\cup \mathcal X_s$ (by Proposition \ref{prop-VReay-modularCompletion}(a)), and since $\mathcal R$ is anchored, it follows  that we can apply Proposition \ref{prop-VReay-Lattice} to conclude that $\mathbf w_{r}(i)-\tilde{\mathbf w}_r^{(t'-1)}(i)$ is constant (and nonzero) for all  sufficiently large $i$.
Thus $\mathbf w_r(i)=x_i=\tilde{ \mathbf  w}^{(t'-1)}_r(i)+a^{(t')}_iu_{t'}+y_i$ with $\tilde{\mathbf w}_r^{(t'-1)}(i)\in \R^\cup\la \partial(\{\mathbf w_r\})\ra=\R^\cup \la \mathcal A\ra\subseteq \R^\cup \la\mathcal C\ra$, $\|y_i\|\in o(a_i^{(t')})$ and $\mathbf w_r(i)-\tilde{\mathbf w}_r^{(t'-1)}(i)=a^{(t')}_iu_{t'}+y_i=\xi\neq 0$ constant, for all sufficiently large $i$. Since $(u_1,\ldots,u_{t'})$ is anchored, we have $a_i^{(t')}\in O(1)$, so that $\|y_i\|\in o(a_i^{(t)})\subseteq o(1)$. But now $y_i\rightarrow 0$ and $\xi=\lim_{i\rightarrow \infty}(a_i^{(t')}u_{t'}+y_i)=(\lim_{i\rightarrow \infty}a^{(t')}_i)u_{t'}\in \R^\cup \la \mathcal C\ra$, with the inclusion since $\mathcal C$ minimally encases $-(u_1,\ldots,u_{t'})$.
 Thus  $x_i=\mathbf w_r(i)=\mathbf w_r^{(t'-1)}(i)+\xi\in \R\la \mathcal C\ra$, ensuring that the limit $\vec u=(u_1,\ldots,u_t)$ of the asymptotically filtered sequence $\{x_i\}_{i=1}^\infty$ has  $u_i\in \R^\cup \la \mathcal C\ra$ for all $i\in [1,t]$. Proposition \ref{prop-orReay-minecase-char} now ensures that $\mathcal C$ not only minimally encases $-(u_1,\ldots,u_{t'})$, but also $-\vec u=(u_1,\ldots,u_t)$.
  Consequently, if $\wt(\mathcal C)\leq 1$, then Item 4 follows taking $\mathcal B=\mathcal C$. Assume by contradiction that $\wt(\mathcal C)\geq 2$. Then there are distinct $\mathbf v_{k_1},\,\mathbf v_{k_2}\in \mathcal C=\mathcal Y_r$, and since $\mathcal R$ is purely virtual, both $\vec u_{\mathbf v_{k_1}}$ and $\vec u_{\mathbf v_{k_2}}$ are fully unbounded.

As noted after the definition of finitary, $\mathcal R''=(\mathcal Y_1\cup\{\mathbf w_1\},\ldots,\big((\mathcal Y_r\cup \{\mathbf w_r\})\setminus \{\mathbf v_{k_1}\}\big)\cup \{\mathbf v_{k_1}\})$ is also a virtual Reay system in $G_0$. As noted several times already in the proof, the strict truncation of any limit defining a half-space $\x$ from $\mathcal R$ is either trivial or fully unbounded. By construction, $\vec u_{\mathbf v_{k_1}}$ is fully unbounded. Consider an arbitrary $\mathbf w_j$ with $j<r$. Then, in view of Proposition \ref{prop-VReay-modularCompletion}(b) (recalling that $\mathcal C_{s_{\ell-1}}^{(\ell-1)}=\mathcal A$), it follows that  either $\vec u_{\mathbf w_j}=(u_1,\ldots,u_{r_j})$ for some $r_j\leq t''$, in which case $\vec u_{\mathbf w_j}$ is  fully unbounded by definition of $t''$, or else $\vec u_{\mathbf w_j}$ is the strict truncation of some defining limit for a  half-space from $\darrow \mathcal C$, and thus also fully unbounded. Consequently, $\mathcal R''$ is purely virtual.
 Thus, since $\R^\cup \la \mathcal C\ra\subseteq \R^\cup \la \mathcal X\ra$ ensures $\dim \R^\cup \la \mathcal C\ra\leq \R^\cup \la \mathcal X\ra$, and since $G_0$ is finitary up to dimension $\dim \R^\cup \la \mathcal X\ra$, it follows that  $\mathcal R''$ is anchored. In particular, $\vec u_\x$ is anchored for every $\x\in \mathcal Y_r\setminus \{\mathbf v_{k_1}\}=\mathcal C\setminus \{\mathbf v_{k_1}\}$, contradicting our assumption that $\mathbf v_{k_2}\in \mathcal C$ is a distinct half-space from
$\mathbf v_{k_1}$ with $\vec u_{\mathbf v_{k_2}}$ fully unbounded. Thus Item 4 is now also established, completing the proof.
\end{proof}

We will now show that $G_0^\diamond\cap -\C(G_0^\diamond)$ being bounded implies that $G_0$ is finitary, giving an important example of finitary sets apart from $G_0$ finite. The converse to this statement is false, as seen by the example $G_0=\{(-1,y):\; y\in \Z_+\}\cup \{(x,-1):\;x\in \Z_+\}\subseteq \Z^2$, which has $G_0^\diamond=G_0$ and $G_0^{\mathsf{lim}}=\{(1,0),(0,1)\}$. On the other hand, the set $$G_0=\{(-1,y):\;y\geq 1, \,y\in \Z\}\cup \{(x,0):\; x\geq 1,\, x\in \Z\}\cup \{(0,-1)\}$$ gives a prototypical example of a basic finitary set in $\Z^2$ with $G_0^\diamond\cap -\C(G_0^\diamond)=\emptyset$, $G_0^{\mathsf{lim}}=\{(1,0),(0,1)\}$ and $G_0^\diamond=\{(-1,y):\;y\geq 1, \,y\in \Z\}\cup\{(0,-1)\}$.

\begin{theorem}\label{thm-keylemmaII}
 Let $\Lambda\subseteq \R^d$ be a full rank lattice, where $d\geq 0$, and let
$G_0\subseteq \Lambda$ be a subset with $\C(G_0)=\R^d$. If $G_0^{\diamond}\cap -\C(G_0^{\diamond})$ is a bounded set, then $G_0$ is finitary. In particular, if  $0\notin \C^*(G_0^\diamond)$, then $G_0$ is finitary.
\end{theorem}

\begin{proof}
 Assume by contradiction that the theorem is false and let $\mathcal R=(\mathcal X_1\cup \{\mathbf v_1\}, \ldots, \mathcal X_s\cup\{\mathbf v_s\})$ be a purely virtual Reay system over $G_0$ which is not anchored having $\dim\R^\cup \la \mathcal X\ra$ minimal, where $\mathcal X=\mathcal X_1\cup\ldots\cup \mathcal X_s$. Then $s\geq 1$, and in view of the minimality of $\dim \R^\cup \la \mathcal X\ra$, we conclude that $G_0$ is finitary up to dimension $\dim \R^\cup \la \mathcal X\ra-1\geq 0$.
 Let $\mathcal R'_0=(\mathcal X_1\cup \{\mathbf v_1\},\ldots,\mathcal X_{s-1}\cup \{\mathbf v_{s-1}\})$. Then we can apply Proposition \ref{prop-finitary-basics} to $\mathcal R'_0$. In particular, $\mathcal R'_0$ is anchored and  $(X_1\cup\ldots\cup X_{s-1})(k)\subseteq G_0^\diamond$ for any tuple $k=(i_\z)$ with all $i_\z$ sufficiently large.

Since $\mathcal R$ is purely virtual,   $\vec u_{\mathbf v_s}$ is fully unbounded. Thus, since $\mathcal R'_0$ is anchored, Lemma \ref{Lemma-VReay-RidRemainders}.1 implies that $\x(i)\in G_0^\diamond$ for every $\x\in \darrow \mathcal A_s\setminus \{\mathbf v_s\}$ and all sufficiently large $i$. Since $\mathcal X_s\subseteq \mathcal A_s$, we now see that $\x(i)\in G_0^\diamond$ for any $\x\in \mathcal X$ and sufficiently large $i$. Moreover, Lemma \ref{Lemma-VReay-RidRemainders}.1 also implies that, for each $\y\in \mathcal X_s$, there is an asymptotically filtered sequence $\{y_{i_\y}^\y\}_{i_\y=1}^\infty$ of terms $y_{i_\y}^\y\in \C(G_0^\diamond)\cap \R^\cup\la\mathcal A_s\ra$ with limit $\vec u_\y$.

Suppose there were some $\z\in \mathcal A_s\setminus \{\mathbf v_s\}$ such that $\vec u_\z$ is also fully unbounded. Then we could apply Lemma \ref{Lemma-VReay-RidRemainders}.1 to conclude that there is also  an asymptotically filtered sequence $\{x_i\}_{i=1}^\infty$ of terms $x_i\in \C(G_0^\diamond)\cap \R^\cup\la\mathcal A_s\ra$ with limit $\vec u_{\mathbf v_s}$. Indeed, $\mathbf v_s(i)\in G_0^\diamond$ for all sufficiently large $i$ and there is a finite set $Z\subseteq G_0^\diamond$ such that every $x_i\in \mathbf v_s(i)+ \C(Z)$ once $i$ is sufficiently large.
By Lemma \ref{Lemma-VReay-RidRemainders}.2,  the set $$Y_k:=\darrow \partial(A_s)(k)\cup \{y_{i_\y}^\y:\;\y\in \mathcal A_s\setminus\{\mathbf v_s\}\}$$ minimally encases $-\vec u_{\mathbf v_s}$  with $\R^\cup \la \mathcal A_s\ra=\R\la Y_k\ra$ for any tuple $k=(i_\z)_{\z\in \darrow \mathcal A_s\setminus \{\mathbf v_s\}}$ with all $i_\z$ sufficiently large.
Fix one such tuple $\kappa=(\iota_\z)_{\z\in  \darrow \mathcal A_s\setminus \{\mathbf v_s\}}$ such that all the above conclusions are  true and let $Y=Y_\kappa$.
Since $\darrow \partial(\mathcal A_s)\subseteq \mathcal X_1\cup\ldots\cup \mathcal X_{s-1}$, we have  $\darrow \partial(A_s)(\kappa)\subseteq G_0^\diamond$, and thus $Y\subseteq \C(G_0^\diamond)$ in view of the definition of the $y_{i_\y}^\y$.
Recall that $\R\la Y\ra=\R^\cup\la \mathcal A_s\ra$. Consequently, applying Proposition \ref{prop-min-encasement-minposbasis} to the sequence
$\{x_i\}_{i=1}^\infty$, we find that $Y\cup \{x_i\}$ is a minimal positive basis for all sufficiently large $i$. Hence, since $x_i\in \mathbf v_s(i)+ \C(Z)$, we conclude that $\mathbf v_s(i)\in -\C(Y\cup Z)\subseteq -\C(G_0^\diamond)$ for all sufficiently large $i$. However, since $\vec u_{\mathbf v_s}$ is fully unbounded (in view of $\mathcal R$ being purely virtual), we see that $\{\mathbf v_s(i)\}_{i=1}^\infty$ is an unbounded sequence of elements contained in $-\C(G_0^\diamond)$ with $\mathbf v_s(i)\in G_0^\diamond$ for all sufficiently large $i$, contradicting the hypothesis that $G_0^\diamond \cap -\C(G_0^\diamond)$ is bounded. So we instead conclude that $\vec u_\z$ is anchored for every $\z\in \mathcal A_s\setminus \{\mathbf v_s\}$. In particular, this ensures that every $\vec u_\y$ with $\y\in \mathcal X_s\subseteq \mathcal A_s\setminus \{\mathbf v_s\}$ is anchored, implying that $\mathcal R$ is anchored (as $\mathcal R'_0$ is anchored), contradicting that $\mathcal R$ was a minimal counterexample.

For the in particular statement, suppose $0\notin \C^*(G_0^\diamond)$. Then, since
$0\notin G_0^\diamond$ (by Proposition \ref{prop-G_0diamond-1st-easy-equiv}.2), it follows that $G_0^\diamond\cap -\C(G_0^\diamond)$ is bounded (indeed, empty), and so  $G_0$ is finitary by the main part of the theorem.
\end{proof}

Proposition \ref{prop-finitary-basics} has some important consequences beyond Theorem \ref{thm-keylemmaII}. Suppose $G_0\subseteq \Lambda\subseteq \R^d$ is finitary with $\C(G_0)=\R^d$. If $X
\subseteq G_0\cup G_0^\infty$ is a minimal positive basis, then each $x\in X$ defines a one-dimensional half-space $\x=\R_+ x$, and letting $\mathcal X$ be the corresponding set of half-spaces, it follows that $((\mathcal X\setminus \{\x\})\cup \{\x\})$ is a virtual Reay system in $G_0$ for any $\x\in \mathcal X$. For a half-space $\z\in \mathcal X$, we either have $z\in G_0^\infty$, in which case $\vec u_\z\in G_0^\infty$ is unbounded, or else $z\in G_0$, in which case $\vec u_\z=z/\|z\|$ is bounded (and in case $z\in G_0$ and $z/\|z\|\in  G_0^\infty$, then both are simultaneously possible). Consequently, if $X\cap G_0^\infty\neq \emptyset$, then we can choose $\x\in \mathcal X$ so that the resulting virtual Reay system $((\mathcal X\setminus \{\x\})\cup \{\x\})$ is purely virtual. Since $G_0$ is finitary, the virtual Reay system is anchored,  i.e., $|X\cap G_0^\infty|\leq 1$. Moreover, if $|X\cap G_0^\infty|=1$, then $X\setminus G_0^\infty\subseteq G_0^\diamond$ by Proposition \ref{prop-finitary-basics}.1. Thus there are limitations on how minimal positive bases can be formed using elements from $G_0\cup G_0^\infty$.

In general, if $\{z_i\}_{i=1}^\infty$ is an asymptotically filtered  sequence of lattice points with limit $u_1$, say with corresponding representation $z_i=a_iu_1+y_i$, then we either have $a_i\rightarrow \infty$ or $a_i\rightarrow C$ for some  $C\geq 0$. If $a_i\rightarrow C$, then $\|y_i\|\in o(a_i)\subseteq o(1)$, ensuring that  $\|y_i\|\rightarrow 0$, which means that $\{z_i\}_{i=1}^\infty$ is a convergent sequence of lattice points with limit $(\lim_{i\rightarrow \infty} a_i)u_1=Cu_1$. Thus, since any convergent sequence of lattice points is eventually constant, it follows that $z_i=(\lim_{i\rightarrow \infty}a_i)u_1=Cu_1$ for all sufficiently large $i$, forcing $C>0$ as $z_i=Cu_1=a_iu_1+y_i\neq 0$ for all sufficiently large $i$ (as $a_i>0$ with $y_i\in \R\la u_1\ra^\bot$ for all $i$). Moreover, since $z_i=Cu_1$ for all sufficiently large $i$, and since $u_1$ and $y_i$ are orthogonal (and thus linearly independent), it follows that $y_i=0$ for all sufficiently large $i$, and if $\{b_i\}_{i=1}^\infty$ is a sequence of positive real numbers with $a_i\in o(b_i)$, then $a_i\rightarrow C>0$ ensures $\{b_i\}_{i=1}^\infty$ is unbounded.
Now suppose  $\{x_i\}_{i=1}^\infty$ is an asymptotically filtered sequence of terms  $x_i\in G_0\subseteq \Lambda$ having limit $\vec u=(u_1,\ldots,u_t)$, say with $x_i=a_i^{(1)}u_1+\ldots+a_i^{(t)}u_t+y_i$, and that $\pi:\R^d\rightarrow \mathcal E^\bot$ is the orthogonal projection for some subspace $\mathcal E=\R\la Z\ra$ generated by a subset of lattice points $Z\subseteq \Lambda$.
Then $\pi(G_0)\subseteq \pi(\Lambda)$ with $\pi(\Lambda)$ a full rank lattice in $\mathcal E^\bot$ by Proposition \ref{Prop-lattice-homoIm}. If $\pi(\vec u)$ is the empty tuple, then $u_i\in \mathcal E$ for all $i$.
Otherwise, there is a minimal index $r_1\in [1,t]$ with $\pi(u_{r_1})\neq 0$, we may consider  $\{x_i\}_{i=1}^\infty$ as an asymptotically filtered sequence with truncated limit $(u_1,\ldots,u_{r_1})$, and Proposition \ref{prop-infinite-limits-proj}.1  implies that the sufficiently large index terms in $\{\pi(x_i)\}_{i=1}^\infty$ form a radially convergent sequence with limit $\overline u_{r_1}:=\pi(u_{r_1})/|\pi(u_{r_1})\|$, say with representation $\pi(x_i)=b_i^{(r_1)}\overline u_{r_1}+\overline y_i$ where $b_i^{(r_1)}\in \Theta(a_i^{(r_1)})$ and $\overline y_i$ is a vector orthogonal to $\overline u_{r_1}$ with $\|\overline y_i\|\in o(b_i^{(r_1)})=o(a_i^{(r_1)})$.
Consequently, if $\{a_i^{(r_1)}\}_{i=1}^\infty$ is bounded, then we can apply the previous observations to conclude that  $\pi(x_i)$ is constant with $\overline y_i=0$ for all sufficiently large $i$, and that $\{a_i^{(j)}\}_{i=1}^\infty$ is unbounded for any $j<r_1$. In such case, $\pi(\vec u)=\overline u_{r_1}$ is a complete limit for $\{\pi(x_i)\}_{i=1}^\infty$ with $\pi(x_i)=C\overline u_{r_1}$ for all sufficiently large $i$, for some $C>0$, and  $a_i^{(j)}\rightarrow \infty$ for any $j<r_1$. This does not imply $r_1=t$, but does ensure $u_i\in \mathcal E+\R u_{r_1}$ for all $i\in [1,t]$. We will need these observations in our discussion below.

Now suppose that
$$\mathcal R=(\mathcal X_1\cup \{\mathbf v_1\},\ldots,\mathcal X_s\cup \{\mathbf v_s\})$$ is a purely virtual Reay system (possibly trivial) in the finitary set $G_0$ and let $$\pi:\R^d\rightarrow \R^\cup\la \mathcal X_1\cup\ldots\cup \mathcal X_s\ra^\bot$$ be the orthogonal projection. Since $G_0$ is finitary, it follows that $\mathcal R$ is anchored, ensuring that $\R^\cup \la \mathcal X_1\cup \ldots\cup \mathcal X_s\ra=\R\la \darrow (X_1\cup \ldots  X_s)(k)\ra$ for any tuple $k=(i_\z)_{\z\in \darrow (\mathcal X_1\cup \ldots\cup \mathcal X_s)}$ with all $i_\z$ sufficiently large (in view of Proposition \ref{prop-orReay-BasicProps}.1 and Proposition \ref{prop-VReay-Lattice}).
Consequently, since $\darrow (X_1\cup \ldots  X_s)(k)\subseteq G_0\subseteq\Lambda$ is a subset of lattice points, it follows from Proposition \ref{Prop-lattice-homoIm} that $\pi(\Lambda)$ is a full rank lattice in $\pi(\R^d)=\R^\cup\la \mathcal X_1\cup\ldots\cup \mathcal X_s\ra^\bot$.
Let $v\in \pi(G_0)\cup \pi(G_0)^\infty$ be nonzero. In view Proposition \ref{prop-infinite-limits-proj}.2, there is an asymptotically filtered sequence $\{x_i\}_{i=1}^\infty$ of terms $x_i\in G_0$ with limit $\vec u=(u_1,\ldots,u_t)$ such that $\{\pi(x_i)\}_{i=1}^\infty$ is an asymptotically filtered sequence of terms with limit $\pi(\vec u)=\pi(u_t)/\|\pi(u_t)\|=v/\|v\|$. For instance, if $v\in \pi(G_0)$, then $\{x_i\}_{i=1}^\infty$ may be taken to be a constant sequence, though there may be other non-constant sequences with $\pi(x_i)=v$ for all $i$ as well which we could choose instead.
  Moreover, we either have $v\in \pi(G_0)^\infty$, in which case $\vec u$ is fully unbounded by Proposition \ref{prop-infinite-limits-proj}.1 (with $r_\ell=r_1=t$ in the notation of Proposition \ref{prop-infinite-limits-proj}.1), or else $v\in \pi(G_0)$, in which case $\vec u$ is anchored (and in case $v\in \pi(G_0)$ and $v/\|v\|\in \pi(G_0)^\infty$, then both are possible).  Let $x_i=a_i^{(1)}u_1+\ldots+a_i^{(t-1)}u_{t-1}+a^{(t)}_iu_t+y_i$ be the representation of $\{x_i\}_{i=1}^\infty$ as an asymptotically filtered sequence with limit $\vec u$. If $a_i^{(t)}\rightarrow \infty$, then $a_i^{(j)}\rightarrow \infty$ for all $j<t$ holds trivially. On the other hand, if $a_i^{(t)}$ is bounded, then the earlier discussion above also ensures that $a_i^{(j)}\rightarrow \infty$ for all $j<t$, so $\vec u^\triangleleft$ is trivial or fully unbounded in all cases.
But now Proposition \ref{prop-finitary-basics}.3 ensures that there is some $\mathcal B_\x\subseteq \mathcal X_1\cup\ldots\cup \mathcal X_s$ which minimally encases $-\vec u^\triangleleft$.
This allows us to define a half-space $\x$ by setting  $\overline \x=\R^\cup\la\mathcal B_\x\ra+\R_+ u_t$, $\partial(\x)=\R^\cup\la \mathcal B_\x\ra$ and  $\partial(\x)\cap \x=\C(\mathcal B_\x)$.
 Now suppose $X\subseteq \pi(G_0)\cup \pi(G_0)^\infty$ is a minimal positive basis and let $\mathcal X$ be the set of half-spaces obtained from the elements $x\in X$ as just described (using some compatible choice of asymptotically  filtered sequences for each $x\in X$).
 Then $\mathcal R'=(\mathcal X_1\cup \{\mathbf v_1\},\ldots,\mathcal X_s\cup \{\mathbf v_s\},(\mathcal X\setminus \{\x\})\cup \{\x\})$ will be a virtual Reay system over $G_0$ for any $\x\in \mathcal X$.
Since $\C(G_0)=\R^d$, we have $\C(\pi(G_0))=\pi(\R^d)=\R^\cup\la \mathcal X_1\cup\ldots\cup \mathcal X_s\ra^\bot$.
Consequently, given any nonzero $x\in \pi(G_0)\cup \pi(G_0)^\infty$, there exists a minimal positive basis $X\subseteq \pi(G_0)\cup \pi(G_0)^\infty$ which contains $x$. If $X\cap \pi(G_0)^\infty$ is nonempty, then we can choose $\x\in \mathcal X$ so that $\vec u_\x$ is fully unbounded, in which case $\mathcal R'$ will be purely virtual. Since $G_0$ is finitary, we conclude that $\mathcal R'$ must be anchored, i.e., $|X\cap \pi(G_0)^\infty|\leq 1$.
Moreover, if $|X\cap \pi(G_0)^\infty|=1$, then $X\setminus \pi(G_0)^\infty\subseteq \pi(G_0)^\diamond$ by Proposition \ref{prop-G_0diamond-1st-easy-equiv}.1. Thus we have generalized our initial observation regarding minimal positive bases over $G_0\cup G_0^\infty$.
Furthermore, we now see that purely virtual Reay systems over a finitary set $G_0$ can be constructed greedily by recursively applying the construction just described with our only being prevented from extending the purely virtual Reay system $\mathcal R$ to a larger one when $\pi(G_0)^\infty=\emptyset$, i.e., when $\pi(G_0)$ is finite (since $\pi(\Lambda)$ is a lattice).

Likewise, if we have a purely virtual Reay system $\mathcal R'=(\mathcal Y_{s+1}\cup\{\mathbf w_{s+1}\},\ldots, \mathcal Y_{s'}\cup \{\mathbf w_{s'}\})$ over $\pi(G_0)$,
then $Y_{s+1}\cup \{w_{s+1}\}$ will be a minimal positive basis, and we can define a virtual set $\mathcal X_{s+1}\cup \{\mathbf v_{s+1}\}$ (using any compatible choice of asymptotically filtered sequence for each $\z\in \mathcal Y_{s+1}\cup \{\mathbf w_{s+1}\}$) as above so that $\mathcal R''=(\mathcal X_1\cup \{\mathbf v_1\},\ldots,\mathcal X_{s+1}\cup \{\mathbf v_{s+1}\})$ is a purely virtual Reay system over $G_0$ with $\pi(\mathcal R'')=(\mathcal X_{s+1}^\pi\cup \{\pi(\mathbf v_{s+1})\})=(\mathcal Y_{s+1}\cup \{\mathbf w_{s+1}\})$.
Then it follows that $\pi'(\mathcal R')=(\pi'(\mathcal Y_{s+2})\cup \{\pi'(\mathbf w_{s+2})\},\ldots, \pi'(\mathcal Y_{s'})\cup \{\pi'(\mathbf w_{s'})\})$ is a purely virtual Reay system over $\pi'(G_0)$ (by Proposition \ref{prop-VReay-modulo}), where $\pi':\R^d\rightarrow \R^\cup \la \mathcal X_1\cup\ldots\cup \mathcal X_{s+1}\ra^\bot$ is the orthogonal projection.
As before, $\pi'(\Lambda)$ a full rank lattice and, for each  $\y\in \mathcal Y_{s+2}\cup \{\mathbf w_{s+2}\}$, we can choose via Proposition \ref{prop-infinite-limits-proj}.2 an asymptotically filtered sequence $\{x_i\}_{i=1}^\infty$ of terms $x_i\in G_0$ with limit $\vec u_\x=(u_1,\ldots,u_t)$ such that $\pi(\vec u_\x)=\vec u_{\y}$ and $\pi(\vec u_\x^\triangleleft)=\pi(\vec u_\x)^\triangleleft=\vec u_{\y}^\triangleleft$, and then
$\pi'(\vec u_\x)=\pi'(\pi(\vec u_\x))=\pi'(\vec u_\y)=\vec u_{\pi'(\y)}$ follows in view of $\ker \pi\leq \ker \pi'$ and Proposition \ref{prop-VReay-modulo}.
Since $\pi'(\y)$ lies at depth $1$ in the oriented Reay system $\pi'(\mathcal R')$, it has trivial boundary, ensuring that $\vec u_{\pi'(\y)}$ is a single unit vector.
Moreover, since $\pi'(\y)$ is a relative half-space in $\pi'(\mathcal R')$, we must have $\pi'(\pi(u_t))=\pi'(u_t)$ being a representative for $\pi'(\y)$ since $\pi(u_t)$ is a representative for $\y$ (in view of $\pi(\vec u_\x)=\vec u_\y$ and  $\pi(\vec u_\x^\triangleleft)=\pi(\vec u_\x)^\triangleleft$), ensuring that $\vec u_{\pi'(\y)}=\pi'(u_t)/\|\pi'(u_t)\|$. (Alternatively, since $\partial(\{\y\})\subseteq \mathcal Y_{s+1}=\pi(\mathcal X_{s+1})$,  $\pi(\vec u_\x)=\vec u_\y$ and  $\pi(\vec u_\x^\triangleleft)=\pi(\vec u_\x)^\triangleleft$, we have $\pi'(u_i)=0$ for all $i<t$, ensuring the same conclusion.)
Thus,
per prior discussion,  $\vec u_\x^\triangleleft$ will be fully unbounded (or trivial),  and thus $-\vec u_\x^\triangleleft$ is minimally encased by some $\mathcal B_\x\subseteq \mathcal X_1\cup\ldots\cup \mathcal X_{s+1}$, allowing us to define a half-space $\x$ with $\partial(\{\x\})=\mathcal B_\x$ and $\overline \x=\R^\cup\la \mathcal B_\x\ra+\R_+u_t$.
Since the support set $\partial(\{\x\})=\mathcal B_\x$ minimally encases $-\vec u_\x^\triangleleft$, Proposition \ref{prop-orReay-minecase-char}.4 and Proposition \ref{prop-orReay-modulo}.6 ensure that the support set  $\partial(\{\x\})^\pi\subseteq \pi(\mathcal X_{s+1})=\mathcal Y_{s+1}$ minimally encases $-\pi(\vec u_\x^\triangleleft )=-\pi(\vec u_\x)^\triangleleft=-\vec u_\y^\triangleleft$, which means $\partial(\{\x\})^\pi=\partial(\{\y\})$.
Hence, since the representative $u_t$ for $\x$ maps to the representative $\pi(u_t)$ for $\y$ (in view of $\pi(\vec u_\x)=\vec u_{\y}$ and $\pi(\vec u_\x^\triangleleft)=\pi(\vec u_\x)^\triangleleft$), we conclude that $\pi(\x)=\y$.
We can iterate this procedure, resulting  in a purely virtual Reay system $\mathcal W=(\mathcal X_1\cup \{\mathbf v_1\},\ldots,\mathcal X_{s'}\cup\{\mathbf v_{s'}\})$ over $G_0$ such that $\pi(\mathcal W)=\mathcal R'$ and $\mathcal X_j^\pi\cup \{\pi(\mathbf v_j)\}=\mathcal Y_j\cup \{\mathbf w_j\}$ for all $j\in [s+1,s']$, and
  $\pi(\vec u_\x)=\vec u_{\pi(\x)}$ and $\pi(\vec u_\x^\triangleleft)=\pi(\vec u_\x)^\triangleleft$ for all $\x\in \mathcal X_{s+1}\cup \{\mathbf v_{s+1}\}\cup \ldots\cup \mathcal X_{s'}\cup \{\mathbf v_{s'}\}$,
  with each $\vec u_\x$  the limit of a compatible asymptotically filtered sequence initially chosen for the half-space $\y=\pi(\x)\in \mathcal Y_{s+1}\cup \{\mathbf w_{s+1}\}\cup \ldots\cup \mathcal Y_{s'}\cup \{\mathbf w_{s'}\}$.
  We call $\mathcal W$ an \textbf{extension} of $\mathcal R$.

We call the purely virtual Reay system $\mathcal R$ over $G_0$ with $\pi(G_0)^\infty=\emptyset$ a \textbf{maximal} purely virtual Reay system over $G_0$. Since $\pi(G_0)$ is a subset of lattice points, this is equivalent to $\pi(G_0)$ being finite as noted earlier. Now further assume that $\mathcal R$ is a maximal purely virtual Reay system over $G_0$, so $\pi(G_0)$ is a finite set of lattice points. Suppose $\mathcal R'=(\mathcal Y_{s+1}\cup\{\mathbf w_{s+1}\},\ldots, \mathcal Y_{s'}\cup \{\mathbf w_{s'}\})$ is a virtual Reay system over $\pi(G_0)$. For instance, we could  use  Proposition \ref{prop-reay-basis-exists} to find an ordinary Reay system $(Y_{s+1}\cup \{w_{s+1}\},\ldots,Y_{s'}\cup \{w_{s'}\})$ (for the entire space $\pi(\R^d)$ if we like), and then replace each element $y\in Y_{s+1}\cup \{w_{s+1}\}\cup \ldots\cup Y_{s'}\cup \{w_{s'}\}$ with the one-dimensional half-space $\y=\R_+ y$ it represents to obtain  $\mathcal R'$.
We can associate to each $\y\in Y_{s+1}\cup \{w_{s+1}\}\cup \ldots\cup Y_{s'}\cup \{w_{s'}\}$ an asymptotically  filtered sequence $\{x_i\}_{i=1}^\infty$ of terms from $G_0$ with limit $\vec u=(u_1,\ldots,u_t)$ such that $\pi(\vec u)=\pi(u_t)/\|\pi(u_t)\|=y/\|y\|$, which ensures $u_i\in \ker \pi=\R^\cup \la \mathcal X_1\cup \ldots\cup \mathcal X_s\ra$ for all $i<t$.
For instance, we could take $\{x_i\}_{i=1}^\infty$ to be a constant sequence, though there may be other possibilities (we simply need $\pi(x_i)$ to eventually be constant).
Per prior discussion, $\vec u^\triangleleft$ must be trivial or fully unbounded.
Thus, since  $u_1,\ldots,u_{t-1}\in\R^\cup\la \mathcal X_1\cup \ldots\cup \mathcal X_s\ra$, Proposition \ref{prop-finitary-basics}.3 ensures that there is some $\mathcal B_\x\subseteq \mathcal X_1\cup\ldots\cup \mathcal X_s$ which minimally encases $-\vec u^\triangleleft$. As before, this means we can  define a half-space $\x$ by setting  $\overline \x=\R^\cup\la\mathcal B_\x\ra+\R_+ u_t$, $\partial(\x)=\R^\cup\la \mathcal B_\x\ra$ and  $\partial(\x)\cap \x=\C(\mathcal B_\x)$, and  replacing each $\y\in \mathcal Y_{s+1}\cup \{\mathbf w_{s+1}\}\cup \ldots\cup \mathcal Y_{s'}\cup \{\mathbf w_{s'}\}$ with the corresponding half-space $\x$ just constructed again results in virtual sets $\mathcal X_{s+1}\cup \{\mathbf v_{s+1}\},\ldots, \mathcal X_{s'}\cup\{\mathbf v_{s'}\}$ such that $\mathcal W=(\mathcal X_1\cup \{\mathbf v_1\},\ldots,\mathcal X_{s'}\cup\{\mathbf v_{s'}\})$ is a virtual Reay system over $G_0$ with $\pi(\mathcal R)=\mathcal R'$, \
 $\mathcal X_j^\pi\cup \{\pi(\mathbf v_j)\}=\mathcal Y_j\cup \{\mathbf w_j\}$ for all $j\in [s+1,s']$, and $\pi(\vec u_\x)=\vec u_{\pi(\x)}$ for all $\x\in \mathcal X_{s+1}\cup \{\mathbf v_{s+1}\}\cup \ldots\cup \mathcal X_{s'}\cup \{\mathbf v_{s'}\}$, with each $\vec u_\x$  the limit of the compatible asymptotically filtered sequence initially chosen for the half-space $\y=\pi(\x)\in \mathcal Y_{s+1}\cup \{\mathbf w_{s+1}\}\cup \ldots\cup \mathcal Y_{s'}\cup \{\mathbf w_{s'}\}$. We also call $\mathcal W$ an extension of $\mathcal R$.


We continue by showing that finitary sets remain finitary modulo the subspace generated by a purely virtual Reay system.

\begin{proposition}\label{prop-finitary-Modulo-Inheritence}
Let $\Lambda\subseteq \R^d$ be a full rank lattice, where $d\geq 0$, let
$G_0\subseteq \Lambda$ be a finitary subset with $\C(G_0)=\R^d$, let $(\mathcal X_1\cup\{\mathbf v_1\},\ldots, \mathcal X_s\cup \{\mathbf v_s\})$ be a purely virtual Reay system over $G_0$, and let $\pi:\R^d\rightarrow \R^\cup \la \mathcal X_1\cup \ldots\cup \mathcal X_s\ra^\bot$ be the orthogonal projection.  Then $\pi(\Lambda)$ is a full rank lattice in $\C(\pi(G_0))=\R^\cup \la \mathcal X_1\cup \ldots\cup \mathcal X_s\ra^\bot$ and $\pi(G_0)$ is finitary.
\end{proposition}

\begin{proof}
As already remarked after Theorem \ref{thm-keylemmaII}, $\pi(\Lambda)$ is a full rank lattice.
Since $\C(G_0)=\R^d$, we have $\C(\pi(G_0))=\R^\cup \la \mathcal X_1\cup \ldots\cup \mathcal X_r\ra^\bot$.
Let $\mathcal R'=(\mathcal Y_{s+1}\cup \{\mathbf w_{s+1}\},\ldots,\mathcal Y_r\cup \{\mathbf w_r\})$ be an arbitrary purely virtual Reay system over $\pi(G_0)$.  Let $\mathcal W=(\mathcal X_1\cup \{\mathbf v_1\},\ldots,\mathcal X_r\cup \{\mathbf v_r\})$ be the extension of $\mathcal R$ using $\mathcal R'$. Then $\mathcal W$ will be purely virtual (as both $\mathcal R$ and $\mathcal R'$ are purely virtual). Thus $\mathcal W$ is anchored since $G_0$ is finitary, implying that $\pi(\mathcal W)=\mathcal R'$ is also anchored (as $\pi(\vec u_\x^\triangleleft)=\pi(\vec u_\x)^\triangleleft$ for all $\x\in \mathcal X_1\cup\ldots\cup \mathcal X_s$, so if $\pi(\vec u_\x)=u_{\pi(\x)}$ were fully unbounded, then $\vec u_\x$ could be fully unbounded too). Since $\mathcal R'$ is an arbitrary purely virtual Reay system over $\pi(G_0)$, we conclude that $\pi(G_0)$ is finitary, as desired.
\end{proof}

\begin{proposition}\label{prop-finitary-diamond-modulo-transfer} Let $\Lambda\subseteq \R^d$ be a full rank lattice, where $d\geq 0$, let
$G_0\subseteq \Lambda$ be a finitary subset with $\C(G_0)=\R^d$, let $\mathcal R=(\mathcal X_1\cup \{\mathbf v_1\},\ldots,\mathcal X_s\cup\{\mathbf v_s\})$ be a purely virtual Reay system over $G_0$,  and let  $\pi:\R^d\rightarrow \R^\cup\la \mathcal X_1\cup\ldots\cup
\mathcal X_r\ra^\bot$ be the orthogonal projection. Then, for any $\pi(g)\in \pi(G_0)^{\diamond}$ with $g\in G_0$, we have  $g\in G_0^\diamond$. In particular, $$\pi(G_0)^\diamond\subseteq \pi(G_0^\diamond).$$
\end{proposition}

\begin{proof}
Let $\mathcal X=\mathcal X_1\cup\ldots\cup \mathcal X_s$. As $G_0$ is finitary,  $\mathcal R$ is anchored with $\pi(\Lambda)$ a full rank lattice by Proposition \ref{prop-finitary-Modulo-Inheritence}.
Let $g\in G_0$ be an arbitrary element with $\pi(g)\in \pi(G_0)^\diamond$. Then, in view of Proposition \ref{prop-G_0diamond-1st-easy-equiv}, there exists a subset $Y\subseteq G_0$ such that  $|Y|=|\pi(Y)|$, $g\in Y$,  and $\pi(Y)$ minimally encases $-\vec v$, where $\vec v$ is a fully unbounded limit of  an asymptotically filtered sequence of terms $\{x'_i\}_{i=1}^\infty$ from $\pi(G_0)$. By Proposition \ref{prop-infinite-limits-proj}.2, we can w.l.o.g. assume there is an asymptotically filtered sequence $\{x_i\}_{i=1}^\infty$ of terms $x_i\in G_0$ with limit $\vec u=(u_1,\ldots,u_t)$ such that $\pi(x_i)=x'_i$ for all $i$ and $\pi(\vec u)=\vec v$, and as $\vec v$ is fully unbounded, we can assume by Proposition \ref{prop-infinite-limits-proj}.1 that $\vec u$ is fully unbounded as well (by choosing such $\vec u$ with $t$ minimal, i.e., with $\pi(\vec u^\triangleleft)=\pi(\vec u)^\triangleleft=\vec v^\triangleleft$). Thus $\pi(Y)$ minimally encases
$-\pi(\vec u)=-\vec v$, and  thus also the equivalent tuple $-(\pi(u_{s_1}),\ldots,\pi(u_{s_{r}}))$, where $1\leq s_1<\ldots<s_r=t$  are the associated indices for $\pi(\vec u)$. Hence
 Proposition \ref{prop-min-encasement-char} implies there is a partition $Y=Y_1\cup\ldots\cup Y_{\ell}$ such that $(\pi(Y_1)\cup \{\pi(u_{r_1})\},\ldots,\pi(Y_s)\cup\{\pi(u_{r_\ell})\})$ is a Reay system for some $r_1<\ldots<r_\ell$ with $\{r_1,\ldots,r_\ell\}\subseteq \{s_1,\ldots,s_{r}\}$ and $r_1=s_1$. By the same argument  described after Theorem \ref{thm-keylemmaII}, we can extend the virtual Reay system $\mathcal R$ to obtain a virtual Reay system
$\mathcal W=(\mathcal X_1\cup \{\mathbf v_1\},\ldots,\mathcal X_s\cup \{\mathbf v_s\},\mathcal Y_1\cup \{\mathbf w_1\},\ldots,\mathcal Y_\ell\cup \{\mathbf w_\ell\})$ over $G_0$ with $\vec u_{\mathbf w_j}=(u_1,\ldots,u_{r_j})$  and $Y_j$ a set of representatives for the one-dimensional half-spaces in $\mathcal Y_j$,  for each $j\in [1,\ell]$, where we use the definition of the $r_j$ as given by Proposition \ref{prop-min-encasement-char}.1, that $u_i\in \R^\cup \la \mathcal X_1\cup\ldots\cup \mathcal X_s\cup \mathcal Y_1\cup\ldots\cup \mathcal Y_{j-1}\ra$ for all $i<r_j$ with $u_{r_j}\notin \R^\cup \la \mathcal X_1\cup\ldots\cup \mathcal X_s\cup \mathcal Y_1\cup\ldots\cup \mathcal Y_{j-1}\ra$, to ensure  each half-space $\mathbf w_j$ is well-defined with $\partial(\{\mathbf w_j\})\subseteq \mathcal X_1\cup\ldots\cup\mathcal X_s\cup \mathcal Y_1\cup\ldots\cup \mathcal Y_{j-1}$ the support set minimally encasing $-\vec u^\triangleleft$ by Proposition \ref{prop-finitary-basics}.3.
Note $\mathcal W$ is purely virtual as $\mathcal R$ is purely virtual and $\vec u$ is fully unbounded. Since $g\in Y=Y_1\cup \ldots\cup Y_\ell$, it follows that $g$ is a representative for some $\y\in \mathcal Y_j$ and $j\in [1,\ell]$. Indeed, we may take $g=\y(i)=\tilde \y(i)$ for all $i\geq 1$. But now Proposition \ref{prop-finitary-basics}.1 implies that $g\in G_0^\diamond$, as desired.
\end{proof}

Our next goal is to give a characterization of $G_0^\diamond$, for $G_0$ finitary,  in terms of $\mathcal A^{\mathsf{elm}}(G_0)$ and linear combinations over $\Q_+$. Towards that end, we continue with the following lemmas.

\begin{lemma}\label{lem-lattice-smithnormal}
Let $\Lambda\subseteq \R^d$ be a full rank lattice, where $d\geq 0$, and let $X\subseteq \Lambda$ be a linearly independent subset. There exists a positive integer $N>0$ such that, for any $g\in \Lambda\cap\R\la X\ra$, we have $Ng\in \Z\la X\ra$. Moreover, if we also have $g\in -\C^\circ(X)$, then $X\cup \{g\}$ is a minimal positive basis (or $X=\emptyset$ with $g=0$) and there is an elementary atom $U$ with $\supp(U)=X\cup \{g\}$ and $1\leq \vp_g(U)\leq N$.

\end{lemma}

\begin{proof}If $X=\emptyset$, then the lemma holds with $N=1$, so we may assume $X$ is nonempty.
Since $X$ is linearly independent,  $\Z\la X\ra$ is a full rank lattice in $\R\la X\ra$. Via the Smith normal form (\cite[Theorem III.7.8]{lang} \cite[Theorem 2]{Handbodk-Gruber}) applied to the sublattice $\Z\la X\ra \leq \Lambda$, we can find a lattice basis $\{e_1,\ldots,e_d\}$ for $\Lambda$ and positive integers $a_1\mid \ldots\mid a_s$, where $1\leq s\leq d$, such that $\{a_1e_1,\ldots,a_se_s\}$ is a lattice basis for $\Z\la X\ra$. Note that $\R\la e_1,\ldots,e_s\ra=\R\la a_1e_1,\ldots,a_se_s\ra=\R\la X\ra$. Thus, since $\{e_1,\ldots,e_d\}$ is a lattice basis for $\Lambda$, and in particular linearly independent, it follows that any $g\in \Lambda\cap \R\la X\ra$ lies in the lattice $\Z\la e_1,\ldots,e_s\ra$. Since $\Z\la X\ra\leq \Z\la e_1,\ldots,e_s\ra$ is a sublattice of full rank, it follows that $\Z\la e_1,\ldots,e_s\ra/\Z\la X\ra\cong G:=\Z/a_1\Z\times\ldots\times \Z/a_s\Z$ is a finite abelian group of exponent $N=a_s>0$. As a result, any $g\in \Lambda \cap \R\la X\ra$ has $Ng\in \Z\la X\ra$.

Now suppose $g\in \Lambda\cap -\C^\circ(X)\subseteq \Lambda \cap \R\la X\ra$. Then Proposition \ref{prop-char-minimal-pos-basis}.4 implies that $X\cup \{g\}$ is a minimal positive basis (since $X$ is linearly independent) and Proposition \ref{prop-char-minimal-pos-basis}.5 implies that there is a strictly positive linear combination $\Summ{x\in X\cup \{g\}}\alpha_xx=0$, so $\alpha_x>0$ for all $x\in X\cup \{g\}$, with the property that $a_g=1$ and, if $\Summ{x\in X\cup \{g\}}\beta_xx=0$ is another linear combination with $\beta_x\in \R$, then the vector $(\beta_x)_{x\in X\cup \{g\}}$ is a real scalar multiple of the vector $(\alpha_x)_{x\in X\cup \{g\}}$.

Since $Ng\in \Z\la X\ra$, it follows that there is a linear combination $\Summ{x\in X\cup \{g\}}\beta_xx=0$ with $\beta_x\in \Z$ for all $x\in X\cup \{g\}$ and $\beta_g=N>0$. Thus, by the previous observation, we must have $\beta_x=N\alpha_x$ for all $x\in X\cup \{g\}$. In particular, since $N>0$ and $a_x>0$ for all $x\in X\cup \{g\}$, it follows that $\beta_x>0$ for all $x\in X\cup \{g\}$, implying that $S=\prod^\bullet_{x\in X\cup \{g\}}x^{[\beta_x]}$ is a zero-sum sequence. As a result, since $X\cup \{g\}$ is a minimal positive basis, Proposition \ref{prop-char-minimal-pos-basis}.8 implies that  there is an elementary atom $U\mid S$ with $\supp(U)=X\cup \{g\}$ and $S=U^{[m]}$ for some integer $m\geq 1$, ensuring that  $1\leq \vp_g(U)\leq \beta_g=N$, which completes the proof.
\end{proof}

\begin{lemma}
\label{lem-finite-nozs} Let $d\geq 1$ and let $X\subseteq \R^d$ be finite. Then $0\notin \C^*(X)$ if and only if there is a co-dimension one subspace $\mathcal H\subseteq \R^d$ such that $X\subseteq \mathcal H_+^\circ$.
\end{lemma}

\begin{proof}
Any set $X\subseteq \R^d$ contained in an open half-space clearly has $0\notin \C^*(X)$. For the reverse inclusion, we proceed by induction  on $d$, with the case $d=1$ clear. By  Proposition \ref{prop-no-zs-char}, $X\subseteq \mathcal H_+$ for some co-dimension $1$ subspace $\mathcal H\subseteq \R^d$. Applying the induction hypothesis to $X\cap \mathcal H$ yields  a subspace $\mathcal H'\subseteq \mathcal H$  such that $X\cap \mathcal H$ lies entirely on one side of the subspace $\mathcal H'$. Since $X$ is finite, the  space $\mathcal H$ can be slightly perturbed by a rotation about $\mathcal H'$ (fixing the points in $\mathcal H'$) to yield the needed subspace.
\end{proof}

 Proposition \ref{prop-G_0diamond-1st-easy-equiv}.2  means $g\in G_0^\diamond$ when there is a fixed linearly independent set $X$ containing $g$ such that $x_i=-\Summ{x\in X}\alpha_i^{(x)}x\in -\C^\circ(X)$ for some $x_i\in G_0$ with $\alpha_i^{(g)}\rightarrow \infty$. By  Propositions \ref{prop-char-minimal-pos-basis}.4  and \ref{prop-char-minimal-pos-basis}.5, this means there are elementary atoms $U_i$ with $\supp(U_i)=X\cup \{x_i\}$ and $\vp_x(U_i)=\vp_{x_i}(U_i)\alpha_i^{(x)}\geq \alpha_i^{(x)}$ for all $x\in X$.
 In particular, $\vp_g(U_i)\rightarrow \infty$ as $\alpha_i^{(g)}\rightarrow \infty$. Thus Theorem \ref{thm-neg-char} implies that, if $\sup\{\vp_g(U):\;U\in \mathcal A^{\mathsf{elm}}(G_0)\}=\infty$, then this supremum can be obtained by restricting to a subfamily of elementary atoms $U_i$ whose supports are each equal apart from one varying element.
 Obtaining the $\Z_+$ linear combination analog of Theorem \ref{thm-neg-char}  using atoms instead of elementary atoms (including the stronger statement regarding the uniform bound $N$) will be one of the main steps in our characterization result, and one which needs further machinery regarding finitary sets developed in the next subsections.

\begin{theorem}\label{thm-neg-char} Let $\Lambda\subseteq \R^d$ be a full rank lattice, where $d\geq 0$, and let
$G_0\subseteq \Lambda$ be a finitary subset with $\C(G_0)=\R^d$.  Then
$$G_0^\diamond=\{g\in G_0:\; \sup\{\vp_{g}(U):\; U\in \mathcal A^{\mathsf{elm}}(G_0)\}=\infty\}.$$ Indeed, there is a bound $N>0$ such that $\vp_g(U)\leq N$ for any $U\in \mathcal A^{\mathsf{elm}}(G_0)$ and $g\in G_0\setminus G_0^\diamond$.
\end{theorem}

\begin{proof}  By Proposition \ref{prop-diamond-basic-containment}, we have the basic inclusion
$$G_0^\diamond\subseteq \{g\in G_0:\; \sup\{\vp_{g}(U):\; U\in \mathcal A^{\mathsf{elm}}(G_0)\}=\infty\}.$$ It remains to establish the reverse inclusion, which follows from the stronger conclusion about the existence of the uniform bound $N$.
To this end, assume by contradiction that $\{x^{(0)}_i\}_{i=1}^\infty $ is a sequence of terms $x_i^{(0)}\in G_0\setminus G_0^\diamond$ and
$\{U_i\}_{i=1}^\infty$ is a sequence of elementary atoms, so  $U_i\in \mathcal A^{\mathsf{elm}}(G_0)$, with $\vp_{x_i^{(0)}}(U_i)\rightarrow \infty$.
Then $|\supp(U_i)|\leq d+1$ for all $i$, so by passing to a subsequence, we can w.l.o.g. assume all $U_i$ have the same cardinality support (say) $s+1$ with $s\in [1,d]$ (if $s=0$, then $U_i=0$ for all $i$, contradicting that $\vp_{x_i^{(0)}}(U_i)\rightarrow \infty$). Let $\supp(U_i)=\{x_i^{(0)},x_i^{(1)},\ldots,x_i^{(s)}\}$ for $i\geq 1$. Since $|\supp(U_i)|=s+1\geq 2$, we have $0\notin\supp(U_i)$ for all $i$. If $\{x_i^{(0)}\}_{i=1}^\infty$ is unbounded, then by passing to a subsequence, we can w.l.o.g. assume $\{x_i^{(0)}\}_{i=1}^\infty$ is radially convergent with unbounded limit $g_0$.
On the other hand, if $\{x_i^{(0)}\}_{i=1}^\infty$ is bounded, then, since $x_i^{(0)}\in G_0\subseteq \Lambda$ and since any bounded subset of lattice points is finite, we can, by passing to a subsequence, w.l.o.g. assume that the sequence $\{x_i^{(0)}\}_{i=1}^\infty$ is constant, say with $x_i^{(0)}=g_0$ for all $i\geq 1$. Repeating this argument for each $\{x_i^{(j)}\}_{i=1}^\infty$ for $j=1,2,\ldots,s$, we can likewise assume that, for every $j\in [0,s]$, either  $\{x_i^{(j)}\}_{i=1}^\infty$ is constant, say with $x_i^{(j)}=g_j$ for all $i\geq 1$, or else $\{x_i^{(j)}\}_{i=1}^\infty$ is radially convergent  with unbounded limit $g_j$. Partition $[0,s]=I_u\cup I_b$, with $I_u$ consisting of all indices $j\in [0,s]$ such that $\{x_i^{(j)}\}_{i=1}^\infty$ is unbounded, and $I_b$ consisting of all indices $j\in [0,s]$ such that $x_i^{(j)}=g_j$ for all $i\geq 1$. Thus $$\{g_j:\; j\in I_b\}\subseteq G_0\quad\und\quad \{g_j:\; j\in I_u\}\subseteq G_0^\infty.$$

If $I_u=\emptyset$, then $\supp(U_i)=\{g_0,\ldots,g_s\}$ is constant for all $i$, in which case the elementary atoms $U_i$ are all equal to one another (in view of Proposition \ref{prop-char-minimal-pos-basis}.8), contradicting that $\vp_{x_i^{(0)}}(U_i)\rightarrow \infty$. Therefore we can assume $I_u\neq \emptyset$.

Suppose $0\notin \C^*(g_0,\ldots,g_s)$. Then Lemma \ref{lem-finite-nozs} gives a co-dimension one subspace $\mathcal H\subseteq \R^d$ such that $g_0,\ldots,g_s\in \mathcal H_+^\circ$.
 In consequence, since the sequences $\{x^{(j)}_i\}_{i=1}^\infty$ are radially convergent with limit $g_j/\|g_j\|$, it follows that, for all sufficiently large $i$, we also have  $x_i^{(0)},\ldots,x_i^{(s)}\in \mathcal H_+^\circ$, ensuring that $0\notin \C^*(x_i^{(0)},\ldots,x_i^{(s)})=\C^*(\supp(U_i))$. However, this contradicts that $U_i\in \mathcal A^{\mathsf{elm}}(G_0)$ is an elementary atom. So we instead conclude that $0\in \C^*(g_0,\ldots,g_s)$. Thus, since $g_j\neq 0$ for all $j\in [0,s]$, there is a minimal positive basis $B\subseteq \{g_0,\ldots,g_s\}$. Let $$X=B\cap \{g_j:\; j\in I_b\}\subseteq G_0.$$

Since $I_u\neq \emptyset$, it follows that $X\subset \supp(U_i)$ is a proper subset of $\supp(U_i)$ for all $i\geq 1$. Thus, since each $\supp(U_i)$ is a minimal positive basis, it follows from Proposition \ref{prop-char-minimal-pos-basis}.3 applied to $U_i$ that $X$ is linearly independent. Thus, since $B$ is a minimal positive basis, it follows from Proposition \ref{prop-char-minimal-pos-basis}.2 applied to $B$ that $X\subset B$ is a proper subset, which then implies  that $|B\cap \{g_j:\; j\in I_u\}|\geq 1$. On the other hand, if $|B\cap \{g_j:\; j\in I_u\}|>1$, then this contradicts that $G_0$ is finitary, which ensures  (as described after Theorem \ref{thm-keylemmaII}) that any minimal positive basis with elements from $G_0\cup G_0^{\infty}$ can involve at most one element from $G_0^\infty$. So we conclude that $|B\cap \{g_j:\; j\in I_u\}|= 1$, and thus $|X|=|B|-1$. Consequently, if $g_0\in X$, then Proposition \ref{prop-G_0diamond-1st-easy-equiv}.1 implies that $g_0=x_i^{(0)}\in G_0^\diamond$ for all $i\geq 1$, contradicting that $x_i^{(0)}\in G_0\setminus G_0^\diamond$ by assumption. So we may instead assume $g_0\notin X$.

Suppose $B=\{g_0,\ldots,g_s\}$. Then $X=\{g_1,\ldots,g_s\}$ as $g_0\notin X$, ensuring that $-x_i^{(0)}\in
\C^\circ(x_i^{(1)},\ldots,x_i^{(s)})=\C^\circ(g_1,\ldots,g_s)$ for all $i\geq 1$ with $g_1,\ldots,g_s\in G_0\subseteq \Lambda$ linearly independent elements (as the sequences $\{x_i^{(j)}\}_{i=1}^\infty$ for $j\in [1,s]$ will all be constant). Applying Lemma \ref{lem-lattice-smithnormal} to the linearly independent set $X=\{g_1,\ldots,g_s\}$, we find that there exists some constant $N>0$ such that $\vp_{x_i^{(0)}}(U_i)\leq N$ for all $i$ (as there is a unique elementary atom with given support by Proposition \ref{prop-char-minimal-pos-basis}.8), contradicting that $\vp_{x_i^{(0)}}(U_i)\rightarrow \infty$. So we instead conclude that $B\subset \{g_0,\ldots,g_s\}$ is a proper subset.

By re-indexing the $g_j$ with $j\in [1,s]$, we can assume $$X=\{g_j:\;j\in [s'+1,s]\}\subseteq G_0\quad\und\quad B=X\cup \{g_t\}\quad\mbox{ for some $t\in [0,s']$ and  $s'\in [1,s-1]$}$$ (we have $s'>0$ in view of $B\subset \{g_0,\ldots,g_s\}$ being proper). Note $t\in I_u$, and $s'\in [1,s-1]$ is only possible when $d\geq s\geq 2$, which means the proof is now complete for $d=1$ (and trivial for $d=0$), allowing us to proceed inductively on $d$.

Each $g_j\in B$ defines a one-dimensional half-space $\mathbf g_j=\R_+ g_j$, and letting $\mathcal B$ be the corresponding set of half-spaces, it follows that
$\mathcal R=((\mathcal B\setminus \{\mathbf g_{t}\})\cup \{\mathbf g_{t}\})$ is a purely virtual Reay system over $G_0$ (in view of $t\in I_u$ and since all other $g_j\in B\setminus \{g_t\}=X\subseteq G_0$). Let $\pi:\R^d\rightarrow \R^\cup \la \mathcal B\ra^\bot=\R\la g_{s'+1},\ldots,g_s\ra^\bot=\R\la X\ra^\bot$ be the orthogonal projection. Proposition \ref{prop-finitary-Modulo-Inheritence}  ensures that $\pi(\Lambda)$ is a full rank lattice in $\C(\pi(G_0))=\R\la X\ra^\bot$ with $\pi(G_0)$ finitary.

 Note that $\pi(x)=0$ for all $x\in X$. On the other hand, if $J\subset [0,s']$ is a proper subset, then (for any $i\geq 1$) the elements $\pi(x_i^{(j)})$ for $j\in J$ must be distinct and linearly independent, as otherwise  $\{x_i^{(j)}:\;j\in J\}\cup  X=\{x_i^{(j)}:\;j\in J\}\cup \{x_i^{(s'+1)},\ldots,x_i^{(s)}\}$ would be a linearly dependent, proper subset of $\supp(U_i)$, contradicting that $\supp(U_i)$ is a minimal positive basis. In particular, since $s'\geq 1$, it follows that $\pi(x_i^{(j)})\neq 0$ for all $j\in [0,s']$ and  $i\geq 1$. Since there is a strictly positive linear combination of the elements from $\supp(U_i)$ equal to zero, this is also true for the elements from $\pi(\supp(U_i))\setminus \{0\}=\pi(\supp(U_i)\setminus X)$.
 Thus it now follows from Proposition \ref{prop-char-minimal-pos-basis}.2 that $\pi(\supp(U_i))\setminus \{0\}=\{\pi(x_i^{(0)}),\ldots,\pi(x_i^{(s')})\}$ is a minimal positive basis of size $s'+1$ for every $i\geq 1$.
 Consequently, Proposition \ref{prop-char-minimal-pos-basis}.8 now further implies that there is a unique atom $A_i\in \mathcal A(\pi(G_0))$ with $\supp(A_i)\subseteq \pi(\supp(U_i))\setminus \{0\}$, and we have $A_i\in \mathcal A^{\mathsf{elm}}(\pi(G_0))$ with $\supp(A_i)=\pi(\supp(U_i))\setminus \{0\}$ for this atom. Hence each $\pi(U_i)=A_i^{[N_i]}\bdot 0^{[M_i]}$ for some integers $N_i,\,M_i>0$, and there is a unique subsequence $V_i\mid U_i$ such that $\pi(V_i)=A_i$ (as $\pi$ is injective on the elements of $\supp(U_i)\setminus X$), meaning $$U_i\bdot W_i^{[-1]}=V_i^{[N_i]},$$ where $W_i\mid U_i$ is the subsequence  consisting of all terms from $X$.
 Since each $U_i$ is zero-sum, and since $\vp_x(W_i)>0$ for every $x\in X$ with $X$ linearly independent, it follows that $$0=\sigma(U_i\bdot W_i^{[-1]})+\sigma(W_i)=N_i \sigma(V_i)+\sigma(W_i)\in N_i\sigma(V_i)+\C^\circ(X).$$ Thus, since $N_i>0$ and $\supp(V_i)\subseteq G_0\subseteq \Lambda$, it follows that $$\sigma(V_i)\in \Lambda\cap -\C^\circ(X)\subseteq \Lambda \cap \R\la X\ra\quad\mbox{ for all $i\geq 1$}.$$
Since $X\neq \emptyset$ is linearly independent, we can apply Lemma \ref{lem-lattice-smithnormal} with $g=\sigma(V_i)$ to conclude there is a positive integer $N>0$ such that $N\sigma(V_i)\in \Z\la X\ra$  and, moreover, there is an elementary atom $S_i$ with $\supp(S_i)=X\cup \{\sigma(V_i)\}$ and $\vp_{\sigma(V_i)}(S_i)\leq N$, for every $i\geq 1$.
 Then $$T_i=V_i^{[\vp_{\sigma(V_i)}(S_i)]}\bdot {\prod}^\bullet_{x\in X}x^{[\vp_x(S_i)]}\in \mathcal F(G_0)$$ is a zero-sum sequence with $\supp(T_i)=\supp(U_i)$. In consequence, since $U_i$ is an elementary atom, it follows from Proposition \ref{prop-char-minimal-pos-basis}.8 that $\vp_g(U_i)\leq \vp_g(T_i)$ for every $g\in \supp(U_i)$. In particular, $\vp_{x_i^{(0)}}(U_i)\leq \vp_{x_i^{(0)}}(T_i)=\vp_{\sigma(V_i)}(S_i)\cdot \vp_{x_i^{(0)}}(V_i)\leq N\cdot \vp_{x_i^{(0)}}(V_i)$. Thus, since $\vp_{x_i^{(0)}}(U_i)\rightarrow \infty$, we conclude that $\vp_{x_i^{(0)}}(V_i)\rightarrow \infty$ as well. However, since $\pi(V_i)=A_i\in \mathcal A^{\mathsf{elm}}(\pi(G_0))$, this implies that $\vp_{\pi(x_i^{(0)})}(A_i)\rightarrow \infty$.
 By the induction hypothesis applied to $\pi(G_0)$, we know there is a bound $N'>0$ such that $\pi(g)\in \pi(G_0)^\diamond $ whenever $\vp_{\pi(g)}(A)>N'$ for some $A\in \mathcal A^{\mathsf{elm}}(\pi(G_0))$. Thus, since $\vp_{\pi(x_i^{(0)})}(A_i)\rightarrow \infty$, it follows that
 $\pi(x_i^{(0)})\in \pi(G_0)^\diamond$ for all sufficiently large $i$. But now Proposition \ref{prop-finitary-diamond-modulo-transfer} implies that $x_i^{(0)}\in G_0^\diamond$ for all sufficiently large $i$, contradicting our assumption that $x_i^{(0)}\in G_0\setminus G_0^\diamond$ for all $i\geq 1$, which  completes the proof.
\end{proof}

We continue with an important property of  purely virtual Reay systems.

\begin{proposition}\label{prop-finitary-diamond-containment-nonmax} Let $\Lambda\subseteq \R^d$ be a full rank lattice, where $d\geq 0$,  let
$G_0\subseteq \Lambda$ be a finitary subset with $\C(G_0)=\R^d$, let $\mathcal R=(\mathcal X_1\cup \{\mathbf v_1\},\ldots,\mathcal X_s\cup\{\mathbf v_s\})$ be a purely virtual Reay system over $G_0$, let $\mathcal X=\mathcal X_1\cup \ldots\cup \mathcal X_s$, and  let $\pi:\R^d\rightarrow \R^\cup \la \mathcal X\ra^\bot$ be the orthogonal projection.
 If $U\in \Fc(G_0)$ with $\pi(U)\in \mathcal A^{\mathsf{elm}}(\pi(G_0))$ and $|\supp(U)|=|\supp(\pi(U))|$, then   $\mathsf{wt}(-g)\leq 1$, where $g=\sigma(U)$.  Moreover, if $\mathsf{wt}(-g)= 1$, then $\supp(U)\subseteq G_0^\diamond$. In particular,   $(G_0\cap \R^\cup \la \mathcal X\ra)\setminus -\C^\cup(\mathcal X)\subseteq G_0^\diamond$.
\end{proposition}

\begin{proof}
Since $\pi(U)\in \mathcal A^{\mathsf{elm}}(\pi(G_0))$ is an atom,  we have $g=\sigma(U)\in \ker \pi=\R^\cup\la \mathcal X\ra$. Furthermore, since $|\supp(U)|=|\supp(\pi(U))|$, we conclude that $\pi$ is injective on $\supp(U)$ with $\supp(\pi(U))$ either $\{0\}$ or a minimal positive basis.
Note $g=\sigma(U)\in \C_\Z(G_0)\subseteq \Lambda$ as $U\in \Fc(G_0)$.
Let $\mathcal B=\supp_\mathcal R(-g)$, in which case $-g\in \C^\cup(\mathcal B)^\circ$. Thus the support set $\mathcal B$ minimally encases the limit $-\vec u:=(-g/\|g\|)$ (cf. the comments after Lemma \ref{lem-orReay-mincase-t=1}), with the encasement trivially urbane as the limit $\vec u$ is composed of a single coordinate. We may assume $g\neq 0$ as the proposition is   trivial for $g=0$.
By
Proposition \ref{prop-VReay-modularCompletion} applied to $\mathcal B$, there is a virtual Reay system $\mathcal R'=(\mathcal Y_1\cup \{\mathbf w_1\},\ldots,\mathcal Y_\ell\cup \{\mathbf w_\ell\})$ over $G_0\cup \{g\}$ with $\mathcal Y_\ell=\mathcal B$ and $\vec u_{\mathbf w_\ell}=\vec u$. Moreover, each $\mathbf w_j$ with $j\in [1,\ell-1]$ is defined by a strict truncation of a limit $\vec u_\y$ for some  $\y\in \darrow \mathcal B$, and is thus fully unbounded by   Proposition \ref{prop-VReay-Lattice} (as $\mathcal R$ is purely virtual and  anchored in view of $G_0$ being finitary).
Let $\mathcal U$ be the set of one-dimensional half-spaces generated by the elements from $\supp(U)$, in which case $\supp(U)$ is a set of representatives for $\mathcal U$.
Let $\pi':\R^d\rightarrow \R^\cup\la \mathcal Y_1\cup \ldots\cup \mathcal Y_{\ell-1}\ra^\bot$ be the orthogonal projection. Then \be\label{sret}\pi'(Y_\ell\cup \{ w_\ell\})=\pi'(B\cup \{g\})\ee is a minimal positive basis of size $|\mathcal B|+1$ by (OR2) with $g=\sigma(U)$.
Thus $0$ is a strictly positive linear combination of the elements from $\pi'(\supp(U)\cup B)$.
Since $\mathcal Y_1\cup \ldots\cup \mathcal Y_{\ell-1}
\subseteq \mathcal X$ (by Proposition \ref{prop-VReay-modularCompletion}(a)),  we have $\ker \pi'\subseteq \R^\cup\la \mathcal X\ra=\ker \pi$, and since
$ \mathcal B\subseteq \mathcal X\cup \{\mathbf v_1,\ldots,\mathbf v_s\}$, we have $\R\la B\ra\subseteq \R^\cup \la \mathcal B\ra\subseteq \R^\cup\la \mathcal X\ra=\ker \pi$.

If $|\supp(U)|=1$, then $U=g$ and $\pi'$ is injective on $\supp(U)\cup B$ by (OR2) (as it is injective on $B\cup \{g\}$). If  $|\supp(U)|>1$, then $\pi(y)\neq 0$ for all $y\in \supp(U)$, while $\pi(x)=0$ for all $x\in B$, ensuring that $\pi'(x)\neq \pi'(y)$ for $x\in B$ and $y\in \supp(U)$ in view of $\ker \pi'\leq \ker \pi$.  Thus, since $\pi'$ is injective on $\supp(U)$ (as $\pi$ is injective on $\supp(U)$ with $\ker \pi'\leq \ker \pi$) and since $\pi'$ is injective on $B$ (by (OR2)), we conclude that $\pi'$ is injective on $B\cup \supp(U)$ in both cases.

Let us show that $\pi'(B\cup \supp(U))$ is a minimal positive basis of size $|\mathcal B|+|\mathcal U|$ (as $\pi'$ is injective on $B\cup \supp(U)$). If $|\supp(U)|=1$, then $U=g$ with $\supp(U)=\{g\}$, in which case this follows from \eqref{sret}. Next assume $|\supp(U)|>1$ and consider an arbitrary linear combination
 $\Summ{x\in B}\alpha_x \pi'(x)+\Summ{y\in \supp(U)}\beta_y \pi'(y)=0$ for some $\alpha_x,\,\beta_y\in \R$.
  Since $\pi'(B\cup \{g\})$ is a minimal positive basis of size $|B|+1$, it follows that $\pi'(B)$ is linearly independent, and thus $\beta_y\neq 0$ for some $y\in \supp(U)$, allowing us to w.l.o.g. assume $\beta_y>0$ for some $y\in \supp(U)$.
  Applying $\pi$ to this linear combination, we find that $\Summ{y\in \supp(U)}\beta_y \pi(y)=0$.
  Thus, since $\supp(\pi(U))$ is a minimal positive basis of size $|\supp(U)|$, it now follows from Proposition \ref{prop-char-minimal-pos-basis}.5 that $\beta_y>0$ for all $y\in \supp(U)$.
Since $\Summ{y\in \supp(U)}\vp_y(U)y=\sigma(U)=g$ with $\pi(g)=0$ (as $g=\sigma(U)\in \ker \pi$), it follows from Proposition \ref{prop-char-minimal-pos-basis}.5 that $\Summ{y\in \supp(U)}\beta_y y$ is a  positive scalar multiple of $g$ (indeed, each $\beta_y=\alpha\vp_y(U)$ for the same $\alpha>0$).
But now $\Summ{x\in B}\alpha_x \pi'(x)=-\Summ{y\in \supp(U)}\beta_y\pi'(y)$ is positive scalar multiple of $-\pi'(g)$, which combined with $\pi'(B\cup \{g\})$ being a minimal positive basis ensures that $\alpha_x>0$ for all $x\in B$ (by Proposition \ref{prop-char-minimal-pos-basis}.3). Combined with the fact that $\beta_y>0$ must also hold for all $y\in \supp(U)$, shown earlier, we conclude that  any proper subset of $\pi'(B\cup \supp(U))$ is linearly independent. Thus, since we  showed earlier that  $0$ is a strictly positive linear combination of the elements from $\pi'(B\cup \supp(U))$, it now follows from  Proposition \ref{prop-char-minimal-pos-basis}.2 that $\pi'(B\cup \supp(U))$ is a minimal positive basis in the case $|\supp(U)|>1$ as well.

 As a result,  $\mathcal R''=(\mathcal Y_1\cup \{\mathbf w_1\},\ldots,\mathcal Y_{\ell-1}\cup \{\mathbf w_{\ell-1}\}, ((\mathcal B\cup \mathcal U)\setminus \{\x\})\cup \{\x\})$
will be a virtual Reay system over $G_0$ for any $\x\in \mathcal B\cup \mathcal U$ with all half-spaces from $\mathcal U$ having trivial boundary with representative sequence a constant term equal to an element from $U$. Consequently, if $\wt(-g)\geq 1$, then we can choose $\x\in \mathcal B$ with $\x=\mathbf v_j$ for some $j\in [1,s]$, thus ensuring that $\vec u_\x$ is fully unbounded (as $\mathcal R$ is purely virtual), in which case $\mathcal R''$ is purely virtual. In this case, since $G_0$ is finitary, we must have $\mathcal R''$ anchored. In particular,  all  limits defining the half-spaces from $\mathcal B\setminus\{\x\}$ are anchored. Thus $\wt(-g)=1$ as $\mathcal R''$ is purely virtual, and Proposition \ref{prop-finitary-basics}.1 applied to $\mathcal R''$ yields $\supp(U)\subseteq G_0^\diamond$,  as desired.
\end{proof}

\begin{lemma}\label{lem-shadow-limits}
Let $\mathcal E\subseteq \R^d$ be a subspace, let $\pi:\R^d\rightarrow \mathcal E^\bot$ be the orthogonal projection, let $G_0\subseteq \R^d$ be a subset such that $\pi(G_0)$ is finite and let $$\wtilde G_0=\{\sigma(U):\; U\in \mathcal F(G_0),\,\pi(U)\in \mathcal A(\pi(G_0))\}.$$ Suppose $C\subseteq \R^d$ is a polyhedral cone such that $C$ encases every $\vec u\in G_0^{\mathsf{lim}}$. Then $C$ also encases every  $\vec  u\in \wtilde G_0^{\mathsf{lim}}$.
\end{lemma}

\begin{proof} If $C$ is empty, then $G_0^{\mathsf{lim}}$ must also be empty, ensuring that $G_0$ is bounded.  Since  $\pi(G_0)$ is finite, it follows by Dickson's Theorem \cite[Theorem 1.5.3, Corollary 1.5.4]{alfredbook} that  $\mathcal A(\pi(G_0))$ is finite, and now $\wtilde G_0$ is also bounded, ensuring that ${\wtilde G_0}^{\mathsf{lim}}$ is empty. In such case, the lemma holds trivially. Therefore we now assume $C$ is nonempty.
Let $\{x_i\}_{i=1}^\infty$ be an asymptotically  filtered sequence of terms from $\wtilde G_0$ with fully unbounded limit $\vec u=(u_1,\ldots,u_t)$, say with each $x_i=\sigma(U_i)$ for some $U_i\in \mathcal F(G_0)$ with $\pi(U_i)\in \mathcal A(\pi(G_0))$. As previously remarked, since  $\pi(G_0)$ is finite, it follows that  $\mathcal A(\pi(G_0))$ is also finite. As a result, by passing to an appropriate subsequence, we can w.l.o.g. assume $\pi(U_i)=V$ is constant and equal to the same atom for all $i$. Let $V=a_1\bdot\ldots\bdot a_\ell$ with $a_i\in \pi(G_0)$. Then each $x_i=y_i^{(1)}+\ldots+y_i^{(\ell)}$ for some $y_i^{(j)}\in G_0$ with $\pi(y_i^{(j)})=a_j$ for all $j\in [1,\ell]$.
By passing to a subsequence, we can w.l.o.g. assume each $\{y_i^{(j)}\}_{i=1}^\infty$ is either bounded or an asymptotically  filtered sequence of terms from $G_0$ with complete fully unbounded limit $\vec u_j$.  By hypothesis, $C$ encases every $\vec u_j$, whence Proposition  \ref{prop-encasementcones-contain-aprox-seq}.3 implies that $\{y_i^{(j)}\}_{i=1}^\infty$ is bound to $C$. This is trivially true for any bounded sequence  $\{y_i^{(j)}\}_{i=1}^\infty$ as $C$ is nonempty, and so $\{y_i^{(j)}\}_{i=1}^\infty$ is bound to $C$ for each $j\in [1,\ell]$.
Thus, for each $j\in [1,\ell]$, there is a constant $N_j\geq 0$ such that, for every  $i\geq 1$, there is some  $z_i^{(j)}\in C$ such that $\|y_i^{(j)}-z_i^{(j)}\|\leq N_j$. Since $C$ is convex, it follows that $z_i=z_i^{(1)}+\ldots+z_i^{(\ell)}\in C$ with $\|x_i-z_i\|\leq \Sum{j=1}{\ell}N_j$ for all $i\geq 1$ (by the triangle inequality). Hence $\{x_i\}_{i=1}^\infty$ is bound to $C$, whence Lemma \ref{lem-weak-nearness-implication} implies that there is an asymptotically filtered sequence $\{x'_i\}_{i=1}^\infty$ of terms $x'_i\in C$ having fully unbounded limit $\vec u$, in which case
Proposition \ref{prop-finite-union--convergence-encasement} ensures that $C$ encases $\vec u$, as desired.
\end{proof}

Next, we derive additional properties regarding the geometry of a finitary set $G_0$ via its  maximal purely virtual Reay systems.
Theorem \ref{thm-finitary-diamond-containment}.2 ensures that a finitary set $G_0$ has a linearly independent subset $X\subseteq G_0^\diamond\subseteq G_0$ such that $G_0$ is bound to $-\C(X)$, meaning the set $G_0$ must be concentrated around the simplicial cone $-\C(X)$.
Theorem \ref{thm-finitary-diamond-containment}.3 then further implies that the \emph{finitely} generated convex cone $\C(X)$  perfectly approximates the (possibly) non-finitely generated convex cone $\C^\cup(\mathcal X)$ in regards to containment of elements from $-\wtilde G_0$, giving our first example of finite-like behavior for finitary sets. Combined with  Proposition \ref{prop-finitary-diamond-containment-nonmax}, we obtain  restrictions on the location of the elements from $G_0$ in relation to the simplicial cone $\C(X)$. Note $G_0\cap \R\la X\ra\subseteq \wtilde G_0$ since the single term equal to $0$ is always an atom. Since we may not have a maximal Reay system spanning the entire space $\R^d$, the set $\wtilde G_0$ as well as Proposition \ref{prop-finitary-diamond-containment-nonmax}  will be used for dealing with elements lying outside $\R\la  X\ra$.

\begin{theorem}\label{thm-finitary-diamond-containment} Let $\Lambda\subseteq \R^d$ be a full rank lattice, where $d\geq 0$,  let
$G_0\subseteq \Lambda$ be a finitary subset with $\C(G_0)=\R^d$, let $\mathcal R=(\mathcal X_1\cup \{\mathbf v_1\},\ldots,\mathcal X_s\cup\{\mathbf v_s\})$ be a  maximal purely virtual Reay system over $G_0$, let $\mathcal X=\mathcal X_1\cup \ldots\cup \mathcal X_{s}$, and let $\pi:\R^d\rightarrow \R^\cup \la \mathcal X\ra^\bot$ be the orthogonal projection.
\begin{itemize}
\item[1.] $-\mathcal X$ encases every  $\vec u\in G_0^{\mathsf{lim}}$. In particular, $G_0$ is bound to $-\C^\cup(\mathcal X)$.
\item[2.] For any tuple $k=(i_\z)_{\z\in \mathcal X}$ with all $i_\z$ sufficiently large, $-\C(X(k))$ encases every  $\vec u\in G_0^{\mathsf{lim}}$. Moreover, $G_0$ is bound to $-\C(X(k))$ with  $X(k)\subseteq G_0^\diamond$ a linearly independent set.


\item[3.] Let $\wtilde G_0=\{\sigma(U):\; U\in \mathcal F(G_0),\,\pi(U)\in \mathcal A(\pi(G_0))\}$ and let $\mathcal X'\subseteq\mathcal X$. Then $$\wtilde G_0\cap -\C (\darrow X'(k))=\wtilde G_0\cap -\C^\cup(\mathcal X')$$ for any tuple $k=(i_\z)_{\z\in \darrow \mathcal B}$ with all $i_\z$ sufficiently large.
\end{itemize}
\end{theorem}

\begin{proof}
1.
Let $\{x_i\}_{i=1}^\infty$ be an arbitrary asymptotically filtered sequence of terms $x_i\in G_0$ with fully unbounded limit $\vec u=(u_1,\ldots,u_t)$. Since $\mathcal R$ is a \emph{maximal} purely virtual Reay system over the finitary set $G_0$, it follows that $u_1,\ldots,u_t\in \R^\cup \la \mathcal X\ra$ (otherwise $\{\pi(x_i)\}_{i=1}^\infty$ would be an unbounded sequence of terms $\pi(x_i)\in \pi(G_0)$,
 contradicting that $\pi(G_0)$ is finite in view of the maximality of $\mathcal R$). Thus Proposition \ref{prop-finitary-basics}.3 implies that $\C^\cup(\mathcal X)$ encases $-\vec u$. Since $\vec u\in G_0^{\mathsf{lim}}$ was arbitrary, and since $\overline{\C^\cup(\mathcal X)}$ is a polyhedral cone by Proposition \ref{prop-orReay-BasicProps}.2, Theorem \ref{thm-nearness-characterization}.4 implies  that $G_0$ is bound to $-\overline{\C^\cup(\mathcal X)}$, and thus also to $-\C^\cup(\mathcal X)$.

2. Since $G_0$ is finitary and $\mathcal R$ is purely virtual, it follows that $\mathcal R$ is anchored. For any tuple $k=(i_\z)_{\z\in \mathcal X}$ with all $i_\z$ sufficiently large, we have $\darrow \tilde X(k)=\darrow X(k)=X(k)$ by Proposition \ref{prop-VReay-Lattice}, in which case $X(k)=\darrow X(k)$ is a set of representatives for $\mathcal X=\mathcal X_1\cup \ldots\cup \mathcal X_s$, ensuring that $X(k)$ is linearly independent by Proposition \ref{prop-reay-basis-properties}.1. Moreover, by Proposition \ref{prop-finitary-basics}.1, we have $ X(k)\subseteq G_0^\diamond$ once all $i_\z$ are sufficiently large. It remains to show there is some $N>0$ so that, so long as all $i_\z\geq N$, then $G_0$ is bound to $-\C(X(k))$. In view of Theorem \ref{thm-nearness-characterization}, this is equivalent to showing any fully unbounded limit $\vec u$ of an asymptotically filtered sequence of terms from $G_0$ has $-\vec u$ encased by $X(k)$. By Item 1, each such $-\vec u$ is encased by $\mathcal X$, and thus must be minimally encased by some support set $\mathcal B\subseteq \mathcal X$, with the encasement urbane in view of $\mathcal B\subseteq \mathcal X$.

 Let $\mathfrak X$ consist of all subsets $\mathcal B\subseteq \mathcal X$ for which there is some $\vec u\in G_0^{\mathsf{lim}}$  with $-\vec u$ minimally encased by $\mathcal B$. Note $\mathfrak X$ is finite as $\mathcal X$ is finite.  Let $\mathcal B\in \mathfrak X$ be arbitrary. Then there is some fully unbounded limit $\vec u=(u_1,\ldots,u_t)$ of an asymptotically  filtered sequence $\{x_i\}_{i=1}^\infty$ of terms from $G_0$ with $-\vec u$ minimally encased by $\mathcal B$. We fix the tuple $\vec u$ for $\mathcal B$ and will show $-\C(\darrow B(k))$ encases every $\vec v=(v_1,\ldots,v_{t'})\in G_0^{\mathsf{lim}}$ with $v_1,\ldots,v_{t'}\in \R^\cup \la \mathcal B\ra$ for any tuple  $k=(i_\z)_{\z\in \mathcal X}$ with all $i_\z\geq N_\mathcal B$, for some $N_\mathcal B$ that depends only on $\mathcal B$ and the fixed tuple $\vec u$ (thus not dependent on the potentially infinite number of varying tuples $\vec v\in G_0^{\mathsf{\lim}}$). Taking $N=\max_{\mathcal B\in \mathfrak X}N_\mathcal B$, which exists as $\mathcal X$ is finite, Item 2 will  follow as every fully unbounded limit of an asymptotically  filtered sequence of terms from $G_0$ is encased by some $\mathcal B\in \mathfrak X$, as noted above the definition of $\mathfrak X$.

 Proposition \ref{prop-min-encasement-char} ensures that if a subset $X\subseteq \R^d$ minimally encases $-\vec w=-(w_1,\ldots,w_{r})$ and $\vec w'=(w_1,\ldots,w_{r},w_{r+1},\ldots,w_{r+r'})$ has $w_i\in \R\la X\ra$ for all $i>r$, then $X$ also minimally encases $-\vec w'$. Indeed, if $\mathcal F=(\mathcal E_1,\ldots,\mathcal E_\ell)$ is the filter given by the application of Lemma \ref{prop-min-encasement-char} to the minimal encasement of $-\vec w$ by $X$, then the assumption $w_i\in \R\la X\ra$ for all $i>r$ ensures that it remains a compatible filter for $-\vec w'$ having the same associated set of indices (note $X$ encasing $-\vec w$ ensures $w_i\in \R\la X\ra$ for all $i\leq r$), and then the oriented Reay system from Lemma \ref{prop-min-encasement-char} applied to the minimal encasement of $-\vec  w$ also shows that $X$ minimally encases $-\vec w'$.
 We use this observation several times below.

 For any tuple $k=(i_\z)_{\z\in \mathcal B}$ with all $i_\z$ sufficiently large, we have $\darrow B(k)=\darrow \tilde B(k)$ by Proposition \ref{prop-VReay-Lattice}, in which case $\darrow B(k)$ is a set of representatives for $\darrow \mathcal B\subseteq \mathcal X_1\cup \ldots\cup \mathcal X_s$. Thus Proposition \ref{prop-orReay-BasicProps}.1 implies that $\R^\cup \la\mathcal B\ra=\R\la \darrow B(k)\ra$. Since $-\vec u$ is minimally encased by $\mathcal B$ urbanely,
  let $1=r_1<\ldots<r_\ell<r_{\ell+1}=t+1$ be the indices and $\emptyset =\mathcal C_0\prec \mathcal C_1\prec\ldots\prec \mathcal C_{\ell-1}\prec \mathcal C_\ell=\mathcal B\subseteq \mathcal X$ the support sets given by Proposition \ref{prop-orReay-minecase-char}.2.
 Then $\mathcal C_1$ minimally encases $-u_{r_1}=-u_1$, so  Lemma \ref{lemma-minencase-Rep} implies that $\darrow C_1(k)\cup \{u_{r_1}\}$ is a minimal positive basis for $\R^\cup\la \mathcal C_1\ra$ with $\darrow \tilde C_1(k)=\darrow C_1(k)$ minimally encasing $-u_1=-u_{r_1}$, and thus also minimally encasing
  $-(u_1,\ldots,u_{r_2-1})$ as $u_i\in \R^\cup \la \mathcal C_1\ra=\R\la \darrow C_1(k)\ra$ for $i<r_2$, so long as all $i_\z$ are sufficiently large.
  Let $\pi_1:\R^d\rightarrow \R^\cup \la \mathcal C_1\ra^\bot$ be the orthogonal projection.  Proposition \ref{prop-orReay-minecase-char} implies that $\mathcal B^{\pi_1}$ minimally encases $-\pi_1(\vec u)$ urbanely with associated support sets $\emptyset \prec \mathcal C_2^{\pi_1}\prec\ldots\prec \mathcal C_{\ell-1}^{\pi_1}\prec \mathcal C^{\pi_1}_\ell=\mathcal B^{\pi_1}\subseteq \mathcal X^{\pi_1}$ (as in the proof of Proposition \ref{prop-VReay-modularCompletion}).
  Note $\darrow \mathcal C_2^{\pi_1}=\pi_1(\darrow \mathcal C_2\setminus \darrow \mathcal C_1)$ in view of Proposition \ref{prop-orReay-BasicProps}.9.
 Proposition \ref{prop-orReay-minecase-char}.2 implies that $\mathcal C_2$ minimally encases $-(u_1,\ldots,u_{r_2})$ urbanely, and then Proposition \ref{prop-orReay-minecase-char}.4 ensures  $\mathcal C_2^{\pi_1}$ minimally encases $-\pi_1((u_1,\ldots,u_{r_2}))=-\pi_1(u_{r_2})/\|\pi_1(u_{r_2})\|$ urbanely, allowing us to apply Lemma \ref{lemma-minencase-Rep} to conclude $(\darrow \wtilde C_2^{\pi_1})(k)$
 minimally encases $-\pi_1(u_{r_2})$,  once all $i_\z$ are sufficiently large.

Since $\mathcal R$ is anchored, Proposition \ref{prop-VReay-modulo} ensures that $\pi_1(\mathcal R)$ is also anchored, while $\pi_1(\Lambda)$ is a lattice by Proposition \ref{Prop-lattice-homoIm} since $\ker \pi_1=\R\la \darrow C_1(k)\ra$ with $\darrow C_1(k)\subseteq G_0\subseteq \Lambda$. Thus Proposition \ref{prop-VReay-Lattice} implies that
$(\darrow \wtilde C_2^{\pi_1})(k)=(\darrow  C_2^{\pi_1})(k)$. Combined with the conclusions from the previous paragraph, we find that  $(\darrow \wtilde C_2^{\pi_1})(k)=(\darrow  C_2^{\pi_1})(k)=\pi_1(\darrow C_2\setminus \darrow C_1)(k)$ minimally encases  $-\pi_1(u_{r_2})$,  once all $i_\z$ are sufficiently large, so  $(\darrow C_1(k)\cup \{u_{r_1}\},(\darrow C_2\setminus \darrow C_1)(k)\cup \{u_{r_2}\})$ is a Reay system. Moreover, since $u_1,\ldots,u_{r_3-1}\in \R^\cup \la \mathcal C_2\ra=\R^\cup \la \darrow C_2(k)\ra$ by Proposition \ref{prop-orReay-minecase-char}.2, it follows that $\darrow C_2(k)$ minimally encases $-(u_1,\ldots,u_{r_3-1})$, once all $i_\z$ are sufficiently large. Iterating this argument, we find that
   $(\darrow C_1(k)\cup \{u_{r_1}\},(\darrow C_2\setminus \darrow C_1)(k)\cup \{u_{r_2}\}, \ldots,(\darrow C_\ell\setminus \darrow C_{\ell-1})(k)\cup \{u_{r_\ell}\})$ is an ordinary Reay system for any tuple $k=(i_\z)_{\z\in \mathcal X}$ with all $i_\z\geq N_\mathcal B$ sufficiently large.
   Moreover, for each $j\in [1,\ell]$, we have $-(u_1,\ldots,u_{r_j})^\triangleleft=-(u_1,\ldots,u_{r_j-1})$ minimally encased by $\darrow C_{j-1}(k)$. Thus replacing each element from $\darrow C_\ell(k)=\darrow B(k)$ with the one-dimensional half-space it defines and using the asymptotically filtered sequence $\{x_i\}_{i=1}^\infty$ as the representative sequence of each $\mathbf w_j$ gives rise to a purely virtual Reay system $\mathcal R_\mathcal B=(\mathcal Y_1\cup \{\mathbf w_1\},\ldots,\mathcal Y_\ell\cup \{
  \mathbf w_\ell\})$ over $G_0$ with $\R^\cup \la \mathcal Y_1\cup\ldots\cup \mathcal Y_\ell\ra=\R\la \darrow B(k)\ra=\R^\cup \la \mathcal B\ra$ and $\vec u_{\mathbf w_j}=(u_1,\ldots,u_{r_j})$ for all $j\in [1,\ell]$ (it is purely virtual as $-\vec u$ is fully unbounded). Applying Proposition \ref{prop-finitary-basics}.3 to $\mathcal R_\mathcal B$, we conclude that, if $\vec v=(v_1,\ldots,v_{t'})\in G_0^{\mathsf{lim}}$  with all $v_i\in \R^\cup \la \mathcal B\ra$, then $-\vec v$ is encased by $\mathcal Y_1\cup \ldots\cup \mathcal Y_\ell$, and thus also by $\darrow C_\ell(k)=\darrow B(k)\subseteq X(k)$ (as all half-spaces in the sets $\mathcal Y_i$ are one-dimensional), completing the proof as remarked earlier. 

3. Note $\darrow \mathcal X=\mathcal X$ since $\mathcal R=(\mathcal X_1\cup \{\mathbf v_1\},\ldots,\mathcal X_s\cup \{\mathbf v_s\})$ is an oriented Reay system. Each element $g\in \wtilde G_0$ has $g=\sigma(U)$ for some  $U\in \Fc(G_0)$, implying $g\in\C_\Z(G_0)\subseteq \Lambda$. Thus $$\wtilde G_0\subseteq \Lambda.$$
 Since $G_0$ is finitary, $\mathcal R$ is purely virtual and also  anchored. Thus,
 by Proposition \ref{prop-VReay-Lattice}, we have $X(k)=\tilde X(k)$ for any tuple $k=(i_\z)_{\z\in \mathcal X}$ with all $i_\z$ sufficiently large. Since $\mathcal R$ is a maximal purely virtual Reay system over $G_0$, it follows that $\pi(G_0)$ is finite. Thus Item 2 and Lemma \ref{lem-shadow-limits} imply $-\C(X(k))$ encases every $\vec u\in \wtilde G_0^{\mathsf{lim}}$ for any tuple $k=(i_\z)_{\z\in \mathcal X}$ with all $i_\z$ sufficiently large. Thus $\wtilde G_0$ is bound to $-\C(X(k))$ by  Theorem \ref{thm-nearness-characterization}.
Fix one tuple $\kappa=(\iota_\z)_{\z\in \mathcal X}$ such that $\wtilde G_0$ is bound to $-\C(X(\kappa))$ and $X(\kappa)=\tilde X(\kappa)$. By Proposition \ref{prop-VReay-RepBasics}.2 (applied to each $\x\in \mathcal B$), we can assume $\C(\darrow B(\kappa))\subseteq \C(\darrow B(k))$ for any tuple $k=(i_\z)_{\z\in \mathcal X}$ with all $i_\z$ sufficiently large and any $\mathcal B\subseteq \darrow\mathcal X'\subseteq \mathcal X$. Now $X(\kappa)\subseteq G_0\subseteq \Lambda$ is a linearly independent set of lattice points.
Since $\wtilde G_0\cap  \R^\cup\la \mathcal B\ra$  is bound to both $\R^\cup\la \mathcal B\ra$ and $-\C(X(\kappa))$, the former trivially, it follows from Corollary \ref{cor-nearness-intersection-prop} that $\wtilde G_0\cap  \R^\cup\la \mathcal B\ra$ is bound to  $\R^\cup\la \mathcal B\ra\cap -\C(X(\kappa))=\R\la \darrow B(\kappa)\ra\cap -\C(X(\kappa))=-\C(\darrow B(\kappa))$, with former equality by Proposition \ref{prop-orReay-BasicProps}.1 and the latter in view of the linear independence of $X(\kappa)$.
We can tile $-\C(\darrow B(\kappa))$ with translates of the fundamental parallelepiped defined using the  linearly independent set $\darrow B(\kappa)$ as the lattice basis  for  $\Lambda_\mathcal B:=\Z\la \darrow B(\kappa)\ra\leq \Lambda$, and then any point $x\in -\C(\darrow B(\kappa))$ will be within distance $M$ of some lattice point from $\Lambda_\mathcal B\cap -\C(\darrow B(\kappa))$, where $M$ is the maximal distance between a point of the fundamental parallelepiped and the set of vertices for that parallelepiped. Thus, since $\wtilde G_0\cap \R^\cup \la\mathcal B\ra$ is bound to $-\C(\darrow B(\kappa))$, it follows that there is a bound $N$ such that any $g\in \wtilde G_0\cap \R^\cup \la \mathcal B\ra$ is within distance $N$ of some lattice point $-x_g\in \Lambda_\mathcal B\cap -\C(\darrow B(\kappa))$. Moreover, as there are only a finite number of possibilities for $\mathcal B$, we can assume the same $N$ suffices for each possible $\mathcal B\subseteq \darrow \mathcal X'\subseteq \mathcal X$.

Any $g\in \wtilde G_0\cap -\C^\cup(\mathcal X')$ has  $-g$ encased by $\mathcal X'$, and thus there is some $\mathcal B\preceq \mathcal X'\subseteq \mathcal X_1\cup\ldots\cup \ldots \mathcal X_s$ which minimally encases $-g$, ensuring that $\mathcal B\subseteq \darrow \mathcal X'$ is a support set with $-g\in \C^\cup(\mathcal B)^\circ\subseteq \R^\cup\la \mathcal B\ra=\R\la \darrow B(\kappa)\ra$, with the final equality in view of  Proposition \ref{prop-orReay-BasicProps}.1. Since $\darrow B(\kappa)\subseteq X(\kappa)$ is linearly independent, it follows that $\darrow B(\kappa)$ is a linear basis for $\R^\cup\la \mathcal B\ra$. For $\x\in \darrow \mathcal B$, let $x=\tilde \x(\iota_\x)=\x(\iota_\x)\in \darrow B(\kappa)$, and let \be\label{onewingduck}\Summ{\x\in \darrow \mathcal B}\alpha_\x x=-g\ee be the representation of $-g\in \R^\cup\la \mathcal B\ra$ as a linear combination of the basis elements $\darrow B(\kappa)$, where $\alpha_\x\in \R$.
Since $\mathcal B$ is a support set, and thus virtual independent,  Proposition \ref{prop-orReay-BasicProps}.4 ensures that any element contained in $\C^\cup(\mathcal B)^\circ$ is a \emph{strictly} positive linear combination of some choice of  representatives from \emph{all} the  half-spaces $\x\in \mathcal B$. Moreover, any representative set $B$ for $\mathcal B$ is linearly independent modulo $\R^\cup\la \partial(\mathcal B)\ra$ by Proposition \ref{prop-orReay-BasicProps}.3.  Thus,
since $-g\in \Lambda\cap \C^\cup(\mathcal B)^\circ$, it follows by considering \eqref{onewingduck} modulo $\R^\cup \la \partial(\mathcal B)\ra$ that  $\alpha_\x>0$ for every $\x\in \mathcal B$ (though not necessarily for all $\x\in \darrow \mathcal B$). By the work above, there is some $x_g\in \Lambda_\mathcal B\cap \C(\darrow B(\kappa))$ with $\mathsf d(-g,x_g)\leq N$. Let $$\Summ{\x\in \darrow \mathcal B}\beta_\x x=x_g$$ be the representation of $x_g\in \C(\darrow B(\kappa))$ as a positive linear combination of the basis elements $\darrow B(\kappa)$, so $\beta_\x\geq 0$ for all $\x\in \darrow \mathcal B$.
Since $x_g\in \Lambda_\mathcal B\cap \C(\darrow B(\kappa))=\C_\Z(\darrow B(\kappa))$, where the latter equality follows in view of the linear independence of $\darrow B(\kappa)$ and $\Lambda_\mathcal B=\Z\la \darrow B(\kappa)\ra$, we have $\beta_\x\in \Z_+$ for all $\x\in \darrow \mathcal B$.

Observe that $\gamma_\y^2\leq \Summ{\x\in \darrow \mathcal B} \gamma_\x^2=\|\Summ{\x\in \darrow \mathcal B}\gamma_\x T(x)\|^2\leq \|T\|^2 \cdot \|\Summ{\x\in \darrow \mathcal B}\gamma_\x x\|^2$ for any $\gamma_\x \in \R$ and $\y\in \darrow \mathcal B$, where $T:\R^\cup \la \mathcal B\ra\rightarrow \R^\cup \la \mathcal B\ra$ is a linear transformation mapping $\darrow B(\kappa)$ to an orthonormal basis and $\|T\|$ is the operator norm of $T$ with respect to the Euclidean $L_2$-norm.
Thus, since $\mathsf d(-g,x_g)\leq N$, it follows that there is some $N'\geq 0$ such that  $|\alpha_\x-\beta_\x|\leq N'$ for all $\x\in \darrow \mathcal B$ and $-g\in \wtilde G_0\cap \C^\cup(\mathcal B)^\circ$ (namely, $N'=\|T\|\cdot N$). For any $\beta_\x\geq \alpha_\x$ with $\x\in \mathcal B$, let $\beta'_\x$ be the greatest integer strictly less than $\alpha_\x$. Note $\beta'_\x\geq 0$ in view of $\alpha_\x>0$. For all other $\x\in \darrow \mathcal B$, let $\beta'_\x=\beta_\x\geq 0$.  Let $$\Summ{\x\in \darrow \mathcal B}\beta'_\x x=x'_g.$$ Since $|\beta_\x-\beta'_\x|\leq |\alpha_\x-\beta_\x|+1\leq N'+1$ for all $\x\in \darrow \mathcal B$, there is some $N''$ such that  $\|x_g-x'_g\|\leq N''$ for all $g\in \wtilde G_0\cap -\C^\cup(\mathcal B)^\circ$ (for instance, in view of the triangle inequality, we could take $N''=(N'+1)\Summ{\x\in \darrow \mathcal B}\|x\|$). Thus, by replacing $N$ with a larger value and using $x'_g$ in place of $x_g$, we obtain the additional conclusion that
\be\label{addconcl} \beta_\x\in \Z_+\quad \mbox{ for all
$\x\in \darrow \mathcal B$} \quad \und\quad \alpha_\x>\beta_\x\geq 0\quad \mbox{ for all $\x\in \mathcal B$},\ee
ensuring that  $$x_g=x'_g\in \C_{\Z}(\darrow B(\kappa))\subseteq \Lambda.$$
Thus, since $g\in \Lambda$ and $\mathsf d(-g,x_g)\leq N$, we have $x_g+y=-g$ with $y\in \Lambda$ and $\|y\|\leq N$. Note
$$y=-g-x_g=\Summ{\x\in \darrow \mathcal B}(\alpha_\x-\beta_\x)x$$ with $\alpha_\x-\beta_\x>0$ for all $\x\in \mathcal B$ in view of \eqref{addconcl}. Thus $y\in \Summ{\x\in \mathcal B}(\x^\circ+\partial(\x))\subseteq \Summ{\x\in \mathcal B}\x\subseteq \C^\cup(\mathcal B)$.
In summary, we now have \be\label{teetime}y\in \Lambda \cap \C^\cup(\mathcal B) \quad\und\quad \|y\|\leq N.\ee

Since any bounded set of lattice points is finite, there are only a finite number of $y$ satisfying \eqref{teetime}.
In consequence, Proposition \ref{prop-VReay-RepBasics}.2 implies that, for any tuple $k=(i_\z)_{\z\in \mathcal X}$ with all $i_\z$ sufficiently large, we have $\C(\darrow B(\kappa))\subseteq \C(\darrow B(k))\subseteq \C(\darrow X'(k))$ and $y\in \C(\darrow B(k))\subseteq \C(\darrow X'(k))$, for every  support set  $\mathcal B\subseteq \darrow \mathcal X'$ and every possible  $y$ satisfying \eqref{teetime}. But then $y\in \C(\darrow X'(k))$ and $x_g\in \C(\darrow B(\kappa))\subseteq \C(\darrow B(k))\subseteq \C(\darrow X'(k))$, for every $g\in \wtilde G_0\cap -\C^\cup(\mathcal X')$, ensuring that $-g=y+x_g\in \C(\darrow X'(k))$ in view of the convexity of $\C(\darrow X'(k))$.
Thus $\wtilde G_0\cap -\C^\cup (\mathcal X')\subseteq \wtilde G_0\cap -\C(\darrow X'(k))$. Since $\C^\cup(\darrow X'(k))=\C^\cup(\darrow \tilde X'(k))\subseteq \C^\cup(\mathcal X')$ for any tuple $k$, the reverse inclusion is trivial, and Item 3 follows, completing the proof.
\end{proof}

\subsection{Series Decompositions and Virtualizations}\label{sec-series-decomp}

Let $\Lambda\leq \R^d$ be a full rank lattice, let $G_0\subseteq \Lambda$ be a finitary subset with $\C(G_0)=\R^d$, and let  $\mathcal R=(\mathcal X_1\cup \{\mathbf v_1\},\ldots,\mathcal X_s\cup \{\mathbf v_s\})$ be a purely virtual Reay system over $G_0$. Since $G_0$ is finitary, any purely virtual Reay system must be anchored. However, if $\mathcal R$ is anchored, then $G_0\subseteq \Lambda$ together with Proposition \ref{prop-VReay-Lattice} ensures that $\tilde \x(i)=\x(i)$  for any $\x\in \bigcup_{i=1}^s\mathcal X_i$ once $i$ is sufficiently large. Suppose $\x\in \mathcal X_j$ with $j\in [1,s]$.  Any $\x(i)$ with $i$ sufficiently large is equal to $\tilde \x(i)$ and is then  an actual lattice point from $G_0$ which is a representative for $\x$. Let $x'\in G_0$ be any element which is a positive  scalar multiple of $u_t$ modulo $\R^\cup\la \bigcup_{i=1}^{j-1}\mathcal X_i)$, where $\vec u_\x=(u_1,\ldots,u_t)$, so $x'\in (\R^\cup\la \bigcup_{i=1}^{j-1}\mathcal X_i\ra+\x)^\circ$. For instance, $x'\in G_0$ could be any representative  for the half-space $\x$, including any $x'=\x(i)$ with $i$  sufficiently large. Define a new half-space $\x'=\R_+x'$ with $\partial(\x')=\{0\}$. Then $\x$ and $\x'$ are both equal modulo
$\R^\cup\la \bigcup_{i=1}^{j-1}\mathcal X_i)$, so  replacing $\x$ by $\x'$ would preserve (OR2) in the definition of an oriented Reay system.
Of course, if $j=1$, then $\x=\x'$ since $\partial(\x)=\{0\}$ in this case, and we have more or less done nothing apart from changing the representative for $\x$.
Suppose, for $j=1,2,\ldots,s$, we replace each $\x\in \bigcup_{i=1}^s\mathcal X_i$ with some half-space $\x'$ with $\partial(\x')=\{0\}$ as just described to result in  $(\mathcal X'_1\cup \{\mathbf v_1\},\ldots,\mathcal X'_s\cup \{\mathbf v_s\})$. Now $(\mathcal X'_1\cup \{\mathbf v_1\})$ is still a virtual Reay system over $G_0$ with $\vec u_{\x}=\vec u_{\x'}$ for all $\x\in \mathcal X_1$, the value of $\vec u_{\mathbf v_1}$ unchanged and $\R^\cup\la \mathcal X_1\ra=\R^\cup\la \mathcal X'_1\ra$.
Since $\mathcal R$ is purely virtual,  $-\vec u_{\mathbf v_2}^\triangleleft$ is fully unbounded or trivial and minimally encased by $\partial(\{\mathbf v_2\})\subseteq \R^\cup \la \mathcal X_1\ra=\R^\cup\la \mathcal X'_1\ra$. Thus, since $G_0$ is finitary and $\mathcal R$ is purely virtual, it follows from Proposition \ref{prop-finitary-basics}.3 applied to $(\mathcal X'_1\cup \{\mathbf v_1\})$ that $-\vec u^\triangleleft_{\mathbf v_2}$ is minimally encased by some subset $\mathcal B_2\subseteq \mathcal X'_1$, allowing us to define a new half-space $\mathbf v'_2$ with $\vec u_{\mathbf v'_2}=\vec u_{\mathbf v_2}$ and $\partial(\{\mathbf v'_2\})=\mathcal B_2$.
This makes $(\mathcal X'_1\cup\{\mathbf v'_1\},\mathcal X'_2\cup \{\mathbf v'_2\})$ into a virtual Reay system, where $\mathbf v'_1:=\mathbf v_1$, with $\partial(\x')=\{0\}$ and $\x'(i)=x'$ constant for $\x'\in \mathcal X'_2$. Moreover, since the value of each $\x'\in \mathcal X'_2$ has not changed modulo $\R^\cup \la \mathcal X_1\ra=\R^\cup \la \mathcal X'_1\ra$, it follows that $\R^\cup\la \mathcal X'_1\cup \mathcal X'_2\ra=\R^\cup\la \mathcal X_1\cup \mathcal X_2\ra$. Iterating this argument, we find that $\mathcal R':=(\mathcal X'_1\cup \{\mathbf v'_1\},\ldots,\mathcal X'_s\cup \{\mathbf v'_s\})$ is a purely Virtual Reay system over $G_0$ with $\R^\cup\la \bigcup_{i=1}^j\mathcal X'_i\ra=\R^\cup\la \bigcup_{i=1}^j\mathcal X_i\ra$ and $\vec u_{\mathbf v'_j}=\vec u_{\mathbf v_j}$ for all $j\in [1,s]$, and with $\partial(\x')=\{0\}$ and $\x'(i)=x'\in G_0$ constant for all $\x'\in \mathcal X'_1\cup\ldots\cup\mathcal X'_s$. In particular, Proposition \ref{prop-finitary-basics}.1 applied to $\mathcal R'$ ensures that any $x'\in G_0$ which lies in the open half-space $(\R^\cup \la \bigcup_{i=1}^{j-1}\mathcal X_i\ra+\x)^\circ$, for $\x\in \mathcal X_j$, must satisfy $x'\in G_0^\diamond$ (assuming $G_0\subseteq \Lambda\subseteq \R^d$ is finitary with $\C(G_0)=\R^d$ and $\mathcal R$ purely virtual).

The above construction of $\mathcal R'$ was done using arbitrary elements $x'\in G_0$ lying in the open half-space $(\R^\cup \la \bigcup_{i=1}^{j-1}\mathcal X_i\ra+\x)^\circ$, for $\x\in \mathcal X_j$, and resulting only in the existence of subsets $\partial(\{\mathbf v'_j\})\subseteq \mathcal X'_1\cup\ldots\cup \mathcal X'_{j-1}$.  Suppose instead we choose each $x'=\x(i_\x)$ for some fixed but sufficiently large $i_\x$. Let $k=(i_\x)_{\x\in \mathcal X}$ be a fixed tuple with all $i_\x$ sufficiently large (as determined below).
For a half-space $\x\in \mathcal X_1\cup\{\mathbf v_1\}\cup\ldots\cup\mathcal X_s\cup \{\mathbf v_s\}$ from $\mathcal R$, let $\mathbf x'\in \mathcal X'_1\cup\{\mathbf v'_1\}\cup\ldots\cup\mathcal X'_s\cup \{\mathbf v'_s\}$ denote the corresponding half-space from $\mathcal R'$, so each $\x'=\R_+\x(i_\x)$ whenever $\x\in \mathcal X_1\cup\ldots\cup \mathcal X_s$, and for a subset $\mathcal X\subseteq \mathcal X_1\cup\{\mathbf v_1\}\cup\ldots\cup\mathcal X_s\cup \{\mathbf v_s\}$, let $\mathcal X'=\{\x':\; \x\in \mathcal X\}\subseteq \mathcal X'_1\cup\{\mathbf v'_1\}\cup\ldots\cup\mathcal X'_s\cup \{\mathbf v'_s\}$.
Assuming the indices $i_\x$ are chosen sufficiently large, Proposition  \ref{prop-VReay-RepBasics}.1 implies  that \be\label{tulet}\partial(\{\mathbf v'_j\})=(\darrow \partial(\{\mathbf v_j\}))'=\darrow \partial(\{v_j\})(k)\quad
\mbox{ for each $j\in [1,s]$}.\ee
Note, for the equality \eqref{tulet} to make sense, we informally identify $\darrow \partial(\{v_j\})(k)$ with the collection of one-dimensional half-spaces generated by the elements from $\darrow \partial(\{v_j\})(k)$, with this convention continued at later points of the discussion.
If $\mathcal A_j$ is the set defined in Proposition \ref{prop-VReay-SupportSet}.1 for $\mathcal R$, then Proposition \ref{prop-finitary-basics}.2 ensures that $\mathcal A_j\setminus \{\mathbf v_j\}\subseteq \mathcal X_1\cup\ldots\cup \mathcal X_j$, in which case Proposition \ref{prop-VReay-Lattice} implies
$(\darrow A_j\setminus \{v_j\})(k)=\wtilde{(\darrow A_j\setminus \{v_j\})}(k)$ so long as all indices  are sufficiently large. In view of Proposition \ref{prop-VReay-SupportSet}.2(e), we have $\R^\cup \la \mathcal A_j\ra=\R^\cup \la \mathcal A_j\setminus \{\mathbf v_j\}\cup\partial(\{\mathbf v_j\})\ra$, while $(A_j\setminus \{v_j\})(k)\subseteq \R^\cup \la \mathcal A_j\setminus \{\mathbf v_j\}\ra$ as the elements of  $(A_j\setminus \{v_j\})(k)$ are representatives for the half-spaces from $\mathcal A_j\setminus \{\mathbf v_j\}$.
This ensures that the set $Y_k$ in Proposition \ref{prop-VReay-SupportSet}.3 (for the element $\mathbf v_j$) is empty, in which case Proposition \ref{prop-VReay-SupportSet}.3  implies that $(\darrow A_j\setminus \{v_j\})(k)$ minimally encases $-\vec u_{\mathbf v_j}$ so long as  all indices are sufficiently large.

In view of \eqref{tulet} and  Proposition \ref{prop-orReay-BasicProps}.1, we have   $$\R^\cup \la \partial(\{\mathbf v'_j\})\ra=\R\la \darrow \partial(\{v_j\})(k)\ra=\R^\cup\la \partial(\{\mathbf v_j\})\ra.$$
Let $\pi_{\mathbf v_j}=\R^d\rightarrow \R^\cup \la \partial(\{\mathbf v_j\})\ra^\bot$ be the orthogonal projection.
By Proposition \ref{prop-orReay-modulo}.1 and Proposition \ref{prop-orReay-BasicProps}.9, $\pi_{\mathbf v_j}$ is injective on $(\mathcal X_1\cup\ldots\cup \mathcal X_s)\setminus \darrow \partial(\{\mathbf v_j\})$ and maps all such elements to non-zero half-spaces.
In consequence, since $(\darrow A_j\setminus \{v_j\})(k)\subseteq \mathcal X'_1\cup\ldots\cup \mathcal X'_s$ minimally encases $-\vec u_{\mathbf v_j}=-\vec u_{\mathbf v'_j}$, it follows from Proposition \ref{prop-orReay-minecase-char}.4 that $\Big((\darrow \mathcal A_j\setminus \{v_j\})(k)\Big)^{\pi_{\mathbf v_j}}=\pi_{\mathbf v_j}\big((\darrow A_j\setminus \darrow v_j)(k)\big)$ minimally encases $-\pi_{\mathbf v_j}(\vec u_{\mathbf v'_j})$. Thus, since $\mathcal A_j\setminus \{\mathbf v_j\}\subseteq \mathcal X$, we conclude that $\pi_{\mathbf v_j}\big((\darrow A_j\setminus \darrow v_j)(k)\big)=\supp_{\pi_{\mathbf v_j}(\mathcal R')}(-\pi_{\mathbf v_j}(\vec u_{\mathbf v'_j}))$, meaning $(\darrow A_j\setminus \darrow v_j)(k)$ is the pull-back of $\supp_{\pi_{\mathbf v_j}(\mathcal R')}(-\pi_{\mathbf v_j}(\vec u_{\mathbf v'_j}))$.
  As a result, it now follows from Propositions \ref{prop-VReay-SupportSet}.1 and \ref{prop-VReay-SupportSet}.2(b) that, letting $\mathcal A'_j$ and $\mathcal A'_{\mathbf v'_j}$ denote the sets given by Proposition \ref{prop-VReay-SupportSet} for $\mathcal R'$, we have
 $\mathcal A'_j=(\darrow A_j\setminus \darrow v_j)(k)\cup \{\mathbf v'_j\}=(\darrow \mathcal A_j\setminus \darrow \partial(\{\mathbf v_j\})\big)'$ and $\mathcal A'_{\mathbf v'_j}=\big(\mathcal A'_j\setminus \{\mathbf v'_j\}\cup \partial(\{\mathbf v'_j\}))^*=((\darrow A_{j}\setminus \{v_j\})(k))^*=(\darrow A_{j}\setminus \{v_j\})(k)$, with the final equality holding since all half-spaces from  $\mathcal X'$ have trivial boundary, and the second in view \eqref{tulet}. In summary,
  \be\nn\partial(\{\mathbf v'_j\})=\big(\darrow \partial(\{\mathbf v_j\})\big)',\quad \mathcal A'_{\mathbf v'_j}=(\darrow \mathcal A_{\mathbf v_j})'\quad\und\quad\mathcal A'_j=\big(\darrow \mathcal A_j\setminus \darrow \partial(\{\mathbf v_j\})\big)' \quad\mbox{ for every $j\in [1,s]$}.\ee This ensures  the virtual Reay system structure associated to each $\mathbf v_j$ is preserved in $\mathbf v'_j$ when passing to the virtual Reay system $\mathcal R'$ (as much as possible given that $\partial(\x')=\{0\}$ for every $\x'\in \mathcal X'_1\cup\ldots\cup \mathcal X'_s$). 

In view of the observations just made, when considering a purely virtual Reay system $\mathcal R=(\mathcal X_1\cup \{\mathbf v_1\},\ldots,\mathcal X_s\cup \{\mathbf v_s\})$ over a finitary subset $G_0\subseteq \Lambda$ with $\C(G_0)=\R^d$, we can often restrict attention to when $\partial(\x)=\{0\}$ for all $\x\in \mathcal X_1\cup\ldots\cup X_s$, as if this fails for $\mathcal R$, then another virtual Reay system $\mathcal R'$ over $G_0$ can be constructed with this property as described above with the values of $\partial(\{\mathbf v_j\})$, $\mathcal A_j$ and $\mathcal A_{\mathbf v_j}$ minimally affected.

Let $\Lambda\leq \R^d$ be a full rank lattice and let $G_0\subseteq \Lambda$ be a finitary subset with $\C(G_0)=\R^d$. We now define the followings sets:
\begin{align*}&\mathfrak X(G_0)=\{\bigcup_{i=1}^s\mathcal X_i:\;\mbox{there is a purely virtual Reay sytem $(\mathcal X_1\cup \{\mathbf v_1\},\ldots,\mathcal X_s\cup \{\mathbf v_s\})$ over $G_0$}\}, \\ & X(G_0)=\{X:\;X\subseteq G_0 \mbox{ is a set of representatives for the half-spaces from some $\mathcal X\in \mathfrak X(G_0)$}\}.\end{align*}
Since the empty tuple is by default a purely virtual Reay system, we always have $\emptyset\in X(G_0)$.
If $X\in X(G_0)$, then there is a purely virtual Reay System $(\mathcal X_1\cup\{\mathbf v_1\},\ldots,\mathcal X_s\cup \{\mathbf v_s\})$ over $G_0$ with $X\subseteq G_0$ being a set of representative for the half-spaces from $\bigcup_{i=1}^s\mathcal X_i$.
Per the discussion above, replacing each $x\in X$ with the half-space $\R_+ x$  results in a purely virtual Reay system $\mathcal R'=(\mathcal X'_1\cup \{\mathbf v'_1\},\ldots,\mathcal X'_s\cup\{\mathbf v'_s\})$ over $G_0$ such that $X$ is a set of representatives for the half-spaces from  $\bigcup_{i=1}^s\mathcal X'_i$ and $\partial(\x')=\{0\}$ with $\x'(i_\x)=x\in G_0$ for all $\x'\in \bigcup_{i=1}^s\mathcal X'_i$ and $i_\x\geq 1$. Thus it can always be assumed that $X\in X(G_0)$ came from a virtual Reay system having these properties. Prior discussion  ensures that $X\subseteq G_0^\diamond$ for any $X\in X(G_0)$, while any $X\in X(G_0)$ is a linearly independent subset of $G_0\subseteq \Lambda$ by Proposition \ref{prop-reay-basis-properties}.1, thus generating a sublattice of $\Lambda$. Let
\begin{align*}
&\mathfrak P_\Z(G_0)=\{\Z\la X\ra:\; X\in X(G_0)\} \und\\
&\mathfrak P_{\R}(G_0)=\{\R^\cup \la \mathcal X\ra:\; \mathcal X\in \mathfrak X(G_0)\}=\{\R\la X\ra:\; X\in X(G_0)\}=\{\R\la \Lambda'\ra:\;\Lambda'\in \mathfrak P_\Z(G_0)\},
\end{align*}
with the second equality for  $\mathfrak P_\R(G_0)$ by Proposition \ref{prop-orReay-BasicProps}.1. Note $\mathfrak P_\Z(G_0)$ consists of sublattices of $\Lambda$ as just discussed.
The linearly independent set $X\in X(G_0)$ can be recovered from $\C_{\Z}(X)$ by considering each ray defined by a vertex  in the polyhedron $\C(X)\cap B_1(0)$ lying on the unit sphere, and then taking the minimal nonzero element of this ray contained in $\C_\Z(X)$. Thus $\C_{\Z}(X)=\C_{\Z}(Y)$ implies $X=Y$, allowing us to  define a partial order $\preceq_\Z$ on $X(G_0)$ by declaring $X\preceq_\Z Y$ when $\C_{\Z}(X)\subseteq \C_{\Z}(Y)$.
 Since each $X\in X(G_0)$ is linearly independent, we have \be\label{cone-integer}\C_{\Z}(X)=\C(X)\cap \Z\la X\ra\quad\mbox{ for $X\in X(G_0)$}.\ee
 Thus $$X\preceq_\Z Y\quad\mbox{if and only if }\quad\C(X)\subseteq \C(Y)\quad\und\quad\Z\la X\ra\leq \Z\la Y\ra.$$ Indeed, if $\C_{\Z}(X)\subseteq \C_{\Z}(Y)$, then $x\in \C_{\Z}(Y)\subseteq \Z\la Y\ra\cap \C(Y)$ for each $x\in X$, ensuring that $\C(X)\subseteq \C(Y)$ and $\Z\la X\ra\leq \Z\la Y\ra$, while if $\C(X)\subseteq \C(Y)$ and $\Z\la X\ra\leq \Z\la Y\ra$, then $\C_{\Z}(X)=\C(X)\cap \Z\la X\ra\subseteq \C(Y)\cap \Z\la Y\ra=\C_{\Z}(Y)$.

\begin{definition}For $X\in X(G_0)$, there is a purely virtual Reay system $\mathcal R=(\mathcal X_1\cup \{\mathbf v_1\},\ldots,\mathcal X_s\cup \{\mathbf v_s\})$ over $G_0$ and ordered partition $X=X_1\cup\ldots\cup X_s$ such that each $X_j$ is a set of representatives for the half-spaces from $\mathcal X_j$, for $j\in [1,s]$. We call any such $\mathcal R$ a \textbf{realization} of $X$ and the ordered partition $X=X_1\cup\ldots\cup X_s$ a \textbf{series decomposition} of $X$. The set $\mathcal X=\bigcup_{i=1}^s\mathcal X_i$ is called a \textbf{virtualization} of $X$.
For $\mathcal X\in \mathfrak X(G_0)$, there is a purely virtual Reay system $\mathcal R=(\mathcal X_1\cup \{\mathbf v_1\},\ldots,\mathcal X_s\cup \{\mathbf v_s\})$ over $G_0$. We also call $\mathcal X=\bigcup_{i=1}^s\mathcal X_i$ a \textbf{series decomposition} of $\mathcal X$ and $\mathcal R$ a \textbf{realization} of $\mathcal X$.
\end{definition}

To explain the name, it may be helpful to view a series decomposition as an ascending chain $\emptyset\subset X_1\subset (X_1\cup X_2)\subset\ldots\subset (X_1\cup\ldots\cup X_s)=X$ or $\emptyset\subset \mathcal X_1\subset (\mathcal X_1\cup \mathcal X_2)\subset\ldots\subset (\mathcal X_1\cup\ldots\cup \mathcal X_s)=\mathcal X$
If $X=X_1\cup\ldots\cup X_s$ is a series decomposition of $X\in X(G_0)$, then $Y_j:=X_1\cup\ldots\cup X_j\subseteq X$ is a subset of $X$ with $Y_j\in X(G_0)$ for any $j\in [1,s]$.
Likewise, if $\mathcal X=\bigcup_{i=1}^\cup \mathcal X_i$ is a series decomposition of $\mathcal X\in \mathfrak X(G_0)$, then every $\mathcal Y_j:=\mathcal X_1\cup\ldots\cup \mathcal X_j\subseteq \mathcal X$ is a subset of $\mathcal X$ with $\mathcal Y_j\in \mathfrak X(G_0)$ for any $j\in [1,s]$.
The converse to both these statements regarding $X(G_0)$ and $\mathfrak X(G_0)$ is also true in the following strong sense.

\begin{lemma}\label{lemma-X0-fact}
Let $\Lambda\leq \R^d$ be a full rank lattice and let $G_0\subseteq \Lambda$ be a finitary subset with $\C(G_0)=\R^d$. Suppose $\mathcal X,\,\mathcal Y\in \mathfrak X(G_0)$ with $\mathcal Y\subseteq \mathcal X$ and that $\mathcal Y=\mathcal Y_1\cup\ldots\cup \mathcal Y_t$ is a series decomposition of $\mathcal Y$. Let  $\pi:\R^d\rightarrow \R^\cup \la \mathcal Y\ra^\bot$ be the orthogonal projection.
\begin{itemize}
\item[1.] There is a series decomposition $\mathcal X=\mathcal X_1\cup\ldots\cup \mathcal X_s$ with $s\geq t$ and $\mathcal X_i=\mathcal Y_i$ for $i\in [1,t]$.
\item[2.] $\mathcal X^\pi=\pi(\mathcal X\setminus \mathcal Y)\in \mathfrak X(\pi(G_0))$.
\item[3.] If $\mathcal X\setminus \mathcal Y=\mathcal X'_1\cup \ldots\cup \mathcal X'_r$ with $\pi(\mathcal X\setminus \mathcal Y)=\pi(\mathcal X'_1)\cup \ldots\cup \pi(\mathcal X'_r)$ a series decomposition, 
     then $\mathcal X=\mathcal Y_1\cup\ldots\cup \mathcal Y_t\cup \mathcal X'_1\cup \ldots\cup\mathcal X'_r$ is a series decomposition.
\end{itemize}
\end{lemma}

\begin{proof}
1. Let $\mathcal R=(\mathcal X_1\cup\{\mathbf v_1\},\ldots,\mathcal X_s\cup \{\mathbf v_s\})$ be a realization of $\mathcal X\in \mathfrak X(G_0)$, so $\mathcal X=\mathcal X_1\cup \ldots\cup \mathcal X_s$,  and let  $\mathcal R_\mathcal Y=(\mathcal Y_1\cup \{\mathbf w_1\},\ldots,\mathcal Y_t\cup \{\mathbf w_t\})$ be a realization of $\mathcal Y\subseteq \mathcal X$ associated to the series decomposition  $\mathcal Y=\mathcal Y_1\cup\ldots\cup \mathcal Y_t$. Thus $\darrow \mathcal Y=\mathcal Y$ when considered as a subset of half-spaces from $\mathcal R_\mathcal Y$, and thus also (by Proposition \ref{prop-orientedReay-halfspace-uniquelydetermines}) when considered as a subset of half-spaces from  $\mathcal R$. In view of Proposition \ref{prop-VReay-Lattice}, we have $\R^\cup \la \mathcal Y\ra=\R \la \darrow Y(k)\ra=\R\la Y(k)\ra$ for any tuple $k=(i_\y)_{\y\in \mathcal Y}$ with all $i_\y$ sufficiently large. By Proposition \ref{prop-VReay-modulo}, $\pi(\mathcal R)=(\mathcal X_j^\pi\cup \{\pi(\mathbf v_j)\})_{j\in J}$ is a Virtual Reay system over $\pi(G_0)$. Since $\mathcal R$ is purely virtual,
 it follows from Proposition \ref{prop-VReay-modulo} that $\pi(\mathcal R)$ is also purely virtual, and Proposition \ref{prop-VReay-modulo}  ensures that \be\label{wheed}\vec u_{\pi(\x)}=\pi(\vec u_\x),\quad \pi(\vec u_\x)^\triangleleft=\pi(\vec u_\x^\triangleleft),\quad\und\quad \partial(\{\pi(\x)\})=\partial(\{\x\})^\pi\ee for all $\x\in \bigcup_{j\in J}(\mathcal X_j\cup \{\mathbf v_j\})$ with  $\pi(\x)\neq 0$.

 Let $j_1<j_2<\ldots<j_r$ be the indices from $J$ and let $\pi(\mathcal R)=(\mathcal C_1\cup \{\mathbf c_1\},\ldots,\mathcal C_r\cup \{\mathbf c_r\})$, so each  $\mathcal C_i=\mathcal X_{j_i}^\pi$ and each $\mathbf c_i=\pi(\mathbf v_{j_i})$ for $i\in [1,r]$. For each $i\in [1,r]$, let $\mathcal B_i=\pi^{-1}(\mathcal C_i)\subseteq \mathcal X_{j_i}$.  Thus Proposition \ref{prop-orReay-BasicProps}.9 and $\darrow \mathcal Y=\mathcal Y$ give  \be\label{YX-equzal}\mathcal Y\cup \bigcup_{i=1}^r\mathcal B_i=\mathcal X,\ee with the union disjoint. Each $\mathbf c_i=\pi(\mathbf v_{j_i})$, for $i\in [1,r]$, is defined by the limit $\vec u_{\pi(\mathbf v_{j_i})}=\pi(\vec u_{\mathbf v_{j_i}})$ and satisfies $\pi(\vec u_{\mathbf v_{j_i}})^\triangleleft=\pi(\vec u_{\mathbf v_{j_i}}^\triangleleft)$ by \eqref{wheed}.
 Moreover, each $-\vec u_{\mathbf v_{j_i}}^\triangleleft$, for $i\in [1,r]$, is minimally encased by $\partial(\{\mathbf v_{j_i}\})\subseteq \bigcup_{k=1}^{j_i-1}\mathcal X_k$, and is thus encased by $\mathcal Y\cup \mathcal B_{\mathbf c_i}$, where $\mathcal B_{\mathbf c_i}:=\partial(\{\mathbf v_{j_i}\})\setminus \mathcal Y$ with $\mathcal B_{\mathbf c_i}\subseteq \mathcal B_1\cup \ldots\cup \mathcal B_{i-1}$. Also, $\pi(\mathcal B_{\mathbf c_i})=\mathcal B_{\mathbf c_i}^\pi=\partial(\{\mathbf v_{j_i}\})^\pi=\partial(\{\pi(\mathbf v_{j_i})\})=\partial(\{\mathbf c_i\})$, where the first equality follows since $\mathcal B_{\mathbf c_i}\subseteq \mathcal B_1\cup \ldots\cup \mathcal B_{i-1}$ with $\pi(\x)\neq \{0\}$ for all $\x\in \mathcal B_1\cup\ldots\cup \mathcal B_r$, the third by \eqref{wheed}, and the other two in view of the definitions of  $\mathcal B_{\mathbf c_i}$ and $\mathbf c_i$.
 But now we have the needed hypotheses to apply Lemma \ref{lemma-VReay-extension} (using $\mathcal R_\mathcal Y$ for  $\mathcal R_\mathcal A$, and $\pi(\mathcal R)$ for $\mathcal R_\mathcal C$) to thereby conclude that $$\mathcal R'=(\mathcal Y_1\cup \{\mathbf w_1\},\ldots,\mathcal Y_t\cup \{\mathbf w_t\},\mathcal B_1\cup \{\mathbf b_1\},\ldots,\mathcal B_r\cup \{\mathbf b_r\})$$ is a virtual Reay system in $G_0$ with each $\mathbf b_i$, for $i\in [1,r]$, defined by the limit $\vec u_{\mathbf v_{j_i}}$. By \eqref{YX-equzal}, we have $\mathcal Y_1\cup\ldots\mathcal Y_t\cup \mathcal B_1\cup \ldots\cup \mathcal B_r=\mathcal X$.
 Now $\mathcal R_\mathcal Y$ and $\pi(\mathcal R)$ are both purely virtual, the former as it is a realization of a set from $\mathfrak X(G_0)$, and the latter as observed above \eqref{wheed}, whence Lemma \ref{lemma-VReay-extension} implies that $\mathcal R'$ is purely virtual as well, completing Item 1.


2. By Item 1, there exists a realization $\mathcal R=(\mathcal X_1\cup\{\mathbf v_1\},\ldots,\mathcal X_s\cup \{\mathbf v_s\})$ of $\mathcal X=\mathcal X_1\cup\ldots\cup \mathcal X_s$ with $\mathcal Y_i=\mathcal X_i$ for $i\in [1,t]$. Proposition \ref{prop-VReay-modulo} implies that $\pi(\mathcal R)=(\pi(\mathcal X_{t+1})\cup \{\pi(\mathbf v_{t+1})\},\ldots,\pi(\mathcal X_s)\cup \{\pi(\mathbf v_s)\})$ is a realization of $\mathcal X^\pi=\pi(\mathcal X)\setminus \{\{0\}\}=\pi(\mathcal X\setminus\mathcal Y)$. Thus $\pi(\mathcal X\setminus \mathcal Y)\in \mathfrak X(\pi(G_0))$, establishing Item 2. Note $\mathfrak X(\pi(G_0))$ is well-defined by Proposition \ref{prop-finitary-Modulo-Inheritence}.

3. Suppose $\mathcal X\setminus \mathcal Y=\mathcal X'_1\cup \ldots\cup \mathcal X'_r$ with $\pi(\mathcal X\setminus \mathcal Y)=\pi(\mathcal X'_1)\cup \ldots\cup \pi(\mathcal X'_r)$ a series decomposition.  Let $\mathcal R=(\mathcal X_1\cup \{\mathbf v_1\},\ldots,\mathcal X_s\cup \{\mathbf v_s\})$ be a realization of $\mathcal X$.
We proceed inductively to show $\mathcal Y_1\cup\ldots\cup \mathcal Y_t\cup \mathcal X'_1\cup\ldots\cup \mathcal \mathcal X'_{r'}$ is a series decomposition, for $r'=1,2,\ldots,r$.  We begin with the case $r'=1$. Let $\mathcal R_C=(\pi(\mathcal X'_1)\cup \{\mathbf c_1\})$ be a realization of $\pi(\mathcal X'_1)\in \mathfrak X(G_0)$.  Note each $\pi(\x)$ with $\x\in \mathcal X'_1$ has trivial boundary, so $\vec u_{\pi(\x)}$ is a tuple consisting of one element, which is a representative for the half-space $\pi(\x)$, and $\partial(\x)\subseteq \ker \pi=\R^\cup \la\mathcal Y\ra$, in turn ensuring that $\pi(\vec u_\x)$ is a tuple consisting of one element which is also a representative for the half-space $\pi(\x)$ (by (VR1) for $\x$ when considered part of the realization $\mathcal R$). Thus $\pi(\vec u_\x)=\vec u_{\pi(\x)}$ for $\x\in \mathcal X_1$.
We aim to use Lemma \ref{lemma-VReay-extension} with the realizations $\mathcal R$ and $\mathcal R_C$ as given above, $\mathcal A=\mathcal Y$ and $\mathcal R_A=(\mathcal Y_1\cup \{\mathbf w_1\},\ldots,\mathcal Y_t\cup \{\mathbf w_t\})$ a realization of $\mathcal Y\in \mathfrak X(G_0)$. By Proposition \ref{prop-infinite-limits-proj}.2, we have $\vec u_{\mathbf c_1}=\pi(\vec u)$ for some  limit $\vec u$ of an asymptotically filtered sequence of terms from $G_0$. Moreover,  $\pi(\vec u^\triangleleft)=\pi(\vec u)^\triangleleft$. Let $\vec u=(u_1,\ldots,u_\ell)$. Since $\partial(\{\mathbf c_1\})=\emptyset$, we have $u_1,\ldots,u_{\ell-1}\in \ker \pi=\R^\cup \la \mathcal Y\ra$. Since $\vec u_{\mathbf c_1}=\pi(\vec u)$ is fully unbounded (as $\mathcal R_\mathcal C$ is purely virtual) with $\pi(\vec u^\triangleleft)=\pi(\vec u)^\triangleleft$, it follows by Proposition \ref{prop-infinite-limits-proj}.1 that $\vec u$ is also fully unbounded. Thus  $\vec u^\triangleleft$ is either trivial or fully unbounded, and so  $-\vec u^\triangleleft$ is encased by $\mathcal Y$ in view of Proposition \ref{prop-finitary-basics}.3
applied to the realization $\mathcal R_A$.
 But now we can apply Lemma \ref{lemma-VReay-extension} to conclude $\mathcal Y_1\cup\ldots\cup \mathcal Y_t\cup \mathcal X'_1$ is a series decomposition (the virtual Reay system given by Lemma \ref{lemma-VReay-extension} will be purely virtual since both $\mathcal R_C$ and $\mathcal R_A$ are purely virtual  by assumption). This completes the base case $r'=1$.
 However, for $r'>1$, the induction hypothesis gives that $\mathcal Y':=\mathcal Y_1\cup \ldots\cup \mathcal Y_t\cup \mathcal X'_1\cup \ldots\cup \mathcal X'_{r'-1}$ is a series decomposition of $\mathcal Y'\subseteq \mathcal X$, showing $\mathcal Y'\in \mathfrak X(G_0)$. Let
  $\pi':\R^d\rightarrow \R^\cup\la \mathcal Y'\ra^\bot$ be the orthogonal projection. Since $\pi(\mathcal X\setminus \mathcal Y)=\pi(\mathcal X'_1)\cup \ldots\cup \pi(\mathcal X'_r)$ is a series decomposition, Proposition \ref{prop-VReay-modulo} and Proposition \ref{prop-orReay-modulo}.1 imply that  $\mathcal X^{\pi'}=\pi'(\mathcal X\setminus \mathcal Y')=\pi'(\mathcal X_{r'})$ is a series decomposition, and applying the base case with $\mathcal Y'$ in place of $\mathcal Y$ completes the induction and proof.
\end{proof}

If $X\in X(G_0)$, then $X$ has a series decomposition $X=X_1\cup \ldots\cup X_s$ and associated realization $\mathcal R=(\mathcal Y_1\cup \{\mathbf v_1\},\ldots,\mathcal Y_s\cup\{\mathbf v_s\})$. If we replace each $x\in X$ with the one-dimensional half-space $\x=\R_+x$ and let $\mathcal X$ and $\mathcal X_j$ for $j\in [1,s]$ be the resulting set of half-spaces represented by  $X$ and $X_j$, respectively, then per the discussion regarding the construction of $\mathcal R'$ at the beginning of Section \ref{sec-series-decomp}, it follows that $\mathcal X\in \mathfrak X(G_0)$ with $\mathcal X=\mathcal X_1\cup\ldots\cup \mathcal X_s$ a series decomposition and $\mathcal X$ a virtualization of $X$. This allows us to translate statements involving $\mathfrak X(G_0)$ to ones regarding $X(G_0)$ by replacing any series decomposition $X=X_1\cup\ldots\cup X_s$ with the series decomposition $\mathcal X=\mathcal X_1\cup\ldots\cup\mathcal X_s$, applying the appropriate result regarding $\mathfrak X(G_0)$ to $\mathcal X$, and then returning to $X(G_0)$ by using that each $X_i$ is a set of representatives for the one-dimensional half-spaces from $\mathcal X_i$. We do this out concretely as an example in the next lemma.

\begin{lemma}\label{lemma-X0-fact-repver}
Let $\Lambda\leq \R^d$ be a full rank lattice and let $G_0\subseteq \Lambda$ be a finitary subset with $\C(G_0)=\R^d$. Suppose $X,\,Y\in X(G_0)$ with $Y\subseteq X$ and that $Y=Y_1\cup\ldots\cup Y_t$ is a series decomposition of $Y$. Let $\pi:\R^d\rightarrow \R\la Y\ra^\bot$ be the orthogonal projection.
\begin{itemize}
\item[1.] There is a series decomposition $X=X_1\cup\ldots\cup  X_s$ with $s\geq t$ and $X_i=Y_i$ for $i\in [1,t]$.
\item[2.] $\pi(X\setminus  Y)=\pi(X)\setminus \{0\}\in  X(\pi(G_0))$.
\item[3.] If $ X\setminus Y= X'_1\cup \ldots\cup  X'_r$ with $\pi(X\setminus  Y)=\pi( X'_1)\cup \ldots\cup \pi(X'_r)$ a series decomposition, then $X= Y_1\cup\ldots\cup Y_t\cup X'_1\cup \ldots\cup X'_r$ is a series decomposition.
\end{itemize}
\end{lemma}

\begin{proof}
Replace each element $x\in X$ with the half-space $\x=\R_+ x$ to define a new set $\mathcal X$. The subset $Y\subseteq X$ then defines a new set $\mathcal Y\subseteq  \mathcal X$ .  Per the discussion regarding the construction of $\mathcal R'$ from the beginning of Section \ref{sec-series-decomp}, we have $\mathcal X,\,\mathcal Y\in \mathfrak X(G_0)$ with $X$ being a set of representatives for $\mathcal X$. Indeed, each subset $Y_i\subseteq Y$ defines a new set $\mathcal Y_i$ such that $\mathcal Y=\mathcal Y_1\cup\ldots\cup \mathcal Y_t$ is a series decomposition  of $\mathcal Y$ with each $Y_i$ a set of representatives for $\mathcal Y_i$. Applying Lemma \ref{lemma-X0-fact}.1 to $\mathcal Y\subseteq \mathcal X$ using the series decomposition  $\mathcal Y=\mathcal Y_1\cup\ldots\cup\mathcal Y_t$, we find a series decomposition  $\mathcal X=\mathcal X_1\cup \ldots\cup \mathcal X_s$ with $s\geq t$ and $\mathcal X_i=\mathcal Y_i$ for $i\in [1,t]$. Since $X$ is a set of representatives for the half-spaces from  $\mathcal X=\mathcal X_1\cup \ldots\cup \mathcal X_s$, we have a partition   $X=X_1\cup \ldots\cup X_s$ with each $X_j$ a set of representatives for the half-spaces from $\mathcal X_j$, for $j\in [1,s]$. Thus $X=X_1\cup \ldots\cup X_s$ is a series decomposition of $X$. Moreover, since $\mathcal Y_j=\mathcal X_j$ for $j\in [1,t]$ with $Y_j\subseteq Y\subseteq X$ the set of representatives for $\mathcal Y_j$, we have $Y_j=X_j$ for $j\in [1,t]$.
This establishes Item 1. Let $\mathcal R=(\mathcal X_1\cup \{\mathbf v_1\},\ldots,\mathcal X_s\cup \{\mathbf v_s\})$ be a realization associated to the series decomposition $\mathcal X=\mathcal X_1\cup\ldots\cup \mathcal X_s$. Then the realization  $\pi(\mathcal R)=(\pi(\mathcal X_{t+1})\cup \{\pi(\mathbf v_{t+1})\},\ldots,\pi(\mathcal X_s)\cup \{\pi(\mathbf v_s)\})$ shows $\pi(\mathcal X\setminus \mathcal Y)\in \mathfrak X(G_0)$ with $\pi(X\setminus Y)$ as a set of representatives, which shows  $\pi(X\setminus Y)\in X(G_0)$, giving Item 2. For item 3, each set $X'_j\subseteq X$ defines a set of half-spaces $\mathcal X'_j$ by replacing the elements of $X'_j$ with the one-dimensional rays they define. As argued in  Item 1, $\pi(\mathcal X\setminus Y)=\pi(\mathcal X'_1)\cup \ldots\cup \pi(\mathcal X'_r)$ is a series decomposition. Applying Lemma \ref{lemma-X0-fact}.3 using this decomposition, we find that  $\mathcal X=\mathcal Y_1\cup\ldots\cup \mathcal Y_t\cup \mathcal X'_1\cup \ldots\cup \mathcal X'_r$ is a series decomposition, and since each $Y_i$ is a set of representatives for $\mathcal Y_i$ and each $X'_i$ is a set of representatives for $\mathcal X_i$, Item 3 follows.
\end{proof}

Let $G_0\subseteq \Lambda\subseteq \R^d$ be a finitary set with $\C(G_0)=\R^d$ and suppose $\mathcal X\in \mathfrak X(G_0)$. Let $\mathcal X=\mathcal X_1\cup\ldots\cup \mathcal X_s$ be a series decomposition of $\mathcal X$. This corresponds to a chain $\emptyset\subset \mathcal X_1\subset (\mathcal X_1\cup \mathcal X_2)\subset \ldots \subset (\mathcal X_1\cup \ldots\cup \mathcal X_s)=\mathcal X$.
A \textbf{refinement} of the series decomposition $\mathcal X=\mathcal X_1\cup\ldots\cup \mathcal X_s$ is a series decomposition $\mathcal X=\mathcal X'_1\cup \ldots\cup \mathcal X'_{r}$ such that each $\mathcal X_j=\mathcal X'_{t_{j-1}+1}\cup \ldots\cup \mathcal X'_{t_j}$ for some $0=t_0<t_1<\ldots<t_s=r$.
Equivalently, this means that the subsets $(\mathcal X_1\cup \ldots\cup \mathcal X_j)$ for $j\in [0,s]$ each occur in the chain $\emptyset \subset \mathcal X'_1\subset (\mathcal X'_1\cup \mathcal X'_2)\subset \ldots\subset (\mathcal X_1'\cup \ldots\cup \mathcal X'_r)=\mathcal X'=\mathcal X$. Such a refinement is proper if $r>s$, and the series decomposition $\mathcal X=\mathcal X_1\cup\ldots\cup \mathcal X_s$ is called a \textbf{maximal} series decomposition  if it has no proper refinements.
We say that $\mathcal X\in \mathfrak X(G_0)$ is \textbf{irreducible} if no proper, nonempty subset of $\mathcal X$ lies in $\mathfrak X(G_0)$.
Note this ensures that the only series decomposition of $\mathcal X$ is $\mathcal X=\mathcal X$, since if $\mathcal X=\mathcal X_1\cup \mathcal X_2$ is a series decomposition, then $\mathcal X_1\in \mathfrak X(G_0)$. Moreover, Lemma \ref{lemma-X0-fact}.1 ensures the converse also holds, meaning $\mathcal X\in \mathfrak X(G_0)$ is irreducible if and only if $\mathcal X=\mathcal X$ is the only series decomposition of $\mathcal X$. The above terms are used with the analogous definitions  for $X\in X(G_0)$ as well.
Using Lemma \ref{lemma-X0-fact}, we now characterize maximal series decompositions.

\begin{proposition}\label{prop-finitary-MaxDecompChar}Let $\Lambda\leq \R^d$ be a full rank lattice, where $d\geq 0$, let  $G_0\subseteq \Lambda$ be a finitary set with $\C(G_0)=\R^d$, and let $\mathcal X=\mathcal X_1\cup \ldots\cup \mathcal X_s$ be a series decomposition of $\mathcal X\in \mathfrak X(G_0)$. Then $\mathcal X=\mathcal X_1\cup \ldots\cup \mathcal X_s$ is maximal if and only if each $\pi_{j-1}(\mathcal X_j)$ is irreducible for $j\in [1,s]$, where $\pi_{j-1}:\R^d\rightarrow \R^\cup \la \mathcal X_1\cup \ldots\cup \mathcal X_{j-1}\ra^\bot$ is the orthogonal projection.
\end{proposition}

\begin{proof}
Suppose the series decomposition $\mathcal X=\mathcal X_1\cup \ldots\cup \mathcal X_s$ is maximal but by contradiction there is some $j\in [1,s]$ such that $\pi_{j-1}(\mathcal X_j)$ is not irreducible. Then there is some series decomposition $\pi_{j-1}(\mathcal X_j)=\pi_{j-1}(\mathcal Y_1)\cup \ldots\cup \pi_{j-1}(\mathcal Y_t)$, where $\mathcal X_j=\mathcal Y_1\cup \ldots\cup \mathcal Y_t$ with $t\geq 2$. By Proposition \ref{prop-VReay-modulo} (applied to a realization $\mathcal R$ of $\mathcal X=\mathcal X_1\cup \ldots\cup \mathcal X_s$), it follows that $\pi_{j-1}(\mathcal X_{j})\cup \ldots\cup \pi_{j-1}(\mathcal X_s)$ and $\pi_j(\mathcal X_{j+1})\cup \ldots\cup \pi_j(\mathcal X_s)$ are  series decompositions.
Using these decomposition and applying Lemma \ref{lemma-X0-fact}.3 in $\pi_{j-1}(G_0)\subseteq \pi_{j-1}(\Lambda)$ taking $\pi_{j-1}(\mathcal X)=\pi_{j-1}(\mathcal X_j)\cup\ldots\cup \pi_{j-1}(\mathcal X_s)$ to be $\mathcal X$, taking $\pi_{j-1}(\mathcal X_j)$ to be $\mathcal Y$, and taking $\pi$ to be $\pi_j$, we conclude that $\pi_{j-1}(\mathcal X\setminus (\mathcal X_1\cup\ldots\cup \mathcal X_{j-1}))=\pi_{j-1}(\mathcal Y_1)\cup \ldots\cup \pi_{j-1}(\mathcal Y_t)\cup \pi_{j-1}(\mathcal X_{j+1})\cup\ldots\cup \pi_{j-1}(\mathcal X_s)$ is a series decomposition, and now a second application of Lemma \ref{lemma-X0-fact}.3 yields that $\mathcal X=\mathcal X_1\cup\ldots\cup \mathcal X_{j-1}\cup \mathcal Y_1\cup \ldots\cup \mathcal Y_t\cup \mathcal X_{j+1}\cup \ldots\cup \mathcal X_s$ is a series decomposition, contradicting that $\mathcal X=\mathcal X_1\cup \ldots\cup \mathcal X_s$ was maximal in view of $t\geq 2$.

Conversely, now suppose each $\pi_{j-1}(\mathcal X_j)$ is irreducible but by contradiction  $\mathcal X=\mathcal X_1\cup \ldots\cup \mathcal X_s$ is not maximal.
Then there exists a proper refinement $\mathcal X=\mathcal Y_1\cup \ldots\cup \mathcal Y_r$.
Let $j\in [1,s]$ be the minimal index such that $\mathcal X_j\neq \mathcal Y_j$, so $\mathcal X_i=\mathcal Y_i$ for $i\in [1,j-1]$ and $\mathcal X_j=\mathcal Y_j\cup \ldots\cup \mathcal Y_{j+t}$ with $t\geq 1$.
Applying Proposition \ref{prop-VReay-modulo} to a realization of the series decomposition  $\mathcal Y_1\cup \ldots\cup \mathcal Y_{j+t}$, we conclude that $\pi_{j-1}(\mathcal X_j)=\pi_{j-1}(\mathcal Y_j)\cup\ldots\cup \pi_{j-1}(\mathcal Y_{j+t})$ is a series decomposition, contradicting that $\pi_{j-1}(\mathcal X_j)$ is irreducible  and $t\geq 1$.
\end{proof}

\begin{definition}
Let $\Lambda\leq \R^d$ be a full rank lattice, where $d\geq 0$, let $G_0\subseteq \Lambda$,  let $\mathcal R=(\mathcal X_1\cup \{\mathbf v_1\},\ldots,\mathcal X_s\cup \{\mathbf v_s\})$ be an anchored virtual Reay system over $G_0$, and let $\pi_\x:\R^d\rightarrow \R^\cup\la \partial(\{\x\})\ra^\bot$ be the orthogonal projection, for $\x\in \mathcal X= \mathcal X_1\cup\ldots\cup\mathcal X_s$.  We say that $\mathcal R$ is \textbf{$\Delta$-pure} if $\x(i_\x)=\tilde \x(i_\x)$ with $\pi_\x(\x(i_\x))$  constant  and $\Z\la X(k)\ra=\Delta$, for all tuples $k=(i_\x)_{\x\in \mathcal X}$ and $\x\in \mathcal X$.
\end{definition}


Let $\Lambda\leq \R^d$ be a full rank lattice, where $d\geq 0$, let $G_0\subseteq \Lambda$,  let $\mathcal R=(\mathcal X_1\cup \{\mathbf v_1\},\ldots,\mathcal X_s\cup \{\mathbf v_s\})$ be an anchored virtual Reay system over $G_0$. Set $\mathcal X=\mathcal X_1\cup\ldots\cup \mathcal X_s$.
Then Proposition \ref{prop-VReay-Lattice} ensures that $\x(i_\x)=\tilde \x(i_\x)$ with $\pi_\x(\x(i_\x))$ constant, for any $\x\in \mathcal X$, once $i_\x$ is sufficiently large, so  the first two conditions in the definition of $\Delta$-pure can always be achieved by discarding the first few terms in each representative sequence $\{\x(i_\x)\}_{i=1}^\infty$. Assume this is the case. For $j\in [0,s]$, let $$\mathcal E_{j}=\R\la \mathcal X_1\cup\ldots\cup \mathcal X_{j}\ra=\R\la (X_1\cup\ldots\cup X_{j})(k)\ra$$ and let  $\pi_{j}:\R^d\rightarrow \mathcal E_{j}^\bot$ be the orthogonal projection. Suppose $\x\in \mathcal X_j$. Then, since $\pi_\x(i_\x)$ is constant with $\partial(\{\x\})\subseteq \mathcal X_1\cup\ldots\cup \mathcal X_{j-1}$, it follows that $\pi_{j-1}(\x(i_\x))$ is also constant, say $\pi_{j-1}(\x(i_\x))=\pi_{j-1}(y)$ for all $i_\x$, for some $y\in G_0\subseteq \Lambda$. This ensures that we have subsets $Y_1,\ldots,Y_s\subseteq G_0\subseteq \Lambda$ such that each  $\pi_{j-1}(X_j(k))=\pi_{j-1}(Y_j)$ is  constant and independent of the tuple $k=(i_\x)_{i_\x\in \mathcal X}$ with $|Y_j|=|\pi_{j-1}(Y_j)|=|\mathcal X_j|$,  for $j\in [1,s]$. By definition of an ordinary Reay system, each $\pi_{j-1}(Y_j)=\pi_{j-1}(X_j(k))$, for $j\in [1,s]$, is a linearly independent set of size $|Y_j|$, ensuring that $Y=Y_1\cup\ldots\cup Y_s$ is linearly independent.

Now consider arbitrary subsets $X_1,\ldots,X_s\subseteq \Lambda$ with $\pi_{j-1}(X_j)=\pi_{j-1}(Y_j)$ and $|X_j|=|\pi_{j-1}(X_j)|=|Y_j|$ for all $j\in [1,s]$, and set $X=X_1\cup\ldots\cup X_s$.
Since $\pi_0$ is the identity map, this forces $X_1=Y_1$, and a short inductive argument then gives $$\R\la X_1\cup\ldots\cup X_j\ra=\R\la Y_1\cup\ldots\cup Y_j\ra=\mathcal E_j\quad\mbox{ for all $j\in [0,s]$}.$$ Thus, since each $\pi_{j-1}(X_j)=\pi_{j-1}(Y_j)$ is a linearly independent set of size $|X_j|$, it follows that $X=X_1\cup\ldots\cup X_s$ is linearly independent.
For $j\in [0,s]$, let $$\Delta_j=\Z\la X_1\cup\ldots\cup X_j\ra,$$ which is a full rank lattice in $\mathcal E_j=\R\la X_1\cup\ldots\cup X_j\ra=\R\la \Delta_j\ra$ (as $X$ is linearly independent).  Since $\Lambda\cap\mathcal E_j$ is also a full rank lattice in $\mathcal E_j$ with $\Delta_j=\Z\la X_1\cup\ldots\cup X_j\ra\leq \Lambda\cap \mathcal E_j$, it follows that $(\Lambda\cap \mathcal E_j)/\Delta_j$ is a finite abelian group, for $j\in [1,s]$.

Let $j\in [1,s]$ be arbitrary, and let $x_1,\ldots,x_r\in X_j$ and $y_1,\ldots,y_r\in Y_j$ be the distinct elements of $X_j$ and $Y_j$ indexed so that $\pi_{j-1}(x_i)=\pi_{j-1}(y_i)$ for all $i\in [1,r]$. For each $i\in [1,r]$, we have  $x_i=y_i+\xi_{x_i}$ for some $\xi_{i}:=\xi_{x_i}\in \ker\pi_{j-1}=\mathcal E_{j-1}$, and  since $x_i,\,y_i\in \Lambda$, it follows that $\xi_i\in  \Lambda\cap \mathcal E_{j-1}$. We now have
\be\label{lattice-splitting}\Delta_j=\Delta_{j-1}+\Z\la X_j\ra=\Delta_{j-1}+\Z\la y_1+\xi_1,\ldots,y_r+\xi_r\ra.\ee
From \eqref{lattice-splitting}, we see that lattice $\Delta_j$ is completely determined by the value of the sublattice $\Delta_{j-1}$ as well as the values $\xi_i\mod \Delta_{j-1}$ for $i=1,\ldots,r$. (Recall that the set $Y=Y_1\ldots\cup Y_s$, and thus also the elements $y_1,\ldots,y_r\in Y_j$, are fixed). Moreover, assuming $\Delta_{j-1}$ is fixed, then distinct possibilities for the values of $\xi_1,\ldots,\xi_r$ modulo $\Delta_{j-1}$ yield distinct lattices $\Delta_j$, which can be seen by noting that $y_i+(\xi_i+\Delta_{j-1})\subseteq \Delta_j$ is precisely the subset of all elements $z\in \Delta_j$ with $\pi_{j-1}(z)=\pi_{j-1}(y_i)$ (in view of $\pi_{j-1}(X_j)=\pi_{j-1}(Y_j)$ being a linearly independent set of size $|X_j|=|Y_j|$). Note, if we apply the above setup with $X$ taken to be $Y$, then  $\xi_y=0$ for all $y\in Y$, meaning $\xi_x\equiv 0\mod \Delta_{j-1}$ for all $x\in X_j$ whenever $\Delta_j=\Z\la X_1\cup\ldots\cup X_j\ra=\Z\la Y_1\cup\ldots\cup Y_j\ra$ for all $j\in [1,s]$.

Since $(\Lambda\cap \mathcal E_{j-1})/\Delta_{j-1}$ is a finite group, there are only a finite number of possibilities for the value of each $\xi_i\in \Lambda\cap \mathcal E_{j-1}$ modulo $\Delta_{j-1}$. Thus an inductive argument on $j=1,\ldots,s$ shows that there are only a finite number of possibilities for each $\Delta_j$. Indeed, since $X_1=Y_1$, there is only one possibility for $\Delta_1$, completing the base of the induction, while  there are only a finite number of possibilities for each $\Delta_j$ with $\Delta_{j-1}$ fixed, one for each choice of the values of $\xi_1,\ldots,\xi_r$ modulo $\Delta_{j-1}$, with only a finite number of possibilities for $\Delta_{j-1}$ by induction hypothesis, leading to only a finite number of possibilities for $\Delta_j$.

All of the above basic observations will be crucial to several of our later arguments, including the following lemma showing how a pure virtual Reay system can be obtained from a given virtual Reay system by passing to appropriate subsequences of the representative sequences.

 \begin{lemma}\label{lem-VReay-lattice-rep}Let $\Lambda\leq \R^d$ be a full rank lattice, where $d\geq 0$, let $G_0\subseteq \Lambda$,  let  $\mathcal R=(\mathcal X_1\cup \{\mathbf v_1\},\ldots\,\mathcal X_s\cup \{\mathbf v_s\})$ be an anchored  virtual Reay system over $G_0$, let $\mathcal X=\bigcup_{i=1}^s\mathcal X_i$. Assume (by  Proposition \ref{prop-VReay-Lattice}) that the first few  terms in each representative sequence $\{\x(i)\}_{i=1}^\infty$ have been discarded so that  $\x(i)=\tilde\x(i)\in G_0$ is constant modulo $\R\la \partial(\{\x\})\ra$ for all $\x\in \mathcal X=\mathcal X_1\cup\ldots\cup \mathcal X_s$ and  $i\geq 1$.

\begin{itemize}
\item[1.]There are only a finite number of possibilities for $\Z\la X(k)\ra$ as we range over all tuples $k=(i_\x)_{\x\in\mathcal X}$.

\item[2.] If $\Delta=\Z\la X(\kappa)\ra$ for some tuple
 $\kappa=(\iota_\x)_{\x\in \mathcal X}$, then either
 \begin{itemize}
 \item[(a)] there is a finite set of indices $I_\Delta$  such that every  tuple $k=(i_\x)_{\x\in \mathcal X}$ with $\Delta=\Z\la X(k)\ra$ has $i_\x\in I_\Delta$ for some $\x\in \mathcal X$, or
 \item[(b)] replacing each representative sequence $\{\x(i_\x)\}_{\x\in \mathcal X}$ by an appropriate subsequence, for $\x\in \mathcal X$, we can obtain $\Z\la X(k)\ra=\Delta$ for all tuples $k=(i_\x)_{\x\in \mathcal X}$,
 \end{itemize}
 with (b) holding for at least one lattice $\Delta$.
 \end{itemize}
 In particular, by discarding the first few terms in each representative sequence $\{\x(i)\}_{i=1}^\infty$, for $\x\in\mathcal X$, it follows that $\mathcal R$ can be made $\Delta$-pure, for any $\Delta=\Z\la X(k)\ra$, by replacing each representative sequence $\{\x(i)\}_{i=1}^\infty$ by an appropriate subsequence.
 \end{lemma}

 \begin{proof}
Let $Y=Y_1\cup \ldots\cup Y_s$ with each $Y_j=X_j(\kappa)$, where $\kappa=(\iota_\x)_{\x\in \mathcal X}$ is an arbitrary tuple.
Per the discussion above Lemma \ref{lem-VReay-lattice-rep} (applied taking $X$ to be $X(k)$ for an arbitrary tuple $k$), for each $j\in [1,s]$,  there are only a finite number of possibilities for $\Z\la (X_1\cup \ldots\cup X_j)(k)\ra$ as we range over all tuples $k$, yielding Item 1. If $s=1$, then $\partial(\x)=\{0\}$ for $\x\in \mathcal X$, so  each $\x(i_\x)$ with $\x\in \mathcal X$ is constant (cf. Proposition \ref{prop-VReay-Lattice}), in which case the lemma holds trivially. Therefore we can assume $s\geq 2$ and proceed by induction on $s$. Let $\x_1,\ldots,\x_r\in \mathcal X_s$ be the distinct half-spaces from $\mathcal X_s$, adapt the abbreviations $i_j:=i_{\x_j}$ and $\iota_j:=\iota_{\x_j}$, and set $y_j=\x_j(\iota_j)$ for $j\in [1,r]$, so $Y_s=\{y_1,\ldots,y_r\}$. Since $\x_j(i_j)$ is constant modulo $\R^\cup\la \partial(\{\x_j\}\ra\subseteq \R^\cup \la \mathcal X_1\cup \ldots\cup \mathcal X_{s-1}\ra$,  we have $\x_j(i_{j})-y_j\in \R^\cup\la \mathcal X_1\cup\ldots\cup \mathcal X_{s-1}\ra=\R\la (X_1\cup\ldots\cup X_{s-1})(k)\ra=\R\la \Delta'(k)\ra$ for all $j\in [1,r]$ and $i_{j}\geq 1$.

 Let $\Delta=\Z\la X(\kappa)\ra$. For a more general  tuple $k=(i_\x)_{i_\x\in \mathcal X}$, let $\Delta(k)=\Z\la X(k)\ra$ and $\Delta'(k)=\Z\la (X_1\cup\ldots\cup X_{s-1})(k)\ra$. For $j\in [1,r]$, let $$\xi_j(k)=\xi_j(i_{j})=\x_j(i_{j})-y_j\in \Lambda\cap \R\la \Delta'(k)\ra.$$
 If there is a finite set of indices $I_\Delta$ such that all tuples $k=(i_\x)_{\x\in\mathcal X}$ with $\Z\la X(k)\ra=\Delta$ have $i_\x\in I_\Delta$ for some $\x\in \mathcal X$, then (a) holds and the induction is complete (note this cannot happen for every one of the finite number of possibilities for $\Delta$, by Item 1, as that would mean we could destroy all tuples $k=(i_\x)_{\x\in \mathcal X}$ by only removing a finite number of terms from the $\{\x(i_\x)\}_{i_\x=1}^\infty$ with $\x\in \mathcal X$, which is absurd). So we can assume otherwise.
 As we range over all tuples $k=(i_\x)_{\x\in \mathcal X}$ with $\Delta(k)=\Delta$,
 there are only a finite number of possibilities for the lattices $\Delta'(k)$ (by Item 1).
 If $k=(i_\x)_{\x\in\mathcal X}$ and $k'=(i'_\x)_{\x\in\mathcal X}$ are two tuples with $\Delta(k)=\Delta(k')=\Delta$ and  $\Delta'(k)=\Delta'(k')=\Delta'$, then the comments after \eqref{lattice-splitting} ensure that $$\x_j(i_{j})-y_j=\xi_j(k)\equiv \xi_j(k')=\x_j(i'_{j})-y_j\mod \Delta'\quad\mbox{ for all $j\in [1,r]$}.$$ 
Indeed, the lattice $\Delta(k)$ is completely determined by $\Delta'(k)$ and the values of $\xi_1(k),\ldots,\xi_r(k)$ modulo $\Delta'(k)$, with different values of $\xi_j(k)\mod \Delta'(k)$ giving rise to different lattices $\Delta(k)$, so there are fixed $\xi_1,\ldots,\xi_r\in \Lambda\cap \R\la \Delta'\ra$ such that a tuple  $k=(i_\x)_{\x\in \mathcal X}$ has $\Delta(k)=\Delta$ and $\Delta'(k)=\Delta'$ precisely when $\Delta'(k')=\Delta'$ and $\xi_j(k)\equiv \xi_j\mod \Delta'$ for all $j\in [1,r]$, where $k'=(i_\x)_{\x\in \mathcal X\setminus \mathcal X_s}$ is the restriction of the tuple $k$ to $\mathcal X\setminus \mathcal X_s$.
Consequently, if we let $\mathcal I_{\Delta'}$ consist of all tuples $k'=(i_\x)_{\x\in \mathcal X\setminus \mathcal X_s}$ with $\Z\la (X\setminus X_s)(k')\ra=\Delta'$ and, for each $j\in [1,r]$, let  $I_j$ consist of all indices $i_{j}$ with $\xi_j(i_{j})\equiv \xi_j\mod \Delta'$, then
\be\label{debecrossprod}\{k=(i_\x)_{\x\in \mathcal X}:\; (i_\x)_{\mathcal X\setminus \mathcal X_s}\in \mathcal I_{\Delta'} \und i_{j}\in I_j\mbox{ for all $j\in [1,r]$}\}=
\mathcal I_{\Delta'}\times I_1\times\ldots\times I_r\ee is the set of all tuples $k$ with $\Delta(k)=\Delta$ and $\Delta'(k)=\Delta'$.

As we range over all tuples $k$ with $\Delta(k)=\Delta$, there are only a finite number of possibilities for $\Delta'(k)$. If, for some possibility $\Delta'$, it is possible to remove a finite number of terms from each $\{\x(i_\x)\}_{i_\x=1}^\infty$ with $\x\in \mathcal X$ and destroy all occurrences where $\Delta'(k)=\Delta'$ from  among tuples the remaining tuples $k$, then  do so (replacing each representative sequence $\{\x(i_\x)\}_{i_\x=1}^\infty$ for $\x\in \mathcal X$ with an appropriate subsequence of sufficiently large indexed terms). Thus we may w.l.o.g. assume this is not possible for every  $\Delta'(k)$ that occurs among tuples $k$ with $\Delta(k)=\Delta$. Note, we cannot have destroyed all such lattices $\Delta'(k)$ by such a procedure, as this would contradict our assumption about there not existing a finite set of indices $I_\Delta$ such that every $k=(i_\x)_{\x\in \mathcal X}$ with $\Delta(k)=\Delta$ has $i_\x\in I_\Delta$ for some $\x\in \mathcal X$. Let $\Delta'$ be one possible lattice $\Delta'(k)$ that has survived.
Then we can apply the induction hypothesis to $(X_1\cup \ldots\cup X_{s-1})(k)$ allowing us to replace each $\{\x(i_\x)\}_{i_\x=1}^\infty$ for $\x\in \mathcal X\setminus \mathcal X_s$ with appropriate subsequences resulting in $\Delta'(k)=\Delta'$ for all tuples $k$. Each $I_j$ is infinite, else removing all indices from the finite set $I_j$ would destroy all tuples $k$ with $\Delta'(k)=\Delta'$ and $\Delta(k)=\Delta$, contrary to assumption. But now, in view of \eqref{debecrossprod}, we can replace each $\{\x_j(i_{j})\}_{i_{j}=1}^\infty$ with the infinite subsequence $\{\x_j(i_{j})\}_{i_{j}\in I_j}$, for $j\in [1,r]$, and thereby attain (b), which completes the induction and proof.
 \end{proof}

Let $\Lambda\leq \R^d$ be a full rank lattice and let $G_0\subseteq \Lambda$ be a finitary set with $\C(G_0)=\R^d$. Let $\mathcal X\in \mathfrak X(G_0)$ with  $\mathcal R=(\mathcal X_1\cup \{\mathbf v_1\},\ldots,\mathcal X_s\cup \{\mathbf v_s\})$  a purely virtual Reay system over $G_0$ realizing $\mathcal X$. Then Lemma \ref{lem-VReay-lattice-rep} ensures that, assuming $\mathcal X$ is fixed, there are only a finite number of  lattices $\Delta=\Z\la X(k)\ra\in \mathfrak P_\Z(G_0)$ as we range over all tuples $k=(i_\x)_{\x\in \mathcal X}$ with each $i_\x$ sufficiently large,
 and by refining the representative sequences $\{\x(i_\x)\}_{\x\in\mathcal X}$, it can be assumed via Lemma \ref{lem-VReay-lattice-rep} that $\mathcal R$ is $\Delta$-pure, for any $\Delta=\Z\la X(k)\ra$ with $k=(i_\x)_{\x\in \mathcal X}$ so long as  all $i_\x$ are sufficiently large.
We will see below that $\mathfrak P_\Z(G_0)$ is finite without need to restrict to a fixed $\mathcal X\in \mathfrak X(G_0)$.
If $\{\x(i_\x)\}_{i_\x=1}^\infty$ is a representative sequence for $\x\in \mathcal X$ in $\mathcal R$, then Proposition \ref{prop-orientedReay-halfspace-uniquelydetermines} ensures  it can also be substituted for a representative sequence for $\x$ in any other realization $\mathcal R'$ of $\mathcal X$. Thus the possibilities for $\Delta=\Z\la X(k)\ra$ depend only on $\mathcal X$ and not on the particular series decomposition $\mathcal X=\bigcup_{i=1}^{s}\mathcal X_i$ nor realization $\mathcal R$, and if $\mathcal X=\mathcal X_1\cup\ldots\cup \mathcal X_s$ is a series decomposition having a realization that is $\Delta$-pure, then any series decomposition of $\mathcal X$ has a realization that is $\Delta$-pure.
Furthermore, if $\Delta=\Z\la X\ra\in \mathfrak P_\Z(G_0)$, where $X\in X(G_0)$, then per earlier discussions, there is a realization $\mathcal R=(\mathcal X_1\cup \{\mathbf v_1\},\ldots,\mathcal X_s\cup \{\mathbf v_s\})$ with each $\x(i)\in X\subseteq G_0$ constant, for $\x\in \mathcal X:=\mathcal X_1\cup\ldots\cup \mathcal X_s$, ensuring there is a virtualization $\mathcal X$ of $X$ having a $\Delta$-pure realization.

\begin{definition}
For $\mathcal X\in \mathfrak X(G_0)$, let  $\mathfrak P_\Z(\mathcal X)$  denote all  $\Delta\in \mathfrak P_\Z(G_0)$ for which $\mathcal X$ has a $\Delta$-pure realization. This is independent of series decomposition, and every $\Delta\in \mathfrak P_\Z(G_0)$ has  $\Delta\in \mathfrak P_\Z(\mathcal X)$ for some $\mathcal X\in \mathfrak X(G_0)$, as explained above.
We let $\mathfrak X(G_0,\Delta)\subseteq \mathfrak X(G_0)$ consist of all $\mathcal X\in \mathfrak X(G_0)$ with $\Delta\in \mathfrak P_\Z(\mathcal X)$, and we let $X(G_0,\Delta)$ consist of all $X\in X(G_0)$ with $\Z\la X\ra=\Delta$, where $\Delta\in \mathfrak P_\Z(G_0)$.\end{definition}

\subsection{Finiteness Properties of Finitary Sets}

The sets $X(G_0)$, $\mathfrak X(G_0)$ and $\mathfrak X(G_0,\Delta)$ may be infinite. Our next goal is to show that, nonetheless, they still exhibit some finite-like behavior. Theorems \ref{thm-finitary-FiniteProps-I}, \ref{thm-finitary-FiniteProps-II} and \ref{thm-finitary-FiniteProps-III}
 contain some of the important finite-like properties possessed by a  finitary set,  explaining the choice of name. We begin with Theorem \ref{thm-finitary-FiniteProps-I}, which contains  the essential finitary property of $\mathfrak X(G_0)$ and $X(G_0)$ as well as  the finiteness of $\mathfrak P_\Z(G_0)$ and $\mathfrak P_{\R}(G_0)$.

\begin{theorem}\label{thm-finitary-FiniteProps-I}
Let $\Lambda\subseteq \R^d$ be a full rank lattice, where $d\geq 0$, and let $G_0\subseteq \Lambda$ be a finitary subset with $\C(G_0)=\R^d$.
\begin{itemize}
\item[1.]  There are only a finite number of irreducible sets  $\mathcal X\in \mathfrak X(G_0)$ and $X\in X(G_0)$.
\item[2.]  $\mathfrak B_\Z(G_0)$ and $\mathfrak P_\R(G_0)$ are both  finite.
\end{itemize}
\end{theorem}

\begin{proof}
1. If $\mathcal X\in \mathfrak X(G_0)$ is irreducible, then any realization $\mathcal R=(\mathcal X_1\cup\{\mathbf v_1\},\ldots,\mathcal X_s\cup \{\mathbf v_s\})$ of $\mathcal X$ must have $s=1$ with $\mathcal X=\mathcal X_1$, in which case every $\x\in \mathcal X$ is a one-dimensional half-space spanned by any representative. Thus to show there are only a finite number of irreducible sets $\mathcal X\in \mathfrak X(G_0)$, it suffices to show there are only a finite number of irreducible sets $X\in X(G_0)$. Assume by contradiction that this fails and let $\{X_i\}_{i=1}^\infty$ be a sequence of distinct irreducible sets $X_i\in X(G_0)$. Now $X\in X(G_0)$ being irreducible implies there exists an unbounded limit $u\in G_0^\infty$ of a radially convergent sequence of terms from $G_0$ such that $X$ minimally encases $-u$, which is equivalent to $X\cup \{u\}$ being a minimal positive basis. Since each $X_i$ is linearly independent, we have $|X_i|\leq d$, and so, by passing to a subsequence of $\{X_i\}_{i=1}^\infty$, we can assume all $X_i$ have the same size, say $|X_i|=s$ with  $X_i=\{x_i^{(1)},\ldots,x_i^{(s)}\}$ for all $i$. For each $i$, let $x_i^{(s+1)}\in G_0^\infty$ be such that $X_i\cup \{x_i^{(s+1)}\}$ is a minimal positive basis. By successively passing to an appropriate subsequence of each $\{x_i^{(j)}\}_{i=1}^\infty$, for $j=s+1,s,\ldots,1$, we can w.l.o.g. assume each $\{x_i^{(j)}\}_{i=1}^\infty$, for $j\in [1,s+1]$, is a radially convergent sequence with limit (say) $u_j$.
Since $x_i^{(s+1)}\in G_0^\infty$ for all $i$, it follows from Lemma \ref{lemma-radpoints-closed} that $u_{s+1}\in G_0^\infty$. For $j\in [1,s]$, we have  $x_i^{(j)}\in G_0\subseteq \Lambda$ with $\Lambda$ a lattice, so each sequence $\{x_i^{(j)}\}_{i=1}^\infty$ either eventually stabilizes to a nonzero constant value or else $u_j\in G_0^\infty$. Since the $X_i$ are all distinct, the latter cannot occur for all $j\in [1,s]$.

If $0\notin \C^*(u_1,\ldots,u_s,u_{s+1})$, then Lemma \ref{lem-finite-nozs} implies there is an open half-space $\mathcal E^\circ_+$ containing all $u_j$ for $j\in [1,s+1]$. Consequently, since $x_i^{(j)}/|x_i^{(j)}\|\rightarrow u_j$ for each $j\in [1,s+1]$, it follows that each $X_i\cup \{x_i^{(s+1)}\}$, with $i$ sufficiently large, will also be contained in the open half-space $\mathcal E^\circ_+$, ensuring $0\notin \C^*(X_i\cup \{x_i^{(s+1)}\})$. However, this contradicts that each $X_i\cup \{x_i^{(s+1)}\}$ is a minimal positive basis. Therefore we instead conclude that $0\in \C^*(u_1,\ldots,u_s,u_{s+1})$.
Thus there is some nonempty subset $J\subseteq [1,s+1]$ such that $U=\{u_j:\;j\in J\}$ is a minimal positive basis.
If $u_j\notin G_0^\infty$ for all $j\in J$, then $J\subseteq [1,s]$ and $x_i^{(j)}/\|x_i^{(j)}\|=u_j$ for all sufficiently large $i$ and $j\in J$.
In such case, $\{x_i^{(j)}:j\in J\}\subseteq X_i$ will be a minimal positive basis for all sufficiently large $i$, contradicting that each $X_i\cup\{x_i^{(s+1)}\}$ is a minimal positive basis with $J\subseteq [1,s]$.
Therefore we instead conclude that there is some $t\in J$ such that $u_t\in G_0^\infty$, and in view of $G_0$ being finitary, it then follows that $J\cap G_0^\infty=\{t\}$ (cf. the comments after Theorem \ref{thm-keylemmaII}).
In particular, $J\setminus \{t\}\subseteq [1,s]$ as $u_{s+1}\in G_0^\infty$.
But now $\{x_i^{(j)}:j\in J\setminus \{t\}\}\cup \{u_t\}$ is a minimal positive basis for all sufficiently large $i$ with $u_t\in G_0^\infty$, implying $\{x_i^{(j)}:j\in J\setminus \{t\}\}\in X(G_0)$.
However, since $\{x_i^{(j)}:j\in J\setminus \{t\}\}\subseteq X_i$ with $X_i$ irreducible, this is only possible if $J\setminus \{t\}=[1,s]$ and $t=s+1$.
Hence every $\{x_i^{(j)}\}_{i=1}^\infty$ for $j\in [1,s]$ eventually stabilizes, contrary to what we concluded in the previous paragraph. This establishes Item 1.

2. Since every subspace from $\mathfrak P_\R(G_0)$ is linearly spanned by a lattice from $\mathfrak P_\Z(G_0)$,  it suffices to show  $\mathfrak P_\Z(G_0)$ is finite.
For $s\in [1,d]$, let $X_s(G_0)\subseteq X(G_0)$ consist of all $X\in X(G_0)$ having a maximal series decomposition of length $s$, say $X=X_1\cup\ldots\cup X_s$. We proceed by induction on $s=1,2,\ldots,d$ to show that the number of lattices generated by subsets from $X_1(G_0)\cup\ldots\cup X_s(G_0)$ is finite. Note $X_1(G)$ is precisely the subset of all irreducible subsets  $X\in X(G_0)$, so  $X_1(G_0)$ is  finite by Item 1, ensuring the number of lattices generated by subsets $X\in X_1(G_0)$ is also finite. Thus the base $s=1$ of the induction is complete, and we assume $s\geq 2$.

Let $X\in X_s(G_0)$ be arbitrary and let $X=X_1\cup\ldots\cup X_s$ be a maximal series decomposition of $X$. Then $X\setminus X_s=X_1\cup\ldots\cup X_{s-1}$ is a maximal series decomposition by  Proposition \ref{prop-finitary-MaxDecompChar}. Let $\Delta'=\Z\la X\setminus X_s\ra$. By induction hypothesis, there are only a finite number of possibilities for the lattice $\Delta'$, so it suffices to show that there are only a finite number of lattices generated by $X'\in X_s(G_0)$ having a maximal series decomposition $X'=X'_1\cup\ldots\cup X'_{s}$ with $\Z\la X'\setminus X'_s\ra=\Delta'$, for each possible $\Delta'$. To this end, fix an arbitrary possible lattice $\Delta'$, let $\mathcal E=\R\la \Delta'\ra$, let $\pi:\R^d\rightarrow \mathcal E^\bot$ be the orthogonal projection.
Proposition
\ref{prop-finitary-Modulo-Inheritence} implies $\pi(G_0)$ is a finitary subset of the lattice $\pi(\Lambda)$.   Proposition \ref{prop-finitary-MaxDecompChar} implies  $\pi(X)\setminus \{0\}=\pi(X_s)\in X(\pi(G_0))$ is irreducible. As a result, Item 1 implies that there are only a finite number of possibilities for $\pi(X_s)$. Consequently, it suffices to show, for each possible irreducible set  $\pi(Y_s)=\pi(X_s)$, that there are only a finite number of lattices generated by a set $X\in X(G_0)$   having a maximal series decomposition $X=X_1\cup\ldots\cup X_s$ with $\Z\la X_1\cup \ldots\cup X_{s-1}\ra=\Delta'$ and $\pi(X_s)=\pi(Y)$. However, this follows from the discussion given after \eqref{lattice-splitting}, completing the proof.
\end{proof}

\begin{lemma}\label{lemma-radconv-super}Let $\mathcal R=(\mathcal X_1\cup\{\mathbf v_1\},\ldots,\mathcal X_s\cup \{\mathbf v_s\})$ be an oriented Reay system in $\R^d$, where $d\geq 0$, let $\mathcal B\subseteq \mathcal X_1\cup\ldots\cup \mathcal X_s$, and let $\pi:\R^d\rightarrow \R^\cup \la \mathcal B\ra^\bot$ be the orthogonal projection. For $\x\in \mathcal B$, let $\mathcal B_\x=\mathcal B\setminus \{\x\}\cup \partial(\{\x\})$ and
 let $\pi_\x:\R^d\rightarrow \R^\cup\la \mathcal B_\x\ra^\bot$ be the orthogonal projection.  Let
 $\{x_i\}_{i=1}^\infty$ be an asymptotically filtered sequence of terms $x_i\in \R^d$ with fully unbounded limit $\vec u=(u_1,\ldots,u_t)$.

  \begin{itemize}
\item[1.]
If  $-\vec u$ is minimally encased by $\mathcal B$, then $\|\pi_\x(x_i)\|\rightarrow \infty$ for all $\x\in \mathcal B$.

\item[2.] If $\vec u$ is a complete fully unbounded limit and $-\vec u$ is encased by $\mathcal B$ with $\{\pi_\x(x_i)\}_{i=1}^\infty$ unbounded for all $\x\in \mathcal B$, then $\mathcal B$ minimally encases $-\vec u$.
\item[3.]  If $-\vec u$ is minimally encased by $\mathcal B$ and $\{y_i\}_{i=1}^\infty$ is a  sequence of terms $y_i\in \R^d$ such that $\|y_i\|\in o(\|\pi_\x(x_i)\|)$ for all $\x\in \mathcal B$, then the sequence $\{x_i+y_i\}_{i=1}^\infty$ is an asymptotically  filtered sequence with fully unbounded limit $(u_1,u_2,\ldots,u_{r_\ell})$ (after discarding the first few terms), where $1=r_1<\ldots<r_\ell<r_{\ell+1}=t+1$ are the indices given by Proposition \ref{prop-orReay-minecase-char} regarding the minimal encasement of $-\vec u$ by $\mathcal B$.
\end{itemize}
\end{lemma}

\begin{proof}
Let $x_i=a_i^{(1)}u_1+\ldots+a_i^{(t)}u_t+\varepsilon_i$ be the representation of $\{x_i\}_{i=1}^\infty$ as an asymptotically filtered sequence with fully unbounded limit $\vec u$, so $a_i^{(j)}\rightarrow \infty$ for all $j\in [1,t]$, \ $a_i^{(j)}\in o(a_i^{(j-1)})$ for all $j\in [2,t]$, and $\|\varepsilon_i\|\in o(a_i^{(t)})$.

1.  Since $\mathcal B\subseteq \mathcal X_1\cup\ldots\cup \mathcal X_s$ minimally encases $-\vec u$, it must do so urbanely and be a support set.
Let $\x\in \mathcal B$ be arbitrary. Then $\mathcal B_\x\prec \mathcal B\subseteq \mathcal X_1\cup\ldots\cup \mathcal X_s$ and  $\pi_\x(x_i)=a_i^{(1)}\pi_\x(u_1)+\ldots+a_i^{(t)}\pi_\x(u_t)+\pi_\x(\varepsilon_i)$.
Proposition \ref{prop-orReay-minecase-char}.4 implies that $\mathcal B^{\pi_\x}$ minimally encases $-\pi_\x(\vec u)$.
If $\pi_\x(u_i)=0$ for all $i\in [1,t]$, then $\pi_\x(\vec u)$ is the empty tuple, in which case $\mathcal B^{\pi_\x}$ minimally encasing $-\pi_\x(\vec u)$ is only possible if $\mathcal B^{\pi_\x}=\emptyset$. Hence $\y\in \R^\cup \la \mathcal B_\x\ra$ for all $\y\in \mathcal B$, whence Proposition \ref{prop-orReay-BasicProps}.9 implies that $\mathcal B\subseteq \darrow \mathcal B_\x$. However, since $\mathcal B$ is a support set, we have $\mathcal B^*=\mathcal B$, in which case $\x\in \mathcal B\setminus \darrow \mathcal B_\x$, contradicting that $\mathcal B\subseteq \darrow \mathcal B_\x$. Therefore we must instead have $\pi_\x(u_j)\neq 0$ for some $j\in [1,t]$, and we may assume $r\in [1,t]$ is the minimal index with $\pi_\x(u_r)\neq 0$. Since  $a_i^{(j)}\in o(a_i^{(j-1)})$ for all $j\in [2,t]$ and $\|\varepsilon_i\|\in o(a_i^{(t)})$,  we have $\|\pi_\x(x_i)\|\sim a_i^{(r)}\|\pi_\x(u_r)\|$, and thus $\|\pi_\x(x_i)\|\rightarrow \infty$ as $a_i^{(r)}\rightarrow \infty$, completing Item 1.

2.
 Since $\vec u$ is a \emph{complete} fully unbounded limit, we have $\|\epsilon_i\|$ bounded.
 Since $\mathcal B$ encases $-\vec u$, it follows that there is some $\mathcal A\preceq \mathcal B$ that minimally encases $-\vec u$. Assume be contradiction that $\mathcal A\neq \mathcal B$. Then $\mathcal A\preceq \mathcal B_\x$ for some $\x\in \mathcal B$. Let $\tau:\R^d\rightarrow \R^\cup \la \mathcal A\ra^\bot$ be the orthogonal projection. Since $\mathcal A$ encases $-\vec u$, we have $\tau(u_j)=0$ for all $j$, and thus $\pi_\x(u_j)=0$ for all $j$ (since $\mathcal A\preceq \mathcal B_\x$ implies $\mathcal A\subseteq \darrow \mathcal B_\x$). But then $\{\pi_\x(x_i)\}_{i=1}^\infty=\{\pi_\x(\varepsilon_i)\}_{i=1}^\infty$ is a bounded sequence (as $\|\epsilon_i\|$ is bounded), contrary to hypothesis.

3.  Since $\mathcal B\subseteq \mathcal X_1\cup\ldots\cup \mathcal X_s$ minimally encases $-\vec u$, it is a support set (so $\mathcal B^*=\mathcal B$) and  must do so urbanely, and since $\vec u$ is fully unbounded, this ensures that $\mathcal B\neq \emptyset$.
Let $$\emptyset =\mathcal B_0\prec \mathcal B_1\prec\ldots\prec \mathcal B_\ell=\mathcal B$$ be the support sets and  let $1=r_1<\ldots<r_\ell<r_{\ell+1}=t+1$ be the indices given by  Proposition \ref{prop-orReay-minecase-char}. For $j\in [0,\ell]$, let $\pi_j:\R^d\rightarrow \R^\cup \la \mathcal B_j\ra^\bot$ be the orthogonal projection. For $j\in [1,\ell]$,  let
$\overline u_j=\pi_{j-1}(u_{r_j})/\|\pi_{j-1}(u_{r_j})\|$. Let $\mathcal F=(\R^\cup\la \mathcal B_1\ra,\ldots,\R^\cup\la \mathcal B_\ell\ra)$. In view of Proposition \ref{prop-orReay-minecase-char}.3 and Proposition \ref{prop-compat-filter}, let $$x_i=(b^{(1)}_i\overline u_1+w^{(1)}_i)+\ldots+(b_i^{(\ell)}\overline u_\ell+w_i^{(\ell)})+\varepsilon'_i$$
be a representation of $\{x_i\}_{i=1}^\infty$ as an $\mathcal F$-filtered sequence, so $b_i^{(j)}\in \Theta(a_i^{(r_j)})$ for $j\in [1,r]$.  Now $\pi_{\ell-1}(u_{r_{\ell}})/\|\pi_{\ell-1}(u_{r_{\ell}})\|=\overline u_\ell$ and
$\supp_{\pi_{\ell-1}(\mathcal R)}(-\overline u_\ell)=\mathcal B^{\pi_{\ell-1}}$ (by Proposition \ref{prop-orReay-minecase-char}.2). Thus $ B^{\tau}\cup \{\tau(\overline u_\ell)\}$ is a minimal positive basis of size $|\mathcal B^{\pi_{\ell-1}}|+1$ (by Propositions \ref{prop-orReay-BasicProps}.3 and \ref{prop-orReay-BasicProps}.4), where $\tau:\R^d\rightarrow \R^\cup\la \mathcal B_{\ell-1}\cup \partial(\mathcal B)\ra^\bot$ is the orthogonal projection, whence  $\overline u_{\ell}\notin \R\la \darrow B\setminus \{x\}\ra$ for any $\x\in\mathcal B\setminus\darrow \mathcal B_{\ell-1}$. (Indeed, if we write $-\overline u_{\ell}$ as a linear combination of elements from $\darrow B$ and apply $\tau$ to this linear combination, then we obtain a linear combination of $-\tau(\overline u_{\ell})$ using the elements from $\tau(B)$, which, as $B^{\tau}\cup \{\tau(\overline u_\ell)\}$ is a minimal positive basis of size $|\mathcal B^{\pi_{\ell-1}}|+1$, is only possible if the coefficient of each element $x\in B$ is strictly positive.)
Since $\mathcal B_{\ell-1}\prec \mathcal B_\ell=\mathcal B\neq \emptyset$, we have $\mathcal B_{\ell-1}\preceq \mathcal B_\x$ for some $\x\in \mathcal B$, and thus there exists some $\x\in \mathcal B\setminus\darrow \mathcal B_{\ell-1}$ (as $\mathcal B^*=\mathcal B$).
For this $\x$, we have $\pi_\x(\overline u_\ell)\neq 0$, for otherwise  $\overline u_\ell\in \R^\cup\la \mathcal B_\x\ra=\R\la \darrow B_\x\ra\subseteq \R\la \darrow B\setminus \{x\}\ra$, contrary to what was just noted. Since $\mathcal B_{\ell-1}\preceq \mathcal B_\x$ and $\pi_\x(\overline u_\ell)\neq 0$, we have $$\pi_\x(x_i)=b_i^{(\ell)}\pi_\x(\overline u_\ell)+\pi_\x(w_i^{(\ell)})+\pi_\x(\varepsilon'_i),$$ whence $\|\pi_\x(x_i)\|\in \Theta(b_i^{(\ell)})=\Theta(a_i^{(r_{\ell})})$ (as $\|\varepsilon'_i\|,\,\|w_i^{(\ell)}\|\in o(b_i^{(\ell)})$ in view of the $\mathcal F$-filtration representation for $\{x_i\}_{i=1}^\infty$). Thus the hypothesis $\|y_i\|\in o(\|\pi_\x(x_i)\|)$ implies that $\|y_i\|\in o(b_i^{(\ell)})$ and $\|y_i\|\in o(a_i^{(r_\ell)})$. It follows that $\{x_i+y_i\}_{i=1}^\infty$ is an asymptotically filtered sequence with fully unbounded limit $(u_1,\ldots,u_{r_\ell})$ (after discarding the first few terms), completing the proof.
\end{proof}

Suppose $\mathcal X,\,\mathcal X'\in \mathfrak X(G_0)$. We define a partial order on $\mathfrak X(G_0)$ (and thus also on each $\mathfrak X(G_0,\Delta)$), by declaring $$\mathcal X\preceq_\cup \mathcal X'$$ when the half-spaces from $\mathcal X$ are in bijective correspondence with the half-spaces from $\mathcal X'$, say with $\x\in \mathcal X$ corresponding to $\x'\in \mathcal X'$, such that, for every $\x\in \mathcal X$, we have  $\x\subseteq \x'$ and $\x^\circ\subseteq (\x')^\circ$, the latter meaning  any representative for $\x$ is also one for $\x'$. For $\mathcal A\subseteq \mathcal X$, let $\mathcal A'=\{\x':\;\x\in \mathcal A\}\subseteq \mathcal X'$ denote the  image of $\mathcal A$ under the bijection $\x\mapsto \x'$.   The relation $\preceq_\cup$ is clearly transitive and reflexive. We first make some observations regarding the defining condition before giving the argument that $\preceq_\cup$ is anti-symmetric.

Suppose $\x,\,\y\in \mathcal X$ with $\y\subseteq \x$. Then
$\y\subseteq \x\subseteq \x'$, so that $\x'$ contains a representative $y$ for $\y$, and thus also one for $\y'$ as a representative for $\y$ is also one for $\y'$. If $\y'\notin \darrow \x'$, then the representative $y$ for $\y'$ would be linearly independent from the representatives $\darrow x'$ for $\darrow \x'$  (as any set of representatives for $\mathcal X'$ is linearly independent), contradicting that $y\in \x'\subseteq \R^\cup\la \darrow \x'\ra=\R \la \darrow x'\ra$. Therefore we must have $\y'\in \darrow \x'$, i.e., $\y'\subseteq \x'$. Moreover, if $\y\subset \x$, then the injectivity of the map $\z\mapsto \z'$ ensures $\y'\neq \x'$. In summary, this shows that \be\label{zebra}\y\subset \x\quad\mbox{ implies }
\quad\y'\subset \x'\quad\mbox{ for $\x,\,\y\in \mathcal X$}.\ee
In particular, we  have \be\label{darowprime}(\darrow \mathcal A)'\subseteq \darrow (\mathcal A')=\darrow \mathcal A'\quad\und\quad (\mathcal A')^*\subseteq (\mathcal A^*)'\quad \mbox{for $\mathcal A\subseteq \mathcal X$}.\ee To see the second inclusion in \eqref{darowprime}, note that $(\mathcal A')^*$ and $(\mathcal A^*)'$ are both subsets of $\mathcal A'$. Thus, if the inclusion were to fail, then there would be some $\y\in \mathcal A\setminus \mathcal A^*$ with $\y'\in (\mathcal A')^*$. But then, $\y\in \mathcal A\setminus \mathcal A^*$ ensures there is some $\x\in \mathcal A$ with $\y\subset \x$, in turn implying $\y'\subset \x'$ by \eqref{zebra}, which contradicts that $\y'\in (\mathcal A')^*$.
Next, we observe that   \be\label{boundary-prime}(\darrow \partial(\{\x\}))'\subseteq \darrow \partial(\{\x'\})\quad\mbox{ for all $\x\in \mathcal X$}.\ee
Indeed, if $\y\in \darrow \partial(\{\x\})$, then $\y\subset \x$, which implies $\y'\subset \x'$ by  \eqref{zebra}, yielding \eqref{boundary-prime}.

The fact that we only have inclusions in \eqref{darowprime}  prompts us to define $\darrow \overline{\mathcal A},\,\mathcal A^{\circ *}\subseteq \mathcal X$ to be the subsets such that
$$(\darrow \overline{\mathcal A})'=\darrow \mathcal A'\quad\und\quad(\mathcal A^{\circ *})'=(\mathcal A')^*.$$ In view of \eqref{darowprime}, we have $$\nn\darrow \mathcal A\subseteq \darrow \overline{\mathcal A}\quad\und\quad\mathcal A^{\circ *}\subseteq \mathcal A^*.$$
Note $\darrow (\darrow \overline{\mathcal A})=\darrow \overline{\mathcal A}$, for if  $\x\subseteq \y\in \darrow \overline{\mathcal A}$, then $\y'\in \darrow \mathcal A'$, implying $\y'\subseteq \z'$ for some $\z'\in \mathcal A'$, whence \eqref{zebra} implies $\x'\subseteq \y'\subseteq \z'\in \mathcal A'$, which implies $\x'\in \darrow \mathcal A'$, and thus $\x\in \darrow \overline{\mathcal A}$.
We also have $(\mathcal A^{\circ *})^*=\mathcal A^{\circ *}$. Indeed, $(\mathcal A^{\circ *})^*\subseteq \mathcal A^{\circ *}$ holds by definition, while  if $\x\in \mathcal A^{\circ *}\setminus (\mathcal A^{\circ *})^*$, then  $\x\subset \y$ for some $\y\in \mathcal A^{\circ *}$, whence $\x'\subset \y'$ by \eqref{zebra} with $\x',\,\y'\in (\mathcal A')^*$ (by definition of $\mathcal A^{\circ *}$), which contradicts the definition of $*$.  Since $(\darrow \overline{\mathcal A^{\circ *}})'=\darrow (\mathcal A^{\circ *})'=\darrow((\mathcal A')^*)=\darrow \mathcal A'=(\darrow \overline{\mathcal A})'$, the injectivity of the map $\z\mapsto \z'$ implies $\darrow \overline{\mathcal A^{\circ *}}=\darrow \overline{\mathcal A}$.
Since $((\darrow \overline{\mathcal A})^{\circ *})'=((\darrow \overline{\mathcal A})')^*=(\darrow \mathcal A')^*=(\mathcal A')^*=(\mathcal A^{\circ *})'$, the injectivity of the map $\z\mapsto \z'$ implies $(\darrow \overline{\mathcal A})^{\circ *}=\mathcal A^{\circ *}$. In summary,
\be\label{base-prime-closureinterior} \darrow (\darrow \overline{\mathcal A})=\darrow \overline{\mathcal A},\quad (\mathcal A^{\circ *})^*=\mathcal A^{\circ *}, \quad \darrow \overline{\mathcal A^{\circ *}}=\darrow \overline{\mathcal A}, \quad\und\quad (\darrow \overline{\mathcal A})^{\circ *}=\mathcal A^{\circ *}.\ee Thus the arrow closure and star interior defined above behave similar to the closure and interior operations on convex sets.

The condition that any representative for $\x$ also be a representative for $\x'$ may be replaced by the alternate condition that  $\partial(\x)\subseteq \partial(\x')$ and $\x\nsubseteq \partial(\x')$, as the following argument shows. Suppose $\x\subseteq \x'$ and every representative for $\x$ is also a representative for $\x'$.
Since $x\in \x$ is a representative for $\x'$, we have $x\notin \partial(\x')$, ensuring $\x\nsubseteq \partial(\x')$.
If $\partial(\x)\subseteq \partial(\x')$ fails, then $\partial(\x)\subseteq \overline \x\subseteq \overline{\x'}$ ensures that there is some representative $x'$ for $\x'$ with $x'\in \partial(\x)$.
Note that $\{x'\}\cup \darrow \partial(\{x'\})$ is a basis for $\R^\cup\la \darrow \x'\ra$, with the representatives for $\x'$ those elements of $\R^\cup\la \darrow \x'\ra$ lying strictly on the positive side of the hyperplane $\R^\cup \la\darrow \partial(\{x'\})\ra$.
If $x=\alpha x'+\Summ{\y'\in \darrow \partial(\{\x'\})}\alpha_\y y'$ is a representative for $\x$, where $\alpha,\,\alpha_\y\in \R$ (possible as $\x\subseteq \x'$), then so too is $x-(\alpha+1)x'$ in view of $x'\in \partial(\x)$.
However, the coefficient of  $x'$ when expressing $x-(\alpha+1)x'$ as a linear combination of the elements from the basis $\{x'\}\cup \darrow \partial(\{x'\})$ is negative, meaning $x-(\alpha+1)x'$ lies on the \emph{negative} side of the hyperplane  $\R^\cup \la\darrow \partial(\{x'\})\ra$, which contradicts that the representative $x-(\alpha+1)x'$ for $\x$ should also be a representative for $\x'$.
Therefore we instead conclude that $\partial(\x)\subseteq \partial(\x')$, as desired.
Next suppose $\x\subseteq \x'$, \, $\partial(\x)\subseteq \partial(\x')$ and $\x\nsubseteq \partial(\x')$.
Let $x$ be an arbitrary representative for $\x$.
Then $x\in \x\subseteq \x'$, so if $x$ is not a representative for $\x'$, then $x\in \partial(\x')$.
Thus $\x\subseteq \partial(\x)+\R x\subseteq \partial(\x')+\partial(\x')=\partial(\x')$, contrary to assumption. Therefore we instead conclude that every representative $x$ for $\x$ is also a representative for $\x'$, and the equivalence of the two stated conditions follows.


Finally, to see that $\preceq_\cup$ is anti-symmetric, suppose $\mathcal X\preceq_\cup \mathcal X'$ and $\mathcal X'\preceq_\cup \mathcal X$ and let $\mathcal X=\mathcal X_1\cup \ldots\cup \mathcal X_s$  be a series decomposition of $\mathcal X$.  Since $\mathcal X\preceq_\cup \mathcal X'$, we have  $\partial(\x)\subseteq \partial(\x')$ for $\x\in \mathcal X$, implying  $\dim \partial(\x)\leq \dim\partial(\x')$. Thus $\Summ{\x\in \mathcal X}\dim \partial(\x)\leq \Summ{\x\in \mathcal X}\dim \partial(\x')$ when $\mathcal X\preceq_\cup \mathcal X'$. Since we also have $\mathcal X'\preceq_\cup \mathcal X$,  the reverse inequality likewise holds, yielding $\Summ{\x\in \mathcal X}\dim \partial(\x)=\Summ{\x\in \mathcal X}\dim \partial(\x')$. As $\dim \partial(\x)\leq \dim\partial(\x')$, this is only possible if $\dim \partial(\x)= \dim\partial(\x')$ for all $\x\in \mathcal X$. Consequently, since $\partial(\x)\subseteq \partial(\x')$ are subspaces, and since  $\x$ and $\x'$ share a common representative by definition of  $\mathcal X\preceq_\cup \mathcal X'$, it then follows  that $\partial(\x)=\partial(\x')$ and  $\overline\x=\overline \x'$ for $\x\in \mathcal X$.  Now $\partial(\x)$ is a subspace of dimension $|\darrow \partial(\{\x\})|$ by Proposition \ref{prop-orReay-BasicProps}.1. Likewise, $\partial(\x')$ is a subspace of dimension $|\darrow \partial(\{\x'\})|$. As a result, $\partial(\x)=\partial(\x')$ implies $|\darrow \partial(\{\x\})|=|\darrow \partial(\{\x'\})|$. By \eqref{boundary-prime}, we know  $(\darrow \partial(\{\x\}))'\subseteq \darrow \partial(\{\x'\})$ is subset with cardinality $|(\darrow \partial(\{\x\}))'|=|\darrow \partial(\{\x\})|=|\darrow \partial(\{\x'\})|$, whence equality must hold in \eqref{boundary-prime}, i.e., $$(\darrow \partial(\{\x\}))'= \darrow \partial(\{\x'\}).$$ We can now show $\x=\x'$ by induction on $j$, where $x\in \mathcal X_j$. For $\x\in \mathcal X_1$, we have   $\partial(\x')=\partial(\x)=\{0\}$, ensuring that $\partial(\x)\cap \x=\{0\}=\partial(\x')\cap \x'$, which forces  $\x=\x'$  as $\x=\overline \x'$.
Thus we can assume $j\geq 2$, completing the base of the induction.  Applying the induction hypothesis to all half-spaces from $ \darrow \partial(\{\x\})$ combined with equality holding in \eqref{boundary-prime} yields  $\partial(\x)\cap\x=\C^\cup \big(\darrow \partial(\{\x\})\big)=\C^\cup \big((\darrow \partial(\{\x\}))'\big)=\C^\cup \big(\darrow \partial(\{\x'\})\big)=\partial(\x')\cap \x'$. Combined with $\overline \x=\overline \x'$, it follows that $\x=\x'$, completing the induction,  which shows  $\preceq_\cup$ is anti-symmetric. Note, it is not initially  as evident that the bijection from $\mathcal X'$ to $\mathcal X$ given by $\mathcal X'\preceq_\cup \mathcal X$ should be the inverse of that from $\mathcal X$ to $\mathcal X'$, which is why we have argued above without using the assumption $(\x')'=\x$. Of course, the above argument establishes this, since it gives $(\x')'=\x'=\x$.
Note this argument shows that $\Summ{\x\in \mathcal X}\dim \partial(\x)=\Summ{\x\in \mathcal X}\dim \partial(\x')$ implies $\x=\x'$ for all $\x\in \mathcal X$.

Next we show  that  if $\mathcal B\subseteq \mathcal X$ minimally encases $-\vec u$, then   $(\mathcal B')^*$  minimally encases $-\vec u$.

\begin{proposition}\label{prop-finitary-finiteabove-align}
Let $\Lambda\subseteq \R^d$ be a full rank lattice, where $d\geq 0$, and let $G_0\subseteq \Lambda$ be a finitary subset with $\C(G_0)=\R^d$. Suppose $\mathcal X,\,\mathcal X'\in \mathfrak X(G_0)$ with $\mathcal X\preceq_\cup \mathcal X'$ and each $\x\in \mathcal X$ in bijective correspondence with the half-space $\x'\in \mathcal X'$. Let $\vec u=(u_1,\ldots,u_t)$ be a tuple of orthonormal vectors from $\R^d$. If  $\mathcal B\subseteq \mathcal X_1\cup\ldots\cup \mathcal X_s$ minimally encases $\vec u$, then $(\mathcal B')^*$ minimally encases $\vec u$.
\end{proposition}

\begin{proof}
Since $\x\subseteq \x'$ for all $\x\in \mathcal X$ (by definition of $\mathcal X\preceq_\cup \mathcal X'$), and since $\mathcal B$ minimally encases $\vec u$,  we conclude that $\mathcal B'$ encases $\vec u=(u_1,\ldots,u_t)$. Thus there is some  $\mathcal A'\preceq \mathcal B'$ which minimally encases $\vec u$, where $\mathcal A\subseteq \mathcal X$. It remains to show $\mathcal A'=(\mathcal B')^*$.

Since $\mathcal A'$ minimally encases $\vec u=(u_1,\ldots,u_t)$, we have $u_1,\ldots,u_t\in \R^\cup \la\mathcal A'\ra=\R\la \darrow A'\ra=\R\la \darrow \overline{A}\ra=\R^\cup \la \darrow \overline{\mathcal A}\ra$, with the first equality in view of Proposition \ref{prop-orReay-BasicProps}.1,  the second in view of the definition of $\darrow \overline{\mathcal A}$ and the fact that any representative for a half-space $\x\in \darrow \overline{\mathcal A}$ is a also a representative for the half-space $\x'\in (\darrow \overline{\mathcal A})'=\darrow \mathcal A'$, and the third  in view of the first equality in \eqref{base-prime-closureinterior} and Proposition \ref{prop-orReay-BasicProps}.1.
Thus, since $\mathcal B$ minimally encases $\vec u$, Proposition \ref{prop-orientedReay-ExtraUnique} and the first equality in \eqref{base-prime-closureinterior} imply that $\mathcal B\subseteq \darrow \overline{\mathcal  A}$, in turn implying $\mathcal B'\subseteq (\darrow \overline{\mathcal A})'=\darrow \mathcal A'$. Thus $(\mathcal B')^*\preceq \mathcal A'$ by \eqref{lassie3}.
On the other hand, $\mathcal A'\preceq \mathcal B'$ implies $\mathcal A'\subseteq \darrow \mathcal B'= \darrow (\mathcal B')^*$, whence $(\mathcal A')^*\preceq (\mathcal B')^*$ follows again by \eqref{lassie3}. Since $\mathcal A'=(\mathcal A')^*$ as $\mathcal A'$ minimally encases $\vec u$, we see that $(\mathcal B')^*\preceq \mathcal A'=(\mathcal A')^*\preceq (\mathcal B')^*$, yielding the desired conclusion $\mathcal A'=(\mathcal B')^*$.
\end{proof}

\begin{definition}
Let $\mathfrak X^*(G_0)$ be the set consisting of all \emph{maximal} series decompositions of the sets $\mathcal X\in \mathfrak X (G_0)$. We informally write the elements of  $\mathfrak X^*(G_0)$ in the form $\bigcup_{i=1}^s\mathcal X_i$, where $\mathcal X=\bigcup_{i=1}^s\mathcal X_i$ is a maximal series decomposition of some $\mathcal X\in \mathfrak X(G_0)$. For $\Delta\in \mathfrak P_\Z(G_0)$, let $\mathfrak X^*(G_0,\Delta)$ be the set consisting of all maximal series decompositions of the sets $\mathcal X\in \mathfrak X (G_0,\Delta)$.
\end{definition}

If $\bigcup_{i=1}^s\mathcal X_i\in \mathfrak X^*(G_0,\Delta)$ and we set $\mathcal X=\bigcup_{i=1}^s\mathcal X_i\in \mathfrak X(G_0,\Delta)$, then there is a  realization $\mathcal R=(\mathcal X_1\cup\{\mathbf v_1\},\ldots,\mathcal X_s\cup\{\mathbf v_s\})$ of $\mathcal X$ as well as a $\Delta$-pure realization $\mathcal R'$ of $\mathcal X$, in which case the comments above the definition of $\mathfrak P_\Z(\mathcal X)$ allow us to assume (by exchanging the representative sequences in $\mathcal R$ for those from $\mathcal R'$), that $\mathcal X$ has a $\Delta$-pure realization $(\mathcal X_1\cup\{\mathbf v_1\},\ldots,\mathcal X_s\cup\{\mathbf v_s\})$. We defined a partial order $\preceq_\cup$ on $\mathfrak X(G_0)$ earlier.
We now define a related partial order $\preceq_\cup^*$ on $\mathfrak X^*(G_0)$ as follows. Let $\bigcup_{i=1}^{s}\mathcal X_i,\,\bigcup_{i=1}^s\mathcal X'_i\in \mathfrak X^*(G_0)$ and set  $\mathcal X=\bigcup_{i=1}^{s}\mathcal X_i$ and $\mathcal X'=\bigcup_{i=1}^s\mathcal X'_i$.
Then we declare $\bigcup_{i=1}^{s}\mathcal X_i\preceq_\cup^*\bigcup_{i=1}^s\mathcal X'_i$ if there is bijection between $\mathcal X$ and $\mathcal X'$ showing $\mathcal X\preceq_\cup \mathcal X'$ which restricts to a bijection between $\mathcal X_j$ and $\mathcal X'_j$ for all $j\in [1,s]$.
Thus, if we let $\x\mapsto \x'$ denote the bijection showing $\mathcal X\preceq_\cup \mathcal X'$, then each $\mathcal X'_i=\{\x':\; \x\in \mathcal X_i\}$, so that the notation agrees with our previous definition of  $\mathcal A'$ for  $\mathcal A\subseteq \mathcal X$.
By definition, $\bigcup_{i=1}^{s}\mathcal X_i\preceq_\cup^*\bigcup_{i=1}^s\mathcal X'_i$ implies $\mathcal X\preceq_\cup \mathcal X'$, in which case anti-symmetry for $\preceq_\cup^*$ follows from that for $\preceq_\cup$ (recall that we showed $\x'=\x$ when $\mathcal X\preceq_\cup \mathcal X'$ and $\mathcal X'\preceq_\cup \mathcal X$), while reflexivity and transitivity are immediate. Thus $\preceq_\cup^*$ is a partial order on $\mathfrak X^*(G_0)$, and so restricts to one on  $\mathfrak X^*(G_0,\Delta)$ as well. The relation between the partial orders $\preceq_\cup$ and $\preceq_\cup^*$ is given by the following proposition.

\begin{proposition}\label{prop-VReay-poset-equiv}
Let $\Lambda\subseteq \R^d$ be a full rank lattice, where $d\geq 0$, let $G_0\subseteq \Lambda$ be a finitary subset with $\C(G_0)=\R^d$, and let $\mathcal X,\,\mathcal X'\in \mathfrak X(G_0)$. Then $\mathcal X\preceq_\cup \mathcal X'$ if and only if there are maximal series decompositions $\mathcal X=\bigcup_{i=1}^s\mathcal X_i$ and $\mathcal X'=\bigcup_{i=1}^s\mathcal X'_i$ such that $\bigcup_{i=1}^s\mathcal X_i\preceq_\cup^* \bigcup_{i=1}^s\mathcal X'_i$.
\end{proposition}

\begin{proof}
If $\bigcup_{i=1}^s\mathcal X_i\preceq_\cup^* \bigcup_{i=1}^s\mathcal X'_i$, then $\mathcal X\preceq_\cup \mathcal X'$ follows by definition of $\preceq_\cup^*$, so one direction is trivial. Assume $\mathcal X\preceq_\cup \mathcal X'$ and let $\mathcal X'=\bigcup_{i=1}^s\mathcal X'_i$ be any maximal series decomposition of $\mathcal X'$. It suffices to show $\mathcal X=\bigcup_{i=1}^s\mathcal X_i$ is also a maximal series decomposition.
Since any representative of $\x$ is also a representative for $\x'$, Proposition \ref{prop-orReay-BasicProps}.1 implies $$\R^\cup \la \mathcal X_1\cup\ldots\cup \mathcal X_{j}\ra=\R\la X_1\cup\ldots\cup X_{j}\ra=\R\la X'_1\cup\ldots\cup X'_{j}\ra=\R^\cup \la \mathcal X'_1\cup\ldots\cup \mathcal X'_{j}\ra$$ for all $j\in [0,s]$.
For $j\in[1,s]$, let $\pi_{j-1}:\R^d\rightarrow \R^\cup \la \mathcal X_1\cup\ldots\cup \mathcal X_{j-1}\ra^\bot$ be the orthogonal projection.

Let $\mathcal R'=(\mathcal X'_1\cup \{\mathbf v'_1\},\ldots,\mathcal X'_s\cup \{\mathbf v'_s\})$ be a realization of $\mathcal X'=\bigcup_{i=1}^s\mathcal X'_i$. Since each $\x\in \mathcal X'_1$ has trivial boundary, it follows that $\x=\x'$ for all $\x\in \mathcal X_1$ (as $\x\subseteq\x'$ forces this in view of $\x'$ being one-dimensional). Thus $(\mathcal X'_1\cup \{\mathbf v'_1\})=(\mathcal X_1\cup \{\mathbf v'_1\})$ is a purely virtual Reay system with $\mathcal X'_1=\mathcal X_1$ irreducible (by Proposition \ref{prop-finitary-MaxDecompChar} applied to the maximal series decomposition $\bigcup_{i=1}^s\mathcal X'_i$), showing that $\mathcal X_1=\mathcal X'_1\in \mathfrak X^*(G_0)$. We  proceed by induction on $j$ to show $\bigcup_{i=1}^{j}\mathcal X_i\in \mathfrak X^*(G_0)$ with $\pi_{j-1}(\mathcal X_j)=\pi_{j-1}(\mathcal X'_j)$, having just completed the  base case $j=1$. By induction hypothesis, $\bigcup_{i=1}^{j-1}\mathcal X_i\in \mathfrak X^*(G_0)$, so there is realization  $\mathcal R_{j-1}=(\mathcal X_1\cup\{\mathbf v_1\},\ldots,\mathcal X_{j-1}\cup \{\mathbf v_{j-1}\})$.
 Let $\vec u_{\mathbf v'_j}=(u_1,\ldots,u_t)$. Then $-\vec u_{\mathbf v'_j}^\triangleleft$ is a fully unbounded (or trivial) limit of an asymptotically filtered sequence of terms from $G_0$ which is minimally encased by $\partial(\{\mathbf v'_j\})\subseteq \mathcal X'_1\cup \ldots\cup \mathcal X'_{j-1}$, thus ensuring that $u_i\in \R^\cup \la \mathcal X'_1\cup\ldots\cup  \mathcal X'_{j-1}\ra= \R^\cup \la \mathcal X_1\cup\ldots\cup \mathcal X_{j-1}\ra$ for all $i<t$. Applying Proposition \ref{prop-finitary-basics}.3 to $\mathcal R_{j-1}$, we conclude that there is a support set $\mathcal B\subseteq \mathcal X_1\cup\ldots\cup \mathcal X_{j-1}$ which minimally encases $-\vec u^\triangleleft_{\mathbf v'_j}$, allowing us to define a half-space $\mathbf v_j$ with $\partial(\{\mathbf v_j\})=\mathcal B$ using the limit $\vec u_{\mathbf v_j}=\vec u_{\mathbf v'_j}$ with $\mathbf v_j(i)=\mathbf v'_j(i)$ for all $i$. By \eqref{boundary-prime}, we have $\partial(\{\x\})'\subseteq \darrow \partial(\{\x'\})\subseteq \mathcal X'_1\cup\ldots\cup \mathcal X'_{j-1}$ for all $\x\in \mathcal X_j$, ensuring via the bijection $\z\mapsto \z'$ that $\partial(\{\x\})\subseteq \mathcal X_1\cup\ldots\cup\mathcal X_{j-1}$ for all $\x\in \mathcal X_j$.
 Thus, since
$ \R^\cup \la \mathcal X'_1\cup\ldots\cup  \mathcal X'_{j-1}\ra= \R^\cup \la \mathcal X_1\cup\ldots\cup \mathcal X_{j-1}\ra$, it follows that $\pi_{j-1}(\mathcal X_j\cup \{\mathbf v_j\})=\pi_{j-1}(\mathcal X'_j\cup \{\mathbf v'_j\})$ with
 $(\mathcal X_1\cup\{\mathbf v_1\},\ldots,\mathcal X_j\cup \{\mathbf v_j\})$ a purely virtual Reay system, showing $\bigcup_{i=1}^j\mathcal X_i\in \mathfrak X^*(G_0)$. This
 completes the induction. The case $j=s$ shows $\mathcal X=\bigcup_{i=1}^s\mathcal X_i$ is a series decomposition. Since  $ \R^\cup \la \mathcal X'_1\cup\ldots\cup  \mathcal X'_{j-1}\ra= \R^\cup \la \mathcal X_1\cup\ldots\cup \mathcal X_{j-1}\ra$, applying Proposition \ref{prop-finitary-MaxDecompChar}  to the maximal decomposition $\bigcup_{i=1}^s\mathcal X'_i$ gives that  $\pi_{j-1}(\mathcal X'_j)=\pi_{j-1}(\mathcal X_j)$ is irreducible for all $j\in [1,s]$, and then a second application of Proposition \ref{prop-finitary-MaxDecompChar} to $\bigcup_{i=1}^s\mathcal X_i$ implies that  $\mathcal X=\bigcup_{i=1}^s\mathcal X_i$ is maximal.
\end{proof}

Theorem \ref{thm-finitary-FiniteProps-II} and Corollary \ref{cor-finitary-FiniteProps-II} contain the finiteness from above property for the partial orders $\preceq_\cup$ and $\preceq_\cup^*$ associated to a finitary set $G_0$.

\begin{theorem}\label{thm-finitary-FiniteProps-II}
Let $\Lambda\subseteq \R^d$ be a full rank lattice, where $d\geq 0$, and let $G_0\subseteq \Lambda$ be a finitary subset with $\C(G_0)=\R^d$. Then $\mathsf{Max}\big(\mathfrak X^*(G_0,\Delta),\preceq^*_\cup)$ is finite for every $\Delta\in \mathfrak P_\Z(G_0)$.
\end{theorem}

\begin{proof}
For $s\in [1,d]$, let $\mathfrak X^*_s(G_0,\Delta)\subseteq \mathfrak X^*(G_0,\Delta)$ consist of all maximal series decompositions $\mathcal X=\bigcup_{i=1}^s\mathcal X_i$ of some $\mathcal X\in \mathfrak X(G_0,\Delta)$ having  length $s$.  Note $\mathsf{Max}\big(\mathfrak X^*(G_0,\Delta)\big)=\bigcup_{s=1}^d \mathsf{Max}\big(\mathfrak X^*_s(G_0,\Delta)\big)$ since comparable series decompositions  under $\preceq_\cup^*$ must have the same length.
 We proceed by induction on $s=1,2,\ldots,d$ to show that $\mathsf{Max}\big(\mathfrak X^*_s(G_0,\Delta)\big)$ is finite for every $\Delta\in \mathfrak P_\Z(G_0)$.  Note $\mathfrak X^*_1(G_0)$ is in one-to-one correspondence with the subset of all irreducible   $\mathcal  X\in \mathfrak X(G_0)$, so  $ \mathsf{Max}(\mathfrak X_1^*(G_0,\Delta))$ is finite by Theorem \ref{thm-finitary-FiniteProps-I}.1. Thus the base $s=1$ of the induction is complete, and we  assume $s\geq 2$. Fix  $\Delta\in \mathfrak P_\Z(G_0)$.

 Let $$\bigcup_{i=1}^s\mathcal X_i\in \mathsf{Max}(\mathfrak X_s^*(G_0,\Delta))$$ be arbitrary, set $\mathcal X=\bigcup_{i=1}^s\mathcal X_i$, and let $$\mathcal R=(\mathcal X_1\cup\{\mathbf v_1\},\ldots,\mathcal X_s\cup\{\mathbf v_s\})$$ be a $\Delta$-pure realization.
By Lemma \ref{lem-VReay-lattice-rep}, we can assume (by passing to subsequences of the representative sequences as need be) that  $\Z\la (X\setminus X_s)(k)\ra=\Delta'$ is constant for all tuples $k=(i_\x)_{\x\in \mathcal X}$, while
$\mathcal X\setminus \mathcal X_s=\bigcup_{i=1}^{s-1}\mathcal X_i$ is  a maximal series decomposition by Proposition \ref{prop-finitary-MaxDecompChar}. Hence \be\label{s-1-go}\mathcal X_1\cup\ldots\cup \mathcal X_{s-1}\in \mathfrak X^*_{s-1}(G_0,\Delta').\ee 

\subsection*{Claim A}  $\mathcal X_1\cup\ldots\cup \mathcal X_{s-1}\in \mathsf{Max}(\mathfrak X_{s-1}^*(G_0,\Delta')).$

\begin{proof}Suppose by contradiction $\mathcal X_1\cup\ldots\cup \mathcal X_{s-1}\notin \mathsf{Max}(\mathfrak X_{s-1}^*(G_0,\Delta'))$.
Then, in view of \eqref{s-1-go},  there exists some $\bigcup_{i=1}^{s-1}\mathcal X'_i\in \mathfrak X^*(G_0,\Delta')$ with $\mathcal X_1\cup\ldots\cup \mathcal X_{s-1}\prec^*_\cup \mathcal X'_1\cup\ldots\cup \mathcal X'_{s-1}$.
 Let $\x'\in (\mathcal X\setminus \mathcal X_s)'$ be the image of $\x\in \mathcal X\setminus \mathcal X_s$ under the bijection showing $\mathcal X\setminus \mathcal X_s\prec_\cup (\mathcal X\setminus \mathcal X_s)'$. Let $\mathcal R'=(\mathcal X'_1\cup \{\mathbf w'_1\},\ldots,\mathcal X'_{s-1}\cup \{\mathbf w'_{s-1}\})$ be a $\Delta'$-pure realization of $(\mathcal X\setminus \mathcal X_s)'=\bigcup_{i=1}^{s-1}\mathcal X'_i$. Since $\mathcal X\setminus \mathcal X_s\prec_\cup (\mathcal X\setminus \mathcal X_s)'$ ensures that a set of representatives for $\mathcal X\setminus \mathcal X_s$ is also a set of representatives for $(\mathcal X\setminus \mathcal X_s)'$,  we must have $\R^\cup \la \mathcal X\setminus \mathcal X_s\ra=\R\la (X\setminus X_s)(k)\ra=\R\la (X\setminus X_s)'(k)\ra=\R^\cup \la (\mathcal X\setminus \mathcal X_s)'\ra$ by Proposition \ref{prop-orReay-BasicProps}.1.
 For each $\y\in \mathcal X_s\cup \{\mathbf v_s\}$, we have $-\vec u_\y^\triangleleft$ minimally encased by $\partial(\{\y\})\subseteq \mathcal X_1\cup\ldots\cup \mathcal X_{s-1}=\mathcal X\setminus \mathcal X_s$.
Thus Proposition \ref{prop-finitary-finiteabove-align} (applied to $\mathcal X\setminus \mathcal X_s\prec_\cup (\mathcal X\setminus\mathcal X_s)'$) implies that $(\partial(\{\y\})^{\circ *})'=(\partial(\{\y\})')^*$ minimally encases $-\vec u_\y$.
This allows us to define a new half-space $\y'$ such that $\partial(\{\y'\})=(\partial(\{\y\})')^*\subseteq \partial(\{\y\})'\subseteq (\mathcal X\setminus \mathcal X_s)'=\mathcal X'_1\cup\ldots\cup \mathcal X'_{s-1}$ with $\vec u_{\y'}=\vec u_\y$ and  $\y(i)=\y'(i)$ for all $i$,  making $$\mathcal R''=(\mathcal X'_1\cup \{\mathbf w'_1\},\ldots,\mathcal X'_{s-1}\cup \{\mathbf w'_{s-1}\}, \mathcal X'_{s}\cup \{\mathbf w'_{s}\}),$$ where $\mathcal X'_{s}=\{\y':\;\y\in \mathcal X_s\}$, $\mathbf w'_{s}=\mathbf v'_s$ and $\vec u_{\mathbf w'_{s}}=\vec u_{\mathbf v_s}$, a virtual Reay system over $G_0$ in view of $\R^\cup \la \mathcal X\setminus \mathcal X_s\ra=\R^\cup \la (\mathcal X\setminus \mathcal X_s)'\ra$, which ensures that  \be\label{chococake}\pi(\mathcal X'_s\cup \{\mathbf w'_s\})=\pi(\mathcal X_s\cup \{\mathbf v_s\})\ee with $\pi:\R^d\rightarrow \R^\cup \la \mathcal X\setminus \mathcal X_s\ra^\bot$ the orthogonal projection.
For $\y\in \mathcal X_s\cup\{\mathbf v_s\}$, we have
 $$\partial(\y)=\R^\cup \la  \darrow \partial(\{\y\})\ra\subseteq
 \R^\cup\la \darrow \overline{\partial(\{\y\})}\ra=\R^\cup\la \darrow \partial(\{\y\})'\ra=\R^\cup\la (\partial(\{\y\})')^*\ra=\partial(\y')$$
 and $\partial(\y)\cap \y=\C^\cup\big(\partial(\{\y\})\big)\subseteq \C^\cup\big(\partial(\{\y\})'\big)\subseteq \C^\cup \big(\partial(\{\y'\})\big)=\partial(\y')\cap \y'$ (with the second inclusion by \eqref{boundary-prime}, and the first as $\x\subseteq \x'$ for all $\x$).
 Hence, since $\y$ and $\y'$ share a common representative, it follows that $\y\subseteq \y'$ for all $\y\in \mathcal X_s$, and that any representative for $\y$ is also a representative for $\y'$.
Since $\mathcal R'$ is purely virtual, and since $\vec u_{\mathbf w'_{s}}=\vec u_{\mathbf v_s}$ is fully unbounded (as $\mathcal R$ is purely virtual),   it follows that $\mathcal R''$ is  purely virtual, showing that $\bigcup_{i=1}^s\mathcal X'_i$ is a series decomposition.
Since $\mathcal X=\bigcup_{i=1}^s\mathcal X_i$ is a maximal series decomposition, Proposition \ref{prop-finitary-MaxDecompChar}  ensures that $\pi(\mathcal X_s)=\pi(\mathcal X'_s)$ is irreducible (the equality follows from \eqref{chococake}), and now two further applications of Proposition \ref{prop-finitary-MaxDecompChar}, one  to the maximal series decomposition $(\mathcal X\setminus \mathcal X_s)'=\bigcup_{i=1}^{s-1}\mathcal X'_i$ and one to the series decomposition $\mathcal X'=\bigcup_{i=1}^s\mathcal X'_i$, implies that   $\mathcal X'=\bigcup_{i=1}^s\mathcal X'_i$ is maximal.
By construction, each representative of $\x\in \mathcal X$ is a representative of $\x'\in \mathcal X'_1\cup\ldots\cup \mathcal X'_{s}$, while $\x\subseteq \x'$ for $\x\in \mathcal X\setminus \mathcal X_s$ (since $\mathcal X\setminus \mathcal X_s\prec_\cup (\mathcal X\setminus \mathcal X_s)'$) as well as for $\x\in \mathcal X_s$ (as just argued).
Thus  \be\label{welltown}\mathcal X_1\cup\ldots\cup\mathcal X_{s-1}\cup  \mathcal X_s\prec^*_\cup \mathcal X'_1\cup\ldots\cup\mathcal X'_{s-1}\cup \mathcal X'_s,\ee with the strict inequality following as $\x\neq \x'$ for some $\x\in \mathcal X\setminus \mathcal X_s$ in view of $\mathcal X\setminus \mathcal X_s\prec_\cup (\mathcal X\setminus \mathcal X_s)'$.
Since the representative sequences for $\mathcal X'_{s}$ are the same as those from $\mathcal X_s$ with $\Z\la (X\setminus X_s)(k)\ra=\Z\la (X\setminus X_s)'(k)\ra=\Delta'$, it follows that
$\Delta=\Z\la X(k)\ra=\Z\la (X\setminus X_s)(k)\ra+\Z\la X_s(k)\ra=\Delta'+\Z\la X_s(k)\ra=\Delta'+\Z\la X'_{s}(k)\ra=\Z\la X'(k)\ra$, for any tuple $k$, ensuring that $\Delta\in \mathfrak P_\Z(\mathcal X')$. Thus $\bigcup_{i=1}^s\mathcal X'_i\in  \mathfrak X_{s}^*(G_0,\Delta)$, which combined with \eqref{welltown} contradicts that $\bigcup_{i=1}^s\mathcal X_i\in \mathsf{Max}(\mathfrak X_s^*(G_0,\Delta))$. This completes Claim A.\end{proof}

By induction hypothesis, $\mathsf{Max}\big(\mathfrak X_{s-1}^*(G_0,\Delta')\big)$ is finite for any of the finite choices for $\Delta'\in \mathfrak P_\Z(G_0)$.
Thus, since $\bigcup_{i=1}^s\mathcal X_i\in \mathsf{Max}(\mathfrak X_s^*(G_0,\Delta))$ was arbitrary, it follows in view of  \eqref{s-1-go} and Claim A that every $\bigcup_{i=1}^s\mathcal X_i\in  \mathsf{Max}(\mathfrak X_s^*(G_0,\Delta))$ has  $\bigcup_{i=1}^{s-1}\mathcal X_i \in  \mathsf{Max}(\mathfrak X_{s-1}^*(G_0,\Delta'))$ for some $\Delta'\in \mathfrak P_\Z(G_0)$.
Thus to show $\mathsf{Max}(\mathfrak X_s^*(G_0,\Delta))$ is finite, it suffices to show there are only a finite number of $\bigcup_{i=1}^s\mathcal X_i\in  \mathsf{Max}(\mathfrak X_s^*(G_0,\Delta))$ extending each possible maximal series decomposition $\bigcup_{i=1}^{s-1}\mathcal X_i\in\mathsf{Max}(\mathfrak X_{s-1}^*(G_0,\Delta'))$. To this end, fix $\Delta'\in \mathfrak P_\Z(G_0)$ and $\bigcup_{i=1}^{s-1}\mathcal X_i\in  \mathsf{Max}(\mathfrak X^*_{s-1}(G_0,\Delta'))$ and let $\bigcup_{i=1}^s\mathcal X_i\in \mathsf{Max}(\mathfrak X^*_s(G_0,\Delta))$ be arbitrary having $\mathcal R=(\mathcal X_1\cup\{\mathbf v_1\},\ldots,\mathcal X_s\cup\{\mathbf v_s\})$ as a $\Delta$-pure realization with $(\mathcal X_1\cup\{\mathbf v_1\},\ldots,\mathcal X_{s-1}\cup\{\mathbf v_{s-1}\})$ a $\Delta'$-pure realization. Set $\mathcal X=\bigcup_{i=1}^s\mathcal X_i$ and $\mathcal Y=\bigcup_{i=1}^{s-1}\mathcal X_i$, so $\mathcal Y$ is fixed and  $\mathcal X_s$ is arbitrary (subject to the constraints described).

Let $\mathcal E=\R^\cup \la \mathcal X_1\cup\ldots\cup \mathcal X_{s-1}\ra=\R^\cup \la \mathcal Y\ra=\R\la \Delta'\ra$ and let $\pi:\R^d\rightarrow \mathcal E^\bot$ be the orthogonal projection. Since $\mathcal R$ is $\Delta$-pure, we have  $\x(i)=\tilde \x(i)$ for all $i\geq 1$ and $\x\in \mathcal X_s$, with $\pi(\x(i))$ constant. Thus $\pi(X_s(k))$ is a fixed set for all tuples $k$.
By Proposition \ref{prop-finitary-MaxDecompChar}, the set of one-dimensional half-spaces $\pi(\mathcal X_s)\in \mathfrak X(\pi(G_0))$ is irreducible, and thus  $\pi(X_s(k))\in X(\pi(G_0))$ is also irreducible. As a result, Theorem \ref{thm-finitary-FiniteProps-I}.1 implies  there are only a finite number of possibilities for $\pi(\mathcal X_s)$ and $\pi(X_s(k))$. Thus it suffices to show there are only a finite number of possibilities for $\mathcal X$ when $\pi(\mathcal X_s)=\mathcal Y_s$ and also $\pi(X_s(k))=Y_s$ are fixed, say for the fixed subsets $Y_s\subseteq \pi(G_0)\subseteq \mathcal E^\bot$ and $\mathcal Y_s\in \mathfrak X(\pi(G_0))$.  Since there are only a finite number of possible choices for the sets $\partial(\{\x\})\subseteq \mathcal X_1\cup\ldots\cup \mathcal X_{s-1}=\mathcal Y$ for $\x\in \mathcal X_s$, we can likewise assume these sets are also fixed for $\mathcal X$.

\medskip

To summarize, we have arbitrary fixed values for $\Delta$, $\Delta'$, $\mathcal Y=\bigcup_{i=1}^{s-1}\mathcal X_i$, $\mathcal Y_s$, $Y_s$ and $\mathcal B_{\y}\subseteq \mathcal X_1\cup\ldots\cup \mathcal X_{s-1}$ for $\y\in \mathcal Y_s$ (we will later drop the subscript $\y$ when the corresponding set $\mathcal B$ is needed for only one $\y\in \mathcal Y_s$). Set $\mathcal E=\R^\cup \la\mathcal X_1\cup\ldots\cup\mathcal X_{s-1}\ra=\R\la \Delta'\ra$ and let $\pi:\R^d\rightarrow \mathcal E^\bot$ be the orthogonal projection. Let $\Omega$ consist of all maximal series decompositions $$\bigcup_{i=1}^{s-1}\mathcal X_i\cup \mathcal X_s\in \mathsf{Max}(\mathfrak X_s^*(G_0,\Delta))$$ having a $\Delta$-pure realization $\mathcal R=(\mathcal X_1\cup \{\mathbf v_1\},\ldots,\mathcal X_s\cup \{\mathbf v_s\})$ of $\mathcal X=\bigcup_{i=1}^{s-1}\mathcal X_i\cup \mathcal X_s$ with $(\mathcal X_1\cup \{\mathbf v_1\},\ldots,\mathcal X_{s-1}\cup \{\mathbf v_{s-1}\})$
a $\Delta'$-pure realization such that  $\pi(\mathcal X_s)=\mathcal Y_s$,  $\pi(X_s(k))=Y_s$ for all tuples $k$, and   $\partial(\{\x\})=\mathcal B_{\pi(\x)}$ for every $\x\in\mathcal X_s$.  Our goal is to show $\Omega$ is finite for the arbitrary choice of fixed values used to define $\Omega$.

\subsection*{Claim B} Suppose $\bigcup_{i=1}^{s-1}\mathcal X_i\cup \mathcal X_s\in \Omega$ and $\bigcup_{i=1}^{s-1}\mathcal X_i\cup \mathcal Z_s\in \Omega$. If $\x\in \mathcal X_s$ and $\z\in \mathcal Z_s$ with $\pi(\x)=\pi(\x)$, then $\bigcup_{i=1}^{s-1}\mathcal X_i\cup (\mathcal X_s\setminus\{\x\}\cup \{\z\})\in \Omega$.

\begin{proof} Let $\bigcup_{i=1}^{s-1}\mathcal X_i\cup \mathcal X_s\in \Omega$ and $\bigcup_{i=1}^{s-1}\mathcal X_i\cup \mathcal Z_s\in \Omega$, and
let $\mathcal R_\mathcal X=(\mathcal X_1\cup \{\mathbf v_1\},\ldots,\mathcal X_{s-1}\cup \{\mathbf v_{s-1}\},\mathcal X_s\cup \{\mathbf v_s\})$ and $\mathcal R_\mathcal Z=(\mathcal X_1\cup \{\mathbf w_1\},\ldots,\mathcal X_{s-1}\cup \{\mathbf w_{s-1}\},\mathcal Z_s\cup \{\mathbf w_s\})$ be $\Delta$-pure realizations of $\mathcal X$ and $\mathcal Z$, respectively, with $(\mathcal X_1\cup \{\mathbf v_1\},\ldots,\mathcal X_{s-1}\cup \{\mathbf v_{s-1}\})$ and $(\mathcal X_1\cup \{\mathbf w_1\},\ldots,\mathcal X_{s-1}\cup \{\mathbf w_{s-1}\})$ each $\Delta'$-pure realizations of $\mathcal Y=\mathcal X\setminus \mathcal X_s=\mathcal Z\setminus \mathcal Z_s$.
Since each $\mathbf w_j$ is a maximal element in the poset of half-spaces from $\mathcal R_\mathcal Z$, we can swap each $\mathbf w_j$ for $\mathbf v_j$ (for $j\in [1,s-1]$), allowing us to assume $\mathbf w_j=\mathbf v_j$ for all $j\in [1,s-1]$, so that  $\mathcal R_\mathcal Z=(\mathcal X_1\cup \{\mathbf v_1\},\ldots,\mathcal X_{s-1}\cup \{\mathbf v_{s-1}\},\mathcal Z_s\cup \{\mathbf w_s\})$.
 For $\x\in \mathcal X_s$,  there is some nonzero $y\in Y_s\subset \pi(G_0)$ such that $\pi(\x(i))=y$ for all $i$. Then $y=\pi(\overline y)$ for some $\overline y\in G_0\subseteq \Lambda$, meaning each $\x(i)=\overline y+\xi_{\x}(i)$ for some $\xi_\x(i)\in \Lambda\cap \mathcal E$ (as $\ker \pi=\mathcal E$ and $\x(i),y\in \Lambda$).
 As discussed after \eqref{lattice-splitting}, the value of $\Z\la X\la k\ra\ra=\Z\la Y(k)\ra+\Z\la X_k\ra=\Delta'+\Z\la X_s(k)\ra$ depends solely on the values of $\xi_\x(i)\mod \Delta'$ for $\x\in \mathcal X_s$, with distinct possibilities for the $\xi_\x(i)$ modulo $\Delta'$ giving rise to distinct lattices.
 Thus, since $\Delta=\Z\la X\la k\ra\ra$ and $\Delta'=\Z\la Y(k)\ra$ are fixed for all tuples $k$, it follows that there are $\xi_\y\in \Lambda\cap \mathcal E$ for $\y\in \mathcal Y_s$ such that $\xi_\x(i)\equiv \xi_\y\mod \Delta'$ for all $i$ and $\x\in \mathcal X_s$ with $\pi(\x)=\y$.
 Applying these same arguments to $\mathcal Z_s$ instead of $\mathcal X_s$, we likewise conclude that $\xi_\z(i)\equiv \xi_\y\mod \Delta'$ for all $i$ and $\z\in\mathcal Z_s$ with $\pi(\z)=\y$. As a result, if $\pi(\x)=\pi(\z)=\y\in \mathcal Y_s$, then $\Z\la (X\setminus \{x\}\cup \{z\})(k)\ra=\Z\la X(k)\ra=\Delta$ for all tuples $k$, as the values of $\Delta'$ and $\xi_\y$ modulo $\Delta'$ for $\y\in\mathcal Y_s$ completely determine the  lattice $\Delta$.

 Let $\x\in \mathcal X_s$ and $\z\in \mathcal Z_s$ with $\pi(\x)=\pi(\z)$. Then $\pi(\mathcal X_s\setminus \{\x\}\cup \{\z\})=\pi(\mathcal X_s)=\mathcal Y_s$ with
 $\mathcal R'_\mathcal X=(\mathcal X_1\cup \{\mathbf v_1\},\ldots,\mathcal X_{s-1}\cup \{\mathbf v_{s-1}\}, (\mathcal X_s\setminus \{\x\}\cup \{\z\})\cup \{\mathbf v_s\})$ a purely virtual Reay system over $G_0$, $\partial(\{\z\})=\partial(\{\x\})=\mathcal B_{\pi(\x)}=\mathcal B_{\pi(\z)}$ and $\pi(X_s\setminus \{x\}\cup\{z\})(k)=Y_s$ for all tuples $k$. In view of the conclusion of the previous paragraph, we see that $\mathcal R'_\mathcal X$ is $\Delta$-pure, while $(\mathcal X_1\cup \{\mathbf v_1\},\ldots,\mathcal X_{s-1}\cup \{\mathbf v_{s-1}\})$ is $\Delta'$-pure by assumption.
 Applying Proposition \ref{prop-finitary-MaxDecompChar} to the maximal decomposition $\bigcup_{i=1}^s\mathcal X_i$ and using that $\pi(\mathcal X_s\setminus \{\x\}\cup \{\z\})=\mathcal Y_s=\pi(\mathcal X_s)$ shows that the series decomposition $\bigcup_{i=1}^{s-1}\mathcal X_i\cup (\mathcal X_s\setminus\{\x\}\cup \{\z\})$ is maximal. It follows that $\bigcup_{i=1}^{s-1}\mathcal X_i\cup (\mathcal X_s\setminus\{\x\}\cup \{\z\})\in \mathfrak X^*_s(G_0,\Delta)$. It remains to show $\bigcup_{i=1}^{s-1}\mathcal X_i\cup (\mathcal X_s\setminus\{\x\}\cup \{\z\})\in \mathsf{Max}\big(\mathfrak X^*_s(G_0,\Delta)\big)$ in order to complete the claim.

\medskip

 Assume by contradiction that $\bigcup_{i=1}^{s-1}\mathcal X'_i\cup (\mathcal X_{s}\setminus \{\x\}\cup \{\z\})'$ is a maximal $\Delta$-pure series decomposition with $\bigcup_{i=1}^{s-1}\mathcal X_i\cup (\mathcal X_s\setminus\{\x\}\cup \{\z\})\prec_\cup^*\bigcup_{i=1}^{s-1}\mathcal X'_i\cup (\mathcal X_{s}\setminus \{\x\}\cup \{\z\})'$. Let $$\mathcal R'=(\mathcal X'_1\cup\{\mathbf v'_1\},\ldots,\mathcal X'_{s-1}\cup \{\mathbf v'_{s-1}\},(\mathcal X_s\setminus \{\x\}\cup \{\z\})'\cup \{\mathbf v'_s\})$$ be a $\Delta$-pure realization.
 Since each representative of a half-space $\y$ is also one for $\y'$, it follows from Proposition \ref{prop-orReay-BasicProps}.1 that  \be\label{kernelprimo}\ker \pi=\R^\cup\la \mathcal X_1\cup\ldots\cup\mathcal X_{s-1}\ra=
 \R^\cup\la \mathcal X'_1\cup\ldots\cup\mathcal X'_{s-1}\ra.\ee
 For $\y\in \mathcal X_s\setminus \{\x\}\cup\{\z\}$, the condition  $\y\subseteq \y'$ implies $\pi(\y)\subseteq \pi(\y')$ with both $\pi(\y')$ and $\pi(\y)$ one-dimensional half-spaces (which follows by applying the definition of oriented Reay system  to $\mathcal R_\mathcal X$, $\mathcal R_\mathcal Z$ and  $\mathcal R'$), whence
\be\label{frog-prime-hop}\pi(\x)=\pi(\z)=\pi(\z')\quad\und\quad \pi(\y)=\pi(\y')\;\mbox{ for all $\y\in \mathcal X_s\setminus\{\x\}$}.\ee

 Since $\mathcal R'$ is  $\Delta$-pure and $\mathcal E=\R^\cup \la \mathcal X_1\cup\ldots\cup \mathcal X_{s-1}\ra=\R^\cup \la \mathcal X'_1\cup\ldots\cup \mathcal X'_{s-1}\ra=\R\la (X\setminus X_s)'(k)\ra$ for all tuples $k$ (by Proposition \ref{prop-orReay-BasicProps}.1), we have
 \be\label{chop1}\Delta\cap \mathcal E=\Big(\Z\la (X\setminus X_s)'(k)\ra+\Z\la (X_s\setminus \{x\}\cup\{z\})'(k)\ra\Big)\cap \mathcal E=\Z\la (X\setminus X_s)'(k)\ra,\ee for all tuples $k$, where the final equality above follows since $\pi((X_s\setminus \{x\}\cup\{z\})'(k))$ is a linearly independent set of size $|\mathcal X_s|$ by (OR2) for $\mathcal R'$. Since $\mathcal R_\mathcal X=(\mathcal X_1\cup\{\mathbf v_1\},\ldots,\mathcal X_s\cup \{\mathbf v_s\})$ is $\Delta$-pure and
 $(\mathcal X_1\cup\{\mathbf v_1\},\ldots,\mathcal X_{s-1}\cup \{\mathbf v_{s-1}\})$ is $\Delta'$-pure, we have \be\label{chop2}\Delta\cap \mathcal E=\Big(\Z\la (X\setminus X_s)(k)\ra+\Z\la X_s(k)\ra\Big)\cap \mathcal E=\Z\la (X\setminus X_s)(k)\ra=\Delta'\ee for all tuples $k$, with the second equality above following since $\ker \pi=\mathcal E=\R^\cup \la \mathcal X_1\cup\ldots\cup \mathcal X_{s-1}\ra=\R\la (X\setminus X_s)(k)\ra$ (by Proposition \ref{prop-orReay-BasicProps}.1) with $\pi(X_s(k))$  a linearly independent set of size $|X_s|$ (by (OR2) for $\mathcal R_\mathcal X$). Combining \eqref{chop1} and \eqref{chop2}, we find that \be\nn\Delta'=\Delta\cap \mathcal E=\Z\la (X\setminus X_s)'(k)\ra\ee for all tuples $k$, ensuring that $(\mathcal X'_1\cup \{\mathbf v'_1\},\ldots,\mathcal X'_{s-1}\cup \{\mathbf v'_{s-1}\})$ is $\Delta'$-pure.

We have \begin{align}\label{legofrog}\Z\la Y_s\ra&=\Z\la \pi\big((X_s)(k)\big)\ra=\pi\Big(\Z\la (X\setminus X_s)(k)\ra+\Z\la X_s(k)\ra\Big)=\pi(\Delta)\\&=\pi\Big(\Z\la (X\setminus X_s)'(k)\ra+\Z\la (X_s\setminus \{x\}\cup \{z\})(k)\ra\Big)=\Z\la\pi\Big((X_s\setminus \{x\}\cup \{z\})'(k)\Big)\ra,\nn\end{align}
for all tuples $k$, with the first equality since $\bigcup_{i=1}^s\mathcal X_i\in \Omega$, the third since $\mathcal R_\mathcal X$ is $\Delta$-pure, the fourth since $\mathcal R'$ is $\Delta$-pure, and the second and fifth since $\ker\pi=\R^\cup \la \mathcal X_1\cup\ldots\cup \mathcal X_{s-1}\ra=\R\la (X\setminus X_s)(k)\ra=\R\la (X\setminus X_s)'(k)\ra=\R^\cup \la \mathcal X'_1\cup\ldots\cup \mathcal X'_{s-1}\ra$.
 In view of \eqref{frog-prime-hop}, we have $\pi(\z')=\pi(\x)$, which is a one-dimensional ray \emph{positively} spanned by an element from $Y_s$ (per definition of $\Omega$). Likewise, each $\y\in  \mathcal X_s\setminus\{\x\}$ has $\pi(\y')=\pi(\y)$ being a one-dimensional ray \emph{positively} spanned by an element from $Y_s$. Thus, since $Y_s$ is linearly independent, it follows that only way \eqref{legofrog} can hold is if $\pi\big((X_s\setminus\{x\}\cup \{z\})'(k)\big)=Y_s$ for all tuples $k$.

 As a result,  $\mathcal R'$ now satisfies all the same hypotheses as  $\mathcal R_\mathcal X$ and $\mathcal R_\mathcal Z$ needed to define analogously the elements $\xi_{\y'}(i)$ for $\y\in \bigcup_{i=1}^{s-1}\mathcal X_i\cup (\mathcal X_s\setminus \{\x\}\cup \{\z\})$  for $\mathcal R'$ as they were defined for $\mathcal R_\mathcal X$ and $\mathcal R_\mathcal Z$. By \eqref{frog-prime-hop} and the discussion after \eqref{lattice-splitting}, it again follows that \be\label{align-xi}\xi_{\y'}(i)\equiv \xi_{\pi(\y')}=\xi_{\pi(\y)}\mod \Delta'\ee for all $\y\in \mathcal X_s\setminus\{\x\}\cup\{\z\}$ and $i\geq 1$

 Suppose $\y=\y'$ for all  half-spaces $\y\in \bigcup_{i=1}^{s-1}\mathcal X_i\cup (\mathcal X_s\setminus\{\x\})$. Then we must have $\z\subset \z'$ since $\bigcup_{i=1}^{s-1}\mathcal X_i\cup (\mathcal X_s\setminus\{\x\}\cup \{\z\})\prec_\cup^*\bigcup_{i=1}^{s-1}\mathcal X'_i\cup (\mathcal X_{s}\setminus \{\x\}\cup \{\z\})'$. In such case, since $\pi(\z')=\pi(\z)$ (as follows by \eqref{frog-prime-hop}) and $\partial(\{\z'\})\subseteq \bigcup_{i=1}^{s-1}\mathcal X'_i=\bigcup_{i=1}^{s-1}\mathcal X_i$ (as $\y=\y'$ for all other half-spaces), it follows that $\bigcup_{i=1}^{s-1}\mathcal X_i\cup (\mathcal Z_s\setminus \{\z\}\cup\{\z'\})$ is a series decomposition.
  Since $\pi(\mathcal Z_s\setminus \{\z\}\cup \{\z'\})=\pi(\mathcal Z_s)$ (by \eqref{frog-prime-hop}), it follows from  Proposition \ref{prop-finitary-MaxDecompChar} applied to the maximal series decomposition $\bigcup_{i=1}^{s-1}\mathcal X_i\cup \mathcal Z_s$, and then to $\bigcup_{i=1}^{s-1}\mathcal X_i\cup (\mathcal Z_s\setminus \{\z\}\cup\{\z'\})$, that $\bigcup_{i=1}^{s-1}\mathcal X_i\cup (\mathcal Z_s\setminus \{\z\}\cup\{\z'\})$ is maximal.
   Moreover, $\bigcup_{i=1}^{s-1}\mathcal X_i\cup (\mathcal Z_s\setminus \{\z\}\cup\{\z'\})$ is $\Delta$-pure since $\mathcal R_\mathcal Z$ is $\Delta$-pure with $\xi_{\z'}(i)\equiv \xi_{\pi(\z')}=\xi_{\pi(\z)}\mod \Delta'$ by \eqref{align-xi}. Thus $\bigcup_{i=1}^{s-1}\mathcal X_i\cup \mathcal Z_s\prec_{\cup}^*\bigcup_{i=1}^{s-1}\mathcal X_i\cup (\mathcal Z_s\setminus \{\z\}\cup\{\z'\})$, which contradicts that $\bigcup_{i=1}^{s-1}\mathcal X_i\cup \mathcal Z_s\in \mathsf{Max}\big(\mathfrak X_{s}^*(G_0,\Delta)\big)$.
 So we can instead assume $\y_0\subset \y_0'$ for some  $\y_0\in \bigcup_{i=1}^{s-1}\mathcal X_i\cup (\mathcal X_s\setminus\{\x\})$.

 Let $\vec u_\x=(u_1,\ldots,u_t)$. Since $\partial(\{\x\})\subseteq \mathcal X_1\cup\ldots\cup \mathcal X_{s-1}$ encases $-\vec u^\triangleleft_\x$, we have $u_1,\ldots,u_{t-1}\in \R^\cup\la \partial(\{\x\})\ra\subseteq \R^\cup \la \mathcal X_1\cup\ldots\cup \mathcal X_{s-1}\ra=\R^\cup \la \mathcal X'_1\cup\ldots\cup \mathcal X'_{s-1}\ra$ (by \eqref{kernelprimo}).  Thus Proposition \ref{prop-finitary-basics}.3 applied to a realization of $\bigcup_{i=1}^{s-1}\mathcal X'_i$ implies that there is a subset $\mathcal C\subseteq \mathcal X'_1\cup\ldots\cup \mathcal X'_{s-1}$ that minimally encases $-\vec u_\x^\triangleleft$, allowing us to define a half-space $\x'$ with $\partial(\{\x'\})=\mathcal C$, $\vec u_{\x'}=\vec u_\x$ and $\x'(i)=\x(i)$ for all $i$, so $\overline{\x'}=\R^\cup\la \mathcal C\ra+\R_+u_t$.
 Then $\pi(\x')=\pi(\z')=\pi(\x)=\R_+\pi(u_t)$ is a one-dimensional half-space (as $\ker \pi=\R^\cup\la  \mathcal X_1\cup\ldots\cup \mathcal X_{s-1} \ra$) with $$\pi((\mathcal X_s\setminus\{\x\}\cup \{\z\})')=\pi(\mathcal X_s)=\pi(\mathcal X'_s)=\mathcal Y_s$$
 by \eqref{frog-prime-hop}. Combined with  \eqref{kernelprimo} and Proposition \ref{prop-finitary-MaxDecompChar}, we find
   $\bigcup_{i=1}^{s-1}\mathcal X'_i\cup \mathcal X'_s$ is a  maximal series decomposition since $\bigcup_{i=1}^{s-1}\mathcal X'_i\cup (\mathcal X_s\setminus \{\x\}\cup \{\z\})'$ is a maximal series decomposition.
   We have  $\xi_{\z'}(i)\equiv \xi_{\pi(\z')}=\xi_{\pi(\z)}=\xi_{\pi(\x)}\equiv \xi_\x(i)\mod \Delta'$ for all $i$ by \eqref{align-xi} and \eqref{frog-prime-hop}. Thus, since $\x'(i)=\x(i)$ with  $\bigcup_{i=1}^{s-1}\mathcal X'_i\cup (\mathcal X_s\setminus\{\x\}\cup \{\z\})'$ $\Delta$-pure, it follows that  $\bigcup_{i=1}^{s-1}\mathcal X'_i\cup \mathcal X'_s$ is also $\Delta$-pure. We have now established that $\bigcup_{i=1}^s\mathcal X'_i\in \mathfrak X_s^*(G_0,\Delta)$.

 Since $\mathcal C\subseteq \mathcal X'_1\cup\ldots\cup\mathcal X'_{s-1}$ and $\partial(\{\x\})\subseteq \mathcal X_1\cup\ldots\cup\mathcal X_{s-1}$ both minimally encase $-\vec u^\triangleleft_{\x}=-\vec u^\triangleleft_{\x'}$, it follows from Proposition \ref{prop-finitary-finiteabove-align} and Proposition \ref{prop-orReay-minecase-char}.1 that $(\partial(\{\x\})')^*=\mathcal C$, whence $\partial(\{\x\})'\subseteq \darrow \mathcal C$, ensuring \begin{align*}\partial(\x)=\R^\cup \la \partial(\{\x\})\ra&\subseteq \R^\cup \la \partial(\{\x\})'\ra\subseteq \R^\cup\la \darrow \mathcal C\ra=\R^\cup \la \darrow\partial(\{\x'\})\ra=\partial(\x')\quad\und\\ \partial(\{\x\})\cap \x=\C^\cup (\partial(\{\x\}))&\subseteq \C^\cup (\partial(\{\x\})')\subseteq \C^\cup (\darrow \mathcal C)=\C^\cup(\darrow \partial(\{\x'\}))=\partial(\x')\cap \x',\end{align*} with the first inclusion in both lines above since $\y\subseteq \y'$ for all $\y\in \partial(\{\x\})$.
     Combined with the fact that $u_t$ is representative for both $\x$ and $\x'$, it follows that $\x\subseteq \x'$ with every representative for $\x$ one for $\x'$.  Since $\bigcup_{i=1}^{s-1}\mathcal X_i\cup (\mathcal X_s\setminus\{\x\}\cup \{\z\})\prec_\cup^*\bigcup_{i=1}^{s-1}\mathcal X'_i\cup \big((\mathcal X_{s}\setminus \{\x\})'\cup \{\z'\}\big)$), we have $\y\subseteq \y'$ and $\y^\circ\subseteq (\y')^\circ$ for all $\y\in \bigcup_{i=1}^{s-1}\mathcal X_i\cup (\mathcal X_s\setminus \{\x\})$. It follows that  $\bigcup_{i=1}^{s-1}\mathcal X_i\cup \mathcal X_s\prec^*_\cup \bigcup_{i=1}^{s-1}\mathcal X'_i\cup \mathcal X'_s$, with the relation strict since $\y_0\subset \y_0'$ for some half-space $\y_0\in \bigcup_{i=1}^{s-1}\mathcal X_i\cup (\mathcal X_s\setminus\{\x\})$, which contradicts  that $\bigcup_{i=1}^{s-1}\mathcal X_i\cup \mathcal X_s\in \mathsf{Max}(\mathfrak X_s^*(G_0,\Delta))$, completing Claim B.
 \end{proof}

We continue assuming $\bigcup_{i=1}^s\mathcal X_i\in \Omega$ with all notation as defined earlier. In particular, $\mathcal R=(\mathcal X_1\cup \{\mathbf v_1\},\ldots,\mathcal X_s\cup \{\mathbf v_s\})$ is $\Delta$-pure realization of $\mathcal X=\bigcup_{i=1}^s\mathcal X_i$.
Let $\y\in\mathcal Y_s$ be arbitrary. Then there is a support set $\mathcal B=\mathcal B_\y\subseteq \mathcal X_1\cup\ldots\cup\mathcal X_{s-1}=\mathcal Y$ and some representative $y\in Y_s\subseteq \pi(G_0)$ for $\y$, such that, if  $\x\in \mathcal X_s$ is any potential half-space with $\pi(\x)=\y$,  then $\partial(\{\x\})=\mathcal B$  and $\pi(\x(i))=y\in \mathcal E^\bot$ for all $i\geq 1$. Note $y\neq 0$ since it is a representative of a half-space from $\mathcal Y_s\in \mathfrak X(\pi(G_0))$.
Let   $$\pi_\mathcal B:\R^d\rightarrow \R^\cup\la \mathcal B\ra^\bot$$ be the orthogonal projection.
Letting $\vec u_\x=(u^{(1)}_\x,\ldots,u_\x^{(t_\x)})$,
we have (by Proposition \ref{prop-VReay-Lattice}) \be\label{repso}\x(i)=(a_{i,\x}^{(1)}u_\x^{(1)}+\ldots+a_{i,\x}^{(t_\x-1)}u_\x^{(t_\x-1)}
+w_{i,\x})+\xi_\x+y\ee for some $a_{i,\x}^{(j)}>0$, \ $w_{i,\x}\in \R^\cup \la \mathcal B\ra$ and $\xi_\x\in \mathcal E\cap \R^\cup\la \mathcal B\ra^\bot$, with $a_{i,\x}^{(j)}\rightarrow \infty$ for $j\in [1,t_\x-1]$, $a_{i,\x}^{(j)}\in o(a_{i,\x}^{(j-1)})$ for $j\in [2,t_\x-1]$, and $\|w_{i,\x}\|\in o(a_{i,\x}^{(t_\x-1)})$. Note $\xi_\x+y$ is a representative for the half-space $\x$. Since $\xi_\x\in \mathcal E\cap \R^\cup\la \mathcal B\ra^\bot$ and  $y\in \mathcal E^\bot$ is nonzero, it follows that distinct values for $\xi_\x$ determine distinct half-spaces $\x$. Also, we have $\pi_{\mathcal B}(\x(i))=\xi_\x+y\in \pi_\mathcal B(G_0)\subseteq \pi_\mathcal B(\Lambda)$, while $\pi_{\mathcal B}(\Lambda)$ is  a lattice by Proposition \ref{Prop-lattice-homoIm}. The set $\mathcal B$ is empty precisely when $t_\x=1$. For each $\z\in \mathcal B$, let $$\mathcal B_\z=\mathcal B\setminus\{\z\}\cup \partial(\{\z\})$$ and let $\pi_\z:\R^d\rightarrow\R^\cup \la \mathcal B_\z\ra^\bot$ be the orthogonal projection.

If, for each $\y\in \mathcal Y_s$, there are only a finite number of possible half-spaces $\x $ with $\pi(\x)=\y$ that occur for the  $\mathcal X\in \Omega$, then $\Omega$ will be finite, as desired. Therefore we can assume there is some $\y\in \mathcal Y_s$ for which this fails, instead having an infinite sequence of distinct half-spaces $\{\x_j\}_{j=1}^\infty$ with $\pi(\x_j)=\y$ for all $j\geq 1$, and each $\x_j$ lying in some
$\mathcal X_s^{(j)}$ with $\bigcup_{i=1}^{s-1}\mathcal X_i\cup \mathcal X_s^{(j)}\in \Omega$. We use the notation of the previous paragraph with each $\x_j$ except that we replace everywhere the subscript $\x=\x_j$ by $j$
(to lighten notation some). We also write $\x_j(i):=\x(i,j)$ for $i,\,j\geq 1$, to emphasize the dependence of the representation in \eqref{repso} on both the parameters $i$ and $j$.
Let $$\x_\mathcal B(i,j)=a_{i,j}^{(1)}u_j^{(1)}+\ldots+a_{i,j}^{(t_j-1)}u_j^{(t_j-1)}
+w_{i,j}\in \R^\cup \la \mathcal B \ra, $$
so $\x(i,j)=\x_\mathcal B(i,j)+\xi_j+y$.

In view of Claim B, we can assume the other half-spaces in $\mathcal X^{(j)}_s$ remain fixed and consider the half-space $\x_j\in \mathcal X_s$, corresponding to $\y$,  as varying. Then, each choice of a half-space $\x_j$ with $\pi(\x_j)=\y$, gives rise to a set $\mathcal X_s^{(j)}$ with  $\bigcup_{i=1}^{s-1}\mathcal X_i\cup \mathcal X_s^{(j)}\in  \Omega$ and $\mathcal X^{(j)}:=\bigcup_{i=1}^{s-1}\mathcal X_i\cup \mathcal X_s^{(j)}$, where $\mathcal X^{(j)}\setminus\{\x_j\}$ is a fixed set. Note, we will for most of the argument drop the super-scripts $(j)$ as this information will be irrelevant to all but the final arguments.

Since the half-spaces $\x_j$ are distinct with $\partial(\{\x_j\})=\mathcal B$  fixed,  the  representatives  $\xi_j+y\in \pi_\mathcal B(G_0)\subseteq \pi_\mathcal B(\Lambda)$ are also all distinct, and since these are lattice points, this means the sequence $\{\xi_j\}_{j=1}^\infty$ is unbounded. Thus, by passing to a subsequence, we can assume $$\|\xi_j\|\rightarrow \infty\quad\mbox{ with }\quad \|\xi_j\|\geq 1\quad\mbox{ for all $j$},$$ and that $\{\xi_j\}_{j=1}^\infty$ is an asymptotically filtered sequence of terms with fully unbounded limit.

\subsection*{Claim C} If $\vec v=(v_1,\ldots,v_{t})\in G_0^{\mathsf{lim}}$ with $v_1,\ldots,v_t\in \R^\cup\la \mathcal B\ra$, then $\mathcal B$ encases $-\vec v$.

\begin{proof}The set $\mathcal B=\partial(\{\x\})$ minimally encases the fully unbounded (or trivial) limit $-\vec u_\x^\triangleleft$ (by Proposition \ref{prop-VReay-Lattice}), and since $\mathcal B\subseteq \mathcal X_1\cup \ldots\cup \mathcal X_{s-1}$, the encasement is urbane with $\mathcal B$ a support set. Proposition \ref{prop-VReay-modularCompletion} implies there is a virtual Reay system $\mathcal R_\mathcal B=(\mathcal Z_1\cup \{\z_1\},\ldots,\mathcal Z_{s'}\cup \{\z_{s'}\})$ over $G_0$ with $\darrow \mathcal B=\bigcup_{i=1}^{s'}\mathcal Z_i$. If $\mathcal B=\emptyset$, then $\mathcal R_\mathcal B$ is the trivial (empty) virtual Reay system. Otherwise,  $\vec u_\x^\triangleleft$ is fully unbounded, and the strict truncation of any $\vec u_\y$ with $\y\in \darrow \mathcal B\subseteq \mathcal X_1\cup\ldots\cup \mathcal X_{s-1}$ is also either trivial or fully unbounded (by Proposition \ref{prop-VReay-Lattice}). Thus, per the comments above Proposition \ref{prop-VReay-modularCompletion}, $\mathcal R_\mathcal B$ is purely virtual. Claim C now follows from Proposition \ref{prop-finitary-basics}.3 applied to $\mathcal R_\mathcal B$.\end{proof}

We aim to contradict the maximality of $\bigcup_{i=1}^s\mathcal X_i$ by  using the infinite sequence of distinct half-spaces  $\{\x_j\}_{j=1}^\infty$ to construct a new,  strictly larger half-space $\mathbf u$. To define $\mathbf u$, we will need to construct an appropriate representation sequence $\{\mathbf u(j)\}_{j=1}^\infty$. We define $\mathbf u(j)=\x\big(i(j),j\big)\in G_0$, where $i:\Z_+\setminus\{0\}\rightarrow \Z_+$ is a function that associates to each $j\geq 1$ a sufficiently large index $i$, constructed as follows. If $\mathcal B=\emptyset$, then simply take $i(j)=j$ for all $j\geq 1$.
Otherwise, for each fixed $j\geq 1$,  the sequence $\{\x_\mathcal B(i,j)\}_{i=1}^\infty$ is an asymptotically  filtered sequence of terms  with  limit $\vec u_j^\triangleleft$ such that  $-\vec u_j^\triangleleft$ is fully unbounded (by Proposition \ref{prop-VReay-Lattice}) and  minimally encased by $\mathcal B=\partial(\{\x_j\})$. Thus Lemma \ref{lemma-radconv-super}.1 implies that $\|\pi_\z(\x_\mathcal B(i,j))\|\rightarrow \infty$ for all $\z\in \mathcal B$ (as $i\rightarrow\infty$). Since $\xi_j$ is a fixed element, it is thus possible to choose $i(j)$ to be sufficiently large so that \be\label{lintin}\|\pi_\z\big(\x_\mathcal B\big(i,j\big)\big)\|\geq 2^j\|\xi_j\|\geq 2^j\quad\mbox{for all $i\geq i(j)$}\ee and all $\z\in \mathcal B$. For $j\geq 1$ and $\mathcal B\neq \emptyset$, define $i(j)$ to be any sufficiently large index such that \eqref{lintin} holds. Note that, even if we pass to a subsequence of $\{\x_j\}_{j=1}^\infty$, then \eqref{lintin} remains true under the new re-indexing of the remaining terms from $\{\x_j\}_{j=1}^\infty$.
With the function $i:\Z_+\setminus \{0\}\rightarrow \Z_+$ fixed, we set
$$\mathbf u_\mathcal B(j)=\x_\mathcal B\big(i(j),j\big)\in \R^\cup\la \mathcal B\ra \quad\und\quad
\mathbf u(j)=\x\big(i(j),j\big)=\mathbf u_\mathcal B(j)+\xi_j+y\in G_0.$$
If $\mathcal B=\emptyset$, then $\mathbf u_\mathcal B(j)=0$  and $\mathbf u(j)=\xi_j+y$  for all $j\geq 1$. Otherwise, \eqref{lintin} ensures that \be \label{order-stuff}\|\pi_\z(\mathbf u_\mathcal B(j))\|\rightarrow \infty\quad\und \quad\|\xi_j\|\in o(\|\pi_\z(\mathbf u_\mathcal B(j))\|)\quad \mbox{ for all $\z\in \mathcal B$}.\ee In particular, since $\|\mathbf u_\mathcal B(j)\|\geq \|\pi_\z(\mathbf u_\mathcal B(j))\|$, we have
\be\label{order-more} \|\mathbf u_\mathcal B(j)\|\rightarrow \infty\quad\und \quad\|\xi_j\|\in o(\|\mathbf u_\mathcal B(j)\|)\quad \mbox{ when  $\mathcal B\neq\emptyset$.}\ee

If $\mathcal B\neq \emptyset$, then \eqref{order-more} gives $\|\mathbf u_\mathcal B(j)\|\rightarrow \infty$. Thus, by passing to an appropriate subsequence of $\{\x_j\}_{j=1}^\infty$, we can assume $\{\mathbf u_\mathcal B(j)\}_{j=1}^\infty$ is an asymptotically filtered sequence with complete fully unbounded limit $\vec v=(v_1,\ldots,v_t)$. Let $$\mathbf u_\mathcal B(j)=\x_\mathcal B\big(i(j),j\big)=\alpha_j^{(1)}v_1+\ldots+\alpha_j^{(t)}v_t+\varepsilon_j$$ be the representation of
$\mathbf u_\mathcal B(j)$ as an asymptotically filtered sequence with limit $\vec v$,
so $\alpha_j^{(r)}>0$ and $\alpha_j^{(r)}\rightarrow \infty$ for $r\in [1,t]$, \ $\alpha_j^{(r)}\in o(\alpha_j^{(r-1)})$ for $r\in [2,t]$, and $\|\varepsilon_j\|$ is bounded.
Since $\mathbf u_\mathcal B(j)\in \R^\cup\la \mathcal B\ra$ for all $j$, we must have $v_1,\ldots,v_t\in \R^\cup \la \mathcal B\ra$.  When $\mathcal B=\emptyset$, we instead set $\vec v$ to be the empty tuple, so $t=0$ and $\varepsilon_j=0$ for all $j$.

\subsection*{Claim D} $-\vec v$ is encased by $\mathcal B$.

\begin{proof} If $\mathcal B=\emptyset$, then the claim is trivial, so we may assume $\mathcal B\neq \emptyset$.
Note that $\mathbf u(j)=\mathbf u_\mathcal B(j)+\xi_j+y\in G_0$ for all $j$.
Since $\alpha_j^{(r)}\rightarrow \infty$ for all $r\in [1,t]$, we have the constant sequence  $y\in o(\alpha_j^{(r)})$ for all $r\in [1,t]$.
Since $\alpha_j^{(1)}\sim \|\mathbf u_\mathcal B(j)\|$, \eqref{order-more} implies that  $\|\xi_j\|\in o(\alpha_j^{(1)})$.
 Thus there is a maximal index $t'\in [1,t]$ such that
 $$\|\xi_j\|\in o(\alpha_j^{(r)})\mbox{ for all $r\in [1,t']$}.$$ In such case, $\{\mathbf u(j)\}_{j=1}^\infty$ is an asymptotically filtered sequence of terms from $G_0$ with fully unbounded limit $\vec v'=(v_1,\ldots,v_{t'})$, in which case $-\vec v'$ is encased by $\mathcal B$ in view of Claim C (note we need  $\mathbf u(j)\in G_0$ to use Claim C, so we cannot directly apply Claim C to $\{\mathbf u_\mathcal B(j)\}_{j=1}^\infty$).
 Let $\mathcal B'\preceq \mathcal B$ be such that $\mathcal B'$ minimally encases $-\vec v'$. Since $\mathcal B'\preceq \mathcal B\subseteq \mathcal X_1\cup\ldots\cup \mathcal X_{s-1}$, the encasement of $-\vec v'$ by $\mathcal B'$ is urbane with $\mathcal B'$ a support set. Thus we can apply Proposition \ref{prop-orReay-minecase-char} and let
 $$\emptyset =B'_0\prec \mathcal B'_1\prec \ldots\prec\mathcal B'_\ell=\mathcal B'\preceq  \mathcal B$$ be the support sets and $1=r_1<\ldots<r_{\ell}<r_{\ell+1}=t'+1$ be the indices given by Proposition \ref{prop-orReay-minecase-char} applied to the encasement of $-\vec v'$ by $\mathcal B'$.
 If $\mathcal B'=\mathcal B$, then $\mathcal B=\mathcal B'$ encases $-\vec v$ in view of $v_1,\ldots,v_t\in \R^\cup\la \mathcal B\ra$ and Proposition \ref{prop-orReay-minecase-char}.2 (particularly, part (c)), yielding the claim. Therefore we may assume $\mathcal B'\prec \mathcal B$, whence  $\mathcal B'\preceq \mathcal B_\z$ for some $\z\in \mathcal B$, implying $\mathcal B'\subseteq \darrow \mathcal B_\z$. Hence $\pi_\z(v_i)=0$ for all $i\in [1,t']$.

If $\pi_\z(v_i)=0$ for all $i\in [1,t]$, then $\|\pi_\z(\mathbf u_\mathcal B(j))\|=\|\pi_\z(\varepsilon_j)\|$, which is bounded as $\vec v$ is a \emph{complete} fully unbounded limit for $\mathbf u_\mathcal B(j)$, contradicting \eqref{order-stuff}. Therefore there must be some minimal $t''\in [t'+1,t]$ such that $\pi_\z(u_{t''})\neq 0$. But then $\|\pi_\z(\mathbf u_\mathcal B(j))\|\in \Theta(\alpha_j^{(t'')})$, which combined with  $\|\xi_j\|\in o(\|\pi_\z(\mathbf u_\mathcal B(j))\|)$ from \eqref{order-stuff} yields  $\|\xi_j\|\in o(\alpha_j^{(t'')})$. In consequence, since $\alpha_j^{(r)}\in o(\alpha_j^{(r-1)})$ for all $r\in [2,t]$, we have $\|\xi_j\|\in o(\alpha_j^{(r)})$ for all $r\in [1,t'']$, implying by the maximality in the definition of $t'$ that $t'\geq t''$, contradicting that $t''\in [t'+1,t]$ as shown above. This establishes Claim D.\end{proof}

If $\mathcal B\neq \emptyset$, then Claim D implies that   $\{-\mathbf u_\mathcal B(j)\}_{j=1}^\infty$ is an asymptotically filtered sequence of terms with complete fully unbounded limit $-\vec v$ encased by $\mathcal B$, while $\|\pi_\z(\mathbf u_\mathcal B(j))\|$ is unbounded by \eqref{order-stuff} for all $\z\in \mathcal B$. Thus Lemma \ref{lemma-radconv-super}.2 implies that $\mathcal B$ minimally encases $-\vec v$, a fact which is trivially true when $\mathcal B=\emptyset$ as well. Since $\mathcal B\subseteq \mathcal X_1\cup\ldots\cup\mathcal X_{s-1}$, the encasement is urbane.  Let
 $$\emptyset =\mathcal B_0\prec \mathcal B_1\prec \ldots\prec\mathcal B_\ell=\mathcal B$$ be the support sets and $1=r_1<\ldots<r_{\ell}<r_{\ell+1}=t+1$ be the indices given by Proposition \ref{prop-orReay-minecase-char} applied to the encasement of $-\vec v$ by $\mathcal B$.
By \eqref{order-stuff},  we have  $\|\xi_j\|\in o(\|\pi_\z(\mathbf u_\mathcal B(j))\|)$ for all $\z\in \mathcal B$, and thus also $\|\xi_j+y\|\in o(\|\pi_\z(\mathbf u_\mathcal B(j))\|)$ for all $\z\in \mathcal B$ as $\|\pi_\z(\mathbf u_\mathcal B(j))\|\rightarrow \infty$ by \eqref{order-stuff}. Thus Lemma \ref{lemma-radconv-super}.3 (applied with $j$ replacing $i$, with $x_i$ taken to be $\mathbf u_\mathcal B(j)$, and with $y_i$ taken to be $\xi_j+y$) implies that $\{\mathbf u(j)\}_{j=1}^\infty$ is an asymptotically filtered sequence of terms from $G_0$ with fully unbounded limit $(v_1,\ldots,v_{r_\ell})$ (after discarding the first few terms), when $\mathcal B\neq \emptyset$.
However, when $\mathcal B=\emptyset$, we have $\mathbf u(j)=\xi_j+y$, which is also an asymptotically filtered sequence with fully unbounded limit, which we can set equal to $(v_1,\ldots,v_{r_\ell})$.
 By passing to a subsequence of $\{\mathbf u(j)\}_{j=1}^\infty$, we may assume $\{\mathbf u(j)\}_{j=1}^\infty$ is an asymptotically filtered sequence of terms from $G_0$ with \emph{complete} fully unbounded limit $\vec w=(v_1,\ldots,v_{r_\ell},w_{1},\ldots,w_{n})$.

Since $\{\mathbf u(j)\}_{j=1}^\infty$ is  asymptotically  filtered  with complete \emph{fully unbounded} limit $\vec w$, so too is the sequence $\{\mathbf u(j)-y\}_{j=1}^\infty$ obtained by translating all terms by a fixed constant (after discarding the first few terms). As a result, since  $\mathbf u(j)-y=\x_\mathcal B((i(j),j)+\xi_j\in \R^\cup \la \mathcal X_1\cup\ldots\cup \mathcal X_{s-1}\ra=\mathcal E$, we conclude that $v_1,\ldots,v_{r_\ell},w_1,\ldots,w_n\in \mathcal E$. Consequently, since $\{\mathbf u(j)\}_{j=1}^\infty$ is an asymptotically filtered sequence of \emph{terms from $G_0$} with fully unbounded limit $\vec w$, it follows from Proposition \ref{prop-finitary-basics}.3 (applied to $(\mathcal X_1\cup \{\mathbf v_1\},\ldots,\mathcal X_{s-1}\cup \{\mathbf v_{s-1}\})$) that $\mathcal Y=\mathcal X_1\cup\ldots\cup \mathcal X_{s-1}$ encases $-\vec w$. Let $\mathcal C\subseteq \mathcal Y$ be a subset which minimally encases $-\vec w$. Then, since $\mathcal B\subseteq \mathcal Y$ minimally encases $-(v_1,\ldots,v_{r_\ell})$, it follows in view of Proposition \ref{prop-orReay-minecase-char} that $\mathcal B\preceq \mathcal C$, with equality only possible if $w_1,\ldots, w_n\in \R^\cup \la \mathcal B\ra$. However, if $w_1,\ldots, w_n\in \R^\cup \la \mathcal B\ra$, then  $\{\pi_\mathcal B(\mathbf u(j))\}_{j=1}^\infty$ will be a bounded sequence (as $\vec w$ is a \emph{complete} fully unbounded limit and $v_1,\ldots,v_{r_\ell}\in \R^\cup \la \mathcal B\ra$), contradicting that $\pi_\mathcal B(\mathbf u(j))=\pi_\mathcal B\big(\x((i(j),j))\big)=\xi_j+y$ with $\|\xi_j\|\rightarrow \infty$. Therefore we conclude that $$\mathcal B\prec \mathcal C.$$ Let $\pi_\mathcal C:\R^d\rightarrow \R^\cup\la \mathcal C\ra^\bot$ be the orthogonal projection.

Now $\mathbf u(j)=\mathbf u_\mathcal B(j)+\xi_j+y$ with $\mathbf u_\mathcal B(j)\in \R^\cup \la \mathcal B\ra\subseteq \R^\cup\la \mathcal C\ra$ and  $y\in  \mathcal E^\bot=\R^\cup\la \mathcal Y\ra^\bot\subseteq \R^\cup \la \mathcal C\ra^\bot$. Thus  $\pi_\mathcal C(\mathbf u(j))=\pi_\mathcal C(\xi_j)+y$ with $\{\pi_\mathcal C(\mathbf u(j))\}_{j=1}^\infty$ a bounded sequence (as $\vec w$ is a \emph{complete} fully unbounded limit) of lattice points from $\pi_\mathcal C(G_0)$ (in view of Proposition \ref{Prop-lattice-homoIm}).
By passing once more to a subsequence, we can thus assume $\pi_\mathcal C(\mathbf u(j))=\pi_\mathcal C(\xi_j)+y$ is constant (as a bounded set of lattice points is finite), say with $\pi_\mathcal C(\xi_j)=\xi\in \R^\cup \la\mathcal C\ra^\bot$ for all $j$. As $y\neq 0$ with $y\in \mathcal E^\bot$ and $\xi_j\in \mathcal E$, we have $\xi+y\neq 0$.
Let $$\mathbf u(j)=\beta_j^{(1)}v_1+\ldots+\beta_j^{(r_\ell)}v_{r_\ell}+\gamma_j^{(1)}w_1+\ldots+\gamma_j^{(n)}w_n
+z_j$$
be the representation of $\{\mathbf u(j)\}_{j=1}^\infty$ as an asymptotically filtered sequence with complete fully unbounded limit $\vec w$.
Since $\pi_{\mathcal C}(\mathbf u(j))=\xi+y\neq 0$, by passing to a subsequence
we can assume $\{\mathbf u(j)\}_{j=1}^\infty$ is an asymptotically filtered sequence with $(v_1,\ldots,v_{r_\ell},w_1,\ldots,w_n,w_{n+1},\ldots,w_{n+r})$ as limit, for some $w_{n+1},\ldots,w_{n+r-1}\in \R^\cup \la \mathcal C\ra$ and $w_{n+r}\notin \R^\cup \la \mathcal C\ra$, say with representation
$$\mathbf u(j)=\beta_j^{(1)}v_1+\ldots+\beta_j^{(r_\ell)}v_{r_\ell}+\gamma_j^{(1)}w_1+\ldots+
\gamma_j^{(n)}w_n
+\gamma_j^{(n+1)}w_{n+1}+\ldots+\gamma_j^{(n+r)}w_{n+r}+z'_j.$$
Then $\xi+y=\pi_\mathcal C(\mathbf u(j))=\gamma_j^{(n+r)}\pi_\mathcal C(w_{n+r})+\pi_\mathcal C(z'_j)$ is constant. Thus, since $\|z'_j\|\in o(\gamma_j^{(n+r)})$, we cannot have $\gamma_j^{(n+r)}\rightarrow 0$, for this would imply  $\xi+y=\pi_\mathcal C(\mathbf u(j))\rightarrow 0$, contradicting that $\xi+y\neq 0$. Hence, since $\{\gamma_j^{(n+r)}\}_{j=1}^\infty$ is bounded, we must have $\gamma_j^{(n+r)}\rightarrow a$ for some $a>0$, and now $\gamma_j^{(n+r)}\in o(\gamma_j^{(n+r-1)})$ ensures that $\gamma_j^{n+r-1}\rightarrow \infty$, which is only possible if $r=1$. Therefore
$$\mathbf u(j)=\beta_j^{(1)}v_1+\ldots+\beta_j^{(r_\ell)}v_{r_\ell}+\gamma_j^{(1)}w_1+\ldots+
\gamma_j^{(n)}w_n
+\gamma_j^{(n+1)}w_{n+1}+z'_j$$
with $\gamma_j^{(n+1)}\pi_\mathcal C(w_{n+1})+\pi_\mathcal C(z'_j)=\xi+y$ for all $j\geq 1$. Since $\|z'_j\|\in o(\gamma_j^{(n+1)})$, we have $\|z'_j\|\rightarrow 0$, whence $\xi+y=\lim_{j\rightarrow \infty} (\gamma_j^{(n+1)}\pi_\mathcal C(w_{n+1})+\pi_\mathcal C(z'_j))=a\pi_\mathcal C(w_{n+1})$, ensuring that $y=\pi(\xi+y)=a\pi(w_{n+1})$.

Define a new half-space $\mathbf u$ with $\overline{\mathbf u}=\R^\cup\la \mathcal C\ra+\R_+ w_{n+1}$ and $\partial(\mathbf u)\cap \mathbf u=\C^\cup (\mathcal C)$. Defining the limit
$\vec u_\mathbf u=(v_1,\ldots,v_{r_\ell},w_1,\ldots,w_n,w_{n+1})$, we see that we then have $\mathcal C\subseteq \mathcal X_1\cup\ldots\cup \mathcal X_{s-1}$ minimally encasing $-\vec u_\mathbf u^\triangleleft=-(v_1,\ldots,v_{r_\ell},w_1,\ldots,w_n)$, allowing us to set $\partial(\{\mathbf u\})=\mathcal C$.
As just shown above,  $\xi+y=a\pi_\mathcal C(w_{n+1})$, ensuring that $\pi_\mathcal C(\xi_r+y)=\xi+y=a\pi_\mathcal C(w_{n+1})$ for all $r$.
Thus the representative $\xi_r+y$ for the half-space $\x_r$ will also be a representative for the half-space $\mathbf u$. Moreover, $a\pi(w_{n+1})=\pi(\xi_r+y)=\pi(\xi+y)=y$, ensuring that the representative for the half-space $\mathbf u$ modulo $\mathcal E$ is equal to a positive multiple of  the representative for the half-space $\x_r$ modulo $\mathcal E$. Consequently, if we fix some $\x_r\in \mathcal X_s$ for defining $\mathcal X$, say $\x_r\in \mathcal X_s^{(r)}$ with $\mathcal X=\mathcal X^{(r)}$,  and then replace $\x_r\in \mathcal X^{(r)}_s$  with $\mathbf u$, we obtain a new virtual Reay system  $$\mathcal R'=(\mathcal X_1\cup \{\mathbf v_1\},\ldots,\mathcal X_{s-1}\cup \{\mathbf v_{s-1}\},(\mathcal X^{(r)}_s\setminus \{\x_r\}\cup \{\mathbf u\})\cup \mathcal \{\mathbf v_s\})$$ over $G_0$ with  $\mathcal X':=\mathcal Y\cup \mathcal X^{(r)}_s\setminus \{\x_r\}\cup \{\mathbf u\}$.
Set $\x_r'=\mathbf u$ and $\y'=\y$ for all other half-spaces $\y$ from $\mathcal R^{(r)}$.
Since $\pi(\mathcal X_s^{(r)}\setminus \{\x_r\}\cup \{\mathbf u\})=\pi(\mathcal X_s^{(r)})=\mathcal Y_s$, it follows from Proposition \ref{prop-finitary-MaxDecompChar} applied to the maximal series decomposition $\bigcup_{i=1}^{s-1}\mathcal X_i\cup \mathcal X_i^{(r)}$ that $\bigcup_{i=1}^{s-1}\mathcal X_i\cup (\mathcal X_s^{(r)}\setminus \{\x_r\}\cup\{\mathbf u\})$ is a maximal series decomposition.
Since $\mathcal B\prec \mathcal C$, we have $\partial(\x_r)=\R^\cup \la \mathcal B\ra\subseteq \R^\cup \la \mathcal C\ra=\partial(\mathbf u)$. Thus, since we have some representative for $\x_r$ that is also one for $\mathbf u$, it follows that any representative for $\x_r$ is also one for $\mathbf u$. Since $\R^\cup \la \mathcal B\ra\subseteq \R^\cup\la \mathcal C\ra$, Proposition \ref{prop-orReay-BasicProps}.9 implies $\mathcal B\subseteq \darrow \mathcal C$, whence $\x_r\cap \partial(\x_r)=\C^\cup (\mathcal B)\subseteq \C^\cup(\mathcal C)=\mathbf u\cap \partial(\mathbf u)$.  If $\R^\cup \la \mathcal C\ra=\R^\cup \la \mathcal B\ra$, then Proposition \ref{prop-orReay-BasicProps}.9 ensures that $\mathcal C\subseteq \darrow \mathcal B$, implying $\mathcal C=\mathcal C^*\preceq \mathcal B$, which would contradict that $\mathcal B\prec \mathcal C$ (note $\mathcal C^*=\mathcal C$ follows as $\mathcal C$ minimally encases the limit $-\vec w$).
In consequence, since any representative for $\x_r$ is one for $\mathbf u$, it follows that $\x_r\subset \mathbf u$, whence $\bigcup_{i=1}^{s-1}\mathcal X_i\cup \mathcal X_s^{(r)}\prec^*_\cup \bigcup_{i=1}^{s-1}\mathcal X_i\cup (\mathcal X_s^{(r)}\setminus \{\x_r\}\cup\{\mathbf u\})$.
Finally,  the representative sequence $\{\mathbf u(j)\}_{j=1}^\infty$ for $\mathbf u$ has $\mathbf u(j)=\x\big((i(j),j\big)$, which is a term from the representative sequence $\x_j(i)$. By construction, any such term, together with the remaining representatives for the other half-spaces from $\mathcal X^{(r)}\setminus \{\x_r\}=\mathcal X'\setminus\{\mathbf u\}$, generates the lattice $\Delta$ (as argued in Claim B).
But this ensures $\Delta\in \mathfrak P_\Z(\mathcal X')$, which together with $\bigcup_{i=1}^{s-1}\mathcal X_i\cup \mathcal X_s^{(r)}\prec^*_\cup \bigcup_{i=1}^{s-1}\mathcal X_i\cup (\mathcal X_s^{(r)}\setminus \{\x_r\}\cup\{\mathbf u\})$ contradicts that $\bigcup_{i=1}^{s-1}\mathcal X_i\cup\mathcal X_s^{(r)}\in \mathsf{Max}(\mathfrak X^*_s(G_0,\Delta))$, completing the proof.
\end{proof}

\begin{corollary}\label{cor-finitary-FiniteProps-II}
Let $\Lambda\subseteq \R^d$ be a full rank lattice, where $d\geq 0$, and let $G_0\subseteq \Lambda$ be a finitary subset with $\C(G_0)=\R^d$.  Then $\mathsf{Max}\big(\mathfrak X(G_0),\preceq_\cup)$, $\mathsf{Max}\big(\mathfrak X^*(G_0),\preceq^*_\cup)$, $\mathsf{Max}\big(\mathfrak X(G_0,\Delta),\preceq_\cup)$ and  $\mathsf{Max}\big(\mathfrak X^*(G_0,\Delta),\preceq^*_\cup)$ are all  finite,  with $\darrow \mathsf{Max}\big(\mathfrak X(G_0)\big)=\mathfrak X(G_0)$, $\darrow \mathsf{Max}\big(\mathfrak X^*(G_0)\big)=\mathfrak X^*(G_0)$, $\darrow \mathsf{Max}\big(\mathfrak X(G_0,\Delta)\big)=\mathfrak X(G_0,\Delta)$ and $\darrow \mathsf{Max}\big(\mathfrak X^*(G_0,\Delta)\big)=\mathfrak X^*(G_0,\Delta)$,  for any $\Delta\in \mathfrak P_\Z(G_0)$.
\end{corollary}

\begin{proof}
Let $\mathcal X\in \mathfrak X(G_0)$ and let $\mathcal X=\mathcal X_1\cup\ldots\cup \mathcal X_s$ be a maximal series decomposition. If $\mathcal X'\in \mathfrak X(G_0)$ with $\mathcal X\prec_\cup \mathcal X'$, then $\Summ{\x\in \mathcal X}\dim \partial(\x)<\Summ{\x'\in \mathcal X'}\dim \partial(\x')$ (as argued when showing anti-symmetry for $\preceq_\cup$). Thus, since $0\leq \Summ{\x\in \mathcal X}\dim \partial(\x)\leq d^2$, we see that there can be no infinite ascending (or descending) chains in $(\mathfrak X(G_0),\preceq_\cup)$.
 Since $\bigcup_{i=1}^s\mathcal X_i\prec_\cup^*\bigcup_{i=1}^s\mathcal X'_s$ implies $\mathcal X\prec_\cup \mathcal X'$, there are also no infinite ascending (or  descending) chains  in $(\mathfrak X^*(G_0),\preceq^*_\cup)$.
 It thus suffices to show $\mathsf{Max}(\mathfrak X(G_0),\preceq_\cup)$, $\mathsf{Max}(\mathfrak X^*(G_0),\preceq^*_\cup)$, $\mathsf{Max}(\mathfrak X(G_0,\Delta),\preceq_\cup)$ and $\mathsf{Max}(\mathfrak X^*(G_0,\Delta),\preceq^*_\cup)$ are all finite, for any $\Delta\in \mathfrak P_\Z(G_0)$. Since each $\mathcal X\in \mathsf{Max}(\mathfrak X(G_0))$ with $\Delta\in \mathfrak P_\Z(\mathfrak X)$ is also an element of $\mathsf{Max}(\mathfrak X(G_0,\Delta))$, and since $\mathfrak P_\Z(G_0)$ is finite by Theorem \ref{thm-finitary-FiniteProps-I}.2, it suffices to show $\mathsf{Max}(\mathfrak X(G_0,\Delta))$ and $\mathsf{Max}(\mathfrak X^*(G_0,\Delta))$ are finite for an arbitrary $\Delta\in \mathfrak P_\Z(G_0)$. If $\mathcal X\in \mathsf{Max}(\mathfrak X(G_0,\Delta))$ and $\mathcal X=\bigcup_{i=1}^s\mathcal X_i$ is any maximal series decomposition, then  $\bigcup_{i=1}^s\mathcal X_i\in \mathsf{Max}(\mathfrak X^*(G_0,\Delta))$ (since $\bigcup_{i=1}^s\mathcal X_i\prec_\cup^*\bigcup_{i=1}^s\mathcal X'_i$ would imply $\mathcal X\prec_\cup \mathcal X'=\bigcup_{i=1}^s\mathcal X'_i$). Thus it suffices to prove $\mathsf{Max}(\mathfrak X^*(G_0,\Delta))$ is finite for an arbitrary $\Delta\in \mathfrak P_\Z(G_0)$, as all other parts follow routinely from this as just explained, and this was proved in Theorem \ref{thm-finitary-FiniteProps-II}.
\end{proof}

We continue
with another finite-like property of finitary subsets, showing that, even though a finitary set $G_0$ can have an infinite number of atoms, every atom must nonetheless contain some term from a fixed \emph{finite} subset of  $G_0$.

\begin{proposition}\label{prop-finitary-FiniteDeletion}
Let $\Lambda\subseteq \R^d$ be a full rank lattice, where $d\geq 0$, and let $G_0\subseteq \Lambda$ be a finitary subset with $\C(G_0)=\R^d$. Let $$\mathcal A_\bullet(G_0)=\{U\in \mathcal A(G_0):\;\mbox{ there exists $X\subseteq \supp(U)$ with $\emptyset\neq X\in X(G_0)$\}}$$
Then the following hold.
\begin{itemize}
\item[1.] $\mathcal A(G_0)\setminus \mathcal A_\bullet(G_0)$ is finite.
\item[2.] There exist finite sets $\wtilde X\subseteq G_0^\diamond$ and $\wtilde Y\subseteq G_0$ with $\mathcal A(G_0\setminus \wtilde X)$ finite and $\mathcal A(G_0\setminus \wtilde Y)=\emptyset$.
\end{itemize}
\end{proposition}

\begin{proof} 1. Let $\mathcal A_\bullet^{\mathsf{elm}}(G_0)=\mathcal A^{\mathsf{elm}}(G_0)\cap \mathcal A_\bullet(G_0)$.
We begin by showing that $\mathcal A^{\mathsf{elm}}(G_0)\setminus \mathcal A_\bullet^{\mathsf{elm}}(G_0)$ is finite. Assume by contradiction that $\{U_i\}_{i=1}^\infty$ is an infinite sequence of distinct elementary atoms $U_i\in \mathcal A^{\mathsf{elm}}(G_0)$ such that no nonempty subset of $\supp(U_i)$ lies in $X(G_0)$, for all $i\geq 1$. By passing to a subsequence of $\{U_i\}_{i=1}^\infty$, we can w.l.o.g. assume $|U_i|=s\geq 2$ for all $i$, say with each $U_i=\{x_i^{(1)},\ldots,x_i^{(s)}\}$. Since the $U_i$ are elementary atoms with $|U_i|\geq 2$, we have $0\notin \supp(U_i)$ with $\supp(U_i)$ a minimal positive basis, for all $i$. Thus, by again passing to a subsequence of $\{U_i\}_{i=1}^\infty$, we can further assume each sequence $\{x_i^{(j)}\}_{i=1}^\infty$ is a radially convergent sequence of terms from $G_0$ with limit $u_j$, for $j\in [1,s]$. Moreover, since $G_0\subseteq \Lambda$ is a set of lattice points (so that any bounded subset of $G_0$ is finite), we either have $u_j$ unbounded or w.l.o.g. $\{x_i^{(j)}\}_{i=1}^\infty$ constant, for each $j\in [1,s]$. Since the $U_i$ are all distinct, the latter cannot occur for all $j\in [1,s]$.

If $0\notin \C^*(u_1,\ldots,u_s)$, then Lemma \ref{lem-finite-nozs} implies  there will be an open half-space $\mathcal E^\circ_+$ containing all $u_j$ for $j\in [1,s]$. Consequently, since $x_i^{(j)}/|x_i^{(j)}\|\rightarrow u_j$ for each $j\in [1,s]$, it follows that each $\supp(U_i)$, with $i$ sufficiently large, will also be contained in the open half-space $\mathcal E^\circ_+$, ensuring $0\notin \C(\supp(U_i))$, which contradicts that each $U_i$ is an atom. Therefore we instead conclude that $0\in \C^*(u_1,\ldots,u_s)$. Thus there is some nonempty subset $J\subseteq [1,s]$ such that $U=\{u_j:\;j\in J\}$ is a minimal positive basis.
If $u_j\notin G_0^\infty$ for all $j\in J$, then $J\subset [1,s]$ is a proper, nonempty subset (since at least one $u_j$ with $j\in[1,s]$ is unbounded).
 In such case, $\{x_i^{(j)}:j\in J\}\subset \supp(U_i)$ will be a minimal positive basis, contradicting that each $\supp(U_i)$ is a minimal positive basis, thus ensuring any proper subset is linearly independent. Therefore we instead conclude that there is some $t\in J$ such that $u_t\in G_0^\infty$, and in view of $G_0$ being finitary, it then follows that $J\cap G_0^\infty=\{t\}$. But now $\{x_i^{(j)}:\;j\in J\setminus \{t\}\}\cup \{u_t\}$ is a minimal positive basis. Thus $\{x_i^{(j)}:\;j\in J\setminus \{t\}\}$ minimally encases $-u_t\in G_0^\infty$, ensuring that each $X_i=\{x_i^{(j)}:\;j\in J\setminus \{t\}\}\subset \supp(U_i)$ is a subset with $\emptyset\neq X_i\in X(G_0)$, contrary to assumption. So we instead conclude that $\mathcal A^{\mathsf{elm}}(G_0)\setminus \mathcal A_\bullet^{\mathsf{elm}}(G_0)$ is finite.

 By Theorem \ref{thm-carahtheodory-elm-atom}, every atom $U\in \mathcal A(G_0)$ has a factorization $U=\prod_{i\in [1,\ell]}^\bullet U_i^{[\alpha_i]}$ with each $U_i\in \mathcal A^{\mathsf{elm}}(G_0)$ and each $\alpha_i\in \Q$ with $0<\alpha_i\leq 1$. In particular, $\supp(U)=\bigcup_{i=1}^\ell \supp(U_i)$. Thus, if $U\notin\mathcal A_\bullet(G_0)$, then $U_i\notin \mathcal A_\bullet^{\mathsf{elm}}(G_0)$ for every $i\in [1, \ell]$ as well, implying that \be\supp(U)\subseteq \bigcup_{V\in \mathcal A^{\mathsf{elm}}(G_0)\setminus \mathcal A^{\mathsf{elm}}_\bullet(G_0)}\supp(V).\label{cree}\ee  However, as just shown, $\mathcal A^{\mathsf{elm}}(G_0)\setminus \mathcal A^{\mathsf{elm}}_\bullet(G_0)$ is finite, ensuring that  $\bigcup_{V\in \mathcal A^{\mathsf{elm}}(G_0)\setminus \mathcal A^{\mathsf{elm}}_\bullet(G_0)}\supp(V)$ is finite,
 implying via Dickson's Theorem \cite[Theorem 1.5.3, Corollary 1.5.4]{alfredbook} that the set of atoms $\mathcal A\big(\bigcup_{V\in \mathcal A^{\mathsf{elm}}(G_0)\setminus \mathcal A^{\mathsf{elm}}_\bullet(G_0)}\supp(V)\big)$ is also finite. Hence there can only be a finite number of atoms $U$ satisfying \eqref{cree}, which implies $\mathcal A(G_0)\setminus \mathcal A_\bullet(G_0)$ is finite, completing the proof of Item 1.

2. Every $U\in \mathcal A_\bullet(G_0)$ has some nonempty $X\in X(G_0)$ with $X\subseteq \supp(U)$. Since any $X\in X(G_0)$ contains some $X_1\subseteq X$ with $X_1\in X(G_0)$ irreducible, we conclude that every  $U\in \mathcal A_\bullet(G_0)$ has some (nonempty) irreducible $X\in X(G_0)$ with $X\subseteq \supp(U)$. By Theorem \ref{thm-finitary-FiniteProps-I}.1, there are only a finite number of irreducible sets in $X(G_0)$. Let $\wtilde X$ be obtained by including one element from each irreducible set from $X(G_0)$. Then $\wtilde X$ is finite and $\mathcal A(G\setminus \wtilde X)\subseteq \mathcal A(G_0)\setminus \mathcal A_\bullet(G_0)$ is finite by Item 1. Moreover, since $X\subseteq G_0^\diamond$ for any $X\in X(G_0)$ (as remarked at the start of Section \ref{sec-series-decomp}), we have $\wtilde X\subseteq G_0^\diamond$. Let $\wtilde Y$ be obtained by taking $\wtilde X$ and including one element from the support of each atom from $\mathcal A(G\setminus \wtilde X)$. Since $\mathcal A(G\setminus \wtilde X)$ and $\wtilde X$ are finite, so too is $\wtilde Y$, and $\mathcal A(G_0\setminus \wtilde Y)=\emptyset$ by construction, completing the proof of Item 2.
\end{proof}

Next, we give the finiteness from below property for the partial order $\preceq_\Z$. We will strengthen the implied consequence that $\mathsf{Min}\big(X(G_0),\preceq_\Z)$ is finite later in Proposition \ref{prop-finitary-mintype}.

\begin{theorem}\label{thm-finitary-FiniteProps-III}
Let $\Lambda\subseteq \R^d$ be a full rank lattice, where $d\geq 0$, and let $G_0\subseteq \Lambda$ be a finitary subset with $\C(G_0)=\R^d$. Then there are neither infinite descending chains nor infinite antichains in $(X(G_0),\preceq_\Z)$. In particular, $\mathsf{Min}\big(X(G_0),\preceq_\Z)$ is finite and $\uarrow \mathsf{Min}\big(X(G_0),\preceq_\Z)=X(G_0)$.
\end{theorem}

\begin{proof}
We need only prove $X(G_0)$ contains no infinite anti-chains nor infinite descending chains, as the remainder of the theorem  then follows from Proposition \ref{prop-poset-conseq}. In view of Theorem \ref{thm-finitary-FiniteProps-I}.2, there are only a finite number of lattices $\Delta\in \mathfrak P_\Z(G_0)$.  Thus it suffices to prove $X(G_0,\Delta)$ contains no infinite anti-chains nor infinite descending chains for an arbitrary $\Delta\in \mathfrak P_\Z(G_0)$, which we now fix. Recall $X\preceq_\Z Y$ for $X,\,Y\in X(G_0,\Delta)$ is equivalent to $\C(X)\subseteq \C(Y)$.
Let $X_s(G_0,\Delta)\subseteq X(G_0,\Delta)$ consist of all $X\in X(G_0,\Delta)$, so $\Z\la X\ra=\Delta$, having a maximal series decomposition of length $s$, say $X=X_1\cup \ldots\cup X_s$.  We proceed by induction on $s=1,2\ldots,d$ to show that $X_s(G_0,\Delta)$ contains no infinite anti-chains nor infinite descending chains, which will complete the proof as there are only a finite number of possible values for $s\in [1,d]$. Note $X_1(G_0,\Delta)$ consists precisely of all irreducible sets $X\in X(G_0,\Delta)$. Thus $X_1(G_0,\Delta)$ is finite by Theorem \ref{thm-finitary-FiniteProps-I}.1, in which case it trivially can contain no infinite descending chains nor infinite anti-chains, completing the base $s=1$ of the induction. So we now assume $s\geq 2$.

Let $X\in X_s(G_0,\Delta)$ be arbitrary with maximal series decomposition $X=X_1\cup \ldots\cup X_s$. Then $Y=X_1\cup \ldots\cup X_{s-1}\in X_{s-1}(G_0)$ by Proposition \ref{prop-finitary-MaxDecompChar}.  Let $\Delta'=\Z\la Y\ra\in \mathfrak P_\Z(G_0)$. As there are only a finite number of possibilities for $\Delta'$ by Theorem \ref{thm-finitary-FiniteProps-I}.2, it suffices to show there are no infinite anti-chains nor infinite descending chains involving sets $X\in X_s(G_0,\Delta)$ having a series decomposition $X=X_1\cup \ldots\cup X_s$ with $\Delta'=\Z\la X_1\cup\ldots\cup X_{s-1}\ra$, for an arbitrary $\Delta'\in \mathfrak P_\Z(G_0)$ which we now fix. Let $\mathcal E=\R\la Y\ra=\R\la \Delta'\ra$ and let  $\pi:\R^d\rightarrow \mathcal E^\bot$ be the orthogonal projection. Observe that $\Delta'$ and $\Delta\cap \mathcal E$ are both full rank lattices in $\mathcal E$ with $\Delta'\leq \Delta\cap \mathcal E$. Thus $(\Delta\cap \mathcal E)/\Delta'$ is finite, ensuring there is a finite set $\overline \Delta\subseteq \Delta\cap \mathcal E$ of coset representatives for $\Delta'$.

Proposition \ref{prop-finitary-MaxDecompChar} implies that $\pi(X_s)\in X(\pi(X_0))$ is irreducible. Proposition \ref{prop-finitary-Modulo-Inheritence} implies that $\pi(G_0)$ is finitary with $\pi(\Lambda)$ a lattice. Thus Theorem \ref{thm-finitary-FiniteProps-I}.1 implies that there are only a finite number of possibilities for $\pi(X_s)$, and so we may add the condition that there is a fixed subset $Y_s\subseteq G_0\cap \Delta$ with $|\pi(Y_s)|=|Y_s|$ so that   $\pi(X_s)=\pi(Y_s)$ is the same  fixed irreducible set for all $X$ under consideration. Note $Y\cup Y_s$ is a lattice basis for $\Delta$.

If $x\in X_s$ has $\pi(x)=\pi(y)$ with $y\in Y_s$, then $x-y\in \Delta\cap \ker \pi=\Delta\cap \mathcal E$. It follows that $x=y+\xi_x-\varpi_x$ for some  $\varpi_x\in \Delta'$ and $\xi_x\in \overline \Delta$. Since there are only a finite number of possibilities for $\xi_x\in \overline \Delta$, we can add the condition that, for each $y\in Y_s$, there is a fixed $\xi_y\in \overline \Delta$, such that, for all $X$ under consideration,  the element $x\in X_s$ with $\pi(x)=\pi(y)$ always has $x=y+\xi_y-\varpi_x$ for some $\varpi_x\in \Delta'$.
Let $Y_s=\{y_1,\ldots,y_r\}$ be the distinct elements of $Y_s$ and adapt the abbreviation $\xi_j=\xi_{y_j}$ for $j\in [1,r]$. Thus, if $X=\{x_1,\ldots,x_r\}$ are the elements of $X$ indexed so that $\pi(x_j)=\pi(y_j)$ for $j\in [1,r]$, then $x_j=y_j+\xi_j-\varpi_{x_j}$ with $\varpi_{x_j}\in \Delta'$.

By induction hypothesis and Proposition \ref{prop-poset-conseq}, we have $\uarrow \mathsf{Min}(X_{s-1}(G_0,\Delta'))=X_{s-1}(G_0,\Delta')$ with  $\mathsf{Min}(X_{s-1}(G_0,\Delta'))$ finite. Thus every $Y\in X_{s-1}(G_0,\Delta')$ has some $Y_\emptyset\in \mathsf{Min}(X_{s-1}(G_0,\Delta'))$ with $Y_\emptyset\preceq_\Z Y$. As there are only a finite number of possibilities for $Y_\emptyset$, we may add the condition that there is some fixed subset  $Y_\emptyset\in \mathsf{Min}(X_{s-1}(G_0,\Delta'))$ with $Y_\emptyset\preceq_\Z Y=X_1\cup\ldots\cup X_{s-1}$ for all $X$ under consideration.

In summary, the above work means there are fixed $\Delta,\,\Delta'\in \mathfrak P_\Z(G_0)$, a fixed subset $Y_s=\{y_1,\ldots,y_r\}\subseteq G_0\cap \Delta$ with $\pi(Y_s)\in X(\pi(G_0))$ irreducible, a fixed set of representatives $\overline{\Delta}$ for the lattice $\Delta\cap \mathcal E$ modulo $\Delta'$,  where $\mathcal E=\R\la \Delta'\ra$, fixed $\xi_1,\ldots, \xi_r\in \overline{\Delta}\subseteq \Delta\cap \mathcal E$, and a fixed $Y_\emptyset\in \mathsf{Min}(X_{s-1}(G_0,\Delta'))$ such that, letting  $\Omega\subseteq X(G_0)$ consist of all $X\in X(G_0)$ having a maximal series decomposition $X=X_1\cup\ldots\cup X_s$ such that $\Z\la X\ra=\Delta$, \ $\Z\la Y\ra=\Delta'$, \ $\R\la Y\ra=\mathcal E$, \ $\mathcal Y_\emptyset \preceq_\Z Y$, and $X_s=\{x_1,\ldots,x_r\}$ with each  $$x_j=y_j+\xi_j-\varpi_{x_j}$$ for some $\varpi_{x_j}\in \Delta'$, where $j\in [1,r]$ and $Y=X_1\cup\ldots\cup X_{s-1}$, then it suffices to show there are no infinite anti-chains nor infinite descending chains in $\Omega$, for the arbitrary fixed parameters indicated, to complete the induction and thus also the proof.

For $j\in [1,r]$, let $\wtilde X_j=\{\varpi_{x_j}:\: X\in \Omega\}\subseteq \Delta'\subseteq \mathcal E$.
Let $\mathcal R_\emptyset$ be a realization of $Y_\emptyset\in X(G_0,\Delta)$ in which all half-spaces corresponding to $Y_\emptyset$ are one-dimensional.
If $\{z_i\}_{i=1}^\infty$ is an asymptotically filtered sequence of terms from $\wtilde X_j$ with fully  unbounded limit $\vec u=(u_1,\ldots,u_t)$, then $u_1,\ldots,u_t\in \mathcal E=\R\la \Delta'\ra=\R\la Y_\emptyset\ra$ (as $\wtilde X_j\subseteq \mathcal E$), and then  $\{-z_i+\xi_j+y_j\}_{i=1}^\infty$ will be an asymptotically filtered sequence of terms from $G_0$ with  fully  unbounded limit $-\vec u$ (after discarding the first few terms).
Thus Proposition \ref{prop-finitary-basics}.3 applied to $\mathcal R_\emptyset$ ensures that $\vec u$ is encased by $Y_\emptyset$, in which case Theorem \ref{thm-nearness-characterization}.4 applied to $\C(Y_\emptyset)$ implies that $\wtilde X_j$ is bound to $\C(Y_\emptyset)$ for all $j\in [1,r]$, and thus $\wtilde X_j$ is also bound  to $\C_{\Z}(Y_\emptyset)$ for all $j\in [1,r]$. This means there is some fixed ball $B\subseteq \mathcal E$ such that every $\varpi\in \bigcup_{j=1}^r\wtilde X_j$ has some $z\in \C_{\Z}(Y_\emptyset)$ such that $\varpi\in B+z$. As $\varpi,\,z\in \Delta'$, we are assured  that $\varpi-z\in B\cap \Delta'$, which is a bounded set of lattice points. Thus each $\varpi_j \in \wtilde X_j$ has $\varpi_j-\xi'_j\in \C_{\Z}(Y_\emptyset)$ for some $\xi'_j\in B\cap \Delta'$.
As there are only a finite number of choices for $\xi'_j$, we can make one last restriction by only considering $X\in \Omega$ for which the   same fixed value of $\xi'_j$ occurs for all $\varpi_j\in \wtilde X_j$, for each $j\in [1,r]$. Let $\Omega'\subseteq \Omega$ be resulting subset. This means  it suffices to show there are no infinite anti-chains nor infinite descending chains in $\Omega'$ to complete the proof. Moreover, replacing the fixed coset representative $\xi_j$ with the alternative fixed coset representative $\xi_j-\xi'_j$ for each $j\in [1,r]$ (thus replacing $\varpi_j$ by $\varpi_j-\xi'_j$), we obtain $$\varpi_{x_j}\in \C_{\Z}(Y_\emptyset)$$ rather than simply $\varpi_{x_j}\in \Delta'=\Z\la Y_\emptyset\ra$, for each $j\in [1,r]$.

Each  $X\in \Omega'$ corresponds to the tuple $$\varphi(X)=(Y,\varpi_{1},\ldots,\varpi_{r})\in X_{s-1}(G_0,\Delta')\times \underbrace{\C_{\Z}(Y_\emptyset)\times \ldots\times \C_{\Z}(Y_\emptyset)}_r.$$ Indeed,  $Y\subseteq X$ consists of all $x\in X$ with $\pi(x)=0$, while, since  $\pi(Y_s)$ consists of $|Y_s|=r$ distinct elements, there is a unique indexing of $X_s=X\setminus Y=\{x_1,\ldots,x_r\}$ such that $\pi(x_j)=\pi(y_j)$ for all $j\in [1,r]$, and then $\varpi_j=\varpi_{x_j}=-x_j+y_j+\xi_j$ as defined above.
Moreover, if $X'\in \Omega'$ is a set with corresponding tuple $\varphi(X')=(Y',\varpi'_{1},\ldots,\varpi'_{r})=(Y,\varpi_{1},\ldots,\varpi_{r})=\varphi(X)$, then $X=X'$. Thus we see that the elements $X\in \Omega'$ are in bijective correspondence with a subset $$\overline {\Omega'}\subseteq X_{s-1}(G_0,\Delta')\times \underbrace{\C_{\Z}(Y_\emptyset)\times \ldots\times \C_{\Z}(Y_\emptyset)}_r.$$
Given a set $Y\in X_{s-1}(G_0,\Delta')$, we define a partial order $\preceq_Y$ on $\C_{\Z}(Y_\emptyset)$ by declaring $$\varpi\preceq_Y \varpi'\quad\mbox{ when } \quad  \varpi'\in \varpi+\C(Y),$$ for $\varpi,\,\varpi'\in \C_{\Z}(Y_\emptyset)$.
Note $\preceq_Y$ is transitive as $\C(Y)$ is convex, is reflexive as $0\in \C(Y)$, and is antisymmetric since  $Y$ is linearly independent (which ensures that $\C(Y)\cap -\C(Y)=\{0\}$).
Also, since $\varpi,\,\varpi'\in\Delta'=\Z\la Y\ra$ for all $Y\in X_{s-1}(G_0,\Delta')$ with $Y$ linearly independent,  we have $$\varpi'\in \varpi+\C(Y)\quad\mbox{ if and only if }\quad\varpi'\in \varpi+\C_{\Z}(Y).$$

\subsection*{Claim A} For $X,\,X'\in \Omega'$ with $\varphi(X)=(Y,\varpi_1,\ldots,\varpi_r)$ and $\varphi(X')=(Y',\varpi'_1,\ldots,\varpi'_r)$, we have $X\preceq_\Z X'$ if and only if $Y\preceq_\Z Y'$ and $\varpi_j\preceq_{Y'} \varpi'_j$ for all $j\in [1,r]$.

\begin{proof} Suppose $X\preceq_\Z X'$. Then $\C_{\Z}(X)\subseteq \C_{\Z}(X')$. Let  $\Summ{y\in Y}\alpha_yy\in \C_{\Z}(Y)$ be an arbitrary element, where $\alpha_y\in \Z_+$. Then $\Summ{y\in Y}\alpha_yy\in \C_{\Z}(Y) \subseteq \C_{\Z}(X)\subseteq \C_{\Z}(X')$,  whence
$\Summ{y\in Y}\alpha_yy=\Summ{x'\in X'}\beta_{x'}x'$ for some $\beta_{x'}\in \Z_+$. Thus, since $\R\la Y\ra=\R\la Y'\ra=\mathcal E=\ker \pi$ with $\pi(X'\setminus Y')$ a set of $|X'\setminus Y'|$ linearly independent elements, it follows that $\beta_{x'}=0$ for all $x'\in X'\setminus Y'$, whence
$\Summ{y\in Y}\alpha_yy=\Summ{x'\in Y'}\beta_{x'}x'\in \C_{\Z}(Y')$, which shows $\C_{\Z}(Y)\subseteq \C_{\Z}(Y')$, and thus that $Y\preceq_\Z Y'$. Let $X\setminus Y=X_s=\{x_1,\ldots,x_r\}$ and $X'\setminus Y'=X'_s=\{x'_1,\ldots,x'_r\}$ with $$x_j=y_j+\xi_j-\varpi_{j}\quad\und\quad x'_j=y_j+\xi_j-\varpi'_{j} \quad\mbox{ for $j\in [1,r]$}.$$ Let $j\in [1,r]$ be arbitrary. Then $x_j\in \C_{\Z}(X)\subseteq \C_{\Z}(X')$, implying that $x_j=\Summ{x'\in X'}\beta_{x'}x'$ for some $\beta_{x'}\in\Z_+$. Since $\pi(x_j)=\pi(y_j)$ with $\pi(X'\setminus Y')=\{\pi(y_1),\ldots,\pi(y_r)\}$ linearly independent, it follows that $\beta_{x'}=0$ for all $x'\in X'_s\setminus\{x'_j\}$ and $\beta_{x'_j}=1$. Hence $$y_j+\xi_j-\varpi_j=x_j=x'_j+\Summ{y'\in Y'}\beta_{y'}y'=y_j+\xi_j-\varpi'_j+\Summ{y'\in Y'}\beta_{y'}y'\in y_j+\xi_j-\varpi'_j+\C_{\Z}(Y'),$$ implying that $\varpi'_j\in \varpi_j+\C_{\Z}(Y')$, and thus that $\varpi_j\preceq_{Y'} \varpi'_j$. This establishes one direction of the claim.

Next suppose that $Y\preceq_\Z Y'$ and $\varpi_j\preceq_{Y'} \varpi'_j$ for all $j\in [1,r]$.
Let $\Summ{x\in X}\alpha_xx\in \C_{\Z}(X)$ be arbitrary, where $\alpha_x\in \Z_+$. Since $Y\preceq_\Z Y'$, we have $\C_{\Z}(Y)\subseteq \C_{\Z}(Y')$, ensuring that $\Summ{y\in Y}\alpha_yy\in \C_{\Z}(Y')$. Now $X_s=X\setminus Y=\{x_1,\ldots,x_r\}$ with $x_j=y_j+\xi_j-\varpi_j$ and $X'_s=X'\setminus Y'=\{x'_1,\ldots,x'_r\}$ with $x'_j=y_j+\xi_j-\varpi'_j$. Thus, since $\varpi_j\preceq_{Y'} \varpi'_j$, which ensures $-\varpi_j\in -\varpi'_j+\C_{\Z}(Y')$, we have $x_j\in x'_j+\C_{\Z}(Y')$, for $j\in [1,r]$. As a result, since $\C_{\Z}(Y')$ is a convex lattice cone (closed under addition and positive scalar multiplication by integers), we have  $$\Summ{x\in X}\alpha_xx=\Summ{y\in Y}\alpha_yy+\Sum{j=1}{r}\alpha_{x_j}x_j\in \Sum{j=1}{r}\alpha_{x_j}x'_j+\C_{\Z}(Y')\subseteq \C_{\Z}(X'_s)+\C_{\Z}(Y')=\C_{\Z}(X').$$ Since $\Summ{x\in X}\alpha_xx\in \C_{\Z}(X)$ was arbitrary, this shows $\C_{\Z}(X)\subseteq \C_{\Z}(X')$, implying $X\preceq_\Z X'$, which completes the other direction in Claim A.
\end{proof}

 If $\C_{\Z}(Y_1)\subseteq \C_{\Z}(Y_2)$, where $Y_1,\,Y_2\in X_{s-1}(G_0,\Delta')$, then $\varpi\preceq_{Y_1} \varpi'$ implies $\varpi\preceq_{Y_2} \varpi'$. In particular, $\varpi\preceq_{Y_\emptyset} \varpi'$ implies $\varpi \preceq_{Y} \varpi'$ for any $Y$ occurring as the first coordinate of an element from $\overline{\Omega'}$
(recall that $\mathsf C_{\Z}(Y_\emptyset)\subseteq \C_{\Z}(Y)$ since $Y_\emptyset \preceq_\Z Y$ by construction of $\Omega$).
We make $\overline{\Omega'}$ into a poset by declaring $(Y,\varpi_1,\ldots,\varpi_r)\preceq(Y',\varpi'_1,\ldots,\varpi'_r)$ when  $Y\preceq_\Z Y'$ and $\varpi_j\preceq_{Y'} \varpi'_j$ for all $j\in [1,r]$, with transitivity guaranteed by the observation just noted.
Claim A then ensures that $\varphi$ gives an isomorphism between the posets $\Omega'$ and  $\overline{\Omega'}$.
Moreover, it is readily seen (since $Y_\emptyset$ is linearly independent) that there is an isomorphism of  posets $(\C_{\Z}(Y_\emptyset),\preceq_{Y_\emptyset})\cong \Z_+^{|Y_\emptyset|}$ using the product order on $\Z_+^{|Y_\emptyset|}$. Consequently, in view of Claim A, we see that we can remove relations from the poset $ \Omega'\cong \overline{\Omega'}$ to result in a poset isomorphic to a subset of $X_{s-1}(G_0,\Delta')\times \underbrace{\Z_+^{|Y_\emptyset|}\times \ldots\times \Z_+^{|Y_\emptyset|}}_r$.
By induction hypothesis, there are no infinite anti-chains nor infinite descending chains in $X_{s-1}(G_0,\Delta')$, and there are also no  infinite anti-chains nor infinite descending chains in $\Z_+^{|Y_\emptyset|}$ (recall the discussion in Section \ref{sec-poset-lattices}). Therefore iterated application of Proposition \ref{prop-poset-prod} implies there are no
 infinite anti-chains in $X_{s-1}(G_0,\Delta')\times \underbrace{\Z_+^{|Y_\emptyset|}\times \ldots\times \Z_+^{|Y_\emptyset|}}_r$, and thus also no infinite anti-chains in $\Omega'$, as this lattice can be obtained from  $X_{s-1}(G_0,\Delta')\times \underbrace{\Z_+^{|Y_\emptyset|}\times \ldots\times \Z_+^{|Y_\emptyset|}}_r$ by adding additional relations and removing elements (which both can only reduce the size of an anti-chain). It remains to show that $\Omega'$ also has no infinite descending chains.

Suppose by contradiction that $\{X_i\}_{i=1}^\infty$ is a strictly descending sequence of $X_i\in \Omega'$ and let $\varphi(X_i)=(Y_i,\varpi_i^{(1)},\ldots,\varpi_i^{(r)})$ for $i\geq 1$. Then, in view of Claim A, it follows that $\{Y_i\}_{i=1}^\infty$ is a  descending sequence of $Y_i\in X_{s-1}(G_0,\Delta')$. As a result, since $X_{s-1}(G_0,\Delta')$ contains no infinite descending chain by induction hypothesis, it follows that, by discarding the first few terms in $\{X_i\}_{i=1}^\infty$, we obtain $Y_i=Y$ for all $i$.  Since $\{X_i\}_{i=1}^\infty$ is a descending sequence,  Claim A now ensures that  $\{\varpi_i^{(j)}\}_{i=1}^\infty$ is a descending sequence in $(\mathsf C_{\Z}(Y_{\emptyset}),\preceq_{Y})$ for each $j\in [1,r]$.
Since $\C_\Z(Y_\emptyset)\subseteq \C_\Z(Y)$ in view of $Y_\emptyset\preceq_\Z Y$, it follows that $(\mathsf C_{\Z}(Y_{\emptyset}),\preceq_{Y})$ is a sub-poset of $(\mathsf C_{\Z}(Y),\preceq_{Y})$, which is isomorphic to $\Z_+^{|Y|}$ as $Y$  is linearly independent. However, the poset $\Z_+^{|Y|}$ contains no infinite descending chains (as discussed in Section \ref{sec-poset-lattices}), meaning neither the isomorphic poset $(\mathsf C_{\Z}(Y),\preceq_{Y})$ nor the sub-poset $(\mathsf C_{\Z}(Y_\emptyset),\preceq_{Y})$ contain infinite descending chains. It follows that each chain $\{\varpi_i^{(j)}\}_{i=1}^\infty$ for $j\in [1,r]$ eventually stabilizers, contradicting that $\{X_i\}_{i=1}^\infty$ is a strictly descending chain. With this contradiction,  we instead conclude that $\Omega'$ contains no infinite descending chains,  completing the induction and proof.
\end{proof}

\subsection{Interchangeability and the Structure of $X(G_0)$}

The goal of this subsection is study the structure of $X(G_0)$. One of our main aims is to partition $X(G_0)$ into a \emph{finite} number of subsets, each consisting of various $X\in X(G_0)$ having the same type, such that sets of the same type  posses common regularity properties. We begin with the basic definitions.

Let $\Lambda\leq \R^d$ be a full rank lattice and let $G_0\subseteq \Lambda$ be a finitary subset with $\C(G_0)=\R^d$. Let $X\in X(G_0)$, let $X=X_1\cup\ldots\cup X_s$ be a maximal series decomposition of $X$. For $j\in [0,s]$, let $\Delta_j=\Z\la X_1\cup\ldots\cup
X_j\ra$ and let $\pi_j:\R^d\rightarrow \R\la X_1\cup \ldots\cup X_j\ra^\bot$ be the orthogonal projection.  Note $\Delta_j=\R\la X_1\cup \ldots\cup X_j\ra\cap \Delta$, for $j\in [0,s]$, as $X$ is linearly independent, where $\Delta:=\Delta_s=\Z\la X\ra$. Also, if $x\in X_j$, then $j$ is the minimal index such that $\pi_j(x)=0$.
Suppose $X'\in X(G_0)$ is another set with maximal series decomposition  $X'=X'_1\cup\ldots\cup X'_s$ such that $\pi_{j-1}(X_j)=\pi_{j-1}(X'_j)$ for $j=1,2,\ldots,s$. Then $\R\la X_1\cup \ldots\cup X_j\ra=\R\la X'_1\cup \ldots\cup X'_j\ra$ for $j=0,1,2\ldots,s$.  Moreover,  each $x\in X$ identifies uniquely with some $x'\in X'$. Namely, $x\in X_j$ with $j\in [1,s]$ the minimal index such that $\pi_j(x)=0$, and then $x'\in X'_j$ is the unique element with $\pi_{j-1}(x)=\pi_{j-1}(x')$. This allows us to define a bijection between and $X$ and $X'$ by sending $x\mapsto x'$ as just defined. In such case, we say that $X$ and $X'$ have the same \textbf{linear type}.
If we additionally have $$\Z\la X'_1\cup\ldots\cup X'_{j}\ra=\Z\la X_1\cup\ldots\cup X_{j}\ra\quad\mbox{ for every $j\in [1,s]$},$$ then we say that $X$ and $X'$ have the same  \textbf{lattice type}.
For each possible lattice type, we can fix a representative set $Z\in X(G_0)$ with series decomposition $Z=Z_1\cup\ldots\cup Z_s$ having this lattice type. Then, given any other $X\in X(G_0)$ having the same lattice type, say with associated maximal series decomposition $X=X_1\cup\ldots\cup X_s$, there is a bijection $\phi:X\rightarrow Z$ such that $\Z\la X_1\cup \ldots\cup X_{j}\ra=\Z\la \phi(X_1)\cup \ldots \cup \phi(X_{j})\ra$ with $\phi(X_j)=Z_j$ for $j\in [1,s]$.  If $X'\in X(G_0)$ is another set with the same lattice type, we likewise have a bijection $\phi':X'\rightarrow Z$. Moreover, if $x\in X\cap X'$, then $\phi(x)=\phi'(x)$ by how the bijections $\phi$ and $\phi'$ are defined above. Thus we can extend the domain of $\phi$ to obtain a map $\phi:\bigcup_X X\rightarrow Z$, where the union runs over all $X\in X(G_0)$ having the same lattice type as $Z$, and we  identify $\phi$ as the lattice type of $Z$ itself, allowing us to say a set $X\in X(G_0)$ has lattice type $\phi$.
Note a given $X\in X(G_0)$ may have multiple types. However, each maximal series decomposition $X=X_1\cup \ldots\cup X_s$ corresponds to precisely one. The type $\phi$ associated to the empty set is called the trivial type, which corresponds to when $s=0$. Only the empty set has trivial type. We remark that the lattice type $\phi$ depends on the associated linear type and the lattices $\Delta_j$ for $j=1,\ldots,s$. However, once the linear type is fixed, in turn fixing the values of the subspaces $\R\la X_1\cup \ldots\cup X_j\ra$ for $j\in [0,s]$, the value of $\Delta_s=\Z\la X\ra$ then completely determines the other values  $\Delta_j$ with $j<s$ since
$\Z\la Z_1\cup\ldots\cup Z_j\ra=\Z\la Z_1\cup \ldots\cup Z_s\ra\cap \R\la Z_1\cup\ldots\cup Z_j\ra=\Delta_s\cap \R\la X_1\cup \ldots\cup X_j\ra=\Z\la X_1\cup \ldots\cup X_s\ra\cap \R\la X_1\cup \ldots\cup X_j\ra=\Z\la X_1\cup \ldots\cup X_j\ra$ as $Z$ and $X$ are linearly independent.

The linear type of $X\in X(G_0)$ associated to the maximal  series decomposition $X=X_1\cup\ldots\cup X_s$ is determined by the values of  $\pi_{j-1}(X_j)$ for $j\in [1,s]$, where $\pi_{j-1}:\R^d\rightarrow \R\la X_1\cup\ldots\cup X_{j-1}\ra^\bot$ is the orthogonal projection. The sets $\phi_{j-1}(X_j)$ are each  irreducible by Proposition \ref{prop-finitary-MaxDecompChar}, while Theorem \ref{thm-finitary-FiniteProps-I}.1 ensures that each finitary set $\pi_{j-1}(G_0)$ (by Proposition \ref{prop-finitary-Modulo-Inheritence}) only has a finite number of irreducible sets as possibilities for $\phi_{j-1}(X_j)$. It follows that there are only a finite number of linear types for $G_0$. The discussion after \eqref{lattice-splitting} ensures that each linear type only has a finite number of associated lattice types. Thus the number of lattice types for $G_0$ is also finite.
Indeed, per the discussion after \eqref{lattice-splitting},  the value of the lattice $\Delta_j=\Z\la X_1\cup\ldots\cup X_j\ra$ is uniquely determined by $\Delta_{j-1}$ and the values of $\xi_x\mod \Delta_{j-1}$ for $x\in X_j$, where $$\xi_x:=x-\phi(x)\in \Lambda\cap \R\la \Delta_{j-1}\ra.$$
Note $\xi_{y}=0$ for all $y\in Y$, ensuring $\xi_{\phi(x)}=0$ for all $x\in X$.
Thus, given a lattice type $\phi$ having  $Z=Z_1\cup\ldots\cup Z_s$ as the associated maximal series decomposition of its codomain, there are elements $\xi_z\in \Lambda\cap \R\la \Delta_{j-1}\ra$ for $z\in Z$ such that, if  $X\in X(G_0)$ has the same linear type as $\phi$, then $X$ also has the same lattice type as $\phi$ precisely when  \be\label{gargoyle}\xi_x\equiv \xi_{\phi(x)}=0\mod \Delta_{j-1}\quad\mbox{ for every $x\in X_j$ and $j\in [1,s]$},\ee where $X=X_1\cup\ldots\cup X_s$ is the associated maximal series decomposition and $\Delta_{j-1}=\Z\la Z_1\cup\ldots\cup Z_{j-1}\ra$.

We let
$\mathfrak T(G_0)$ consist of the set of possible lattice types $\phi$ for the sets  $X\in X(G_0)$.
Given $\phi\in \mathfrak T(G_0)$, let $$\mathfrak X(\phi)\subseteq X(G_0)$$ consist of all $X\in X(G_0)$ that have lattice type $\phi$ and  set $$\mathfrak X^\cup(\phi)=\bigcup_{X\in \mathfrak X(\phi)}X,$$ so each $\phi\in \mathfrak T(G_0)$ is a map $\phi:\mathfrak X^\cup (\phi)\rightarrow Z$, where $Z\in \mathfrak X(\phi)$ is a distinguished set with type $\phi$ having associated series decomposition $Z=Z_1\cup\ldots\cup Z_s$. Note the lattice types associated to $X\in X(G_0)$ are in bijective correspondence with the maximal series decompositions of $X$.

Suppose $X,\,X'\in \mathfrak X(\phi)$ and $x\in X$ and $x'\in X'$ are elements with $\phi(x)=\phi(x')$. Define a new set $X''=X\setminus\{x\}\cup \{x'\}\subseteq G_0$. Then $X''\in \mathfrak X(\phi)$ as well, as the following short argument shows.
We must have $x\in X_j$ and $x'\in X'_j$ for some $j\in [1,s]$, where $X=X_1\cup \ldots\cup X_s$ and $X'=X'_1\cup \ldots\cup X'_s$ are  the associated maximal series decompositions.
Since $X$ and $X'$ both have the same type, we have $\Z\la X_1\cup \ldots\cup X_{j-1}\ra=\Delta_{j-1}=\Z\la X'_1\cup \ldots\cup X'_{j-1}\ra$, and since $\phi(x)=\phi(x')$, it follows that $\xi_x\equiv \xi_{x'}\equiv 0\mod \Delta_{j-1}$, for $j\in [1,s]$, ensuring that  \be\label{rike}\Delta_j=\Z\la X_1\cup\ldots\cup X_j\ra=\Z\la X_1\cup \ldots\cup X_{j-1}\cup (X_j\setminus\{x\}\cup \{x'\})\ra\quad\mbox{ for $j\in [1,s]$},\ee since the values of $\Delta_{j-1}$ and each $\xi_x\mod \Delta_{j-1}$ for $x\in X_j$ determine the lattice $\Delta_j$.
 Let $\mathcal R=(\mathcal X_1\cup \{\mathbf v_1\},\ldots,\mathcal X_s\cup \{\mathbf v_s\})$ be a realization of $X\in X(G_0)$ associated to the maximal series decomposition $X=X_1\cup \ldots\cup X_s$ with all $\x\in \mathcal X=\mathcal X_1\cup\ldots\cup \mathcal X_s$ having dimension one.
 Note $\mathcal R$ exists per the discussion at the beginning of Section \ref{sec-series-decomp}, and by this same discussion, replacing the half-space $\x=\R_+ x$ with $\x'=\R_+ x'$ results in another virtual Reay system (after possibly modifying the $\mathbf v_i$ for $i\in [j+1,s]$), showing that $X_1\cup\ldots \cup X_{j-1}\cup (X_j\setminus \{x\}\cup \{x'\})\cup X_{j+1}\cup \ldots\cup X_s$ is a series decomposition of $X''$, which is maximal in view of Proposition \ref{prop-finitary-MaxDecompChar} since $X=X_1\cup\ldots\cup X_s$ was maximal and $\pi_{j-1}(x)=\pi_{j-1}(x')$.
 In view of \eqref{rike}, it now follows that $X''$ also has type $\phi$, as claimed.

What the above important observation means is, if $X,\,X'\in \mathfrak X(\phi)$ have the same type $\phi$, then it is possible to exchange elements between $X$ and $X'$ equal under the mapping $\phi$ and have the resulting set remain in the class $\mathfrak X(\phi)$. Indeed, if $Z_\phi\in X(\phi)$ is the representative for the type $\phi$ and  $X\subseteq \mathfrak X^\cup(\phi)=\bigcup_{Z\in \mathfrak X(\phi)}Z$ is \emph{any} subset consisting of precisely one element $x\in X$ with $\phi(x)=z$ for each $z\in Z_\phi$, then $X\in \mathfrak X(\phi)$ (as $X$ can be constructed by a process of at most $|Z_\phi|-1$ exchanges just described). We call this  the \textbf{interchangeability property} of $\mathfrak X(\phi)$, which essentially amounts to saying each $\mathfrak X^\cup(\phi)$ has the structure of a direct product. Thus $$X(G_0)=\bigcup_{\phi\in \mathfrak T(G_0)}\mathfrak X(\phi)$$ is a decomposition of $X(G_0)$ into a finite  number of exchange closed subsets (with the union not necessarily disjoint).

A lattice type corresponds to  a map $\phi:\mathfrak X^\cup(\phi)=\bigcup_{X\in \mathfrak X(\phi)}X\rightarrow Z_\phi$, where $Z_\phi\in \mathfrak X(\phi)$ is a distinguished representative.  The choice of representative $Z_\phi$ is rather arbitrary. A map $$\varphi:\bigcup_{X\in \mathfrak X_m(\varphi)} X\rightarrow Z_\varphi$$ obtained from $\phi$ by first changing the representative set used for the codomain from $Z_\phi$ to some $Z_\varphi\in \mathfrak X(\phi)$ and then restricting the domain to a union of sets from  some subset  $\mathfrak X_m(\varphi)\subseteq \mathfrak X(\phi)$ which contains the new codomain, so $Z_\varphi \in \mathfrak X_m(\varphi)$,  will be called a \textbf{refinement} of $\phi$. In such case, we set $$\mathfrak X^\cup(\varphi)=\bigcup_{X\in \mathfrak X_m(\varphi)} X$$ and note $\varphi:\mathfrak X^\cup(\varphi)\rightarrow Z_\varphi$.

In view of Proposition \ref{prop-finitary-mintype} below, we refer to the sets $X\in \mathfrak X_m(\varphi)$ as  having the same \textbf{minimal type} $\varphi$. Note, if $X$ and $X'$ both have the same minimal type $\varphi$, then they must also have the same lattice type, as the notion of minimal type refines that of lattice type, and by a small abuse of notation, we will often use $\varphi$ for this lattice type as well.
In what follows, we say an element $x\in X\in \mathfrak X_m(\varphi)$ is at \textbf{depth} $j\in [1,s]$ when $\varphi(x)\in Z_j$, where $Z_\varphi=Z_1\cup\ldots\cup Z_s$ is the maximal series decomposition associated to the codomain $Z_\varphi$ of the minimal type $\varphi$, and we call $$\mathsf{dep}(\varphi)=s$$ the \textbf{depth} of the minimal type $\varphi$. We likewise define all such terms, as well as $\mathsf{dep}(\phi)$, for a lattice type $\phi\in \mathfrak T(G_0)$.
The argument for the existence of the finite set $\mathfrak T_m(G_0)$ is similar to that used for Theorem \ref{thm-finitary-FiniteProps-III}.

\begin{proposition}\label{prop-finitary-mintype}
Let $\Lambda\subseteq \R^d$ be a full rank lattice, where $d\geq 0$, and let $G_0\subseteq \Lambda$ be a finitary subset with $\C(G_0)=\R^d$. Then there is a finite collection  $\mathfrak T_m(G_0)$ with each $\varphi\in \mathfrak T_m(G_0)$ being a map  $\varphi:\mathfrak X^\cup_m(\varphi)=\bigcup_{X\in \mathfrak X_m(\varphi)} X\rightarrow Z_\varphi$ satisfying the following properties:
\begin{itemize}
\item[(a)] $\mathfrak X(G_0)=\bigcup_{\varphi\in \mathfrak T_m(G_0)}\mathfrak X_m(\varphi)$.
\item[(b)] $\varphi$ is a refinement of some lattice type for $G_0$.
\item[(c)] $Z_\varphi\preceq _\Z X$ for every $X\in \mathfrak X_m(\varphi)$, i.e., $\C_\Z(Z_\varphi)\subseteq \C_\Z(X)$.
\item[(d)] For any $X\in \mathfrak X_m(\varphi)$ and $x\in X$ at depth $j\in [1,s]$, we have $\varphi(x)\in x+\C_\Z(Z_1\cup\ldots\cup Z_{j-1})$, where $Z_\varphi=Z_1\cup\ldots\cup Z_s$ is the maximal series decomposition associated to $Z_\varphi$.
\item[(e)] The interchangeability property holds for $\mathfrak X_m(\varphi)$, meaning, given any $X,\,X'\in \mathfrak X_m(\varphi)$ and $x\in X$ and $x'\in X'$ with $\varphi(x)=\varphi(x')$, we have $X\setminus \{x\}\cup \{x'\}\in \mathfrak X_m(\varphi)$.
\end{itemize}
\end{proposition}

\begin{proof}Since $\mathfrak T(G_0)$ is finite and $\mathfrak X(G_0)=\bigcup_{\phi\in \mathfrak T(G_0)}\mathfrak X(\phi)$, it suffices to show each  $\mathfrak X(\phi)$, for $\phi\in \mathfrak T(G_0)$,  can be written as a finite union of sets satisfying (b)--(e).
Let $\phi\in \mathfrak T(G_0)$ be an arbitrary lattice type with codomain $Y\in \mathfrak X(\phi)$ and associated maximal series decomposition $Y=Y_1\cup\ldots\cup Y_s$, so $$\phi:\mathfrak X^\cup(\phi)\rightarrow Y_1\cup \ldots\cup Y_{s-1}\cup Y_s.$$ We proceed by induction on $s=0,1,\ldots,d$. If $s=0$, then (b)--(e) hold with $\varphi=\phi$ the trivial lattice type and $\mathfrak X_m(\varphi)=\mathfrak X(\phi)=\{\emptyset\}$. If $s=1$, then (b)--(e) hold with $\varphi=\phi$,  $\mathfrak X_m(\varphi)=\mathfrak X(\phi)=\{Y\}=\{Y_1\}$ and $Z_\varphi=Y$. Therefore we may assume $s\geq 2$.  Let $X\in \mathfrak X(\phi)$ be arbitrary with associated maximal series decomposition $X=X_1\cup\ldots\cup X_s$.
By definition of a lattice type, $$\Delta:=\Z\la Y_1\cup\ldots\cup Y_s\ra=\Z\la X_1\cup \ldots\cup X_s\ra\quad \und \quad \Delta':=\Z\la Y_1\cup\ldots\cup Y_{s-1}\ra=\Z\la X_1\cup \ldots\cup X_{s-1}\ra$$ are fixed (with  $\Delta=\Delta_s$ and $\Delta'=\Delta_{s-1}$  in the notation from the start of Section \ref{sec-finitary}.4).

By induction hypothesis, there are only a finite number of possible minimal types $\varphi'$ for the set $X_1\cup\ldots\cup X_{s-1}\in X(G_0)$ that are refinements of the lattice type defined for $Y\setminus Y_s$ by the maximal  series decomposition $Y\setminus Y_s=Y_1\cup\ldots\cup Y_{s-1}$ (maximality follows by Proposition \ref{prop-finitary-MaxDecompChar} since $Y=Y_1\cup\ldots\cup Y_s$ is maximal). Note any $X\in \mathfrak X(\varphi)$ with associated maximal series decomposition $X=X_1\cup \ldots\cup X_s$ has  $X\setminus X_s=X_1\cup\ldots\cup X_{s-1}$ of the same lattice type as $Y\setminus Y_s=Y_1\cup\ldots\cup Y_{s-1}$, as $X$ and $Y$ both have lattice type $\phi$. For any such minimal type $\varphi'$ (refining the lattice type defined by $Y_1\cup\ldots\cup Y_{s-1}$),  let $\mathfrak X_{\varphi'}\subseteq \mathfrak X(\phi)$ consist of all $X\in \mathfrak X(\phi)$ having associated maximal series decomposition $X=X_1\cup \ldots\cup X_s$ such that $X\setminus X_s\in \mathfrak X_m(\varphi')$.
Note the  condition that $\varphi'$ refine the lattice type defined by $Y_1\cup\ldots\cup Y_{s-1}$, and thus also that defined by $X_1\cup\ldots\cup X_{s-1}$, guarantees the subtle condition that  $\varphi'(x)=\varphi'(z)$ if and only if $\phi(x)=\phi(z)$, for any $x$ and $z$ in the common domain of both $\varphi'$ and $\phi$.
Now $\mathfrak X(\phi)=\bigcup_{\varphi'}\mathfrak X_{\varphi'}$ is a finite union of sets $\mathfrak X_{\varphi'}$ with all $X\in \mathfrak X_{\varphi'}$  having $X\setminus X_s$ of the same minimal type $\varphi'$,  so it suffices to show each $\mathfrak X_{\varphi'}$  can be written as a finite union of sets satisfying (b)--(e). To this end, let $\varphi'$ be arbitrary and suppose $Z_{\varphi'}=Z_1\cup\ldots\cup Z_{s-1}$ is the codomain of $\varphi'$, so $$\varphi':\mathfrak X_m^\cup (\varphi')\rightarrow Z_1\cup\ldots\cup Z_{s-1}.$$
In view of the interchangeability properties for $\mathfrak X(\phi)$ and $\mathfrak X(\varphi')$ and that guaranteed in (e) by induction for $\mathfrak X_m(\varphi')$, it follows that $\mathfrak X_{\varphi'}$ also has the interchangeability property. Note that $$\Z\la Z_{\varphi'}\ra=\Z\la Z_1\cup \ldots\cup Z_{s-1}\ra=\Delta'$$ as $\varphi'$ refined the lattice type associated to $Y_1\cup\ldots\cup Y_{s-1}$.

If $Z'_1\cup \ldots\cup Z'_{s-1}\in \mathfrak X_m(\varphi')\subseteq \mathfrak X(\varphi')$ and $X'_1\cup \ldots\cup X'_s\in \mathfrak X(\phi)$, then we have $\Z\la Z'_1\cup\ldots\cup Z'_{s-1}\ra=\Delta'=\Z\la Y_1\cup\ldots\cup Y_{s-1}\ra=\Z\la X'_1\cup \ldots\cup X'_{s-1}\ra$, with the first equality as $\varphi'$ refines the lattice type associated to $Y_1\cup \ldots\cup Y_{s-1}$.
Thus \be\label{working}\Delta=\Z \la X'_1\cup \ldots\cup X'_{s-1}\ra+\Z\la X'_s\ra=\Z\la Z'_1\cup \ldots\cup Z'_{s-1}\ra +\Z\la X'_s\ra=\Z\la Z'_1\cup \ldots\cup Z'_{s-1}\cup X'_s\ra.\ee Since $\varphi'$ refines the lattice type associated to $X_1\cup \ldots\cup X_{s-1}$ and $X'_1\cup \ldots\cup X'_s\in \mathfrak X(\phi)$, it follows that $Z'_1\cup \ldots\cup Z'_{s-1}\cup X'_s$ has the same linear type as $\phi$, which combined with \eqref{working} ensures that $Z'_1\cup \ldots\cup Z'_{s-1}\cup X'_s$ has the same lattice type as $\phi$.
Hence $Z'_1\cup \ldots\cup Z'_{s-1}\cup X'_s\in \mathfrak X(\phi)$ with $Z'_1\cup \ldots\cup Z'_{s-1}\in \mathfrak X_m(\varphi')$, ensuring by definition of $\mathfrak X_{\varphi'}$ that $Z'_1\cup \ldots\cup Z'_{s-1}\cup X'_s\in \mathfrak X_{\varphi'}$. Applying this with $Z'_1\cup \ldots\cup Z'_{s-1}=Z_1\cup \ldots\cup Z_{s-1}$ and $X'_1\cup \ldots\cup X'_s=Y_1\cup \ldots\cup Y_s$, we conclude that
 \be\label{YZ-in}Y_{\varphi'}:=Z_1\cup\ldots\cup Z_{s-1}\cup Y_s\in \mathfrak X_{\varphi'}.\ee

Since all $X\in \mathfrak X_{\varphi'}\subseteq \mathfrak X(\phi)$ have the same lattice type,
it follows by \eqref{gargoyle} that  any $x\in \mathfrak X^\cup(\phi)$ with $\mathsf{dep}(x)=s$ has  \be\label{wrench}\xi_x:=x-\phi(x)\equiv 0\mod \Delta'.\ee
 For $y\in Y_s$, define $$\wtilde X_y=\{\xi_x: \mbox{ there is some }X\in \mathfrak X_{\varphi'}\und x\in X\mbox{ with }\phi(x)=y\}\subseteq \Delta'.$$
 By definition of $\xi_x$ and $\wtilde X_y$, we find that  \be\label{transfer}\xi\in \wtilde X_y\quad\mbox{ implies }\quad \xi+y=\xi_x+\phi(x)=x\in\mathfrak X^\cup(\varphi)\subseteq G_0.\ee Let $\mathcal R_{\varphi'}$ be a realization of $Z_{\varphi'}=Z_1\cup\ldots\cup Z_{s-1}$ in which all anchored half-spaces have dimension one.
If $\{-\xi_i\}_{i=1}^\infty$ is an asymptotically filtered sequence of terms from $-\wtilde X_y$ with fully unbounded limit $\vec u=(u_1,\ldots,u_t)$, then $u_1,\ldots,u_t\in \R\la \Delta'\ra=\R\la Z_{\varphi'}\ra$, and then \eqref{transfer} ensures $\{\xi_i+y\}_{i=1}^\infty$ is  an asymptotically  filtered sequence of terms from $G_0$ with fully unbounded limit $-\vec u$.
Thus Proposition \ref{prop-finitary-basics}.3  applied to $\mathcal R_{\varphi'}$ ensures that $\vec u$ is encased by $Z_{\varphi'}$, in which case Theorem \ref{thm-nearness-characterization}.4 applied to $\C(Z_{\varphi'})$ implies that $-\wtilde X_y$ is bound to $\C(Z_{\varphi'})$, and thus $-\wtilde X_y$ is also bound  to $\C_{\Z}(Z_{\varphi'})$.
This means there is some fixed ball $B_y\subseteq \mathcal \R\la Z_{\varphi'}\ra$ such that every $\xi\in \wtilde X_y$ has some $z\in \C_{\Z}(Z_{\varphi'})$ with  $-\xi\in B_y+z$. As $\xi,\,z\in \Delta'$, we are assured  that $-\xi-z\in \Delta'\cap B_y$, which is a bounded set of lattice points, and thus finite. Consequently, \be\mbox{each $\xi \in \wtilde X_y$ has $\varpi+\xi\in \C_{\Z}(-Z_{\varphi'})$ for some $\varpi\in \Delta'\cap  B_y$}.\label{texttry}\ee
For each $y\in Y_s$, there are only a finite number of choices for an element $\varpi_y\in \Delta'\cap B_y$. For each fixed choice of elements $\varpi_y\in \Delta'\cap B_y$ for each $y\in Y_s$, we can define a subset $\Omega\subseteq \mathfrak X_{\varphi'}$ consisting of all $X\in \mathfrak X_{\varphi'}$ such that, whenever $x\in X$ with $\phi(x)=y\in Y_s$, then $\varpi_y+\xi_x\in \C_\Z(- Z_{\varphi'})$. In view of \eqref{texttry}, we see that $\mathfrak X_{\varphi'}$ is the union of all possible $\Omega$, as we range over the finite number of possible choices for the $\varpi_y$. It thus suffices to show (b)--(e) hold for each possible $\Omega\subseteq \mathfrak X_{\varphi'}$.
Fix one such arbitrary possibility for $\Omega$. For $y\in Y_s$, define
$$\wtilde X_y^\Omega=\{\xi_x: \mbox{ there is some }X\in \Omega\und x\in X\mbox{ with }\phi(x)=y\}\subseteq \wtilde X_y.$$  Then, by definition of $\Omega$, we are assured that, for all $y\in Y_s$, we have  \be\label{transferII}\varpi_y+\wtilde X_y^\Omega\subseteq \C_\Z(-Z_{\varphi'})\;\und\; y+\wtilde X_y^\Omega\subseteq \{x\in G_0:\; x\in X\mbox{ for some } X\in \Omega\und \phi(x)=y\},\ee with the latter inclusion above following by noting that $\xi\in \wtilde X_y^\Omega$ implies $\xi=x-\phi(x)=x-y$ for some $x\in X$  with $X\in \Omega$ (by \eqref{wrench} and definition of $\wtilde X_y^\Omega$).

 The restriction added to a set $X=X_1\cup \ldots\cup X_s$ when passing from $\mathfrak X_{\varphi'}$ to $\Omega\subseteq \mathfrak X_{\varphi'}$ only applies to conditions involving $X_s$, and the restrictions imposed for each $x\in X_s$ are independent of each other. Thus, since the  interchangeability property holds for  $\mathfrak X_{\varphi'}$, it follows
that it also holds for $\Omega$, i.e., if $X,\,X'\in \Omega$  and $x\in X$ and $x'\in X'$ with $\phi(x)=\phi(x')$, then $X\setminus \{x\}\cup \{x'\}\in \Omega$. Moreover, if $X'=X'_1\cup\ldots\cup X'_s\in \mathfrak X_{\varphi'}$ and $X=X_1\cup \ldots\cup X_s\in \Omega$, then $X'_1\cup \ldots\cup X'_{s-1}\cup X_s\in \Omega$.
In particular, the case $X'= Y_{\varphi'}\in \mathfrak X_{\varphi'}$ (by \eqref{YZ-in}) tell us  that,
if $X=X_1\cup\ldots\cup X_{s-1}\cup X_s\in \Omega$, then $Z_1\cup\ldots\cup Z_{s-1}\cup X_s\in \Omega$. If $\xi\in \wtilde X^\Omega_y$, then \eqref{transferII} implies there is some $X\in \Omega$ with $$y+\xi=\phi(x)+\xi=x\in X=X_1\cup\ldots\cup X_{s-1}\cup X_s,$$ and  by the  previous conclusion, we can w.l.o.g. assume $X=Z_1\cup \ldots\cup Z_{s-1}\cup X_s$ (note $x\in X_s$ as $\phi(x)=y\in Y_s$). As a result, if, for each $y\in Y_s$, we have some $\xi_y\in \wtilde X_y^\Omega$, then swapping the elements $y+\xi_y$ into the series decomposition $X=Z_1\cup \ldots\cup Z_{s-1}\cup X_s$, one by one, yields a series decomposition  $Z_1\cup\ldots\cup Z_{s-1}\cup Z_s\in \Omega$ with $Z_s=\{y+\xi_y:\;y\in Y_s\}$.

The natural partial order on $\C_\Z(-Z_{\varphi'})$, by declaring $a\preceq b$ for $a,b\in \C_\Z(-Z_{\varphi'})$ when $b\in a+\C_\Z(-Z_{\varphi'})$, makes $\C_\Z(-Z_{\varphi'})$ into a partially ordered set isomorphic to $\Z_+^{|Z_{\varphi'}|}$ (since $Z_{\varphi'}$ is a linearly independent set). Since $\varpi_y+\wtilde X^\Omega_y\subseteq \C_\Z(-Z_{\varphi'})$ by \eqref{transferII}, Hilbert's Basis Theorem (see Section \ref{sec-poset-lattices}) ensures that $\mathsf{Min} (\varpi_y+\wtilde X^\Omega_y)$ is finite with $\uarrow \mathsf{Min} (\varpi_y+\wtilde X^\Omega_y)=\varpi_y+\wtilde X^\Omega_y$, for each $y\in Y_s$.
For each $y\in Y_s$, we can choose some minimal element $\varpi_y+\xi_y\in \mathsf{Min}(\varpi_y+\wtilde X^\Omega_y)$ and then set $$z_y=y+\xi_y.$$ Letting $Z_s=\{z_y:\;y\in Y_s\}$, the discussion of the previous paragraph ensures that $$Z_\varphi:=Z_1\cup\ldots\cup Z_{s-1}\cup Z_s\in \Omega.$$
 As $\mathsf{Min}(\varpi_y+\wtilde X^\Omega_y)$ is finite, there are only a finite number of possibilities for how to construct the set $Z_s$. For each possible $Z_s$, we can define a subset $\mathfrak X_m(\varphi)\subseteq \Omega$ consisting of all $X\in \Omega$ such that, for every $y\in Y_s$ and $x\in X$ with $\phi(x)=y$, we have $\varpi_y+\xi_y\preceq \varpi_y+\xi_x$. Changing the codomain of $\phi$ to $Z_\varphi$ and restricting to the domain $\bigcup_{X\in \mathfrak X_m(\varphi)}X$ gives rise to a refinement $\varphi:\bigcup_{X\in \mathfrak X_m(\varphi)}X\rightarrow Z_\varphi$ of $\phi$. We now have $\Omega$ written as a finite union of sets $\mathfrak X_m(\varphi)$ such that (b) holds.

Let $\varphi$ be an arbitrary possible restriction defined above, let $X\in \mathfrak X_m(\varphi)$ be arbitrary and let $x\in X$. If $x$ is at depth $j<s$, then, since $X\in \mathfrak X_m(\varphi)\subseteq \Omega\subseteq X_{\varphi'}$, we have $X_1\cup \ldots\cup X_{s-1}\in \mathfrak X_m(\varphi')$ satisfying (d), meaning $\varphi(x)=\varphi'(x)\in x+\C_\Z(Z_1\cup\ldots\cup Z_{j-1})$. If $x$ is at depth $s$, then $\varpi_{y}+\xi_y\preceq \varpi_y+\xi_x$, where $y=\phi(x)$. Thus $\xi_x\in \xi_{y}-\C_{\Z}(Z_1\cup\ldots\cup Z_{s-1})$, in turn implying $\varphi(x)=z_y=\xi_y+y\in \xi_x+y+\C_{\Z}(Z_1\cup \ldots\cup Z_{s-1})=x+\C_{\Z}(Z_1\cup \ldots\cup Z_{s-1})$, with the final equality by \eqref{wrench}.  In particular, (d) holds, and we have $z\in \C_\Z(Z_1\cup \ldots\cup Z_{s-1})+\C_\Z(X_s)$ for every $z\in Z_{\varphi}$ (as $\varphi(X_s)=Z_s$). Since
$X\in \mathfrak X_m(\varphi)\subseteq \Omega\subseteq X_{\varphi'}$, we have $X_1\cup \ldots\cup X_{s-1}\in \mathfrak X_m(\varphi')$ satisfying (c), meaning $\C_{\Z}(Z_1\cup \ldots\cup Z_{s-1})\subseteq \C_\Z(X_1\cup\ldots\cup  X_{s-1})$, whence $z\in \C_{\Z}(Z_1\cup \ldots\cup Z_{s-1})+\C_\Z(X_s)\subseteq \C_\Z(X)$ for every $z\in Z_\varphi$. Thus $Z_\varphi\subseteq \C_\Z(X)$, implying $\C_\Z(Z_\varphi)\subseteq \C_\Z(X)$, showing that
 (c) holds for $\mathfrak X_m(\varphi)$.
 It remains to show that the interchangeability property holds for $\mathfrak X_m(\varphi)$. However,
the restriction added to a set $X=X_1\cup \ldots\cup X_s$ when passing from $\Omega$ to  $\mathfrak X_m(\varphi)\subseteq \Omega$ only applies to conditions involving $X_s$, and the restrictions imposed for each $x\in X_s$ are independent of each other. Thus, since the  interchangeability property holds for  $\Omega$, it follows
that it also holds for $\mathfrak X_m(\varphi)$, i.e., if $X,\,X'\in \mathfrak X_m(\varphi)$  and $x\in X$ and $x'\in X'$ with $\varphi(x)=\varphi(x')$, and hence $\phi(x)=\phi(x')$ too, then $X\setminus \{x\}\cup \{x'\}\in \mathfrak X_m(\varphi)$. Thus (e) holds, which completes the induction and proof.
\end{proof}

\section{Factorization Theory}\label{sec-fact}

\subsection{Lambert Subsets and Elasticity} In this final section, we apply the machinery regarding Convex Geometry and finitary sets developed in prior sections to derive some striking consequences regarding the behavior of factorizations.
Our first goal, which is one of the most difficult and crucial steps in the characterization of finite elasticities, is to obtain a multi-dimensional generalization of a result of Lambert for subsets of $\Z$. After we have achieved this, we then characterize when $\rho_{d+1}(G_0)$ is finite. Recall that $\vp_X(S)=\Summ{x\in X}\vp_x(S)$ for $S\in \Fc(G_0)$ and $X\subseteq G_0$.

\begin{definition}Let $G$ be an abelian group and  let $G_0\subseteq G$ be a subset. We say that  $X\subseteq G_0$ is a \textbf{Lambert} subset if there exists a bound $N\geq 0$ such that $\vp_X(U)=\Summ{x\in X}\vp_x(U)\leq N$ for all $U\in \mathcal A(G_0)$.
\end{definition}

Lambert \cite{Lambert} (see also \cite[Lemma 4.3]{Gerold-lambert-rankone} and \cite[Theorem 3.2]{Wolfgang-lambert}) showed that, if $G_0\subseteq \Z$ is a subset with $G_0\cap -\Z_+$ finite, then $\Summ{x\in G_0\cap \Z_+}\vp_x(U)\leq N$ for every $U\in \mathcal A(G_0)$, where $N=|\min (G_0\cap -\Z_+)|$, which,
in terms of the notation just introduced,  means $G_0\cap \Z_+\subseteq G_0$ is a Lambert subset. His proof was a clever adaptation the well-known argument for obtaining the basic upper bound $\mathsf D(G)\leq |G|$ for the Davenport Constant \cite[Theorem 10.2]{Gbook}.   The hypothesis $|G_0\cap -\Z_+|<\infty$ readily implies that $G_0$ is finitary. Indeed, it is essentially a characterization of being a finitary set in $\Z$ (either $G_0\cap -\Z_+$ or $G_0\cap \Z_+$ must be finite for a subset of $\Z$ to be finitary). Moreover, if $G_0$ is infinite, then $G_0^\diamond =(G_0\cap -\Z_+)\setminus\{0\}$, so Lambert's conclusion can be reworded as saying $G_0\setminus G_0^\diamond\subseteq G_0\subseteq \Z$ is a Lambert subset.

Consider a more general subset $G_0\subseteq \Lambda$, where $\Lambda\leq \R^d$ is a lattice and $d\geq 0$.
If $\C(G_0)\neq \R^d$, then terms from $G_0\setminus \mathcal E$, where $\mathcal E=\C(G_0)\cap -\C(G_0)$ is the lineality space of $\C(G_0)$, are contained in no atom, and can essentially be discarded. Assuming $\C(G_0)=\R^d$, then  Proposition \ref{prop-diamond-basic-containment} implies $X\subseteq G_0\setminus G_0^\diamond$ for any Lambert subset $X\subseteq G_0$. We will later in this subsection see that, when $G_0$ is finitary, then $G_0\setminus G_0^\diamond$ is the unique maximal Lambert subset of $G_0$, giving  a multi-dimensional analog of Lambert's result. However, we
begin first with two propositions entailing some of the basic relationships between finitary sets and finite elasticities.

\begin{proposition}\label{prop-lambert-easy}
Let $G$ be an abelian group and let  $G_0\subseteq G$ be a subset. Suppose $X\subseteq G_0$ is a subset with $\mathcal A(X)=\emptyset$ and $G_0\setminus X\subseteq G_0$ is a Lambert subset with bound $N$. Then $\rho(G_0)\leq N<\infty$. In particular, $\rho_k(G)\leq \rho(G_0)k\leq Nk<\infty$ for any $k\geq 1$.
\end{proposition}

\begin{proof}
Let $U_1,\ldots,U_k\in\mathcal A(G_0)$ be atoms. Suppose $U_1\bdot\ldots\bdot U_k=V_1\bdot\ldots\bdot V_\ell$ for some $V_1,\ldots,V_\ell\in \mathcal A(G)$. Since $\mathcal A(X)=\emptyset$, it follows that each $V_i$ must contain a term from $G_0\setminus X$. However, since $G_0\setminus X\subseteq G_0$ is a  Lambert subset with bound $N$, there are at most $kN$ terms in $U_1\bdot \ldots\bdot U_k$ from $G_0\setminus X$. Thus $\ell\leq kN$, implying $\ell/k\leq N$, and the proposition follows by definition of $\rho(G_0)$.
\end{proof}

\begin{proposition}\label{prop-pre-rho-char}
Let $G$ be an abelian group with torsion-free rank $d\geq 0$ and let  $G_0\subseteq G$ be a subset. Then, regarding the statements below, we have the implications $1.\Rightarrow 2.\Rightarrow 3. \Rightarrow 4.$
\begin{itemize}
\item[1.] There exists a subset $X\subseteq G_0$ such that $\mathcal A(X)=\emptyset$ and $G_0\setminus X\subseteq G_0$ is a Lambert subset.
\item[2.] $\rho(G_0)<\infty$.
\item[3.] $\rho_k(G_0)<\infty$ for all $k\geq 1$.
\item[4.] $\rho_{d+1}(G_0)<\infty$.
\end{itemize}
\end{proposition}

\begin{proof}
The implication $1.\Rightarrow 2.$ follows by Proposition \ref{prop-lambert-easy}, while the implications
$2.\Rightarrow 3.$ and $3.\Rightarrow 4.$ follow from the definition of $\rho(G_0)$.
\end{proof}

\begin{proposition}\label{prop-pre-rho-char-diamond}
Let $\Lambda\leq \R^d$ be a full rank lattice in $\R^d$, where $d\geq 0$, and  let $G_0\subseteq \Lambda$ be a subset with $\C(G_0)=\R^d$. If $\rho_{d+1}(G_0)<\infty$, then $0\notin \C^*(G_0^\diamond)$. In particular, $G_0$ is finitary.
\end{proposition}

\begin{proof}
Note $0\notin G_0^\diamond$ by Proposition \ref{prop-G_0diamond-1st-easy-equiv}.2. Assume by contradiction that $0\in \C^*( G_0^\diamond)$. Then Carath\'eordory's Theorem implies that there is a minimal positive basis $X\subseteq G_0^\diamond$. Let $X=\{x_1,\ldots,x_s\}$ be the distinct elements of $X$, where $2\leq s\leq d+1$, let $V\in \mathcal A^{\mathsf{elm}}(G_0)$ be the unique elementary atom with $\supp(V)=X$ (by Proposition \ref{prop-char-minimal-pos-basis}), and let $N=\max\{\vp_{x_j}(V):\;j\in [1,s]\}$.  Since each $x_j\in G_0^\diamond$, Proposition \ref{prop-diamond-basic-containment} implies that, for each $j\in [1,s]$, there is a sequence $\{U_i^{(j)}\}_{i=1}^\infty$ of atoms $U^{(j)}_i\in\mathcal A^{\mathsf{elm}}(G_0)$ with $\vp_{x_j}(U_i^{(j)})\rightarrow \infty$. Let $M_i=\min\{\vp_{x_j}(U_i^{(j)}):\; j\in [1,s]\}$ and observe that $M_i\rightarrow \infty$ since each $\vp_{x_j}(U_i^{(j)})\rightarrow \infty$.  Consider the product $W_i=U_i^{(1)}\bdot\ldots\bdot U_i^{(s)}$ for $i\geq 1$ and let $m_i= \lfloor \frac{M_i}{N}\rfloor$.
Then $W_i=V^{[m_i]}\bdot W'_i$ is a factorization of $W_i$, where $W'_i=W_i\bdot V^{[-m_i]}\in \mathcal B(G_0)$ and $V\in \mathcal A(G_0)$. Moreover, since $N$ is fixed and $M_i\rightarrow \infty$, it follows that $m_i\rightarrow \infty$, showing that $\rho_s(G_0)=\infty$. However, since $s\leq d+1$, we have $\rho_s(G_0)\leq \rho_{d+1}(G_0)$ by \eqref{rho-ascend-chain}, forcing $\rho_{d+1}(G_0)=\infty$ as well, contrary to hypothesis. Note $0\notin \C^*(G_0^\diamond)$ implies $G_0$ is finitary  by Theorem \ref{thm-keylemmaII}.
\end{proof}

We will need the following lemma for our generalization of Lambert's result.

\begin{lemma}\label{lem-finitary-diamond-containment-nonmax} Let $\Lambda\subseteq \R^d$ be a full rank lattice, where $d\geq 0$,  let
$G_0\subseteq \Lambda$ be a finitary subset with $\C(G_0)=\R^d$, let $\mathcal R=(\mathcal X_1\cup \{\mathbf v_1\},\ldots,\mathcal X_s\cup\{\mathbf v_s\})$ be a purely virtual Reay system over $G_0$, let $\mathcal X=\mathcal X_1\cup \ldots\cup \mathcal X_s$,  and let $\pi:\R^d\rightarrow \R^\cup \la \mathcal X\ra^\bot$ be the orthogonal projection.
 If $U\in \Fc_{\mathsf{rat}}(G_0)$ with $\pi(U)=W^{[\alpha]}$ for some rational number $\alpha>0$ and $W\in \mathcal A^{\mathsf{elm}}(\pi(G_0))$,  and there exists some $w\in \supp(W)$ such that $\pi^{-1}(w)\cap \supp(U)\subseteq G_0\setminus G_0^\diamond$, then $-\sigma(U)\in\C^\cup (\mathcal X)$.
\end{lemma}

\begin{proof}
Suppose $U$ is a counter-example with $|\supp(U)|$ minimal. By replacing $U$ with $U^{[1/\alpha]}$, we can w.l.o.g. assume $\alpha=1$. Let $w_1,\ldots,w_t\in \supp(W)$ be the distinct elements of $\supp(W)$ and let $\alpha_i=\vp_{w_i}(W)\in \Z_+$ for $i\in [1,t]$. Since $W\in \mathcal A^{\mathsf{elm}}(\pi(G_0))$, it follows that $\{w_1,\ldots,w_t\}$ is either $\{0\}$ or  a minimal positive basis with $\alpha_1w_1+\ldots+\alpha_t w_t=0$.  By hypothesis, we have some $w_i$, say $w_1$, such that $\pi^{-1}(w_1)\cap \supp(U)\subseteq G_0\setminus G_0^\diamond$, i.e., every $u\in  \supp(U)$ with $\pi(u)=w_1$ has $u\in G_0\setminus G_0^\diamond$. Partition $\supp(U)=U_1\cup \ldots\cup U_t$ such that each $U_i$, for $i\in [1,t]$, consists of those $u\in \supp(U)$ with $\pi(u)=w_i$. Note $U_1\subseteq G_0\setminus G_0^\diamond$.

Let $X\subseteq \supp(U)$ be an arbitrary subset such that $|X\cap U_i|=1$ for every $i\in [1,t]$, say with $X\cap U_i=\{x_i\}$ for $i\in [1,t]$. Then $U_X:=\prod^\bullet_{i\in [1,t]}x_i^{[\alpha_i]}\in \Fc(G_0)$ with $\pi(U_X)=W$ and $|\supp(U_X)|=|\supp(W)|$. As a result, since $x_1\in U_1\subseteq G_0\setminus G_0^\diamond$, it follow from Proposition \ref{prop-finitary-diamond-containment-nonmax} applied to $U_X$ that $\mathsf{wt}(-\sigma(U_X))=0$, i.e., \be\label{seeit}-\sigma(U_X)\in \C^\cup (\mathcal X).\ee Let $$\beta=\min\{\frac{\vp_{x_i}(U)}{\vp_{x_i}(U_X)}:\;i\in [1,t]\}=\min\{\frac{\vp_{x_i}(U)}{\alpha_i}:\;i\in [1,t]\}>0.$$
Since $\vp_{x_i}(U_X)=\alpha_i=\Summ{u\in U_i}\vp_u(U)\geq \vp_{x_i}(U)$ for each $i\in [1,t]$ (as $\pi(U)=W$), we have $\beta\leq 1$ with equality only possible if $U_X^{[\beta]}=U_X=U$.

In view of the definition of $\beta$, we have $U_X^{[\beta]}\mid U$ with $\vp_{x_i}(U_X^{[\beta]})=\vp_{x_i}(U)$ for every $i\in [1,t]$ attaining the minimum in the definition of $\beta$. Thus $U_X^{[\beta]}\in \Fc_{\mathsf{rat}}(G_0)$ with $$|\supp(U\bdot U_X^{[-\beta]})|<|\supp(U)|\quad\und\quad \pi(U_X^{[\beta]})=W^{[\beta]}.$$
 By \eqref{seeit}, we have $-\sigma(U_X^{[\beta]})=-\beta\sigma(U_X)\in \C^\cup (\mathcal X)$ (since  $\C^\cup (\mathcal X)$ is a convex cone by Proposition \ref{prop-orReay-BasicProps}.2). Consequently, if $\beta=1$, so that $U_X^{[\beta]}=U$, then  $-\sigma(U)=-\sigma(U_X^{[\beta]})\in \C^\cup(\mathcal X)$, as desired. On the other hand, if $\beta<1$, then we have
 $U\bdot U_X^{[-\beta]}\in \Fc_{\mathsf{rat}}(G_0)$ with $\pi(U\bdot U_X^{[-\beta]})=W^{[1-\beta]}$. Thus the hypotheses hold using $U\bdot U_X^{[-\beta]}$, in which case the minimality of  $|\supp(U)|$ for the counter-example $U$ ensures that $-\sigma(U\bdot U_X^{[-\beta]})\in \C^\cup (\mathcal X)$. Combined with \eqref{seeit} and the convexity of $\C^\cup(\mathcal X)$ (Proposition \ref{prop-orReay-BasicProps}.2), we conclude that $-\sigma(U)=-\sigma(U_X^{[\beta]})-\sigma(U\bdot U_X^{[-\beta]})\in \C^\cup (\mathcal X)$, as desired.
\end{proof}

We now come to the key generalization of Lambert's result. Theorem \ref{thm-structural-Lambert} requires the hypothesis $0\notin\C^*(G_0^\diamond)$, which is slightly stronger than being finitary (in view of Theorem \ref{thm-keylemmaII}). In exchange, we actually attain a decomposition of an arbitrary atom $U\in \mathcal A(G_0)$, reminiscent of Theorem \ref{thm-carahtheodory-elm-atom} (Carath\'eordory's Theorem), which implies that $G_0\setminus G_0^\diamond\subseteq G_0$ is a Lambert subset. A simplification of the argument used to prove Theorem \ref{thm-structural-Lambert} works when we only have $G_0$ being finitary, though it only yields that $G\setminus G_0^\diamond\subseteq G_0$ is a Lambert subset and not the additional decomposition result for  atoms. We deal with this in more  detail afterwards in Corollary \ref{cor-structural-Lambert}. The proof of Theorem \ref{thm-structural-Lambert} is algorithmic and yields recursively defined values for the constants $N_S$ and $N_T$, which are quite large and dependent on the structure of the individual set $G_0$. We have not  optimized the estimates for $N_S$ and $N_T$, instead opting for arguments that simplify the recursive definitions and presentation.  We remark that the proof of Theorem \ref{thm-structural-Lambert} does not require minimal types, that is:

\smallskip

\emph{Theorem \ref{thm-structural-Lambert} remains valid  when $\varphi_1,\ldots,\varphi_\mathfrak s\in \mathfrak T(G_0)$ are the distinct nontrivial lattice types, rather than the distinct nontrivial minimal types}.

\smallskip

\noindent In such case $\mathfrak s+1=|\mathfrak T(G_0)|$ rather than $|\mathfrak T_m(G_0)|$, which is in general smaller, thus requiring less iterations of the algorithm, and so yielding improved bounds for  $N_T$ and $N_S$. However, we will need the added refinements provided by  minimal types later for Theorem \ref{thm-structural-char}.

\begin{theorem}\label{thm-structural-Lambert}
Let $\Lambda\leq \R^d$ be a full rank lattice, where $d\geq 0$,  and let $G_0\subseteq \Lambda$ be a subset with  $\C(G_0)=\R^d$. Suppose $0\notin\C^*(G_0^\diamond)$ and let $\varphi_1,\ldots,\varphi_\mathfrak s\in \mathfrak T_m(G_0)$ be the distinct nontrivial minimal types $\varphi_j:\mathfrak X_m^\cup(\varphi_j)\rightarrow Z_{\varphi_j}$ indexed so that $\mathsf{dep}(\varphi_i)\leq \mathsf{dep}(\varphi_j)$ whenever
$i\leq j$.
Then there are bounds  $N_S\geq 0$ and $N_T\geq 0$ such that any atom $U\in \mathcal A(G_0)$ has a factorization $$U=A_0\bdot A_1\bdot\ldots\bdot A_\mathfrak s,\quad\mbox{ where }\quad A_i\in \mathcal B_{\mathsf{rat}}(G_0),$$ such that,  for every $j\in[ 1,\mathfrak s]$,  \begin{itemize}
\item[(a)] no subset $X\subseteq \supp(A_0\bdot\ldots\bdot A_{j-1})$ with $X\in X(G_0)$ has minimal type $\varphi_j$;
\item[(b)]  $|\supp(\{A_0\})|\leq N_S$ and $|\supp(\{A_j\})|\leq N_S$;
\item[(c)] $|A_0|\leq N_T$ and $\vp_{G_0\setminus \mathfrak X^\cup_m(\varphi_j)}(A_j)=\Summ{x\in G_0\setminus \mathfrak X^\cup_m(\varphi_j)}\vp_x(A_j)\leq N_T$; and
\item[(d)]  $\Summ{x\in G_0\setminus \mathfrak X^\cup_m(\varphi_j)}\vp_x(A_j)x\in-\C(Z_{\varphi_j})$.
\end{itemize}
\end{theorem}

\begin{proof} If $d=0$, then there are no nontrivial types and $U=0$ is the only possible atom, in which case the theorem holds with $\mathfrak s=0$, \ $N_S=0$ and $N_T=1$. Therefore we may assume $d\geq 1$.
Since $0\notin \C^*(G_0^\diamond)$, it follows from  Theorem \ref{thm-keylemmaII} that $G_0$ is finitary.  We construct the bounds $N_S$ and $N_T$ as well as the subsequences $A_j$ inductively for $j=\mathfrak s,\ldots,1,0$. Assume, for some $\mathfrak t\in [1,\mathfrak s+1]$, we have already constructed bounds $N_S\geq 0$ and $N_T\geq 0$ so that, for any atom $U\in \mathcal A(G_0)$, we can find a  factorization $$U=V\bdot A_{\mathfrak t}\bdot\ldots\bdot A_\mathfrak s,\quad\mbox{ where }\quad V,\,A_i\in \mathcal B_{\mathsf{rat}}(G_0),$$ such that, for every $j\in [\mathfrak t,\mathfrak s]$, the following hold:

\begin{itemize}
\item[(a)] no subset $X\subseteq \supp(V\bdot A_{\mathfrak t}\ldots\bdot A_{j-1})$ with $X\in X(G_0)$ has minimal type $\varphi_j$;
\item[(b)] $|\supp(\{V\})|\leq N_S$ and $|\supp(\{A_j\})|\leq N_S$;
\item[(c)] $\vp_{G_0\setminus \mathfrak X_m^\cup(\varphi_j)}(A_j)\leq N_T$; and
\item[(d)] $\Summ{x\in G_0\setminus \mathfrak X_m^\cup(\varphi_j)}\vp_x(A_j)x\in -\C(Z_{\varphi_j})$.
\end{itemize}
For instance, (a)--(d) hold trivially   when $\mathfrak t=\mathfrak s+1$ taking $V=U$ and $N_S=N_T=0$, completely the base of the inductive argument.

\smallskip

Suppose the inductive process is finished, that is, (a)--(d) hold with  $\mathfrak t=1$. Set $A_0=V$.  Then the theorem holds apart from showing $|V|=|A_0|\leq N_T$. From (a) we conclude that there is no nonempty subset $X\subseteq \supp(V)$ with $X\in  X(G_0)$ (as any such $X$ has some minimal type $\varphi_j$). In view of  Proposition \ref{prop-finitary-FiniteDeletion}, there are only a finite number of atoms $W\in \mathcal A(G_0)$ such that there is no subset $X\subseteq \supp(W)$ with $\emptyset \neq X\in X(G_0)$. In particular, there are only a finite number of such elementary atoms, say  $W_1,\ldots,W_\ell\in \mathcal A^{\mathsf{elm}}(G_0)$, such that there is no nonempty subset $X\subseteq \supp(W_j)$ with $X\in  X(G_0)$, for $j\in [1,\ell]$.
By Theorem \ref{thm-carahtheodory-elm-atom}, any zero-sum rational sequence has a factorization as a product of rational powers of elementary atoms. For the rational zero-sum $V\in \mathcal B_{\mathsf{rat}}(G_0)$, only the elementary atoms $W_1,\ldots,W_\ell$ can occur in such a factorization (in view of our prior work). Thus $V=\prod^\bullet_{j\in [1,\ell]}W_j^{[\alpha_j]}$ for some rational numbers $\alpha_j\geq 0$. As $V\mid U$ with $U\in \mathcal A(G_0)$, we must have $\alpha_j\leq 1$ for all $j\in [1,\ell]$. Hence $|V|=\Sum{j=1}{\ell}\alpha_j|W_j|\leq \Sum{j=1}{\ell}|W_j|$. Observing that $\Sum{j=1}{\ell}|W_j|$ is a fixed constant independent of $U$, we can now replace $N_T$ by $\max\{N_T,\Sum{j=1}{\ell}|W_j|\}$ and thereby obtain the final remaining conclusion $|V|\leq N_T$, which would complete the proof. So we may instead assume $\mathfrak t>1$ and proceed with the construction of  $A_{\mathfrak t-1}$.

\smallskip

Let $\varphi=\varphi_{\mathfrak t-1}$, let $Z_\varphi=Z_1\cup\ldots\cup Z_s$  be the codomain of $\varphi$, let $$Z_\varphi=\{z^{(\varphi)}_1,\ldots,z^{(\varphi)}_n\},$$ with $z^{(\varphi)}_1,\ldots,z^{(\varphi)}_n\in \Z_\varphi$ the  $n\leq d$ distinct element in $Z_\varphi$, and let $$\pi:\R^d\rightarrow \R\la Z_\varphi\ra^\bot$$ be the orthogonal projection. Since $\varphi$ is nontrivial, we have $Z_\varphi\neq \emptyset$, whence $\dim \R\la Z_\varphi\ra^\bot<d$.
 Let $\mathfrak X=\mathfrak X_m^\cup(\varphi)=\bigcup_{X\in \mathfrak X_m(\varphi)}X\subseteq G_0^\diamond$ (the inclusion follows by the remarks from the beginning of Section \ref{sec-series-decomp}) and, for each $i\in [1,n]$, let $\mathfrak X_i\subseteq \mathfrak X$ consist of all $x\in \mathfrak X$ with $\varphi(x)=z_i^{(\varphi)}$. Thus $$\mathfrak X=\bigcup_{i=1}^n\mathfrak X_i\subseteq G_0^\diamond.$$
In view of Proposition \ref{prop-finitary-Modulo-Inheritence}, $\pi(G_0)$ is also finitary.   Let $U\in \mathcal A(G_0)$ be an arbitrary atom. By induction hypothesis, we have bounds $N_S\geq 0$ and $N_T\geq 0$ and a factorization $U=V\bdot A_\mathfrak t\bdot\ldots\bdot A_\mathfrak s$ such that  (a)--(d) hold.
Let $V_\mathfrak X\mid V$ be the rational subsequence consisting of all terms from $\mathfrak X$ and let  $V_{G_0\setminus \mathfrak X}\mid V$ be the rational subsequence consisting of all terms from $G_0\setminus \mathfrak X$. Thus $$V=V_\mathfrak X\bdot V_{G_0\setminus \mathfrak X}\quad\mbox{ with $V_\mathfrak X\in \Fc_{\mathsf{rat}}(\mathfrak X)$}\; \und\;\mbox{ $V_{G_0\setminus \mathfrak X}\in \Fc_{\mathsf{rat}}(G_0\setminus \mathfrak X)$}.$$
Note, since $\supp(V_\mathfrak X)\cap \supp(V_{G_0\setminus \mathfrak X})=\emptyset$, that
$$|\{V\}|=|\{V_\mathfrak X\}|+|\{V_{G_0\setminus \mathfrak X}\}|.$$
Since any $X\in \mathfrak X_m(\varphi)$ has $\R\la X\ra=\R\la Z_\varphi\ra=\ker \pi$, we have  $\pi(x)=0$ for all $x\in \mathfrak X$ by definition of $\mathfrak X$. Thus $V\in \mathcal B_{\mathsf{rat}}(G_0)$ implies  $$\pi(V_{G_0\setminus \mathfrak X})\in \mathcal B_{\mathsf{rat}}(\pi(G_0)).$$
If there is no subset $X\subseteq \supp(V)$ with $X\in X(G_0)$ having minimal type $\varphi=\varphi_{\mathfrak t-1}$, then we may take $A_{\mathfrak t-1}$ to be the trivial sequence and find that (a)--(d) hold for the factorization $U=V\bdot A_{\mathfrak t-1}\bdot\ldots\bdot A_\mathfrak s$, completing the induction and proof. Therefore we may assume otherwise that there is some $Y\in \mathfrak X_m(\varphi)$ with $Y\subseteq \supp(V)$, say with corresponding maximal series decomposition $Y=Y_1\cup \ldots\cup Y_s$. In particular, $Y\subseteq \supp(V_\mathfrak X)$ with $V_\mathfrak X$ nontrivial (as $\varphi=\varphi_{\mathfrak t-1}$ is nontrivial).

If there were some nonempty subset $Y_{s+1}\subseteq \supp(V)$ such that $|\pi(Y_{s+1})|=|Y_{s+1}|$ and $\pi(Y_{s+1})\in X(\pi(G_0))$ is irreducible, then we could greedily extend the realization associated to $Y=Y_1\cup\ldots\cup Y_s$ by that associated to $\pi(Y_{s+1})$ to conclude $Y'=Y_1\cup \ldots\cup Y_s\cup Y_{s+1}\in X(G_0)$ with $Y'\subseteq \supp(V)$. By Proposition \ref{prop-finitary-MaxDecompChar}, $Y'=Y_1\cup \ldots\cup Y_s\cup Y_{s+1}$ is a maximal series decomposition.
Letting $\varphi'$ be a minimal type associated to $Y'=Y_1\cup \ldots\cup Y_s\cup Y_{s+1}$, it follows that $\mathsf{dep}(\varphi')=s+1>s=\mathsf{dep}(\varphi)=\mathsf{dep}(\varphi_{\mathfrak t-1})$, which
 in view of the choice of indexing for the $\varphi_i$ forces $\varphi'=\varphi_j$ for some $j\geq \mathfrak t$,  contrary to the conclusion of  (a). Therefore we can instead assume no nonempty subset  $X\subseteq \supp(\pi(V))$ has $X\in  X(\pi(G_0))$, as any such set must contain an irreducible subset.  In consequence, since $\supp(\pi(V))=\supp(\pi(V_{G_0\setminus \mathfrak X}))\cup \{0\}$, we also find that no nonempty subset  $X\subseteq \supp(\pi(V_{G_0\setminus \mathfrak X}))$ has $X\in  X(\pi(G_0))$.

\smallskip

\textbf{Definition and Properties of the $T'_k$.}
Since $\pi(G_0)$ is finitary, it follows from Proposition \ref{prop-finitary-FiniteDeletion} that there are only a finite number of atoms $W\in \mathcal A(\pi(G_0))$ such that there is no nonempty subset $X\subseteq \supp(W)$ with $X\in X(\pi(G_0))$.
In particular, there are only a finite number of such elementary atoms, say  $W_0,W_1\ldots,W_{\ell_\varphi}\in \mathcal A^{\mathsf{elm}}(\pi(G_0))$, where $W_0$ is the zero-sum consisting of a single term equal to $0$, such that there is no nonempty subset $X\subseteq \supp(W_j)$ with $X\in X(\pi(G_0))$, for $j\in [0,\ell_\varphi]$. By Theorem \ref{thm-carahtheodory-elm-atom},
any zero-sum rational sequence has a factorization as a product of rational powers of  elementary atoms. For the rational zero-sum $\pi(V_{G_0\setminus \mathfrak X})\in \mathcal B_{\mathsf{rat}}(\pi(G_0))$, only the elementary atoms $W_0,\ldots,W_{\ell_\varphi}$ can occur in such a factorization (in view of our prior comments). Thus $\pi(V_{G_0\setminus \mathfrak X})=\prod^\bullet_{j\in [0,\ell_\varphi]}W_j^{[w_i]}$ for some rational numbers $w_j\geq 0$.
Since $\pi(U)$ is not an atom, we cannot conclude  $w_j\leq 1$. However, letting $m_{j}=\lfloor w_{j}\rfloor$ and $\varepsilon_j=w_j-m_j\in [0,1)\cap \Q$, we have
$$\pi(V_{G_0\setminus \mathfrak X})=\Big({\prod}^\bullet_{j\in [0,\ell_\varphi]} W_j^{[m_{j}]}\Big)\bdot\Big({\prod}^\bullet_{j\in [0,\ell_\varphi]} W_j^{[\varepsilon_j]}\Big).$$
Let $\ell'_U=\Sum{i=0}{\ell_\varphi}m_{i}+\Sum{i=0}{\ell_\varphi}\lceil \varepsilon_i \rceil$. Then there must be a factorization $$V_{G_0\setminus \mathfrak X}=T'_1\bdot\ldots\bdot T'_{\ell'_U},\quad\mbox{ where each $T'_k\in \Fc_{\mathsf{rat}}(G_0)$},$$  such that
\begin{itemize}
\item for each $k\in [1,\Sum{i=0}{\ell_\varphi}m_i]$, we have $\pi(T'_k)=W_j$ for some $j\in [0,\ell_\varphi]$, and
\item for each $k\in [\Sum{i=0}{\ell_\varphi}m_i+1,\ell'_U]$, say with $k=\Sum{i=0}{\ell_\varphi}m_i+\alpha$ and  $\alpha\in [1,\Sum{i=0}{\ell_\varphi}\lceil \varepsilon_j\rceil]\subseteq [1,\ell_\varphi+1]$,  we have $\pi(T'_k)=W_j^{[\varepsilon_j]}$ for the $\alpha$-th largest  $j\in [0,\ell_\varphi]$ with $\varepsilon_j>0$.
\end{itemize}
Indeed, the $T'_k$ can be sequentially constructed for $k=1,2,\ldots,\ell'_U$ by always first attempting to only include terms in $T'_{k}$ from $V_{G_0\setminus\mathfrak X}\bdot(T'_{1}\bdot\ldots\bdot T'_{k-1})^{[-1]}$ with integer multiplicities until this is no longer possible, after which we attempt to only include terms in $T'_k$ from
 $V_{G_0\setminus\mathfrak X}\bdot(T'_{1}\bdot\ldots\bdot T'_{k-1})^{[-1]}$ with their full remaining multiplicity in $V_{G_0\setminus\mathfrak X}\bdot(T'_{1}\bdot\ldots\bdot T'_{k-1})^{[-1]}$ until this is no longer possible, and then finally (potentially) including one last term from $V_{G_0\setminus\mathfrak X}\bdot(T'_{1}\bdot\ldots\bdot T'_{k-1})^{[-1]}$ with  multiplicity less than one, which we need only do at most once for each element  of $\supp(\pi(T'_k))$. Additionally, when choosing a term $y$ to include in $T'_k$ according to the previous guidelines, always choose one with $y\in G_0\setminus G_0^\diamond$ whenever possible. Assume the $T'_k$ have been constructed according to such restrictions.

 Let us analyse how the fractional subsequence $\{V_{G_0\setminus \mathfrak X}\}$ compares to $\{V_{G_0\setminus \mathfrak X}\bdot (\prod_{k\in I}^\bullet T'_k)^{[-1]}\}$, for a subset $I\subseteq [1,\ell'_\varphi]$. For each $T'_k$, we have $\pi(T'_k)=W_j$ or $W_j^{[\varepsilon_j]}$ for some $j\in [0,\ell_\varphi]$. The process of constructing $T'_k$ first takes available terms with integer multiplicities in $V_{G_0\setminus \mathfrak X}\bdot (T'_1\bdot\ldots\bdot T'_{k-1})^{[-1]}$ to construct $\lfloor T'_k\rfloor$. Removing these terms cannot create a new fractional term. The process then includes terms with their full available fractional multiplicity. Including these terms actually removes a term with fractional multiplicity. Finally, only the last term $y$ included into $T'_k$, for each $w\in \supp(W_j)$ (so $\pi(y)=w$), can actually create a new fractional term that was not already present in $V_{G_0\setminus \mathfrak X}\bdot (T'_1\bdot\ldots\bdot T'_{k-1})^{[-1]}$, and this can only occur if $y$ has multiplicity at least one in $V_{G_0\setminus \mathfrak X}\bdot (T'_1\bdot\ldots\bdot T'_{k-1}\bdot \lfloor T'_k\rfloor)^{[-1]}$. Let $L_k\subseteq G$ consist of all $y$ included as a last term into $T'_k$, potentially one for each $w\in \supp(W_j)$, such that the term $y$ is included with non-integer multiplicity in $T'_k$.
 Note $|\supp(W_j)|\leq d$, since $W_j\in \mathcal A^{\mathsf{elm}}(\pi(G_0))$ is an elementary atom with $\dim (\im \pi)=\dim \R\la Z_\varphi\ra^\bot<d$, so $|L_k|\leq |\supp(W_j)|\leq d$.
 Note that $L_{\ell'_U}\subseteq \supp(\{V_{G_0\setminus \mathfrak X}\})\cup L_1\cup \ldots \cup L_{\ell'_U-1}$, since every term is included into $T'_{\ell'_U}$ with its full remaining multiplicity.
 Moreover, if $j=0$, then we obtain the improved estimate $|\supp(W_0)|=|W_0|=1$ as an upper bound rather than $d$.


 For $k\in [1,\Sum{i=0}{\ell_\varphi}m_i]$, we have $\pi(T'_k)=W_j$, with the multiplicity of each term in $W_j$ being an integer (since $W_j$ is an elementary atom). If $T'_k\in \Fc(G_0)$, then $T'_k\mid \lfloor V_{G_0\setminus \mathfrak X}\rfloor$, and removing $T'_k$ creates no new fractional term. If $T'_k\notin \Fc(G_0)$, then there is some  $w\in \supp(W_j)$ such that $\vp_w(W_j)>\vp_w(\pi(\lfloor T'_k\rfloor))$. However, since both $\vp_w(W_j)$ and $\vp_w(\pi(\lfloor T'_k\rfloor))$ are integers, we have $\vp_w(W_j)-\vp_w(\pi(\lfloor T'_k\rfloor))\geq 1$. Thus all remaining  terms $y$ to be included into $T'_k$, each with  $\pi(y)=w$ for some  $w\in \supp(W_j)$ satisfying $\vp_w(W_j)>\vp_w(\pi(\lfloor T'_k\rfloor))$, must come from $\{V_{G_0\setminus \mathfrak X}\}$ and have their full remaining multiplicity in $V_{G_0\setminus \mathfrak X}\bdot (T'_1\bdot\ldots\bdot T'_{k-1}\bdot \lfloor T'_k\rfloor)^{[-1]}$ strictly less than one (else $y$ could be included with integer multiplicity in view of $\vp_w(W_j)-\vp_w(\pi(\lfloor T'_k\rfloor))\geq 1$).
 In particular, including the last term   into $T'_k$, for each $w\in \supp(W_j)$, does not create a new fractional term, meaning $L_k\subseteq \supp(\{V_{G_0\setminus \mathfrak X}\})$ for $k\in [1,\Sum{i=1}{\ell_\varphi}m_i]$.  Moreover, for each $w\in \supp(W_j)$ satisfying $\vp_w(W_j)>\vp_w(\pi(\lfloor T'_k\rfloor))$, of which there is at least one in view of $T'_k\notin \Fc(G_0)$, the  next term $y$ included into $T'_k$ with $\pi(y)=w$ is from $\{V_{G_0\setminus \mathfrak X}\}$ and included with its full remaining multiplicity as $\vp_y(V_{G_0\setminus \mathfrak X}\bdot (T'_1\bdot\ldots\bdot T'_{k-1}\bdot \lfloor T'_k\rfloor)^{[-1]})<1\leq \vp_w(W_j)-\vp_w(\pi(\lfloor T'_k\rfloor))$. Thus we lose one fractional term from $\{V_{G_0\setminus \mathfrak X}\}$ for each $k\in [1,\Sum{i=1}{\ell_\varphi}m_i]$ with $T'_k\notin \Fc(G_0)$.

  The results of the previous two paragraphs can now be combined to derive the following summary statements.
 For each $T'_{k}$ with $k\in [1,\Sum{i=0}{\ell_\varphi}m_i]$, we either have $T'_{k}\in \Fc(G_0)$  or  else \be\label{T'suppDescends}|\supp(
  \Big\{ V_{G_0\setminus \mathfrak X}\bdot (T'_1\bdot\ldots\bdot T'_k)^{[-1]}\Big\})|<|\supp(
 \Big\{V_{G_0\setminus \mathfrak X}\bdot(T'_{1}\bdot\ldots\bdot T'_{k-1})^{[-1]}\Big\})|.\ee
 Consequently, (b) ensures there are at most $N_S+\Sum{i=0}{\ell_\varphi}\lceil \varepsilon_i\rceil\leq N_S+\ell_\varphi+1$ indices $k\in [1,\ell'_U]$ with $T'_k\notin \Fc(G_0)$. 
 Moreover, for any $I\subseteq [1,\ell'_U]$, we have
 \be\label{Lsupp}\supp(\Big\{V_{G_0\setminus \mathfrak X}\bdot ({\prod}^\bullet_{k\in I}T'_k)^{[-1]}\Big\})\subseteq \supp(\{V_{G_0\setminus \mathfrak X}\})\cup \bigcup_{k=1}^{\max I} L_k.\ee
 %
We have $|L_k\setminus \supp(\{V_{G_0\setminus \mathfrak X}\})|=0$ for $k\in [1,\Sum{i=0}{\ell_\varphi}m_i]$; we have $|L_k|\leq d$ for each of the $\Sum{i=0}{\ell_\varphi}\lceil \varepsilon _i\rceil\leq \ell_\varphi+1$ values $k\in [\Sum{i=0}{\ell_\varphi}m_i+1,\ell'_U]$; and we have $L_{\ell'_U}\subseteq \supp(\{V_{G_0\setminus \mathfrak X}\})\cup L_1\cup \ldots \cup L_{\ell'_U-1}$. Moreover, if equality holds in the estimate $\Sum{i=0}{\ell_\varphi}\lceil \varepsilon _i\rceil\leq \ell_\varphi+1$, then  $\pi(T'_k)=W_0^{[\varepsilon_0]}$ for $k=\Sum{i=0}{\ell_\varphi}m_i+1$ with corresponding improved estimate $|L_k|\leq 1$.
Thus \eqref{Lsupp} applied to ${\prod}^\bullet_{k\in I}T'_k=V_{G_0\setminus \mathfrak X}\bdot\Big( {\prod}^\bullet_{k\in [1,\ell'_U]\setminus I}T'_k\Big)^{[-1]}$, combined with (b), implies (since $d\geq 1$) that \be
 \label{support-myst} |\supp(\Big\{ {\prod}^\bullet_{k\in I}T'_k\Big\})|\leq N_S+\max\{0,(\ell_\varphi-1) d+1\}\leq N_S+\ell_\varphi d,
 \quad\mbox{ for any $I\subseteq [1,\ell'_U]$}.
 \ee
 Each $T'_k$ with $k\in [1,\ell'_U]$ has $\pi(T'_k)=W_j$ or $W_j^{[\varepsilon_j]}$ for some $j\in [0,\ell_\varphi]$, with $\varepsilon _j\leq 1$. Thus $|\lfloor T'_k\rfloor |\leq |T'_k|\leq |W_j|$. Let $$M_W=\max_{j\in [0,\ell_\varphi]}|W_j|.$$ Then, for every $k\in [1,\ell'_U]$,   \eqref{support-myst} applied with $I=\{k\}$ yields
 \be\label{conq1}
 |\supp(\{T'_k\})|\leq N_S+\max\{0,(\ell_\varphi-1) d+1\}\quad\und\quad
 |\lfloor T'_k\rfloor| \leq |T'_k|\leq M_W.\ee

\smallskip

\textbf{Definition and Properties of the $T_k$.}
We next proceed to define a new factorization  $$V_{G_0\setminus \mathfrak X}=T_1\bdot\ldots\bdot T_{\ell_U},\quad\mbox{ where each $T_k\in \Fc_{\mathsf{rat}}(G_0)$},$$ by modifying the sequences $T'_k$ with $k\in [1,\ell'_U]$ as follows. Let $$\Delta=\Z\la Z_\varphi\ra\quad\und\quad \Lambda'=\Lambda\cap \R\la Z_\varphi\ra.$$ Then $\Delta$ and $\Lambda'$ are both full rank lattices in $\R\la Z_\varphi\ra$ with $\Delta\leq \Lambda'$. It follows that $\Lambda'/\Delta$ is a finite abelian group with Davenport Constant $D_\Delta:=\mathsf D(\Lambda'/\Delta)\leq |\Lambda'/\Delta|$ (see Section \ref{sec-intro-factorization}).

Consider an arbitrary sequence $T'_k$ with $k\in [1,\ell'_U]$. Then  $\pi(T'_k)=W_j$ or $W_j^{[\varepsilon_j]}$ for some $j\in [0,\ell_\varphi]$. Let us call  $T'_k$ \emph{pure} if there is some $w\in \supp(W_j)$ such that $\pi^{-1}(w)\cap \supp(T'_k)\subseteq G_0\setminus G_0^\diamond$.

By definition of the $T'_k$, we have $\sigma(T'_k)\in \ker \pi=\R\la Z_\varphi\ra$ for every $k\in [1,\ell'_U]$. Thus, if $T'_k\in \Fc(G_0)$, then we have $\sigma(T'_k)\in \R\la Z_\varphi\ra\cap \Lambda=\Lambda'$. If $I\subseteq [1,\ell'_U]$ is any subset of indices $k\in [1,\ell'_U]$ with $T'_k\in \Fc(G_0)$ and $|I|\geq D_\Delta$, then the definition of $D_\Delta$ ensures that there is some minimal nonempty subset $I'\subseteq I$ such that  $|I'|\leq D_\Delta$ and  $\sigma(\prod_{k\in I'}^\bullet T'_k)\in \Delta$.

The rational sequences $T_k$ and sets $I_k\subseteq [1,\ell'_U]$ such that $T_1\bdot\ldots\bdot T_k=\prod^\bullet_{i\in I_k}T'_i$ will be constructed sequentially for $k=1,2,\ldots,\ell_U$ (though we will sometimes need to define multiple $T_k$ at the same time). Set $I_0=\emptyset$ and take $T_0$ to be the trivial sequence. Assume we have already constructed $T_1\bdot\ldots\bdot T_k=\prod^\bullet_{i\in I_k}T'_i$. Then we construct
$T_{k+1}$ (or possibly $T_{k+1}$ and $T_{k+2}$ simultaneously, in which case  only $I_{k+2}$ is defined and not $I_{k+1}$) according to the following three possibilities.

(i) Suppose there is a nonempty subset $I\subseteq [1,\ell'_U]\setminus I_k$ of indices $i\in [1,\ell'_U]\setminus I_k$ such that $T'_i\in \Fc(G_0)$ is pure for every $i\in I$ and $\sigma(\prod_{i\in I}^\bullet T'_i)\in \Delta$. Then choose a minimal such subset $I$ and set $T_{k+1}=\prod_{i\in I}^\bullet T'_i$ and $I_{k+1}=I_k\cup I$. Note $|I|\leq D_\Delta$ in such case, as remarked  previously. Since each $T'_i\in \Fc(G_0)$, for $i\in I$, it follows that $\pi(T'_i)\in \mathcal B(\pi(G_0))$ is a zero-sum \emph{sequence} (not just rational sequence) with $\supp(\pi(T'_i))=\supp(W)$ for some elementary atom $W\in \mathcal A^{\mathsf{elm}}(\pi(G_0))$, which forces  $\pi(T'_i)$ to be an \emph{integer} power of $W$ (by definition of an elementary atom), and thus $\pi(T'_i)=W$ by construction of the $T'_i$.


(ii) Suppose the hypotheses of (i)  fail but there is some $i\in [1,\ell'_U]\setminus I_k$ with $\supp(T'_i)\nsubseteq G_0^\diamond$. We have $\pi(T'_i)=W^{[\varepsilon]}$ for some elementary atom $W\in\mathcal A^{\mathsf{elm}}(\pi(G_0))$ and rational number $\varepsilon\in (0,1]$. For each $w\in \supp(W)$, write $\alpha_w:=\vp_w(W^{[\varepsilon]})=\alpha^{\overline \diamond}_w+\alpha^\diamond_w$, where $$\alpha^{\overline \diamond}_w=\underset{x\in \supp(T'_i)\setminus G_0^\diamond}{\Summ{\pi(x)=w}}\vp_x(T'_i)\quad\und\quad\alpha^\diamond_w=\underset{x\in \supp(T'_i)\cap G_0^\diamond}{\Summ{\pi(x)=w}}\vp_x(T'_i).$$
Let $$\varepsilon'=\varepsilon\cdot \max\{\alpha_w^{\overline \diamond}/\alpha_w:\:w\in \supp(W)\}\in (0,\epsilon]\cap \Q\subseteq (0,1]\cap \Q.$$ Note that $\varepsilon'>0$ in view of  $\supp(T'_i)\nsubseteq G_0^\diamond$ and that $\varepsilon'=\varepsilon$ precisely when $T'_i$ is pure.  Now $W^{[\varepsilon']}$ has $W^{[\varepsilon']}\mid W^{[\varepsilon]}$ and $\vp_w(W^{[\varepsilon']})\geq \alpha_w^{\overline \diamond}$ for every $w\in \supp(W)$, with equality holding for each $w$ obtaining the maximum in the definition of $\varepsilon'$. Define a new rational subsequence $T_{k+1}\mid T'_i$ with $\pi(T_{k+1})=W^{[\varepsilon']}$ as follows. Include in $T_{k+1}$ all terms $x\in \supp(T'_i)\setminus G_0^\diamond$ with their full multiplicity from $T'_i$, which is possible in view of $\vp_w(W^{[\varepsilon']})\geq \alpha_w^{\overline \diamond}$. Continue to include terms in $T_{k+1}$ from $x\in \supp(T'_i)\cap G_0^\diamond$ with their full multiplicity from $T'_i$, so long as this is possible, and finally (potentially) add one last term to $T_{k+1}$, for each $w\in \supp(W)$, with only part of its multiplicity from $T'_i$.
Note $T_{k+1}=T'_i$ when $T'_i$ pure. In such case, we set $I_{k+1}=I_k\cup \{i\}$. If $T'_i$ is not pure, then  $T_{k+2}=T'_i\bdot T_{k+1}^{[-1]}$ is a nonempty sequence, and in such case, we set $I_{k+2}=I_k\cup \{i\}$, defining  $T_{k+1}$ and $T_{k+2}$ simultaneously  under these circumstances. For any $w\in \supp(W)$ attaining the maximum in the definition of $\varepsilon'$, we have \be\label{ii-purity}\pi^{-1}(w)\cap \supp(T_{k+1})\subseteq G_0\setminus G_0^\diamond.\ee
Also, $\supp(T_{k+2})\subseteq G_0^\diamond$ since all terms from $\supp(T'_i)\setminus G_0^\diamond$ were included in $T_{k+1}$ with their full multiplicity from $T'_i$.

(iii) Suppose the hypotheses of (i) and (ii) both fail, i.e., $\supp(T'_i)\subseteq G_0^\diamond$ for every $i\in [1,\ell'_U]\setminus I_k$. Then take any remaining $i\in [1,\ell'_U]\setminus I_k$ and set $T_{k+1}=T'_i$ and $I_{k+1}=I_k\cup \{i\}$. If no such $i$ exists, i.e., if $I_k=[1,\ell'_U]$, then set $\ell_U=k$.

We partition the interval $$[1,\ell_U]=I_\Z\cup I_\Q\cup I_\diamond,$$ where $I_\Z\subseteq [1,\ell_U]$ consists of all $k+1\in [1,\ell_U]$ for which $T_{k+1}$ was constructed under condition (i), where $I_\Q$ consists of all $k+1\in [1,\ell_U]$ for which  $T_{k+1}$ was constructed under  condition (ii), and where $I_\diamond$  consists of all $k+1\in [1,\ell_U]$ for which  $T_{k+1}$ was constructed under  condition   (iii) as well as all $k+2\in [1,\ell_U]$ for which  $T_{k+2}$ was constructed under condition (ii).
By construction, $I_\Z$ consists of the first $|I_\Z|$ elements from $[1,\ell_U]$.
We have the following  observations.

If $k\in I_\Z$, then there is a nonempty subset $I\subseteq [1,\ell'_U]$ with $|I|\leq D_\Delta$ and $T_k={\prod}_{i\in I}^\bullet T'_i$ such that, for every $i\in I$, we have  \begin{align}\label{Iz-purity}  T'_i\in \Fc(G_0),\quad  \pi(T'_i)\in\mathcal A^{{\mathsf{elm}}}(\pi(G_0)),\quad\und\\ \label{Iz-purity-ii} \pi^{-1}(w)\cap \supp(T'_i)\subseteq G_0\setminus G_0^\diamond \quad\mbox{ for some $w\in \supp(\pi(T'_i))$}.\end{align}
 Moreover,
\be\label{Iz-integral}|T_k|\leq D_\Delta M_W,\quad T_k\in \Fc(G_0)\quad\und \quad \sigma(T_k)\in \Delta,\quad\mbox{ for every $k\in I_\Z$},
\ee where the first inequality follows in view of \eqref{conq1} and $|I|\leq D_\Delta$.

If $k\in I_\Q\cup I_\diamond$, then $T_k\mid T'_i$ for some $i\in [1,\ell'_U]$, and $\pi(T_k^{[\varepsilon_k]})\in\mathcal A^{{\mathsf{elm}}}(\pi(G_0))$ for some rational $\epsilon_k\geq 1$.
Thus \eqref{conq1} implies
\be\label{TkQBound}|T_k|\leq |T'_i|\leq M_W\quad\mbox{ for all $k\in I_\Q\cup I_\diamond$}.\ee  If $k\in I_\Q$, then \eqref{ii-purity} implies that
\begin{align}\label{Iq-purity}    \pi^{-1}(w)\cap \supp(T_k)\subseteq G_0\setminus G_0^\diamond \quad\mbox{ for some $w\in \supp(\pi(T_k))$}.\end{align}
If $k\in I_\diamond$, then (per the hypotheses of Condition (iii) and the remarks regarding $T_{k+2}$ in Condition (ii)) \be\label{T-diamondpart}\supp(T_k)\subseteq G_0^\diamond.\ee


 \subsection*{Claim A} $|I_\Q|\leq D_\Delta+N_S+2\ell_\varphi$.

\begin{proof}
If $k\in I_\Q$, then there was some $i\in [1,\ell'_U]$ such that $T_k\mid T'_i$ with $\supp(T'_i)\nsubseteq G_0^\diamond$. By construction, the same $i$ does not occur for two distinct $k\in I_\Q$. We can assume there are at most $D_\Delta-1$ such $T'_i$ with $T'_i\in \Fc(G_0)$ and $T'_i$ pure, else the algorithm for producing $T_k$ would have used condition (i) rather than (ii).  As noted after \eqref{T'suppDescends}, there are at most $N_S+\ell_\varphi+1$ indices $i\in [1,\ell'_U]$ with $T'_i\notin \Fc(G_0)$, and this estimate includes all indices from $[\Sum{\alpha=0}{\ell_\varphi}m_\alpha+1,\ell'_U]$.
It remains to find an upper bound for the number of $i\in [1,\Sum{\alpha=1}{\ell_\varphi}m_\alpha]$ such that $\supp(T'_i)\nsubseteq G_0^\diamond$, \  $T'_i\in\Fc(G_0)$  but  $T'_i$ is not pure.
 For any such $T'_i$, we have $\pi(T'_i)=W_j$ for some $j\in [0,\ell_\varphi]$ (since $i\in [1,\Sum{\alpha=1}{\ell_\varphi}m_\alpha]$). Since $T'_i\in \Fc(G_0)$, it follows that $T'_i$ was filled entirely with terms from $\lfloor V\rfloor$ (cf. \eqref{RatSeq6}).
  Since $\supp(T'_{i})\nsubseteq G_0^\diamond$, there is at least one $x\in \supp(T'_{i})$ with $x\in G_0\setminus G_0^\diamond$. Thus we must have $j\neq 0$, else $\supp(T'_{i})=\{x\}\subseteq G_0\setminus G_0^\diamond$ (since $|W_0|=1$), contradicting that $T'_i$ is not pure. Since $T'_i$ is not pure, we have $\pi^{-1}(w)\cap \supp(T_i)\cap G_0^\diamond\neq \emptyset$ for every $w\in \supp(W_j)=\supp(\pi(T'_i))$.

Suppose $T'_{i'}$ is another sequence with $i'\in [1,\Sum{\alpha=0}{\ell_\varphi}m_\alpha]$ such that $$T'_i\in\Fc(G_0)\quad\und\quad \pi(T'_{i'})=W_j=\pi(T'_i).$$
Since $T'_i,\,T'_{i'}\in \Fc(G_0)$, it follows that both $T'_i$ and $T'_{i'}$ were constructed by only including terms with integer multiplicities. Moreover, when there was a choice between a term from $G\setminus G_0^\diamond$ and another
from $G_0^\diamond$ (with available integer multiplicity), those from $G\setminus G^\diamond$ always had preference.
If  $i'<i$, then,
 since there is some  $x\in \supp(T'_{i})$ with $x\in G_0\setminus G_0^\diamond$, it follows that, when choosing terms $y$ to place into $T'_{i'}$ with $\pi(y)=\pi(x)$, we did not exhaust all possible choices from $G_0\setminus G_0^\diamond$.
 Thus, due to our preference to include terms from $G_0\setminus G_0^\diamond$ when possible when constructing the $T'_i$, we conclude that every term $y\in \supp(T'_{i'})$ with $\pi(y)=\pi(x)$ satisfies $y\in G_0\setminus G_0^\diamond$, ensuring that $T'_{i'}$ is pure.
 On the other hand, if $i'>i$, then, since every $w\in  \supp(W_j)$ has $\pi^{-1}(w)\cap \supp(T'_i)\cap G_0^\diamond \neq \emptyset$, we must have exhausted all possible terms from $G_0\setminus G_0^\diamond$ that could every be included in $T'_{i'}$ (again, in view of our preference to choose terms from $G_0\setminus G_0^\diamond$). Hence $\supp(T'_{i'})\subseteq G_0^\diamond$. In conclusion, we see that, for each $j\in [1,\ell_\varphi]$, there is at most one $i\in [1,\Sum{\alpha=0}{\ell_\varphi}m_\alpha]$ such that $\supp(T'_{i})\nsubseteq G_0^\diamond$, \ $T'_i\in \Fc(G_0)$ but $T'_i$ is not pure. Combined with our previous estimates, we now conclude that
$|I_\Q|\leq (D_\Delta-1)+(N_S+\ell_\varphi+1)+\ell_\varphi=D_\Delta+N_S+2\ell_\varphi$,
completing Claim A.
\end{proof}

\subsection*{Claim B} For each  $k\in I_\Z\cup I_\Q$,  we have $-\sigma(T_{k})\in \C(X)$ for any $X\in \mathfrak X_m(\varphi)$.  For each  $k\in I_\Z$,  we have $-\sigma(T_{k})\in \C_\Z(X)$ for any $X\in \mathfrak X_m(\varphi)$. In particular, both these statements hold with $X=Z_\varphi$.

\begin{proof}
Let $\mathcal R_X=(\mathcal X_1\cup \{\mathbf v_1\},\ldots,\mathcal X_{s}\cup \{\mathbf v_{s}\})$ be a realization of $X$ (associated to the minimal type $\varphi$) with all half-spaces from $\bigcup_{j=1}^{s}\mathcal X_j$ having dimension one, in which case $\R^\cup \la \mathcal X\ra=\R\la Z_\varphi\ra=\ker \pi$ and $\Z\la X\ra=\Z\la Z_\varphi\ra=\Delta$. Since $\varphi\in \mathfrak X_m(G_0)$, $\mathcal R_X$ is purely virtual.  Since all half-spaces in $\mathcal R_X$ have dimension one with the elements from $X$ being representatives for the half-spaces from $\mathcal X$, we have $\C^\cup (\mathcal X)=\C(X)$.
If $k\in I_\Q$, then \eqref{Iq-purity} allows us to apply Lemma \ref{lem-finitary-diamond-containment-nonmax} to $T_k$ and thereby conclude $-\sigma(T_k)\in \C^\cup (\mathcal X)=\C(X)$, as desired. If $k\in I_\Z$, then  $T_k=\prod_{i\in I}T'_i$ for some $I\subseteq [1,\ell'_U]$. Then \eqref{Iz-purity} implies that each $T'_i\in \Fc(G_0)$ with $\pi(T'_i)=W$ for some $W\in \mathcal A^{\mathsf{elm}}(\pi(G_0))$. Hence, in view of   \eqref{Iz-purity-ii}, we can apply Lemma \ref{lem-finitary-diamond-containment-nonmax} to each $T'_i$ to conclude  $-\sigma(T'_i)\in \C^\cup(\mathcal X)=\C(X)$,  for $i\in I$. Thus, in view of the convexity of $\C(X)$, we have $-\sigma(T_k)=-\Summ{i\in I}\sigma(T'_i)\in \C(X)$. In view of \eqref{Iz-integral}, we have $\sigma(T_k)\in \Delta$. Thus $-\sigma(T_k)\in \C(X)\cap \Delta=\C_\Z(X)$, with the latter equality in view of the linear independence of $X$ and $\Z\la X\ra=\Delta$, completing Claim B.
\end{proof}

Each $k\in I_\Q$ has $T_k\mid T'_{i(k)}$ for some $i(k)\in [1,\ell'_U]$, whence $\supp(\prod_{k\in I_\Q}^\bullet T_k)\subseteq \supp(\prod_{k\in I_\Q}^\bullet T'_{i(k)})$. By \eqref{support-myst}, we have $|\supp(\Big\{\prod_{k\in I_\Q}^\bullet T'_{i(k)}\Big\})|\leq N_S+\ell_\varphi d$.   Claim A and \eqref{conq1} imply  $|\lfloor \prod_{k\in I_\Q}^\bullet T'_{i(k)}\rfloor|\leq \Summ{k\in I_\Q}|T'_{i(k)}|\leq |I_\Q|M_W\leq (D_\Delta+N_S+2\ell_\varphi)M_W$. As a result of all these calculations, we now find that (via \eqref{RatSeq2})
\begin{align}\nn
|\supp({\prod}_{k\in I_\Q}^\bullet T_k)|\leq |\supp({\prod}_{k\in I_\Q}^\bullet T'_{i(k)})|&\leq |\lfloor {\prod}_{k\in I_\Q}^\bullet T'_{i(k)}\rfloor|+|\supp(\Big\{{\prod}_{k\in I_\Q}^\bullet T'_{i(k)}\Big\})|
\\\label{findfound} &\leq (D_\Delta+N_S+2\ell_\varphi)M_W+N_S+\ell_\varphi d.
\end{align}

Recall the definition of $\mathfrak X=\bigcup_{i=1}^n\mathfrak X_i$. In view of the interchangeability property of $\mathfrak X_m(\varphi)$ (Proposition \ref{prop-finitary-mintype}), the sets $X\in \mathfrak X_m(\varphi)$ with $X\subseteq \supp(V_\mathfrak X)$ are precisely those obtained by choosing one element  $x_i\in \mathfrak X_i\cap \supp(V_\mathfrak X)$ for each $i\in [1,n]$. Moreover, since $Y\subseteq \supp(V_\mathfrak X)$ with $Y\in \mathfrak X_m(\varphi)$,  the sets $\mathfrak X_i\cap \supp(V_\mathfrak X)$ are all nonempty, for $i\in [1,n]$.

Let  $X\in \mathfrak X_m(\varphi)$ be arbitrary, say with $X=\{x_1,\ldots,x_n\}$ the distinct elements $x_i\in X$, and let $g\in \Delta\cap -\C(X)$. Then $X\subseteq G_0^\diamond$ is a linearly independent set  with $\Z\la X\ra=\Delta$ (cf. the remarks at the beginning of Section \ref{sec-series-decomp}). It follows that \be\label{commetsA}\C_{\Z}(X)=\C(X)\cap \Delta.\ee As a result, if $-g=\alpha_1x_1+\ldots+\alpha_n x_n$ is a linear combination with $\alpha_j\in \R_+$ expressing that $-g\in \C(X)$, then it follows by the linear independence of $X$ (and recalling that $g\in \Delta=\Z\la X\ra$) that \be\label{ZintegersA}\alpha_j\in \Z_+\quad\mbox{ for all $j\in [1,n]$}.\ee 

\smallskip

\textbf{The Outer Loop.} We proceed to recursively construct fractional sequences $B_k\in \Fc_{\mathsf {rat}}(G_0)$, for $k=0,1,\ldots$, with $$B_k\mid V_\mathfrak X\bdot (B_0\bdot\ldots\bdot B_{k-1})^{[-1]}.$$  We refer to the process which constructs the $B_k$ as the \emph{Outer Loop}, and initially set $B_0$ to be the trivial sequence.  We view
$$V_k:=V_\mathfrak X\bdot (B_0\bdot\ldots\bdot B_{k})^{[-1]}=V_\mathfrak X\bdot (B_1\bdot\ldots\bdot B_{k})^{[-1]}$$ as the current state resulting after $k$ iterations of the outer loop. We view each $j\in [1,n]$ (indexing $\mathfrak X_j$) as a box which contains some of the damaged elements $$D_{k}:=\supp(\{V_k\})\subseteq G_0^\diamond.$$ A box $j\in [1,n]$ is \emph{depleted} if $\mathfrak X_j\cap \supp(\lfloor V_k\rfloor )=\emptyset$, which means every $x\in\supp(V_k)\cap \mathfrak X_j$ has $\vp_x(V_k)<1$. All elements contained in a depleted box are  damaged.  A box  $j\in [1,n]$ with  $\mathfrak X_j\cap \supp(\lfloor V_k\rfloor )\neq \emptyset$, meaning there is some $x\in \supp(V_k)\cap \mathfrak X_j$ with $\vp_x(V_k)\geq 1$, is called \emph{undepleted}. An undepleted box contains at least one undamaged element.
 Let $J_{k}^d\subseteq [1,n]$ consist of all depleted boxes $j\in [1,n]$ after $k$ iterations, and let $J_{k}^u\subseteq [1,n]$ consist of the remaining undepleted boxes $j\in [1,n]$ after $k$ iterations. Note, once a box becomes depleted, it remains depleted in all further iterations of the outer loop.
If there is some box $j\in [1,n]$ with $\mathfrak X_j\cap \supp(V_{k-1})=\emptyset$, that is, an empty box, then the outer loop process halts after $(k-1)$ iterations. Assuming this is not the case, then the next sequence $B_k$, for $k\geq 1$, must satisfy the following properties.
\begin{itemize}
\item[1.] $|D_k\cap \mathfrak X_j|\leq |D_{k-1}\cap \mathfrak X_j|+1$ for all $j\in J_{k-1}^u$.
\item[2.] $|D_k\cap \mathfrak X_j|\leq |D_{k-1}\cap \mathfrak X_j|$ for all $j\in J^d_{k-1}$.
\item[3.] Either $|J^d_{k}|>|J^d_{k-1}|$ or else  $|D_{k}\cap \mathfrak X_j|<|D_{k-1}\cap \mathfrak X_j|$ for some  $j\in J^d_k=J^d_{k-1}$.
\item[4.] $\sigma(B_k)=-\beta_k g_k$ for some $\beta_k\in \Q_+$ with $0<\beta_k\leq 1$, where
    $$g_k=\left\{
            \begin{array}{ll}
              \sigma(T_k), & \hbox{for $k\leq |I_\Z|$;} \\
              \Sum{i=1}{|I_\Z|}(1-\beta_i)g_i+\Summ{i\in I_\Q}\sigma(T_i), & \hbox{for $k=|I_\Z|+1$;} \\
              g_k=(1-\beta_{k-1})g_{k-1}, & \hbox{for $k\geq |I_\Z|+2$.}
            \end{array}
          \right.$$ Moreover, $\beta_k=1$ is only possible if either $k\geq |I_\Z|+1$ or else $k=1=|I_\Z|$, \ $|I_\Q|=|I_\diamond|=0$ and $V=U$, and in either case, the Outer Loop process  then halts after $k$ iterations.
\end{itemize}
Recall that $I_\Z=[1,|I_\Z|]$ by construction.

Suppose it is possible to construct such sequences $B_k$ as described above (the algorithm for their construction will be given afterwards). We proceed to give some basic properties which must then hold, as well as an estimate for how many iterations the outer loop process can run.

First observe that Claim B combined with  $\C(X)$ being a convex cone ensures that
\begin{align}\label{gk-house}&-g_k\in\C(X) \quad\mbox{for all $X\in \mathfrak X_m(\varphi)$ and $k\geq 1$,}\quad \und \\ \nn &-g_k\in\C_\Z(X) \quad\mbox{for all $X\in \mathfrak X_m(\varphi)$ and $k\leq |I_\Z|$.}\end{align}

If $g_k=0$ with $k\leq |I_\Z|$, then $\sigma(T_k)=0$ with $T_k\mid U$ and  $T_k\in \Fc(G_0)$ (by \eqref{Iz-integral}), in which case $U\in \mathcal A(G_0)$ being an atom forces $U=T_k\mid V_{G_0\setminus \mathfrak X}$, contradicting that $V_\mathfrak X\mid U$ is nontrivial. Therefore we instead conclude that $g_k\neq 0$ for all $k\leq |I_\Z|$.
We next show $g_k\neq 0$ for $k=|I_\Z|+1$. In this case, $g_k=\Sum{i=1}{|I_\Z|}(1-\beta_i)\sigma(T_i)+\Summ{i\in I_\Q}\sigma(T_i)$ with $\sigma(T_i)\in -\C(X)$ for all $i\in I_\Z\cup I_\Q$ (in view of Claim B).
Since $0\notin \C^*(X)$ (as each $X\in \mathfrak X_m(\varphi)$ is linearly independent), a sum of elements from $\C(X)$ can only equal zero if all elements in the sum are themselves zero.
Thus, assuming by contradiction that  $g_k=0$, it follows that  $(1-\beta_i)\sigma(T_i)=0$ for all $i\in I_\Z$ and $\sigma(T_i)=0$ for all $i\in I_\Q$.
As we have already established that $\sigma(T_i)=g_i\neq 0$ for each $i\in I_\Z$, this forces $\beta_i=1$ for all $i\leq |I_\Z|$.
However, per Item 4, we know the Outer Loop halts immediately at step $i$ if every some $\beta_i=1$. Thus, if $I_\Z\neq \emptyset$, then the Outer Loop must have halted at step $i=1$, meaning $k= |I_\Z|+1\geq 2$ does not exist. We are left to conclude that $g_k=0$ for $k=|I_\Z|+1$ is only possible if $I_\Z= \emptyset$ and $\sigma(T_i)=0$ for all $i\in I_\Q$. As a result, since $V\in \mathcal B_{\mathsf{rat}}(G_0)$ has sum zero, it follows that $\sigma(V_\mathfrak X\bdot \prod^\bullet_{i\in I_\diamond}T_i)=\sigma(V\bdot (\prod^\bullet_{i\in I_\Q\cup I_\Z}T_i)^{[-1]})=0$. However, $\supp (V_\mathfrak X\bdot \prod^\bullet_{i\in I_\diamond}T_i)\subseteq G_0^\diamond$ by \eqref{T-diamondpart} and $\mathfrak X\subseteq G_0^\diamond$. Thus our hypothesis $0\notin \C^*(G_0^\diamond)$ (note this is the first time we have used the hypothesis $0\notin \C^*(G_0^\diamond)$ rather than the weaker consequence that $G_0$ is finitary) forces $V_\mathfrak X\bdot \prod^\bullet_{i\in I_\diamond}T_i$ to be the trivial sequence, contradicting that $V_\mathfrak X$ is nontrivial.
So we instead conclude that $g_k\neq 0$ also for $k=|I_\Z|+1$. As a consequence, a simple inductive argument now shows  $g_k\neq 0$ for all $k\geq |I_\Z|+1$ (since Item 4 ensures that $\beta_{k-1}\neq 1$, else the  Outer Loop halts immediately at step $k-1$). Summarizing:
\be\label{gk-notzero} g_k\neq 0 \quad\mbox{ for all $k$}.\ee

Next observe that $|D_0\cap \mathfrak X_j|\leq N_S$ for all $j\in [1,n]$ (in view of (b)), which combined with Items 1 and 2 ensures that \be\label{GoodBoxMaxA}|D_k\cap \mathfrak X_j|\leq k+N_S\quad\mbox{ for all $j\in [1,n]$}.\ee The value $|J_k^d|\in [0,n]$ is nondecreasing with $k$. There are therefore at most $n'\leq n$ `jump' values of $k$ where $|J_k^d|>|J_{k-1}^d|$, say occurring for $1\leq \kappa_1<\ldots<\kappa_{n'}$. Let $\kappa_0=0$ and $\kappa_{n'+1}=\infty$. Note that Item 3 ensures that $|J_1^d|\geq 1$.

We wish to estimate how long the outer loop process can run without halting.
Items 2 and  3 ensure that $|D_k\cap \bigcup_{j\in J_k^d}\mathfrak X_j|$ decreases by at least one after each iteration except when $|J_k^d|>|J_{k-1}^d|$,  in which case it instead increases by at most $(|J_k^d|-|J_{k-1}^d|)(k+N_S)$ in view of \eqref{GoodBoxMaxA}. Thus \be\label{lemonpieA}|D_k\cap \bigcup_{j\in J_k^d}\mathfrak X_j|-|D_{k-1}\cap \bigcup_{j\in J_{k-1}^d}\mathfrak X_j|\leq (|J_k^d|-|J_{k-1}^d|)(k+1+N_S)-1\quad\mbox{ for all $k\geq 1$}.\ee
Consequently, since $|D_k\cap \bigcup_{j\in J_k^d}\mathfrak X_j|\geq 0$ must hold for all $k\geq 0$, it follows that every $$k\leq \kappa_{n'}+|D_{\kappa_{n'}}\cap \bigcup_{j\in J_{\kappa_{n'}}^d}\mathfrak X_j|.$$
Indeed,  apart from (possibly) $k=0$ and the final value of $k$, we must have $|D_k\cap \bigcup_{j\in J_k^d}\mathfrak X_j|\geq |J_k^d|\geq |J_1^d|\geq 1$, as otherwise the mechanism for halting the outer loop process is triggered. In particular, either $\kappa_{n'}=1$ (in which case $n'=1$) or $|D_{\kappa_{n'}-1}\cap \bigcup_{j\in J_{\kappa_{n'}-1}^d}\mathfrak X_j|\geq |J_{\kappa_{n'}-1}^d|\geq |J_1^d|\geq 1$, which is slightly better then the estimate that this quantity be non-negative.
Thus we likewise obtain estimates \be\label{itt1A}\kappa_t\leq \kappa_{t-1}+|D_{\kappa_{t-1}}\cap \bigcup_{j\in J_{\kappa_{t-1}}^d}\mathfrak X_j|\quad\mbox{ for all $t\geq 2$}\ee and  $\kappa_1\leq 1+|D_0\cap \bigcup_{j\in J^d_0}\mathfrak X_j|$ (with equality only possible if $\kappa_1=1$ and $|D_0\cap \bigcup_{j\in J^d_0}\mathfrak X_j|=0$).
For   $k<\kappa_{t+1}$, we have passed $t'\leq t$ jump values, which add at most $$\Sum{\alpha=1}{t'}(|J_{\kappa_\alpha}^d|-|J_{\kappa_{\alpha-1}}^d|)(\kappa_\alpha+1+N_S)-t'\leq \Sum{\alpha=1}{t}(|J_{\kappa_\alpha}^d|-|J_{\kappa_{\alpha-1}}^d|)(\kappa_\alpha+1+N_S)-t'$$ damaged elements into the depleted boxes (in view of \eqref{lemonpieA}), on top of the initial number  $|D_0\cap \bigcup_{j\in J_0^d}\mathfrak X_j|$ of damaged elements. At the same time, for each of the $k-t'$ non-jump steps, the number of damaged elements decreases by at least one (as remarked above \eqref{lemonpieA}). Thus  \be\label{eertA}|D_k\cap \bigcup_{j\in J_k^d}\mathfrak X_j|\leq |D_0\cap \bigcup_{j\in J_0^d}\mathfrak X_j|-k+\Sum{\alpha=1}{t}(|J_{\kappa_\alpha}^d|-|J_{\kappa_{\alpha-1}}^d|)
(\kappa_\alpha+1+N_S)\quad\mbox{ for $k<\kappa_{t+1}$}.\ee
The above estimates gets larger when the values of $\kappa_\alpha$ are each as large as possible, that is, delaying when the jumps $|J_k^d|>|J_{k-1}^d|$ occur can only increase the estimates. Also, increasing the number of times that we have jumps $|J_k^d|>|J_{k-1}^d|$ only increases these estimates, since this  breaks multiple simultaneous jumps into several jumps spaced out, with the later jumps having larger contributing factor. Thus we can obtain an upper bound for how long the outer loop process will run by taking $n'=n$ and delaying each jump as long as possible per the estimates above.
In particular, \eqref{eertA} now yields \be\label{itt2A}|D_{\kappa_t}\cap \bigcup_{j\in J_{\kappa_t}^d}\mathfrak X_j|\leq |D_0\cap \bigcup_{j\in J_0^d}\mathfrak X_j|-\kappa_t+\Sum{\alpha=1}{t}
(\kappa_\alpha+1+N_S).\ee
%
%
Note  $|D_0 \cap \bigcup_{j\in J_0^d}\mathfrak X_j|\leq N_S$ (in view of (b)). Thus $\kappa_1\leq 1+N_S$ and $|D_{\kappa_1}\cap \bigcup_{j\in J_{\kappa_1}^d}\mathfrak X_j|\leq |D_0 \cap \bigcup_{j\in J_0^d}\mathfrak X_j|-\kappa_1+(\kappa_1+1+N_S)\leq 2(1+N_S)$. A simple inductive argument, using \eqref{itt1A}, \eqref{itt2A} and the estimate  $|D_0 \cap \bigcup_{j\in J_0^d}\mathfrak X_j|\leq 1+N_S$, now  gives $\kappa_t\leq (2^{t}-1)(N_S+1)$ and  $|D_{\kappa_t}\cap \bigcup_{j\in J_{\kappa_t}^d}\mathfrak X_j|\leq 2^t(N_S+1)$ for all $t\in [1,n]$. We thus find  every $k\leq \kappa_n+|D_{\kappa_n}\cap \bigcup_{j\in J^d_{\kappa_n}}\mathfrak X_j|\leq (2^n-1)(1+N_S)+2^n(1+N_S)=(2^{n+1}-1)(1+N_S)$, meaning the outer loop process must halt after at most $$(2^{n+1}-1)(1+N_S)$$ steps, which is independent of $U$.

\medskip

Suppose we run the above outer loop process and it stops at step $k\leq (2^{n+1}-1)(1+N_S)$. For each $i\in [1, k]$ with $i\leq |I_\Z|$, let $C_i=B_i\bdot T_i^{[\beta_i]}$. If $\beta_k=1=k=|I_\Z|$ and $|I_\Q|=|I_\diamond|=0$, then $V_{G_0\setminus \mathfrak X}=T_1\mid C_1$, ensuring $\supp(V\bdot C_1^{[-1]})\subseteq\mathfrak X\subseteq G_0^\diamond$.  If $k\geq |I_\Z|+1$, then we have
\ber\nn V\bdot(C_1\bdot\ldots\bdot C_{|I_\Z|})^{[-1]}&=&
\Big(V_\mathfrak X\bdot (B_1\bdot\ldots\bdot B_{|I_\Z|})^{[-1]}\Big)\bdot\Big(V_{G_0\setminus \mathfrak X}\bdot T_1^{[-\beta_1]}\bdot \ldots\bdot T_{|I_\Z|}^{[-\beta_{|I_\Z|}]}\Big)\\\nn
&=&\Big(V_\mathfrak X\bdot (B_1\bdot\ldots\bdot B_{|I_\Z|})^{[-1]}\Big)\bdot\Big({\prod}^\bullet_{i\in I_\Z}T_i^{[1-\beta_i]}\bdot  {\prod}^\bullet_{i\in I_\Q}T_i\bdot{\prod}^\bullet_{i\in I_\diamond}T_i\Big).\eer
Let $S_{|I_\Z|+1}= {\prod}_{i\in I_\Z}^\bullet T_i^{[1-\beta_i]}\bdot {\prod}^\bullet_{i\in I_\Q}T_i$ and $C_{|I_\Z|+1}=B_{|I_\Z|+1}\bdot S_{|I_\Z|+1}^{[\beta_{|I_\Z|+1}]}$. For $i\in [|I_\Z|+2,k+1]$, let $S_i=S_{i-1}^{[1-\beta_{i-1}]}$. For  $i\in [|I_\Z|+2,k]$, set $C_i=B_i\bdot S_i^{[\beta_i]}$. For $i\geq |I_\Z|+1$, we have $g_i=\sigma(S_i)$,  while $S_j$ is the subsequence obtained from $\prod^\bullet_{i\in I_\Z}T_i\bdot\prod_{i\in I_\Q}^\bullet T_i$ by removing the terms from the subsequence $\prod_{i\in I_\Z}^\bullet T_i^{[\beta_i]}\bdot \prod^\bullet_{i\in [|I_\Z|+1,j-1]}S_i^{[\beta_i]}$, as can be seen by a short inductive argument on $j=|I_\Z|+1,\ldots, k+1$. Thus
  $C_1\bdot\ldots\bdot C_k\mid V$ and
$\prod_{i\in I_\Z}^\bullet T_i^{[\beta_i]}\bdot \prod^\bullet_{i\in [|I_\Z|+1,k]}S_i^{[\beta_i]}\mid \prod^\bullet_{i\in I_\Z}T_i\bdot\prod_{i\in I_\Q}^\bullet T_i$
with
equality holding precisely when $\beta_k=1$. In particular, combined with the previous observation for what happens when $\beta_k=1=k=|I_\Z|$ with $|I_\Q|=|I_\diamond|=0$, we conclude that
\be \label{diamondgondi}\supp(V\bdot (C_1\bdot\ldots\bdot C_k)^{[-1]})\subseteq G_0^\diamond\quad\mbox{ when $\beta_k=1$}.\ee
 In view of Item 4, we have  $C_i\in \mathcal B_{\mathsf{rat}}(G_0)$ for all $i\in [1,k]$. Define $$A_{\mathfrak t-1}=C_1\bdot\ldots\bdot C_k\quad\und\quad V'=V\bdot (C_1\bdot\ldots\bdot C_k)^{[-1]}=V\bdot A_{\mathfrak t-1}^{[-1]}$$
 Since the $V,\,C_1,\ldots,C_k\in \mathcal B_{\mathsf{rat}}(G_0)$, it follows that $V'\in \mathcal B_{\mathsf{rat}}(G_0)$.
 Let us show that there are bounds $N'_S\geq N_S$ and $N'_T\geq N_T$ such that (a)--(d) hold for $U=V'\bdot A_{\mathfrak t-1}\bdot A_\mathfrak t\bdot\ldots\bdot A_\mathfrak s$. Note, since (a) held originally for $U=V\bdot A_\mathfrak t\bdot\ldots\bdot A_\mathfrak s$, it suffices to show $\supp(V')$ contains no subset $X\in\mathfrak X_m(\varphi)=\mathfrak X_m(\varphi_{\mathfrak t-1})$ in order to show (a) holds for the new factorization.

If $k\leq |I_\Z|$, then  $\prod_{i\in [1,k]}^\bullet T_i^{[\beta_i]}$ is the subsequence of $A_{\mathfrak t-1}=C_1\bdot\ldots\bdot C_k$ consisting of all terms from $G_0\setminus \mathfrak X$, and $\Sum{i=1}{k}\sigma(T_i^{[\beta_i]})=\Sum{i=1}{k}\beta_i\sigma(T_k)\in -\C(Z_\varphi)$ in view of Step B and the convexity of $\C(Z_\varphi)$. If $k\geq |I_\Z|+1$, then $\sigma(S_{|I_\Z|+1})=\sigma({\prod}_{i\in I_\Z}^\bullet T_i^{[1-\beta_i]}\bdot {\prod}^\bullet_{i\in I_\Q}T_i)\in -\C(Z_\varphi)$ in view of Claim B and the convexity of $\C(Z_\varphi)$. A short inductive argument now shows $\sigma(S_i)\in -\C(Z_\varphi)$ for all $i\in [|I_\Z|+1,k]$, as $\C(Z_\varphi)$ is a convex cone. Now
$\prod_{i\in I_\Z}^\bullet T_i^{[\beta_i]}\bdot \prod^\bullet_{i\in [|I_\Z|+1,k]} S_i^{[\beta_i]}$ is  the subsequence of $A_{\mathfrak t-1}=C_1\bdot\ldots\bdot C_k$ consisting of all terms from $G_0\setminus \mathfrak X$. Thus, since $\sigma(S_i)\in -\C(Z_\varphi)$ for all $i\in [|I_\Z|+1,k]$, the convexity of the convex cone $\C(Z_\varphi)$ ensures that $\sigma(\prod_{i\in I_\Z}^\bullet T_i^{[\beta_i]}\bdot \prod^\bullet_{i\in [|I_\Z|+1,k]} S_i^{[\beta_i]})\in -\C(Z_\varphi)$. We conclude that (d) holds for the factorization $U=V'\bdot A_{\mathfrak t-1}\bdot A_\mathfrak t\bdot\ldots\bdot A_\mathfrak s$.

There are two ways the outer loop process can halt. First, we may have $\beta_k=1$.
In this case, \eqref{diamondgondi} implies $\supp(V')\subseteq G_0^\diamond$, while $V'\in \mathcal B_{\mathsf{rat}}(G_0)$ as already remarked above.  In consequence, since $0\notin \C^*(G_0^\diamond)$ by hypothesis, we conclude that $V'$ is the trivial sequence. (Here, we have again used the hypothesis $0\notin \C^*(G_0^\diamond)$ rather than weaker consequence that $G_0$ is finitary.) In such case,  the remaining part of (a)   holds trivially. Moreover, $V'$ being trivial  forces $V_k$ to be trivial, as $V_k=V_\mathfrak X\bdot (B_1\bdot\ldots\bdot B_k)^{[-1]}$ consists of all terms of $V'$ from $\mathfrak X$. In particular, $D_k=\emptyset$ and $J_k^d=[1,n]$, so that
\be\label{particular}\mbox{Item 4 holding with $\beta_k=1$ implies Items 1--3 trivially hold in the Outer Loop}.\ee

The second way the outer loop process can halt is if there is an empty box $j\in [1,n]$ with $\mathfrak X_j\cap \supp V_k=\emptyset$. In view of the definition of $V'$, we have $V'_\mathfrak X=V_\mathfrak X\bdot (B_1\bdot\ldots\bdot B_k)^{[-1]}=V_k$, where $V'_\mathfrak X\mid V'$ is the subsequence of terms from $\mathfrak X$. Thus the  remaining part of (a) holds in this case as well.

If $k\leq |I_\Z|$, then \eqref{Iz-integral} implies \be\nn\vp_{G_0\setminus \mathfrak X}(A_{\mathfrak t-1})=\Sum{i=1}{k}\beta_i|T_i|\leq \Sum{i=1}{k}|T_i|\leq kD_\Delta M_W\leq (2^{n+1}-1)(1+N_S)D_\Delta M_W.\ee
By Claim A and \eqref{TkQBound}, we have   $\Summ{k\in I_\Q}|T_k|\leq  |I_\Q|M_W\leq (D_\Delta+N_S+2\ell_\varphi)M_W$. As a result,
if $k\geq |I_\Z|+1$, then $|I_\Z|\leq k-1\leq (2^{n+1}-1)(1+N_S)-1$, and combined with \eqref{Iz-integral} we now obtain
\begin{multline*}\vp_{G_0\setminus \mathfrak X}(A_{\mathfrak t-1})\leq\Summ{i\in I_\Z}|T_i|+\Summ{i\in I_\Q}|T_i|\leq  ((2^{n+1}-1)(1+N_S)-1)D_\Delta M_W\\+(D_\Delta+N_S+2\ell_\varphi)M_W= ((2^{n+1}-1)(1+N_S)D_\Delta+N_S+2\ell_\varphi)M_W.
\end{multline*}
Thus, setting \be\label{recursive-nT} N'_T=\max\{N_T,\,((2^{n+1}-1)(1+N_S)D_\Delta+N_S+2\ell_\varphi)M_W\},\ee which is independent of $U$, we see that   (c) holds.


If $k\leq |I_\Z|$, then $V'=V\bdot (C_1\bdot\ldots\bdot C_k)^{[-1]}=\Big(V_\mathfrak X\bdot( B_1\bdot \ldots\bdot B_k)^{[-1]}\Big)\bdot \Big (V_{G_0\setminus \mathfrak X}\bdot  (T_1^{[\beta_1]}\bdot\ldots\bdot T_k^{[\beta_k]})^{[-1]}\Big)$ with $V_k=V_\mathfrak X\bdot( B_1\bdot \ldots\bdot B_k)^{[-1]}$ and $T_i\in \Fc(G_0)$ for all $i\leq k$ (by \eqref{Iz-integral}).
In view of \eqref{GoodBoxMaxA}, we have $|D_k\cap \mathfrak X_j|\leq k+N_S$ for all $j\in [1,n]$;
in view of (b) and $\supp(V_\mathfrak X)\cap \supp(V_{G_0\setminus \mathfrak X})=\emptyset$, we have $|\supp(\{V_{G_0\setminus \mathfrak X}\})|=|\supp(\{V\})|-|\supp(\{V_\mathfrak X\})|\leq |\supp(\{V\})|\leq N_S$;
in view of $T_i\in \Fc(G_0)$ and \eqref{Iz-integral}, we have $|\supp(\{T_i^{[\beta_i]}\})|\leq |\supp(T_i)|\leq |T_i|\leq D_\Delta M_W$.
Thus \eqref{RatSeq4} and \eqref{RatSeq5} yield
\ber\nn
|\supp(\{V'\})|&\leq& |\supp(\{V_k\})|+|\supp(\{V_{G_0\setminus \mathfrak X}\})|+\Sum{i=1}{k}|\supp(T_i)|\\\nn &\leq&\Sum{j=1}{n}|D_k\cap \mathfrak X_j|+|\supp(\{V_{G_0\setminus \mathfrak X}\})|+\Sum{i=1}{k}|\supp(T_i)|\\
\nn &\leq & n(k+N_S)+N_S+k D_\Delta M_W\\\nn
&\leq &
(n+D_\Delta M_W)(2^{n+1}-1)(1+N_S)+(n+1)N_S.
\eer
If $k\geq |I_\Z|+1$,  then $V'=V\bdot (C_1\bdot\ldots\bdot C_k)^{[-1]}=\Big(V_\mathfrak X\bdot( B_1\bdot \ldots\bdot B_k)^{[-1]}\Big)\bdot \Big (V_{G_0\setminus \mathfrak X}\bdot  \prod^\bullet_{i\in I_\Z\cup I_\Q}T_i^{[-\gamma'_i]}\Big)$ for some $\gamma'_i\leq 1$ with $V_k=V_\mathfrak X\bdot( B_1\bdot \ldots\bdot B_k)^{[-1]}$. In view of \eqref{GoodBoxMaxA}, we have $|D_k\cap \mathfrak X_j|\leq k+N_S$ for all $j\in [1,n]$; in view of (a) and $\supp(V_\mathfrak X)\cap \supp(V_{G_0\setminus \mathfrak X})=\emptyset$, we have $|\supp(\{V_{G_0\setminus \mathfrak X}\})|=|\supp(\{V\})|-|\supp(\{V_\mathfrak X\})|\leq |\supp(\{V\})|\leq N_S$; in view of $T_i\in \Fc(G_0)$ for $i\in I_\Z$ and \eqref{Iz-integral}, we have $|\supp(T_i)|\leq |T_i|\leq D_\Delta M_W$ for $i\in I_\Z$;
in view of $k\geq |I_\Z|+1$, we have $|I_\Z|< k\leq (2^{n+1}-1)(1+N_s)$. Combining these estimates with \eqref{findfound}, we obtain (via \eqref{RatSeq4} and \eqref{RatSeq5}) \ber\nn
|\supp(\{V'\})|&\leq& |\supp(\{V_k\})|+|\supp(\{V_{G_0\setminus \mathfrak X}\})|+\Summ{i\in I_\Z}|\supp(T_i)|+|\supp({\prod}_{i\in I_\Q}^\bullet T_i)|\\\nn
&\leq&\Sum{j=1}{n}|D_k\cap \mathfrak X_j|+|\supp(\{V_{G_0\setminus \mathfrak X}\})|+\Summ{i\in I_\Z}|\supp(T_i)|+|\supp({\prod}_{i\in I_\Q}^\bullet T_i)|\\
\nn &\leq & n(k+N_S)+N_S+|I_\Z| D_\Delta M_W+|\supp({\prod}_{i\in I_\Q}^\bullet T_i)|\\\nn
&\leq &
n(k+N_S)+N_S+kD_\Delta M_W+(D_\Delta+N_S+2\ell_\varphi)M_W+N_S+\ell_\varphi d
\\\nn
&\leq &(n+D_\Delta M_W)(2^{n+1}-1)(1+N_S)+(n+2)N_S+\ell_\varphi d+(D_\Delta+N_S+2\ell_\varphi)M_W.
\eer
Setting $N''_S=(n+D_\Delta M_W)(2^{n+1}-1)(1+N_S)+(n+2)N_S+\ell_\varphi d+(D_\Delta+N_S+2\ell_\varphi)M_W\geq N_S$, which is independent of $U$, we see that the first bound in (b) holds.

Finally, \eqref{RatSeq5} implies $\supp(\{A_{\mathfrak t-1}\})\subseteq \supp(\{V'\})\cup \supp(\{V'\bdot A_{\mathfrak t-1}\})=\supp(\{V'\})\cup \supp(\{V\})$. Thus, in view of the above work and  (b) for the original factorization $U=V\bdot V_{\mathsf t}\bdot\ldots\bdot V_\mathfrak s$, we have $|\supp(\{A_{\mathfrak t-1}\})|\leq |\supp(\{V'\})|+| \supp(\{V\})|\leq N''_S+N_S$. Thus, letting  \be\label{recur-NS}N'_S:=(n+D_\Delta M_W)(2^{n+1}-1)(1+N_S)+(n+3)N_S+\ell_\varphi d+(D_\Delta+N_S+2\ell_\varphi)M_W\geq N_S,\ee we see that the second bound in (b) also holds, which completes the induction. It remains only to show that it is indeed possible to construct the sequences $B_k$ of the outer loop with  the desired list of properties, and then the proof will be complete.


\smallskip

\textbf{The Inner Loop.}
Assume that  the rational sequences $B_1,\ldots,B_{k-1}\in \Fc_{\mathsf{rat}}(G_0)$ of the Outer Loop have already been constructed, for $k\geq 1$, and that $\mathcal X_j \cap \supp(V_{k-1})\neq \emptyset$  for all $j\in[1,n]$, and $0<\beta_i<1$ for all $i<k$, so that the Outer Loop process has not terminated. Let $$V'=V_{k-1}=V_\mathfrak X\bdot(B_1\bdot\ldots\bdot B_{k-1})^{[-1]}.$$ We then construct the sequence $B_k\mid V'$ by a separate recursive process which we refer to as the Inner Loop.
 Since $\mathcal X_j \cap \supp(V')\neq \emptyset$ for all $j\in[1,n]$,  we can select a fixed element  $z_j\in\mathcal X_j \cap \supp(V')$ from each depleted box $j\in J_{k-1}^d$. Let $X=\{z_j:\;j\in J^d_{k-1}\}$. Observe that \be\label{xdepleteA}X\subseteq D_{k-1}\cap \bigcup_{j\in J_{k-1}^d}\mathfrak X_j\quad\und\quad X\subseteq\supp(\{V'\})\ee by its definition and that of a depleted box. Let $$W_X={\prod}^\bullet_{x\in X} x^{[\vp_{x}(V')]}\mid \{V'\}.$$ Since all the terms in $W_X$ are depleted in $V_{k-1}=V'$, it follows that $\supp(\lfloor V'\rfloor )\cap \supp(W_X)=\emptyset$.

Assume we have already constructed a sequence  $C\in\mathcal F_{\mathsf{rat}}(G_0)$ which  satisfies
$$C\mid \lfloor V'\rfloor\bdot W_X$$ such that
$\sigma(C)=-\beta g_k$ for some $\beta\in \Q_+$ with $0\leq \beta\leq 1$, and such that \be\label{Z'no1A}Z'=\supp(C)\cap \supp(\lfloor V'\rfloor \cdot W_X\bdot {C}^{[-1]})\ee is  a subset $Z'\subseteq \mathfrak X$ satisfying
\begin{itemize}
\item[(a)] $|Z'\cap \mathfrak X_j|\leq 1$ for all $j\in [1,n]$;
\item[(b)] $Z'\cap \bigcup_{j\in J_{k-1}^d}\mathfrak X_j\subseteq X$;
\item[(c)] If $\beta=0$, then $C$ is trivial; and
\item[(d)] If $\beta=1$, then either $k\geq |I_\Z|+1$ or else $k=1=|I_\Z|$, $|I_\Q|=|I_\diamond|=0$ and $U=V$.
 \end{itemize}
 For instance, we could initially start with $C$ taken to be the trivial sequence, in which case $Z'=\emptyset$ and $\beta=0$. In view of \eqref{Z'no1A}, we see that the terms  $z\in \supp(C)\setminus Z'$  are those which are completely removed from  $\lfloor V'\rfloor \bdot W_X$ when we remove the sequence $C$. Thus, \be\label{latergater}\mbox{if $z\in \supp(C)\setminus Z'$ and $z\notin X$, then $\vp_z(V'\bdot C^{[-1]})=\vp_z(\{V'\})<1$},\ee and if $z\in \supp(C)\setminus Z'$ and $z\in X$, then $\vp_z(V'\bdot C^{[-1]})=0<\vp_z(V')<1$ (as the elements from $X$ come from depleted boxes).
As a result, if we were to halt the Inner Loop process and set $B_k=C$ and $V_k=V'\bdot C^{[-1]}$, then  $D_k$ would be obtained from $D_{k-1}$ by removing all elements from $(X\cap \supp(C))\setminus Z'$ and including (possibly) some of the elements from $Z'$, meaning \be\label{DgogoA}D_{k}\subseteq \big(D_{k-1}\setminus (X\cap \supp(C))\big)\cup Z'.\ee

If we have  $X\nsubseteq \supp(V'\bdot C^{[-1]})$, which is equivalent to $X\nsubseteq \supp(\lfloor V'\rfloor \bdot W_X\bdot C^{[-1]})$ in view of the definition of $W_X$ and the fact that all elements of $X$ are from depleted boxes, then  the Inner loop halts and we set $B_{k}=C$ and $V_k=V'\bdot C^{[-1]}$.
 Note, since $X\subseteq \supp(V')$, that this is only possible if $C$ is nontrivial. Thus Item 4 of the Outer Loop holds with $\beta_k=\beta$ by definition of $C$, (c) and  (d). It follows in view of \eqref{DgogoA} and (a)  that Item 1  of the Outer Loop holds.
 It follows in view of \eqref{DgogoA},  (b), and  \eqref{xdepleteA} that Item 2 from the Outer Loop holds.
Since $X\nsubseteq \supp(\lfloor V'\rfloor\bdot W_X\bdot C^{[-1]})$ and $X=\supp(W_X)$, it follows that there must be some term from $X\cap \supp(C)$ completely removed from $\lfloor V'\rfloor \bdot W_X$, whence $X\cap \supp(C)\nsubseteq Z'$.
 If there is some $j\in J_k^d\setminus J_{k-1}^d$, then Item 3 of the Outer Loop holds. Otherwise, we have $J_k^d=J_{k-1}^d$. In this case, any element $x\in (X\cap \supp(C))\setminus Z'$, which exists as we just observed $(X\cap \supp(C))\nsubseteq Z'$, is an element of $D_{k-1}$ (by \eqref{xdepleteA}) not in $D_k$ (by \eqref{DgogoA}), which combined with the already established Item 2 of the Outer Loop yields Item 3 of the Outer Loop. Thus $B_{k}=C$ satisfies all conditions for the Outer Loop, as required. Therefore   instead assume $X\subseteq \supp(V'\bdot C^{[-1]})$, equivalent to \be\label{Xcont}X\subseteq \supp(\lfloor V'\rfloor \bdot W_X\bdot C^{[-1]}).\ee

If $\beta=1$, then (d) ensures either  $k\geq |I_\Z|+1$ or else $k=1=|I_\Z|$, $|I_\Q|=|I_\diamond |=0$ and $U=V$. In this case,  the Inner loop requirement that $\sigma(C)=-\beta g_k=-g_k$ ensures that Item 4 of the Outer Loop holds with $\beta_k=\beta=1$ taking $B_k=C$, and then Items 1--3 do as well by \eqref{particular}, meaning we can set $B_k=C$ for the Outer Loop, as required. Therefore we may instead assume \be\label{betaless1}\beta<1.\ee

If $x\in \mathfrak X_j$ with $j\in J_{k-1}^u$, then $x\notin X=\supp(W_X)$, as all terms from $X$ are from depleted boxes. Thus $C\mid \lfloor V'\rfloor \bdot W_X$ ensures that $\vp_x(C)\leq \vp_x(\lfloor V'\rfloor)$. Hence either $\vp_x(C)=\vp_x(\lfloor V'\rfloor)$ or $x\in \mathfrak X_j\cap \supp(\lfloor V'\rfloor \bdot W_X\bdot C^{[-1]}\rfloor)$. Now suppose there were some  $x\in \mathfrak X_j\cap \supp(\lfloor V'\bdot C^{[-1]}\rfloor)$ with $x\notin \mathfrak X_j\cap \supp(\lfloor V'\rfloor \bdot W_X\bdot C^{[-1]})$, where $j\in J_{k-1}^u$. Then $\vp_x(V')-\vp_x(C)\geq 1$ and $\vp_x(C)=\vp_x(\lfloor V'\rfloor )$, ensuring that $\vp_x(\{V'\})=\vp_x(V')-\vp_x(\lfloor V'\rfloor)=\vp_x(V')-\vp_x(C)\geq 1$, contrary to the definition of $\{V'\}$. So we instead conclude that
\be\label{supp-identity}\mathfrak X_j\cap \supp(\lfloor V'\bdot C^{[-1]}\rfloor)\subseteq \mathfrak X_j\cap \supp(\lfloor V'\rfloor \bdot W_X\bdot C^{[-1]})\quad\mbox{ for every $j\in J_{k-1}^u$}.\ee

If there is some $j\in J_{k-1}^u$ with $\mathfrak X_j\cap \supp(\lfloor V'\bdot C^{[-1]}\rfloor)=\emptyset$, then the Inner loop halts and we set $B_{k}=C$ and $V_k=V'\bdot C^{[-1]}$. Note this ensures that $|J_{k}^d|>|J_{k-1}^d|$ (as $j\in J_{k}^d\setminus  J_{k-1}^d$), so Item 3 of the Outer Loop holds. Since $\mathfrak X_j\cap \supp(\lfloor V'\rfloor)\neq \emptyset$ for every $j\in J_{k-1}^u$ by definition of an undepleted box, we must have $C$ nontrivial. Thus Item 4 of the Outer Loop holds with $\beta_k=\beta$ by definition of $C$, (c) and \eqref{betaless1}. By the same arguments used to establish \eqref{Xcont}, it follows from  (a), (b), \eqref{xdepleteA} and \eqref{DgogoA} that Items 1 and 2 of the Outer loop hold. Thus $B_{k}=C$ satisfies all conditions for the Outer Loop, as required. Therefore   instead assume every $j\in J_{k-1}^u$ has some $z_j\in \mathfrak X_j\cap \supp(\lfloor V'\bdot C^{[-1]}\rfloor)\subseteq \mathfrak X_j\cap \supp(\lfloor V'\rfloor \bdot W_X\bdot C^{[-1]})$, with the inclusion  in view of \eqref{supp-identity}.


In view (b), for every $j\in J_{k-1}^d$, the set $Z'$ either contains $z_j\in X\cap \mathfrak X_j$ or no element from $\mathfrak X_j$ at all. In view of \eqref{Xcont}, every $j\in J_{k-1}^d$ has $z_j\in \supp(\lfloor V'\rfloor \bdot W_X\bdot C^{[-1]})$.
In view of (a), for every $j\in J_{k-1}^u$, the set $Z'$ contains at most one element from $\mathfrak X_j$.
In view of the conclusion of the previous paragraph, for every $j\in J_{k-1}^u$ for which $Z'\cap \mathfrak X_j=\emptyset$, there is some $z_j\in \mathfrak X_j\cap \supp(\lfloor V'\bdot C^{[-1]}\rfloor)\subseteq \mathfrak X_j\cap \supp(\lfloor V'\rfloor \bdot W_X\bdot C^{[-1]})$.
As a result of all these observations, \eqref{Z'no1A} and the Interchangeability Property of $\mathfrak X_m(\varphi)$ (Proposition \ref{prop-finitary-mintype}), we can find some subset (say) $$Z=\{z_1,\ldots,z_n\}$$ with   $Z\in \mathfrak X_m(\varphi)$ and
\be\label{wehoveA}X\cup Z'\subseteq Z\subseteq \supp(\lfloor V'\rfloor \bdot W_X\bdot C^{[-1]})\subseteq \supp(V'\bdot C^{[-1]}),\ee where the latter inclusion follows by recalling that $W_X\mid \{V'\}$. In particular, since $X=\supp(W_X)$, we have
\be\label{ZminusX-inclu}Z\setminus X\subseteq \supp(\lfloor V'\rfloor).\ee


In view of the definition of the $g_i$, it follows that $g_k$ (indeed, every $g_i$) is a positive rational linear combination of the $\sigma(T_i)$ with $i\in I_\Z\cup I_\Q$. Let $\beta<1$ be the rational number $\beta\in \Q_+$ from the definition of $C$.
Since $\beta<1$, it follows that $(1-\beta)g_k$ is a positive scalar multiple of $g_k$, and thus $-(1-\beta)g_k\in \C(Z)$ in view of Claim B, \ $Z\in \mathfrak X_m(\varphi)$, and the convexity of $\C(Z)$.
Since $Z\in \mathfrak X_m(\varphi)$ is linearly independent, let $\alpha_{1},\ldots,\alpha_n\in \R_+$ be the unique real numbers such that $$\alpha_{1}z_1+\ldots+\alpha_{n}z_n=-(1-\beta)g_k.$$
Since $g_k$  is a positive rational linear combination of the $\sigma(T_i)$ with $i\in I_\Z\cup I_\Q$, with each $\sigma(T_i)$ a positive rational linear combination of terms from $G_0\subseteq \Lambda$ (since $T_i\in \Fc_{\mathsf{rat}}(G_0)$), and since  $\beta\in \Q_+$,  it follows that $-m'(1-\beta)g_k\in \Lambda\cap \C(Z)\subseteq \Lambda \cap \R\la \Delta\ra=\Lambda'$, for some integer $m'\geq 1$ (recall that $\Z\la Z\ra=\Z\la Z_\varphi\ra=\Delta$ since $Z,\,Z_\varphi\in \mathfrak X_m(\varphi)$).
Thus, since $\Lambda'/\Delta$ is a finite abelian group, it follows that there is some integer $m\geq 1$ such that $-m(1-\beta)g_k\in \Delta\cap \C(Z)=\C_\Z(Z)$, with the equality in view of \eqref{commetsA}. Consequently, \eqref{ZintegersA} ensures that
 $m\alpha_i\in\Z_+$ for all $i\in [1,n]$, showing that  $\alpha_1,\ldots,\alpha_n\in \Q_+$.

\subsection*{Claim C} Either $\vp_{z_i}(\lfloor V'\rfloor \bdot W_X\bdot C^{[-1]})<\alpha_{i}$ for some $i\in [1,n]$, or $k\geq |I_\Z|+1$, or $k=1=|I_\Z|$, $|I_\Q|=|I_\diamond|=0$ and $U=V$.

\begin{proof}
Suppose $k\leq |I_\Z|$ and  $\vp_{z_i}(\lfloor V'\rfloor \bdot W_X\bdot C^{[-1]})\geq \alpha_{i}$ for all $i\in [1,n]$. We will show that $k=1=|I_\Z|$, $|I_\Q|=|I_\diamond|=0$ and $U=V$.
Since $k\leq |I_\Z|$, we have $$-g_k=-\sigma(T_k)\in \C_\Z(Z)\subseteq \Z\la Z\ra=\Delta$$ from \eqref{gk-house} and $Z\in \mathfrak X_m(\varphi)$. 
Also, $W_X\mid \{V'\}$ ensures that $\lfloor V'\rfloor \bdot W_X\mid V'=V_{k-1}$, which is a subsequence of $V_\mathfrak X$, thus disjoint from $T_k\mid V_{G_0\setminus \mathfrak X}$.
Hence, since $\sigma(C)=-\beta g_k=-\beta\sigma(T_k)$ by hypothesis of the Inner Loop process, we conclude that  $$T:=T_{k}\bdot C\bdot {\prod}^\bullet_{i\in [1,n]}z_i^{[\alpha_i]}\in \mathcal B_{\mathsf{rat}}(G_0)\quad\und\quad  T\mid U.$$ In particular, $T$ is zero-sum.
By \eqref{Iz-integral}, we have $$T_k\in \Fc(G_0).$$
By \eqref{latergater}, we have $\vp_x(V'\bdot C^{[-1]})=\vp_x(\{V'\})$ for all  $x\in \supp(C)\setminus (Z'\cup X)$. Thus
\ber\nn\vp_x(C)&=&\vp_x(V')-\vp_x(V'\bdot C^{[-1]})\\\label{churntA}&=&\vp_x(V')-\vp_x(\{V'\})=\vp_x(\lfloor V'\rfloor )\in \Z_+\quad\mbox{ for all $x\in \supp(C)\setminus( Z'\cup X)$}.\eer
  Since $T$ is zero-sum and $\supp(T_k)\cap Z\subseteq \supp(V_{G_0\setminus \mathfrak X})\cap \mathfrak X=\emptyset$, we have $$\Sum{i=1}{n}\vp_{z_i}(T)z_i=-\Summ{x\in \supp(C)\setminus Z}\vp_x(C)x-\sigma(T_k)=-\Summ{x\in \supp(C)\setminus Z}\vp_x(C)x-g_k.$$
As the elements from any $\mathfrak X_j$ are contained in some set from  $\mathfrak X_m(\varphi)$, with every such set being a lattice basis for $\Delta$, it follows that
$\supp(C)\subseteq \mathfrak X\subseteq \Delta$. Thus, since $g_k\in \Delta$ and $\vp_x(C)\in \Z_+$ for all $x\in \supp(C)\setminus Z\subseteq \supp(C)\setminus (Z'\cup X)$ (in view of \eqref{churntA} and \eqref{wehoveA}), it follows that  $\Sum{i=1}{n}\vp_{z_i}(T)z_i\in \Delta.$ Since $\vp_{z_i}(T)=\vp_{z_i}(C)+\alpha_i\geq 0$, we have $\Sum{i=1}{n}\vp_{z_i}(T)z_i\in \C(Z)$.
But now  $\Sum{i=1}{n}\vp_{z_i}(T)z_i\in \Delta\cap \C(Z)=\C_{\Z}(Z)$  by \eqref{commetsA}, whence $\vp_{z_i}(T)\in \Z_+$ for all $i\in [1,n]$ by \eqref{ZintegersA}. Combined with \eqref{churntA}, $X\cup Z'\subseteq Z$ (in view of \eqref{wehoveA}) and $T_{k}\in \Fc(G_0)$, we conclude that $T\in\Fc(G_0)$ is an ordinary zero-sum sequence. However, since $T\mid U$ is nontrivial (as each $T_{k}$ is nontrivial) with $U\in \mathcal A(G_0)$ an atom, this is only possible if $T=U$, which, in turn, is only possible if $|I_\Z|=1=k$, $|I_\Q|=|I_\diamond|=0$, $V_{G_0\setminus \mathfrak X}=T_1$, and $U=V$. This establishes Claim C.
\end{proof}

Let $$\gamma=\min(\{1\}\cup \{\vp_{z_i}(\lfloor V'\rfloor \bdot W_X\bdot C^{[-1]})/\alpha_i :\; i\in [1,n],\,\alpha_i\neq 0\})\leq 1.$$
We cannot have
 $\alpha_i= 0$ for all $i\in [1,n]$ in view of $\beta<1$ and \eqref{gk-notzero}. Step C ensures that $\gamma=1$ is only possible if $k\geq |I_\Z|+1$ or else $k=1=|I_\Z|$, \ $|I_\Q|=|I_\diamond|=0$ and $U=V$.
 In view of \eqref{wehoveA}, $\vp_{z_i}(\lfloor V'\rfloor \bdot W_X\bdot C^{[-1]})>0$ for all $i\in [1,n]$, implying $\gamma>0$.

By definition, we have $\vp_{z_i}(\lfloor V'\rfloor \bdot W_X\bdot C^{[-1]})\geq \gamma \alpha_i$ for all $i\in [1,n]$, with equality holding for any $i\in [1,n]$ attaining the minimum in the definition of $\gamma$.
Moreover, since $V',\,C\in \Fc_{\mathsf{rat}}(G_0)$, we have  $\alpha_i,\,\vp_{z_i}(\lfloor V'\rfloor \bdot W_X\bdot C^{[-1]})\in \Q_+$, implying  $\gamma\in \Q_+$ is a positive rational number, and thus  $\gamma\alpha_i\in \Q_+$ for every  $i\in[1,n]$.  Define $$C'=C\bdot {\prod}_{i\in [1,n]}^\bullet z_i^{[\gamma \alpha_i]}.$$ Note $C'$ is nontrivial in view of $\gamma>0$ and not all $\alpha_i=0$. By construction, $C'\mid \lfloor V'\rfloor \bdot W_X$  and $C'\in \Fc_{\mathsf{rat}}(G_0)$ (since $C\in \Fc_{\mathsf{rat}}(G_0)$ with $\gamma\alpha_i\in \Q_+$ for all $i\in [1,n]$).
Note $$\sigma(C')=\sigma(C)-\gamma (1-\beta)g_k=-\beta g_k-\gamma (1-\beta)g_k=-(\beta+\gamma(1-\beta))g_k.$$ Since $\beta,\,\gamma\in \Q_+$ with $0<\gamma\leq 1$ and $0\leq \beta<1$, we have $\beta':=\beta+\gamma(1-\beta)\in \Q_+$ with
$0<\beta'\leq 1$. Furthermore, $\beta'=1$ if and only if $\gamma=1$. Since $Z'\subseteq Z=\{z_1,\ldots,z_n\}$ and $C\mid C'$, it follows from \eqref{Z'no1A} that $$Z''=\supp(C')\cap \supp(\lfloor V'\rfloor \bdot W_X\bdot  {C'}^{[-1]})\subseteq Z.$$
Since $Z\subseteq \mathfrak X$ with $|Z\cap \mathfrak X_j|=1$ for all $j\in [1,n]$ (as $Z\in \mathfrak X_m(\varphi)$), it follows that the subset $Z''\subseteq Z$ satisfies $Z''\subseteq \mathfrak X$ with $|Z''\cap \mathfrak X_j|\leq 1$ for all $j\in [1,n]$.
 Note $Z''\cap \bigcup_{j\in J_{k-1}^d}\mathfrak X_j\subseteq Z\cap \bigcup_{j\in J_{k-1}^d}\mathfrak X_j=X$, with the equality in view of $X\subseteq Z\in \mathfrak X_m(\varphi)$ (from \eqref{wehoveA}). Finally, if $\beta'=1$, then $\gamma=1$, whence either $k\geq |I_\Z|+1$ or else $k=1=|I_\Z|$, \ $|I_\Q|=|I_\diamond|=0$ and $U=V$ (as remarked after the definition of $\gamma$).  Thus the rational sequence $C'$ satisfies all the requirements for our recursive construction of the Inner Loop. Hence we may replace $C$ by $C'$ and repeat the Inner Loop process just described once more.

If at any point while iterating the Inner Loop Process we obtain $\gamma=1$, then $\beta'=1$ follows, and the Inner Loop immediately halts and outputs $C'=B_k$ as described earlier (when arguing that we could assume $\beta<1$). Otherwise, there will always be an index  $i\in [1,n]$   attaining the minimum in the definition of $\gamma$, in which case   $z_i\in Z\subseteq  \supp(\lfloor V'\rfloor \bdot W_X\bdot C^{[-1]})$ but $z_i\notin \supp(\lfloor V'\rfloor \bdot W_X\bdot {C'}^{[-1]})$, thus ensuring  $|\supp(\lfloor V'\rfloor \bdot W_X\bdot C^{[-1]})|>|\supp(\lfloor V'\rfloor \bdot W_X\bdot {C'}^{[-1]})|$ has strictly decreased.
As a result, we cannot iterate the Inner Loop process indefinitely, as the non-negative integer $|\supp(\lfloor V'\rfloor \bdot W_X\bdot C^{[-1]})|$ strictly decreases after each iteration, ensuring that the process must halt for one of the three possibilities described earlier, all of which lead to the construction of the next sequence $B_k$ in the outer loop process. As this shows it is always possible to construct the sequence $B_k$ with the needed properties for the Outer Loop process, the proof is now complete.
\end{proof}

\begin{corollary}\label{cor-structural-Lambert}
Let $\Lambda\leq \R^d$ be a full rank lattice, where $d\geq 0$,  and let $G_0\subseteq \Lambda$ be a finitary subset with $\C(G_0)=\R^d$. Then $G_0\setminus G_0^\diamond\subseteq G_0$ is a Lambert subset.
\end{corollary}

\begin{proof}
Suppose $0\notin \C^*(G_0)$ (which implies $G_0$ is finitary by Theorem \ref{thm-keylemmaII}). Let  $U\in \mathcal A(G_0)$ be arbitrary. Per the remarks above Theorem \ref{thm-structural-Lambert}, apply the lattice type version of Theorem \ref{thm-structural-Lambert} to the atom $U$, so $\mathfrak s$ equals the number of nontrivial lattice types rather than minimal types.  Letting $\mathfrak X^\cup (\varphi_j)=\bigcup_{X\in \mathfrak X(\varphi_j)}X$, observe that $\mathfrak X^\cup(\varphi_j)\subseteq G_0^\diamond$, for all $j\in [1,\mathfrak s]$, since $\mathfrak X^\cup(\varphi_j)$ is a union of sets $X\in X(G_0)$ with each $X\subseteq G_0^\diamond$, as discussed at the beginning of Section \ref{sec-series-decomp}. Thus $G_0\setminus G_0^\diamond\subseteq  G_0\setminus \mathfrak X^\cup(\varphi_j)$ for all $j\in [1,\mathfrak s]$, whence Theorem \ref{thm-structural-Lambert}(c) ensures $\vp_{G_0\setminus G_0^\diamond}(U)=\Sum{i=0}{\mathfrak s}\vp_{G_0\setminus G_0^\diamond}(A_i)\leq (\mathfrak s+1)N_T$, showing that $G_0\setminus G_0^\diamond\subseteq G_0$ is a Lambert subset with bound $(\mathfrak s+1)N_T$, where $N_T$ is bound given in Theorem \ref{thm-structural-Lambert}.

If $G_0$ is only finitary, with $0\in \C^*(G_0)$, then the corollary does not follow directly from Theorem \ref{thm-structural-Lambert}. However, most parts of the proof of Theorem \ref{thm-structural-Lambert} only require the hypothesis that $G_0$ is finitary, not the stronger hypothesis that $0\in \C^*(G_0)$. Indeed, this stronger hypothesis is only used twice in the proof: first, in the paragraph after \eqref{gk-house}, where it is used to show $g_k\neq 0$ for $k>|I_\Z|$, which is needed for the Inner Loop to construct $B_k$, and second, in the paragraph above \eqref{particular}, where it is used in the case when $\beta_k=1$, which is needed to show (a) holds for the next iteration  $U=V'\bdot A_{\mathfrak t-1}\bdot A_\mathfrak t\bdot\ldots\bdot A_{\mathfrak s}$. Note the latter only occurs when $k\geq |I_\Z|+1$ or when $k=1=|I_\Z|$, $|I_\Q|=|I_\diamond|=0$ and $V=U$.

One can modify the Outer Loop  by requiring it to prematurely halt during the construction of $A_{\mathfrak t-1}$ if ever $k=|I_\Z|$ is reached. In this way, both cases where the hypothesis $0\notin \C^*(G_0)$ is used are never encountered. To avoid confusion, let $N_T^{(\mathfrak t)}$ and $N_S^{(\mathfrak t)}$ denote the recursively defined constants for which (a)--(d) holds for $U=V\bdot A_t\bdot\ldots\bdot A_{\mathfrak s}$ (assuming the process did not prematurely halt before $A_t$ could be constructed), and let $N_T$ and $N_S$ denote the final values, so $N_S=N_S^{(1)}$ and $N_T=\max\{N_T^{(1)}, \Sum{j=1}{\ell}|W_j|\}$, where $W_1,\ldots,W_\ell\in \mathcal A(G_0)$ are the distinct atoms having no nonempty $X\subseteq \supp(W_j)$ with $X\in X(G_0)$, of which there are a finite number as noted at the beginning of the proof of Theorem \ref{thm-structural-Lambert} (which only requires the assumption that $G_0^\diamond$ is finitary, not that $0\notin \C^*(G_0^\diamond)$). Note $N_T$ and $N_S$ still exist, with their values independent of $U$, even if we cannot complete the construction of all $A_{\mathfrak t-1}$ for some atom $U$. Indeed, $N_T^{(\mathfrak t)}$ and $N_S^{(\mathfrak t)}$ are simply  defined by the dual recursion  given by \eqref{recursive-nT} and \eqref{recur-NS}, so $N_T^{(\mathfrak t-1)}=N'_T$ and $N_S^{(\mathfrak t-1)}=N'_S$ with $N'_T$ and $N'_S$ defined by \eqref{recursive-nT} and \eqref{recur-NS} using $N_S^{(\mathfrak t)}$ and $N_T^{(\mathfrak t)}$ in place of  $N_S$ and $N_T$ in the formulas \eqref{recursive-nT} and \eqref{recur-NS} (and initial values $N_T^{(\mathfrak s+1)}=N_S^{(\mathfrak s+1)}=0$). Since all constants used in the recursive formulas \eqref{recursive-nT} and \eqref{recur-NS} depend only on the individual lattice types $\varphi$ for the finitary subset $G_0$, the resulting values of $N_T$ and $N_S$ are well-defined.

If, for the atom $U$, the Outer Loop process never prematurely halts during the construction of any $A_{t-1}$, then we obtain $\vp_{G_0\setminus G_0^\diamond}(U)\leq (\mathfrak s+1)N_T$, as we did when $0\notin \C^*(G_0^\diamond)$. Now instead suppose, for the atom $U$, that the process prematurely halts during the construction of $A_{\mathfrak t-1}$ for some $\mathfrak t\in [2,\mathfrak s+1]$, so we have $U=V\bdot \ A_{\mathfrak t}\bdot\ldots\bdot A_{\mathfrak s}$ satisfying (a)--(d)  but fail to construct $A_{\mathfrak t-1}$. There are two ways this failure can arise. First, we may have completed the Outer Loop process with it ending when $\beta_k=1$,
 $k=1=|I_\Z|$, $|I_\Q|=|I_\diamond|=0$ and $V=U$, leaving us unable to conclude (a) holds for $A_{\mathfrak t-1}$. In such case,  $V_{G_0\setminus \mathfrak X}=T_1$ with $\vp_{G_0\setminus G_0^\diamond}(U)=\vp_{G_0\setminus G_0^\diamond}(V)\leq |V_{G_0\setminus \mathfrak X}|\leq |T_1|\leq D_\Delta M_W\leq N^{(\mathfrak t-1)}_T\leq N_T$ by \eqref{Iz-integral} and \eqref{recursive-nT}, as desired.
  Second, the Outer Loop process has constructed $B_1,\ldots,B_{|I_\Z|}$ and has not finished, so the terminal value $k>|I_\Z|$, yet $g_{|I_\Z|+1}=0$, leaving us unable to construct the next sequence $B_{|I_\Z|+1}$ using the Inner Loop. In this case, since $\supp(T_i)\subseteq G_0^\diamond$ for $i\in I_\diamond$ by \eqref{T-diamondpart}, and since $\supp(V_\mathfrak X)\subseteq \mathfrak X\subseteq G_0^\diamond$, we obtain \begin{align}\nn\vp_{G_0\setminus G_0^\diamond}(V)&\leq \Summ{i\in I_\Z\cup I_\Q}|T_i|\leq |I_\Z|D_\Delta M_W+|I_\Q|M_W\leq (k-1) D_\Delta M_W+(D_\Delta+N^{(t)}_S+2\ell_\varphi)M_W\\
&\leq ((2^{n+1}-1)(1+N^{(t)}_S)D_\Delta+N^{(t)}_S+2\ell_\varphi)M_W\leq  N^{(\mathfrak t-1)}_T\leq N_T,\label{cutoff}
\end{align}
with the second inequality in view of \eqref{Iz-integral} and \eqref{TkQBound}, with the third inequality in view of $k>|I_\Z|$ and Claim A,  with fourth in view of the upper bound $k\leq (2^{n+1}-1)(1+N^{(\mathfrak t)}_S)$, estimating how long the Outer Loop Process can run, and with the fifth in view of \eqref{recursive-nT}. Thus $\vp_{G\setminus G_0^\diamond}(U)=\vp_{G\setminus G_0^\diamond}(V)+\Sum{i=\mathfrak t}{\mathfrak s}\vp_{G\setminus G_0^\diamond}(A_i)\leq N_T^{(\mathfrak t-1)}+(\mathfrak s-\mathfrak t+1)N^{(\mathfrak t)}_T\leq (\mathfrak s+1)N_T$ now follows from (c) holding for $U=V\bdot A_\mathfrak t\bdot \ldots\bdot A_\mathfrak s$.
In all cases,  $G_0\setminus G_0^\diamond\subseteq G_0$ is a Lambert subset with bound $(\mathfrak s+1)N_T$, completing the proof.
\end{proof}

We can now extend Theorem \ref{thm-neg-char} to characterize the set $G_0^\diamond$ in terms of $\Z_+$-linear combinations and atoms, rather than $\Q_+$-linear combinations and elementary atoms.

\begin{corollary}\label{cor-G_0diamond-Zequiv}
Let $\Lambda\leq \R^d$ be a full rank lattice, where $d\geq 0$,  and let $G_0\subseteq \Lambda$ be a finitary subset with $\C(G_0)=\R^d$. Then \
\ber\nn G_0^\diamond&=&\{g\in G_0:\; \sup\{\vp_{g}(U):\; U\in \mathcal A^{\mathsf{elm}}(G_0)\}=\infty\}\\\nn &=&\{g\in G_0:\; \sup\{\vp_{g}(U):\; U\in \mathcal A(G_0)\}=\infty\}.\eer
\end{corollary}

\begin{proof}
The equality $G_0^\diamond=\{g\in G_0:\; \sup\{\vp_{g}(U):\; U\in \mathcal A^{\mathsf{elm}}(G_0)\}=\infty\}$ holds by Theorem \ref{thm-neg-char}. Since $\mathcal A^{\mathsf{elm}}(G_0)\subseteq \mathcal A(G_0)$,  the inclusion $\{g\in G_0:\; \sup\{\vp_{g}(U):\; U\in \mathcal A^{\mathsf{elm}}(G_0)\}=\infty\}\subseteq \{g\in G_0:\; \sup\{\vp_{g}(U):\; U\in \mathcal A(G_0)\}=\infty\}$ is trivial. To see the reverse inclusion, let $g\in G_0$ with $\sup\{\vp_{g}(U):\; U\in \mathcal A(G_0)\}=\infty$. Then Corollary \ref{cor-structural-Lambert} ensures $g_0\in G_0^\diamond=\{g\in G_0:\; \sup\{\vp_{g}(U):\; U\in \mathcal A^{\mathsf{elm}}(G_0)\}=\infty\}$.
\end{proof}

We can now achieve one of our mains goals: the characterization of finite elasticities for subsets $G_0\subseteq \Z^d$. Note that Theorem \ref{thm-rho-char}.1 characterizes $\rho_{d+1}(G_0)<\infty$ in terms of a basic, combinatorial property of $G_0$, one which trivially implies $\rho_{d+1}(G_0)<\infty$ (cf. Proposition \ref{prop-pre-rho-char}), while  Theorem \ref{thm-rho-char}.5 characterizes $\rho_{d+1}(G_0)<\infty$ in terms of the subset $G_0^\diamond\subseteq G_0$, which can be defined purely in terms of Convex Geometry (cf. Proposition \ref{prop-G_0diamond-1st-easy-equiv}).

\begin{theorem}\label{thm-rho-char}
Let $\Lambda\leq \R^d$ be a full rank lattice in $\R^d$, where $d\geq 0$, and let $G_0\subseteq \Lambda$ be a subset. Then the following are equivalent.
\begin{itemize}
\item[1.] There exists a subset $X\subseteq G_0$ such that $\mathcal A(X)=\emptyset$ and $G_0\setminus X\subseteq G_0$ is a Lambert subset.
\item[2.] $\rho(G_0)<\infty$
\item[3.] $\rho_k(G_0)<\infty$ for all $k\geq 1$
\item[4.] $\rho_{d+1}(G_0)<\infty$.
\item[5.] $0\notin \C^*( G_0^\diamond)$.
\end{itemize}
\end{theorem}

\begin{proof} The implications $1.\Rightarrow 2.\Rightarrow 3.\Rightarrow 4.$ follow by Proposition \ref{prop-pre-rho-char}.   Let $\mathcal E=\C(G_0)\cap -\C(G_0)$ be the lineality space of $\C(G_0)$. Let $\wtilde G_0=\{g\in G_0:\; g\in \supp(U)\mbox{ for some $U\in \mathcal A(G_0)$}\}$. Corollary \ref{cor-nondegen} implies  \be\label{begit}G_0\cap \mathcal E=\wtilde G_0\quad\und\quad\C(\wtilde G_0)=\C(G_0)\cap -\C(G_0)=\mathcal E.\ee Thus, by definition of $G_0^\diamond$, we have $G_0^\diamond=\wtilde G_0^\diamond$.
Note that  $\wtilde \Lambda =\Lambda\cap\mathcal E\leq \Lambda$ is a sublattice. Since $\wtilde G_0\subseteq \wtilde \Lambda$ is a subset which linearly generates $\mathcal E$, it follows that $\wtilde \Lambda\leq \mathcal E$ is a full rank lattice in $\mathcal E$. Let $\tilde d=\dim \mathcal E$.
Now suppose we knew Theorem \ref{thm-rho-char} held for $\wtilde G_0\subseteq \wtilde \Lambda\leq \mathcal E$. Then $\rho_{d+1}(G_0)<\infty$ would imply $\rho_{\tilde d+1}(\wtilde G_0)\leq \rho_{d+1}(\wtilde G_0)=\rho_{d+1}(G_0)<\infty$, with the first inequality in view of $\tilde d\leq d$ and \eqref{rho-ascend-chain}, in turn implying $0\notin \C^*(\wtilde G_0^\diamond)=\C^*(G_0^\diamond)$, in turn implying there is a subset $\wtilde X\subseteq \wtilde G_0$ with $\mathcal A(\wtilde X)=\emptyset$ and $\wtilde G_0\setminus \wtilde X\subseteq \wtilde G_0$  a Lambert subset.
But then clearly $X=\wtilde X\cup (G_0\setminus \wtilde G_0)\subseteq G_0$ is a subset with $\mathcal A(X)=\emptyset$ and $G_0\setminus X=\wtilde G_0\setminus \wtilde X\subseteq G_0$ a Lambert subset.
Thus the remaining implications follow from the case when $G_0=\wtilde G_0$, which we now assume. In this case, the implication $4.\Rightarrow 5.$ follows by Proposition \ref{prop-pre-rho-char-diamond} applied to $G_0=\wtilde G_0$, which we can apply in view of \eqref{begit} and  $\rho_{\tilde d+1}(G_0)\leq \rho_{d+1}(G_0)<\infty$ (by \eqref{rho-ascend-chain}).  It remains to prove $5.\Rightarrow 1.$
To this end, suppose $0\notin \C^*(G_0)$. Then $G_0$ is finitary by Theorem \ref{thm-keylemmaII} and \eqref{begit}, whence Corollary \ref{cor-structural-Lambert} and \eqref{begit} imply $G_0\setminus G_0^\diamond\subseteq G_0$ is a Lambert subset. Since $0\notin \C^*(G_0^\diamond)$, Proposition \ref{prop-rational-atoms} implies $\mathcal A(G_0^\diamond)=\emptyset$. Thus Item 1 holds with $X=G_0^\diamond$, completing the proof.
\end{proof}

Let $G_0\subseteq G$ be a subset of an abelian group $G$. For $k\geq 1$, we define the elementary elasticity $\rho_k^{\mathsf{elm}}(G_0)$ to be the minimal integer $N$ such that, if $$U_1\bdot\ldots\bdot U_k=V_1\bdot\ldots\bdot V_\ell$$ with $U_1,\ldots,U_k\in \mathcal A^{\mathsf{elm}}(G_0)$ and $V_1,\ldots,V_\ell\in \mathcal A(G_0)$, then $\ell\leq N$. If no such $N$ exists, then we set $\rho_k^{\mathsf{elm}}(G_0)=\infty$.
It is readily noted (by the same simple argument establishing \eqref{rho-ascend-chain}) that $\rho_1^{\mathsf{elm}}(G_0)\leq \rho_2^{\mathsf{elm}}(G_0)\leq \ldots$. Indeed, if $U_1\bdot \ldots\bdot U_{k-1}=V_1\bdot\ldots\bdot V_{\ell}$ with $U_1,\ldots,U_{k-1}\in \mathcal A^{\mathsf{elm}}(G_0)$ and $V_1,\ldots,V_{\ell}\in \mathcal A(G_0)$, then the factorization $U_1\bdot \ldots\bdot U_{k-1}\bdot U_1=V_1\bdot\ldots\bdot V_{\ell}\bdot U_1$ shows $\rho^{\mathsf{elm}}_k(G_0)\geq \ell+1$. If $\rho_{k-1}^{\mathsf{elm}}(G_0)$ is finite, then we may take $\ell=\rho_{k-1}^{\mathsf{elm}}(G_0)$, and if  $\rho_{k-1}^{\mathsf{elm}}(G_0)$ is infinite, then we may take $\ell$ to be arbitrarily large. In either case,  $\rho_k^{\mathsf{elm}}(G_0)\geq \rho_{k-1}^{\mathsf{elm}}(G_0)$ follows.
 Trivially, we also have $$\rho_k^{\mathsf{elm}}(G_0)\leq \rho_k(G_0).$$ Thus if $\rho_k(G_0)$ is finite, then so is $\rho_k^{\mathsf{elm}}(G_0)$.
As a consequence of Theorem \ref{thm-rho-char}, we have the following converse.

\begin{corollary}
\label{cor-rhoelem} Let $\Lambda\subseteq \R^d$ be a full rank lattice, where $d\geq 0$, and let $G_0\subseteq \Lambda$ be a subset. Then $\rho_{d+1}(G_0)<\infty$ if and only if $\rho_{d+1}^{\mathsf{elm}}(G_0)<\infty$.
\end{corollary}

\begin{proof}We may w.l.o.g. suppose every $g\in G_0$ is contained in an atom,
whence $\C(G_0)=\mathcal E:=\C(G_0)\cap -\C(G_0)$ is the lineality space of $\C(G_0)$ by Corollary \ref{cor-nondegen}. Let $\wtilde d=\dim \mathcal E\leq d$.
As already remarked, one direction of the corollary is trivial. Suppose $\rho_{d+1}^{\mathsf{elm}}(G_0)<\infty$. Then $\rho_{\wtilde d+1}^{\mathsf{elm}}(G_0)\leq \rho_{d+1}^{\mathsf{elm}}(G_0)<\infty$.
If $0\notin \C^*(G_0^\diamond)$, then $\rho_{d+1}(G_0)<\infty$ follows by Theorem \ref{thm-rho-char}.3, as desired. Therefore we may assume $0\in \C^*(G_0^\diamond)$, in which case Carth\'edory's Theorem implies that there is a minimal positive basis $X\subseteq G_0^\diamond$ with $s=|X|\leq d+1$.  Let $V\in \mathcal A^{\mathsf{elm}}(X)$ be the unique elementary atom with support $X=\{g_1,\ldots,g_s\}\subseteq G_0^\diamond$ (by Proposition \ref{prop-char-minimal-pos-basis}) and let $m=\max\{\vp_{g_j}(V):\;j\in [1,s]\}$.  By Proposition \ref{prop-diamond-basic-containment}, every $g\in G_0^\diamond$ has $\supp\{\vp_g(U):\;U\in \mathcal A^{\mathsf{elm}}(G_0)\}=\infty$. Thus, for every $i\geq 0$ and $j\in [1,s]$, we can find an elementary atom $U_i^{(j)}\in \mathcal A^{\mathsf{elm}}(G_0)$ with $\vp_{g_j}(U_i^{(j)})\geq i$. Then $V^{[\lfloor i/m\rfloor]}\mid U_i^{(1)}\bdot\ldots U_i^{(s)}$ for every $i\geq 1$. Hence, for every $i\geq 1$,  we have a factorization $$U_i^{(1)}\bdot\ldots\bdot U_i^{(s)}=\Big({\prod}_{j=1}^{\lfloor i/m\rfloor} V\Big)\bdot B_i$$ for some $B_i\in \mathcal B(G_0)$. Since $i\rightarrow \infty$ and $m$ is a fixed constant, it follows that $\rho^{\mathsf{elm}}_s(G_0)=\infty$. However, since $s\leq d+1$, we have $\rho_s^{\mathsf{elm}}(G_0)\leq \rho_{d+1}^{\mathsf{elm}}(G_0)$, so that $\rho_{d+1}^{\mathsf{elm}}(G_0)=\infty$ as well, contrary to assumption.
\end{proof}

\subsection{The Structure of Atoms and Arithmetic Invariants}

In this section, we give a weak structural description of the atoms $U\in \mathcal A(G_0)$ when $0\notin \C^*(G_0)$, equivalently, when $\rho_{d+1}(G_0)<\infty$. We then use this result, along with our characterization of finite elasticities and prior machinery, to show that having finite elasticities implies nearly all  standard invariants of Factorization Theory are also finite, meaning the hypothesis $\rho_{d+1}(G_0)<\infty$ is sufficient to guarantee factorizations are as well-behaved as could be hoped.

We begin by describing what a set $X$ with minimal type $\varphi$ looks like when expressed using the lattice basis $Z_\varphi$.

\begin{lemma}\label{lem-finitary-minrep}
Let $\Lambda\leq \R^d$ be a full rank lattice, where $d\geq 0$,  let $G_0\subseteq \Lambda$ be a finitary subset with  $\C(G_0)=\R^d$,  let $\varphi\in \mathfrak X_m(G_0)$ be a minimal type with codomain $Z_\varphi=Z_1\cup \ldots\cup Z_s$, and let $X\in \mathfrak X_m(\varphi)$ be a subset with minimal type $\varphi$ and associated maximal series decomposition $X=X_1\cup \ldots\cup X_s$. For $k\in [1,s]$, let $z_1^{(k)},\ldots z_{t_k}^{(k)}\in Z_k$ and $x_1^{(k)},\ldots,x_{t_k}^{(k)}\in X_k$ be the distinct elements of $Z$ and $X$, respectively, indexed so that $\varphi(x_i^{(k)})=z_i^{(k)}$ for all $k\in [1,s]$ and $i\in [1,t_k]$. Then, for every $k\in [1,s]$ and $i\in [1,t_k]$, we have
\be\label{minrep}x_i^{(k)}=z_i^{(k)}-
\Sum{\kappa=1}{k-1}\Sum{\iota=1}{t_\kappa}\xi_\iota^{(\kappa)} z_\iota^{(\kappa)}
\quad\mbox{ for some }\quad \xi_\iota^{(\kappa)}\in \Z_+.\ee
\end{lemma}

\begin{proof}
For $k\in [0,s]$, let $\Delta_k=\Z\la Z_1\cup \ldots\cup Z_k\ra$ and let $\pi_k:\R^d\rightarrow \R\la \Delta_k\ra^\bot$ be the orthogonal projection. Set $\Delta=\Delta_s$.  Since $X$ has minimal type $\varphi$, we have $\Z\la X_1\cup \ldots\cup X_j\ra=\Delta_j$ for all $j\in [1,s]$. Let $k\in [1,s]$ and $i\in [1,t_k]$ be arbitrary.
Let $Y=Z_\varphi\setminus \{z_i^{(k)}\}\cup \{x_i^{(k)}\}$. By the interchangeability property of $\mathfrak X_m(\varphi)$ (Proposition \ref{prop-finitary-mintype}), we have $Y\in \mathfrak X_m(\varphi)$. Thus $Y$ also has minimal type $\varphi$, and so $Z_\varphi\subseteq \C_\Z(Z_\varphi)\subseteq \C_\Z(Y)$. In particular, we can write  $z_i^{(k)}$ as a positive integer combination of the elements from $Y$, say $$z_i^{(k)} =\Sum{\kappa=1}{k-1}\Sum{\iota=1}{t_\kappa}\xi_\iota^{(\kappa)} z_\iota^{(\kappa)}+\xi_i^{(k)}x_i^{(k)}+\underset{\iota\neq i}{\Sum{\iota=1}{t_k}}\xi_{\iota}^{(k)}z_i^{(k)}+
\Sum{\kappa=k+1}{s}\Sum{\iota=1}{t_\kappa}\xi_\iota^{(\kappa)} z_\iota^{(\kappa)}\quad\mbox{ for some $\xi_\iota^{(\kappa)}\in \Z_+$}.$$ The coefficients of $z_\iota^{(\kappa)}$ with $\kappa\geq k+1$ are all zero in view of the linear independence of $Z_{k+1}\cup\ldots\cup Z_s$ modulo $\R\la \Delta_k\ra$ (cf. Proposition \ref{prop-reay-basis-properties}). Since $X$ and $Z_\varphi$ share the same minimal type $\varphi$ with $\varphi(x_i^{(k)})=z_i^{(k)}$, we have $z_i^{(k)}-x_i^{(k)}\in \R\la \Delta_{k-1}\ra\cap \Delta=\Delta_{k-1}$. In particular, $x_i^{(k)}$ and $z_i^{(k)}$ are equal modulo $\R\la \Delta_{k-1}\ra$, in which case the linear independence of $Z_k$ modulo $\R\la \Delta_{k-1}\ra$ ensures that
the coefficients of $z_\iota^{(k)}$ with $\iota\neq i$ are also zero, and that the coefficient of $x_i^{(k)}$ is one. The result now follows.
\end{proof}

For the next proposition, we need to view our sequences as indexed, so that two distinct terms of a sequence $S\in \Fc(G_0)$ that are equal as elements can still be viewed as  distinct terms in the sequence $S$. We follow the notation introduced in Section \ref{sec-intro-factorization}.

\begin{proposition}\label{prop-swapping}
Let $\Lambda\leq \R^d$ be a full rank lattice, where $d\geq 0$,  let $G_0\subseteq \Lambda$ be a finitary subset with  $\C(G_0)=\R^d$,
let $\varphi\in \mathfrak X_m(G_0)$ be a minimal type with codomain $Z_\varphi=Z_1\cup \ldots\cup Z_s$, let $\mathfrak X=\mathfrak X_m^\cup (\varphi)=\bigcup_{X\in \mathfrak X_m(\varphi)}X$, let $S=g_1\bdot\ldots\bdot g_\ell\in \Fc(\mathfrak X)$ be a sequence, and  for every $k\in [1,s]$, let $I_k\subseteq [1,\ell]$ be the subset of all $a\in [1,\ell]$ with $g_a$ at depth $k$. Suppose $\sigma(S)\in \C(Z_\varphi)$. Then
there exists a system of subsets $T_a\subseteq [1,\ell]$ for $a\in [1,\ell]$  such that the following hold.
\begin{itemize}
\item[1.]   For every $a\in [1,\ell]$, say $a\in I_k$, there is some $\partial(a)\subseteq I_1\cup\ldots\cup I_{k-1}$ such that $T_a=\{a\}\cup \bigcup_{b
  \in \partial(a)}T_b$ is a disjoint union.
\item[2.] If $b\notin T_a$ and $a\notin T_b$, then $T_a\cap T_b=\emptyset$, for every $a,\,b\in [1,\ell]$.
 \item[3.] $\sigma(S(T_a))=\varphi(g_a)$, for every $a\in [1,\ell]$.
 \end{itemize}
In particular, for any system of sets $T_a\subseteq [1,\ell]$ satisfying Items 1--3, the following hold.
\begin{itemize}
\item[4.] If $a\in I_k$, then $T_a\setminus \{a\}\subseteq I_1\cup \ldots\cup I_{k-1}$.
\item[5.] If $b\in T_a$, then $T_b\subseteq T_a$.
\item[6.] For  every $k\in [1,s]$, there exists a subset $J_k\subseteq I_1\cup \ldots\cup I_k$ such that $\bigcup_{i\in J_k}T_i=I_1\cup\ldots\cup I_k$ is a disjoint union, ensuring that $${\prod}^\bullet_{i\in [1,k]}S(I_i)={\prod}^\bullet_{i\in J_k} S(T_i)$$ is a factorization of the subsequence of all terms in $S$ with depth at most $k$ into a product of subsequences whose sums each lie in  $Z_\varphi$.
\end{itemize}
\end{proposition}

\begin{proof}
Item 4 follows from  Item 1 and a quick inductive argument on $k=1,2,\ldots,s$. Likewise, Item 5 follows from Item 1  and a quick inductive argument on $k$, where $a\in I_k$.  To show Item 6 holds, let  $k\in [1,s]$. The set $J_k$ is then constructed recursively in $k$ steps. Set $J_k^{(k)}=I_k$. Assuming $J_k^{(j)}$ has been constructed, with $j\in [2,k]$, define $J_k^{(j-1)}=J_k^{(j)}\cup \Big(I_{j-1}\setminus \bigcup_{c\in J_k^{(j)}}T_c\Big)\subseteq I_{j-1}\cup I_{j}\cup  \ldots\cup  I_k$. Set $J_k=J_k^{(1)}$. In view of Item 4, any $a\in I_j$ with $j\leq k$ has $T_a\subseteq I_1\cup \ldots\cup I_k$. Thus $J_k=J_k^{(1)}\subseteq I_1\cup\ldots \cup I_k$ and $\bigcup_{i\in J_k}T_i=I_1\cup\ldots\cup I_k$ by construction.
 By Item 4, we  have $b\notin T_a$ and $a\notin T_b$ for any distinct $a,\,b\in I_j$ and $j\in [1,s]$.  Hence Item 2 ensures that $\bigcup_{i\in I_j}T_i$ is a disjoint union for any $j\in [1,s]$.
In particular, $\bigcup_{i\in J_k^{(k)}}T_i$ is a disjoint union. Also, if $a\in J_k^{(j)}$ and $b\in I_{j-1}\setminus \bigcup_{c\in J_k^{(j)}}T_c$, then $b\notin T_a$. Since $b\in I_{j-1}$ ensures that $T_b\subseteq I_1\cup\ldots\cup I_{j-1}$ (by Item 4), while  $a\in J_k^{(j)}\subseteq I_j\cup \ldots\cup I_k$, we also have $a\notin T_b$.  Thus Item 2 ensures that $T_a$ is disjoint from $T_b$, for every $a\in J_k^{(j)}$. We also have $T_b$ disjoint from any other $T_{b'}$ with $b'\in I_{j-1}$ as already mentioned. An inductive argument now shows that $\bigcup_{i\in J_k^{(j)}}T_i$ is a disjoint union for $j=k,k-1,\ldots,1$. The case $j=1$ implies $\bigcup_{i\in J_k^{(1)}}T_i=\bigcup_{i\in J_k}T_i$ is a disjoint union, and we now see that Item 6 follows. It remains to establish Items 1--3.

For $k\in [1,s]$, let $z^{(k)}_1,\ldots z^{(k)}_{t_k}\in Z_k$ be the distinct elements of $Z_k$.
Let $\mathfrak X=\bigcup_{k=1}^s\bigcup_{i=1}^{t_k}\mathfrak X_{k,i}$, with $\mathfrak X_{k,i}\subseteq \mathfrak X$ consisting of all $x\in \mathfrak X$ with $\varphi(x)=z_i^{(k)}$.
Given any $x\in \R\la Z_\varphi\ra$, we have $$x=\Sum{k=1}{s}\Sum{i=1}{t_k}\xi_{k,i}(x)z_i^{(k)}\quad\mbox{ for some uniquely defined $\xi_{k,i}(x)\in \R$},$$ as $Z_\varphi\in \mathfrak X_m(\varphi)$ is linearly independent. If $x\in \mathfrak X$, then Lemma \ref{lem-finitary-minrep} ensures that $\xi_{k,i}(x)\in \Z$ for all $k$ and $i$. Moreover, if $\varphi(x)=z_\iota^{(\kappa)}$, then $\xi_{k,i}(x)=0$  whenever $k\geq\kappa$, apart from the value $\xi_{\kappa,\iota}(x)=1$, and $\xi_{k,i}(x)\leq 0$ whenever $k<\kappa$.

We construct the subsets $\partial(a)$ and $T_a$ recursively, first for all $a\in I_1$, then for all $a\in I_2$, and so forth. For $a\in I_1$, define $\partial(a)=\emptyset$ and $T_a=\{a\}$. Then $\partial(a)\subseteq I_1\cup \ldots\cup I_0=\emptyset$, and Items 1 and 2 trivially hold. Moreover, since $g_a\in X$ for some $X\in \mathfrak X_m(\varphi)$ with $\varphi(g_a)\in Z_1$ (since $a\in I_1$), it follows from the definition of lattice type (each minimal type is a refinement of a lattice type) that $g_a=\varphi(g_a)$. Thus Item 3 holds. Assume the $T_a$ and $\partial(a)$ have been constructed for all $a\in I_1\cup \ldots\cup I_{\kappa-1}$, where $\kappa\geq 2$, such that Items 1--3 hold.
Let $J_{\kappa-1}\subseteq I_1\cup \ldots\cup I_{\kappa-1}$ be a subset such that $\bigcup_{i\in J_{\kappa-1}}T_i=I_1\cup \ldots\cup I_{\kappa-1}$ is a disjoint union with $${\prod}^\bullet_{i\in [1,\kappa-1]}S(I_i)={\prod}^\bullet_{i\in J_{\kappa-1}} S(T_i)$$ a factorization of the subsequence of all terms in $S$ with depth at most $\kappa-1$ into a product of subsequences whose sums each lie in  $Z_\varphi$, which exists by the argument used to derive Item 6 applied to the subsequence of $S$ consisting of all terms having depth at most $\kappa-1$.
For $k\in [1,\kappa-1]$ and $j\in [1,t_k]$, let $J_{\kappa-1}^{(k,j)}\subseteq J_{\kappa-1}$ be all those $b\in J_{\kappa-1}$ with $\sigma(S(T_b))=z_j^{(k)}$. Then $J_{\kappa-1}=\bigcup_{k=1}^{\kappa-1}\bigcup_{j=1}^{t_k}J_{\kappa-1}^{(k,j)}$ is a disjoint union.
We proceed to construct $\partial(a)$ and $T_a$ for $a\in I_\kappa$ as follows.

For each $k\in [1,\kappa-1]$ and $j\in [1,t_k]$, we have $\xi_{k,j}(g_a)\leq 0$, for $a\in I_\kappa$, in view of Lemma \ref{lem-finitary-minrep} as remarked earlier. Since $I_k\subseteq \bigcup_{i\in J_{\kappa-1}}T_i$ (as $k\leq \kappa-1$), it follows that ${\prod}^\bullet_{i\in J_{\kappa-1}} S(T_i)$ contains all terms $x$ from $S$ with $\xi_{k,j}(x)>0$. Thus,
  since $\sigma(S)\in \C(Z_\varphi)$ by hypothesis with $Z_\varphi$ linearly independent, we must have $\Summ{a\in I_\kappa}|\xi_{k,j}(g_a)|\leq |J_{\kappa-1}^{(k,j)}|$. Consequently, it is possible to find  disjoint subsets $D^{(k,j)}_a\subseteq J_{\kappa-1}^{(k,j)}$, for $a\in I_\kappa$, such that $|D^{(k,j)}_a|=|\xi_{k,j}(g_a)|$  for each $a\in I_\kappa$. Note $D_a^{(k,j)}\subseteq J_{\kappa-1}^{(k,j)}\subseteq J_{\kappa-1}\subseteq I_1\cup\ldots\cup I_{\kappa-1}$.
  Define $\partial(a)=\bigcup_{k=1}^{\kappa-1}\bigcup_{j=1}^{t_k}D_a^{(k,j)}\subseteq J_{\kappa-1}\subseteq I_1\cup\ldots\cup I_{\kappa-1}$ and $T_a=\{a\}\cup \bigcup_{c\in \partial(a)}T_c$.
  By construction, Items 1 and 3 both hold for all $a\in I_\kappa$, with disjointness in Item 1 following from the disjointness of the union $\bigcup_{i\in J_{\kappa-1}}T_i$.
   Since \be\label{labeltray}\bigcup_{i\in J_{\kappa-1}}T_i=I_1\cup \ldots\cup I_{\kappa-1}\ee is a disjoint union, and since $\bigcup_{a\in I_\kappa} \partial(a)\subseteq J_{\kappa-1}$ is also a disjoint union (as the $D_a^{(k,j)}$ are disjoint), it follows that $T_a\cap T_b=\emptyset$ for all distinct $a,\,b\in I_\kappa$.
  If $a\in I_\kappa$ and $b\in I_1\cup\ldots\cup I_{\kappa-1}$ with $b\notin T_a$, then $b\in T_c$ for some $c\in J_{\kappa-1}\setminus \partial(a)$ (in view of \eqref{labeltray}), in which case the disjointness of  the union in \eqref{labeltray} ensures that $T_a\cap T_c=\emptyset$.
  However, since Item 4 holds for $I_1\cup \ldots \cup I_{\kappa-1}$, it follows that $b\in T_c$ implies $T_b\subseteq T_c$.
  Hence $T_a\cap T_b=\emptyset$ follows in view of  $T_a\cap T_c=\emptyset$. Finally, if $a,\,b\in I_1\cup\ldots\cup I_{\kappa-1}$ with $a\notin T_b$ and $b\notin T_a$, then Item 2 holding for $I_1\cup \ldots\cup I_{\kappa-1}$ ensures that $T_a\cap T_b=\emptyset$. Hence Item 2 holds for $I_1\cup \ldots\cup I_{\kappa}$, and iterating the construction for $\kappa=1,2,\ldots, s$ completes the proof.
\end{proof}

Item 2 in Proposition \ref{prop-swapping} implies that either $b\in T_a$ (implying $T_b\subseteq T_a$ by Item 5) or $a\in T_b$ (implying $T_1\subseteq T_a$ by Item 5) or $T_a\cap T_b=\emptyset$. Thus any system of sets $T_a\subseteq [1,\ell]$ satisfying Items 1--3 in Proposition \ref{prop-swapping} also satisfies
\begin{itemize}
\item[$1'$.] $a\in T_a$ and $T_a\setminus \{a\}\subseteq I_1\cup\ldots\cup I_{k-1}$, where $a\in I_k$,
\item[$2'$.] either $T_a\subseteq T_b$ or $T_b\subseteq T_a$ or $T_a\cap T_b=\emptyset$,  and
\item[$3'$.]  $\sigma(S(T_a))=\varphi(g_a)$,
\end{itemize}
for any $a,\,b\in [1,\ell]$.
However, a system of sets $T_a\subseteq [1,\ell]$ satisfying Items $1'$--$3'$ must also satisfy Items 1--3, meaning these are equivalent defining conditions for  the set system $T_a\subseteq [1,\ell]$---as the following short argument shows. If $a\notin T_b$, then Item $1'$ ensures $T_a\nsubseteq T_b$. Likewise, if $b\notin T_a$, then Item $1'$ ensures $T_b\nsubseteq T_a$. As a result, if $a\notin T_b$ and $b\notin T_b$, then Item $2'$ implies $T_a\cap T_b=\emptyset$, yielding Item 2. Item 3 is the same as Item $3'$. Let $a\in [1,\ell]$.
 If $b\in T_a\setminus \{a\}$, then Item $1'$ ensures $\mathsf{dep}(g_c)\leq \mathsf{dep}(g_b)<\mathsf{dep}(g_a)$ for all $c\in T_b$, implying $a\notin T_b$. Thus $T_a\nsubseteq T_b$ (as $a\in T_a$ but $a\notin T_b$) and $T_a\cap T_b\neq \emptyset$ (as $b\in T_a\cap T_b$), in which case Item $2'$ yields $T_b\subseteq T_a$, showing Item 5 holds.  We conclude (by Items $1'$ and 5) that $b\in T_b\subseteq T_a$ for all $b\in T_a$, implying $T_a=\{a\}\cup \bigcup_{b\in T_a\setminus\{a\}}T_b$. Let  $\partial(a)\subseteq T_a\setminus \{a\}$
 consist of all $b\in T_a\setminus \{a\}$ such that there does not exist any $c\in T_a\setminus \{a,b\}$ with $T_b\subseteq T_{c}$. By Item $1'$, we have $\partial(a)\subseteq I_1\cup\ldots\cup I_{k-1}$, where $a\in I_k$. By definition of $\partial(a)$, we have $T_a=\{a\}\cup \bigcup_{b\in T_a\setminus\{a\}}T_b=\{a\}\cup \bigcup_{b\in \partial(a)}T_b$. If $b,\,c\in \partial(a)$ are distinct, then the definition of $\partial(a)\subseteq T_a\setminus \{a\}$ ensures that $T_b\nsubseteq T_c$ and $T_c\nsubseteq T_b$, so that Item $2'$ implies $T_b$ and $T_c$ are disjoint. By Item $1'$, any $b\in \partial(a)\subseteq I_1\cup \ldots\cup I_{k-1}$ has $T_b\subseteq I_1\cup \ldots\cup I_{k-1}$, ensuring that $a\notin T_b$ (as $a\in I_k$). It follows that the union $T_a=\{a\}\cup \bigcup_{b\in \partial(a)}T_b$ is disjoint, yielding Item 1, which completes the equivalence of the conditions.

It is entirely possible for the size of $|\supp(U)|$, for an atom $U\in \mathcal A(G_0)$, to be arbitrarily large even under the assumption $\rho_{d+1}(G_0)<\infty$. This means that, if we partition the terms of the sequence $U$ according to whether they are equal as elements from $G_0$, then we cannot hope to achieve a global bound on the number of partition classes we must use. However, the following theorem shows that there is a less restrictive notion of support that \emph{does} have the property that any atom $U\in \mathcal A(G_0)$ can have its terms partitioned into  at most $N+\Summ{\varphi\in \mathfrak T_m(G_0)}|Z_\varphi|\leq N+d|\mathfrak T_m(G_0)|$ types of elements, with elements of the same type behaving in the same essential manner as described by Proposition \ref{prop-swapping}. Theorem \ref{thm-structural-char} can be viewed as a weak structural description of the atoms $U\in \mathcal A(G_0)$, allowing us to effectively simulate globally bounded finite support.

\begin{theorem}\label{thm-structural-char}
Let $\Lambda\leq \R^d$ be a full rank lattice, where $d\geq 0$,  and let $G_0\subseteq \Lambda$ be a subset with  $\C(G_0)=\R^d$. Suppose $0\notin \C^*(G^\diamond_0)$. Then there exists a bound $N\geq 0$ such that any atom $U\in \mathcal A(G_0)$ has a factorization $$U=R\bdot {\prod}^\bullet_{\varphi\in \mathfrak T_m(G_0)} S_{\varphi},\quad
\mbox{ with } \quad R\in \Fc(G_0)\quad \und\quad S_\varphi\in \Fc(\mathfrak X_m^\cup (\varphi)),$$
such that $|R|\leq N$ and $\sigma(S_\varphi)\in \C_\Z(Z_\varphi)$ for all $\varphi\in \mathfrak T_m(G_0)$, where $Z_\varphi$ is the codomain of $\varphi$. In particular, all conclusions of Proposition \ref{prop-swapping} hold for each $S_\varphi$.
\end{theorem}

\begin{proof} Since $0\notin \C^*(G_0)$, Theorem \ref{thm-keylemmaII} ensures that  $G_0$ is finitary. Let $U\in \mathcal A(G_0)$ be an arbitrary atom and
let $U=A_0\bdot\prod_{\varphi\in \mathfrak T_m(G_0)}A_\varphi$ be the factorization of $U$ given by Theorem \ref{thm-structural-Lambert}, with bounds $N_S\geq 0$ and $N_T\geq 0$ and $A_0,\,A_\varphi\in \mathcal B_{\mathsf{rat}}(G_0)$ (with $A_\varphi$ the trivial sequence when $\varphi\in \mathfrak T_m(G_0)$ is the trivial type). For each $\varphi\in \mathfrak T_m(G_0)$, we define the subsequence $S_\varphi\mid \lfloor A_\varphi\rfloor$ as follows.

For the trivial type $\varphi$, we set $S_\varphi$ to be the trivial sequence.
Let  $\varphi\in \mathfrak T_m(G_0)$ be an arbitrary, nontrivial type, let $Z_\varphi=Z_1\cup\ldots\cup Z_s$ be the codomain of $\varphi$ and associated maximal series decomposition, and let $\mathfrak X=\mathfrak X_m^\cup(\varphi)=\bigcup_{X\in \mathfrak X_m(\varphi)}X$.  For $k\in [1,s]$, let $z_1^{(k)},\ldots z_{t_k}^{(k)}\in Z_k$ be the distinct elements of $Z_k$.
Let $\mathfrak X=\bigcup_{k=1}^s\bigcup_{i=1}^{t_k}\mathfrak X_{k,i}$, with $\mathfrak X_{k,i}\subseteq \mathfrak X$ consisting of all $x\in \mathfrak X$ with $\varphi(x)=z_i^{(k)}$. Let $\mathfrak X_k=\bigcup_{i\in [1,t_k]}\mathfrak X_{k,i}$, $$A=A_\varphi\quad\und\quad A=A_\mathfrak X\bdot A_{G_0\setminus \mathfrak X},$$ where $A_\mathfrak X\mid A$ is the rational subsequence consisting of all terms from $\mathfrak X$, and $A_{G_0\setminus \mathfrak X}\mid A$ is the rational subsequence consisting of all terms from $G_0\setminus \mathfrak X$.
For $k\in [1,s]$ and $i\in [1,t_k]$, let $A_{\mathfrak X_{k,i}}\mid A_\mathfrak X$ be the rational subsequence consisting of all terms $x$ with $\varphi(x)=z_i^{(k)}$, i.e., all terms $x\in \mathfrak X_{k,i}$, and let $A_{\mathfrak X_k}$ be the rational subsequence consisting of all terms from $\mathfrak X_k$.  By Theorem \ref{thm-structural-Lambert}(d), we have $\sigma(A_{G_0\setminus \mathfrak X})\in -\C(Z_\varphi)$, whence
\be\label{corun} \sigma(A_\mathfrak X)\in \C(Z_\varphi)=\C\big(\{z_i^{(k)}:\; k\in [1,s],\,i\in [1,t_k]\}\big)\ee follows in view of  $A=A_\varphi\in \mathcal B_{\mathsf{rat}}(G_0)$ being zero-sum.

Given any $x\in \R\la Z_\varphi\ra$, which includes any $x\in \mathfrak X$, we have \be\label{explanation}x=\Sum{k=1}{s}\Sum{i=1}{t_k}\xi_{k,i}(x)z_i^{(k)}\quad\mbox{ for some uniquely defined $\xi_{k,i}(x)\in \R$},\ee as $Z_\varphi$ is linearly independent.
If $x\in \mathfrak X$, then Lemma \ref{lem-finitary-minrep} ensures that $\xi_{k,i}\in \Z$ for all $k$ and $i$. Moreover, if $x\in \mathfrak X_{\kappa,\iota}$, then $\xi_{k,i}(x)=0$  whenever $k\geq\kappa$, apart from the value $\xi_{\kappa,\iota}(x)=1$, and $\xi_{k,i}(x)\leq 0$ whenever $k<\kappa$. Consequently, $\xi_{k,i}(x)\leq 0$ for all $x\in \supp(A_\mathfrak X)$, apart from those terms $x\in \mathfrak X_{k,i}$, for which we instead have $\xi_{k,i}(x)=1$.

We proceed to iteratively define rational subsequences $D_r\mid A_\mathfrak X$, for $r=0,1,\ldots,s-1$, such that $D_0\mid D_1\mid  \ldots\mid D_{s-1}$ and
\begin{itemize}
\item[(a)] $A_\mathfrak X \bdot D_r^{[-1]}\in \Fc(\mathfrak X)$,
\item[(b)] $\xi_{k,i}(\sigma(A_\mathfrak X\bdot D_r^{[-1]}))\in \Z_+$ for all $k\in [1,r]$ and $i\in [1,t_k]$,
\item[(c)] $D_r^{(j)}=D_{k}^{(j)}$ for all $k\in [0,r-1]$ and $j\in [1,k+1]$,
\item[(d)] $|D_r|\leq N_S+(|D^{(1)}_0|+|D^{(2)}_1|+\ldots+|D^{(r)}_{r-1}|)$, and
\item[(e)] $|D^{(k)}_{k-1}|\leq N_S$ for $k\in [1,r]$,
\end{itemize}
where $D_k^{(j)}\mid D_k$ denotes the subsequence of all terms from $\mathfrak X_j$, for $k\in [0,s-1]$ and $j\in[1,s]$.
Set $D_0=\{A_\mathfrak X\}$. Then (a)--(e) hold with $|D_0|=|\{A_\mathfrak X\}|\leq |\supp(\{A_\mathfrak X\})|\leq N_S$ in view of Theorem \ref{thm-structural-Lambert}(b) and \eqref{RatSeq1}. Assume we have already constructed the sequences $D_0,\ldots,D_{r-1}$, with $r\geq 1$. Then we construct $D_r$ as follows.

For $k\in [1,s]$ and $i\in [1,t_k]$, let $D_{r-1}^{(k,i)}\mid D_{r-1}$ denote the rational subsequence of all terms from $\mathfrak X_{k,i}$. The following makes implicit use of the comments after \eqref{explanation}.  All terms $x\in \supp(A_{\mathfrak X})$ with depth less than $r$ have $\xi_{r,i}(x)=0$, while all terms $x\in \supp(A_{\mathfrak X})$ with depth greater than $r$ either have $\xi_{r,i}(x)=0$ or $\xi_{r,i}(x)\leq -1$ (as they must be integer values), for $i\in [1,t_r]$.
Thus, for each $i\in [1,t_r]$, we have $$|D_{r-1}^{(r,i)}|=\xi_{r,i}(\sigma(D_{r-1}^{(r,i)}))\geq \xi_{r,i}(\sigma(D_{r-1})).$$
For $i\in [1,t_r]$, let $B_i\mid  A_\mathfrak X\bdot D_{r-1}^{[-1]}$ be a minimal length sequence $B_i\in \Fc(G_0)$ such that either $\xi_{r,i}(\sigma(B_i))\leq -\xi_{r,i}(\sigma(D_{r-1}))$ (if such $B_i$ exists) or else let $B_i\mid A_\mathfrak X\bdot D_{r-1}^{[-1]}$ consist of all terms $x$ with $\xi_{r,i}(x)<0$.
Note, if $\xi_{r,i}(\sigma(D_{r-1}))\leq 0$, then $B_i$ is the trivial sequence.

Regardless of which case holds in the definition of
$B_i$, we have $\xi_{r,i}(x)\in \Z$ with $\xi_{r,i}(x)\leq -1$ for all $x\in \supp(B_i)$, ensuring $B_i$ only contains terms with depth at least $r+1$ and that $|B_i|\leq |\xi_{r,i}(\sigma(B_i))|=-\xi_{r,i}(\sigma(B_i))$. If the latter case holds in the definition of $B_i$ (but not the former), then $\xi_{r,i}(\sigma(B_i))>-\xi_{r,i}(\sigma(D_{r-1}))$, implying $|B_i|\leq -\xi_{r,i}(\sigma(B_i))\leq  \xi_{r,i}(\sigma(D_{r-1}))\leq |D_{r-1}^{(r,i)}|$. If the former case holds in the definition of $B_i$ and $\xi_{r,i}(\sigma(D_{r-1}))\geq 0$, then,  we  have $|B_i|\leq |\xi_{r,i}(\sigma(D_{r-1}))|=\xi_{r,i}(\sigma(D_{r-1}))\leq |D_{r-1}^{(r,i)}|$ in view of $\xi_{r,i}(x)\leq -1$ for all $x\in \supp(B_i)$ and the minimality of $|B_i|$.
Finally, if the former case holds in the definition of $B_i$ and $\xi_{r,i}(\sigma(D_{r-1}))\leq 0$, then $B_i$ is trivial, whence $|B_i|\leq |D_{r-1}^{(r,i)}|$ holds trivially.
Consequently, in all cases, $|B_i|\leq |D_{r-1}^{(r,i)}|$ for $i\in [1,t_r]$.
Thus, letting $B=\lcm(B_1,\ldots,B_{t_{r}})\in \Fc(G_0)$ (which is the smallest subsequence of $A_\mathfrak X\bdot D_{r-1}^{[-1]}$ containing each $B_i$ as a subsequence), we have \be\label{goatchease}|B|\leq |B_1|+\ldots+|B_{t_{r}}|\leq \Sum{i=1}{t_r}|D_{r-1}^{(r,i)}|=|D_{r-1}^{(r)}|.\ee
Set $D_r=D_{r-1}\bdot B$. Since $B\in \Fc(G_0)$ and $A_\mathfrak X\bdot D_{r-1}^{[-1]}\in \Fc(\mathfrak X)$ (in view of (a) for $D_{r-1}$), it follows that $A_\mathfrak X\bdot D_r^{[-1]}\in \Fc(\mathfrak X)$, whence (a) holds for $D_r$. For $i\in [1,t_r]$, we either have $\xi_{r,i}(\sigma(D_r))=
\xi_{r,i}(\sigma(D_{r-1}))+\xi_{r,i}(\sigma(B))\leq \xi_{r,i}(\sigma(D_{r-1}))+\xi_{r,i}(\sigma(B_i))\leq 0$ (with the first inequality in view of all terms  $x\in\supp(B)$ having depth at least $r+1$, which ensures $\xi_{r,i}(x)\leq 0$), or else $B_i$ contains all terms $x$ with $\xi_{r,i}(x)<0$, in which case $
\xi_{r,i}(\sigma( A_\mathfrak X\bdot D_r^{[-1]}))\geq 0$. However, we can also obtain this latter inequality in the former case by noting that \eqref{corun} combined with $\xi_{r,i}(\sigma(D_r))\leq 0$ implies  $$\xi_{r,i}(\sigma( A_\mathfrak X\bdot D_r^{[-1]}))=\xi_{r,i}(\sigma(A_\mathfrak X))-\xi_{r,i}(\sigma(D_r))
\geq \xi_{r,i}(\sigma(A_\mathfrak X))\geq 0.$$
By the comments after \eqref{explanation}, $\xi_{r,i}(x)\in \Z$ for each $x\in \supp(A_\mathfrak X\bdot D_r^{[-1]})$, whence $$\xi_{r,i}(\sigma(A_\mathfrak X\bdot D_r^{[-1]}))=\Summ{x\in \supp(A_\mathfrak X\bdot D_r^{[-1]})}\vp_x(A_\mathfrak X\bdot D_r^{[-1]})\xi_{r,i}(x)\in \Z,$$ with the final inclusion since (a) for $D_r$ holds.
Thus
$\xi_{r,i}(\sigma( A_\mathfrak X\bdot D_r^{[-1]}))\in \Z_+$ for $i\in [1,t_r]$. For $k<r$ and $i\in [1,t_k]$, since $B\in \mathcal F(\mathfrak X)$ only contains terms with depth at least $r+1>k$, we have $-\xi_{k,i}(\sigma(B))\in \Z_+$, which combined with (b) holding for $D_{r-1}$ implies $\xi_{k,i}(\sigma( A_\mathfrak X\bdot D_r^{[-1]}))=\xi_{k,i}(\sigma( A_\mathfrak X\bdot D_{r-1}^{[-1]}))-\xi_{k,i}(\sigma(B))\in \Z_+$. Thus (b) holds for $D_r$. As $B$ contains only terms with depth at least $r+1$, it follows that $D_r^{(j)}=D_{r-1}^{(j)}$ for all $j\in [1,r]$, whence (c) holds for $D_r$ (as (c) held for $D_{r-1}$).
We have $|D_r|=|D_{r-1}|+|B|\leq N_S+(|D_{0}^{(1)}|+|D^{(2)}_1|+\ldots+|D^{(r-1)}_{r-2}|)+|D^{(r)}_{r-1}|$ in view of \eqref{goatchease} and (d) for $D_{r-1}$, whence (d) holds for $D_r$. Finally, \ber \nn |D_{r-1}^{(r)}|&\leq& |D_{r-1}|-(|D_{r-1}^{(1)}|+|D_{r-1}^{(2)}|+\ldots+|D_{r-1}^{(r-1)}|)\\\nn &=&
|D_{r-1}|-(|D_{0}^{(1)}|+|D_{1}^{(2)}|+\ldots+|D_{r-2}^{(r-1)}|)\\\nn &\leq &(N_S+|D^{(1)}_0|+\ldots+|D^{(r-1)}_{r-2}|)-
(|D_{0}^{(1)}|+|D_{1}^{(2)}|+\ldots+|D_{r-2}^{(r-1)}|)=
N_S,\eer
with the first equality  in view of (c) and the second inequality in view of (d), both for $D_{r-1}$. Thus (e) holds for $D_r$ as well, completing the construction, which shows that the $D_{r}$ exist.

\bigskip

Set $S_\varphi=A_\mathfrak X\bdot D_{s-1}^{[-1]}$. Since $D_{s-1}\mid A_\mathfrak X$, it follows in view of (a) that $S_\varphi\in \Fc(\mathfrak X)$. In view of (b), we have $\xi_{k,i}(\sigma(S_\varphi))\in \Z_+$ for all $k\in [1,s-1]$ and $i\in [1,t_k]$. However, since $\xi_{s,i}(x)\in \Z_+$ for all $x\in \mathfrak X$ (per the comments after \eqref{explanation}), we trivially have $\xi_{s,i}(\sigma(S_\varphi))\in \Z_+$ for all $i\in [1,t_s]$. Thus $\sigma(S_\varphi)\in \C_\Z(Z_\varphi)$. Set $R=U\bdot \Big(\prod_{\varphi\in \mathfrak X_m(G_0)}S_\varphi\Big)^{[-1]}$. Since $U\in \Fc(G_0)$ and every $S_\varphi\in \Fc(G_0)$, it follows that $R\in \Fc(G_0)$. It remains only to bound $|R|$ independent of $U$. To this end, we can combine (d) and (e) (applied with $r=s-1$) to obtain the estimate
$|D_{s-1}|\leq s N_S\leq d N_S$.
We also have $|A_0|\leq N_T$ and $|A_{G_0\setminus \mathfrak X}|\leq N_T$ by Theorem \ref{thm-structural-Lambert}(c), while $|\mathfrak T_m(G_0)|<\infty$ in view of Proposition \ref{prop-finitary-mintype}. Thus $|R|\leq |A_0|+\sum_{\varphi}(|A_{G_0\setminus \mathfrak X}|+|D_{s-1}|)\leq N$ with $N=N_T+(N_T+d N_S)(|\mathfrak T_m(G_0)|-1)<\infty$ (with the summation over all nontrivial minimal types $\varphi$), which is a finite bound independent of the atom $U$, completing the proof.
\end{proof}

The tame degree (see \cite{alfredbook}) is an invariant of factorization theory whose finiteness implies the finiteness of numerous other factorization invariants.
The following theorem shows that having finite elasticities ensures a weaker tameness property holds in  $\mathcal B(G_0)$, though one which is just sufficiently strong to  still deduce the desired finiteness for the other invariants. We let $\mathsf t_w(G_0)$ be the minimal integer $N\geq 1$ for which Theorem \ref{thm-tame-pseudo} holds, which we call the \textbf{weak tame degree} of $G_0$. The proof of Theorem \ref{thm-tame-pseudo} illustrates how Theorem \ref{thm-structural-char} can be used to simulate finite support in an argument.

\begin{theorem}\label{thm-tame-pseudo}
Let $\Lambda\leq \R^d$ be a full rank lattice, where $d\geq 0$,  and let $G_0\subseteq \Lambda$ be a subset with  $\C(G_0)=\R^d$ and $0\notin \C^*(G_0^\diamond)$. Then there exists an integer $N\geq 1$ such that, given any  $U_1,\ldots,U_k,V_1,\ldots,V_\ell\in \mathcal A(G_0)$ with $$U_1\bdot\ldots\bdot U_k=V_1\bdot\ldots\bdot V_\ell,$$ where $k,\,\ell\geq 1$,  there exist atoms $W_1,\ldots,W_{\ell'}\in \mathcal A(G_0)$,  $r\in [1,k]$ and $\mathcal I\subseteq [1,\ell']$ such that $$U_1\bdot\ldots\bdot U_k=W_1\bdot\ldots\bdot W_{\ell'},$$ $\ell'\geq \ell$  and $U_r\mid \prod_{x\in \mathcal I}^\bullet W_x$ with $|\mathcal I|\leq N$.
\end{theorem}

\begin{proof}
Since $0\notin \C^*(G_0^\diamond)$,  Theorem \ref{thm-rho-char} and  Proposition \ref{prop-lambert-easy} imply  there is an integer $N_\rho\geq 1$ such that $$\rho_k(G_0)\leq N_\rho k\quad\mbox{ for all $k\geq 1$}.$$ Since $0\notin\C^*(G_0^\diamond)$, it follows from Theorem \ref{thm-keylemmaII} that $G_0$ is finitary. In view of Proposition \ref{prop-finitary-mintype}, there are only a finite number of minimal types. Let $\varphi_1,\ldots,\varphi_m\in \mathfrak T_m(G_0)$ be the distinct nontrivial minimal types for $G_0$. For each $j\in [1,m]$, let $$Z_{\varphi_j}=Z_1^{(j)}\cup \ldots\cup Z_{s_j}^{(j)}$$ be the codomain of $\varphi_j$ with $s_j\leq d$ and $|Z_{\varphi_j}|\leq d$.

Consider  $U_1,\ldots,U_k,V_1,\ldots,V_\ell\in \mathcal A(G_0)$ with $$S:=U_1\bdot\ldots\bdot U_k=V_1\bdot\ldots\bdot V_\ell,$$ where $k,\,\ell\geq 1$. By definition of $\rho_k(G_0)$, we have \be\label{ellsmall} \ell\leq \rho_k(G_0)\leq N_\rho k.\ee Let $S=g_1\bdot\ldots\bdot g_{|S|}$ be an indexing of the terms of $S$.
Let $I_1\cup \ldots\cup I_k=[1,|S|]=J_1\cup \ldots\cup J_\ell$ be disjoint partitions such that $$S(I_i)=U_i\quad \und\quad S(J_j)=V_j\quad\mbox{ for all $i\in[1,k]$ and $j\in [1,\ell]$}.$$
Since $0\notin \C^*(G_0^\diamond)$ and $\C(G_0)=\R^d$, we can apply Theorem \ref{thm-structural-char} to each atom $V_i$ for $i\in [1,\ell]$. Let $N_R\geq 0$ be the global bound from Theorem \ref{thm-structural-char} (which we can assume is an integer) and let $V_i=R_i\bdot \prod_{j=1}^mS_i^{(j)}$, for $i\in [1,\ell]$, be the resulting factorization given by Theorem \ref{thm-structural-char}, with $S_i^{(j)}\in \Fc(G_0)$ corresponding to the minimal type $\varphi_j$.
Then $\sigma(S_i^{(j)})\in \C_\Z(\Z_{\varphi_j})$, meaning we can  apply Proposition \ref{prop-swapping} to each $S_i^{(j)}$.
Let $J_i^{(j)}\subseteq [1,|S|]$ be disjoint subsets such that $$S(J_i^{(0)})=R_i\quad\und \quad S(J_i^{(j)})=S_i^{(j)},\quad\mbox{ for $i\in [1,\ell]$ and $j\in [1,m]$}.$$
 Moreover, for each $i\in [1,\ell]$, $j\in [1,m]$ and $n\in [1,s_j]$, let $J_i^{(j,n)}\subseteq J_i^{(j)}$ be the subset of all $x\in J_i^{(j)}$ with $\varphi_j(x)\in  Z^{(j)}_n$ at depth $n$.
Then  $$J_i=\bigcup_{j=0}^mJ_i^{(j)}=J_i^{(0)}\cup \bigcup_{j=1}^m\bigcup_{n=1}^{s_j}J_i^{(j,n)}\quad\mbox{  for every $i\in [1,\ell]$.}$$
Let \begin{align*}X_0=\bigcup_{i=1}^\ell J_i^{(0)},\quad
\Omega_0=\{(0,g_x):\; x\in  X_0\}\quad\und\quad
\Omega_\diamond=\{(j,z):\;j\in [1,m],\,z\in Z_{\varphi_j}\}.
\end{align*}
Moreover, partition $$\Omega_\diamond=\Omega_1\cup \ldots\cup \Omega_d$$ such that $\Omega_n$ consists of all $(j,z)\in \Omega_\diamond$ with $z\in Z_n^{(j)}$ at depth $n$.

Since Theorem \ref{thm-structural-char} implies $|J_i^{(0)}|=|R_i|\leq N_R$ for all $i\in [1,\ell]$, and since $|Z_{\varphi_j}|\leq d$ for all $j\in [1,m]$, we have \be\label{finitestuff}|X_0|\leq \ell N_R\quad\und\quad |\Omega_\diamond|=\Sum{j=1}{m}|Z_{\varphi_j}|\leq md.\ee
We view $\Omega:=\Omega_0\cup \Omega_\diamond=\Omega_0\cup \Omega_1\cup \ldots\cup \Omega_d$ as the set of \emph{support types} for $S=V_1\bdot\ldots\bdot V_\ell$. A support type $\tau\in \Omega_n$ is said to be at depth $n$. Note, if $\tau=(j,z)$ with $j\geq 1$, then the depth of $\tau$ equals the depth of $z\in Z_{\varphi_j}$. For each $x\in [1,|S|]$, we have $x\in J_i^{(j)}$ for some unique $i\in [1,\ell]$ and $j\in [0,m]$, allowing us to define $\mathsf s(x)=(0,g_x)\in \Omega_0$ when $j=0$  and $\mathsf s(x)=(j,\varphi_j(g_x))\in \Omega_\diamond$ when  $j\geq 1$. For $I\subseteq [1,|S|]$, $\mathsf s(I)\in \Fc(\Omega)$ is a sequence of support types from $\Omega$.  We associate the depth of $\mathsf s(x)$ (defined above) as the depth of $x\in [1,|S|]$.

Let
$$\alpha=\min\left\{|X_0\cap I_i|+\Summ{\tau\in \Omega_\diamond}\Big(\frac{\ell \vp_\tau(\mathsf s(I_i))}{\vp_\tau\big(\mathsf s\big([1,|S]]\big)\big)}+1\Big):\;i\in [1,k]\right\}.$$ Technically, we exclude any terms $\tau\in \Omega_\diamond$ in the sum defining $\alpha$ with  $\vp_\tau\big(\mathsf s\big([1,|S]]\big)\big)=0$. Then
\ber\nn k\alpha&\leq& \Sum{i=1}{k}|X_0\cap I_i|+\Sum{i=1}{k}\Summ{\tau\in \Omega_\diamond}\frac{\ell \vp_\tau(\mathsf s(I_i))}{\vp_\tau\big(\mathsf s\big([1,|S]]\big)\big)}+k|\Omega_\diamond|
 \\\nn &=& |X_0|+\Summ{\tau\in \Omega_\diamond}\Sum{i=1}{k}\frac{\ell \vp_\tau(\mathsf s(I_i))}{\vp_\tau\big(\mathsf s\big([1,|S]]\big)\big)}+k|\Omega_\diamond|\leq |X_0|+(\ell+k) |\Omega_\diamond|\\
 &\leq & \ell(N_R+md)+kmd\leq kN_\rho(N_R+md)+kmd,\label{firstthere}\eer with the first inequality in \eqref{firstthere} in view of  \eqref{finitestuff}, and the second in view of \eqref{ellsmall}. Thus  $$\alpha\leq N:=N_\rho(N_R+md)+md ,$$ which is a global bound independent of the $U_i$ and $V_j$.

 Let $r\in [1,k]$ be an index attaining the minimum in the definition of $\alpha$. Then $|X_0\cap I_r|\leq \alpha\leq N$, ensuring that there is some subset $\mathcal I_0\subseteq [1,\ell]$ with \be\label{Y0go}X_0\cap I_r\subseteq \bigcup_{i\in \mathcal I_0}J_i^{(0)}\quad \und\quad |\mathcal I_0|\leq |X_0\cap I_r|\leq N.\ee Likewise, letting $$n_\tau=\left\lceil \frac{\vp_{\tau}(\mathsf s(I_r))}{\vp_\tau\big(\mathsf s \big([1,|S|]\big)\big)/\ell}\right\rceil< \frac{\ell \vp_\tau(\mathsf s(I_r))}{\vp_\tau\big(\mathsf s\big([1,|S]]\big)\big)}+1\leq \ell+1\quad\mbox{ for $\tau\in \Omega_\diamond $},$$ we have \be\label{waterpool}\Summ{\tau\in \Omega_\diamond }n_\tau\leq \alpha-|X_0\cap I_r|\leq N-|\mathcal I_0|.\ee We interpret $n_\tau=0$ when $\vp_\tau\big(\mathsf s \big([1,|S|]\big)\big)=0$.

 We now describe how the $W_1,\ldots,W_{\ell'}\in \mathcal B(G_0)$ can be constructed.
 The idea is as follows. An index set $I\subseteq [1,|S|]$ indexes a sequence $S(I)\mid S$, but it also indexes a sequence $\mathsf s(I)\in \Fc(\Omega)$, obtained by replacing each indexed term in the sequence $S(I)$ with its corresponding support type from $\Omega$, so $\mathsf s(I)=\prod^\bullet_{x\in I}\mathsf s(x)$.  When $I\subseteq [1,|S|]\setminus X_0$, we have $\mathsf s(I)\in \Fc(\Omega_\diamond)$ with $\Omega_\diamond$ a fixed, finite set independent of $S$. Let $\tau\in \Omega_\diamond$.  If we select a subset  $\mathcal I_\tau\subseteq [1,\ell]$ with $|\mathcal I_\tau|=n_\tau$ such that the $\vp_\tau(\mathsf s(J_x))$, for $x\in I_\tau$, are the $n_\tau$ largest values occurring over all  $\vp_{\tau}(\mathsf s(J_i))$ with $i\in [1,\ell]$, then the definition of $n_\tau$ ensures that $$\vp_\tau\Big(\mathsf s\Big({\bigcup}_{z\in \mathcal I_\tau}J_z\Big)\Big)\geq
 n_\tau\Big(\vp_\tau\big(\mathsf s \big([1,|S|]\big)\big)/\ell\Big)\geq
 \vp_\tau(\mathsf s(I_r)),$$ with the first inequality holding since the sum of the $n_\tau$ largest terms in a sum of $\ell$ non-negative terms is always at least $n_\tau$ times the average value of all terms being summed.
 As a result, $\mathsf s(I_r)\mid \prod^\bullet_{z\in \mathcal I}\mathsf s(J_z)$, where $\mathcal I=\mathcal I_0\cup \bigcup_{\tau\in \Omega_\diamond}\mathcal I_\tau$, with $|\mathcal I|\leq N$ in view of \eqref{waterpool}.
 However, since the map $\mathsf s$ is not injective, this does not guarantee that the associated sequence $U_r=S(I_r)$ is a subsequence of the associated sequence $\prod^\bullet_{z\in \mathcal I}V_z=\prod^\bullet_{z\in \mathcal I}S(J_z)$.
 We do, however, have $S(X_0\cap I_r)\mid \prod^\bullet_{z\in \mathcal I_0} S(J_z)=\prod^\bullet_{z\in \mathcal I_0}V_z$ in view of \eqref{Y0go}.
 To deal with the terms from $\Omega_\diamond$, we must  use the sequences $S(T_x)$ given by Proposition \ref{prop-swapping} to exchange terms between the $V_i$.

 If there are terms $x\in J_i^{(j)}$ and $y\in J_{i'}^{(j)}$ with $\mathsf s(x)=\mathsf s(y)$, $i\neq i'$ and $j\geq 1$, say with $g_x$ and $g_y$ at depth $n$, then Proposition  \ref{prop-swapping} implies that there are subsets $T_x\subseteq J_i^{(j)}$ and $T_y\subseteq J_{i'}^{(j)}$ such that $\sigma(S(T_x))=\varphi_j(g_x)=\varphi_j(g_y)=\sigma(S(T_y))$. If we exchange these sets, defining $$K_i^{(j)}=(J_i^{(j)}\setminus T_x)\cup T_y\quad\und\quad K_{i'}^{(j)}=(J_{i'}^{(j)}\setminus T_y)\cup T_x,$$ and correspondingly define $$K_i=(J_i\setminus T_x)\cup T_y\quad\und\quad K_{i'}=(J_{i'}\setminus T_y)\cup T_x,$$
 then the new sequences $W_i=S(K_i)$ and $W_{i'}=S(K_{i'})$ (replacing $V_i$ and $V_{i'}$)  will still be zero-sum, though we do not guarantee that they remain atoms.
  However, since $y\in K_i$ and $x\in K_{i'}$, they are non-empty. Consequently, if either $W_i$ or $W_{i'}$ is not an atom, then we can re-factor them to write $V_i\bdot V_{i'}=W_i\bdot W_{i'}=V'_1\bdot\ldots\bdot V'_\omega$ as a product of $\omega\geq 3$ atoms. This leads to a factorization $U_1\bdot\ldots\bdot U_k=V_1\bdot\ldots\bdot V_{\ell}\bdot V_i^{[-1]}\bdot V_{i'}^{[-1]}\bdot V'_1\bdot\ldots \bdot V'_\omega$ into $\ell'=\ell-2+\omega>\ell$ atoms. In this case, we begin from scratch using this factorization in place of the original one $U_1\bdot\ldots\bdot U_k=V_1\bdot\ldots\bdot V_\ell$. As $\ell'\leq |S|<\infty$, we cannot start from scratch endlessly, meaning eventually we will never encounter this problem, allowing us to  w.l.o.g. assume $W_i=S(K_i)$ and $W_{i'}=S(K_{i'})$ are always atoms (where $\ell'=\ell$ may have increased in size from the original $\ell$ given in the hypotheses).
 Furthermore, we still have $\sigma(S(K^{(j)}_i))=\sigma(S(J_i^{(j)}))\in \C(Z_{\varphi_j})$ and $\sigma(S(K^{(j)}_{i'}))=\sigma(S(J_{i'}^{(j)}))\in \C(Z_{\varphi_j})$, with the inclusions originating from our application of Theorem \ref{thm-structural-char} at the start of the proof. Thus Proposition \ref{prop-swapping} can still be applied to $K_i^{(j)}$ and $K_{i'}^{(j)}$ if we later wish to continue with further such swaps between these sets.
  Proposition \ref{prop-swapping} guarantees that the set $T_x$ contains \emph{no} terms with depth greater than $n$ (the depth of $x$), and has $x$ as the \emph{unique} $a\in T_x$  having depth equal to $n$. Likewise for $T_y$. Thus when swapping the sets $T_x$ and $T_y$, we leave unaffected all terms in $J_i^{(j)}$ and $J_{i'}^{(j)}$  with depth at least $n$, apart from the exchanging of $x$ for $y$.
 This ensures that terms previously swapped but at a higher or equal depth will remain unaffected by exchanging $T_x$ and $T_y$. The value of $\mathsf s(a)$ does not change whether regarding $a\in [1,|S|]$ with respect to the original factorization $S=V_1\bdot\ldots\bdot V_\ell$ or to the one obtained after swapping $T_x$ and $T_y$ and replacing $V_{i}$ and $V_{i'}$ by $W_i$ and $W_{i'}$.
 In particular, the value of $\vp_\tau\big(\mathsf s\big([1,|S|]\big)\big)$ remains unchanged, ensuring that the value of $\alpha$ is unaffected when replacing $J_i$ and $J_{i'}$ by $K_i$ and $K_{i'}$,  and that $r\in [1,k]$ remains an index attaining the minimum in the definition of $\alpha$ (note, the numerators in the definition of $\alpha$ depend on the $U_i$, not the $V_i$).
 Swapping the elements $x$ and $y$ in this fashion leaves all elements from $X_0$, as well as any $J_c^{(b)}$ with $b\neq j$, unaltered, and the sequences $W_i$ and $W_{i'}$ remain nontrivial, as $W_i$ contains $g_y$, and $W_{i'}$ contains $g_x$. Since, apart from $x$ and $y$, only terms with depth less than $n$ are affected by the swap, it follows that the sets $\mathcal I_{\tau'}$, corresponding to any type $\tau'\in \Omega_\diamond$ with depth at least $n$, still have the property that they index the $\vp_\tau(\mathsf s(J_z))$, for $z\in I_{\tau'}$, with  the $n_{\tau'}$ largest values occurring over all  $\vp_{\tau'}(\mathsf s(J_i))$ with $i\in [1,\ell]$.

With these observations in mind, we can now describe how the zero-sums $V_i$ must be modified. Begin with any type $\tau\in \Omega_\diamond$ having maximal available depth. Construct the subset $\mathcal I_\tau$ for the current factorization $S=V_1\bdot\ldots\bdot V_\ell$ as described above. Then $\vp_\tau
\big(\mathsf s(\bigcup_{z\in \mathcal I_\tau}J_z)\big)\geq \vp_\tau(\mathsf s(I_r))$. If $\bigcup_{z\in \mathcal I_\tau}J_z$ contains all elements from $I_r$ having type $\tau$, then nothing need be done, we discard $\tau$ from the list of available types from $\Omega_\diamond$,  we select the next available type from $\Omega_\diamond$ with maximal depth, and continue  once more. On the other hand, if there is some $x\in I_r$ with type $\tau$ not contained in $\bigcup_{z\in \mathcal I_\tau}J_z$, then $\vp_\tau
\big(\mathsf s(\bigcup_{z\in \mathcal I_\tau}J_z)\big)\geq \vp_\tau(\mathsf s(I_r))$ ensures that there must be some $y\in \bigcup_{z\in \mathcal I_\tau}J_z$ having type $\tau$ with $y\notin I_r$. In this case, use Proposition \ref{prop-swapping} to define the sequences $T_x$ and $T_y$, perform the swap of $T_x$ and $T_y$ described above, and redefine our factorization $V_1\bdot\ldots\bdot V_\ell$ by replacing $J_i^{(j)}$ and $J_{i'}^{(j)}$ by $K_i^{(j)}$ and $K_{i'}^{(j)}$, where $x\in J_i^{(j)}$ and $y\in J_{i'}^{(j)}$, and correspondingly replacing $V_i$ and $V_{i'}$ by $W_i$ and $W_{i'}$. To simplify notation, redefine  $V_i$, $V_{i'}$, $J_i^{(j)}$ and $J_{i'}^{(j)}$ accordingly so as to reflect the new current state that now has $x\in \bigcup_{z\in \mathcal I_\tau}J_z$ and $y\notin \bigcup_{z\in \mathcal I_\tau}J_z$. If we now have $\bigcup_{z\in \mathcal I_\tau}J_z$ containing all elements from $I_r$ having type $\tau$, then nothing need be done, we discard $\tau$ from the list of available types in $\Omega_\diamond$ and carry on as before. If this is not the case, we again find a new term $x'\in I_r$ with type $\tau$ not contained in $\bigcup_{z\in \mathcal I_\tau}J_z$, find a new term $y'\notin \bigcup_{z\in \mathcal I_\tau}J_z$ with type $\tau$ and swap $x'$ and $y'$ as before by use of Proposition \ref{prop-swapping}. Since the depth of $x$ and $x'$ are the same, we will not swap $x$ back out of $\bigcup_{z\in \mathcal I_\tau}J_z$ when doing so, nor indeed any other element from  $\bigcup_{z\in \mathcal I_\tau}J_z$ having type $\tau$.
Thus, iterating such a procedure, we will eventually obtain that $\bigcup_{z\in \mathcal I_\tau}J_z$ contains every element of $I_r$ having type $\tau$,
in which case we move on to the next available type $\tau'$ with maximal available depth. We repeat the same for procedure for $\tau'$ as we did for $\tau$ (and, later, as we did for any support types previously discarded before selecting $\tau'$).
We first construct the subset $\mathcal I_{\tau'}$ for  the current state for $S=V_1\bdot\ldots\bdot V_\ell$, and then swap elements into $\bigcup_{z\in \mathcal I_{\tau'}}J_z$ until it contains all elements from $I_r$ having type $\tau'$. While doing so, since we always first choose support types with maximal available depth, we are assured that any type $\tau''$ that has already been discarded had depth at least that of $\tau'$, and thus no  element of type $\tau''$ will be moved when swapping at the later stage for $\tau'$, ensuring that prior work cannot be undone. Continue until all support types from $\Omega_\diamond$ have been exhausted. Once the process ends, we now have a new factorization $U_1\bdot\ldots\bdot U_k=S=W_1\bdot\ldots\bdot W_\ell$, where $W_i$ reflects the final state of $V_i$ after running the above process, such that $U_r=S(I_r)\mid \prod_{z\in \mathcal I}^\bullet W_i$, where $\mathcal I=\mathcal I_0\cup\bigcup_{\tau\in \Omega_\diamond}\mathcal I_\tau$, with $|\mathcal I|\leq |X_0\cap I_r|+\Summ{\tau\in \Omega_\diamond}n_\tau\leq N$, completing the proof.
\end{proof}

With the aid of Theorem \ref{thm-tame-pseudo}, it is now possible to  establish that both the set of distances $\Delta(G_0)$ and Catenary degree $\mathsf c(G_0)$ are finite, and that there can be no arbitrarily large jumps in the elasticities $\rho_k(G_0)$, which in turn implies that  the Structure Theorem for Unions  holds in $\mathcal B(G_0)$.

\begin{theorem}\label{thm-Delta-finite}
Let $\Lambda\leq \R^d$ be a full rank lattice, where $d\geq 0$,  and let $G_0\subseteq \Lambda$ be a subset with   $\rho_{d+1}(G_0)<\infty$. Let $N=\max\{2,\mathsf t_w(G_0)\}$.
\begin{itemize}
\item[1.]  $\rho_k(G_0)-\rho_{k-1}(G_0)\leq N<\infty$ for all $k\geq 2$.
\item[2.] $\max \Delta(G_0)\leq \rho_N(G_0)-N<\infty$. In particular,  $\Delta(G_0)$ is finite.
\item[3.] The Structure Theorem for Unions  holds in $\mathcal B(G_0)$.
\item[4.] The catenary degree $\mathsf c(G_0)\leq \rho_N(G_0)<\infty$ is finite.
\end{itemize}
\end{theorem}

\begin{proof}
Removing elements from $G_0$ that are contained in no atom does not affect the quantities $\rho_k(G_0)$ nor $\Delta(G_0)$,  $\mathsf c(G_0)$, $\mathsf t_w(G_0)$ or $\mathcal B(G_0)$. Thus we may w.l.o.g. assume every $g\in G_0$ is contained in some atom, in which case Proposition \ref{prop-notrivialG_0} implies  $\C(G_0)=\R\la G_0\ra$ is a subspace, and $\Lambda\cap \R\la G_0\ra$ is a full rank lattice in $\R\la G_0\ra$. Thus we may w.l.o.g. assume $\C(G_0)=\R^d$. Then, in view of $\rho_{d+1}(G_0)<\infty$ and Theorem \ref{thm-rho-char}, we have $0\notin \C^*(G_0^\diamond)$ and  can apply Proposition \ref{prop-lambert-easy} to conclude that there is a constant $N_\rho\geq 1$ such that $$\rho_k(G_0)\leq N_\rho k<\infty\quad\mbox{ for all $k\geq 1$}.$$ Since $0\notin\C^*(G_0^\diamond)$ and $\C(G_0)=\R^d$, it follows from Theorem \ref{thm-keylemmaII} that $G_0$ is finitary.

1. Apply Theorem \ref{thm-tame-pseudo} to $G_0$ and let $N=\mathsf t_w(G_0)\in \Z_+$ be the resulting constant. Let $k\geq 2$ and let $U_1,\ldots,U_k,V_1,\ldots,V_\ell\in \mathcal A(G_0)$ be atoms with $$U_1\bdot\ldots\bdot U_k=V_1\bdot\ldots\bdot V_\ell\quad\und\quad \ell=\rho_k(G_0).$$ Then Theorem \ref{thm-tame-pseudo} implies that there are atoms $W_1,\ldots, W_{\ell'}\in \mathcal A(G_0)$, where $\ell'\geq \ell$, and $r\in [1,k]$ and $\mathcal I\subseteq [1,\ell']$ such that $$U_1\bdot\ldots\bdot U_k=W_1\bdot\ldots\bdot W_{\ell'}$$
and $U_r\mid \prod_{x\in\mathcal I}^\bullet W_x$ with $|\mathcal I|\leq N$. Since $\ell=\rho_k(G_0)$, the definition of $\rho_k(G_0)$ ensures $\ell'\leq \ell$, whence $\ell'=\ell'$. By re-indexing, we may w.l.o.g assume $r=k$ and $\mathcal I=[1,m]$ with $$m\leq N.$$  Since $U_k\mid \prod_{i\in [1,m]}^\bullet W_i$, there exists a factorization $U_k\bdot W'_2\bdot\ldots\bdot W_{m'}=W_1\bdot\ldots\bdot W_{m}$ with $W'_i\in \mathcal A(G_0)$ for all $i\in [2,m']$. Note $m'\geq 1$. But now $$U_1\bdot\ldots\bdot U_{k-1}=W'_2\bdot\ldots\bdot W'_{m'}\bdot W_{m+1}\bdot\ldots\bdot W_\ell$$ with $W'_i,W_j\in \mathcal A(G_0)$ for all $i$ and $j$. Thus
$$\rho_{k-1}(G_0)\geq \ell-m+m'-1\geq \ell-m=\rho_k(G_0)-m\geq \rho_k(G_0)-N,$$ implying  $\rho_k(G_0)-\rho_{k-1}(G_0)\leq N=\mathsf t_w(G_0)$, as desired.

2. To show $\Delta(G_0)$ is finite, it suffices to show $\max \Delta (G_0)$ is finite. Apply Theorem \ref{thm-tame-pseudo} to $G_0$ and let $N\geq 2$ be the resulting constant (if $N=1$, replace it with $N=2$). We will show that \be\label{Delta-max}\max \Delta(G_0)\leq \rho_N(G_0)-N.\ee

For $S\in \mathcal B(G_0)$, let $\mathsf L(S)\subseteq [1,|S|]$ be the set of lengths for $S$, which consists of all $k\in [1,|S|]$ for which there exists a factorization $S=U_1\bdot\ldots \bdot U_k$ with $U_i\in \mathcal A(G_0)$ for all $i$. Then $\delta\in \Delta(G_0)$ means there is some $S\in \mathcal B(G_0)$ with $\delta=\ell-k$ for two consecutive elements $k,\,\ell\in \mathsf L(S)$ with $k<\ell$. In other words,  there must exist some $S\in \mathcal B(G_0)$ and atoms $U_1,\ldots,U_k,V_1,\ldots,V_\ell\in \mathcal A(G_0)$ such that $$U_1\bdot\ldots\bdot U_k=S=V_1\bdot\ldots\bdot V_\ell$$ with $k<\ell$ and $\ell-k=\delta$, and there cannot exist a factorization $W_1\bdot\ldots\bdot W_r=S$ with $W_i\in \mathcal A(G_0)$ and $k<r<\ell$. Note $\Delta(G_0)=\bigcup_{S\in \mathcal B(G_0)}\Delta(\mathsf L(S))$, where $\Delta(\mathsf L(S))$ consists of all $\delta$ for which there exist consecutive elements $k,\,\ell\in \mathsf L(S)$ with $k<\ell$ and $\ell-k=\delta$.

If \eqref{Delta-max} fails, then there must be some $S\in \mathcal B(G_0)$ and $k,\,\ell\in \mathsf L(S)$ with $k<\ell$ consecutive elements of $\mathsf L(S)$ and \be\label{DeltaBig}\ell-k\geq \rho_N(G_0)-(N-1)\geq \rho_n(G_0)-(n-1)\quad \mbox{ for all $n\in [1,N]$},\ee where the second inequality in \eqref{DeltaBig} follows in view of iterated application of the basic inequality $\rho_{k+1}(G_0)>\rho_{k}(G_0)$ (see \eqref{rho-ascend-chain} in  Section \ref{sec-transferkrull}). Choose such a counter example with $\ell$ minimal. Note $k\geq 2$, else $\ell=k$, contradicting that $\ell-k\geq \rho_N(G_0)-N+1\geq 1$ by \eqref{DeltaBig}. Thus \eqref{DeltaBig} ensures $\ell\geq \rho_N(G_0)-N+1+k\geq k+1\geq 3$. Let \be\label{factorization}U_1\bdot\ldots\bdot U_k=S=V_1\bdot\ldots\bdot V_\ell\ee be factorizations exhibiting that $k,\,\ell\in \mathsf L(S)$, where $U_1,\ldots,U_k,V_1,\ldots,V_\ell\in \mathcal A(G_0)$.
Apply Theorem \ref{thm-tame-pseudo} to the factorization given in \eqref{factorization} with the roles of the $U_i$ and the $V_j$ swapped.
Then there are atoms  $W_1,\ldots,W_{k'}\in \mathcal A(G_0)$, where $k'\geq k$, and  $r\in [1,\ell]$ and $\mathcal I\subseteq [1,k']$ such that
$$W_1\bdot\ldots\bdot W_{k'}=V_1\bdot\ldots \bdot V_\ell$$ with $V_r\mid \prod^\bullet_{x\in \mathcal I}W_x$ and \be\label{matmostN}m:=|\mathcal I|\leq N.\ee By re-indexing, we may w.l.o.g. assume $r=\ell$ and $\mathcal I=[1,m]$.

The algorithm which constructs the sequences $W_i$ given in Theorem \ref{thm-tame-pseudo} proceeds by successively taking two zero-sums $U_i$ and $U_{i'}$ in the factorization $U_1\bdot\ldots\bdot U_k$ and replacing them with a re-factorization of  $U_i\bdot U_{i'}=W_i\bdot W_{i'}$ into $\omega\geq 2$ atoms. We begin will all $U_j\in \mathcal A(G_0)$ being atoms. If, at some point during the process of constructing the $W_i$, we take two atoms and find that their replacement zero-sum sequences $W_i$ and $W_{i'}$ are not both themselves atoms, that is, $\omega\geq 3$,  then the first time that this occurs, we can re-factor $S$ by replacing these two zero-sums with a factorization of length $\omega\in [3,\rho_2(G_0)]$, to thereby find that $S$ has a factorization of length $k-2+\omega$ with $k<k-2+\omega\leq k-2+\rho_2(G_0)<\ell$, where the last inequality follows from \eqref{DeltaBig} and $N\geq 2$. However, this contradicts that $k<\ell$ are \emph{consecutive} elements of $\mathsf L(S)$. Therefore, we instead conclude that we never need to re-factor the sequences $W_i\bdot W_{i'}$ when applying the algorithm for  Theorem \ref{thm-tame-pseudo}, whence $k'=k$ follows.

Since $V_\ell\mid  W_1\bdot\ldots\bdot W_m$, we have a factorization $V_\ell\bdot W'_2\bdot\ldots\bdot W'_{m'}=W_1\bdot\ldots\bdot W_m$ with $W'_i\in \mathcal A(G_0)$ for all $i\in [2,m']$.
Note \be\label{m'bound} 1\leq m'\leq \rho_m(G_0).\ee But now $$W'_2\bdot\ldots\bdot W'_{m'}\bdot W_{m+1}\bdot\ldots\bdot W_k=S\bdot V_\ell^{[-1]}=V_1\bdot\ldots\bdot V_{\ell-1},$$ showing that $$k'-1,\,\ell-1\in \mathsf L(S'),$$ where $$S'=S\bdot V_\ell^{[-1]}\in\mathcal B(G_0)\quad\und\quad k'=k-m+m'.$$
Observe that any factorization of length $t$ for $S'$ gives a factorization of $S$ of length $t+1$ by concatenating the atom $V_\ell$ onto the end of the factorization of $S'$. Consequently, since $k<\ell$ are consecutive elements of $\mathsf L(S)$, it follows that \be\label{intervalClear}[k,\ell-2]\cap \mathsf L(S')=\emptyset.\ee
%
In view of
 \eqref{m'bound}, \eqref{matmostN} and \eqref{DeltaBig},  we have $k'=k-m+m'\leq k+\rho_m(G_0)-m\leq \ell-1$. Thus, since $k'-1\in \mathsf L(S')$, we conclude from \eqref{intervalClear} that  $k'-1\leq k-1$. Hence, in view of \eqref{intervalClear} again, let $r\in [k'-1,k-1]$ be the maximal element of $\mathsf L(S')$ less than $\ell-1\in \mathsf L(S')$.
 Since $r,\,\ell-1\in \mathsf L(S')$ with $r\leq k-1<\ell-1$,   the minimality of $\ell$ ensures that $$\rho_N(G_0)-N\geq (\ell-1)-r\geq (\ell-1)-(k-1)=\ell-k\geq \rho_N(G_0)-N+1,$$ with the final inequality in view of \eqref{DeltaBig}, which is a contradiction. Thus \eqref{Delta-max} is now established, and since $\rho_N(G_0)<\infty$, we conclude that $\max \Delta(G_0)$, and thus also $\Delta(G_0)$ itself, are both finite, completing Item 2.

3. This follows from Items 1 and 2 and \cite[Theorem 4.2]{Gao-Ger-phok}.

4. Let $N\geq 2$ be a constant given by Theorem \ref{thm-tame-pseudo} applied to $G_0$. We will show that $$\mathsf c(G_0)\leq \rho_N(G_0)<\infty.$$ Let $U_1,\ldots,U_k,V_1,\ldots,V_\ell\in \mathcal A(G_0)$ be arbitrary atoms with $$U_1\bdot \ldots\bdot U_k=V_1\bdot\ldots\bdot V_\ell.$$
By re-indexing the $U_i$ and $V_i$, we may collect together all the atoms which occur in both factorizations, say let  $I\subseteq [1,k]\cap [1,\ell]$ consist of all $i$ such that $U_i=V_i$.
Apply Theorem \ref{thm-tame-pseudo} to the factorization $\prod_{i\in [1,k]\setminus I}^\bullet U_i=\prod_{i\in [1,\ell]\setminus I}^\bullet V_i$.
Then there exists a factorization $\prod_{i\in [1,\ell']\setminus I}^\bullet W_i=\prod_{i\in [1,k]\setminus I}^\bullet U_i$ with $\ell'\geq \ell$ and $W_i\in \mathcal A(G_0)$ satisfying the conclusion of Theorem  \ref{thm-tame-pseudo}. Moreover the algorithm from the proof of Theorem \ref{thm-tame-pseudo} constructs the $W_i$ by sequentially modifying pairs of atoms $W_i$ and $W_{i'}$ by replacing such pairs with a re-factorization into $\omega\in [2,\rho_2(G_0)]$ new atoms. Thus each successive factorization in the algorithm constructing the $W_i$ differs from the prior one by at most $\rho_2(G_0)$ factors.
Theorem \ref{thm-tame-pseudo} guarantees that there is some $r\in [1,k]\setminus I$ with $U_r$ a subsequence of a product of at most $N$ of the atoms $W_i$, say $U_r\mid \prod_{i\in J}^\bullet W_i$ with $|J|\leq N$ and $J\subseteq [1,\ell']\setminus I$. We may re-factor $\prod_{i\in J}^\bullet W_i$ into a product of at most $\rho_N(G_0)$ atoms that includes the atom $U_r$. We thus obtain a factorization $U_1\bdot\ldots\bdot U_k=V'_1\bdot\ldots\bdot V'_{\ell''}$ that now has one more shared factor among the $U_i$ and $V'_i$ than the original one, and which can be constructed sequentially with each new factorization differing from the prior one by at most $\rho_N(G_0)\geq \rho_2(G_0)\geq 2$ factors. Iterating this process at most $k$ times thus transforms the factorization $V_1\bdot \ldots\bdot V_\ell$ into the factorization $U_1\bdot\ldots\bdot U_k$ with each successive factorization differing by at most $\rho_N(G_0)$ factors from the previous one, yielding the desired bound for the catenary degree.
\end{proof}

\subsection*{Summary} We can now summarize our results regarding what we have shown regarding $\mathcal B(G_0)$ under an assumption of finite elasticities.
Let $\Lambda\leq \R^d$ be a full rank lattice, where $d\geq 0$,  and let $G_0\subseteq \Lambda$ be a subset. There is little loss of generality to assume every element $g\in G_0$ occurs in some atom (else we can pass to the subset of $G_0$ having this property). Then Proposition \ref{prop-notrivialG_0} implies that $\C(G_0)=\R\la G_0\ra$ with $\Lambda\cap \R\la G_0\ra\leq \R\la G_0\ra$ a full rank lattice. Hence, replacing $\R^d$ with $\R\la G_0\ra$, we can w.l.o.g. assume $\C(G_0)=\R^d$, which is simply a normalization hypothesis to avoid trivial degeneracies. Under these assumptions, we now summarize some of the key results.
\begin{itemize}
\item[1.] There exists a minimal $s\in [1,d+1]$ such that $\rho_s(G_0)<\infty$ implies $\rho_k(G_0)<\infty$ for all $k\geq 1$.
\end{itemize}
Item 1 follows in view of Theorem \ref{thm-rho-char}, which also shows that $$\rho_{d+1}(G_0)<\infty\mbox{ is equivalent to } \rho(G_0)<\infty.$$ We remark that it would be interesting to know whether the estimate $s\leq d+1$ is tight or can be improved; we have focussed primarily on its existence.

Corollary \ref{cor-rhoelem} ensures that $\rho_{d+1}^{\mathsf{elm}}(G_0)<\infty$ implies $\rho_{d+1}(G_0)<\infty$, meaning it is sufficient to know no product of $d+1$ elementary atoms can be re-factored into an arbitrarily large number of atoms. Theorem \ref{thm-rho-char} characterizes when Item 1 occurs either in  terms of  a basic combinatorial property of the atoms $\mathcal A(G_0)$, or the geometric property $0\notin \C^*(G_0^\diamond)$, which involves $\R_+$-linear combinations of elements of $G_0$ rather than $\Z_+$-linear combinations.
Assuming additionally that $\rho_{d+1}(G_0)<\infty$, so that the conclusion of  Item 1 holds, we obtain the following properties for  $\mathcal B(G_0)$ and subset $G_0\subseteq \Lambda$.
\begin{itemize}
\item[2.] $G_0\subseteq \R^d$ is finitary (by Theorem \ref{thm-keylemmaII}). In particular, all the results of Section \ref{sec-finitary} are available for studying the set $G_0$, including Theorems \ref{thm-finitary-diamond-containment}, \ref{thm-finitary-FiniteProps-I}, \ref{thm-finitary-FiniteProps-II} and \ref{thm-finitary-FiniteProps-III}.

\item[3.] Besides the three equivalent defining definitions of the subset $G_0^\diamond\subseteq G_0$ (given in Proposition \ref{prop-G_0diamond-1st-easy-equiv}), we also have (by Corollary \ref{cor-G_0diamond-Zequiv})
\ber\nn G_0^\diamond&=&\{g\in G_0:\; \sup\{\vp_{g}(U):\; U\in \mathcal A^{\mathsf{elm}}(G_0)\}=\infty\}\\\nn &=&\{g\in G_0:\; \sup\{\vp_{g}(U):\; U\in \mathcal A(G_0)\}=\infty\}.\eer
  \item[4.] There is a finite subset $X\subseteq G_0^\diamond$ such that $\mathcal A(G_0\setminus X)$ is finite, and thus also a finite subset $Y\subseteq G_0$ such that $\mathcal A(G_0\setminus Y)=\emptyset$ (by Proposition \ref{prop-finitary-FiniteDeletion}).
  \item[5.] $G_0\setminus G_0^\diamond \subseteq G_0$ is a Lambert subset, indeed, the unique maximal Lambert subset. (Corollaries \ref{cor-structural-Lambert} and \ref{cor-G_0diamond-Zequiv})
  \item[6.] The Weak Tame Degree (as defined in Theorem \ref{thm-tame-pseudo}) is finite: $\mathsf t_w(G_0)<\infty$, and we set $N:=\max\{2,\mathsf t_w(G_0)\}$.
  \item[7.] The elasticities of $G_0$ do not contain arbitrarily large gaps: $\rho_k(G_0)-\rho_{k-1}(G_0)\leq N<\infty$ for all $k\geq 2$. In particular, $\rho_k(G_0)\leq Nk$ grows linearly. (Theorem \ref{thm-Delta-finite}.1)
  \item[8.] $\max \Delta(G_0)\leq \rho_N(G_0)-N<\infty$ is finite. (Theorem \ref{thm-Delta-finite}.2)
  \item[9.] The Set of Distances is  finite: $|\Delta(G_0)|<\infty$. (Theorem \ref{thm-Delta-finite}.2)
  \item[10.] The Catenary Degree is finite: $\mathsf c(G_0)\leq \rho_N(G_0)<\infty$. (Theorem \ref{thm-Delta-finite}.4)
  \item[11.] The Structure Theorem for Unions holds for $\mathcal B(G_0)$. (Theorem \ref{thm-Delta-finite}.3)
  \item[12.] A weak structure theorem holds for the atoms, effectively allowing simulation of globally bounded finite support for the atoms $U\in \mathcal A(G_0)$. (Theorem \ref{thm-structural-char})
\end{itemize}

\subsection{Subsets of Finitely Generated Abelian Groups}\label{sec-transfer-monoid-results}

In this final subsection, we show how our prior results regarding subsets $G_0\subseteq \Lambda \leq \R^d$ can be extended to cover (Transfer) Krull Monoids $H$ over a subset $G_0\subseteq G$, where $G=\Z^d\oplus G_T$ is a finitely generated abelian group with torsion-free rank $d\geq 0$ and torsion subgroup $G_T\leq G$.

When $G=G_T$ is a finite group, all factorization invariants we have encountered are trivially finite. It is natural to suppose their finiteness, for subset of a more general finitely generated abelian group $G$, is thus principally affected by the torsion-free portion of $G$, which naturally embeds as a lattice $\Z^d\leq \R^d$. To this end, we have the following basic relation between $\mathcal A(G_0)$ and $\mathcal A(\pi(G_0))$, where $\pi:G\rightarrow \Z^d$ is the natural projection homomorphism with kernel $G_T$.

\begin{proposition}\label{prop-fg-atoms-basic-transfer}
Let $G=\Z^d\oplus G_T$ be a finitely generated abelian group with torsion subgroup $G_T\leq G$, where $d\geq 0$, let $\pi:G\rightarrow \Z^d$ be the projection homomorphism with $\ker \pi=G_T$,  let $\mathsf D(G_T)$ be the Davenport constant for $G_T$, let $m=\exp(G_T)$, and let $G_0\subseteq G$.
\begin{itemize}
\item[1.] If $U\in \mathcal A(G_0)$ is an atom and $\pi(U)=W_1\bdot\ldots\bdot W_\ell$ is a factorization of $\pi(U)$ with $W_i\in \mathcal A(\pi(G_0))$ for all $i\in [1,\ell]$, then $\ell\leq \D(G_T)$.
\item[2.]  If $U\in \Fc(G_0)$ with $\pi(U)\in \mathcal A^{\mathsf{elm}}(\pi(G_0))$, then $U^{[m]}\in \mathcal B(G_0)$. Moreover, if  $U^{[m]}=V_1\bdot\ldots\bdot V_\ell$ is a factorization  with $V_i\in \mathcal A(G_0)$ for all $i\in [1,\ell]$, then $\ell \leq m=\exp(G_T)$.
\end{itemize}
\end{proposition}

\begin{proof}
1. Factor $U=V_1\bdot\ldots\bdot V_\ell$, with the $V_i\in \Fc(G_0)$ such that $\pi(V_i)=W_i$ for all $i\in [1,\ell]$.
Since  $U\in \mathcal A(G_0)$ is an atom,  $\pi(U)\in \mathcal B(\pi(G_0))$ is a zero-sum sequence. Since $\pi(V_i)=W_i\in \mathcal A(\pi(G_0))$, we have   $\sigma(V_i)\in \ker \pi=G_T$ for each $i\in [1,\ell]$. If $\ell>\mathsf D(G_T)$, then applying the definition of $\D(G_T)$ to the sequence $\prod_{i\in [1,\ell]}^\bullet \sigma(V_i)\in \Fc(G_T)$, we can find some nontrivial, proper subset $I\subseteq [1,\ell]$ such that $\prod_{i\in I}^\bullet V_i\in \Fc(G_0)$ is a nontrivial, proper zero-sum subsequence of $U$, contradicting that $U\in \mathcal A(G_0)$ is an atom, which completes Item 1.

2. Since  $\pi(U)\in \mathcal A^{\mathsf{elm}}(G_0)$ is an elementary atom, we have $\sigma(U)\in \ker \pi=G_T$, whence $\sigma(U^{[m]})=m\sigma(U)=\exp(G_T)\sigma(U)=0$, ensuring $U^{[m]}\in \mathcal B(G_0)$.
If $\pi(U)$ is the subsequence consisting a single term equal to $0$, then $|U^{[m]}|=m$, in which case $\ell\leq m$ is trivial. Therefore we may instead assume $0\notin \supp(\pi(U))$. As a result,
since  $\pi(U)\in \mathcal A^{\mathsf{elm}}(G_0)$ is an elementary atom, it follows by Proposition \ref{prop-char-minimal-pos-basis} that
 $X=\supp(\pi(U))$ is a minimal positive basis, and $\pi(U)$ is the unique atom whose support is contained in $X$. Since each $V_i\in \mathcal A(G_0)$ with $V_i\mid U^{[m]}$, it follows that $\pi(V_i)\in \mathcal B(X)$. However, as already remarked, the unique atom with support contained in $X$ is $\pi(U)$, ensuring that each $\pi(V_i)=\pi(U)^{[m_i]}$ for some $m_i\geq 1$. But now
 $$\pi(U)^{[m]}=\pi(U^{[m]})=\pi(V_1)\bdot\ldots\bdot \pi(V_\ell)=\pi(U)^{[m_1]}\bdot\ldots\bdot \pi(U)^{[m_\ell]},$$ implying $\ell\leq \Sum{i=1}{\ell}m_i=m$, with the inequality since  $m_i\geq 1$ for all $i$, completing  the proof.
\end{proof}

The next proposition  gives a correspondence between finite elasticities in $G_0\subseteq G$ and finite elasticities in $\pi(G_0)\subseteq \Z^d$. In particular, Theorem \ref{thm-rho-char} also characterizes when $\rho_{(d+1)m}(G_0)<\infty$ is finite by applying it to the set $\pi(G_0)\subseteq \Z^d$.

\begin{proposition}\label{prop-fg-rho-transfer}
Let $G=\Z^d\oplus G_T$ be a finitely generated abelian group with torsion subgroup $G_T\leq G$, where $d\geq 0$, let $\pi:G\rightarrow \Z^d$ be the projection homomorphism with $\ker \pi=G_T$,  let $\mathsf D(G_T)$ be the Davenport constant for $G_T$, let $m=\exp(G_T)$, and let $G_0\subseteq G$.
\begin{itemize}
\item[1.] If $\rho_{d+1}(\pi(G_0))<\infty$, then $\rho_k(G_0)\leq \rho_{k\mathsf D(G_T)}(\pi(G_0))<\infty$  for all $k\geq 1$.
\item[2.] If $\rho_{(d+1)m}(G_0)<\infty$,   then $\rho_k(\pi(G_0))<\infty$ for all $k\geq 1$.
\end{itemize}
In particular, $\rho_k(G_0)<\infty$ for all $k\geq 1$ if and only if $\rho_k(\pi(G_0))<\infty$ for all $k\geq 1$. Moreover, if $\rho_{(d+1)m}(G_0)<\infty$, then $\rho_k(G_0)<\infty$ for all $k\geq 1$.
\end{proposition}

\begin{proof}
We may w.l.o.g. discard terms from $G_0$ contained in no atom, as such terms have no bearing on $\rho_k(G_0)$ nor $\rho_k(\pi(G_0))$ (as $G_T$ has finite exponent), and thereby assume every $g\in G_0$ is contained in some atom. We may also w.l.o.g. assume $\la G_0\ra=G$ in view of \eqref{rho-ascend-chain}, and may embed $\Z^d\leq \R^d$, which is a full rank lattice in $\R^d$. Since $\la G_0\ra=G$, we have $\R\la \pi(G_0)\ra=\R^d$. Since every  $g\in G_0$ is contained in some atom, it follows that every $\pi(g)\in \pi(G_0)$ is also contained in some atom, whence Proposition \ref{prop-notrivialG_0} implies  $\C(\pi(G_0))=\R\la \pi(G_0)\ra= \R^d$ with $\pi(G_0)\subseteq \Z^d$.

1. Suppose $\rho_{d+1}(\pi(G_0))<\infty$. Then Theorem \ref{thm-rho-char} implies that   \be\nn\rho_k(\pi(G_0))<\infty\quad\mbox{ for all $k\geq 1$}.\ee Consider an arbitrary  factorization $$U_1\bdot\ldots\bdot U_k=V_1\bdot\ldots\bdot V_{\ell}$$ with  $U_j,\,V_i\in \mathcal A(G_0)$ for all $i$ and   $j$.
Re-factor each $\pi(U_j)=\prod^\bullet_{x\in I_{j}}W^{(j)}_x$  into a product of $|I_{j}|\leq \mathsf D(G_T)$ atoms $W^{(j)}_x\in \mathcal A(\pi(G_0))$ via Proposition \ref{prop-fg-atoms-basic-transfer}.1.  Then $\prod_{j\in [1,k]}^\bullet \prod^\bullet_{x\in I_{j}}W^{(j)}_x$ is a product of at most $k\D(G_T)$ atoms. As a result, since $${\prod}_{j\in [1,k]}^\bullet {\prod}^\bullet_{x\in I_{j}}W_x^{(j)}=\pi(U_1)\bdot\ldots\bdot \pi(U_k)=\pi(V_1)\bdot \ldots\bdot \pi(V_{\ell})$$ with each $\pi(V_j)\in \mathcal B(\pi(G_0))$ nontrivial, it follows in view of \eqref{rho-ascend-chain} that $\ell\leq \rho_{k\mathsf D( G_T)}(\pi(G_0))<\infty$, as desired.

2. Suppose $\rho_{(d+1)m}(G_0)<\infty$ but, by way of contradiction, that there is some $k\geq 1$ such that $\rho_k(\pi(G_0))=\infty$. If $\rho_{d+1}(\pi(G_0))<\infty$, then Theorem \ref{thm-rho-char} ensures that $\rho_k(\pi(G_0))<\infty$, contrary to assumption. Therefore we can assume $\rho_{d+1}(\pi(G_0))=\infty$. Thus Corollary \ref{cor-rhoelem} implies that $\rho^{\mathsf{elm}}_{d+1}(\pi(G_0))=\infty$. Consequently, there is a sequence of factorizations
$$\pi(U_{i,1})\bdot\ldots\bdot \pi(U_{i,d+1})=\pi(V_{i,1})\bdot\ldots\bdot \pi(V_{i,\ell_i})\quad\mbox{for $i=1,2,\ldots$}$$ such that
$$U_{i,j},\,V_{i,j}\in\Fc(G_0), \quad \pi(U_{i,j})\in \mathcal A^{\mathsf{elm}}(\pi(G_0)),\quad \und\quad \pi(V_{i,j})\in\mathcal A(\pi(G_0)),\quad\mbox{ for all $i$ and $j$,}$$ with \be\label{equal}U_{i,1}\bdot\ldots\bdot U_{i,d+1}=V_{i,1}\bdot\ldots\bdot V_{i,\ell_i}\quad\mbox{ for every $i\geq 1$}, \quad\und \quad \ell_i\rightarrow \infty.\ee Since the $U_{i,j}$ and $V_{i,j}$ are zero-sum modulo $G_T$ with $m=\exp(G_T)$, we have  $U_{i,j}^{[m]},\,V_{i,j}^{[m]}\in \mathcal B(G_0)$ for all $i$ and $j$.
Re-factor each $U^{[m]}_{i,j}=\prod^\bullet_{x\in I_{i,j}}W^{(i,j)}_x$  into a product of $|I_{i,j}|\leq m$ atoms $W^{(i,j)}_x\in \mathcal A(G_0)$ via Proposition \ref{prop-fg-atoms-basic-transfer}.2. Then each ${\prod}^\bullet_{j\in [1,d+1]}\prod_{x\in I_{i,j}}^\bullet W^{(i,j)}_x$ is a product of at most $(d+1)m$ atoms from $\mathcal A(G_0)$, for every $i\geq 1$.
In view of \eqref{equal}, we have $${\prod}^\bullet_{j\in [1,d+1]}{\prod}_{x\in I_{i,j}}^\bullet W^{(i,j)}_x=U_{i,1}^{[m]}\bdot\ldots\bdot U_{i,d+1}^{[m]}=V_{i,1}^{[m]}\bdot\ldots\bdot V_{i,\ell_i}^{[m]}$$ with $V_{i,j}^{[m]}\in \mathcal B(G_0)$ a nontrivial zero-sum sequence for every $i\geq 1$ (since $\pi(V_{i,j})\in \mathcal A(\pi(G_0))$ with $m=\exp(G_T)$). In consequence, since  $\ell_i\rightarrow \infty$, it follows via \eqref{rho-ascend-chain} that $\rho_{(d+1)m}(G_0)=\infty$, contrary to assumption, which completes Item 2.
\end{proof}

We next extend Proposition \ref{prop-lambert-easy}.

\begin{corollary}\label{cor-fg-rho-transfer} Let $H$ be a Transfer Krull Monoid over a subset $G_0$ of a finitely generated abelian group $G=\Z^d\oplus G_T$ with torsion subgroup $G_T\leq G$, where $d\geq 0$. If
 $\rho_{(d+1)\exp(G_T)}(H)<\infty$,   then there is a constant $N_\rho\geq 1$ such that  $\rho_k(H)\leq N_\rho k<\infty$ for all $k\geq 1$.
\end{corollary}

\begin{proof}Since $\rho_k(H)=\rho_k(G_0)$ for all $k\geq 1$ by definition of a Transfer Krull Monoid, it suffices to prove the corollary when $H=\mathcal B(G_0)$, which we now assume.
Let $\pi:G\rightarrow \Z^d$   be the projection homomorphism with  kernel $G_T$. In view of \eqref{rho-ascend-chain}, we may embed $\Z^d\leq \R^d$ and  w.l.o.g. assume $\la G_0\ra=G$ and that  every $g\in G_0$ is contained in some atom, whence $\C(\pi(G_0))=\R^d$ (as in the proof of Proposition  \ref{prop-fg-rho-transfer}). Since $\rho_{(d+1)\exp(G_T)}(G_0)<\infty$, Proposition  \ref{prop-fg-rho-transfer} implies that $\rho_k(\pi(G_0))<\infty$ and $\rho_k(G_0)\leq \rho_{k\mathsf D(G_T)}(\pi(G_0))<\infty$ for all $k\geq 1$.
As a result, Theorem \ref{thm-rho-char} implies that we can apply Proposition \ref{prop-lambert-easy} to conclude there is a constant $N_\rho\geq 1$ (namely $\mathsf D(G_T)$ times the constant given by Proposition \ref{prop-lambert-easy}) such that $\rho_k(G_0)\leq \rho_{k\mathsf D(G_T)}(\pi(G_0))\leq N_\rho k<\infty$  for all $k\geq 1$, as desired.
\end{proof}

We now extend Theorem \ref{thm-tame-pseudo} to show the \textbf{weak tame degree} $\mathsf t_w(H)$ (the minimal integer $N\geq 1$ such that the conclusion of Main Theorem \ref{thm-tame-pseudo-tor} holds)  is  finite when $\rho_{(d+1)\exp(G_T)}(H)<\infty$. The definition is  a variation on the $\omega$ constant often used in conjunction with the tame degree \cite{Alfred-Ruzsa-book} \cite{Gerold-omega}. The proof is a variation on that used to prove Theorem \ref{thm-tame-pseudo}. We have opted to first give a proof of the simpler  Theorem \ref{thm-tame-pseudo}, to first present a version of the algorithm that avoids the delicate (and distracting) technicalities needed to get the algorithm to work in the general setting as stated in Theorem \ref{thm-tame-pseudo-tor}.

\begin{mtheorem}\label{thm-tame-pseudo-tor} Let $H$ be a  Krull Monoid  over a subset $G_0$ of a finitely generated abelian group $G=\Z^d\oplus G_T$ with torsion subgroup $G_T\leq G$, where $d\geq 0$. Suppose   $\rho_{(d+1)\exp(G_T)}(H)<\infty$. Then there exists an integer $N\geq 1$ such that, given any atoms $U_1,\ldots,U_k,V_1,\ldots,V_\ell\in \mathcal A(H)$ with $$U_1\cdot\ldots\cdot U_k=V_1\cdot\ldots\cdot V_\ell,$$ where $k,\,\ell\geq 1$,  then there exist atoms $W_1,\ldots,W_{\ell'}\in \mathcal A(H)$,  $r\in [1,k]$ and $\mathcal I\subseteq [1,\ell']$ such that $$U_1\cdot\ldots\cdot U_k=W_1\cdot\ldots\cdot W_{\ell'},$$  $\ell'\geq \ell$ and $U_r\mid \prod_{x\in \mathcal I} W_x$ with $|\mathcal I|\leq N$.
\end{mtheorem}

\begin{proof}
Let us first show that the theorem reduces to the case when $H=\mathcal B(G_0)$. To this end, assume the theorem holds for $\mathcal B(G_0)$ with bound $N\geq 1$.
We can w.l.o.g. assume $H$ is reduced, so we replace $H$ by $H_{\mathsf{red}}=H/H^\times$, and let $\varphi:H\rightarrow \Fc(P)$ be a divisor homomorphism and
 $\theta:H\rightarrow \mathcal B(G_0)$ the associated transfer homomorphism, with $G_0\subseteq G$ and $[\cdot]$ as defined in Section \ref{sec-prelim} (above \eqref{magic-krull-divisors}), so $\theta(S)=[\varphi(S)]$ for $S\in H$.
Then $\theta(U_1)\bdot\ldots\bdot \theta(U_k)=\theta(U_1\cdot\ldots\cdot U_k)=\theta(V_1\cdot\ldots\cdot V_\ell)=\theta(V_1)\bdot\ldots\bdot \theta(V_\ell)$ with all $\theta(U_i),\,\theta(V_j)\in \mathcal A(G_0)$. Applying the theorem to this factorization, we find atoms $W''_1,\ldots,W''_{\ell'}\in \mathcal A(G_0)$ with \be\label{forkknive}\theta(U_1\cdot\ldots\cdot U_k)=\theta(U_1)\bdot\ldots\bdot \theta(U_k)=W''_1\bdot\ldots\bdot W''_{\ell'}\ee and $\mathcal I\subseteq [1,\ell']$ such that $\ell'\geq \ell$ and $\theta(U_r)\mid \prod^\bullet_{x\in \mathcal I}W''_x$ with $|\mathcal I|\leq N$, for some $r\in [1,k]$. For $i\in [1,\ell']$, let $W'_i\in \Fc(P)$ be a sequence with $[W'_i]=W''_i$. Then using the definition of the composition map $\theta$ in   $\theta(U_r)\mid \prod^\bullet_{x\in \mathcal I}W''_x$ and \eqref{forkknive}, we find
\be\label{varphi-prd}[\varphi(U_r)]\mid {\prod}_{x\in \mathcal I}^\bullet [W'_x]\quad\und\quad [\varphi(U_1)]\bdot\ldots\bdot [\varphi(U_k)]=[W'_i]\bdot\ldots\bdot [W'_{\ell'}].\ee
Consequently, we can choose the pre-image sequences $W'_i\in \Fc(P)$ such that $$\varphi(U_r) \mid {\prod}_{x\in \mathcal I}^\bullet W'_x\quad\und\quad \varphi(U_1)\bdot\ldots\bdot \varphi(U_k)=W'_1\bdot\ldots\bdot W'_{\ell'}.$$
Since each $W'_i\in \Fc(P)$ with $[W'_i]=W''_i\in \mathcal A(G_0)\subseteq \mathcal B(G_0)$, for $i\in [1,\ell']$, it follows from \eqref{magic-krull-divisors} that $W'_i\in \varphi(H)$ for all $i\in  [1,\ell']$. Thus there are $W_i\in H$ with $\varphi(W_i)=W'_i$ for $i\in [1,\ell']$. We now have $\varphi(U_r)\mid {\prod}_{x\in \mathcal I}^\bullet \varphi(W_x)$ and $$\varphi(U_1\cdot\ldots\cdot U_k)=\varphi(U_1)\bdot\ldots\bdot \varphi(U_k)=\varphi(W_1)\bdot \ldots\bdot \varphi(W_{\ell'})=\varphi(W_1\cdot\ldots\cdot W_{\ell'}).$$ As a result, since $\varphi:H\rightarrow \varphi(H)$ is an isomorphism (by \eqref{magic-krull-divisors}), it follows that $U_1\cdot\ldots\cdot U_k=W_1\cdot\ldots\cdot W_{\ell'}$ with $W_i\in H$ for all $i$,  and since $\varphi$ is a divisor homomorphism, it follows from $\varphi(U_r)\mid {\prod}_{x\in \mathcal I}^\bullet \varphi(W_x)$ that $U_r\mid {\prod}_{x\in \mathcal I}^\bullet W_x$. Finally, since $\theta$ is a transfer homomorphism with  $\theta(W_i)=[\varphi(W_i)]=[W'_i]=W''_i\in \mathcal A(G_0)$, it follows that $W_i\in \mathcal A(H)$ for all $i$,  showing the theorem holds for $H$. It remains to prove the theorem when  $H=\mathcal B(G_0)$, which we now assume.

Let $\pi:G\rightarrow \Z^d$  and $\pi_T:G\rightarrow G_T$ be the projection homomorphisms with respective kernels $G_T$ and $\Z^d$. In view of \eqref{rho-ascend-chain}, we may embed $\Z^d\leq \R^d$ and  w.l.o.g. assume $\la G_0\ra=G$ and that  every $g\in G_0$ is contained in some atom, whence $\C(\pi(G_0))=\R^d$ (by Proposition  \ref{prop-notrivialG_0}). Since $\rho_{(d+1)\exp(G_T)}(G_0)<\infty$, Proposition \ref{prop-fg-rho-transfer} implies that $\rho_k(\pi(G_0))<\infty$ and $\rho_k(G_0)\leq \rho_{k\mathsf D(G_T)}(\pi(G_0))<\infty$ for all $k\geq 1$.
As a result, Theorem \ref{thm-rho-char} implies that $0\notin \C^*(\pi(G_0)^\diamond)$, while Corollary \ref{cor-fg-rho-transfer} implies there is an integer $N_\rho\geq 1$ such that $$\rho_k(G_0)\leq  N_\rho k\quad\mbox{ for all $k\geq 1$}.$$ Since $0\notin\C^*(\pi(G_0)^\diamond)$, it follows from Theorem \ref{thm-keylemmaII} that $\pi(G_0)$ is finitary. In view of Proposition \ref{prop-finitary-mintype}, there are only a finite number of minimal types. Let $\varphi_1,\ldots,\varphi_m\in \mathfrak T_m(\pi(G_0))$ be the distinct nontrivial minimal types for $\pi(G_0)$. For each $j\in [1,m]$, let $$Z_{\varphi_j}=Z_1^{(j)}\cup \ldots\cup Z_{s_j}^{(j)}$$ be the codomain of $\varphi_j$ with $s_j\leq d$ and $|Z_{\varphi_j}|\leq d$.

Consider  $U_1,\ldots,U_k,V_1,\ldots,V_\ell\in \mathcal A(G_0)$ with $$S:=U_1\bdot\ldots\bdot U_k=V_1\bdot\ldots\bdot V_\ell,$$ where $k,\,\ell\geq 1$. Now \be\label{ellsmall-tor} \ell\leq \rho_k(G_0)\leq N_\rho k.\ee Let $S=g_1\bdot\ldots\bdot g_{|S|}$ be an indexing of the terms of $S$.
Let $I_1\cup \ldots\cup I_k=[1,|S|]=J_1\cup \ldots\cup J_\ell$ be disjoint partitions such that $$S(I_i)=U_i\quad \und\quad S(J_j)=V_j\quad\mbox{ for all $i\in[1,k]$ and $j\in [1,\ell]$}.$$
By  Proposition \ref{prop-fg-atoms-basic-transfer}.1, each $\pi(V_i)$ for $i\in [1,\ell]$ factors into a product of $$d_i\leq\mathsf D(G_T)$$ atoms modulo $G_T$, say $$V_i={\prod}_{t\in [1,d_i] }^\bullet V_{i,t}\quad\mbox{ with every $\pi(V_{i,t})\in\mathcal A(\pi(G_0))$}.$$ We thus have partitions $J_i=J_{i,1}\cup \ldots\cup J_{i,d_i}$ such that $$S(J_{i,t})=V_{i,t}\quad \mbox{ for every $i\in [1,\ell]$ and $t\in [1,d_i]$}.$$
Since $0\notin \C^*(\pi(G_0)^\diamond)$, we can apply Theorem \ref{thm-structural-char} to each atom $\pi(V_{i,t})$ for $i\in [1,\ell]$ and $t\in [1,d_i]$. Let $N_R\geq 0$ be the global bound from Theorem \ref{thm-structural-char} (which we can assume is an integer) and let $V_{i,t}=R_{i,t}\bdot \prod_{j=1}^mS_{i,t}^{(j)}$, for $i\in [1,\ell]$ and $t\in [1,d_i]$, be the resulting factorization given by Theorem \ref{thm-structural-char}, with $\pi(S_{i,t}^{(j)})\in \Fc(\pi(G_0))$ corresponding to the minimal type $\varphi_j$, so $\sigma(\pi(S_{i,t}^{(j)}))\in \C_\Z(Z_{\varphi_j})$.  We can then apply Proposition \ref{prop-swapping} to each $\pi(S_{i,t}^{(j)})$.
Let $J_{i,t}^{(j)}\subseteq [1,|S|]$ be disjoint subsets with $$S(J_{i,t}^{(0)})=R_{i,t}\quad\und \quad S(J_{i,t}^{(j)})=S_{i,t}^{(j)},\quad\mbox{ for $i\in [1,\ell]$, $t\in [1,d_i]$ and $j\in [1,m]$}.$$
 Moreover, for each $i\in [1,\ell]$, $t\in [1,d_i]$, $j\in [1,m]$ and $n\in [1,s_j]$, let $J_{i,t}^{(j,n)}\subseteq J_{i,t}^{(j)}$ be the subset of all $x\in J_{i,t}^{(j)}$ with $\varphi_j(x)\in  Z^{(j)}_n$ at depth $n$.
Then  $$J_{i,t}:=\bigcup_{j=0}^mJ_{i,t}^{(j)}=J_{i,t}^{(0)}\cup \bigcup_{j=1}^m\bigcup_{n=1}^{s_j}J_{i,t}^{(j,n)}\quad\mbox{  for every $i\in [1,\ell]$ and $t\in [1,d_i]$.}$$
Let \begin{align*}&&X_0=\bigcup_{i=1}^\ell \bigcup_{t=1}^{d_i}J_{i,t}^{(0)},\quad
\Omega_0=\{(0,\pi(g_x),\pi_T(g_x)):\; x\in  X_0\}\quad\und\quad\\
&&\Omega_\diamond=\{(j,z,a):\;j\in [1,m],\,z\in Z_{\varphi_j},\,a\in G_T\}.
\end{align*}
Moreover, partition $$\Omega_\diamond=\Omega_1\cup \ldots\cup \Omega_d$$ such that $\Omega_n$ consists of all $(j,z,a)\in \Omega_\diamond$ with $z\in Z_n^{(j)}$ at depth $n$. 
Apply Proposition \ref{prop-swapping} to each $\pi(S_{i,t}^{(j)})$ (for $i\in [1,\ell]$, $t\in [1,d_i]$ and $j\in [1,m]$)  and fix a system of subsets $T_x\subseteq J^{(j)}_{i,t}$, for each  $x\in J^{(j)}_{i,t}$, such that the conclusions of Proposition \ref{prop-swapping} hold.

Since Theorem \ref{thm-structural-char} implies $|J_{i,t}^{(0)}|=|R_{i,t}|\leq N_R$ for all $i\in [1,\ell]$ and $t\in [1,d_i]$, since $d_i\leq \mathsf D(G_T)$ for all $i\in [1,\ell]$, and since $|Z_{\varphi_j}|\leq d$ for all $j\in [1,m]$, we have \be\label{finitestuff-tor}|X_0|\leq \ell \mathsf D(G_T)N_R\quad\und\quad |\Omega_\diamond|=|G_T|\cdot \Sum{j=1}{m}|Z_{\varphi_j}|\leq md|G_T|.\ee
We view $\Omega=\Omega_0\cup \Omega_\diamond=\Omega_0\cup \Omega_1\cup \ldots\cup \Omega_d$ as the set of \emph{support types} for $S=V_1\bdot\ldots\bdot V_\ell$. A support type $\tau\in \Omega_n$ is said to be at depth $n$.
Note, if $\tau=(j,z,a)$ with $j\geq 1$, then the depth of $\tau$ equals the depth of $z\in Z_{\varphi_j}$. For each $x\in [1,|S|]$, we have $x\in J_{i,t}^{(j)}$ for some unique $i\in [1,\ell]$, $t\in [1,d_i]$ and $j\in [0,m]$, allowing us to define \begin{align*}&\mathsf s(x)=(0,\pi(g_x),\pi_T(g_x))\in \Omega_0\quad\mbox{ when $j=0$,}\quad \und \\&\mathsf s(x)=(j,\varphi_j(\pi(g_x)),\pi_T\Big(\sigma(S(T_x))\Big))\in \Omega_\diamond\quad\mbox{ when  $j\geq 1$}.\end{align*}
Note Proposition \ref{prop-swapping} implies $\varphi_j(\pi(g_x))=\sigma(\pi(S)(T_x))=\pi\Big(\sigma(S(T_x))\Big)$.
For $I\subseteq [1,|S|]$, $\mathsf s(I)\in \Fc(\Omega)$ is a sequence of support types from $\Omega$.  We associate the depth of $\mathsf s(x)$ (defined above) as the depth of $x\in [1,|S|]$.

Let
$$\alpha=\min\left\{|X_0\cap I_i|+\Summ{\tau\in \Omega_\diamond}\Big(\frac{\ell \vp_\tau(\mathsf s(I_i))}{\vp_\tau\big(\mathsf s\big([1,|S]]\big)\big)}+1\Big):\;i\in [1,k]\right\}.$$ Technically, we exclude any terms $\tau\in \Omega_\diamond$ in the sum defining $\alpha$ with  $\vp_\tau\big(\mathsf s\big([1,|S]]\big)\big)=0$.  Then
\ber\nn k\alpha&\leq& \Sum{i=1}{k}|X_0\cap I_i|+\Sum{i=1}{k}\Summ{\tau\in \Omega_\diamond}\frac{\ell \vp_\tau(\mathsf s(I_i))}{\vp_\tau\big(\mathsf s\big([1,|S]]\big)\big)}+k|\Omega_\diamond|
 \\\nn &=& |X_0|+\Summ{\tau\in \Omega_\diamond}\Sum{i=1}{k}\frac{\ell \vp_\tau(\mathsf s(I_i))}{\vp_\tau\big(\mathsf s\big([1,|S]]\big)\big)}+k|\Omega_\diamond|\leq |X_0|+(\ell+k) |\Omega_\diamond|\\
 &\leq & \ell(N_R\mathsf D(G_T)+md|G_T|)+kmd|G_T|\leq kN_\rho(N_R\mathsf D(G_T)+md|G_T|)+kmd|G_T|,\label{firstthere-tor}\eer with the first inequality in \eqref{firstthere-tor} in view of  \eqref{finitestuff-tor}, and the second in view of \eqref{ellsmall-tor}. Thus  $$\alpha\leq N:=N_\rho(N_R\mathsf D(G_T)+md|G_T|)+md|G_T| ,$$ which is a global bound independent of the $U_i$ and $V_j$.

 Let $r\in [1,k]$ be an index attaining the minimum in the definition of $\alpha$. Then $|X_0\cap I_r|\leq \alpha\leq N$, ensuring that there is some subset $\mathcal I_0\subseteq [1,\ell]$ with \be\label{Y0go-tor}X_0\cap I_r\subseteq \bigcup_{i\in \mathcal I_0}J_i^{(0)}\quad \und\quad |\mathcal I_0|\leq |X_0\cap I_r|\leq N.\ee Likewise, letting $$n_\tau=\left\lceil \frac{\vp_{\tau}(\mathsf s(I_r))}{\vp_\tau\big(\mathsf s \big([1,|S|]\big)\big)/\ell}\right\rceil< \frac{\ell \vp_\tau(\mathsf s(I_r))}{\vp_\tau\big(\mathsf s\big([1,|S]]\big)\big)}+1\leq \ell+1\quad\mbox{ for $\tau\in \Omega_\diamond $},$$ we have \be\label{waterpool-tor}\Summ{\tau\in \Omega_\diamond }n_\tau\leq \alpha-|X_0\cap I_r|\leq N-|\mathcal I_0|.\ee
 We interpret $n_\tau=0$ when $\vp_\tau\big(\mathsf s \big([1,|S|]\big)\big)=0$.

 We now describe how the $W_1,\ldots,W_{\ell'}\in \mathcal B(G_0)$ can  be constructed.
 The idea is as follows. An index set $I\subseteq [1,|S|]$ indexes a sequence $S(I)\mid S$, but it also indexes a sequence $\mathsf s(I)\in \Fc(\Omega)$, obtained by replacing each indexed term in the sequence $S(I)$ with its corresponding support type from $\Omega$, so $\mathsf s(I)=\prod^\bullet_{x\in I}\mathsf s(x)$.  When $I\subseteq [1,|S|]\setminus X_0$, we have $\mathsf s(I)\in \Fc(\Omega_\diamond)$ with $\Omega_\diamond$ a fixed, finite set independent of $S$. Let $\tau\in \Omega_\diamond$.  If we select a subset  $\mathcal I_\tau\subseteq [1,\ell]$ with $|\mathcal I_\tau|=n_\tau$ such that the $\vp_\tau(\mathsf s(J_x))$, for $x\in I_\tau$, are the $n_\tau$ largest values occurring over all  $\vp_{\tau}(\mathsf s(J_i))$ with $i\in [1,\ell]$,
 then the definition of $n_\tau$ ensures that $$\vp_\tau\Big(\mathsf s\Big({\bigcup}_{z\in \mathcal I_\tau}J_z\Big)\Big)\geq
 n_\tau\Big(\vp_\tau\big(\mathsf s \big([1,|S|]\big)\big)/\ell\Big)\geq \vp_\tau(\mathsf s(I_r)),$$
 with the first inequality holding since the sum of the $n_\tau$ largest terms in a sum of $\ell$ non-negative terms is always at least $n_\tau$ times the average value of all terms being summed.
 As a result, $\mathsf s(I_r)\mid \prod^\bullet_{z\in \mathcal I}\mathsf s(J_z)$, where $\mathcal I=\mathcal I_0\cup \bigcup_{\tau\in \Omega_\diamond}\mathcal I_\tau$, with $|\mathcal I|\leq N$ in view of \eqref{waterpool-tor}. However, since the map $\mathsf s$ is not injective, this does not guarantee that the associated sequence $U_r=S(I_r)$ is a subsequence of the associated sequence $\prod^\bullet_{z\in \mathcal I}V_z=\prod^\bullet_{z\in \mathcal I}S(J_z)$. We do, however, have $S(X_0\cap I_r)\mid \prod^\bullet_{z\in \mathcal I_0}S(J_z)=\prod^\bullet_{z\in \mathcal I_0}V_z$ in view of \eqref{Y0go-tor}. To deal with the terms from $\Omega_\diamond$, we must  use the sequences $S(T_x)$ given by Proposition \ref{prop-swapping} to exchange terms between the $V_i$.

 If there are terms $x\in J_{i,t}^{(j)}$ and $y\in J_{i',t'}^{(j)}$ with $\mathsf s(x)=\mathsf s(y)$, $i\neq i'$ and $j\geq 1$, say with $g_x$ and $g_y$ at depth $n$, then $\pi_T\Big(\sigma(S(T_x))\Big)=\pi_T\Big(\sigma(S(T_y))\Big)$, while  Proposition  \ref{prop-swapping} implies that $$\pi\Big(\sigma(S(T_x))\Big)=\sigma(\pi(S)(T_x))=\varphi_j(\pi(g_x))=\varphi_j(\pi(g_y))=
 \sigma(\pi(S)(T_y))=\pi\Big(\sigma(S(T_y))\Big).$$ Combined with $\pi_T\Big(\sigma(S(T_x))\Big)=\pi_T\Big(\sigma(S(T_y))\Big)$, we conclude that
 \be\label{mainted}\sigma(S(T_x))=\sigma(S(T_y)).\ee
 If we exchange these sets, defining $$K_{i,t}^{(j)}=(J_{i,t}^{(j)}\setminus T_x)\cup T_y\quad\und\quad K_{i',t'}^{(j)}=(J_{i',t'}^{(j)}\setminus T_y)\cup T_x,$$ and correspondingly define $$K_i=(J_i\setminus T_x)\cup T_y\quad\und\quad K_{i'}=(J_{i'}\setminus T_y)\cup T_x,$$
 then the new sequences $W_i=S(K_i)$ and $W_{i'}=S(K_{i'})$ will still be zero-sum, though we do not guarantee that they remain atoms. However, since $y\in K_i$ and $x\in K_{i'}$, they are non-empty. Consequently, if either $W_i$ or $W_{i'}$ is not an atom, then we can re-factor them to write $V_i\bdot V_{i'}=W_i\bdot W_{i'}=V'_1\bdot\ldots\bdot V'_\omega$ as a product of $\omega\geq 3$ atoms. This leads to a factorization $U_1\bdot\ldots\bdot U_k=V_1\bdot\ldots\bdot V_{\ell}\bdot V_i^{[-1]}\bdot V_{i'}^{[-1]}\bdot V'_1\bdot\ldots \bdot V'_\omega$ into $\ell'=\ell-2+\omega>\ell$ atoms. In this case, we begin from scratch using this factorization in place of the original one $U_1\bdot\ldots\bdot U_k=V_1\bdot\ldots\bdot V_\ell$. As $\ell'\leq |S|<\infty$, we cannot start from scratch endlessly, meaning eventually we will never encounter this problem, allowing us to  w.l.o.g. assume $W_i=S(K_i)$ and $W_{i'}=S(K_{i'})$ are always atoms (where $\ell'=\ell$ may have increased in size from the original $\ell$ given in the hypotheses).
 Furthermore, we still have $\sigma(\pi(S)(K^{(j)}_{i,t}))=\sigma(\pi(S)(J_{i,t}^{(j)}))\in \C(Z_{\varphi_j})$ and $\sigma(\pi(S)(K^{(j)}_{i',t'}))=\sigma(\pi(S)(J_{i',t'}^{(j)}))\in \C(Z_{\varphi_j})$ by \eqref{mainted}, with the inclusions originating from our application of Theorem \ref{thm-structural-char} at the start of the proof. Thus Proposition \ref{prop-swapping} can still be applied to $K_{i,t}^{(j)}$ and $K_{i',t'}^{(j)}$ if we later wish to continue with further such swaps between these sets (though we do not guarantee nor need that $\pi(S)(K_{i,t}^{(j)})$ and $\pi(S)(K_{i',t'}^{(j)})$ remain atoms). Proposition \ref{prop-swapping} guarantees that the set $T_x$ contains \emph{no} terms with depth greater than $g_x$, and has $x$ as the \emph{unique} $a\in T_x$ with $g_a$ having depth equal to that of $x$.
 Likewise for $T_y$. Thus when swapping the sets $T_x$ and $T_y$, we leave unaffected all terms in $J_{i,t}^{(j)}$ and $J_{i',t'}^{(j)}$  with depth at least $n$, apart from the exchanging of $x$ for $y$.
 This ensures that terms previously swapped but at a higher or equal depth will remain unaffected by exchanging $T_x$ and $T_y$.

 As in the proof of Theorem \ref{thm-tame-pseudo}, the sequences $\mathsf s(K_{i,t}^{(j)})=\mathsf s(J_{i,t}^{(j)})$ and $\mathsf s(K_{i',t'}^{(j)})=\mathsf s(J_{i',t'}^{(j)})$ remain unchanged, though this  now requires a short argument. The value of $\mathsf s(z)$ for $z\in [1,|S|]$ depends upon the sequence $T_z$. After swapping $T_x$ and $T_y$, the sets $T_z$ need to be redefined via Proposition \ref{prop-swapping}. We must show that this can be done in such a way that $\sigma(S(T_z)))=\sigma(S(T'_z))$ for all $z$, where $T'_z$ is the new sequence associated to $z$ after performing the swap. 
 Unless $z\in J_{i,t}^{(j)}$ or $z\in J_{i',t'}^{(j)}$, the values for $T_z$ can be left unaffected, so $T_z=T'_z$. Let us consider the case when $z\in J_{i,t}^{(j)}$. The case $z\in J_{i',t'}^{(j)}$ will then follow by an analogous argument. If $T_z\cap T_x=\emptyset$, then the value of $T_z$ can also be left unaffected.
 Otherwise, Proposition \ref{prop-swapping}.2 implies that either $z\in T_x$ or $x\in T_z$. If $z\in T_x$, the Proposition \ref{prop-swapping}.1 implies that $T_z\subseteq T_x$, and we can again leave the value of $T_z$ unchanged, though note that $T'_z\subseteq K_{i',t'}^{(j)}$ while $T_z\subseteq J_{i,t}^{(j)}$. In the remaining case $x\in T_z$ but $z\notin T_x$ (so $x\neq z$), we have $T_x\subseteq T_z$ by Proposition \ref{prop-swapping}.1 with $\sigma(S(T_x))=\sigma(S(T_y))$ by \eqref{mainted}.
 This allows us to define $T'_z=T_z\setminus T_x\cup T_y\subseteq K_{i,t}^{(j)}$, which then satisfies $\sigma(S(T'_z))=\sigma(S(T_z))-\sigma(S(T_x))+\sigma(S(T_y))=\sigma(S(T_z))$. The equivalent defining conditions $1'$ and $2'$ for Proposition \ref{prop-swapping} also hold for the newly defined set system $T'_z$. This shows it is possible to adjust the sequences $T_z$ after swapping $T_x$ for $T_y$ in such a way that the value $\mathsf s(z)$ remains unaffected.
 In particular, the value of $\vp_\tau\big(\mathsf s\big([1,|S|]\big)\big)$ remains unchanged, ensuring that the value of $\alpha$ is unaffected when replacing $J_i$ and $J_{i'}$ by $K_i$ and $K_{i'}$,  and that $r\in [1,k]$ remains an index attaining the minimum in the definition of $\alpha$ (note, the numerators in the definition of $\alpha$ depend on the $U_i$, not the $V_j$).
 Swapping the elements $x$ and $y$ in this fashion leaves all elements from $X_0$, as well as any $J_{c,d}^{(b)}$ with $b\neq j$, unaltered, and the sequences $W_i$ and $W_{i'}$ remain nontrivial, as $W_i$ contains $g_y$, and $W_{i'}$ contains $g_x$. Since, apart from $x$ and $y$, only terms with depth less than $n$ are affected by the swap, it follows that the sets $\mathcal I_{\tau'}$, corresponding to any type $\tau'\in \Omega_\diamond$ with depth at least $n$, still have the property that they index the $\vp_\tau(\mathsf s(J_z))$, for $z\in I_{\tau'}$, with  the $n_{\tau'}$ largest values occurring over all  $\vp_{\tau'}(\mathsf s(J_i))$ with $i\in [1,\ell]$.

With these observations in mind, we can now describe how the zero-sums $V_i$ must be modified. Begin with any type $\tau\in \Omega_\diamond$ having maximal available depth. Construct the subset $\mathcal I_\tau$ for the current factorization $S=V_1\bdot\ldots\bdot V_\ell$  as described above. Then $\vp_\tau
\big(\mathsf s(\bigcup_{z\in \mathcal I_\tau}J_z)\big)\geq \vp_\tau(\mathsf s(I_r))$. If $\bigcup_{z\in \mathcal I_\tau}J_z$ contains all elements from $I_r$ having type $\tau$, then nothing need be done, we discard $\tau$ from the list of available types from $\Omega_\diamond$,  we select the next available type from $\Omega_\diamond$ with maximal depth, and continue  once more. On the other hand, if there is some $x\in I_r$ with type $\tau$ not contained in $\bigcup_{z\in \mathcal I_\tau}J_z$, then $\vp_\tau
\big(\mathsf s(\bigcup_{z\in \mathcal I_\tau}J_z)\big)\geq \vp_\tau(\mathsf s(I_r))$ ensures that there must be some $y\in \bigcup_{z\in \mathcal I_\tau}J_z$ having type $\tau$ with $y\notin I_r$. In this case,  perform the swap of $T_x$ and $T_y$ described above, and redefine our factorization $V_1\bdot\ldots\bdot V_\ell$ by replacing $J_{i,t}^{(j)}$ and $J_{i',t'}^{(j)}$ by $K_{i,t}^{(j)}$ and $K_{i',t'}^{(j)}$, where $x\in J_{i,t}^{(j)}$ and $y\in J_{i',t'}^{(j)}$, and correspondingly replacing $V_i$ and $V_{i'}$ by $W_i$ and $W_{i'}$. Also, adjust the values of the auxiliary sets $T_z$ as described above. To simplify notation, redefine  $V_i$, $V_{i'}$, $J_{i,t}^{(j)}$, $J_{i',t'}^{(j)}$ and $T_z$ accordingly so as to reflect the new current state that now has $x\in \bigcup_{z\in \mathcal I_\tau}J_z$ and $y\notin \bigcup_{z\in \mathcal I_\tau}J_z$. If we now have $\bigcup_{z\in \mathcal I_\tau}J_z$ containing all elements from $I_r$ having type $\tau$, then nothing need be done, we discard $\tau$ from the list of available types in $\Omega_\diamond$ and carry on as before. If this is not the case, we again find a new term $x'\in I_r$ with type $\tau$ not contained in $\bigcup_{z\in \mathcal I_\tau}J_z$, find a new term $y'\notin \bigcup_{z\in \mathcal I_\tau}J_z$ with type $\tau$ and swap $x'$ and $y'$ as before by use of Proposition \ref{prop-swapping}. Since the depth of $x$ and $x'$ are the same, we will not swap $x$ back out of $\bigcup_{z\in \mathcal I_\tau}J_z$ when doing so, nor indeed any other element from  $\bigcup_{z\in \mathcal I_\tau}J_z$ having type $\tau$.
Thus, iterating such a procedure, we will eventually obtain that $\bigcup_{z\in \mathcal I_\tau}J_z$ contains every element of $I_r$ having type $\tau$, in which case we move on to the next available type $\tau'$ with maximal available depth. We repeat the same for procedure for $\tau'$ as we did for $\tau$ (and, later, as we did for any support types previously discarded before selecting $\tau'$). We first construct the subset $\mathcal I_{\tau'}$ for the current state for $S=V_1\bdot\ldots\bdot V_\ell$, and then swap elements into $\bigcup_{z\in \mathcal I_{\tau'}}J_z$ until it contains all elements from $I_r$ having type $\tau'$.
While doing so, since we always first choose support types with maximal available depth, we are assured that any type $\tau''$ that has already been discarded had depth at least that of $\tau'$, and thus no  element of type $\tau''$ will be moved when swapping at the later stage for $\tau'$, ensuring that prior work cannot be undone.
Continue until all support types from $\Omega_\diamond$ have been exhausted.
Once the process ends, we now have a new factorization $U_1\bdot\ldots\bdot U_k=S=W_1\bdot\ldots\bdot W_\ell$, where $W_i$ reflects the final state of $V_i$ after running the above process, such that $U_r\mid \prod_{z\in \mathcal I}^\bullet W_i$, where $\mathcal I=\mathcal I_0\cup\bigcup_{\tau\in \Omega_\diamond}\mathcal I_\tau$, with $|\mathcal I|\leq |X_0\cap I_r|+\Summ{\tau\in \Omega_\diamond}n_\tau\leq N$, completing the proof.
\end{proof}

The following basic proposition shows that having finite elasticities implies a Krull Monoid is always locally tame.

 \begin{proposition}\label{prop-localtame}
 Let $H$ be a Krull Monoid with divisor homomorphism $\varphi:H\rightarrow \Fc(P)$ and let $U\in \mathcal A(H)$. Then $\mathsf t(H,U)\leq \rho_{|\varphi(U)|}(H)$.
 \end{proposition}

\begin{proof}
Let $U_1,\ldots,U_k\in \mathcal A(H)$ with $U\mid U_1\bdot\ldots\bdot U_k$ and let $I\subseteq [1,k]$ be a minimal subset with $U\mid \prod_{i\in I}U_i$. Then $\varphi(U)\mid \varphi(\prod_{i\in I}U_i)=\prod_{i\in I}\varphi(U_i)$. Each $\varphi(U_i)$ is a nontrivial sequence as each $U_i\in \mathcal A(H)$. Thus we trivially have $\varphi(U)\mid \prod_{i\in J}\varphi(U_i)$ for some $J\subseteq I$ with $|J|\leq |\varphi(U)|$. By definition of a divisor homomorphism, it follows that $U\mid \prod_{i\in J}U_i$, and now the minimality of $I$ ensures $I=J$ with $|I|=|J|\leq |\varphi(U)|$. Since $U\mid \prod_{i\in I}U_i$, we have $\prod_{i\in I}U_i=U\cdot V_2\cdot\ldots\cdot V_r$ for some $V_2,\ldots,V_r\in \mathcal A(H)$, and by definition of the elasticities and \eqref{rho-ascend-chain}, we must have $r\leq \rho_{|I|}(H)\leq \rho_{|\varphi(U)|}(H)$.
The result now follows.
\end{proof}

Having now established Theorem \ref{thm-tame-pseudo-tor}, we can immediately extend  Theorem \ref{thm-Delta-finite} to the more general finitely generated group setting.

\begin{mtheorem}\label{thm-Delta-finite-tor}
Let $H$ be a Transfer Krull Monoid over a subset $G_0$ of a finitely generated abelian group $G=\Z^d\oplus G_T$ with torsion subgroup $G_T\leq G$, where $d\geq 0$. Suppose    $\rho_{(d+1)\exp(G_T)}(H)<\infty$ and let $N=\max\{2,\mathsf t_w(G_0)\}$.
\begin{itemize}
\item[1.]  $\rho_k(H)-\rho_{k-1}(H)\leq N<\infty$ for all $k\geq 2$.
\item[2.] $\max \Delta(H)\leq \rho_N(H)-N<\infty$. In particular,  $\Delta(H)$ is finite.
\item[3.] The Structure Theorem for Unions holds in $H$.
\item[4.] If $H$ is also a Krull Monoid, the catenary degree $\mathsf c(H)\leq \rho_N(H)<\infty$ is finite.
\item[5.] If $H$ is also a Krull Monoid, then $H$ is locally tame.
\end{itemize}
\end{mtheorem}

\begin{proof}
Since $H$ is a Transfer Krull Monoid over $G_0$, we have $\rho_k(H)=\rho_K(G_0)$  and $\mathcal U_k(H)=\mathcal U_k(G_0)$ for all $k\geq 1$, and $\Delta(H)=\Delta(G_0)$. Moreover, when $H$ is a Krull Monoid, we have $\mathsf c(H)\leq \max\{\mathsf c(G_0),2\}$ (see Section \ref{sec-prelim}). It thus suffices to prove the theorem when $H=\mathcal B(G_0)$, which we now assume.
By Corollary \ref{cor-fg-rho-transfer}, there is a constant $N_\rho\geq 1$ such that $\rho_k(G_0)\leq N_\rho k<\infty$  for all $k\geq 1$.  Items 1--4 now follows by the identical arguments given in Theorem \ref{thm-Delta-finite}, simply replacing the use of Theorem \ref{thm-tame-pseudo} with Theorem \ref{thm-tame-pseudo-tor}. Item 5 follows from Proposition \ref{prop-localtame} and Item 1 (which inductively shows $\rho_k(H)=\rho_k(G_0)<\infty$ for all $k$).
\end{proof}

Consider an infinite abelian group $G$ and let  $G_0 \subseteq G$ be a subset such that there is some  $m\geq 0$ so that  every element $g \in G$ is the sum of at most $m$ elements from $G_0$ (e.g., this assumption holds for integrally closed finitely generated algebras over perfect fields). Then $\Delta (G_0) $ is infinite by \cite[Theorem 1.1]{Ha02c}, so  Main Theorem \ref{thm-Delta-finite-tor} further implies that $\rho_{(d+1)\exp(G_T)} (G_0)=\infty$ must also  be infinite.

We now extend the definition of $G_0^\diamond \subseteq G_0$.

\begin{definition}
Let $G= \Z^d\oplus G_T$ be a finitely generated abelian group with torsion subgroup  $G_T\leq G$, where $d\geq 0$, and let $\pi:G\rightarrow \Z^d\leq \R^d$ be the projection homomorphism with kernel $G_T$. For $G_0\subseteq G$, we define $$G_0^\diamond=\{g\in G_0:\;\pi(g)\in \pi(G_0)^\diamond\}.$$
\end{definition}

The following extends Corollary \ref{cor-G_0diamond-Zequiv}

\begin{proposition}\label{prop-fg-diamondZ}
Let $G= \Z^d\oplus G_T$ be a finitely generated abelian group with torsion subgroup  $G_T\leq G$, where $d\geq 0$, and let $G_0\subseteq G$ be a subset. Suppose   $\rho_{(d+1)\exp(G_T)}(G_0)<\infty$. Then
\ber\nn G_0^\diamond&=&\{g\in G_0:\; \sup\{\vp_{g}(U):\; U\in \mathcal A^{\mathsf{elm}}(G_0)\}=\infty\}\\\nn &=&\{g\in G_0:\; \sup\{\vp_{g}(U):\; U\in \mathcal A(G_0)\}=\infty\}.\eer
\end{proposition}

\begin{proof} Since $G_T$ has finite exponent, an element $g\in G_0$ is contained in an atom if and only if $\pi(g)\in \pi(G_0)$ is contained in an atom. Thus Corollary \ref{cor-nondegen} and the definition of $G_0^\diamond$ ensure that  removing elements from $G_0$ contained in no atom does not affect $G_0^\diamond$, allowing us to w.l.o.g. assume every $g\in G_0$ is contained in an atom. We may embed $\Z^d\leq \R^d$. Let $\pi:G\rightarrow \Z^d$ be the projection homomorphism with kernel $G_T$, and let $m=\exp(G_T)$. In view of \eqref{rho-ascend-chain}, we may w.l.o.g. assume $\la G_0\ra=G$.
Since $\la G_0\ra=G$ and  every $g\in G_0$ is contained in an atom, Corollary \ref{cor-nondegen} yields $\C(\pi(G_0))=\R^d$. Since  $\rho_{(d+1)\exp(G_T)}(G_0)<\infty$,  Proposition \ref{prop-fg-rho-transfer} implies $\rho_{d+1}(\pi(G_0))<\infty$, whence Theorem \ref{thm-rho-char} implies that $0\notin \C^*(\pi(G_0)^\diamond)$, in turn implying that $\pi(G_0)$ is finitary (by Theorem \ref{thm-keylemmaII}).

Let $g\in G_0^\diamond$ be arbitrary. Then $\pi(g)\in \pi(G_0)^\diamond$, whence Corollary \ref{cor-G_0diamond-Zequiv} implies that there is a sequence of elementary atoms $V_i\in \mathcal A^{\mathsf{elm}}(\pi(G_0))$, for $i=1,2,\ldots,$ with $\vp_{\pi(g)}(V_i)\rightarrow \infty$. Moreover, since $0\notin \pi(G_0)^\diamond$ (by Proposition \ref{prop-G_0diamond-1st-easy-equiv}.2), we have $\pi(g)\neq 0$. For $i\geq 1$, let $U_i\in \Fc(G_0)$ be a sequence with $\pi(U_i)=V_i$, \ $\vp_g(U_i)=\vp_{\pi(g)}(V_i)$ and $|\supp(U_i)|=|\supp(V_i)|$. Thus $\vp_g(U_i^{[m]})\rightarrow \infty$.
Since $V_i$ is an elementary atom containing the nonzero element $\pi(g)$, Proposition \ref{prop-char-minimal-pos-basis} implies that $\supp(V_i)$ is a minimal positive basis, in which case $\mathcal A(Y)=\emptyset$ for any proper subset $Y\subset \supp(V_i)$. Thus, since $|\supp(U_i)|=|\supp(V_i)|$, it follows that $\mathcal A(\pi(Y))=\emptyset$ for any proper subset $Y\subset \supp(U_i)$, in turn implying $\mathcal A(Y)=\emptyset$ as well.
In view of Proposition \ref{prop-fg-atoms-basic-transfer}.2, each  $U_i^{[m]}$ factors into a product of at most $m$ atoms. Thus, letting $W_i\mid U_i^{[m]}$ be a subsequence which is an atom with $\vp_g(W_i)$ maximal, it follows that $\vp_g(W_i)\rightarrow \infty$ in view of $\vp_g(U_i^{[m]})\rightarrow \infty$. Moreover, since $\mathcal A(Y)=\emptyset$ for any proper subset $Y\subset \supp(U_i)$, it follows that each $W_i\in \mathcal A^{\mathsf{elm}}(G_0)$ is an elementary atom with $\vp_g(W_i)\rightarrow \infty$, which establishes the inclusion $G_0^\diamond\subseteq \{g\in G_0:\; \sup\{\vp_{g}(U):\; U\in \mathcal A^{\mathsf{elm}}(G_0)\}=\infty\}$. The inclusion $\{g\in G_0:\; \sup\{\vp_{g}(U):\; U\in \mathcal A^{\mathsf{elm}}(G_0)\}\subseteq \infty\}\subseteq \{g\in G_0:\; \sup\{\vp_{g}(U):\; U\in \mathcal A(G_0)\}=\infty\}$ holds trivially in view of $\mathcal A^{\mathsf{elm}}(G_0)\subseteq \mathcal A(G_0)$.

To establish the final reverse inclusion, let $g\in G_0$ be an element with $\sup\{\vp_{g}(U):\; U\in \mathcal A(G_0)\}=\infty$. Let $U_i\in \mathcal A(G_0)$, for $i=1,2,\ldots,$ be a sequence of atoms with $\vp_g(U_i)\rightarrow \infty$. By Proposition \ref{prop-fg-atoms-basic-transfer}.1, each $\pi(U_i)$ factors as a product of at most $\mathsf D(G_T)$ atoms. Thus, letting $V_i\mid \pi(U_i)$ be an atom $V_i\in \mathcal A(\pi(G_0))$ with $\vp_{\pi(g)}(V_i)$ maximal, it follows that $\vp_{\pi(g)}(V_i)\rightarrow \infty$ in view of $\vp_g(U_i)\rightarrow \infty$. Hence $\pi(g)\in \{x\in \pi(G_0):\;\sup\{\vp_{x}(V):\; V\in \mathcal A(\pi(G_0))\}=\infty\}=\pi(G_0)^\diamond$, with the equality in view of Corollary \ref{cor-G_0diamond-Zequiv}, which implies $g\in G_0^\diamond$ by definition of $G_0^\diamond$. This establishes the reverse inclusion $\{g\in G_0:\; \sup\{\vp_{g}(U):\; U\in \mathcal A(G_0)\}=\infty\}\subseteq G_0^\diamond$, completing the proof.
\end{proof}

Next, we extend Corollary \ref{cor-structural-Lambert}.

\begin{proposition}\label{prop-lambert-tor-a}
Let $G= \Z^d\oplus G_T$ be a finitely generated abelian group with torsion subgroup  $G_T\leq G$, where $d\geq 0$, and  let $G_0\subseteq G$ be a subset. Suppose   $\rho_{(d+1)\exp(G_T)}(G_0)<\infty$. Then $G_0\setminus G_0^\diamond\subseteq G_0$ is a Lambert subset with $\mathcal A(G_0^\diamond)=\emptyset$.
\end{proposition}

\begin{proof} In view of \eqref{rho-ascend-chain} and Corollary \ref{cor-nondegen}, we may w.l.o.g. assume $\la G_0\ra=G$ and that every $g\in G_0$ is contained in an atom. We may embed $\Z^d\leq \R^d$. Let $\pi:G\rightarrow \Z^d$ be the projection homomorphism with kernel $G_T$, and let $m=\exp(G_T)$.
As argued in Proposition \ref{prop-fg-diamondZ}, we have  $\C(\pi(G_0))=\R^d$ and  $0\notin \C^*(\pi(G_0)^\diamond)$ with  $\pi(G_0)$ finitary. Note  $\pi(G_0^\diamond)=\pi(G_0)^\diamond$ by definition of $G_0^\diamond$. Thus Proposition \ref{prop-rational-atoms} and $0\notin \C^*(\pi(G_0)^\diamond)$ together imply
$\mathcal A(\pi(G_0^\diamond))=\mathcal A (\pi(G_0)^\diamond)=\emptyset$, in turn implying
 $\mathcal A(G_0^\diamond)=\emptyset$. By Corollary \ref{cor-structural-Lambert}, $\pi(G_0)\setminus \pi(G_0)^\diamond\subseteq \pi(G_0)$ is a Lambert subset, say with bound $N\geq 1$.
By definition of $G_0^\diamond$, we have $\pi(G_0\setminus G_0^\diamond)=\pi(G_0)\setminus \pi(G_0)^\diamond$. Let $U\in \mathcal A(G_0)$ be an arbitrary atom.
Proposition \ref{prop-fg-atoms-basic-transfer}.1 implies that $\pi(U)=W_1\bdot\ldots\bdot W_\ell$ for some atoms $W_1,\ldots,W_\ell\in \mathcal A(\pi(G_0))$ with $\ell\leq \mathsf D(G_T)$. Thus $$\vp_{G_0\setminus G_0^\diamond}(U)\leq \Sum{i=1}{\ell}\vp_{\pi(G_0\setminus G_0^\diamond)}(W_i)=\Sum{i=1}{\ell}\vp_{\pi(G_0)\setminus \pi(G_0)^\diamond}(W_i)\leq \ell N\leq \mathsf D(G_T)N,$$ with the second inequality in view of $\pi(G_0)\setminus \pi(G_0)^\diamond\subseteq \pi(G_0)$ being  a Lambert subset with bound $N\geq 1$, which shows that $G_0\setminus G_0^\diamond\subseteq G_0$ is a Lambert subset with bound $\mathsf D(G_T)N<\infty$.
\end{proof}

The following extends Proposition \ref{prop-finitary-FiniteDeletion} to our current setting.

\begin{proposition}
\label{prop-tor-deletion}
Let $G= \Z^d\oplus G_T$ be a finitely generated abelian group with torsion subgroup  $G_T\leq G$, where $d\geq 0$, and  let $G_0\subseteq G$ be a subset. Suppose   $\rho_{(d+1)\exp(G_T)}(G_0)<\infty$. Then there are finite subsets $X\subseteq G_0^\diamond$ and $Y\subseteq G_0$ such that $\mathcal A(G_0\setminus X)$ is finite and $\mathcal A(G_0\setminus Y)=\emptyset$.
\end{proposition}

\begin{proof}
In view of \eqref{rho-ascend-chain} and Corollary \ref{cor-nondegen}, we may w.l.o.g. assume $\la G_0\ra=G$ and that every $g\in G_0$ is contained in an atom. We may embed $\Z^d\leq \R^d$. Let $\pi:G\rightarrow \Z^d$ be the projection homomorphism with kernel $G_T$. 
As argued in Proposition \ref{prop-fg-diamondZ}, we have  $\C(\pi(G_0))=\R^d$ and  $0\notin \C^*(\pi(G_0)^\diamond)$ with  $\pi(G_0)$  finitary. Thus Proposition \ref{prop-finitary-FiniteDeletion} implies that there is a finite subset $\wtilde X\subseteq \pi(G_0)^\diamond$ such that $\mathcal A(\pi(G_0)\setminus \wtilde X)$ is finite. Hence, since $\ker \pi=G_T$ is finite, it follows that $X:=\pi^{-1}(\wtilde X)\cap G_0\subseteq G_0$ is finite with $\pi(X)=\wtilde X\subseteq \pi(G_0)^\diamond$, ensuring $X\subseteq G_0^\diamond$ by definition of $G_0^\diamond$. It follows that $\mathcal A(\pi(G_0)\setminus \wtilde X)=\mathcal A(\pi(G_0\setminus X))$ is finite.

Consider an arbitrary atom $U\in\mathcal A(G_0\setminus X)$. Then Proposition \ref{prop-fg-atoms-basic-transfer}.1 ensures that $\pi(U)=W_1\bdot\ldots\bdot W_\ell$ for some atoms $W_1,\ldots,W_\ell\in \mathcal A(\pi(G_0\setminus X))$ with $\ell\leq \mathsf D(G_T)$. Since $\mathcal A(\pi(G_0\setminus X))$ is finite, there are at most $|\mathcal A(\pi(G_0\setminus X))|\cdot\ell\leq |\mathcal A(\pi(G_0\setminus X))|\cdot \mathsf D(G_T)<\infty$ possibilities for $\pi(U)$, for our arbitrary atom $U\in \mathcal A(G_0\setminus X)$. Thus, since $\ker \pi=G_T$ is finite, it follows that there are only finitely many possibilities for $U\in \mathcal A(G_0\setminus X)$, meaning  $\mathcal A(G_0\setminus X)$ is finite. Including in $X$ one element from each of the finite number of atoms from $\mathcal A(G_0\setminus X)$ then yields a subset $Y\subseteq G_0$ with $\mathcal A(G_0\setminus Y)=\emptyset$, completing the proof.
\end{proof}

We conclude with the extension of Theorem \ref{thm-rho-char}.

\begin{mtheorem}\label{thm-rho-char-tor}
Let $H$ be a Transfer Krull Monoid over a subset $G_0$ of a finitely generated abelian group $G=\Z^d\oplus G_T$ with torsion subgroup $G_T\leq G$, where $d\geq 0$. Then the following are equivalent.
\begin{itemize}
\item[1.] $\rho(H)<\infty$.
\item[2.] $\rho_k(H)<\infty$ for all $k\geq 1$.
\item[3.] $\rho_{(d+1)\exp(G_T)}(H)<\infty$.
\item[4.] There exists a subset $X\subseteq G_0$ such that $\mathcal A(X)=\emptyset$ and $G_0\setminus X\subseteq G_0$ is a Lambert subset.
\item[5.] $\mathcal A(G_0^\diamond)=\emptyset$.
\item[6.] $0\notin \C^*(\pi(G_0)^\diamond)$, where $\pi:G\rightarrow \Z^d\leq \R^d$ is the projection  with kernel $G_T$.
\end{itemize}
\end{mtheorem}

\begin{proof}Since $H$ is a Transfer Krull Monoid, we have $\rho(H)=\rho(G_0)$ and $\rho_k(H)=\rho_k(G_0)$ for all $k\geq 1$. It thus suffices to prove the theorem when $H=\mathcal B(G_0)$, which we now assume.
The implications $4.\Rightarrow 1.\Rightarrow 2.$ follows from proposition \ref{prop-pre-rho-char}, while the implication $2.\Rightarrow 3.$ is trivial.
The implication $3.\Rightarrow 5.$ follows by Proposition \ref{prop-lambert-tor-a}.
 In view of Proposition \ref{prop-rational-atoms}, Item 6 is equivalent to $\mathcal A(\pi(G_0)^\diamond)=\emptyset$. By definition of $G_0^\diamond$, we have $\pi(G_0^\diamond)= \pi(G_0)^\diamond$. Since $G_T$ has finite exponent, it follows that $\mathcal A(G_0^\diamond)=\emptyset$ if and only if $\mathcal A(\pi(G_0)^\diamond)=\mathcal A(\pi(G_0^\diamond))=\emptyset$. Thus Items 5 and 6 are equivalent.
It remains to establish the implication
$6.\Rightarrow 4.$ Suppose $0\notin \C^*(\pi(G_0)^\diamond)$. Then Theorem \ref{thm-rho-char} implies that $\rho_{d+1}(\pi(G_0))<\infty$, whence Proposition \ref{prop-fg-rho-transfer}.1 implies that $\rho_{(d+1)\exp(G_T)}(G_0)<\infty$, allowing us to apply  Proposition \ref{prop-lambert-tor-a}, which shows that Item 4 holds with $X=G_0^\diamond$.
\end{proof}

\subsection*{Summary}
We can now summarize our results regarding what we have shown  under an assumption of finite elasticities. Let $H$ be a Transfer Krull Monoid over a subset $G_0$ of a finitely generated abelian group $G= \Z^d\oplus G_T$  with torsion subgroup  $G_T\leq G$, where $d\geq 0$. We may embed $\Z^d\leq \R^d$. For example, $H=\mathcal B(G_0)$. Let $\pi:G\rightarrow \Z^d$ be the projection homomorphism with kernel $G_T$. There is little loss of generality to assume every element $g\in G_0$ occurs in some atom (else we can pass to the subset of $G_0$ having this property) and that $\la G_0\ra=G$ (else we can replace $G$ with $\la G_0\ra$), which we now do.
Note, since $G_T$ has finite exponent, that $g\in G_0$ is in an atom from $\mathcal A(G_0)$ if and only if $\pi(g)\in \pi(G_0)$ is contained in an atom from $\mathcal A(\pi(G_0))$. Thus  Proposition \ref{prop-notrivialG_0} implies that $\C(\pi(G_0))=\R\la G_0\ra=\R^d$.  Under these assumptions, we now summarize some of the key results.
\begin{itemize}
\item[1.] There exists a minimal $s\in [1,(d+1)\exp(G_T)]$ such that $\rho_s(H)<\infty$ implies $\rho_k(H)<\infty$ for all $k\geq 1$. (Corollary \ref{cor-fg-rho-transfer})
\end{itemize}
 As for the torsion-free case, it would be interesting to know if the estimate $s\leq (d+1)\exp(G_T)$ is tight or can be improved.
  Main Theorem \ref{thm-rho-char-tor} characterizes when Item 1 occurs either in terms of  a basic combinatorial property of the atoms $\mathcal A(G_0)$, or the geometric property $0\notin \C^*(\pi(G_0)^\diamond)$. It also shows that
   $$\rho_{(d+1)\exp(G_T)}(H)<\infty\mbox{ is equivalent to } \rho(H)<\infty.$$
Assuming additionally that $\rho_{(d+1)\exp(G_T)}(H)<\infty$, so that the conclusion of  Item 1 holds, we obtain the following properties.
\begin{itemize}
\item[2.] $0\notin\C(\pi(G_0)^\diamond)$. In particular, $\pi(G_0)\subseteq \R^d$ is finitary and  all the results of Section \ref{sec-finitary} are available for studying the set $\pi(G_0)$, including Theorems \ref{thm-finitary-diamond-containment}, \ref{thm-finitary-FiniteProps-I}, \ref{thm-finitary-FiniteProps-II} and \ref{thm-finitary-FiniteProps-III}. (Main Theorem \ref{thm-rho-char-tor} and Theorem \ref{thm-keylemmaII})

\item[3.] Besides the  defining definition of the subset $G_0^\diamond\subseteq G_0$ (stated earlier  in Section \ref{sec-transfer-monoid-results} combined with Proposition \ref{prop-G_0diamond-1st-easy-equiv}), we also have (by Proposition \ref{prop-fg-diamondZ})
\ber\nn G_0^\diamond&=&\{g\in G_0:\; \sup\{\vp_{g}(U):\; U\in \mathcal A^{\mathsf{elm}}(G_0)\}=\infty\}\\\nn &=&\{g\in G_0:\; \sup\{\vp_{g}(U):\; U\in \mathcal A(G_0)\}=\infty\}.\eer
  \item[4.] There is a finite subset $X\subseteq G_0^\diamond$ such that $\mathcal A(G_0\setminus X)$ is finite, and thus also a finite subset $Y\subseteq G_0$ such that $\mathcal A(G_0\setminus Y)=\emptyset$. (Proposition \ref{prop-tor-deletion}).
  \item[5.] $G_0\setminus G_0^\diamond \subseteq G_0$ is a Lambert subset with $\mathcal A(G_0^\diamond)=\emptyset$, indeed, $G_0\setminus G_0^\diamond \subseteq G_0$ is the unique maximal Lambert subset. (Propositions \ref{prop-lambert-tor-a} and \ref{prop-fg-diamondZ})
  \item[6.] If $H$ is a  Krull Monoid, the Weak Tame Degree (as defined in Main Theorem \ref{thm-tame-pseudo-tor}) is finite: $\mathsf t_w(H)<\infty$. In particular,  $N:=\max\{2,\mathsf t_w(G_0)\}<\infty$ is finite.
  \item[7.] The elasticities of $H$ do not contain arbitrarily large gaps: $\rho_k(H)-\rho_{k-1}(H)\leq N<\infty$ for all $k\geq 2$. In particular, $\rho_k(H)\leq Nk$ grows linearly. (Main Theorem \ref{thm-Delta-finite-tor}.1)
  \item[8.] $\max \Delta(H)\leq \rho_N(H)-N<\infty$ is finite. (Main Theorem \ref{thm-Delta-finite-tor}.2)
  \item[9.] The Set of Distances is  finite: $|\Delta(H)|<\infty$. (Main Theorem \ref{thm-Delta-finite-tor}.2)
  \item[10.] If $H$ is a Krull Monoid, the catenary degree is finite: $\mathsf c(H)\leq \rho_N(H)<\infty$. (Main Theorem \ref{thm-Delta-finite-tor}.4)
  \item[11.] If $H$ is a Krull Monoid, then $H$ is locally tame. (Main Theorem \ref{thm-Delta-finite-tor}.5)
  \item[12.] The Structure Theorem for Unions  holds for $H$. (Main Theorem \ref{thm-Delta-finite-tor}.3)
\end{itemize}

We remark that there are a few even stronger regularity properties of factorization that are \emph{not} implied by the finiteness of the elasticities, including the finiteness of  the monotone catenary degree  \cite[Section 7]{Gerold-lambert-rankone}, the finiteness of the (global) tame degree \cite[Theorem 4.2]{Ge-Ka10a}, and that  the Structure Theorem for Sets of Lengths holds  \cite[Section 6]{Gerold-lambert-rankone}.

\end{document}